\documentclass[12pt]{amsbook}

\usepackage{latexsym}
\usepackage[all]{xy}
\usepackage{amsmath, amsthm }
\usepackage{amssymb}
\usepackage{amsfonts}
\usepackage{color}
\usepackage{mathtools}
\usepackage{stmaryrd}
\usepackage{cancel}
\usepackage{enumitem}

\usepackage{nameref}
\usepackage{hyperref}
\usepackage{xcolor}
\hypersetup{colorlinks=true,urlcolor=blue,linkcolor=blue,citecolor=blue}

\usepackage{cite}
\usepackage{varioref}

\numberwithin{section}{chapter}
\numberwithin{subsection}{section}

\newcommand {\Top} {\mathbf{Top}}
\newcommand {\FMP} {\mathbf{FMP}}
\newcommand {\Perf} {\mathbf{Perf}}
\newcommand {\APerf} {\mathbf{APerf}}
\newcommand {\QCoh} {\mathbf{QCoh}}
\newcommand {\Coh} {\mathbf{Coh}}
\newcommand {\Vect} {\mathbf{Vect}}

\newcommand{\Gl}{\mathbf{Gl}}
\newcommand{\CAlg}{\mathbf{CAlg}}

\newcommand{\VV}{\mathbb{V}}
\newcommand{\CC}{\mathbb{C}}

\newcommand{\QQ}{\mathbb{Q}}
\newcommand{\ZZ}{\mathbb{Z}}
\newcommand{\NN}{\mathbb{N}}
\newcommand{\RR}{\mathbb{R}}
\newcommand{\LL}{\mathbb{L}}
\newcommand {\TT} {\mathbb{T}}
\newcommand{\Gm}{\mathbb{G}_{m}}
\newcommand{\Ga}{\mathbb{G}_{a}}
\newcommand {\Cs} {\mathcal{C}^\infty}
\newcommand {\C} {\mathcal{C}}
\newcommand {\F} {\mathcal{F}}
\newcommand {\G} {\mathcal{G}}
\newcommand {\A} {\mathcal{A}}
\newcommand {\M} {\mathcal{M}}
\newcommand{\D}{\mathcal{D}}
\newcommand{\B}{\mathcal{B}}
\newcommand{\N}{\mathcal{N}}
\newcommand{\cH}{\mathcal{H}}
\newcommand{\cA}{\mathcal{A}}
\newcommand{\cN}{\mathcal{N}}
\newcommand{\cM}{\mathcal{M}}
\newcommand{\cL}{\mathcal{L}}
\newcommand{\ccL}{\widehat{\mathcal{L}}}
\newcommand{\cJ}{\mathcal{J}}
\newcommand{\cT}{\mathcal{T}}
\newcommand{\cQ}{\mathcal{Q}}

\newcommand {\Map} {\mathbf{Map}}
\newcommand {\uMap} {\underline{\mathbf{Map}}}
\newcommand {\Parf} {\mathsf{Perf}}
\newcommand {\rh} {\mathbb{R}\underline{Hom}}
\newcommand {\rch} {\mathbb{R}\underline{\mathcal{H}om}}
\newcommand {\OO} {\mathcal{O}}
\newcommand {\DR}{\mathbf{DR}}
\newcommand {\Fol}{\mathcal{F}ol}

\newcommand  {\CDR}     {\widehat{C}_{DR}}

\newcommand {\Spec} {\mathbf{Spec}}
\newcommand {\Spf} {\mathsf{Spf}}
\newcommand  {\dgart}     {\mathbf{dgart}}
\newcommand  {\dg}     {\mathbf{dg}}
\newcommand  {\fdg}     {\mathbf{dg}^{fil}}
\newcommand  {\cfdg}     {\widehat{\mathbf{dg}}^{fil}}
\newcommand  {\grdg}     {\mathbf{dg}^{gr}}
\newcommand  {\grdga}     {\mathbf{dga}^{gr}}
\newcommand  {\edg}     {\epsilon-\mathbf{dg}}
\newcommand  {\egrdg}     {\epsilon-\mathbf{dg}^{gr}}
\newcommand {\scat} {\mathbf{Cat}_{\infty}}
\newcommand  {\egrcdga}     {\epsilon-\mathbf{cdga}^{gr}}
\newcommand  {\grcdga}     {\mathbf{cdga}^{gr}}

\newcommand  {\cdgamod}     {\mathbf{cdgamod}}
\newcommand  {\cdga}     {\mathbf{cdga}}
\newcommand  {\dglie}     {\mathbf{dgLie}}

\newcommand  {\grdglie}     {\mathbf{dgLie}^{gr}}

\newcommand  {\sCR}     {\mathbf{sCR}}
\newcommand  {\CR}     {\mathbf{CR}}

\newcommand  {\St}     {\mathbf{St}}
\newcommand  {\Aff}     {\mathbf{Aff}}

\newcommand  {\dAff}     {\mathbf{dAff}}

\newcommand  {\ncdga}     {\mathbf{cdga}^{\leq 0}}

\newcommand  {\holim}   {\mathrm{Holim}}

\newcommand  {\llim}   {\mathrm{Lim}}
\newcommand  {\colim}   {\mathrm{Colim}}
\newcommand  {\dSt}   {\mathbf{dSt}}

\newcommand{\s}{\infty}
\newcommand{\HH}{\mathbb{H}}

\newcommand{\Alg}{\mathbf{Alg}}
\newcommand{\MC}{\mathcal{MC}}

\newcommand{\csMod}{\mathsf{csMod}}

\newcommand{\Cinf}{\widehat{C}_{\mathrm{inf}}}

\newcommand {\T} {\mathbb{T}}

\newcommand{\WW}{\mathbb{W}}

\newcommand  {\csCR}     {\mathbf{csCR}}

\newcommand  {\medg}     {\epsilon-\mathbf{dg}^{gr}}

\newcommand  {\mecdga}     {\epsilon-\mathbf{cdga}^{gr}}

\newcommand  {\dChAff}     {\mathbf{dChAff}}

\newcommand{\uHom}{\underline{Hom}}

\theoremstyle{plain}
\newtheorem{thm}{Theorem}[subsection]
\newtheorem{df}[thm]{Definition}
\newtheorem{prop}[thm]{Propositon}
\newtheorem{rmk}[thm]{Remark}
\newtheorem{cor}[thm]{Corollary}
\newtheorem{ex}[thm]{Example}
\newtheorem{lem}[thm]{Lemma}
\newtheorem{slem}[thm]{Sub-Lemma}

\makeindex

\title{Derived Foliations}
\author{Bertrand To\"en and Gabriele Vezzosi}
\date{August 2025}

\begin{document}

\maketitle

\tableofcontents

\section*{Introduction}

One important feature of derived algebraic geometry (in the sense of \cite{hagII,toenems,pv})
is the observation that singularities of varieties and schemes encountered in real life 
are most of the time tamed by natural derived structures. This is particularly 
visible in moduli 
theory. In many instances, classical/underived moduli spaces are very singular, 
but their derived counterparts typically have finite dimensional 
cotangent complexes, and in the best scenarios, local complete intersection singularities. 
As a result, derived moduli spaces carry important information of geometric origin which 
can not be seen directly on their classical truncations. Typical examples are: 
virtual fundamental classes 
(see \cite{MR3341464, khan2019virtualfundamentalclassesderived}), 
and shifted symplectic structures (see \cite{ptvv}).

The main subject of the present book is guided by a similar principle: foliations encountered in real 
life often come equipped with natural \emph{derived structures} which tame their singularities.
More precisely, there exists a natural notion of \emph{derived foliation}, with a natural
truncation functor to classical, possibly singular, foliations, and this notion helps when dealing with 
foliation singularities in general. The purpose of this book is to make precise the notion
of derived foliation, and to prove results that support this principle. 

\subsection*{From classical foliations to graded mixed algebras}

Before presenting the main definition of this book, we would like to contemplate the
classical situation. For a smooth algebraic complex variety $X$, a (algebraic) foliation $\F$  on $X$
is usually defined as a sub-bundle\footnote{We warn the reader here that all along this book
the expression \emph{sub-bundle} of a vector bundle $V$ means a sub-coherent sheaf $W \subset V$
such that both coherent sheaves $W$ and $V/W$ are in fact vector bundles.} $\T_\F \subset \T_X$ of the tangent bundle $\T_X$, which 
is furthermore stable by the Lie bracket. To avoid possible 
confusions, we will say that such a $\T_\F \subset \T_X$ is a \emph{classical smooth foliation on $X$}.
By duality, the sub-bundle $\T_\F$ defines a quotient bundle
$$ p :\Omega_X^1 \twoheadrightarrow \Omega^1_\F = \T_\F^{\vee}$$
where $\Omega_X^1$ is the sheaf of Kähler differential forms on $X$. Using the following standard formula 
relating the Lie bracket $[-,-]$ of vector fields and the de Rham differential $dR$ 
(for $\omega$ a $1$-form, and $u$ and $v$ vector fields on $X$)
$$dR(\omega)(u,v) = u(\omega(v)) - v(\omega(u)) + \omega([u,v])$$
we see that the kernel 
$K_\F \subset \Omega^1_X$ of the above quotient map $p$ is a differential ideal: the (sheaf) image 
of $dR : K_\F \to \Omega^2_X$ is contained in $K_\F \wedge \Omega_X^1$, where $K_\F \wedge \Omega_X^1$ denotes  the (sheaf) image of the composition 
$$\xymatrix{K_\F \otimes \Omega_X^1
\ar[r] & \Omega^1 \otimes \Omega_X^1 \ar[r]^-{\wedge} &  \Omega^2_X.}$$
Conversely, any differential ideal $K_\F$, i.e. a sub-bundle with the property that 
$dR(K_\F) \subset K_\F \wedge \Omega^1_X \subset \Omega_X^2$, is of the form 
$\T_\F^\vee$ for a sub-bundle $T_\F \subset \T_X$ stable by the Lie bracket. More importantly for our purposes, the differential ideal 
property of $K_\F \subset \Omega_X^1$ is also equivalent to the property that the de Rham 
differential $dR$ on the graded ring of forms $\Omega^*_X$, descends to a differential $dR_\F$
on the quotient graded ring
$$\Omega_\F^* = (\Omega_X^*/K_\F.\Omega^*_X).$$
It is easy to check the association $\F \mapsto (K_\F, dR_\F)$ provides a one-to-one
correspondence between classical smooth foliations $\F$ on $X$ and 
sub-bundles $K_\F \subset \Omega^1_X$ together with a multiplicative 
differential $dR_\F$ on $(\Omega_X^*/K_\F.\Omega^*_X)$ compatible with the de
Rham differential via the quotient map $\Omega^*_X \twoheadrightarrow \Omega^*_\F$.

We now want to consider the sheaf $\Omega^1_X$ as sitting in cohomological degree $-1$, 
and we thus consider $\Omega_X^1[*] = \oplus_p \Omega^p_X[p]$. This unusual 
convention, as opposed to $\Omega^1_X$ being considered in cohomological degrees $1$, 
is explained by the geometric interpretation of forms in terms of \emph{derived loop spaces}
and for which we refer to our sections \S \ref{subsec:Derivedfoliationsasequivariantderivedlinearstacks} 
and  \S \ref{derfolgradedloop}. The sheaf $\Omega_\F^*[*]$ is now a sheaf of 
bi-graded commutative rings, were the bi-degree is given on the one side by 
the degree of differential forms, and on the other side by the cohomological degree.
The differential $dR_\F$ on $\Omega_\F^*[*]$ is now 
a multiplicative differential of bi-degree $(1,-1)$: it increases the degree of forms 
by $1$ but decreases the cohomological degree by $1$. Moreover, with these conventions we have
$$\Omega_\F^*[*]\simeq Sym_{\OO_X}(\Omega^1_X[1]),$$
making the $\Omega_\F^*[*]$ into the free bi-graded commutative ring over $\Omega_\F^1[1]$.
Again, the association $\F \mapsto (\Omega_\F^*[*],dR_\F)$ can be seen to define a one-to-one
correspondence between classical smooth foliations $\F$ on $X$ and 
quotient bi-graded commutative rings $\Omega_X^*[*] \twoheadrightarrow \Omega_\F^*[*]$ 
of the form $Sym_{\OO_X}(\Omega^1_\F[1])$ such that the de Rham differential
descends to a compatible differential on $\Omega_\F^*[*]$.

This somehow complicated alternative point of view on classical smooth foliations on $X$ 
has the advantage of introducing (and motivating) one of the major actor of this book: \emph{graded mixed
commutative differential graded algebras}, or in short \emph{graded mixed cdga's}.
These are bi-graded commutative $\CC$-algebras $A^*(*)$, together with 
two differentials
$$d : A^*(*) \to A^{*+1}(*) \qquad \epsilon : A^*(*) \to A^{*-1}(*+1)$$
satisfying natural compatibilities (see Definition \ref{dI-3}). 
In the example above arising from a smooth classical foliation on $X$, $\Omega_\F^*[*]$ is a graded mixed cdga
with $d=0$ and $\epsilon = dR_\F$. Clearly, smooth classical 
foliations on $X$ can be understood as quotients of $\Omega_X^*[*]$
in the category of graded mixed cdga's. This observation, which is not very helpful
when dealing only with classical smooth foliations, is central in the definition of 
derived foliations in general.

\subsection*{The general notion of derived foliations}

We turn to a more general situation where $X$ is an affine derived scheme, say over some
base field of characteristic zero $k$. Such a derived scheme $X=\Spec\, A$ 
is the spectrum of a connective $k$-linear cdga $A$. The cdga $A$ possesses a cotangent
complex $\LL_{A/k}$ (relative to $k$) which comes equipped with a universal
derivation $dR : A \to \LL_{A/k}$ (see Appendix \S \ref{secapp:Derivedschemesandstacks}).
This derivation can be extended uniquely to 
an extra differential $\epsilon$ on the graded cdga $Sym_A(\LL_{A/k}[1])$, making
it into a graded mixed cdga. This graded mixed cdga is denoted by $\DR(X/k)$ or $\DR(A/k)$
and is called the \emph{de Rham graded mixed cdga} associated to $X$. The underlying
bi-graded alegbra is $Sym_A(\LL_{A/k}[1])$ where the bi-degrees are the cohomological
degrees on one hand, and the natural power degree in the symmetric algebra
on the other. The two differentials $d$ and $\epsilon$ are respectively induced by  
the natural cohomological differential of $A$ and by the universal de Rham
derivation $dR : A\to \LL_{A/k}$.

We can now state the major definition of this book
(see Definition \ref{dII-1}).\\

\noindent \textbf{Definition}. A \emph{derived foliation $\F$ on $X$
(relative to $k$)} consists of a graded mixed cdga $\DR(\F)$ together with a morphism 
of graded mixed cdga's $\DR(X/k)\to \DR(\F)$, such that the underlying graded cdga of  $\DR(\F)$ 
is quasi-isomorphic to $Sym_A(\LL_\F[1])$, for some $A$-dg-module $\LL_\F$.\\

In the definition above the $A$-dg-module $\LL_\F$ is called the \emph{cotangent complex} 
of the derived foliation $\F$; it comes together with an anchor
map $\LL_{X/k} \to \LL_\F$ and plays the role of the quotient 
$\Omega^1_X \twoheadrightarrow \Omega_\F^1$ previously discussed 
in the setting of classical smooth foliations. The mixed structure on $\DR(\F)$ together with the quasi-isomorphism $Sym_A(\LL_\F[1]) \simeq \ DR(\F)$ of \emph{graded} cdga's, 
 define on $Sym_A(\LL_\F[1])$ a de Rham-like differential
$$\epsilon : \Lambda^i_A \LL_{\F} \to  \Lambda^{i+1}_A \LL_{\F}$$
which is well defined as a morphism in the $\s$-category of $k$-dg-modules.
This de Rham differential carries (an infinite number of) homotopy coherent structures
expressing that it squares to $0$, that is multiplicative and compatible 
with the de Rham differential on $\DR(X/k)$. These homotopy coherent structures
are encapsulated and hidden in the full graded mixed cdga structure $\DR(\F)$
and are implicit in our approach, as the de Rham differential
$\epsilon$ exists by definition on $\DR(\F)$ (which is only quasi-isomophic
to $Sym_A(\LL_\F[1])$), and does not exist strictly speaking 
on $Sym_A(\LL_\F[1])$ itself. We refer to Corollary \ref{cpI-5} for more on the
shape of homotopy coherences involved, but the point of our definition of derived foliations,
and more generally the point of view of this book is to avoid, as much as possible, to work
with explicit homotopy coherences and try to keep them implicit in the picture.

This being said, the derived foliations on $X$ as defined above 
form an $\s$-category $\Fol(X/k)$. The assignment $X \mapsto \Fol(X/k)$ is functorial in $X$ in the sense
that a derived foliation $\F$ on $X$ can be pulled back along any morphism
$f : Y \to X$ of derived affine schemes. This pull-back operation is of fundamental importance, as even when 
$Y$ and $X$ are smooth schemes, and $\F$ is a classical smooth foliation on $X$, 
its pull-back $f^*(\F)$ is a derived foliation on $Y$ which is, in general,
not a classical foliation anymore. As a result, pull-backs offer a way to 
\emph{create} interesting examples of derived foliations starting from a classical 
situation. This is, of course, parallel to the fact that the (derived) fiber product of
smooth schemes usually is a non-trivial derived scheme. 

The functoriality $X \mapsto \Fol(X/k)$, on derived affine schemes $X$, can be used in order to provide 
a definition of derived foliation on any \emph{derived stack} $X$, simply by imposing the usual 
descent condition 
$$\Fol(X/k) := \lim_{\Spec\, A \to X}\Fol(\Spec\, A/k).$$
However, the behaviour of derived foliations on general derived stacks brings some
difficulties, as for instance cotangent complexes of derived foliations are not stable by pull-backs. As a result the existence of well defined global cotangent complex
of a derived foliation over a general derived stack is not obvious. In that direction 
we prove several descent results showing that, in spite of that, cotangent complexes still have 
a reasonable meaning when working with derived Artin stacks.

Finally, a given derived foliation $\F$ has an attached \emph{de Rham cohomology}, that 
is the generalisation of the well known de Rham cohomology along the leaves, or
\emph{foliated de Rham cohomology}. Within our approach, this foliated de Rham cohomology is very easy and natural 
to define by means of the \emph{Tate-realization functor}, relating 
graded mixed cdga's and complete filtered cdga's (see Corollary \ref{clI-4}). The image
of $\DR(\F)$ by this realization, $\CDR^*(\F)$, is by definition the commutative dg-algebra
of \emph{de Rham cohomology along} $\F$, and the filtration is here of course
the natural Hodge filtration. This de Rham cohomology is however only 
a piece of a refined cohomological invariant called the \emph{Hodge filtration}
associated to a derived foliation $\F$. The terminology is here unfortunate, 
as it should not be confused with the filtration on $\CDR^*(\F)$. The Hodge filtration 
associated with $\F$ is indeed a filtration on the de Rham cohomology of $X$
(suitably derived and completed) $\CDR^*(X)$, whose first graded piece recovers
$\CDR^*(\F)$. More generally, the higher graded pieces can also be described 
as foliated de Rham cohomology with coefficients in the
powers of the normal complex (see \S\,\ref{sec:Hodgefiltration}). The existence
of this Hodge filtration seems to us more important than the coarser notion of
foliated de Rham cohomology, and it will turn out to be very helpful in the study of 
characteristic classes associated to derived foliations (see Corollary 
\ref{cor-chernvanishing}).
  
\subsection*{Formal aspects of derived foliations}

One specific and important aspect of derived foliations is their
\emph{formal integrability} property. This is a major difference with the classical
situation, where certain singular foliations might not be formally integrable
(see Remark \ref{rdII-14}). It can be shown, under very mild assumptions, that
any derived foliation $\F$ on a derived Artin stack $X$ is formally integrable 
around each point in $X$ (see Corollary \ref{ctII-2}). In particular, this implies
that for a field valued point $x\in X(K)$, there is a well defined 
notion of \emph{formal leaf passing through $x$}, as well as a \emph{formal leaf
space at $x$}. As a result, 
it may happen that a given classical singular foliation cannot be enhanced
to a derived foliation (i.e. it is not the truncation of a derived foliation). This can be seen as either a positive result, 
or a negative result, depending of the situation. Nevertheless, this shows that the
situation diverges from what is happening in derived algebraic geometry, where
any scheme can be endowed with the trivial derived structure.

The precise connection between classical singular foliations and derived foliations
can be made clearer. For a smooth scheme $X$ we have a notion of 
\emph{quasi-smooth and rigid} derived foliation on $X$. This is an important class
of derived foliations as they can be considered as the derived foliations which are the
closest to be classical. These are given by $\F\in \Fol(X/k)$ for which the
cotangent complex sits in an exact triangle of perfect complexes on $X$
$$\xymatrix{
\N_\F^* \ar[r] & \Omega_X^1 \ar[r] & \LL_\F
}$$
where $\N_\F^*$ is a vector bundle on $X$. In many situation $\N_\F^* \to \Omega^1_X$
is a monomorphism, and thus $\F$ is a classical smooth foliation on the 
open $U \subset X$ where $\LL_\F$ is a vector bundle. The closed complement of 
$U$ is by definition the \emph{singular set} of $\F$. The derived foliation $\F$
admits a classical truncation which is the image $\N_\F^*$ in $\Omega_X^1$
and which defines a differential ideal $K_\F$ on $X$. As we 
have already mentioned, not all differential ideals arise 
by this construction. In fact we can show that the differential ideals
obtained this way are the one for which an infinite \emph{Godbillon-Vey} sequence exists
(see \S\ref{sec:Existenceofderivedenhancements}). That of a Godbillon-Vey sequence
is a classical notion in foliation theory and therefore this result sheds some
light of how derived foliations are related to classicial singular foliations.

\subsection*{Derived foliations in the holomorphic setting}

Most of this book is written in the algebraic setting. However, 
except in some extremely specific situations, classical foliations \emph{cannot} be
integrated algebraically, and \emph{analytic} integration is in most situations the best
one can hope. It is therefore important to introduce 
derived foliations also in the holomorphic context. Our definition of
derived foliation as certain graded mixed cdga's can be easily transported
in the holomorphic setting, but with an extra technical component
in order to keep track of the holomorphic structure. Nevertheless 
derived foliations on a derived analytic Artin stack makes sense, and formally 
behave as in the algebraic setting. We develop the strict minimum 
of the holomorphic theory with two objectives in sight: 
the \emph{GAGA theorem}, and the existence of an \emph{analytic leaf space} under appropriate  assumptions.

The GAGA theorem for derived foliations states that for
a proper derived Deligne-Mumford stack $X$, the analytification functor induces
an equivalence between the $\s$-categories of algebraic (and almost perfect) 
derived foliations on $X$ and the $\s$-category of holomorphic derived foliations
on the derived analytic stack $X^h$ associated to $X$ (see Theorem \ref{GAGAperfectfol}). 
This is not very surprising as 
derived foliations are defined in terms of coherent data via the derived de Rham 
algebra of forms, and comparison results relating algebraic and holomorphic cotangent
complexes are well known. 

As in the algebraic situation, holomorphic derived foliations are always formally
integrable and it is thus natural to ask if they can be integrated locally 
in the analytic topology. For this, we can use classical results
from \cite{mal1,mal2} to provide interesting conditions ensuring that 
quasi-smooth and rigid holomorphic derived foliations are indeed locally analytically integrable.
This is the starting point for a construction of a global holomorphic leaf space
and for the global analytic integrability. Indeed, under certain natural 
assumption, it is possible to show that a quasi-smooth holomorphic 
derived foliation $\F$ on a smooth complex variety $X$ admits 
a holomorphic leaf space $\pi : X \to X\sslash\F$ (see Theorem \ref{tVI-1}).
Here $X\sslash \F$ is a complex (highly non-separated) orbifold, $\pi$ is a flat 
surjective holomorphic map, and $\F$ is realized as the 
derived foliation induced by $\pi$. The situation is thus very similar 
to the case of smooth classical foliation for which $X\sslash \F$ is
represented by the well known \emph{holonomy groupoid} of $\F$. The novelty here is
that $\pi$ is only flat and thus can have singular fibers, corresponding to
the possibly singular leaves of $\F$. 

It is interesting to combine together the GAGA theorem and the existence
of the analytic leaf space when we start with a derived foliation
$\F$ on a smooth and proper complex variety $X$. As a first application, 
we prove a \emph{Riemann-Hilbert correspondence} which provides an equivalence
between the category of algebraic vector bundles with integrable connections
along $\F$, and the category of relative local systems along the fibers
of $\pi$ (see Theorem \ref{tVI-2}). This is a generalization to the singular
setting of the Riemann-Hilbert correspondence appearing in \cite{del}. 
In another direction, we can exploit the notion of the \emph{holonomy groups}
deduced from the existence of the analytic leaf space, 
which are, by definition, the stabilizers of the orbifold $X\sslash \F$. 
Classically holonomy groups have been useful in stability results, such as the
famous \emph{Reeb stability} theorem (see \cite{MR55692}). 
In our setting, it is in fact possible to prove that $\F$
is globally algebraic integrable on $X$, if and only if $\F$ possesses at least one
compact leaf with finite holonomy groups (see Theorem \ref{tVI-3}). This is again
an extension to the singular setting of previously known results
for classical smooth foliations (see \cite{MR1913291}).

\subsection*{Crystals along derived foliations}

We have already mentioned that a given derived foliation has 
an associated foliated de Rham cohomology. This can be extended by the introduction of a category 
of coefficients for foliated de Rham cohomology called \emph{crystals along
a derived foliation $\F$}. Morally speaking these are 
quasi-coherent complexes endowed with partial flat connections along $\F$. 
Concretely, they can be described in terms of graded mixed $\DR(\F)$-dg-modules
(see Definition \ref{defqcohcrys}). They also have another description in terms
of modules over a sheaf of associative dg-algebras $\D_\F$, called the \emph{sheaf 
of differential operators along} $\F$ (see Theorem \ref{pdmod}). 
The advantage of the latter description is the existence of a natural filtration
on $\D_\F$ induced by the order of differential operators. This is useful 
to define the notion of \emph{good filtration} on a given crystal, analogue to the
well known notion of a good filtration on a classical $\D$-module. This notion 
can then be used in order to define the \emph{singular support}, as well
as the \emph{characteristic cycle} of a crystal endowed with a good filtration, which is
a $K$-theory class inside the total space of the foliated cotangent complex
of $\F$ (see Definition \ref{dsing1}). Even thoughtthe notion of good
filtration has some drawbacks for derived foliations with singularities (see \S\, 
\ref{sect-inde}), 
it is enough to establish a \emph{Grothendieck-Riemann-Roch} theorem and
\emph{index} theorem for crystals along derived foliations (see Theorem \ref{tgrr}
and Corollary \ref{cindex}).

\subsection*{Derived foliations over bases of arbitrary characteristics}

Classical foliations in positive characteristics is a rich subject 
(see for instance \cite{MR927978,MR1863738}),
and we have therefore decided to include a last chapter on derived foliations 
over arbitrary base ring. It turns out that there are 
two, non-equivalent, 
possible extensions of our notion of derived foliations outside of characteristic
zero, which roughly speaking correspond to the two existing notion 
of differential operators: \emph{crystalline differential operators} and \emph{Grothendieck differential
operators}. These are also related to two different extension of de Rham 
cohomology, namely \emph{crystalline cohomology} and \emph{infinitesimal cohomology}. 

The best way to understand the difference between these two possible definitions of derived foliations
in arbitrary characteristics is to adopt a geometric point of view on derived
foliations based on group actions on derived linear stacks (see Definition \ref{defgeoderfol}).
We introduce two group stacks $\cH$ and $\cH_\pi$. The first one, $\cH$ is the semi-direct
product of $\Gm$ with $B\Ga^\sharp$, where $\Ga^\sharp$ is the additive group
with divided powers. The second one, $\cH_\pi$ is the semi-direct product of $\Gm$ with 
$\Ga[-1]=\Omega_0\Ga$ the additive group in degree $-1$. Both of these groups
control graded mixed complexes, in the sense that the symmetric
monoidal $\s$-categories $\QCoh(B\cH)$ and $\QCoh(B\cH_\pi)$ are 
both equivalent to the symmetric monoidal $\s$-category of graded mixed 
complexes. However, the two groups
$\cH$ and $\cH_\pi$ carry different informations. A complex $E$ together
with an action of $\cH$ on the corresponding derived linear stack $\VV(E)$
is called a \emph{derived foliation} (see Definition \ref{dVII-4}). 
On the other hand, a complex $E$ together
with an action of $\cH_\pi$ on the corresponding derived linear stack $\VV(E)$
is called an \emph{infinitesimal derived foliation} (see Definition \ref{dinfol}).
Derived foliations are directly related to (derived and Hodge completed) 
de Rham cohomology, and thus to crystalline cohomology as well (see 
Proposition \ref{compareDRCrys}). On the other hand, infinitesimal derived foliations
are closely related to infinitesimal cohomology, suitably derived and Hodge completed,
if necessary (see Corollary \ref{c1}).

Both, derived foliations and infinitesimal derived foliations have 
a generic behaviour similar to derived foliations in characteristic zero: they
have pull-backs, can be defined over general derived stacks, have an associated 
foliated cohomology theory, etc. However, their difference can be seen in applications. 
Derived foliations, due to their relation with crystalline cohomology, are
useful to study characteristic classes. For instance they can be used to prove
a characteristic $p>0$ version of the vanishing theorem of \emph{Baum-Bott} (and its
generalization on the existence of \emph{residues}, see \cite{baumbott}). 
On the other hand infinitesimal derived foliations are always formally integrable, as
opposed to derived foliations in general (see Corollary \ref{cint}). This 
is a generalization of the fact that not all classical smooth foliations in characteristic
$p>0$ can be formally integrated, because of the well known obstruction 
associated to $p$-curvature (see e.g. \cite{katz, MR927978}). In a way, infinitesimal derived
foliations in characteristic $p>0$ should be thought as \emph{$p$-restricted
derived foliation}, even though we did not try to make this notion precise.

\subsection*{Relations to other works}

We finish this introduction by mentioning several relations with other works. \\

To our knowledge the notion of derived foliation first appears in a 
talk by Tony Pantev in 2014 (see \cite{tony}). He explained in particular
how derived foliations appear naturally in the setting of shifted symplectic and Poisson 
structures of \cite{ptvv,cptvv}, and in particular how they are useful 
to reconstruct shifted potentials in order to write a $n$-shifted symplectic derived stack
globally as the derived critical locus of a function of degree $n+1$. The present book
has pursued this direction, and owns a lot to the conversations among the authors
of \cite{ptvv,cptvv} on the subject. The reconstruction of potential 
is not part of this book as we have decided to focus on slightly different aspects
of derived foliation theory, but has been carried out in details recently in 
\cite[\S 2.4]{hennion2025gluinginvariantsdonaldsonthomastype}.

In a similar direction, an idea conveyed by Damien Calaque during the last decade 
states that shifted Poisson structures in the sense of \cite{cptvv}
can be considered as derived foliations with shifted symplectic structures
along the leaves (in short a derived symplectic foliation). This is the derived 
analogue of the fact that a classical Poisson structure provides a symplectic foliation. 
This has been realized recently in 
\cite{tomić2025shiftedlagrangianthickeningsshifted},
where it is shown that 
a shifted Poisson structure is always given by a Lagrangian thickening, which can be identified with 
the formal leaf space of the wanted symplectic derived foliation
(see also \cite{calaque2024shiftedcotangentbundlessymplectic}).

The formal integrability 
properties of derived foliations studied in this book has been generalized
recently in \cite{brantner2025formalintegrationderivedfoliations}, including 
the non-zero characteristic situations. Our results of 
\S \ref{ctII-2} and \S \ref{cint} are therefore already subsumed by the results of 
\cite{brantner2025formalintegrationderivedfoliations}. Moreover, 
the results of the above work also provide a very general construction of a
\emph{global formal leaf space} of a derived foliation in the form
of a formal thickening. These type of quotients were already mentioned in 
\cite{tony} in the characteristic zero case, but are not included in this book.

Many of the filtrations appearing in this book are obtained by means 
of the equivalence given by the Tate-realization functor, which 
identifies graded mixed complexes and complete filtered complexes (see 
\S \ref{clI-4}). Therefore, all the notions of this book could also have been 
written in terms of complete filtered objects, such as complete filtered
cdga or complete filtered $LSym$-algebras. We are aware that 
Ben Antieau has also developed a theory of derived foliation 
in the topological context based on this notion of filtered algebra, and that  
probably some of the definitions of this book also appear in his work.

Mauro Porta has been keeping telling us that 
graded mixed cdga which are not necessarly graded free (and thus
do not correspond, strictly speaking, to derived foliations) are nonetheless interesting. More recently, 
Philippe Eyssidieux has convinced the authors that 
there could be indeed an interesting derived version of the notion of 
\emph{exterior differential system} (EDS for short). By definition, these
derived EDS on a variety $X$ are precisely 
the general graded mixed cdga under the de Rham algebra $\DR(X)$. It is pretty clear
that some of the notions and results of this book extend to this
more general setting of derived EDS. However, we have not included any 
content in that direction.

\section*{Contents of individual chapters}

This book is divided in seven chapters and one appendix. Chapter $1$ is devoted
to the algebraic preliminaries necessary to approach the study of derived foliations. 
We have included sections on graded mixed objects, filtered objects
and the Tate-realization relating these two notions. These notions will be used throughout the book. The two last sections of this
chapter include one incarnation of graded mixed complexes as representations
of the group stack $B\cH_0$. This geometric interpretation will turn out to be useful, 
typically to construct natural graded mixed objects as push-forwards 
along a morphism to $B\cH_0$. The last section is devoted to a (partial)
description of the homotopy coherences encoded in a graded mixed structure. 
This is a technical result that will be used occasionally, and does not quite follow
the general spirit of the general approach of this book which is to make these coherences implicit as much as possible.

Chapter $2$ contains the main definitions and notions on the subject. We define the $\s$-category $\Fol(X/k)$ of 
derived foliations, first on an affine derived $k$-schemes $X$, then by descent on general
derived $k$-stacks $X$.  We provide examples, and compare derived foliations in particular with the theory 
of dg-Lie algebroids already existing in the literature. We present
also the two important universal foliations: the initial foliation, also 
called the zero foliation, and the final foliation, also called the de Rham 
or tautological foliation. The third section is devoted to the formal
properties of derived foliations. We have restricted ourselves to the behaviour
formally at a field valued point only, which is a non-trivial restriction. However, 
we have treated the general case of almost perfect derived foliations on general 
derived Artin stacks locally of finite presentation. We explain in particular, 
how derived foliations on the formal completion at a point are governed purely in terms
of dg-Lie algebras. This section is central in the sense that it allows to make sense of
the notion of \emph{formal leaf} and \emph{formal leaf} space at a given point, whereas
those notions are not totally obvious when working over general derived stacks. 
The subsequent section concerns formal integrability in the more restrictive context
of quasi-smooth and rigid derived foliations on smooth varieties. It contains
in particular a comparison with classical singular foliations, and 
provide interpretations of classical notions (such as Godbillon-Vey sequences and
transversal jet spaces) in terms of derived foliations. The last section
is taken from the work of \cite{alf}, and presents a useful 
constuction tool of derived foliations by direct images (or Weil restriction). 
The geometric interpretation of derived foliations which is behind this direct image construction will be useful 
to understand and motivate Chapter 7 about derived foliations in arbitrary characteristics.

In Chapter $3$ we introduce and study foliated de Rham cohomology. We start by 
the introduction of the categories of crystals along a derived foliations,
and their basic functorialities by pull-backs. They give rise to 
the notion of \emph{foliated de Rham cohomology with coefficients in a crystal}
that will be used quite often in the rest of the book. The third section of this chapter introduces
the Hodge filtration. This is a filtration on the de Rham cohomology of a derived
stack $X$ induced by the datum of a derived foliation $\F \in \Fol(X/k)$. It is obtained
by the geometric point of view on derived foliations explained in \S\, \ref{sec:Directimages}
together with a very general deformation to the normal bundle construction. The final section
of this chapter provides some elements of transversal geometry by constructing 
the formal jets transversal to a given derived foliation. We explain how these
jets are endowed with a crystal structure, and how they give rise to 
an algebraic notion of holonomy and monodromy by means of a Tannakian approach. These
include higher homotopical invariants captured in the form of a schematic
homotopy type in the sense of \cite{chaff}, and closely related to 
\emph{Galois-differential homotopy theory}.

Chapter $4$ is devoted to the study of operations on crystals. We have not included
a full $6$ operations formalism, which seems to us out of reach at the moment, but 
have studied direct images and direct images with proper support when they are defined. 
The first section introduces the sheaf of differential operators along a derived
foliation $\F$, which is a sheaf of filtered dg-algebras whose associated graded object  
is of the form $Sym_{\OO_X}(\TT_\F)$, the symmetric algebra over the tangent 
complex of the derived foliation $\F$. Modules over $\D_\F$ are equivalent to crystals
along $\F$, but the existence of the filtration on $\D_\F$ can be used in order
to define good filtrations on crystals and the corresponding notion of characteristic cycles.
We prove in particular that characteristic cycles are stable by push-forward, i.e. 
 a foliated version of the Grothendieck-Riemann-Roch formula.

The analytic theory of derived foliation is covered in chapter $5$. In its first section
we remind the basics of derived analytic geometry based on the notion of holomorphic
cdga. In the second section we introduce the \emph{holomorphic graded mixed cdga},
which can be used to defined holomorphic derived foliations along the same
lines as in the algebraic setting. The two last sections are devoted to the proof
of two main results. First of all we study the analytic integrability, locally 
in the analytic topology, of quasi-smooth and rigid derived holomorphic foliations.
This allow us to make comparisons with previously known integrability results
for holomorphic classicial singular foliations, notably the results from \cite{mal1,mal2}.
In the last section of the chapter we prove a GAGA theorem for derived foliations
on proper derived Deligne-Mumford stacks. 

Chapter $6$ contains the notion of analytic leaves for holomorphic derived foliations.
By restricting to quasi-smooth and rigid derived foliations, we provide 
simple assumptions ensuring the existence of a nice analytic leaf space whose construction
uses techniques arising in the study of \emph{étendues} from $SGA4$, and is a generalization 
of the construction of the holonomy groupoid of a classical smooth holomorphic foliation.
The possible singularities bring new features concerning the leaf space:  
one instance is the non-triviality of the \emph{relative
homotopy type} that encodes the variations of the homotopy types of the leaves. This relative homotopy type is not locally constant (even at the level of germs because of the existence
of vanishing cycles). In the last two sections we give two main applications
of the existence of an analytic leaf space and of our previously established GAGA theorem. The first application
is a Riemann-Hibert correspondence, and the second application is an algebraic integrability
result based on the classical Reeb stability theorem.

The last chapter of the book is chapter $7$. It contains some results on derived foliations in 
non-zero characteristic situations. Its first two sections introduce 
an important group stack $\cH$, that enables us to define derived foliations
as linear derived stacks endowed with an action of $\cH$. In the next
section we explain how this notion is related to de Rham cohomology, and introduce extensions
of the notion of crystals and of the Hodge filtration already encountered in characteristic $0$.
The fourth section is devoted to the construction of characteristic classes of derived
foliations over any base ring, and explain how these interact nicely with the Hodge filtration.
We then establish classical vanishing results for the characteristic classes.
In the subsequent section of this chapter, we use these results to prove a characteristic 
$p>0$ version of the existence of residues of foliations introduced by Baum-Bott
in the complex setting. These are proven in crystalline cohomology thanks 
to the comparison between derived de Rham and crystalline cohomology. 
Finally, the last 
section introduces the notion of infinitesimal derived foliations, explains its
relation with infinitesimal cohomology, and how infinitesimal derived foliations can be always integrated
by formal groupoids.

The Appendix contains some quick and basic reminders on topics in $\s$-categories and derived algebraic geometry (global and formal), together with some important descent results that are used in the main text.

\bigskip
\bigskip

\section*{Guide to our notations}
Most of our notations will be defined in the main text. Here we only fix a few of the most 
recurring notations.\\

Throughout the book, $k$ will be a commutative ring, sometimes supposed to be a field: see the 
beginning of each chapter or section to understand what are the standing hypotheses on $k$.\\ 

We will be freely using the language of model categories and of $\s$-categories (see \S 
\ref{secapp:Modelcategoriesandscategories} for a brief recapitulation about these, and for the 
corresponding notations we adopted). 
We will denote by $\mathbf{Top}$ the $\s$-category of spaces (\cite[\S 1.2.16]{htt}), and by 
$\mathbf{Cat}_{\infty}$ the $\s$-category of $\s$-categories  (see \cite[Chapter 3]{htt}).  \\

For a model category (or a category with weak equivalences) 
$\mathsf{C}$, we will denote by $\mathbf{C}$ its associated $\infty$-category obtained by 
inverting in 
the $\infty$-categorical sense the weak equivalences (see \S 
\ref{secapp:Modelcategoriesandscategories} 
for more details). For a functor $F: \mathsf{C} \to \mathsf{D}$ that induces a functor between 
the 
associated $\infty$-categories (e.g. a functor preserving weak equivalences or a left/right 
Quillen 
functor), we will denote with boldface fonts the induced functor $\mathbf{F}:\mathbf{C} \to 
\mathbf{D}$.\\

The $\infty$-category of derived stacks over $k$ (mostly for the \'etale topology) 
will be denoted by $\dSt_k$ (see \S \ref{secapp:Derivedschemesandstacks}), and, for $G$ a 
group scheme or group stack over $k$, $\dSt^G_k$ will be the $\infty$-category of derived 
stacks endowed with a $G$-action.

\chapter{Algebraic preliminaries}\label{chapter:Algebraicpreliminaries}

In this first chapter we gather preliminaries results that we will use all along the book. 
We fix the cohomological conventions, and present the central notion of graded mixed
complexes. We show that it is related to complete filtered dg-modules by means of the
Tate-realization functor. In the last two section we present a stacky interpretation
of graded mixed complexes as quasi-coherent sheaves on the classifying stack
of the group stack $\Gm \ltimes B\Ga$. Finally, we study the moduli spaces of graded mixed structures. \\

All along this chapter we denoted by $k$ a ground commutative ring. It will be assumed 
to be a $\QQ$-algebra starting from Section \ref{sec:GradedmixedcomplexesandderiveddeRhamtheory} on, but 
most of the results of the chapter are independent of the characteristics of the ground ring.

\section{Homological conventions}\label{sec:Homologicalconventions}

In this section $k$ is of arbitrary characteristics. The sole purpose of this section is to fix notations
and the conventions used in this book. 

\subsection{Basics}\label{subsec:Basics}

Complexes of $k$-modules are
by definition cohomologically indexed and their differential $d$ will raise the 
cohomological degree by $1$. More explicitly, a complex of $k$-modules consists
of a $\ZZ$-indexed family of $k$-modules $\{E^n\}_{n\in \ZZ}$, together
with $k$-linear morphisms
$$d^n : E^n \to E^{n+1}$$
such that $d^{n+1}d^n=0$. We will often write symbolically 
such a complex as a direct sum $E=\oplus_n E^n$, and $d$
as a $k$-linear endomorphism $d : E \to E$ satisfying $d^2=0$.

For two complexes $E$ and $F$, a morphism $f : E \to F$ consists
of a family of $k$-linear morphisms $f^n : E^n \to F^n$ which commute
with the differentials of $E$ and $F$
$$d_F^nf^n=f^{n+1}d_E^n,$$ where we have denoted by $d_E$ and $d_F$ the differentials of $E$ and $F$, respectively.
The complexes of $k$-modules and morphisms between them form a category 
denoted by $C(k)$ \index{$C(k)$}. 

For an object $E \in C(k)$, we define its $i$-th cohomology $k$-module as
$$H^i(E):=\frac{(Ker d^i : E^i \to E^{i+1})}{(Im d^{i-1} : E^{i-1}\to E^i)}.$$
This defines functors $H^i(-) : C(k) \longrightarrow k-Mod$, for $i\in \ZZ$. \\

\subsection{Symmetric monoidal structure}\label{subsec:Symmetricmonoidalstructure}

The category $C(k)$ is endowed with a symmetric monoidal structure $\otimes_k$
defined as follows. For two objets $E$ and $F$, we define $E\otimes_k F \in C(k)$
by letting for $n\in \ZZ$
$$(E\otimes_k F)^n:=\bigoplus_{i+j=n}E^i \otimes_k F^j,$$
the differential on $E\otimes_k F$ being defined as the $k$-linear map
sending a tensor $x\otimes y$ to 
$$d(x\otimes y):=d(x)\otimes y + (-1)^{|x|}x\otimes d(y)$$
where $|x|$ is the cohomological degree of $x$, i.e. $x\in E^{|x|}$. In more concrete terms, the differential
$(E\otimes_k F)^n \to (E\otimes_k F)^{n+1}$ is the
sum of all maps
$$(d_E^i \otimes id , (-1)^i id\otimes d_F^j) : E^i\otimes_k F^j \longrightarrow
(E^{i+1}\otimes_k F^j) \oplus (E^i \oplus_k F^{j+1}) \subset \bigoplus_{a+b=n+1}E^a \otimes_k F^b.$$
This defines the functor $-\otimes_k- : C(k)\times C(k) \to C(k)$. It is
furthermore endowed with the natural unity and associativity constraints isomorphisms
induced from the usual unity and associativity constraints of tensor products
of $k$-modules. The symmetry constraint is itself the usual isomorphism

\begin{equation}
\label{eq:volte} V_{E,F} : E \otimes_k F \simeq F\otimes_k E
\end{equation}

\noindent defined by sending $x\otimes y$ to $(-1)^{ij}y\otimes x$
for $x \in E^{i}$ and $y \in F^{j}$. 
With these definitions, $(C(k),\otimes_k)$ is a closed symmetric monoidal 
category.

\subsection{Shifts}\label{subsec:Shifts}

We denote by $k[1]$ the object of $C(k)$ which is defined by 
$$(k[1])^n=\left\{ \begin{array}{cc}
k \,\,\, \mathrm{if} \; n=-1 \\
0 \,\,\, \mathrm{if} \; n\neq -1
\end{array} \right.$$
obviously endowed with the zero differential.
The shift endofunctor on $C(k)$ is defined to be 
$$E \mapsto E[1]:=k[1]\otimes_k E.$$
Unfolding the definitions of the tensor product we see that $E[1]$ is explicitly given,
up to a natural isomorphism, by
$(E[1])^n=E^{n+1}$ with differential 
$$-d^{n+1} : (E[1])^n=E^{n+1} \longrightarrow (E[1])^{n+2}=E^{n+2}.$$
The endofunctor $E \mapsto E[1]$ is an auto-equivalence of the category $C(k)$
and its inverse will be denoted by $E \mapsto E[-1]$. Finally, we set, for
all $i\in \ZZ$ and $E \in C(k)$
$$E[i]:=(k[1])^{\otimes i}\otimes_k E$$
the $n$-th composite of either $E \mapsto E[1]$ if $i\geq 0$, or $E \mapsto E[-1]$ if $i<0$. 
Explicitly, again up to a natural isomorphism, $E[i]$ is given by 
$(E[i])^n=E^{n+i}$ with differential $(-1)^id^{n+i} : (E[i])^n \to (E[i])^{n+1}$. \\

\subsection{Model category structure}\label{subsec:Modelcategorystructure}

We remind that $C(k)$ possesses a model category structure whose fibrations are
the morphisms $f : E \to F$ with $f^n : E^n \to F^n$ surjective for all $n\in \ZZ$, and 
weak equivalences are the quasi-isomorphisms, i.e. 
the morphisms $f : E \to F$ such that for all $i$ the induced
map $H^i(f) : H^i(E) \to H^i(F)$ is bijective. The cofibrations are defined by 
the left lifting property with respect to the class of fibrations that are also quasi-isomorphisms. This model category structure is suitably compatible
with the symmetric monoidal structure and thus makes $C(k)$ \index{$C(k)$} into a 
\emph{symmetric monoidal model category} in the sense of \cite[\S 4]{hov}. 
It satisfies moreover the monoid axiom of \cite{MR1734325}. 

As explained in the Appendix \ref{secapp:Modelcategoriesandscategories}, 
the model category $C(k)$ gives rise to an 
$\s$-category by localization along the weak equivalences.
This $\s$-category is called the \emph{$\s$-category of complexes}, and will 
be denoted by
$$\dg_k:=L(C(k)).$$
A concrete model for $\dg_k$ \index{$\dg_k$} is given by considering the category of cofibrant\footnote{Note that all objects in $C(k)$ are, obviously, fibrant.} objects
$C(k)^c \subset C(k)$. Between two such cofibrant objects $E$ and $F$, we can define
a simplicial set of maps $Hom^{\Delta}(E,F) \in sSets$, whose set of $n$-dimensional simplices
is given by the formula
$$Hom^{\Delta}(E,F)_n:=Hom_{C(k)}(\mathsf{C}(\Delta^n)\otimes_k E,F)$$
where $\mathsf{C}(\Delta^n) \in C(k)$ is the normalized chain complex of homology
of the simplicial set $\Delta^n$. More precisely, if $\mathsf{N}: sMod_k \to C^{-}(k)$ is the Dold-Kan normalization functor (see ...), $C(\Delta^{\bullet}):=\mathsf{N}(\Delta^{\bullet})$ is a cosimplicial object in $C^{-}(k) \subset C(k)$, hence $$Hom^{\Delta}(E,F):=Hom_{C(k)}(\mathsf{C}(\Delta^{\bullet})\otimes_k E,F) \in sSets.$$
It is easy to see that these data organize themselves 
into an $\mathbb{S}$-category (i.e. a category enriched in $sSets$) which is a concrete model 
for $\dg_k$ (see \S \ref{secapp:Modelcategoriesandscategories} for more details).

Finally, the symmetric monoidal model structure on $C(k)$ induces a natural 
symmetric monoidal structure on the $\s$-category $\dg_k$ 
(see \S \ref{secapp:Modelcategoriesandscategories}). It will
be again denoted by $\otimes_k$. Note however, that the localization $\s$-functor
$$C(k) \longrightarrow \dg_k$$
is not a monoidal $\s$-functor, even though its restriction to 
cofibrant objects is so. Therefore, $E \otimes_k F$ is an ambiguous
notation, whose meaning depends on whether $E$ and $F$ are considered as objects
in $C(k)$ or as objects in $\dg_k$. We hope however that the context will provide 
the necessary extra information, each time we will use such a notation in the text.

For an associative and unital dg-algebra $A$, we will also denote by $Mod(A)$ \index{$Mod(A)$}
the category of (left) $A$-dg-modules. It is endowed with an induced model category structure
whose weak equivalences (resp. fibrations) are morphisms which are quasi-isomorphisms
(resp. epimorphisms) on the underlying complexes (see e.g. \cite{MR1734325}). The corresponding $\s$-category 
will be denoted by $\dg_A$ \index{$\dg_A$}. For a morphism of dg-algebras $\varphi: A \to B$, we have
the usual base change functor $B\otimes_A - : Mod(A) \to Mod(B)$, which 
is a left Quillen functor. It induces an adjunction on the corresponding $\s$-categories
$$B\otimes_A - : \dg_A \leftrightarrows \dg_B : (-)_{[\varphi]}$$
When $\varphi$ will be clear from the context, the functor $(-)_{[\varphi]}$ (restriction of scalars) will be simply omitted from our notations.

\section{Graded mixed complexes}\label{sec:Gradedmixedcomplexes}

In this section $k$ is of arbitrary characteristics.

\subsection{Graded complexes}\label{subsec:Gradedcomplexes}

We consider $C(k)^{gr}$ \index{$C(k)^{gr}$} the category of $\ZZ$-graded objects in $C(k)$. An object
in $C(k)^{gr}$ is by definition a collection $E=\{E^{(n)}\}_{n \in \ZZ}$ of complexes
$E^{(n)}$. Morphisms $f : E \to F$ in $C(k)^{gr}$ are defined levelwise as families 
of morphisms of complexes 
$$f^{(n)} : E^{(n)} \longrightarrow F^{(n)}.$$
In order to avoid confusions with the notion of cohomological degree we
refer to $E^{(n)}$ as the \emph{part of weight $n$ of $E$}, and in general 
refer to \emph{weights} for the degree corresponding to the extra
$\ZZ$-grading on objects in $C(k)^{gr}$. \\

The category $C(k)^{gr}$ is endowed with the usual monoidal structure of graded 
objects. 
For two objects $E$ and $F$ of $C(k)^{gr}$ we set 
$$(E\otimes_k F)^{(n)}:=\bigoplus_{i+j=n}E^{(i)}\otimes_k F^{(j)}.$$
This defines a monoidal structure on the category $C(k)^{gr}$. We define a symmetry 
constraints on $C(k)^{gr}$ induced from 
the symmetry constraint of the monoidal structure on $C(k)$. In formula, this is the
family of isomorphisms
$(E\otimes_k F)^{(n)} \simeq (F\otimes_k F)^{(n)}$ which is the 
sum over $i+j=n$ of the
symmetry maps of complexes reminded in \eqref{eq:volte} on page \pageref{eq:volte}
$$V_{E^{(i)},F^{(j)}} : 
E^{(i)} \otimes_k F^{(j)} \longrightarrow F^{(j)} \otimes_k E^{(i)}$$
We warn the reader here that the signs rule only involves the cohomological degrees and
\emph{not} the weights themselves. With these definitions the category $C(k)^{gr}$ is
endowed with a symmetric monoidal structure $\otimes_k$. 

The model category structure on $C(k)$ defines a natural model category structure
on $C(k)^{gr}$ \index{$C(k)^{gr}$} by defining fibrations, equivalences and cofibrations
of graded complexes levelwise on their weight pieces. A morphism 
$f  : E \to F$ in $C(k)^{gr}$ is thus a fibration (resp. cofibration, resp. 
weak equivalence) if for each $n\in \ZZ$ the morphism $f^{(n)} : E^{(n)} \to 
F^{(n)}$ is
so in $C(k)$. This model category structure is compatible with the symmetric monoidal
structure and makes $C(k)^{gr}$ into a symmetric monoidal model category in the sense
of \cite[\S 4]{hov}

\subsection{Weight twists}\label{subsec:Weighttwists}

As for the shift of complexes, we have weight shifts, which will be rather referred to 
as \emph{weight twists} \index{\emph{weight twists}}. We denote the weight twists 
of an object $E$ by $E(i)$, and it is formally defined by the formula
$$(E(i))^{(n)}:=E^{(n-i)}.$$
We warn the reader about the discrepancy between cohomological shifts and weight shifts: 
for $F \in C(k)$ and $E \in C(k)^{gr}$ we have
$$(F[i])^n=F^{n+i} \qquad (E(i))^{n} = E^{(n-i)}.$$
Therefore, with these conventions, if $E$ is pure of weight $0$, that is $E^{(i)}=0$ if $i\neq 0$, then
$E(1)$ is pure of weight $1$. 

\subsection{Mixed structures on graded complexes}\label{subsec:Mixedstructuresongradedcomplexes}

Because of its importance, we display the notion of mixed structure as a definition by its
own. 

\begin{df}\label{dI-1}
A \emph{graded mixed structure} \index{\emph{graded mixed structure}} on a graded complex $E \in C(k)^{gr}$ consists
of a family of morphisms of complexes
$$\epsilon^{(n)} : E^{(n)} \longrightarrow E^{(n+1)}[-1],$$
for $n\in \ZZ$, such that
\begin{equation}\label{display1}(\epsilon^{(n+1)}[-1])\circ \epsilon^{(n)} : E^{(n)} \to E^{(n+2)}[-2]\end{equation}
equals the zero map for all $n$.

A graded complex $E$ endowed with a graded mixed structure is called 
a \emph{graded mixed complex}.
\end{df}

In order to keep notations simple we will most often write 
$\epsilon : E^{(n)} \to E^{(n+1)}[-1]$
instead of $\epsilon^{(n)}$, and the structural equation (\ref{display1}) by $\epsilon^2=0$. 
Note that graded mixed complexes can also be considered as bi-graded
$k$-modules (for the cohomological degree and the weight) 
endowed with two differentials $d$ and $\epsilon$, satisfying the
following conditions.

\begin{enumerate}
\item The differential $d$ is of bidegree $(1,0)$ and $\epsilon$ is of bidegree $(-1,1)$.

\item We have $d^2=\epsilon^2=0$.

\item We have $d\epsilon + \epsilon d=0$.

\end{enumerate}

Therefore, graded mixed complexes can also be considered as bicomplexes, by 
reindexing the weights in order to get a second differential of bidgree $(0,1)$. 
However, we recommend the reader not doing so, as for instance the symmetry constraints
we defined on graded complexes do not match the usual sign rules of 
tensor products of bicomplexes. It is also important, for many purposes in the sequel, to 
consider $\epsilon$ as an \emph{external} differential which should not be treated on the same footing as
the \emph{internal} cohomological differential. \\

Graded mixed complexes form a category in an obvious manner. The morphisms
between two graded mixed complexes $E$ and $F$ simply being the morphisms
of graded complexes which commute with the mixed structures $\epsilon$ on both sides.
More explicitly, a morphism from $E$ to $F$ consists of a family of morphisms
$f(n) : E^{(n)} \to F^{(n)}$ in $C(k)$ such that the following diagram commutes
$$\xymatrix{
E^{(n)} \ar[rr]^-{f^{(n)}} \ar[d]_-{\epsilon^{(n)}} & & F^{(n)} 
\ar[d]^-{\epsilon^{(n)}} \\
E^{(n+1)}[-1] \ar[rr]_-{f^{(n+1)}[-1]} & & F^{(n+1)}[-1].}$$

This category is denoted by $\epsilon-C(k)^{gr}$ \index{$\epsilon-C(k)^{gr}$}. \\

It is useful to present $\epsilon-C(k)^{gr}$ also 
as a category of modules over a certain algebra in the monoidal category $C(k)^{gr}$ of graded
complexes. For this, we introduce $k[\epsilon] \in C(k)^{gr}$, the object defined as
$$k[\epsilon]:=k(0)\oplus k(1)[1].$$
According to our conventions with shifts and twists (see section \ref{subsec:Weighttwists})
this means that the weight pieces of $k[\epsilon]$ are given by
$$k[\epsilon]^{(n)}= \left\{ \begin{array}{cc}
k & if \, n=0 \\
k[1] & if \, n=1 \\
0 & if \, n\neq 0,1.
\end{array} \right.$$

The graded complex $k[\epsilon]$ possesses a natural monoid structure, with respect to the
tensor product of graded complexes mentioned in \S \ref{subsec:Gradedcomplexes}. It is
given by a morphism of graded complexes
$$k[\epsilon]\otimes_k k[\epsilon] \longrightarrow k[\epsilon]$$
defined by
$$k[\epsilon]\otimes_k k[\epsilon] \simeq
k(0) \oplus k[1](1) \oplus k[1](1) \oplus k[2](2) \longrightarrow k\oplus k[1](1)$$
sending $k[2](2)$ to zero, $k$ isomorphically to $k$ by the identity, and 
each factor $k[1](1)$ isomorphically to $k[1](1)$ by the identity as well. In other words,
as a monoid in $C(k)^{gr}$, $k[\epsilon]$ is the trivial square zero extension of $k$ by 
$k[1](1)$ (see \cite[\S 1.2.1]{hagII}). It is, in particular, an associative, unital and commutative monoid. 

We can now rephrase the definition of graded mixed complexes \ref{dI-1} as follows. A
graded mixed structure on a graded complex $E \in C(k)^{gr}$ consists of the data
of a (unital) $k[\epsilon]$-module structure
$$k[\epsilon] \otimes_k E \longrightarrow E.$$
Indeed, the structure map of such a module structure provides 
$$k[\epsilon] \otimes_k E  \simeq E \oplus E(1)[1] \longrightarrow E,$$
sum of the identity and a morphism 
$E \longrightarrow E(-1)[-1]$ of graded complexes. In degree 
$n$ this
defines a morphism of complexes
$$\epsilon^{(n)} : E^{(n)} \longrightarrow E^{(n+1)}[-1].$$
The fact that $\epsilon$ defined as above satisfies $\epsilon^2=0$ is the 
consequence of the
associativity of the module structure. This construction provides a functor
$$\epsilon-C(k)^{gr} \longrightarrow k[\epsilon]-Mod,$$
which is an isomorphism of categories. Note that the right hand side of this 
equivalence of categories consists of modules inside the monoidal category 
$C(k)^{gr}$, 
which by unfolding the definitions amount to graded
dg-modules over the graded dg-algebra $k[\epsilon]$ (endowed with the zero 
cohomological
differential). \\

The fact that the category $\epsilon-C(k)^{gr}$ can be written as 
a category of modules $k[\epsilon]-Mod$ possesses an important consequence, 
namely the existence of model category structures. Indeed, according to
\cite[Thm. 3.2]{keller} there exists two model category structures on $\epsilon-C(k)^{gr}$, 
the projective model structure and the injective model structure, defined using the forgetful functor
$\epsilon-C(k)^{gr} \longrightarrow C(k)^{gr}$.
The following Proposition summarizes the situation.

\begin{prop}\label{pI-1}
There exists a projective (resp. an injective) model category structure
on $\epsilon-C(k)^{gr}$ whose fibrations (resp. cofibrations) 
and weak equivalences are defined
on the underlying graded complexes obtained by forgetting the mixed structures.
The forgetful functor $\epsilon-C(k)^{gr} \to C(k)^{gr}$ is a left (resp. right)
Quillen functor for the projective (resp. injective) model category structure.
\end{prop}

The $\s$-category associated to the model category $\epsilon-C(k)^{gr}$
will be denoted by $\egrdg_k$. By definition it does not depend on which of the two
model category structure we choose as they both share the same notion of weak 
equivalences. \\

\begin{rmk}
\emph{In the proposition above, the expression} injective model structure \emph{can be confusing. Indeed, 
the injective model structure on $\epsilon-C(k)^{gr}$ is here relative to the projective model structure
on $C(k)^{gr}$. In particular, it is not the same as the usual (absolute) injective model structure 
for which cofibrations are defined to be monomorphisms. One important advantage of the injective
model structure of proposition \ref{pI-1} is the fact that it remains a monoidal model category (which 
is not the case for the absolute injective model structure).}
\end{rmk}

The category $\epsilon-C(k)^{gr}$, of graded mixed complexes, carries a natural 
symmetric monoidal structure denoted by $\otimes_k$ and defined as follows. 
For two objects $E$ and $F$ of $\epsilon-C(k)^{gr}$, we first form the 
tensor product $E\otimes_k F$ as a graded complex. We then 
define a mixed structure on $E\otimes_k F$ whose component of weight $n$ is
given as a morphism of complexes
$$\epsilon^{(n)} : \bigoplus_{i+j=n}E^{(i)} \otimes_k F^{(j)}
\longrightarrow \bigoplus_{a+b=n+1}E^{(a)} \otimes_k F^{(b)}[-1]$$
sending $x\otimes y \in  E^{(i)} \otimes_k F^{(j)}$ to 
$$\epsilon(x)\otimes y + (-1)^{|x|}x\otimes \epsilon(y) \in 
E^{(i+1)} \otimes_k F^{(j)}[-1] \bigoplus E^{(i)} \otimes_k F^{(j+1)}[-1],$$
where as usual $|x|$ denotes the cohomological degree of $x$. It is straightforward to check 
that this defines a graded mixed structure on the graded complex $E\otimes_k F$, and
thus an object, still denoted by $E \otimes_k F$, in the category $\epsilon-C(k)^{gr}$. 
From the formula, it is easy to see that the associativity, unital and symmetry constraints
form the monoidal structure of graded complexes, are all compatible with this definition, and
thus provide an associative, unital and symmetric monoidal structure on $\epsilon-C(k)^{gr}$.
By construction, it is such that the faithful functor
$$F : \epsilon-C(k)^{gr} \longrightarrow C(k)^{gr},$$
obtained by forgetting the mixed structure, is a symmetric monoidal functor
when endowed with the canonical monoidal structure identifying $F(E\otimes_k F)$ with
$F(E)\otimes_k F(F)$ by means of the identity morphism.

The symmetric monoidal structure on $\epsilon-C(k)^{gr}$ can also be described in terms
of a commutative and cocommutative Hopf structure on the object $k[\epsilon] \in C(k)^{gr}$. 
The comultiplication for this Hopf structure is a morphism of monoids in $C(k)^{gr}$
$$\Delta : k[\epsilon] \longrightarrow k[\epsilon]\otimes_k k[\epsilon]$$
obtained by sending $\epsilon$ to $\epsilon \otimes 1 + 1\otimes \epsilon$. It should be noted
here that the multiplication on $k[\epsilon]\otimes_k k[\epsilon]$ is defined as the composition
$$\xymatrix{k[\epsilon]\otimes_k k[\epsilon] \otimes_k k[\epsilon]\otimes_k k[\epsilon] 
\ar[rr]^-{id\otimes V \otimes id} & & k[\epsilon]\otimes_k k[\epsilon] \otimes_k k[\epsilon]\otimes_k
k[\epsilon] \ar[r]^-{\mu \otimes \mu} & k[\epsilon]\otimes_k k[\epsilon],}$$
where $V=V_{k[\epsilon],k[\epsilon]}$ is the volte (see \eqref{eq:volte} on page 
\pageref{eq:volte}), and $\mu$ the multiplication map. 
In particular, $(1\otimes \epsilon).(\epsilon \otimes 1)=
-(\epsilon \otimes1).(1\otimes\epsilon)$, and thus $(\epsilon \otimes 1 + 1\otimes \epsilon)^2=0$
as required. This comultiplication makes $k[\epsilon]$ into a commutative and
cocommutative bialgebra object in 
$C(k)^{gr}$ and thus induces a symmetric monoidal structure on its
category of modules $k[\epsilon]-Mod\simeq \epsilon-C(k)^{gr}$. We leave to the reader that
this monoidal structure is naturally isomorphic to the one already discussed above. \\

Finally, the two model category structures on $\epsilon-C(k)^{gr}$ are
compatible with the monoidal structure, and makes $\epsilon-C(k)^{gr}$ into 
a symmetric monoidal model category, for both the injective and the projective model structures.
As a result, the $\s$-category $\egrdg_k$ of graded mixed complexes over $k$ 
carries a canonical symmetric monoidal structure, again denoted by $\otimes_k$ (see
\ref{secapp:Modelcategoriesandscategories}. 

We would like to finish this section by mentioning that the stable $\s$-category $\egrdg_k$ 
carries a natural t-structure for which an object $E$ is connective if for any $n\in \ZZ$
the part of weight $n$ $E^{(n)}$ is $n$-connective: $H^i(E^{(n)}) = 0$ if $i<-n$. This t-structure
is directly related to the Beilinson's t-structure on filtered complexes by 
the equivalence \ref{clI-4} proved in the next section, and therefore will be referred
as the Beilinson's t-structure on graded mixed complexes as well.

\begin{df}\label{dBeilinsontstructure}
The \emph{Beilinson's t-structure} on the $\s$-category $\egrdg_k$ of graded mixed complexes
is the t-structure for which an object $E$ is connective if $E^{(n)}$ is $n$-connective 
for all $n\in \ZZ$.
\end{df}

\section{Filtered complexes}\label{sec:Filteredcomplexes}

In this section $k$ is of arbitrary characteristics.

\subsection{The homotopy theory of filtered complexes}\label{subsec:Thehomotopytheoryoffilteredcomplexes}
We remind that a filtered complex $E$ (of $k$-modules) consists of 
a sequence of complexes $F^nE$, indexed by $n\in \ZZ$, together with morphisms
of complexes $F^{n}E \to F^{n+1}E$ for all $n$. In another language, 
$E$ consists of a functor
$$\ZZ \to C(k),$$
whose values at $n$ is $F^nE$, and
where $\ZZ$ is the ordered group of integers considered as a category
in the usual manner. Filtered complexes form a category in a natural way, 
for which the morphisms are defined to be the natural transformations
of functors $\ZZ \to C(k)$. We denote by $C(k)^{fil}$ this category, which is
therefore the category of all functors from the ordered group $\ZZ$ to $C(k)$.
Note that a morphism $E \to E'$ consists of a family of morphisms of complexes
$F^nE \to F^nE'$ for $n\in \ZZ$, such that all the diagrams 
$$\xymatrix{
F^nE \ar[r] \ar[d] & F^nE' \ar[d] \\
F^{n+1}E \ar[r] & F^{n+1}E'}$$
are commutative.

\begin{rmk}\label{rI-1}
\emph{Our conventions for filtered complexes differ from the one of \cite{Lurie} and \cite{Moulinos}. In our setting
the filtrations are increasing, whereas the filtrations of these references are decreasing by definition, 
and therefore correspond to functors $\ZZ^{op} \to C(k)$ instead. These two notions are
equivalent by means of the isomorphism of ordered groups $-1 : \ZZ \to \ZZ^{op}$, 
sending $n$ to $-n$. In more concrete terms, we can turn an increasing 
filtered complex $E$ into a decreasing one $E'$ by setting $F^nE':=F^{-n}E$ for
all $n\in \ZZ$.}
\end{rmk}

Being a category of diagrams in $C(k)$, the category $C(k)^{fil}$ comes equipped with a natural
model category structure, for which the fibrations and the weak equivalence 
equivalences are defined levelwise
(the projective model structure on the category of diagrams in a combinatorial 
model category). In a more explicit manner, a morphism $f : E \to E'$ is
defined to be a weak equivalence (resp. a fibration)
if and only if for all $n\in \ZZ$ the morphism of complexes
$F^nE \to F^nE'$ is a quasi-isomorphism (resp. an epimorphism). The corresponding 
$\s$-category (obtained by localization, see \ref{secapp:Modelcategoriesandscategories})
will be denoted by $\fdg_k$. \\

The category $C(k)^{fil}$ is also endowed with a symmetric monoidal structure. It can 
be defined explicitly, for two filtered complexes $E$ and $E'$, by the following formula
$$F^n(E\otimes_k E'):=\colim_{i+j\leq n}(F^iE)\otimes_k (F^jE).$$
The colimit is here taken over the category whose objects are pairs $(i,j)$ of integers
with $i+j\leq n$, and for which there is a unique morphism $(i,j) \to (i',j')$
if and only if $i\leq i'$ and $j\leq j'$. In a more functorial manner, 
this monoidal structure is obtained as the convolution product on the category
of functors $Fun(\ZZ,C(k))$, using that addition in $\ZZ$ to define the
structure of a symmetric monoidal category. For two such functors $E,E' : \ZZ \to C(k)$, 
we can form their external product
$$E\boxtimes_k E' : \ZZ \times \ZZ \to C(k)$$
sending $(i,k)$ to $F^iE \otimes_k F^jE'$. The addition defines 
a functor $+ : \ZZ \times \ZZ \to \ZZ$, and we can thus form the left 
Kan extension along $+$ to get the tensor product
$$E\otimes_k E' := (+)_!(E\boxtimes_k E') \in Fun(\ZZ,C(k)).$$
We refer the reader to \cite{Lurie,Moulinos} for more details on this construction. 

The monoidal structure on $C(k)^{fil}$ comes equipped with natural 
associative, unital and symmetry constraints. Moreover, it is compatible
with the model category structure on $C(k)^{fil}$, and makes it into a symmetric
monoidal model category. As a consequence the $\s$-category 
$\fdg_k$ comes equipped with a canonical symmetric monoidal structure (see 
\ref{secapp:Modelcategoriesandscategories}. \\

We have two important functors out of the category $C(k)^{fil}$, namely the
underlying object and associated graded functors. The associated graded is the functor
$$Gr : C(k)^{fil} \longrightarrow C(k)^{gr}$$
sending a filtered complex $E$ to the graded complex $Gr(E)$ whose component of weight $n$ 
is given by
$$Gr(E)^{(n)}:=Coker(F^{n-1}E \to F^nE) \in C(k).$$
We will often use the more standard notation
$$Gr^n(E):=Gr(E)^{(n)}.$$
This functor is a left adjoint and part of a Quillen adjunction
$$Gr : C(k)^{fil} \rightleftarrows C(k)^{gr} : e.$$
The right adjoint $e$ sends a graded object 
$E$ to the filtered object $e(E)$ with 
$F^n(e(E)):=E^{(n)}$, and for which 
the morphism $E^{(n)} \to E^{(n+1)}$ are all defined to be zero. This adjunction 
is a Quillen adjunction as $e$ obivously preserves fibrations and equivalences. 
The left adjoint $Gr$ is moreover endowed with a canonical 
symmetric monoidal structure
$$Gr(E) \otimes_k Gr(E') \simeq Gr(E \otimes_k E').$$
These functorial isomorphisms are defined by the observation that for all $n\in \ZZ$ and
$a+b < n$, the composite morphism
$$F^a(E) \otimes_k F^b(E') \to F^n(E\otimes_k E') \to Gr^n(E\otimes_k E')$$
is the zero map, which gives rise to natural morphisms
$$Gr^i(E) \otimes_k Gr^j(E') \to Gr^n(E \otimes_k E')$$
for $i+j=n$. We leave the reader to check that the sum of all these morphisms over 
$i+j=n$ is an isomorphism. 
The left derived functor of $Gr$ therefore induces
a well defined symmetric monoidal $\s$-functor (see \ref{secapp:Modelcategoriesandscategories})
$$\mathbf{Gr} : \fdg_k \longrightarrow \grdg_k.$$

\begin{rmk}\label{rI-2}
\emph{Note that the functor $Gr : C(k)^{fil} \to C(k)^{gr}$ and its left derived
$\s$-functor $\mathbf{Gr} : \fdg_k \to \grdg_k$ are related by the following diagram
of $\s$-functors 
$$\xymatrix{
C(k)^{fil} \ar[d] \ar[r]^-{Gr} &  C(k)^{gr} \ar[d] \\
\fdg_k \ar@{=>}[ru]^-{h} \ar[r]_-{\mathbf{Gr}} &  \grdg_k.}
$$
where $h$ is a \emph{non-invertible} natural transformation between $\s$-functors $C(k)^{fil} \to \grdg_k$
(measuring the non-commutativity of $Gr$ with the vertical localization functors). We refer to 
\ref{secapp:Modelcategoriesandscategories} for details on the general situation of deriving
symmetric monoidal $\s$-functors.}
\end{rmk}

The second important functor out of the category $C(k)^{fil}$ is the \emph{underlying 
object functor}. It is the functor 
$$(-)^u : C(k)^{fil} \longrightarrow C(k)$$
sending $E$ to the complex 
$$E^u:=\colim_{n\in \ZZ}F^nE.$$
This functor is again the left adjoint of a Quillen adjunction. The right adjoint 
to this adjunction sends a complex $E$ to the constant diagram $F^nE=E$, where all 
morphisms $F^nE \to F^{n+1}E$ are identities. We will denote this filtered object 
by $E_\s$.
Note that in particular we have $Gr(E_\s)\simeq 0$, whereas $(E_\s)^u\simeq E$ by the
counit of the adjunction. In particular $E \to E_\s$ is fully faithful 
and $(-)^u$ is a localization functor.

The functor $(-)^u$ is again endowed with a natural symmetric monoidal structure
$$(E)^u \otimes_k (E')^u \simeq (E\otimes_k E')^u.$$
Both of these functors preserve quasi-isomorphisms, and 
after localization this adjunction induces an adjunction on the corresponding $\s$-categories
$$(-)^u : \fdg_k \leftrightarrows \dg_k : (-)_\s.$$

\begin{rmk}\label{geominterpfiltobj}
\emph{The $\infty$-category $\fdg_k$ can be canonically identified with the $\s$-category $\mathbf{QCoh}([\mathbb{A}^1/\mathbb{G}_m])$ of quasi-coherent complexes on the (geometric) quotient stack $[\mathbb{A}^1/\mathbb{G}_m]$. This idea goes back to C. Simpson (see \cite{Moulinos} for a recent account). Analogously (and more classically), the $\s$-categoty $\grdg_k$ can be canonically identified with the $\s$-category $\mathbf{QCoh}(B\mathbb{G}_m)$ of quasi-coherent complexes on the (geometric) classifying stack $B\mathbb{G}_m$ of $\mathbb{G}_m$.
Through these identifications, the underlying object functor $ (-)^u: \fdg_k \to \dg_k$ becomes the pullback on quasi-coherent complexes along the inclusion $\Spec \, k \simeq [\Gm/\Gm] \to [\mathbb{A}^1/\mathbb{G}_m] $, while the associated graded object functor $\mathbf{Gr}: \fdg_k \to \grdg_k$ is viewed as the pullback along the map $0: B\mathbb{G}_m \to [\mathbb{A}^1/\mathbb{G}_m]$ corresponding to the $\mathbb{G}_m$-equivariant inclusion of the origin in $\mathbb{A}^1$.
The interested reader is referred to \cite{Moulinos} for details.}
\end{rmk}

\subsection{Graded mixed complexes and
filtered complexes}\label{subsec:Gradedmixedcomplexesandfilteredcomplexes}

We are now ready to compare filtered complexes with graded mixed complexes. This
comparison states the existence of a Quillen adjunction between $C(k)^{fil}$
and $\epsilon-C(k)^{gr}$, which fully faifthfully embeds 
the $\s$-category of graded mixed complexes into the $\s$-category 
of filtered complexes. The image of this embedding consists of 
\emph{complete} filtered complexes, and thus we can state that 
the $\s$-category of graded mixed complexes is equivalent
to the $\s$-category of complete filtered complexes. This equivalence can 
be enhanced into an equivalence of symmetric monoidal $\s$-categories once
the monoidal structure on filtered complexes is suitably completed. \\

We start by describing the \emph{Tate realization functor}
$$|-|^t : \epsilon-C(k)^{gr} \longrightarrow C(k)^{fil}.$$
For a given graded mixed complex $E$ the $i$-th layer of the filtered complex $|E|^t$
is given by 

\begin{equation}\label{eq:Taterealization}
F^i|E|^t:=\prod_{p\geq -i}E^{(p)}[-2p],
\end{equation}

\noindent where the infinite product on the right is endowed with the total differential 
$D:=d+\epsilon$. In formulas, for a family of elements $x:=\{x_p\}$ in 
$\prod_{p\geq -i}E^{(p)}[-2p]$, the element $D(x):=\{D(x)_p\}$ is explicitly given by
$$D(x)_p:=d(x_p)+\epsilon(x_{p-1}),$$
where by convention $x_p=0$ if $p < -i$. The morphism 
$F^i|E|^t \hookrightarrow F^{i+1}|E|^t$ is itself defined to be the natural
inclusion of infinite products 
$$\prod_{p\geq -i}E^{(p)}[-2p] \hookrightarrow \prod_{p\geq -i-1}E^{(p)}[-2p],$$
that sets the first component to be zero.

This defines in a functor 
$$|-|^t : \epsilon-C(k)^{gr} \longrightarrow C(k)^{fil}.$$
This functor commutes with arbitrary limits, and therefore by Freyd adjoint theorem 
(see \cite[Prop. 2.9]{MR3338682}), it admits a left adjoint denoted by $\psi$. The functor
$\psi$ can be made explicit using complicated formula, but we will not need this
and for our purpose its existence will be enough.
Moreover, 
by definition of the various model structures, it is clear that 
$|-|^t$ preserves fibrations. The following lemma shows that it 
furthermore preserves equivalences. 

\begin{lem}\label{lI-1}
The functor 
$|-|^t : \epsilon-C(k)^{gr} \longrightarrow C(k)^{fil}$ 
sends quasi-isomorphisms of graded mixed complexes to 
quasi-isomorphisms of filtered complexes.
\end{lem}

\textit{Proof.} Let $E \to E'$ be a quasi-isomorphism of graded
mixed complexes. We must show that for all $i\in \ZZ$ the induced
morphism of complexes
$$F^i|E|^t \to F^i|E'|^t$$
is a quasi-isomorphism. For this, we use that sequence of subcomplexes
$F^j|E|^t \subset F^{i}|E|^t$, for $j\leq i$. Each successive quotient morphism
$$(F^{i}|E|^t)/F^{j+1}|E|^t) \longrightarrow (F^{i}|E|^t)/(F^j|E|^t)$$
is a fibration of complexes, and by the explicit formula for $|E|^t$ 
(\eqref{eq:Taterealization} on page \pageref{eq:Taterealization}) the natural
map
$$F^i|E|^t \longrightarrow \lim_{j\leq i}(F^{i}|E|^t)/(F^j|E|^t)$$
is an isomorphism of complexes. In particular, the canonical map
$$F^i|E|^t \longrightarrow \holim_{j\leq i}(F^{i}|E|^t)/(F^j|E|^t)$$
is itself an quasi-isomorphism of complexes. As a consequence, 
the lemma will follow if we can prove that for all $j\leq i$ fixed, 
the induced morphisme
$$(F^{i}|E|^t)/(F^j|E|^t) \longrightarrow (F^{i}|E'|^t)/(F^j|E'|^t)$$
is a quasi-isomorphism. But this can be proved by induction on $j$ using the 
short exact sequences of complexes
$$\xymatrix{
0 \ar[r] & (F^{j}|E|^t)/(F^{j-1}|E|^t) \ar[r] & (F^{i}|E|^t)/(F^{j-1}|E|^t) \ar[r] & 
(F^{i}|E|^t)/(F^{j}|E|^t) \ar[r] & 0,}$$
together with the identification of complexes $(F^{j}|E|^t)/(F^{j-1}|E|^t) \simeq
E^{(-j)}[-2j]$. 
\hfill $\Box$ \\

An important consequence of the lemma \ref{lI-1} is that 
$|-|^t$ is a right Quillen functor. In particular, we obtain this way an 
induced adjunction of $\s$-categories
$$\psi : \fdg_k \leftrightarrows \egrdg_k : |-|^t.$$

Note that the associated graded of the realization can be described explicitly as
(beware of the change of signs in the gradings !)
\begin{equation}\label{eq:associatedgraded}
    Gr^p(|E|) \simeq E^{-p}[-2p].
\end{equation}

The importance of the realization $\s$-functor is the following result, that will be
used many times during this book.

\begin{prop}\label{pI-2}
The $\s$-functor
$$|-|^t : \egrdg_k \longrightarrow \fdg_k$$
is fully faithful.
\end{prop}

\textit{Proof.} We must show that for any object $E \in \egrdg_k$, the counit 
morphism 
$$\alpha : \psi(|E|^t) \longrightarrow E$$
is an equivalence in $\egrdg_k$. For this, we combine the adjunction 
$(\psi,|-|^t)$ with the adjunction
$$\mathsf{oub} : \egrdg_k \leftrightarrows \fdg_k : F,$$
where the left adjoint $\mathsf{oub}$ is the forgetful $\s$-functor
sending a graded mixed complex to its underlying graded complex. The fact that 
$\mathsf{oub}$ is a left adjoint follows formally from the description of
$\egrdg_k$ as the $\s$-category of graded dg-modules over $k[\epsilon]$ (see
\S \ref{subsec:Mixedstructuresongradedcomplexes}). The $\s$-functor $\mathsf{oub}$ is conservative, and thus
in order to prove that $\alpha$ is an equivalence, it is enough 
to show the following fact: for any graded complex $M\in \grdg_k$, the induced morphism
$$\Map(E,F(M)) \longrightarrow \Map(\psi(|E|^t),F(M))$$
is an equivalence of spaces. Using the adjunction property this is equivalent to the fact that
for all $M\in \grdg_k$, the natural morphism
$$\Map(\mathsf{oub}(E),M) \longrightarrow \Map(|E|^t,|F(M)|^t)$$
is an equivalence of spaces. 

To prove this last statement on mapping spaces, it will be enough to have an explicit description of
the filtered complex $|F(M)|^t$. As a start, $F(M)$ is the graded mixed complex
whose piece of weight $p$ is $M^{(p)} \oplus M^{(p+1)}[-1]$. The mixed structure
is defined by the obvious map
$$\epsilon : M^{(p)} \oplus M^{(p+1)}[-1] \longrightarrow
M^{(p+1)}[-1] \oplus M^{(p+1)}[-2]$$
which is the identity on the factor $M^{(p+1)}[-1]$ and zero on other factors.
Using that $M$ can be written as an infinite product 
$M \simeq \prod_p M^{(p)}(p)$ inside the $\s$-category $\grdg_k$, we see that
$|F(M)|^t$ is itself given, up to an equivalence, as an infinite product of filtered complexes 
$$|F(M)|^t\simeq \prod_{p\in \ZZ}W_p.$$
Here, each $W_p$ is the Tate realization of the graded mixed complex 
$F(M^{(p)}(p))$. Note that $F(M^{(p)}(p))$ only possesses non-zero 
weight pieces in weight $p$ and $p-1$, respectively $M^{(p)}$ and
$M^{(p)}[-1]$. The mixed structure is itself 
non-zero only in weight $p-1$, for which it is the identity 
$$\xymatrix{
F(M^{(p)}(p))(p-1) = M^{(p)}[-1] \ar[r] & 
F(M^{(p)}(p))(p)[-1] = M^{(p)}[-1].}$$
Using the explicit formula for the Tate realization $|-|^t$ (see \eqref{eq:Taterealization}
on page \pageref{eq:Taterealization}), we deduce that $W_p$ is, up to a natural 
equivalence, the filtered complex
given by 
$$F^iW_p= \left\{
\begin{array}{cc}
0 &  if \, i\neq p \\
M^{(p)}[-2p] &  if \, i=-p.
\end{array}
\right.$$

As a result, for any $F' \in \fdg_k$, we have functorial equivalences in $M \in \grdg_k$
$$\Map_{\fdg_k}(F',|F(M)|^t)\simeq \prod_{p\in \ZZ}\Map_{\dg_k}(Gr^{-p}(F')[2p],M^{(p)}).$$
On the other hand, for $E \in \egrdg_k$ we have canonical equivalences of complexes
(see (\ref{eq:associatedgraded}) above)
$$E^{(-p)}[-2p] \simeq Gr^{p}(|E|^t).$$

Therefore, we have canonical equivalences of spaces
$$\Map_{\grdg_k}(\mathsf{oub}(E),M) \simeq \prod_{p\in \ZZ}\Map_{\dg_k}(E^{(p)},M^{(p)}) \simeq 
\Map_{\fdg_k}(|E|^t,|F(M)|^t).$$
By construction, this equivalence is the natural morphism
$Map(\mathsf{oub}(E),M) \longrightarrow Map(|E|^t,|F(M)|^t)$ induced by adjunction.
This
concludes the proof that the counit $\alpha$ is itself an equivalence, and thus
the proof of the Proposition.
\hfill $\Box$ \\

\subsection{Graded mixed complexes as
complete filtered complexes}\label{subsec:Gradedmixedcomplexesascompletefilteredcomplexes}

The full embedding of Proposition \ref{pI-2} can be refined into an 
equivalence of $\s$-categories, by characterizing the essential 
image of the $\s$-functor $\psi$. For this, we introduce the notion
of \emph{completeness} for objects in $\fdg_k$. 

\begin{df}\label{dI-2}
An object $E \in \fdg_k$ is \emph{complete} if we have
an equivalence in $\dg_k$
$$\lim_{i \leq 0} F^iE \simeq 0.$$
\end{df}

Note that the limit and the equivalence of the above definition 
happen inside the $\s$-category $\dg_k$ of complexes. In particular, 
$\lim_{i \leq 0}F^iE$ is here the homotopy limit of the projective system
of complexes
$$\xymatrix{
\dots \ar[r] & F^iE \ar[r] & F^{i-1}E \ar[r] & \dots \ar[r] & F^{-1}E \ar[r] & 
F^0E.}$$
On the level of cohomology, this homotopy limit also involves the derived functor $\lim^1$.
In particular, if $V$ is a $k$-module endowed with a
filtration by submodules $F^iV\subset V$, being complete in the sense of 
Definition \ref{dI-2} is stronger than the 
naive separateness condition stating that $\cap_i F^iV=0$. Indeed, if we denote
by $F^0\hat{V}:=\lim_{i\leq 0} F^0V/F^iV$ the naive completion of $k$-modules, we have natural isomorphisms
$$H^0(\holim_i F^iV) \simeq \cap_i F^iV \qquad
H^1(\holim_i F^iV)\simeq (F^0\hat{V})/(F^0V),$$
and $H^i(\holim_i F^iV)\simeq 0$ for all $i\neq 0,1$. This implies that 
a filtered $k$-module $V$ is complete in the sense of Definition \ref{dI-2} if and only if 
it satisfies the following two conditions
\begin{itemize}
\item $\cap_{i\leq 0} F^iV=0$

\item the inclusion map $F^0V \to F^0\hat{V}$ is bijective.

\end{itemize}

\begin{rmk}\label{rI-3}
\emph{The completeness condition of Definition \ref{dI-2} can also be stated 
as follows. For a filtered complex $E \in \fdg_k$, and any $j \in \ZZ$ we can 
form the natural morphism of complexes
$$F^jE \longrightarrow \lim_{i\leq j}(F^iE)/(F^jE).$$
Then, $E$ is complete if and only if for all $j$ the above morphism
is a quasi-isomorphism. Indeed, the fiber of the above morphism, taken in $\dg_k$, 
is equivalent to $\lim_{i\leq j}F^iE$. Moreover, if $j\geq 0$, the natural 
projection $\lim_{i\leq j}F^iE \longrightarrow \lim_{i\leq 0}F^iE$
is a quasi-isomorphism, as the inclusion of posets $\{i\leq 0\} \subset \{i\leq j\}$
is final. In the same manner, if $j\leq 0$, the projection
$\lim_{i\leq 0}F^iE \longrightarrow \lim_{i\leq j}F^iE$
is an equivalence.}
\end{rmk}

\begin{prop}\label{pI-3}
The essential image of the full embedding 
$$|-|^t : \egrdg_k \hookrightarrow \fdg_k$$
consists of all complete filtered complexes in the sense
of Definition \ref{dI-2}.
 \end{prop}
 
\textit{Proof.} By the Proposition \ref{pI-2}, the $\s$-functor
$|-|^t$ identifies $\egrdg_k$ with a full reflexive sub-$\s$-category
$\mathcal{C}$ of $\fdg_k$. We have to show that an object $E \in \fdg_k$ lies in this
sub-$\s$-category if and only if it is complete.

As a start, from the explicit form of the Tate realization functor 
(see \eqref{eq:Taterealization} 
on page \pageref{eq:Taterealization}) we see that $|E|^t$ is a complete filtered complex
for any graded mixed complex $E$. Therefore, all objects of $\mathcal{C}$ 
are complete filtered complexes. Conversely, let $E$ be a complete
filtered complex. By definition, $E$ lies in 
$\mathcal{C}$ if and only if the adjuntion unit morphism
$E \longrightarrow |\psi(E)|^t$ is a quasi-isomorphism of filtered complexes. 
During the proof of Proposition \ref{pI-2}, we have seen that 
the composition of $\psi$ with the forgetful $\s$-functor $\mathsf{oub} : \egrdg_k \to \grdg_k$, 
is naturally equivalent to $E \mapsto Gr^{-*}(E)[2*]$. Here, 
$Gr^{-*}(E)[2*]$ denotes the graded complexe whose weight $p$ piece is
$Gr^{-p}(E)[2p]$. We deduce from this the following lemma.

\begin{lem}\label{lpI-3}
The $\s$-functor 
$$\psi : \fdg_k \to \egrdg_k,$$
restricted to the full sub-$\s$-category 
of complete filtered complexes, is conservative.
\end{lem}

\textit{Proof of the lemma.} As the forgetful $\s$-functor $\mathsf{oub}$ is itself
conservative, it is enough to prove the conservativity of the composition 
$\mathsf{oub} \circ \psi : \fdg_k \to \grdg_k,$
when restricted to complete filtered complexes.
By the formula $\mathsf\circ {oub} \psi\simeq Gr^{-*}[2*]$ we see that it is even enough to show that 
the family of $\s$-functors $Gr^i$ for $i\in \ZZ$ is a conservative family of $\s$-functors
when restricted to complete filtered complexes. 

Let $E \to E'$ be a morphism of complete filtered complexes, and assume that 
for all $i\in \ZZ$ the induced morphism $Gr^i(E) \to Gr^i(E')$ is a quasi-isomorphism. 
For all $i \in \ZZ$, we have a canonical morphism of complexes
$$\alpha_i : F^iE \longrightarrow \lim_{j\leq i}(F^iE/F^{j}E).$$
The fiber of this morphism is naturally equivalent to $\lim_{j\leq i}F^jE$, 
which is the zero object by completeness of $E$. Therefore, for all $i$, 
the morphism $\alpha_i$ is an equivalence of complexes. Moreover, for all $j\leq i$ 
we have a commutative diagram of complexes, whose horizontal rows are
fibration sequences of complexes
$$\xymatrix{
\ar[d]  Gr^{-j}(E)[2j] \ar[r] & \ar[d] (F^{i}|E|^t)/(F^{j-1}|E|^t) \ar[r] & 
(F^{i}|E|^t)/(F^{j}|E|^t) \ar[d] \\
  Gr^{-j}(E')[2j] \ar[r] & (F^{i}|E'|^t)/(F^{j-1}|E'|^t) \ar[r] & 
(F^{i}|E|^t)/(F^{j}|E|^t).}$$
We therefore see by induction on $j$, than for all $j\leq i$ the induced morphism
$$(F^{i}|E|^t)/(F^{j}|E|^t) \longrightarrow (F^{i}|E'|^t)/(F^{j}|E'|^t)$$
is a quasi-isomorphism. Passing to the limit along $j\leq i$ and using completeness 
as above we have that the induced morphism
$$F^iE \simeq \lim_{j\leq i}(F^iE/F^{j}E) \longrightarrow
F^iE' \simeq \lim_{j\leq i}(F^iE'/F^{j}E')$$
is a quasi-isomorphism for all $i \in \ZZ$. 
\hfill $\Box$ \\

The lemma implies that complete filtered complexes are objects in $\mathcal{C}$, and thus
finishes the proof of the Proposition.
\hfill $\Box$ \\

\subsection{The completed monoidal structure}\label{subsec:Thecompletedmonoidalstructure}

By Proposition \ref{pI-3} we know that the Tate-realization $\s$-functor
$|-|^t : \egrdg_k \longrightarrow \fdg_k$
is fully faithful, and that its essential image consists of all 
complete filtered complexes in the sense of Definition \ref{dI-2}. We will denote by 
$\cfdg_k \subset \fdg_k$ the full sub-$\s$-category of complete filtered complexes, 
so that $|-|^t$ induces an equivalence of $\s$-categories
$$|-|^t : \egrdg_k \simeq \cfdg_k.$$
We will now study the monoidal behaviour of the above equivalence. For this, we start by 
considering the right Quillen functor
$|-|^t : \epsilon-C(k) \longrightarrow C(k)^{fil},$
and show that it carries a natural lax monoidal structure. For two 
graded mixed complexes $E$ and $F$, we construct a natural morphism
\begin{equation}\label{eq:laxmonoidal}
\alpha_{E,F} : |E|^t \otimes |F|^t \longrightarrow |E\otimes F|^t
\end{equation}

\noindent as follows. Such a morphism is determined by a compatible family of maps
$$\alpha_{i,j} : F^i|E|^t \otimes F^j|F|^t \longrightarrow F^{i+j}|E\otimes F|^t$$
for all integers $i$ and $j$. By the explicit formula for $|-|^t$ the morphism $\alpha_{i,j}$
is a morphism of complexes
$$(\prod_{p\geq -i}E^{(p)}[-2p]) \otimes 
(\prod_{q\geq -j}F^{(q)}[-2q]) \longrightarrow \prod_{n \geq -i-j}
(\bigoplus_{a+b=n}E^{(a)}\otimes F^{(b)}[-2n]).$$
On the level of graded modules, the above morphism is given by a family of maps
$$\alpha_n : (\prod_{p\geq -i}E^{(p)}[-2p]) \otimes 
(\prod_{q\geq -j}F^{(q)}[-2q]) \longrightarrow 
\bigoplus_{a+b=n}E^{(a)}\otimes F^{(b)}[-2n].$$
By definition, the morphism $\alpha_n$ symbolically sends
a tensor of the form $\{x_p\} \otimes \{y_q\}$ to 
the element $\sum_{p+q=n}x_p \otimes y_q$ (note that this sum is finite 
as $p$ and $q$ are bounded below by $-i$ and $-j$). The family of $\alpha_n$
for $n\in \ZZ$ provides the morphism of complexes
$$\alpha_{i,j} : F^i|E|^t \otimes F^j|F|^t \longrightarrow F^{i+j}|E\otimes F|^t.$$
These morphisms are associative and symmetric with respect to the associativity
and unit constraints. Moreover, the image by $|-|^t$ of the unit
object $k(0) \in \epsilon-C(k)^{gr}$ is canonically isomorphic to 
$k$, considered as a filtered complex by 
$$F^ik=\left\{
\begin{array}{cc}
0 & if \, i<0 \\
k & if \, i\geq 0
\end{array}\right.$$
In other words, there is a canonical isomorphism $u : |k(0)|^t\simeq k$, 
where $k$ is the unit object in $C(k)^{fil}$. All together, the $\alpha_{i,j}$ 
and the isomorphism $u$ defines a lax monoidal structure on the
functor $|-|^t$, which is unital, associative and symmetric.

By adjunction, the lax monoidal structure on $|-|^t$ defines a unital, associative
and symmetric colax monoidal structure on the left adjoint functor $\psi : 
C(k)^{fil} \to \epsilon-C(k)^{gr}$ as follows. We start by $E$ and $F$ two 
filtered complexes, and we consider the adjunction morphisms
$$E \to |\psi(E)|^t \qquad F \to |\psi(F)|^t.$$
Using the lax monoidal structure on $|-|^t$ constructed above (see \eqref{eq:laxmonoidal}
on page \pageref{eq:laxmonoidal}), we define a morphism
$$\xymatrix{
E\otimes F \ar[r] & |\psi(E)|^t  \otimes |\psi(F)|^t \ar[r] & |\psi(E)\otimes \psi(F)|^t.}$$
Its image by $\psi$, composed with the counit $\psi |-|^t \Rightarrow id$, 
gives a morphism
$$\xymatrix{
\psi(E\otimes F) \ar[r] & \psi(E) \otimes \psi(F),
}$$
which is a unital, associative and symmetric colax monoidal structure on the functor $\psi$.

\begin{lem}\label{lI-4}
The induced colax monoidal structure on the $\s$-functor
$$\psi : \fdg_k \longrightarrow \egrdg_k$$
is a monoidal structure. In other words, for all $E$ and $F$ objects in $\fdg_k$, the natural
morphism
$$\beta_{E,F} : \psi(E\otimes F) \longrightarrow \psi(E) \otimes \psi(F)$$
is an equivalence of graded mixed complexes.
\end{lem}

\textit{Proof.} To check that $\beta_{E,F}$ is a equivalence in $\egrdg_k$ 
we can use the conservative $\s$-functor $\mathsf{oub} : \egrdg_k \to \grdg_k$ to
graded complexes. Moreover, we know that $\mathsf{oub} \circ \psi$ is naturally 
equivalent to $Gr^{-*}[-2*]$ (see the proof of Proposition \ref{pI-2}). By observation, 
we the image of $\beta_{E,F}$ by $\mathsf{oub}$ induces a morphism
$$Gr^{-*}(E\otimes F)[-2*] \longrightarrow Gr^{-*}(E)[-2*] \otimes Gr^{-*}(F)[-2*]$$
which is nothing else than the standard monoidal structure on the $\s$-functor $Gr^{-*}[-2*]$. 
It is therefore an equivalence, proving that $\beta_{E,F}$ is indeed an equivalence
of graded mixed complexes. \hfill $\Box$ \\

The consequence of the above lemma is that the colax monoidal structure on 
the $\s$-functor
$\psi$ is a symmetric monoidal structure. Therefore, we can use
the equivalence 
$$|-|^t : \egrdg_k \simeq \cfdg_k,$$
in order to transport the monoidal structure of $\egrdg_k$ to 
a symmetric monoidal structure on $\cfdg_k$. It will be denoted by
$\hat{\otimes}$, and called the \emph{completed tensor product}. 
The previous lemma, implies that the $\s$-functor $\psi$ provides 
a symmetric monoidal $\s$-functor
$$\psi : \fdg_k \longrightarrow \cfdg_k\simeq \egrdg_k.$$
This exhibits $\cfdg_k$ has a symmetric monoidal localization of 
$\fdg_k$, obtained by inverting the \emph{graded quasi-isomorphisms} (i.e. morphisms
of filtered complexes inducing quasi-isomorphisms on the associated graded).
The localization $\s$-functor $\psi$ will also be called the \emph{completion $\s$-functor}, 
and denoted by 
$$\widehat{(-)} : \fdg_k \longrightarrow \cfdg_k.$$
Note that by construction we have canonical equivalences of complete filtered complexes
$$\widehat{(E\otimes F)} \simeq \widehat{E} \hat{\otimes} \widehat{F}.$$
Note also, that the completion $\s$-functor naturally commutes with 
forming the associated graded
$$Gr^*(\widehat{E}) \simeq Gr^*(E),$$
as can be seen directly from the identification between $\mathsf{oub} \circ \psi$ and $Gr^*[-2*]$. \\

On the level on notations, $|E|^t$ will always refer to the filtered complex associated
to the graded mixed complex $E$. We will also use the following notation 

\begin{df}\label{nonTaterealizationforgradedmixedmodules} We define the (non Tate) \emph{realization} 
$\s$-functor $$|-| : \egrdg_k \to \dg_k \, , \, E \longmapsto |E|:=F^0|E|^t$$ sending a mixed graded 
module $E$ to the $0$-th layer filtration of its Tate realization.
\end{df}

\begin{rmk}\emph{It is easy to verify that the realization functor $|-|$ is equivalent to 
$\mathbb{R}\underline{Hom}(k(0),-)$, the derived enriched Hom complex from the unit
graded mixed object $k(0)$. Note also that, when $E$ has only non-negative weight pieces, we have}
$$|E| \simeq F^0|E|^t\simeq (|E|^t)^u$$
\emph{so $|E|$ is in fact the underlying complex of the filtered complexe $|E|^t$.}
\end{rmk}

We can subsume our comparison results as follows.

\begin{cor}\label{clI-4}
The Tate realization $\s$-functor
$$|-|^t : \egrdg_k \longrightarrow \cfdg_k$$
is a symmetric monoidal equivalence when $\cfdg_k$ is endowed with the complete
tensor product $\hat{\otimes}_k$. It induces an equivalence on the $\s$-categories of
commutative algebras
$$\CAlg(\egrdg_k) \simeq \CAlg(\cfdg_k).$$
\end{cor}

We note here that the $\s$-categories of commutative algebras of the previous corollary
can be described in simple terms as localizations of model categories of strict commutative
algebra objects, thanks to the rectification result of \cite[Prop. 4.5.4.6]{HA}. In particular, 
$\CAlg(\egrdg_k)$ is canonically equivalent to $\egrcdga_k$, the $\s$-category of graded mixed
cdga's described in the next section. In the same manner, $\CAlg(\cfdg_k)$ can be identified
with the localization of the model category of complete filtered cdgas (i.e. commutative algebra
objects in the model category $C(k)^{fil}$ of filtered complexes, whose underlying 
filtered complex is complete).

\begin{rmk}\label{nonTaterealizationforgradedmixedalgebras} \emph{It is easy to check that the (non Tate) realization functor $|-|$ of Definition \ref{nonTaterealizationforgradedmixedmodules} induces a functor $$|-| : \CAlg(\egrdg_k) \longrightarrow \cdga_k.$$ }
\end{rmk}

\section{Graded mixed cdga's and derived de Rham theory}
\label{sec:GradedmixedcomplexesandderiveddeRhamtheory}

In this section we show how graded mixed complexes provide a setting
for de Rham theory of algebraic varieties and schemes, and more generally 
derived stacks. For this
we assume further more that our base ring $k$ is a $\QQ$-algebra. 
The extension to non-zero characteristics require substantial modification and will be treated separately 
in last chapter of this book (see chapter \ref{chapter:nonzerocharacteristics}).

We remind the
symmetric monoidal model category $\epsilon-C(k)^{gr}$ of graded
mixed complexes, and its corresponding symmetric monoidal $\s$-category 
$\egrdg_k$. We will use in the sequel the category of commutative (unital and
associative) monoids in $\epsilon-C(k)^{gr}$, for which we introduce the following
notation
$$\epsilon-CAlg(k)^{gr} := Comm(\epsilon-C(k)^{gr}).$$

Objects in $\epsilon-CAlg(k)^{gr}$ will be referred to as 
\emph{graded mixed cdga's}.
In a more explicit form, a graded mixed cdga $A$ consists of the data
of a graded mixed complex $A$, together with 
\begin{itemize}
\item (multiplication) morphisms of complexes
$$(.) : A^{(i)} \otimes_k A^{(j)} \longrightarrow A^{(i+j)},$$

\item (unit) a cocycle $1 \in A^{(0)}$ of cohomological degree $0$.

\end{itemize}

The multiplication and the unit are required to satisfy the obvious associativity, 
unity and commutativity axioms. Moreover, the multiplication and the unit 
must be compatible with the mixed structure $\epsilon$ in the following sense:
$$\epsilon(1)=0 \qquad \epsilon(x.y)=\epsilon(x).y+(-1)^{|x|}x.\epsilon(y),$$
for all elements $x$ and $y$ in $A$ with $x$ of cohomological degree $|x| \in \ZZ$. 
We note in particular that the sign appearing the above formula
does not depend on the weights of $x$ and $y$, but only on their cohomological degrees. \\

We can use the forgetful functor $\epsilon-CAlg(k)^{gr} \longrightarrow \epsilon-C(k)^{gr}$, 
which forgets the multiplicative structure, in order to create a model category structure
on $\epsilon-CAlg(k)^{gr}$ for which the fibrations and weak equivalences are defined
in $\epsilon-C(k)^{gr}$. The existence of this model category structure uses in an essential
manner that $k$ is assumed to be a $\QQ$-algebra, and can be deduced for instance 
from the general statement \cite[Prop. 4.5.4.6]{HA}. The homotopy theory of graded mixed cdga's is central
in the definition and study of derived foliations, and we thus set the following 
important definition and notations.

\begin{df}\label{dI-3}
The \emph{model category of graded mixed cdga's} is $\epsilon-CAlg(k)^{gr}$
defined above. The corresponding $\s$-category will be denoted by 
$\egrcdga_k$ and will be called the \emph{$\s$-category of graded mixed cdga's over $k$}.
\end{df}

The model category $\epsilon-CAlg(k)^{gr}$ comes equipped with two important 
Quillen adjunctions, which will be used often all along the text. As a start, we have
the adjunction induced by the forgetful functor
$$Sym^{\epsilon} : \epsilon-C(k)^{gr} \rightleftarrows \epsilon-CAlg(k)^{gr} : \mathsf{oub}.$$
Here the left adjoint sends a graded mixed complex $E$ to 
the usual symmetric algebra object
$$Sym^{\epsilon}(E):=\bigoplus_{n\geq 0} E^{\otimes n}/\Sigma_n,$$
where the sum, tensor products and quotients by symmetric groups are all understood
in the category $\epsilon-C(k)^{gr}$ of graded mixed complexes. We have denoted
this functor by $Sym^{\epsilon}$ to avoid confusions with other $Sym$ constructions
that we will use later.

The second important Quillen adjunction is denoted by
$$DR : CAlg(k) \rightleftarrows \epsilon-CAlg(k)^{gr} : (-)^{(0)},$$
where $CAlg(k)$ denotes the model category of (non-graded and non-mixed) commutative
dg-algebras. The right adjoint of this adjunction 
sends a graded mixed cdga $A$ to its weight $0$ part
$A^{(0)}$, which comes equipped with an induced commutative dg-algebra
structure. The left adjoint is called the \emph{de Rham functor} and will be
studied in detailed in the next two paragraphs. 

The two mentioned above Quillen adjunctions of course induce natural adjunctions of 
$\s$-categories (as usual by \ref{secapp:Modelcategoriesandscategories}
$$Sym^{\epsilon} : \egrdg_k \rightleftarrows \egrcdga_k : \mathsf{oub} \qquad
\DR : \cdga_k \rightleftarrows \egrcdga_k : (-)^{(0)}.$$

\subsection{Derivations and the cotangent complex}\label{subsec:Derivationsandthecotangentcomplex}

For a commutative dg-algebra $A \in CAlg(k)$, and a (left) $A$-dg-module $M$, 
we can form the trivial square zero extension of $A$ by $M$ in the usual
manner. It is a new commutative dg-algebra, denoted by $A\oplus M$ or
$A[M]$, whose underlying complex is the direct sum $A \oplus M$. 
The multiplication on $A\oplus M$ is the morphism
$$\mu : \xymatrix{
(A\oplus M) \otimes (A\oplus M) \simeq (A\otimes A) \oplus (A\otimes M) \oplus(M\otimes A)
\oplus (M\otimes M) \ar[r] & A\oplus M,}$$
defined as follows. On the factor $A\otimes A$ the morphism $\mu$
is the multiplication in $A$, followed by the natural inclusion $A \hookrightarrow A \oplus M$.
On the factor $M\otimes M$ the morphism $\mu$ is the zero map. 
On the factor $A\otimes M$, $\mu$ is the given by the left $A$-dg-module
structure on $M$, followed by the natural inclusion $M \hookrightarrow A\oplus M$. 
Finally, on the factor $M\otimes A$, $\mu$ is given by the right $A$-dg-module
structure on $A$ followed by the inclusion $M \hookrightarrow A\oplus M$. 
We note here that the right $A$-dg-module structure on $M$ is defined using the
symmetry constraints on complexes (see \eqref{eq:volte} on page \pageref{eq:volte}), 
and is given explicitly for $a \in A^i$ and $m \in M^j$ by
$m.a=(-1)^{i+j}a.m$. Unfolding the above description of the multiplication on $A\oplus M$, 
we can also write the following formula, for homogeneous elements 
$(a,m) \in A\oplus M$ and $(a',m') \in A\oplus M$
$$(a,m).(a',m')=(a.a',am'+(-1)^{|a'||m|}a'm).$$
With this definition, a derivation on $A$ with values in $M$ is
a section of the natural projection $A \oplus M \longrightarrow A$ inside
the category of commutative dg-algebras (see \cite[\S 1.2.1]{hagII}). Such a section
is of the form $(id_A,u) : A \longrightarrow (A\oplus M)$, with $u : A \longrightarrow M$
a morphism of complexes. The morphism $u$ must also satisfy the usual condition
$$u(a.b)=u(a).b+a.u(b)=a.u(b)+(-1)^{|a||b|}b.u(a)$$
for homogeneous elements $a$ and $b$ in $A$.

For any commutative dg-algebra $A$, there exists a universal derivation 
$$d : A \longrightarrow \Omega_{A/k}^1.$$
Here $\Omega_{A/k}^1$ is the dg-module of Kähler differentials on the commutative dg-algebra
$A$ (relative to $k$), 
which can be explicitly presented as $I/I^2$ where $I \subset A\otimes A$ is the
kernel of the multiplication map $A \otimes A \longrightarrow A$. The universal
property here states that for any left $A$-dg-module $M$ and any derivation $u : A \to M$, 
there exists a unique morphisme of $A$-dg-modules $f : \Omega_{A/k}^1 \to M$
such that $f\circ d = u$
$$\xymatrix{
A \ar[r]^-d  \ar[rd]_-u &  \Omega^1_{A/k} \ar[d]^-{{\exists !}} \\
& M.}$$
The construction $A \mapsto \Omega^1_{A/k}$ is functorial in $A$, and can 
be seen to be a left adjoint of a Quillen adjunction as follows. We form
the category $CAlgMod(k)$, whose objects are pairs $(A,M)$, consisting of
a commutative dg-algebra $A$ and a left $A$-dg-module $M$. The morphisms
from $(A,M)$ to $(A',M')$ are pairs $(u,f)$, where $u : A \to A'$ is a morphism
of commutative dg-algebras and $f : M \to M'$ is a morphism of $A$-dg-modules
(where $M'$ is considered as a left $A$-dg-module via $u$). The category 
$CAlgMod(k)$ admits a model category structure for which 
a morphism $(u,f)$ as above is a weak equivalence (respectively a fibration)
if $u$ and $f$ are quasi-isomorphisms of complexes (respectively 
an epimorphism of complexes). This model category is a particular example
of a general construction of a relative model category structure, as explained
in \cite[Cor. 6.3.16]{harmat}. 

The construction $(A,M) \mapsto A\oplus M$ defines a functor 
$CAlgMod(k) \longrightarrow CAlg(k)$. This functor is a right Quillen 
functor and its left adjoint sends a commutative dg-algebra $A$ to 
the pair $(A,\Omega_{A/k}^1)$. We obtain this way an induced adjunction on the
corresponding $\s$-categories obtained by inversting quasi-isomorphisms
$$\cdga_k \leftrightarrows \cdgamod_k.$$
The left adjoint sends a commutative dg-algebra $A$ to an object of the form
$(A,\mathbb{L}_{A/k})$, where $\mathbb{L}_{A/k}$ is a left $A$-dg-module called
the \emph{cotangent complex of $A$}. An explicit model 
for the $A$-dg-module $\mathbb{L}_{A/k} \in \dg_A$ is given by 
$$\mathbb{L}_{A/k} \simeq A'\otimes_A \Omega_{A'/k}^1,$$
where $\xymatrix{A' \ar@{->>}[r]^-\sim & A}$
is a cofibrant model of $A$.

Being a left adjoint, the $\s$-functor $A \mapsto (A,\LL_{A/k})$
commutes with arbitrary colimits. If $A=\colim_{i}A_i$ inside
the $\s$-category $\cdga_k$, then we have a canonical equivalence
of $A$-dg-modules
$$\LL_{A/k}\simeq \colim_{i}(A\otimes_{A_i}\LL_{A_i/k}).$$

\subsection{The derived de Rham complex as a left 
adjoint}\label{subsec:ThederiveddeRhamcomplexasaleftadjoint}

We consider the functor $E \mapsto E^{(0)}$, from 
graded mixed complexes to complexes, defined by keeping the sole weight 
$0$ piece. This functor possesses a natural symmetric lax monoidal 
structure and therefore naturally extends to a 
functor $A \mapsto A^{(0)}$ from 
the category of graded mixed cdga's to the category of cdga's. This functor
$CAlg(\epsilon-C(k)^{gr}) \longrightarrow CAlg(C(k))$ is the 
is the right adjoint of a Quillen adjunction whose left adjoint is called the 
de Rham functor
$$DR(-/k) : CAlg(C(k)) \longrightarrow CAlg(\epsilon-C(k)^{gr}).$$
Explicitly, the functor $DR$ sends a commutative dg-algebra $A$ to 
the graded mixed cdga whose weight $n$ piece is $Sym^{n}_{A}(\Omega^1_{A/k}[1])$. 
The multiplicative structure on $DR(A)$ is the standard multiplication induced on
a symmetric algebra. The mixed structure is itself induced from 
the usual de Rham differential $dR : A \longrightarrow \Omega_{A/k}^1$, 
and its natural multiplicative 
extension to the whole symmetric algebra $Sym_A(\Omega^1_{A/k}[1])$.
In short
we write $DR(A/k)=Sym_A(\Omega_{A/k}^1[1])$. We leave to the reader to check that 
this construction $A \mapsto DR(A/k)$ does indeed define 
a left adjoint to $(-)^0$, by observing that for any graded mixed cdga $B$ the
morphism $\epsilon : B^{(0)} \to B{(1)}$ is derivation on $B^{(0)}$ with coefficients 
in the $B^{(0)}$-dg-module $B^{(1)}$, and using the univeral property 
of $\Omega_{A/k}^1$.

The functor $A \mapsto DR(A/k)$ being left Quillen it induces 
an $\s$-functor on the corresponding $\s$-categories (see \S \ref{secapp:Modelcategoriesandscategories})
$$\DR : \cdga_k \longrightarrow \egrcdga_k.$$
This $\s$-functor is left adjoint to $A \mapsto A^{(0)}$, sending
a graded mixed cdga $A$ to its peice of weight $0$. The cotangent complexes
being obtained as the left derived functor of $A \mapsto \Omega_{A/k}^1$ we see that 
on the level of $\s$-categories we have canonical equivalences of graded cdga's
$$\DR(A)\simeq Sym_A(\LL_{A/k}[1])$$
where the symmetric algebra is understood in the symmetric monoidal $\s$-category $\grdg_k$ of
graded complexes.
In more concrete terms a model for $\DR(A/k)$ is $DR(A'/k)$ where $A'$ is a
cofibrant model of $A$.

\begin{df}\label{dI-4}
\begin{enumerate}
\item 
The \emph{(derived) de Rham algebra} of a commutative dg-algebra $A$ (over $k$)
is the graded mixed commutative dg-algebra 
$$\DR(A)\simeq Sym_{A}(\LL_{A/k}^1[1]) \in \egrdg_k$$
defined above.

\item The \emph{(Hodge-completed and derived) de Rham complex} of a commutative dg-algebra $A$ (over $k$)
is the filtered complex 
$$\CDR(A/k):=|\DR(A/k)|^t \in \fdg_k.$$

\end{enumerate}
\end{df}

By definition the derived de Rham complex $\CDR(A/k)$ is 
a complete filtered complex, and its filtration is the Hodge filtration (hence the use of the 
expression \emph{Hodge-completed} in the definition above). Its associated graded object is 
the graded cdga $Sym_A(\LL_{A/k}[-1])$, where $\LL_{A/k}$ sits in weight $-1$. 
It also comes equipped with a canonical commutative multiplicative structure, 
because the functor $|-|^t$ is endowed with a canonical symmetric lax monoidal structure (see 
corollary \ref{clI-4}). 
Therefore, $\CDR(A/k)$ can be promoted to a filtered and complete 
commutative dg-algebra over $k$, but we will use this highly structured version rarely in the 
sequel. \\

In the particular case where $A$ is a smooth $k$-algebra, considered as a dg-algebra
concentrated in degree $0$, $\CDR(A/k)$ turns out to be canonically equivalent to
the usual algebraic de Rham complex of $A$ relative to $k$
$$\CDR(A/k) \simeq (\xymatrix{A \ar[r] & \Omega^1_{A/k} \ar[r] & \Omega_{A/k}^2 
\ar[r] & \dots \ar[r] & \Omega_{A/k}^i \ar[r] & \dots})$$
Indeed, when $A$ is smooth the natural morphism $\LL_{A/k} \to \Omega_{A/k}^1$ is a quasi-isomorphism. Therefore, 
for a cofibrant replacement $\xymatrix{A' \ar@{->>}[r]^-\sim & A}$, the
the canonical projection of graded mixed cdga's
$$Sym_{A'}(\LL_{A'/k}[1]) \to Sym_A(\Omega_{A/k}^1[1])$$
is a quasi-isomoprhism. Passing to the Tate realizations, we get an equivalence of filtered and complete cdga's
$$\CDR(A/k) \simeq C^*_{DR}(A/k)$$
where the right hand side is the usual, \emph{underived} de Rham complex of $A$ over $k$.
This is a bounded complex concentrated in degrees $[0,d]$, where $d$ is the 
relative dimension of $A$ over $k$. This complex is filtered by the usual Hodge
filtration, and is complete because the filtration is finite. 
As, by convention, all our filtrations are \emph{increasing}, 
this Hodge filtration is indexed in a non-standard way 
by non-positive integers $p\leq 0$ as follows
$$F^{p}\CDR(A/k) = (\xymatrix{\Omega^{-p}_{A/k} \ar[r] & \Omega_{A/k}^{-p+1} 
\ar[r] & \dots }) \subset \CDR(A/k).$$
The associated graded $Gr^*(\CDR(A/k))$ is thus the usual Hodge
cohomology $\oplus \Omega^i_{A/k}[-i]$, with $\Omega_{A/k}^i$
of weight $-i$. \\

To finish this part we mention below two well known important properties of 
derived de Rham cohomology. \\

\begin{prop}\label{pI-6}
Let $X=\Spec\, A$ be a derived affine scheme of essentially of finite presentation 
over $\CC$ and $X(\CC)$ its 
space of $\CC$-points endowed with its transcendent topology. Then there exists a natural 
equivalence of cdga's
$$\CDR^*(A/\CC) \simeq C^*(X(\CC),\CC),$$
where the right hand side denotes sheaf cohomology with coefficients in the constant sheaf $\CC$.
\end{prop}

\textit{Proof.} This results is proven in
\cite{bhatt2012completionsderivedrhamcohomology}. We can provide a direct proof as follows. 
As a first statement, the right hand side is well known to be functorially equivalent to 
algebraic de Rham cohomology of $X$. Algebraic de Rham cohomology can be computed 
by chosing a closed embedding $X \to Y=Spec\, B$ for $Y$ a smooth affine scheme over $\CC$, and 
then considering the formal de Rham complex $C_{DR}^*(\hat{Y}_X)$ of $Y$ along $X$. 
The formal de Rham complex is obtained by definition from $DR(Y)$ by completion along $X$
$$C^*_{DR}(\hat{Y}_X) = \xymatrix{\hat{B} \ar[r] & \hat{\Omega}_B^1 \ar[r] & \dots \ar[r] & \hat{\Omega}_B^r
\ar[r] & \dots}$$
where $\hat{B} = lim_{n} B/I^n$ with $I \subset B$ the ideal of definition of $X$, and
$\hat{\Omega}_B^r := \hat{B}\otimes_B \Omega_B^r$. 

From the results reminded in \ref{secapp:Cotangentcomplexesofformalcompletion} 
the de Rham complex $C^*_{DR}(\hat{Y}_X)$ is also given by 
$$C^*_{DR}(\hat{Y}_X) \simeq lim_{n} \CDR(Y_n)$$
where $Y_n = B/I^n$. But the Hodge completed derived de Rham cohomology is invariant
by infinitesimal thickening (this can be see along the same lines as 
\ref{pA-8-1}, see \ref{rem:pA-8-1}), so we have 
$$C^*_{DR}(\hat{Y}_X) \simeq \CDR(t_0(X)) \simeq \CDR(X)$$
via the pull-backs along the closed immersions $t_0(X) \to X \to Y$.
\hfill $\Box$ \\

For the next Proposition, recall that for any connective cdga $A$ over $k$, 
the structural morphism $k \to \DR(A/k)$ is a morphism of graded mixed cdga (where
$k$ is endowed with the zero mixed structure). By passing to realizations we get this way
a natural morphism of cdga
$$k \to |\DR^*(A/k)| \simeq \CDR^*(A/k).$$
Suppose now that $A$ is a discrete commutative $\QQ$-algebra and $I \subset A$ an ideal. 
We can apply the previous morphism to the augmentation $A \to A/I$, and get this way 
a natural morphism of cdga's
$$A \to \CDR^*((A/I)/A).$$
Moreover, as this is obtained by realization from a morphism of graded mixed cdgas, 
this morphism is compatible with the augmentations to $A/I$. In particular, this morphism
is compatible with the filtrations: the $I$-adic filtration on $A$ and the
natural filtration on $\CDR^*((A/I)/A)$. As $\CDR^*((A/I)/A)$ is already known to be complete
it induces a morphism of completed filtered cdga's
$$\hat{A} \to \CDR^*((A/I)/A)$$
where $\hat{A}$ is the $I$-adic completion of $A$.

\begin{prop}\label{pI-7}
Let $A$ be a commutative $k$-algebra of finite type and $I \subset A$ an ideal. Then, 
the canonical morphism
$$\hat{A} \to \CDR^*((A/I)/A)$$
is an equivalence. In particular
$$H^0(\CDR^*((A/I)/A))\simeq \hat{A} \qquad H^i(\CDR^*((A/I)/A)) \simeq 0 \; \forall \; i\neq 0.$$
\end{prop}

\textit{Proof.} See \cite[Thm. 4.10]{bhatt2012completionsderivedrhamcohomology}. \hfill $\Box$ \\

\section{The geometry of graded mixed objects}\label{sec:Thegeometryofgradedmixedobjects}

We recall here some results from \cite{mrt} and \cite{alf} which provide a geometric interpretation
of graded mixed complexes. These will be used for some specific
arguments and are sometimes very useful (see for instance \S \ref{sec:Directimages}). \\

We consider the multiplicative group scheme $\Gm$. It acts by standard weight $-1$ on 
the additive group scheme $\Ga$, and thus on the group stack $B\Ga$. Here the weight of the $\Gm$
action is understood geometrically, the induced action of $\Gm$ on functions $\OO(\Ga)=k[T]$
is such that $T$ has weight $1$.

\begin{df}\label{dI-5}
We define the group stack $\cH_0$ as the semi-direct product
$$\cH_0:=B\Ga\rtimes \Gm$$
\end{df}

We can now consider complexes of $k$-modules endowed with an action of $\cH_0$. They form an $\s$-category
which is formally defined as $\QCoh(B\cH_0)$, the $\s$-category of quasi-coherent complexes over the classifying stack  $B\cH_0$ of $\cH_0$. 
The canonical projection $\cH_0 \to \Gm$ possesses a group section $s : \Gm \to \cH_0$, which
provides a natural morphism $s : B\Gm \to B\cH_0$. The fiber of the projection $\cH_0 \to \Gm$
identifies with $B\Ga$, which is also naturally identified with the fiber of the section $s$.
Therefore, $s$ exhibits $B\Ga$ as a $\cH_0$-equivariant stack.

By writing $B\cH_0$ as the standard colimit
$B\cH_0 = \mathrm{colim}_{[n]\in \Delta^{op}} \cH_0^n,$
we can write 
$$\QCoh(B\cH_0) \simeq \lim_{[n] \in \Delta} \QCoh(\cH_0^n).$$
Now, for each $n$, the stack $\cH_0^n$ is an affine stack with cdga of functions
$\OO(\cH_0^n) \simeq \OO(\cH_0)^{\otimes n} \simeq (k[t,t^{-1}]\oplus k[t,t^{-1}][-1](1))^{\otimes n}$. 
The category of graded mixed complexes $\epsilon-C(k)^{gr}$, before localization by quasi-isomorphisms, 
can be understood as the category of $(k[t,t^{-1}]\oplus k[t,t^{-1}][-1](1))$-dg-comodules, 
where $(k[t,t^{-1}]\oplus k[t,t^{-1}][-1](1))$ is the commutative and 
dg-Hopf algebra obtained by semi-direct product of the multiplicative
Hopf algebra $k[t,t^{-1}]$ with $k\oplus k[-1](1)=\OO(B\Ga)$. In other words, we can write
down an equivalence of genuine categories (where the limit is computed in the $2$-category of categories, ad
is taken before localization by quasi-isomorphisms)
$$\epsilon-C(k)^{gr} \simeq \lim_{[n]\in \Delta} (k[t,t^{-1}]\oplus k[t,t^{-1}][-1](1))^{\otimes n}-dg-mod.$$
When we localize along quasi-isomorphisms, we get a functor $\s$-functor
$$\egrdg_k \longrightarrow \QCoh(B\cH_0).$$
By observation, this $\s$-functor is also endowed with a canonical 
symmetric monoidal structure. In general localizations does not preserve limits, and tere are 
no obious reasons why the previouys $\s$-functor is an equivalence. However
the following result is proven in \cite[Prop. 4.2.3]{mrt}. 

\begin{prop}\label{pI-4}
The $\s$-functor constructed above induces a symmetric monoidal equivalence
$$\egrdg_k \simeq \QCoh(B\cH_0).$$
\end{prop}

A simple corollary of the above Proposition, combined with the rectification results of 
\cite[Prop. 4.5.4.6]{HA}, 
is the following interpretation of graded mixed cdga's.

\begin{cor}\label{cpI-4}
There exists an natural equivalence of $\s$-categories
$$\egrcdga_k \simeq \mathbf{CAlg}(\QCoh(B\cH_0)),$$
where the right hand side is the $\s$-category of commutative algebras in 
$\QCoh(B\cH_0)$ in the sense of \cite{HA}.
\end{cor}

\textit{Proof of the Corollary.} The equivalence of Proposition \ref{pI-4} being 
symmetric monoidal, we get an equivalence on the corresponding 
$\s$-categories of commutative algebra objects
$\mathbf{CAlg}(\egrdg_k) \simeq \mathbf{CAlg}(\QCoh(B\cH_0)).$
To finish the proof we invoke the rectification result of \cite[Prop. 4.5.4.6]{HA},
which can be applied here as $k$ is assumed of characteristic zero (and thus
$\epsilon-C(k)^{gr}$ is a \emph{well powered} symmetric monoidal model category).
We thus get that the natural $\s$-functor
$L(\mathbf{CAlg}(\epsilon-C(k)^{gr}))=\egrcdga_k \to \mathbf{CAlg}(\egrcdga_k)$ is an equivalence. Combining these
two equivalences of $\s$-categories finishes the proof of the Corollary. \hfill $\Box$\\

\begin{rmk}\label{geomintofforgetfulfromepsilongrtogr} 
\emph{Using Proposition \ref{pI-4} and the group section $s: \mathbb{G}_m \to \mathcal{H}$ introduced at the beginning of this subsection, we may realize geometrically the forgetful functor $\egrdg_k \to \grdg_k $ as the pullback on quasi-coherent complexes along the map $Bs: B\mathbb{G}_m \to B\mathcal{H}$.  }
\end{rmk}

The previous description allows for the following easy generalization where $\Spec \, k$ is replaced by an arbitrary derived (pre)stack $X$

\begin{equation}\label{grmixoverX}\egrdg_X := \QCoh(B\cH_0 \times X) \qquad \qquad \egrcdga_X \simeq \mathbf{CAlg}(\QCoh(B\cH_0 \times X))\end{equation}

\noindent (where the indicated products are fibered products over $\Spec\, k$). This equivalence for generally
derived stacks follows directly from corollary \ref{cpI-4} by descent. 

\section{The classifying spaces of graded mixed structures}
\label{sec:Theclassifyingspacesofgradedmixedstructures}

In this section we study the \emph{space of graded mixed structures} on a \emph{fixed} 
graded cdga. We first provide a precisely definition of such a space. 

We consider the forgetful $\s$-functor
$$\mathsf{oub}_\epsilon : \egrcdga_k \to \grcdga_k$$
obtained by forgetting the mixed structure. As $\mathsf{oub}_\epsilon$ is conservative, the fiber
$\s$-category $oub_{\epsilon}^{-1}(B)$ taken at any point $B \in \grcdga_k$ is an $\s$-groupoid, whose
geometric realization will be denoted by $|oub_{\epsilon}^{-1}(B)|$.

\begin{df}\label{dI-6}
For a given graded cdga $B \in \grcdga_k$, \emph{the classifying space of
compatible graded mixed structures on $B$} is the simplicial set 
$$\M_\epsilon(B):=|\mathsf{oub}_{\epsilon}^{-1}(B)| \in \TT.$$
\end{df}

Such spaces of algebraic structures have been heavily studied in the literature. We start by recalling
the well known fact that they can be computed, in our specific case, as mapping spaces 
between graded dg-Lie algebras. For this, let $B \in \grcdga_k$ be a fixed graded cdga over $k$. 
We can form $Der^{gr}(B)$ the graded dg-Lie algebra of graded derivations from $B$ with values in $B$. 
It can be define using a strict cofibrant model $\widetilde{B}$ of $B$, and considering strict graded
derivations of $\widetilde{B}$ into $\widetilde{B}$ as a model for $Der^{gr}(B)$. More explicitly, the
piece of weight $n$ of $Der^{gr}(B)$ is the complex of derivations $\widetilde{B} \to \widetilde{B}(n)$, and the 
Lie bracket is the usual commutator of derivations. The space $\M_\epsilon(B)$ can then be understood
as a mapping space in the $\s$-category $\grdglie_k$ of graded deg-Lie algebras (i.e. Lie algebra object
in $\grdg_k$).

\begin{prop}\label{pI-5}
There is a canonical equivalence $Map_{\grdglie}(k(1)[1],Der^{gr}(B)) \to \M_\epsilon(B)$.
\end{prop}

\textit{Proof.} The proposition is a particular case of a more general result. Let $L$ be any
graded dg-Lie algebra. We have a symmetric monoidal model category $L-Mod^{gr}$ of graded $L$-dg-modules
(for which equivalences and fibrations are defined on the underlying graded complexes), and a corresponding
symmetric monoidal $\s$-category $L-\grdg_k$. Commutative algebras in $L-\grdg_k$ form a new 
$\s$-category $L-\grcdga_k$, whose objects can be called \emph{graded $L$-cdga's}.
It comes equipped with a forgetful $\s$-functor
$$\mathsf{oub}_L : L-\grcdga_k \to \grcdga_k.$$
We claim that, for any fixed $B \in \grcdga_k$, the two spaces
$$\mathsf{oub}_L^{-1}(B) \qquad \Map_{\grdglie_k}(L,Der^{gr}(B))$$
are functorially equivalent. The proposition would follow from the special case where $L=k(1)[1]$,
as $k(1)[1]-\grcdga_k$ is nothing else that the $\s$-category of graded mixed cdga's.

We fix $B \in \grcdga_k$, and chose cofibrant model for it. We thus have two $\s$-functors
$$L \mapsto \mathsf{oub}_L^{-1}(B) \qquad L \mapsto \Map_{\grdglie_k}(L,Der^{gr}(B)).$$
These $\s$-functor can be represented as strict $\s$-functors from the category of cofibrant 
graded dg-Lie algebras to the category of spaces. The first one sends $L$ to the nerve of 
the category $C(L,B)$ whose objects are pairs $(B',u)$, with $B' \in L-\grcdga_k$ and 
$u : B \to \mathsf{oub}_L(B')$ is a quasi-isomorphism. The second functor sends
$L$ to the simplicial set $Hom(c^*(L),Der^{gr}(B))$, where $c^* : \grdglie_k \to (\grdglie_k)^{\Delta}$ 
is a left framing in the sense of \cite[Def. 5.2.7]{hov}.

We can construct a natural transformation between these two $\s$-functors. Indeed, there is a canonical
map from the set $Hom(L,Der^{gr}(B))$ to the nerve of $C(L,B)$, sending a morphism $\rho : L \to Der^{gr}(B)$
to the pair $(B_\rho,id)$, where $B_\rho$ is the object in $L-\grcdga_k$ defined by the moprhism $\rho$
(and whose underlying object equals $B$). This map is functorial in $L$, and thus it defines a morphism of
simplicial sets
$$Hom(c^*(L),Der^{gr}(B)) \to |[n] \to N(C(c^n(L),B))|.$$
Now, as quasi-isomorphisms  between graded dg-Lie induces
Quillen equivalences between the corresponding $\s$-categories of graded cdga's, 
we have that all the morphism in the simplicial diagram $[n] \mapsto N(C(c^n(L),B))$
are weak equivalences, and thus the natural morphism 
$$N(C(L,B)) \to  |[n] \to N(C(c^n(L),B))|$$
is also a weak equivalence. We have thus defined a well defined natural transformation of $\s$-functors
$$\psi_L : Hom(c^*(L),Der^{gr}(B)) \to |[n] \to N(C(c^n(L),B))| \simeq N(C(L,B)).$$

It remains to show that $\psi_L$ is an equivalence for all $L$. For this we use the following steps.
\begin{enumerate}
    \item The induced map 
    $$\psi_L : \pi_0(Hom(c^*(L),Der^{gr}(B))) \to \pi_0(N(C(L,B)))$$
    is bijective for all $L$.
    \item Both $\s$-functors $L \mapsto Hom(c^*(L),Der^{gr}(B))$ and 
    $L \mapsto N(C(L,B))$ sends colimis (in the argument $L$) to limits.
\end{enumerate}

As a first observation, setps $(1)$ and $(2)$ above finishes the proof of the proposition. Indeed, 
we can apply $(1)$ to graded dg-Lie algebras of the form $K \otimes L$, for $K$ an arbitrary simplicial set, 
which by $(2)$ will imply that the morphism $\psi_L$ induces a bijection
$$[K,Hom(c^*(L),Der^{gr}(B))] \to [K,N(C(L,B))].$$
As this is true for all $K$, Yoneda implies that $\psi_L$ is an isomorphism in $Ho(\Top)$ and thus is an 
equivalence.

The point $(1)$ is a direct exercise in model category theory, and is left to the reader. For point $(2)$, 
we use that if $L = colim_i L_i$ is a colimit in the $\s$-category $\grdglie_k$, then the natural 
$\s$-functor
$$L-\grcdga_k \to lim_i L_i-\grcdga_k$$
is an equivalence of $\s$-categories. Indeed, this can be seen easily by first considering the case
of graded dg-modules
$$L-\grdg_k \to lim_i L_i-\grdg_k.$$
The fact that this $\s$-functor is an equivalence follows directly from the equivalence between 
dg-modules over $L$ and dg-modules over the universal enveloping algebra $U(L)$, and the fact that 
$L \mapsto U(L)$ preserves colimits. Indeed, for a colimit in graded dg-algebras $U = colim_i U_i$, 
the $\s$-category of graded $U$-dg-modules is equivalent to the limit of 
the $\s$-categories of graded $U_i$-dg-modules as shown in the next lemma.

\begin{lem}
Let $U = colim_i U_i$ be a colimit in the $\s$-category $\grdga_k$ graded (associative and unital)
dg-algebras over $k$. The canonical $\s$-functor on $\s$-categories of  graded dg-modules
$$U-\grdg_k \to lim_i U_i-\grdg_k$$
is an equivalence.
\end{lem}

\textit{Proof of the lemma.} We have a fibration sequence
$$\Map_{\grdga_k}(U,T(x,x)) \to \Map_{\mathbf{grdgcat_k}}(BU,T) \to \Map_{\mathbf{grdgcat_k}}(\mathbb{1},T),$$
where here $T$ is a graded dg-category, the fiber is taken over the unique morphism $\mathbb{1} \to T$
point an object $x \in T$, and $BU$ is the graded dg-category with a unique object determined by $U$. This is the graded version of the fibration sequence appearing in 
the proof of \cite[Lem. 2.11]{tova}, and can proven the exact same manner using the homotopy theory of
graded dg-categories. This fibration sequence implies the lemma as $\Map(_{\grdga_k}(-,T(x,x))$
transforms colimits in limits.
\hfill $\Box$ \\

The lemma finishes the proof of the point $(2)$ above and thus of the proposition.
\hfill $\Box$ \\ 

Proposition \ref{pI-5} is useful as the mapping space
$Map_{\grdglie}(k(1)[1],Der^{gr}(B))$ can be described using an explicit cofibrant model 
for the graded dg-Lie algebra $k(1)[1]$, introduced in \cite{cptvv}. We let 
$\cL$ be the dg-Lie algebra freely generated by elements $p_i$, for $i>0$, such that 
$p_i$ sits in cohomological degree $1-2i$ and in weight $i$, with the relations
$$d(p_n) + \sum_{a+b=n}[p_q,p_b]=0.$$
We consider the natural projection $\cL \to k(1)[1]$ sending $p_1$ to $\epsilon$ and $p_i$ to 
zero for $i>1$. As explained in \cite[]{ptvv}, this projection is a trivial fibration. As
$\cL$ is clearly cofibrant (as being a cell object), this shows that $\cL$ is a cofibrant model
for $k(1)[1]$. Therefore, we have a diagram of spaces
$$\xymatrix{
Hom_{\grdglie}(\cL,Der^{gr}(B)) \ar[r] & Map_{\grdglie}(\cL,Der^{gr}(B)) & Map_{\grdglie}(k(1)[1],Der^{gr}(B))
\ar[l]_-{\simeq }}$$
where $Hom_{\grdglie}(\cL,Der^{gr}(B))$ is the set of morphisms of graded dg-Lie algebras (which 
depends on the choice of a cofibrant model for $B$). Applying $\pi_0$ provides 
a natural surjective map
$Hom_{\grdglie}(\cL,Der^{gr}(B)) \to \pi_0(Map_{\grdglie}(k(1)[1],Der^{gr}(B)))$. Combined with 
the Proposition \ref{pI-5} we obtain the following corollary.

\begin{cor}\label{cpI-5}
Let $B$ be a fixed cofibrant graded cdga. Then, there exists a natural surjective map
$$Hom_{\grdglie}(\cL,Der^{gr}(B)) \to \pi_0(M_\epsilon(B)).$$
\end{cor}

The previous results can be made a bit more explicit. Indeed, by the definition 
of $\cL$ we have that an element in $Hom_{\grdglie}(\cL,Der^{gr}(B))$ is given by a sequence
of $k$-linear graded derivations
$$d_i : B \to B(i)[1-2i]$$
satisfying the system of equations
$$d(d_n) + \sum_{a+b=n}[d_q,d_b]=0.$$
Suppose moreover that $B$ is explicitly given by $Sym_A(E[1])$, where $A$ is a cofibrant
connective cdga, and $E$ is a cofibrant perfect $A$-dg-module. Then, 
each $d_i$ consists itself as a pair
$$\alpha_i : E \longrightarrow \wedge^{i+1}_A E[1-i] \qquad 
\beta_i : A \longrightarrow \wedge^{i}_A E[1-i],$$
such that $\beta_i$ is a $k$-linear derivation, and $\alpha_i$ is a morphism which 
is $A$-linear with respect to $\beta_i$
$$\alpha_i(a.m)=a\alpha_i(m) + (-1)^? \beta_i(a)\wedge m$$
(where $?$ is the adequate sign). This description provides an explicit 
way to write down elements in the set $Hom_{\grdglie}(\cL,Der^{gr}(B))$, and thus
also an explicit manner to write down all the elements in $\pi_0(M_\epsilon(B))$
by Corollary \ref{cpI-5}. We also refer to \cite{mon2} for another
closely related, but slightly different, way to 
describe these classifying spaces of graded mixed structures.

\chapter{Derived foliations}\label{chapter:Derivedfoliations}

In this chapter we give the basic definitions and study the general properties of derived
foliations. We provide examples, and a comparison with the theory of dg-Lie algebroids.
One of the main result, in the third section,
states that derived foliations are always formally integrable, which is a major difference
with the situation of underived possibly singular foliations. We provide two versions of this
result, a general version for derived foliations on general derived Artin stacks, and 
a more specific version for a class of derived foliations called \emph{quasi-smooth and rigid}
on smooth varieties. This second version is useful to study existence of derived enhancement
of underived singular foliations, and to relate them to the classical notion of Godbillon-Vey sequences. \\

As before, unless other is specified, all structures (algebras, 
modules, schemes \dots) are over our base commutative
ring $k$. We assume that $k$ is furthermore a $\QQ$-algebra. The notion 
of derived foliations in arbitrary characteristics situations will be treated
in a separate later chapter of this book, as it requires some
extra care and more evolved notions (see \S \ref{chapter:nonzerocharacteristics}). 

\section{First definitions and examples}\label{sec:Firstdefinitionsandexamples}

We are now ready to define the notion of derived foliations. We will 
first consider the situation for an affine derived scheme, and then proceed
by gluing in order to define derived foliations on more global objects such
as derived schemes and derived Artin stacks. \\

\subsection{Definition in the affine setting}

Let $X = \Spec\, A$ be a derived affine scheme over $k$. Recall that this
implicitly means that $A$ is a connective cdga over $k$, and thus it 
is cohomologically concentrated in non-positive degrees (see appendix \S \ref{secapp:Derivedschemesandstacks}). 
We will allow ourselves to 
implicitly identify the $\s$-categories of $A$-dg-modules and of
quasi-coherent complexes of $\OO_X$-modules. We will also 
use the derived de Rham algebra $\DR(A/k)$, that will also be denoted 
by $\DR(X/k)$.

A derived foliation $\F$ on $X$ will be given, by definition, a graded mixed
cdga $\DR(\F)$, together with a morphism $u : \DR(A/k) \to \DR(\F)$ of graded mixed cdga's over $k$,
satisfying some extra conditions. In order
to phrase these conditions, note that the morphism $u$ induces
a morphism of graded cdga's
$$A \to \DR(A/k) \to \DR(\F),$$
where the first of this morphism is the inclusion of the weight $0$ part. 
This endows $\DR(\F)$ with a canonical structure of a graded $A$-linear
cdga structure. In particular, the inclusion of the weight $1$ part
$\DR(\F)^{(1)} \subset \DR(\F)$ induces a canonical morphism of graded 
$A$-linear cdga's
$$\phi_u : Sym_{A}(\DR(\F)^{(1)}) \longrightarrow \DR(\F).$$
Note that, the morphism $\phi_u$ above, of course, depends on the original 
morphism $u$.

\begin{df}\label{dII-1}
Let $X=\Spec\, A$ be a derived affine scheme over $k$. A \emph{derived
foliation $\F$ over $X$ relative to $k$} is a pair $(\DR(\F),u)$, 
which consists of 
a graded mixed cdga $\DR(\F)$, and a morphism of graded
mixed cdga's
$$u : \DR(A/k) \longrightarrow \DR(\F),$$
such that the induced morphism
$$\phi_u : Sym_{A}(\DR(\F)^{(1)}) \longrightarrow \DR(\F)$$
is a quasi-isomorphism.
\end{df}

As a first comment, it is important to note that 
the condition that $\phi_u$ being a quasi-isomorphism includes in particular
the condition that the morphism $u$ induces a quasi-isomorphism in weight $0$
$A \simeq \DR(\F)^{(0)}.$
As a second observation, the universal property of the
derived de Rham algebra $\DR(A/k)$ (see Definition \ref{dI-4}), the morphism
$u$ is completely characterized by the morphism of graded cdga's (without
any mixed structure)
$$A \to \DR(A/k) \to \DR(\F).$$
Therefore, the definition of derived foliations above could also be phrased 
as follows: a derived foliation $\F$ on $X$ over $k$ consists of 
a graded mixed cdga $\DR(\F)$, together with a quasi-isomorphism
of cdga's $A \to \DR(\F)^{(0)}$, and such that the canonical 
morphism 
$$ Sym_{A}(\DR(\F)^{(1)}) \longrightarrow \DR(\F)$$ is a quasi-isomorphism. \\

The idea behind Definition \ref{dII-1} is that $\DR(\F)$ is the 
de Rham algebra of forms \emph{along the leaves} of $\F$. Therefore, 
this $A$-linear dg-module $\DR(\F)^{(1)}$ should be understood
as the (shifted by $-1$) cotangent sheaf, or rather cotangent complex,
classifying \emph{differential forms along the leaves of $\F$.} This
complex is one of the fundamental data contained in a derived foliation, 
and we will therefore give a specific name, and introduce a specific 
notation for it.

\begin{df}\label{dII-2}
Let $\F$ be a derived foliation on a derived affine scheme $X=\Spec\, A$.
The \emph{cotangent complex of $\F$ relative to k} is defined to be
$$\LL_{\F/k} := \DR(\F)^{(1)}[-1].$$
When the ground ring $k$ is unambiguous it is simple denoted by 
$\LL_{\F}$.
\end{df}

An important remark is that the cotangent complex $\LL_{\F/k}$ of 
a derived foliation $\F$ always comes equipped with a canonical
morphism of $A$-dg-modules
$$\LL_{A/k} \longrightarrow \LL_{\F/k},$$
which is induced by the (shifted) weight $1$ part of the morphism $u : \DR(A/k) \to \DR(\F)$.
This is often referred to as the \emph{anchor map}.

Using the notation and the terminology above, a derived foliation $\F$ 
on $X$, can be thought as a graded mixed structure $\epsilon$ on 
$Sym_A(\LL_{\F/k}[1])$. This graded mixed structure is itself the
de Rham differential on differential forms along $\F$. However, 
one important point in Definition \ref{dII-1} is that
$\epsilon$ does not strictly exists on $Sym_A(\LL_{\F/k}[1])$, but rather
on $\DR(\F)$ which is only quasi-isomorphic $Sym_A(\LL_{\F/k}[1])$. Trying 
to transport the mixed structure along this quasi-isomorphism provides
a derivation $\epsilon$ on $Sym_A(\LL_{\F/k}[1])$, but 
$\epsilon^2$ will now only be homotopic to zero (as opposed to strictly equal
to zero). We can use the quasi-isomorphism to show that there exists a canonical
such homotopy, which itself should satisfy some higher coherences and so on and
so forth. Despite its naive form, the reader must be warned that 
Definition \ref{dII-1} hides in fact an infinite number of rather
evolved homotopy coherences, very similar to 
the one being found in shifted symplectic and Poisson structures of \cite{cptvv}. 
We will come back to this point later on in this book and will 
make these coherence explicit. The reader can already notice that 
by construction $\epsilon$ is a point in the moduli of compatible graded mixed structures
on $Sym_A(\LL_{\F/k}[1])$ in the sense of Definition \ref{dI-6}, and thus can be described
by the coherence of homotopies given by corollary \ref{cpI-5}.

Keeping this in mind, a derived foliation $\F$ has a
cotangent complex $\LL_{\F/k}$, which comes equipped with
de Rham differentials 
$$A \to \LL_{\F/k} \to \wedge^2_A \LL_{\F/k} \to \dots$$
These differential satisfy all the required formula (such as
$d^2=0$, multiplicativity ect), but only up to natural homotopy coherences. 
As far as we know, it is in general not possible to make these identities
to hold strictly on the noze. \\

Before showing several basic examples of derived foliations, 
we introduce some terminology. Derived foliations will be
characterized according to the \emph{size} of their cotangent 
complexes. 

\begin{df}\label{dII-3}
Let $\F$ be a derived foliation on a derived affine scheme $X=\Spec\, A$.
\begin{enumerate}
\item The derived foliation $\F$ is \emph{perfect} (resp. \emph{almost perfect}) if 
$\LL_{\F/k}$ is a perfect (resp. \emph{almost perfect}) $A$-dg-module.

\item The derived foliation $\F$ is \emph{n-connective}, for some
integer $n\in \ZZ$, if $H^i(\LL_{\F/k})=0$ for all $i >n$
(i.e. if $\LL_{\F/k}[n]$ is connective).

\item The derived foliation $\F$ is \emph{smooth} if it is perfect and if
$\LL_{\F/k}$ is of Tor amplitude contained in $[0,\s[$.

\item The derived foliation $\F$ is \emph{quasi-smooth} if it is perfect and if
$\LL_{\F/k}$ is of Tor amplitude contained in $[-1,\s[$.

\end{enumerate}

\end{df}

\subsection{Basic examples of derived foliations}\label{subsec:basicexamplesofderivedfoliations}

We gather in this subsection the very first examples of derived foliations. More evolved 
examples will be given later all along this book. \\

\begin{ex}\label{ex:classicalfoliations}
\emph{\textbf{Classical foliations on smooth 
affine schemes.} We start by the most basic examples of how usual foliations
define derived foliations in our sense. For this, we let 
$X=\Spec\, A$ be a smooth affine scheme over $k$. Let 
$V \subset T_X$ be a subbundle of the tangent bundle of $X$, 
which is stable by the Lie bracket of $T_X$. We consider the dual
bundle $V^\vee$ and 
$Sym_A(V^\vee[1])$. The inclusion $V \subset T_X$ defines
a quotient map $\Omega_{X/k}^1 \to V^\vee$, and thus a morphism of graded cdga's (with zero 
cohomological differential)
$$\DR(A/k)=Sym_A(\Omega_{A/k}^1[1]) \longrightarrow Sym_A(V^\vee[1]).$$
As $V$ is stable by the Lie bracket of vector fields, 
the de Rham differential $dR : A \to \Omega_{A/k}^1$ descends to 
a de Rham differential $dR : A \to V^\vee$, which can be extended multiplicatively
to a graded mixed structure on the graded mixed cdga $Sym_A(V^\vee[1])$. 
We thus have defined a morphism of graded mixed cdga's
$$\DR(A/k) \longrightarrow Sym_A(V^\vee[1]).$$
This morphism defines a derived foliation $\F$ on $X$ relative to $k$
according to the Definition \ref{dII-1}, such that 
$\DR(\F)=Sym_A(V^\vee[1])$ endowed with its induced de Rham differential.
This derived foliation has $V^\vee$ as its cotangent complex, and thus
is smooth in the sense of Definition \ref{dII-3}.} \end{ex}

\begin{ex}\label{ex:liealgebroids}
\emph{\textbf{Lie algebroids on smooth affine schemes.} The previous
example possesses a slight generalization. We let $X=\Spec\, A$ be
a smooth affine scheme over $k$. Let now $V$ be a vector bundle
on $X$ endowed with the structure of a Lie algebroid 
(see \S \ref{sec:ComparisonwithderivedLiealgebroids} for more
on the relations to Lie algebroids). We remind that 
this consists of a pair $(a,[,])$, where 
$a : V \to T_X$ is a morphism of vector bundles on $X$, and
$[,]$ is a $k$-linear Lie bracket on sections of $V$. These data are furthermore 
assumed to 
satisfy two compatibility conditions. 
First we have the standard Liebniz rule: for any sections $s$ and $t \in V$ and
any function $f\in A$ on $X$, we have
$$[f.s,t]=f.[s,t]+a(t)(f).s,$$
where $a(t)(f)$ is the value of the derivation $a(t) \in T_X$ on the function $f$.
The second condition states that the morphism $a$ is compatible with the Lie
brackets on $T_X$: $a([s,t])=[a(s),a(t)]$ for any pair of section $(s,t)$ of $V$.
Classical foliations on $X$ are therefore Lie algebroids $(V,a,[,])$ 
such that $a$ identifies $V$ with a subbundle of $T_X$.} 

\emph{For a Lie algebroid $(V,a,[,])$ as above, we can form its de Rham
algebra $Sym_A(V^\vee[1])$. This is a graded cdga which comes equipped with
an extension of the de Rham differential
$$\xymatrix{
A \ar[r] &  V^\vee \ar[r] & \wedge_A^2 (V^{\vee}) \ar[r] & \dots \ar[r] & \wedge^r_A(V^\vee) \ar[r] & \dots.}$$
The morphism $A \to V^\vee$ is given by composing the univeral
derivation $A \to \Omega_{A/k}^1$ with the $A$-linear dual of the map 
$a$. The differential $dR : V^\vee \to \wedge_A^2 (V^{\vee})$ is itself
dual to the Lie bracket $[,]$ on $V$ in the following sense: for 
$a$ and $b$ sections of $V$, and $v\in V^\vee$, we have
$$dR(v)(s\wedge t)=a(s)(v(t))-a(t)(v(s))+v([s,t]).$$
This differential extends to a canonical graded mixed structure on 
$Sym_A(V^\vee[1])$, and defines a graded mixed cdga, which we denote by 
$\DR(V)$. By construction,
the $A$-linear dual of $a$ defines a morphism of graded mixed cdga's
$\DR(A/k) \to \DR(V)$. The conditions of Definition \ref{dII-1}
are satisfied and therefore this construction defines a derived foliation on $X$.
We refer to \S \ref{sec:ComparisonwithderivedLiealgebroids} for a more general construction
as well as for a precise comparison between derived foliations and objects of Lie-types.} \end{ex}

\begin{ex}\label{ex:tautologicalfoliations}
\emph{\textbf{Tautological derived foliations.} For a given derived affine scheme 
$X=\Spec\, A$ there are two tautological derived foliations on $X$ of particular
importance. The first one is called the \emph{$0$-foliation}, or 
the \emph{pointwise foliation}, as morally its leaves are points of $X$.
In terms of graded mixed cdga's it is given by the augmentation morphism
$$\DR(A/k) \longrightarrow A,$$
inducing the identity on $A$ and zero on all piece of higher weights. Here $A$ 
is considered as a graded mixed cdga with zero mixed structure, and $A$ is
considered in weight $0$. This augmentation morphism obviously satisfies the
conditions of Definition \ref{dII-1}, and here the cotangent complex
is simply given by $0$. This derived foliation will be denoted by $0_X$, 
and is the initial object of the $\s$-category 
of derived foliations on $X$ (see \S \ref{cpII-5} for a more general statement).}

\emph{The second tautological foliation is the extreme opposite, and is
called the \emph{de Rham foliation}, because of its direct relation 
with the de Rham complex. Morally it has only one leave, namely $X$ itself.
In terms of graded mixed cdga's it is given by the identity morphism
$$\DR(A/k) \longrightarrow \DR(A/k).$$
Here the cotangent complex of this derived foliation 
is $\LL_{X/k}$ and the anchor map is the identity map. It is denoted by $*_X$
and we will see later on that it is the final object in the $\s$-category of
derived foliations over $X$.} \end{ex}

\begin{ex}\label{ex:derhamfoliations}
\emph{\textbf{Derived foliations induced by a morphism.}
We finish by a very last basic example, the relative de Rham derived foliation. 
This is again the derived foliation $*_X$ described above, except that the
base ring $k$ is now replaced by a base cdga $B$ mapping to $A$.\\
Let $f : X=\Spec\, A \longrightarrow Y=\Spec\, B$ be a morphism of derived affine 
schemes. This morphism induces a morphism of cdga $B \to A$ and thus of
graded mixed cdga's $\DR(B/k) \to \DR(A/k)$. We define the graded mixed
cdga's $\DR(A/B)$, of $A$ relative to $B$, by
$$\DR(A/B):=\DR(A/k) \otimes_{\DR(B/k)}B.$$
The tensor product above is taken inside the $\s$-category $\egrcdga_k$ of graded
mixed cdga's over $k$, where the morphism $\DR(B/k) \to B$ is the augmentation
morphism. The natural induced morphism 
$\DR(A/k) \to \DR(A/B)$ defines a derived foliation on $X$ in the sense
of Definition \ref{dII-1}. The cotangent complex of this derived foliation
is here $\LL_{A/B}$ the relative cotangent complex of $A$ over $B$, which
sits in a exact triangle of $A$-dg-modules
$$\xymatrix{\LL_{B/k}\otimes_{B}A \ar[r]& \LL_{A/k} \ar[r] & \LL_{A/B}.}$$
This is the relative de Rham foliation on $X$ over $Y$, and 
is denoted by $*_{X/Y}$. This is the very basic example of algebraically integrable
derived foliations, and we will see later on that it is the pull-back along $f$
of the tautological pointwise foliation $0_Y$ (see \S \ref{subsec:Globaldefinition}).} \end{ex} 

\subsection{Global definition}\label{subsec:Globaldefinition}

We now turn to the global notion of derived foliations on derived Artin stacks
in the sense of \S \ref{secapp:Derivedschemesandstacks}. 
For this, we will use the standard procedure by performing a 
left Kan extension from derived affine schemes to all derived stacks. \\

We start by defining an $\s$-functor
$$\Fol(-/k) : \dAff_k^{op} \longrightarrow \scat,$$
from the $\s$-category of derived affine schemes to the opposite of the $\s$-category of
$\s$-categories. For an affine derived scheme $X=\Spec\, A$, 
$F(X/k)$ is defined to be the $\s$-category of 
graded mixed $\DR(X/k)$-cdga's which satisfies the conditions of Definition 
\ref{dII-1}. For an object $\F \in \Fol(X/k)$, we will 
use the notation $\DR(\F)$ for the corresponding graded mixed $\DR(X/k)$-cdga.

For a morphism of derived affine schemes $f : Y=\Spec\, B \to X=\Spec\, A$,
we define an $\s$-functor
$f^* : \Fol(X/k) \to \Fol(Y/k)$, by sending a graded mixed $\DR(X/k)$-cdga  
$\DR(\F)$ (corresponding to an object $\F \in \Fol(Y/k)$) to 
$\DR(Y/k)\otimes_{\DR(X/k)}\DR(\F)$, which is a graded mixed $\DR(Y/k)$-cdga. 
We note here that $\DR(Y/k)\otimes_{\DR(X/k)}\DR(\F)$ also satisfies the conditions
of Definition \ref{dII-1} and thus defines an object $f^*(\F) \in \Fol(Y/k)$:
if as a graded $A$-linear cdga we have 
$\DR(\F) \simeq Sym_A(\LL[1])$, 
then $\DR(Y/k)\otimes_{\DR(X/k)}\DR(\F)$ is of the form
$Sym_B(\LL'[1])$, where $\LL'$ sits in a push-out diagram of $B$-dg-modules
$$\xymatrix{
B\otimes_A \LL_{A/k}[1] \ar[r] \ar[d] & B\otimes_A \LL[1] \ar[d] \\
\LL_{B/k}[1] \ar[r] & \LL'[1].
}$$
Indeed, after applying the $Sym_B$ $\s$-functor, the above push-out diagram of
$B$-dg-modules is transformed to the push-out diagram of graded cdgas defining 
$\DR(Y/k)\otimes_{\DR(X/k)}\DR(\F)$
$$\xymatrix{
B\otimes_A \DR(A/k) \ar[r] \ar[d] & B\otimes_A \DR(\F) \ar[d] \\
\DR(B/k) \ar[r] & Sym_B(\LL'[1]).
}$$

We thus have defined an $\s$-functor
$$\Fol(-/k) : \dAff_k^{op} \longrightarrow \scat,$$
which sends $X=\Spec\, A$ to $\Fol(X/k)$ and
$f : Y=\Spec\, B \to X=\Spec\, A$ to $f^* : \Fol(X/k) \to \Fol(Y/k)$. It is important to 
notice here that if $\F \in \Fol(X/k)$ is a derived foliation with
$\LL_F$ as a cotangent complex, then $f^*(\F)$ is a derived foliation on $Y$ with 
a cotangent complex $\LL_{f^*(\F)}$ defined by the push-out of $B$-dg-modules
$$\xymatrix{
B\otimes_A \LL_{A/k}[1] \ar[r] \ar[d] & B\otimes_A \LL_\F[1] \ar[d] \\
\LL_{B/k}[1] \ar[r] & \LL_{f^*(\F)}[1].
}$$
In particular, even thought there exists a canonical map $f^*(\LL_\F) \to \LL_{f^*(\F)}$,
of quasi-coherent complexes on $Y$,
this map is not an equivalence in general (unless the map $f$ is itself formally etale). In short: cotangent
complexes of derived foliations are not stable by pull-backs. This fact will make more complicated the 
definition of a global cotangent complex when working with global derived stacks, as the gluing procedure
can not be straightforward.
Note that we
will also sometimes use the notation $\Fol(A/k)$
for $\Fol(X/k)$ (when $X=\Spec\, A$). \\

Our first observation is that the $\s$-functor $\Fol(-/k)$ is a derived stack
for the \'etale topology, showing that derived foliations is a local notion
for the \'etale topology.

\begin{prop}\label{pII-1}
The $\s$-functor $\Fol(-/k) : \dAff_k^{op} \longrightarrow \scat$
is a (hypercomplete) stack for the \'etale topology.
\end{prop}

\textit{Proof.} We start by noticing that $X \mapsto \DR(X/k)$ 
defines a hypercomplete stack for the \'etale topology on $\dAff_k^{op}$. 

\begin{lem}\label{lII-1}
The $\s$-functor $\DR(-/k) : \dAff^{op}_k \to \egrcdga_k$ is a
hypercomplete stack for the \'etale topology.
\end{lem}

\textit{Proof of the lemma.} We consider the forgetful $\s$-functor
$G : \egrcdga_k \longrightarrow \grdg_k$, which forgets the algebra and the mixed structures
and only retains the underlying graded complex. This $\s$-functor is conservative and
commutes with all limits. Therefore, in order to prove the lemma it is
enough to check that the composed $\s$-functor
$$G \circ \DR(-/k) : \dAff^{op}_k \to \grdg_k$$
is a hypercomplete stack. This last statement is equivalent to the fact that for all $n \in \ZZ$,
the weight $n$ part $X \mapsto \DR(X/k)^{(n)}$ is a hypercomplete stack of complexes.
Moreover, this last $\s$-functor sends, by definition, a derived affine scheme $X=\Spec\, A$
to the complex $\wedge_A^n \LL_{A/k}[n]$ (where $\wedge^n$ is declared zero for $n<0$ by 
convention).

For fixed $n \in \ZZ$, and fixed $X=\Spec\, A$, the
$\s$-functor $(\Spec\, A' \to X) \mapsto \wedge_{A'}^n \LL_{A'/k}[n]$
is a quasi-coherent complex of $\OO_X$-modules when restricted to the small
\'etale site of $X$. Indeed, for $\Spec\, A' \to \Spec\, A$ any \'etale morphism, 
the canonical morphism of $A'$-dg-modules
$$\LL_{A/k} \otimes_A A' \to \LL_{A'/k}$$
is a quasi-isomorpism, and thus so are all the induced morphisms on wedge powers
$$(\wedge_{A}^n\LL_{A/k}) \otimes_A A' \to \wedge_{A'}^n\LL_{A'/k}.$$
However, it is well know that quasi-coherent $\OO_X$-modules satisfies hypercomplete
\'etale descent (see \cite[Prop. 2.2.2.13]{hagII}). We thus have proven that $G \circ \DR(-/k)$ is indeed
an \'etale hyperstack, and thus have proven the lemma. \hfill $\Box$ \\

To prove that $\Fol(-/k)$ is a derived stack we use the criterion from \cite[Cor. 1.3.2.4]{hagII}, and thus
we have to check that $\Fol(-/k)$ preserves finite products and 
\'etale hypercoverings of connective cdga's. \\

\textbf{Compatibility with finite products.} Let $A$ be a connective cdga and suppose that 
it can be written as a product of cdga's $A=\prod_{i\in I}A_i$, with $I$ a finite set. 
We can write $X=\Spec\, A \simeq \coprod_i X_i$ with $X_i=\Spec\, A_i$.
We have a canonical $\s$-functor
$$\phi : \Fol(X/k) \longrightarrow \prod_{i\in I}\Fol(X_i/k),$$
simply obtained as the product of the pull-backs $\s$-functors
$j_i^* : \Fol(X/k) \to \Fol(X_i/k)$, where $j_i : X_i \hookrightarrow X$
is the canonical inclusion. The $\s$-functor $\phi$ possesses an inverse defined as
follows. For a family $\F_i \in \Fol(X_i)$, corresponding to graded mixed
$\DR(A_i)$-cdgas $\DR(\F_i)$, we can consider the graded mixed cdga $\prod_{i \in I}\DR(\F_i)$.
It weight $0$ part is $A=\prod_{i\in I} A_i$, and thus it is naturally equipped with 
a structure of a graded mixed $\DR(A/k)$-cdga. Moreover, as each 
graded cdga $\DR(\F_i)$ is of the form $Sym_{A_i}(\LL_i)$, we have an equivalence of graded
cdga's
$$\prod_{i \in I}\DR(\F_i) \simeq Sym_{A}(\prod_{i\in I}\LL_i).$$
Therefore, we have defined an $\s$-functor sending a family $\F_i \in \Fol(X_i)$ 
to the derived foliation on $X$ defined by the graded mixed $\DR(X/k)$-cdga
$\prod_{i \in I}(\DR(\F_i))$. Let us denote by 
$$\psi : \prod_{i\in I}\Fol(X_i/k) \longrightarrow \Fol(X/k)$$
this $\s$-functor.

We claim that $\psi$ and $\phi$ are inverse to each others. Indeed, 
for $\F \in \Fol(X/k)$, the canonical morphism of graded mixed $\DR(X/k)$-cdga's
$$\DR(\F) \longrightarrow \prod_{i\in I}(\DR(X_i/k)\otimes_{\DR(X/k)}\DR(\F))$$
induces a natural morphism $\psi \phi (\F) \to \F$ in $\Fol(X/k)$. Because of the lemma \ref{lII-1} we know that 
$\DR(X/k) \simeq \prod_{i\in I}\DR(X_i/k)$, and thus this natural morphism 
is in fact an equivalence of derived foliations over $X$. This shows that 
$\psi\phi$ is equivalent to the identity $\s$-functor.

In the other direction, for $\{\F_i\} \in \prod_{i\in I}\Fol(X_i/k)$, we have a morphism
$\{\F_i\} \to \phi\psi(\{\F_i\})$ corresponding to the canonical projections
$$(\prod_{i\in I}\DR(\F_i)) \otimes_{\DR(X/k)}\DR(X_{j}/k) \longrightarrow \DR(\F_j)$$
for $j\in I$. Again, as $\DR(X/k) \simeq \prod_{i\in I}\DR(X_i/k)$ this morphism
is an equivalence in $\prod_{i \in I}\Fol(X_i/k)$. This shows that $\phi\psi$ is also
equivalent to the identity $\s$-functor. \\

\textbf{Compatibility with hypercoverings.} Let now 
$A \longrightarrow B_*$ be an etale hypercovering of connective cdga's as in \cite[Cor. 1.3.4]{hagII}.
This defines a coaugmented cosimplicial graded mixed cdga
$\DR(A/k) \longrightarrow \DR(B_*/k)$, which by lemma \ref{lII-1}
induces an equivalence
$\DR(A/k) \simeq \llim_{i \in\Delta}\DR(B_i/k)$. By base change we get 
an $\s$-functor on the $\s$-categories of graded mixed modules
$$\phi : \DR(A/k)-\egrdg \longrightarrow \llim_{i \in\Delta}\DR(B_i/k)-\egrdg.$$
This $\s$-functor possesses a right adjoint
$$\psi : \llim_{i \in\Delta}\DR(B_i/k)-\egrdg \longrightarrow \DR(A/k)-\egrdg.$$
which sends a cosimplicial graded mixed $\DR(B_*/k)$-dg-module $M_*$, to 
$$\psi(M_*)=\llim_{i\in \Delta}M_i \in \DR(A/k)-\egrdg,$$
where each $M_i$ is considered as a graded mixed module over $\DR(A/k)$ via the 
coaugmentation $\DR(A/k) \to \DR(B_i/k)$. 

Let $\F \in \Fol(A/k)$ and let us consider the
adjunction morphism
$$\DR(\F) \longrightarrow \llim_{i\in \Delta}(\DR(B_i/k)\otimes_{DR(A/k)}\DR(\F)).$$
To prove that this morphism is an equivalence, we can forget the mixed and algebra structures, 
as the corresponding forgetful $\s$-functor commutes with limits and is conservative, and thus
only consider the above morphism as a morphism of graded complexes. Moreover, 
as each map $A \to B_i$ is \'etale, we have $\DR(B_i/k) \simeq B_i \otimes_A \DR(A/k)$, 
and thus $\DR(B_i/k)\otimes_{DR(A/k)}\DR(\F)\simeq B_i \otimes_A \DR(\F)$. Therefore,
the above adjunction morphism, when considered as a morphism of graded complexes, 
is equivalent to the natural morphism
$$\DR(\F) \longrightarrow \llim_{i\in \Delta}(B_i\otimes_{A}\DR(\F)).$$
By \'etale descent for modules (see \cite[Lem. 2.2.2.13]{hagII}), this is an equivalence. This shows that 
the counit of the adjunction $\psi\phi \Rightarrow id$ is an equivalence. 

We consider now an object $\F_* \in \llim_{i}\Fol(B_i/k)$, and the unit of the adjunction
$\phi\psi(\DR(\F_*)) \to \DR(\F_*)$, which is
a morphism of cosimplicial graded mixed $\DR(B_*)$-cdga's.
It is a general fact that the projection to the degree $0$ piece
$\llim_{i\in \Delta}\mathcal{C}_i \longrightarrow \mathcal{C}_0$ is conservative 
for any limits of $\s$-categories along $\Delta$. As a consequence, 
to prove that the unit morphism is an equivalence it is enough to show that the natural map
$\phi\psi(\DR(\F_0)) \to \DR(\F_0)$
is an equivalence of graded mixed cdga's. This morphism is given by the canonical map
$$B_0\otimes_A (\llim_{i\in\Delta}\DR(\F_i)) \longrightarrow \DR(\F_0).$$
We can, as above, forget the mixed structures, and consider this as a morphism 
of graded cdga's. But, as a graded cosimplicial $B_*$-cdga, 
$\DR(\F_*)$ is of the form $Sym_{B_*}(\LL_{\F_*})$, where $\LL_{\F_*}$ is
the cosimplicial $B_*$-module defined by the cotangent complex
of $\F_*$. Because all transition maps in $B_*$ are formally \'etale, 
$\LL_{\F_*}$ is a cartesian cosimplicial $B_*$-dg-module, and thus
by \'etale descent (see \cite[Lem. 2.2.2.13]{hagII}) comes from a $A$-dg-module $\LL$. 
Moreover, because the descent equivalence $A-\dg \simeq \llim_i B_i-\dg$
is a symmetric monoidal equivalence, we have that 
the graded $B_*$-dg-module $\DR(\F_*)$ comes 
from $Sym_A(\LL)$ by base change along the coaugmentation 
$A \to B_*$. Moreover, we have
$Sym_A(\LL) \simeq \llim_i(\DR(\F_*))$. In particular, we see that the morphism
$B_0\otimes_A (\llim_{i\in\Delta}\DR(\F_i)) \longrightarrow \DR(\F_0)$ becomes equivalent to
the canonical map
$$B_0 \otimes A Sym_A(\LL) \longrightarrow Sym_{B_0}(B_0 \otimes_A \LL)$$
which is an equivalence. This shows that $\DR(\F_*)$ comes by pull-back from
a graded mixed $\DR(A/k)$-cdga $\DR(\F)$, but also that as a graded $A$-linear cdga $\DR(\F)$ is
of the form $Sym_A(\LL)$. \\

To conclude, we have shown that the adjunction $(\phi,\psi)$ restricts to an equivalence
of $\s$-categories of derived foliations
$$\Fol(A/k) \simeq \llim_{i\in \Delta}\Fol(B_i/k).$$
\hfill $\Box$ \\

The Proposition \ref{pII-1} implies that $\Fol$ can be considered 
itself as a derived stack in $\s$-categories. As such, for any derived stack
$X \in \dSt_k$, we can form the $\s$-category of morphisms from $X$ to $\Fol(-/k)$
\begin{equation}\label{functorFolonderivedstacks}\Fol(X/k):=Map_{\dSt_k}(X,\Fol(-/k)) \in \scat.\end{equation}
In more concrete terms, we have
$$\Fol(X/k) \simeq \llim_{Spec\, A \to X} \Fol(A/k),$$
where the limit is taken over the (opposite) $\s$-category $(\dAff_k/X)^{op}$ of derived affine schemes
over $X$. 

\begin{df}\label{dII-4}
For a derived stack $X$ over $k$, the \emph{$\s$-category of derived foliations of $X$} 
(relative to  $k$) is $\Fol(X/k)$ constructed above. When the base ring $k$ is clear 
from the context, we will simply write $\Fol(X)$ for $\Fol(X/k)$.
\end{df}

When a derived foliation $\F \in \Fol(X/k)$ lives over a derived affine scheme $X=\Spec\, A$, 
its cotangent complex $\LL_{\F}$ is the $A$-dg-module such that 
$\DR(\F) \simeq Sym_A(\LL_{\F}[1])$ as graded cdga's. Using the equivalence between 
quasi-coherent modules over $X$ and $A$-dg-modules (see \cite[Cor. 1.3.7.3]{hagII}), $\LL_\F$ will be considered
an an object in $\QCoh(X)$, the $\s$-category of quasi-coherent complexes over $X$.

When $X$ is a general derived stack, and $\F \in \Fol(X/k)$, the global cotangent 
complex of $\F$ still makes sense but is slightly more complicated to express. 
We first define an $\OO_X$-dg-module as follows. For each $u : \Spec\, A \longrightarrow X$,
we have $u^*(\F)$, the induced derived foliation on $\Spec\, A$, and thus
a quasi-coherent complex $\LL_{u^*(\F)} \in\QCoh(\Spec\, A)$. For 
two morphisms
$$\xymatrix{\Spec\, B \ar[r]^-{f} & \Spec\, A \ar[r]^-{u} & X,}$$
there exists a canonical morphism of quasi-coherent complexes over $\Spec\, B$
$$f^*(\LL_{u^*(\F)}) \longrightarrow \LL_{(uf)^*(\F)}.$$
In general, the above morphism is not an equivalence, unless the map $f$ is formally \'etale. 
Therefore, the rule that associate to $\Spec\, A \to X$ the $A$-dg-module
$\LL_{u^*(\F)}$ defines a complex of $\OO_X$-dg-modules in the sense reminded in \S
\ref{secapp:OXmodulesquasi-coherentandperfectmodules}, 
but which is not quasi-coherent in general. This will be called the 
\emph{big cotangent complex} of $\F$, or the \emph{fake cotangent complex}. We will
denote it by 
$$\LL_{\F}^{big} \in \OO_X-\dg.$$
As reminded  in \S \ref{secapp:OXmodulesquasi-coherentandperfectmodules}, the $\s$-category of quasi-coherent complexes over $X$ sits 
inside $\OO_X-\dg$ as a full sub-$\s$-category, and the natural embedding 
$$\QCoh(X) \hookrightarrow \OO_X-\dg$$
possesses a canonical right adjoint 
$$(-)^{qcoh} : \OO_X-\dg \longrightarrow \QCoh(X)$$
called the \emph{quasi-coherator}. Note that the natural 
morphism $\LL_{X/k}^{big} \to \LL_\F^{big}$ induces
a morphism on the associated quasi-coherent complexes.
The quasi-coherent complex $(\LL_{X/k}^{big})^{qcoh}$ identifies naturally 
with the global cotangent complex $\LL_{X/k}$, as reminded in theorem
\ref{tA-5-1}. Therefore, 
when $X$ is a derived Artin stack and $\F\in \Fol(X/k)$, the 
cotangent complex $\LL_\F$ always comes equipped with a natural map
$$\LL_{X/k} \to \LL_{\F}.$$

\begin{df}\label{dII-5}
For any derived Artin stack $X$ and $\F \in \Fol(X)$, 
the \emph{cotangent complex of $\F$} is defined to be
$$\LL_{\F}:=(\LL_{\F}^{big})^{qcoh} \in \QCoh(X).$$
The \emph{conormal complex of $\F$} is defined to be the fiber 
of the natural morphism
$\LL_X \longrightarrow \LL_\F$, and is denoted by $\mathcal{N}^*_{\F/X}$.
\end{df}

We will see examples of global cotangent complexes in our next paragraph. Before this
we can mention that there is a special case for which the cotangent complex defined
above is easy to understand, precisely when $X$ is a Deligne-Mumford derived stack 
(see \ref{secapp:Derivedschemesandstacks}).
Indeed, under this extra assumption, $\QCoh(X)$ can be described 
using quasi-coherent $\OO_X$-dg-modules on the small \'etale site
$(\dAff/X)_{et}$, of derived affine schemes \'etale over $X$. When restricted 
to this small \'etale site, $\LL_{\F}^{big}$ is already 
quasi-coherent, because all the morphisms in $(\dAff/X)_{et}$ are automatically \'etale. 
Thus,  $\LL_{\F}^{big}=\LL_{\F}$ as a quasi-coherent $\OO_X$-dg-module over the small
\'etale site of $X$. \\

Finally, we introduce the global version of the de Rham algebra $\DR(\F/k)$ of
a derived foliation $\F$. For a derived stack $X$ we have a stack of graded mixed
cdga's on the big \'etale site of $F$
$$\DR_{X/k} : (\dAff_k/X)^{op} \longrightarrow \egrcdga_k,$$
which sends $u : \Spec\, A \to X$ to $\DR(A/k)$. More generally, for
$X$ as above and $\F\in \Fol(X/k)$, 
we have a stack 
$$\DR_{\F/k} : (\dAff_k/X)^{op} \longrightarrow \egrcdga_k,$$
which sends $u : \Spec\, A \to X$ to $\DR(u^*(\F)/k)$. This stack
comes moreover equipped with a canonical
morphism $\DR_{X/k} \to \DR_{\F/k}$ of stacks of graded mixed cdga's, and thus
will be referred to as a \emph{graded mixed $\DR_{X/k}$-algebra}.

\begin{df}\label{dII-6}
Let $X\in\dSt_k$ and $\F\in \Fol(X/k)$. 
\begin{enumerate}
\item The \emph{stack of de Rham algebras of $\F$} (relative to $k$) is
the stack of graded mixed cdga's $\DR_{\F/k}$ defined above. 

\item The \emph{(global) de Rham algera of $\F$} (relative to $k$) is
$$\DR(\F/k):=\Gamma(\DR_{\F/k}) \in \egrcdga_k.$$
\end{enumerate}
\end{df}

With the definition above, it is straightforward to 
reinterpret derived foliations on $X$ as sheaves of graded mixed 
$\DR_{X/k}$-algebras. This reinterpretation simply associates
to a derived foliation $\F$ the sheaf $\DR_{\F/k}$ of the above definition. 
We get this way a full embedding of $\Fol(X/k)$ into the $\s$-category 
$\DR_{X/k}-\egrcdga$ of sheaves of graded mixed $\DR_{X/k}$-algebras
$$\Fol(X/k) \hookrightarrow \DR_{X/k}-\egrcdga.$$
The essential image of this full embedding is easy to describe by simply unwiding the
definitions.

\begin{cor}\label{cdII-6}
The above full embedding identifies $\Fol(X/k)$ with the full sub-$\s$-category of
sheaves of graded mixed $\DR_{X/k}$-algebras $\D$ which are quasi-coherent 
over $\DR_{X/k}$: for any diagram of derived stacks 
$\Spec\, B \to \Spec\, A \to X$ the natural morphism 
$$\DR_{X/k}(\Spec\, B) \otimes_{\DR_{X/k}(\Spec\, B)}\D(\Spec\, A) \to 
\D(\Spec\, B)$$
is an equivalence of graded mixed complexes.
\end{cor}

\begin{rmk}
\emph{We warn the reader that in the statement above $\DR_{X/k}$ is itself not 
quasi-coherent over $\OO_X$, even as merely a sheaf of graded $\OO_X$-cdga's (by forgetting
the mixed structure). In particular, $\D$ is also not quasi-coherent
as a graded $\OO_X$-algebra over $X$. It is rather of the form 
$Sym_{\OO_X}(\LL[1])$ for an $\OO_X$-dg-module $\LL$ which is not necessarly quasi-coherent.}
\end{rmk}

When the mixed structure is forgotten, 
the sheaf $\DR_{\F/k}$ becomes, as a graded $\DR_{X/k}$-cdga, equivalent to 
$Sym_{\OO_X}(\LL_{\F/k}^{big}[1])$. It is important to note that 
the graded cdga $Sym_{\OO_X}(\LL_{\F/k}[1])$ does not carry 
any natural compatible mixed structure, at least not before
passing to global sections (see theorem. \ref{tII-1}). This is already true for the
tautological derived foliation $*_X$ whose de Rham
algebra is $\DR_{X/k}$, as the de Rham differential does not exists
as a morphism $\OO_X \to \LL_{X/k}$ but only as a morphism to $\LL_{X/k}^{big}$. \\

We finish by the important descent statement showing that even thought
the de Rham differential of a derived foliation does not exists
as a morphism $\OO_X \to \LL_{\F/k}$, it does so after having taken global sections.
This result can be seen as an extension of the descent statement for closed
forms of \cite[Prop. 1.14]{ptvv} to the derived foliation setting.

\begin{thm}\label{tII-1}
Let $\F \in \Fol(X/k)$ be a derived foliation  over a derived Artin stack $X$. 
We assume that the conormal complex $\mathcal{N}^*_{\F}$ is a perfect complex on $X$. Then, 
the natural morphism
$\LL_{\F/k} \to \LL_{\F/k}^{big}$ induces an equivalence $\OO_X$-dg-modules
$$\wedge^i_{\OO_X}\LL_{\F/k} \to (\wedge^i_{\OO_X}\LL_{\F/k}^{big})^{qcoh}.$$
In particular, there exists a natural equivalence 
of graded (non-mixed) cdga's
$$\Gamma(X,Sym_{\OO_X}(\LL_{\F/k}[1])) \simeq \Gamma(X,Sym_{\OO_X}(\LL_{\F/k}^{big}[1])) 
\simeq \DR(\F/k).$$
\end{thm}

\textit{Proof.} The last statement of the theorem is a direct consequence of the fact that 
the quasi-coherator, by its definition as the right adjoint of the natural
inclusion $\s$-functor, does not change global sections. Also, for $i=0$ the statement
is trivial, as for $i=1$ for which it is true by definition of $\LL_{\F/k}$.

We first show the following lemma, for which the assumption that $\mathcal{N}_{\F/k}^*$ is
perfect is not needed.

\begin{lem}\label{ltII-1}
The commutative diagram of $\OO_X$-dg-modules
$$\xymatrix{
\LL_{X/k} \ar[r] \ar[d] & \LL_{\F/k} \ar[d] \\
\LL_{X/k}^{big} \ar[r] & \LL_{\F/k}^{big}
}$$
is homotopy cartesian.
\end{lem}

\textit{Proof of the lemma.} In the diagram under consideration, the vertical 
morphisms are the adjunction morphism coming from the fact that top horizontal 
morphism is obtained from the bottom horizontal morphism by 
applying $(-)^{qcoh}$ (see theorem \ref{tA-5-1}). Therefore, as 
$(-)^{qcoh}$ is a stable $\s$-functor, it is enough to show that the cofiber
of $\LL_{X/k}^{big} \to \LL_{\F/k}^{big}$ is already a quasi-coherent $\OO_X$-dg-module.

For $u : U=\,Spec\, A \to X$, the morphism of $A$-dg-modules $\LL_{X/k}^{big}(U) \to 
\LL_{\F/k}^{big}(U)$ is equivalent, by definition, to the natural morphism
$$\LL_{A/k} \to \LL_{u^{*}(\F)/k}.$$
For any $f : V=\Spec\, B \to U$, with $v=uf : V \to X$, the commutative diagram
$$\xymatrix{
B\otimes_A \LL_{A/k} \ar[r] \ar[d] & B\otimes_A \LL_{u^{*}(\F)/k} \ar[d] \\
\LL_{B/k} \ar[r] & \LL_{f^*u^*(\F)/k}
}$$
is cartesian by definition of pull-backs of derived foliations over affine derived schemes.
Therefore, if $C$ denotes the cofiber of $\LL_{X/k}^{big} \to \LL_{\F/k}^{big}$, then the
natural morphism of $B$-dg-modules
$B\otimes_A C(U) \to C(V)$
is an equivalence. This is precisely the statement that $C$ is a quasi-coherent $\OO_X$-dg-module.
\hfill $\Box$ \\

According to the lemma, we have a commutative diagram with exact rows
$$\xymatrix{
\LL_{X/k} \ar[r] \ar[d] & \LL_{\F/k} \ar[d] \ar[r] & C \ar[d]^-{=}\\
\LL_{X/k}^{big} \ar[r] & \LL_{\F/k}^{big} \ar[r] & C.
}$$

For all $i>1$ fixed, 
these exact rows induces finite filtrations on $\wedge^i_{\OO_X}\LL_{\F/k}$ 
and $\wedge^i_{\OO_X}\LL_{\F/k}^{big}$ whose graded pieces of weight $p$ are respectively 
$$\wedge^p_{\OO_X}\LL_{X/k} \otimes_{\OO_X}\wedge^{i-p}_{\OO_X}C \qquad \textrm{and} \qquad
\wedge^p_{\OO_X}\LL^{big}_{X/k} \otimes_{\OO_X}\wedge^{i-p}_{\OO_X}C.$$
These filtrations are moreover compatible with the vertical morphism in the middle.
Therefore, in order to show that the morphism $\wedge^i_{\OO_X}\LL_{\F/k} \longrightarrow
(\wedge^i_{\OO_X}\LL_{\F/k}^{big})^{qcoh}$ is an equivalence, it is enough to show that 
it induces an equivalence of those graded pieces. \\

\begin{lem}\label{ltII-1-2}
For any perfect complex $E$ on $X$ and any 
$E'\in \OO_X-\dg$, the natural morphism
$$E\otimes (E')^{qcoh} \longrightarrow (E\otimes E')^{qcoh}$$
is an equivalence.
\end{lem}

\textit{Proof of the lemma.} We use here that perfect complexes are dualizable objects
in $\QCoh(X)$, and thus also in  the bigger 
symmetric monoidal $\s$-category $\OO_X-\dg$. 
Therefore, for any quasi-coherent complex $E_0$, we have functorial 
equivalences
$$\Map(E_0,E\otimes (E')^{qcoh}) \simeq Map(E_0 \otimes E^\vee,(E')^{qcoh}) \simeq
\Map(E_0\otimes E^\vee,E') \simeq \Map(E_0,E\otimes E').$$
By observation, the composed equivalence $\Map(E_0,E\otimes (E')^{qcoh}) \to \Map(E_0,E\otimes E')$
is induced by the canonical morphism $E\otimes (E')^{qcoh} \longrightarrow (E\otimes E')^{qcoh}$, and
the lemma thus follows from Yoneda.
\hfill $\Box$ \\

To finish the proof of the theorem, we invoke theorem \ref{tA-5-1} and the lemma
\ref{ltII-1-2}. They imply that the natural morphisms
$$\wedge^p_{\OO_X}\LL_{X/k} \otimes_{\OO_X}\wedge^{i-p}_{\OO_X}C \to 
\wedge^p_{\OO_X}\LL^{big}_{X/k} \otimes_{\OO_X}\wedge^{i-p}_{\OO_X}C$$
are equivalences (note that $C$ is perfect by assumption), 
and thus, by what we have already see before, 
that the morphism $\wedge^i_{\OO_X}\LL_{\F/k} \to (\wedge^i_{\OO_X}\LL_{\F/k}^{big})^{qcoh}$ is also 
an equivalence.
\hfill $\Box$ \\

\subsection{The $\s$-categories of derived foliations}\label{subsec:Thescategoriesofderivedfoliations}

By construction, for a derived stack $X \in \dSt_k$, derived foliations on $X$ (relative to 
$k$) form an $\s$-category $\Fol(X/k)$. We study here two basic properties of these
$\s$-categories: existence of finite limits and pull-backs functoriality along arbitrary morphisms. \\

Let $f : X \to Y$ be a morphism of derived stacks. It induces 
an $\s$-functor on the $\s$-categories of maps to the derived stack $\Fol(-/k)$
$$f^* : \Fol(Y/k)=Map(Y,\Fol(-,k)) \longrightarrow  \Fol(X/k)=Map(X,\Fol(-,k)).$$
This defines an $\s$-functor
$$\Fol(-/k) : \dSt_k^{op} \longrightarrow \scat,$$
which is the left Kan extension of $\Fol(-/k) : \dAff_k^{op} \longrightarrow \scat$
along the Yoneda embedding $\dAff_k \hookrightarrow \dSt_k$ (see \cite[\S 2.1]{toenems}).

Note that the big cotangent complex construction of Definition \ref{dII-5} defines an 
$\s$-functor
$$\LL^{big} : \Fol(X/k) \longrightarrow (\OO_X-\dg)^{op}.$$
Composed with the quasi-coherator, it defines another $\s$-functor
$$\LL : \Fol(X/k) \longrightarrow (\OO_X-\dg)^{op} \longrightarrow (\QCoh(X))^{op}.$$
the $\s$-functors $\LL^{big}$ and $\LL$ are not functorial in $X$ and rather 
defines a \emph{lax natural transformations} between derived stacks of $\s$-categories 
$\Fol(-/k) \longrightarrow (\OO_{-}-\dg)^{op}\longrightarrow (\QCoh(-))^{op}.$
In other words, for a morphism of derived stacks $f : X \to Y$,
and $\F \in \Fol(Y/k)$,
there exists a natural morphism of 
$\OO_X$-dg-modules
$u : f^*(\LL_{\F}^{big}) \longrightarrow  \LL_{f^*(\F)}^{big}.$
This morphism $u$ depends functorialy in $f$ and $\F$ in some obvious manner, with 
all the required homotopy coherences. The formal way to express these coherences is by 
stating the existence of an $\s$-functor on the Grothendieck's integral
which commutes with the natural projections to $\dSt_k$ (see \cite[\S 1.3]{MR3338682})
$$\xymatrix{\int_{\dSt_k} \Fol(-/k) \ar[rd] \ar[rr]^-{\LL^{big}} & & 
\int_{\dSt_k} (\OO_{-}-\dg)^{op} \ar[ld] \\
 & \dSt_k. & }$$

The next Proposition shows that $\Fol(X/k)$ admits arbitrary limits for 
any derived stack $X \in \dSt_k$.

\begin{prop}\label{pII-3}
\begin{enumerate}
\item For any derived stack $X\in \dSt_k$, the $\s$-category $\Fol(X/k)$ admits small
limits. 

\item For any morphism of derived stacks $f : X \to Y$, the pull-back 
$\s$-functor
$f^* : \Fol(Y/k) \to \Fol(X/k)$ preserves limits.
\end{enumerate}

\end{prop}

\textit{Proof.} $(1)$ The $\s$-category $\Fol(X/k)$ can 
be written as a limit
$$\Fol(X/k) \simeq \llim_{\Spec\, A \to X}\Fol(A/k).$$
In order to show that $\Fol(X/k)$ admits small limits it is therefore 
enough to show the following two facts.

\begin{enumerate}
\item For all $\Spec\, A \in \dAff_k$, the $\s$-category
$\Fol(A/k)$ admits small limits.

\item For any morphism $f : \Spec\, B \to \Spec\, A$ in $\dAff_k$ the induced
$\s$-functor
$f^* : \Fol(B/k) \to \Fol(A/k)$ preserves small limits. 
\end{enumerate}

For the first of these assertions, let 
$I \to \Fol(A/k)$ be a small diagram. It corresponds to a diagram 
$D : I^{op} : \DR(A/k)-\egrcdga$, of graded mixed $\DR(A/k)$-cdga's.
We consider $C:=\colim_{i\in I^{op}}D(i) \in \DR(A/k)-\egrcdga$. We claim that 
$C$ is a graded mixed $\DR(A/k)$-cdga satisfying the conditions of Definition
\ref{dII-1}. Indeed, $D$ induces a diagram on cotangent complexes
$i \mapsto \LL_{D(i)/k}$, from $I^{op}$ to $A-\dg$, simply by considering 
the $\s$-functor obtained by the weight one piece. This is a diagram 
of $A$-dg-modules under $\LL_{A/k}$. We know that the
forgetful $\s$-functor, from graded mixed cdga's to graded cdga's commutes, 
and even reflects, limits and colimits. Therefore,
$C$ is, as a graded $\DR(A/k)$-cdga, of the form 
$$C\simeq \colim_{i\in I^{op}}Sym_A(\LL_{D(i)/k}[1]) \simeq
Sym_A(\LL[1]),$$
where $\LL=\colim_{i\in I^{op}}\LL_{D(i)/k}$ is the colimit in the $\s$-category
$\LL_{A/k}/A-\dg$ of $A$-dg-modules under $\LL_{A/K}$.
Finally, by definition $\Fol(A/k)$ is a full sub-$\s$-category of the complete $\s$-category
$(\DR(A/k)-\egrcdga)^{op}$, and we just have seen that this sub-$\s$-category
is stable by arbitrary limits. In particular, $\Fol(A/k)$ possesses limits.

The second of the above assertions is also true, as 
$f^*$ is defined, on the level of graded mixed cdga's, by 
the base change $-_{\DR(A/k)}\DR(B/k)$, and thus sends
colimits of graded mixed $\DR(A/k)$-cdga's to colimits of graded
mixed $\DR(B/k)$-cdga's. \\

$(2)$ As before, the statement is reduced to morphisms between affine derived schemes, for
which we already have seen that pull-backs preserve limits.
\hfill $\Box$ \\

As a corollary of the proof of Proposition \ref{pII-3} we have the following formula
for the cotangent complex of a limit of derived foliations. In the statement below
we denote by $\LL_{X/k}^{big}$ the \emph{big} or \emph{fake} cotangent 
complex of a derived stack $X$, defined as an $\OO_X$-dg-module on $\dAff_k/\dSt_k$ by 
sending $\Spec\, A \to X$ to $\LL_{A/k}$. Note that its image by the quasi-coherator
$\OO_X-\dg \to \QCoh(X)$ is the usual cotangent complex $\LL_{X/k}$ 
(see theorem \ref{tA-5-1}).

\begin{cor}\label{cpII-3}
Let $X$ be a derived stack, and $\F$ a derived foliation that 
is expressed as a finite limit $\F\simeq \llim_{i\in I}\F_i$ in $\Fol(X/k)$. Then, the canonical
morphisms
$$\colim_{i\in I^{op}} (\LL^{big}_{\F_i}) \longrightarrow \LL_{\F}^{big} \qquad
\colim_{i\in I^{op}} (\LL_{\F_i}) \longrightarrow \LL_{\F}$$
are equivalences in the $\s$-categories $\LL_{X/k}^{big}/\OO_X-\dg$ and
$\LL_{X/k}/\QCoh(X)$ of objects under $\LL_{X/k}^{big}$ and $\LL_{X/k}$.
\end{cor}

It is also important to note that cotangent complexes of pull-backs of derived foliations
can be computed explicitly.

\begin{cor}\label{cltII-1}
Let $f : X \to Y$ be a morphism of derived Artin stacks, and $\F \in\Fol(Y/k)$. Then, 
there exists a cocartesian square of quasi-coherent complexes on $X$
$$\xymatrix{
f^*(\LL_{Y/k}) \ar[r] \ar[d]  & f^*(\LL_{\F/k}) \ar[d] \\
\LL_{X/k} \ar[r] & \LL_{f^*(\F)/k}.
}$$
\end{cor}

\textit{Proof.} We have tautological equivalences of $\OO_X$-dg-modules
$$f^*(\LL_{Y/k}^{big}) \simeq \LL_{X/k}^{big} \qquad f^*(\LL_{\F/k}^{big}) \simeq 
\LL_{f^*(\F)/k}^{big}.$$
By Lemma \ref{ltII-1} applied to $\F$ we have a cocartesian square of $\OO_Y$-dg-modules
$$\xymatrix{
\LL_{Y/k} \ar[r] \ar[d]  & \LL_{\F/k} \ar[d] \\
\LL_{Y/k}^{big} \ar[r] & \LL_{\F/k}^{big}.}$$
Applying $f^*$ we get another cocartesian square
$$\xymatrix{
f^*(\LL_{Y/k}) \ar[r] \ar[d]  & f^*(\LL_{\F/k}) \ar[d] \\
f^*(\LL_{Y/k}^{big}) \ar[r] & f^*(\LL_{\F/k}^{big}).}$$
The corollary follows by the tautological identifications
$$f^*(\LL_{Y/k}^{big})=\LL_{X/k}^{big} \qquad f^*(\LL_{\F/k}^{big})=\LL_{f^*(\F)/k}^{big},$$
and an application of the $\s$-functor $(-)^{qcoh}$.
\hfill $\Box$ \\

We finish this part with the specific case of derived foliations over Deligne-Mumford
derived stacks. We remind that $X \in \dSt_k$ is called \emph{Deligne-Mumford} if 
it is $n$-geometric for some $n$ in the sense of \cite[Def. 1.3.3.1]{hagII}, and if 
there are \'etale morphisms $\Spec\, A_i \to X$ such that
$\coprod_{i}\Spec\, A_i \to X$ is an epimorphism. Derived schemes are
in particular Deligne-Mumford, as each of the morphism $\Spec\, A_i \to X$ can be
chosen to be a Zariksi open immersion.
For such a Deligne-Mumford derived stack $X$, we can define the small \'etale site
$X_{et}$, which is the full sub-$\s$-category of $\dAff_k/X$ consisting of 
all \'etale morphisms $\Spec\, A \to X$, endowed with the induced \'etale topology.
Stacks on the $\s$-site $X_{et}$ form a sub-$\s$-topos of $(\dAff_k/X)^{\sim,et}$, 
called the small \'etale $\s$-topos of $X$.
One important feature of the small \'etale site is the inclusion of the truncation
$j : t_0X \to X$ induces an equivalence of $\s$-topos
$j^* : X_{et}^{\sim} \simeq (t_0X)_{et}^{\sim}$. In particular the small \'etale topos
does not depend on derived structures, in a similar manner than the small \'etale topos
of a scheme does not depend on non-reduced scheme structures. 

It is immediate to check that the big cotangent complex $\LL_{X/k}^{big}$
of a derived Deligne-Mumford stack $X$ is quasi-coherent when restricted to the small
\'etale site. More precisely, the adjunction morphism of $\OO_X$-dg-modules
$\LL_{X/k} \to \LL_{X/k}^{big}$ restricts on $X_{et}$ to an equivalence. Similarly, 
if $\F \in \Fol(X/k)$ is a derived foliation on $X$, then the canonical 
morphism $\LL_{\F/k} \to \LL_{\F/k}^{big}$ restricts to an equivalence on $X_{et}$. Implementing
on this, we can prove the following proposition, which describes derived foliations in 
terms of sheaf theory on the small \'etale sites. For this we introduce the following notation:
if $\F\in \Fol(X/k)$ we denote by $\DR_{\F/k,et}$ the restriction of the sheaf of graded
mixed cdga's $\DR_{\F/k}$ of Definition \ref{dII-6} to $X_{et}$.

\begin{prop}\label{pII-4}
Let $X\in \dSt_k$ be a derived Deligne-Mumford stack. The $\s$-functor
$$\Fol(X/k) \longrightarrow (\DR_{X/k,et}-\egrcdga)^{op}$$
from derived foliations to stacks of graded mixed $\DR_{X/k,et}$-cdga's on the small
\'etale site of $X$, is fully faithful. Its essential image consists of stacks of
graded mixed $\DR_{X/k,et}$-cdga $D$ satisfying the following two properties.
\begin{enumerate}
\item The weight $1$ piece $D(1)$ of $D$ is a quasi-coherent 
$\OO_X$-dg-module, for the natural $\OO_X$-module structure induced by the graded 
$\DR_{X/k,et}$-algebra
structure on $D$.

\item The natural morphism of stacks of graded $\DR_{X/k,et}$-cdga's
$Sym_{\OO_X}(D(1)) \to D$ is an equivalence.

\end{enumerate}
\end{prop}

\textit{Proof.} This is a formal exercise left to the reader.
\hfill $\Box$ \\

The importance of Proposition \ref{pII-4} is that 
$\DR(\F/k)_{et}$ is, as a stack of graded cdga's, of the form 
$Sym_{\OO_X}(\LL_{\F/k}[1])$. In other words, the de Rham differential
for $\F$ exists as a morphism of stacks of $k$-dg-modules $\OO_X \to \LL_{\F/k}$ globally 
defined on $X_{et}$. This is specific property of Deligne-Mumford stacks, 
for general Artin derived stacks this de Rham differential only exists
as a $k$-linear morphism  $\OO_X \to \LL_{\F/k}^{big}$ and does not factor a priori through the
adjunction morphism $\LL_{\F/k} \to \LL_{\F/k}^{big}$. \\

Finally, the various notions introduced in the definition for derived foliations
of derived affine schemes, extends to derived foliations on 
derived Artin stacks in a very natural manner.

\begin{df}\label{dII-8}
Let $\F$ be a derived foliation on a derived stack $X$.
\begin{enumerate}
\item The derived foliation $\F$ is \emph{perfect} (resp. \emph{almost perfect}) if 
$\LL_{\F/k}$ is a perfect (resp. \emph{almost perfect}) complex on $X$.

\item The derived foliation $\F$ is \emph{n-connective}, for some
integer $n\in \ZZ$, if the cohomology sheaves 
$\underline{H}^i(\LL_{\F/k})=0$ for all $i >n$
(i.e. if $\LL_{\F/k}[n]$ is connective).

\item The derived foliation $\F$ is \emph{smooth} if it is perfect and if
$\LL_{\F/k}$ is of Tor amplitude contained in $[0,\s[$.

\item The derived foliation $\F$ is \emph{quasi-smooth} if it is perfect and if
$\LL_{\F/k}$ is of Tor amplitude contained in $[-1,\s[$.

\end{enumerate}

\end{df}

We finish this part by a general method to compute mapping spaces in the 
$\s$-category $\Fol(X/k)$ of derived foliations on a derived stack $X$, in terms
of their cotangent complexes and actions of the group stack $\cH_0$ described in
\S \ref{sec:Thegeometryofgradedmixedobjects}. 
Let then $X$ be any derived stack (not necessarily assumed to be Artin). We define its 
global cotangent complex using the quas-coherator construstion of \ref{secapp:CotangentcomplexesofderivedArtinstacks}
$$\LL_{X/k}:=(\LL_{X/k}^{big})^{qcoh} \QCoh(X).$$
For a derived foliation $\F \in \Fol(X)$, we also have its big cotangent 
complex $\LL_{\F/k}^{big}$ and its cotangent complex defined the same way
$$\LL_{\F/k}:=(\LL_{\F/k}^{big})^{qcoh} \in \QCoh(X).$$
We consider the exact triangle of $\OO_X$-dg-modules
$$\xymatrix{\cN_{\F}^* \ar[r] & \LL_{X/k}^{big} \ar[r] & \LL_{\F/k}^{big},}$$
where $\cN_{\F}^*$ is by definition the conormal complex of $\F$ (relative to $k$). A simple
observation show that $\cN_{\F}^*$ is quasi-coherent. As a result, the commutative square
$$\xymatrix{
\LL_{X/k}^{big} \ar[r] & \LL_{\F/k}^{big} \\
\LL_{X/k}  \ar[u] \ar[r] & \LL_{\F/k} \ar[u]
}$$
is homotopy cartesian.

For two quasi-coherent complexes $E$ and $E'$ on $X$, we can form
the derived stack of morphism $\uMap(E,E') \in \dSt_k$. By definition, it sends 
a connective cdga $A$ to the simplicial set $\Map_{\QCoh(X)}(E,A \otimes E')$, 
where $A \otimes_k E'$ is the external tensor product of $E'$ by the $k$-linear complex underlying $A$.
It is also the $\s$-functor sending $A$ to $\Map_{\QCoh(X\times \Spec\, A)}(p^*(E),p^*(E'))$, 
where $p : X\times \Spec\, A \to X$ is the natural projection. The stack 
$\uMap(E,E')$ possesses a natural linear structure (i.e. is a stack of $\OO$-dg-modules), but 
we will consider it merely as a space valued stacks and will not keep track of its 
linear structure in what follows. When $E$ and $E'$ both sits under a given object
$a : E_0 \to E$, and $b : E_0 \to E'$, we also have the mapping stack under $E_0$
$\uMap_{E_0/}(E,E')$, of morphisms from $E$ to $E'$ that commute with 
the morphism $a$ and $b$. By definition, we have a fibration sequence of 
derived stacks
$$\uMap_{E_0/}(E,E') \to \uMap(E,E') \to \uMap(E_0,E')$$
where the fiber is taken at the point corresponding to $b$.

Finally, we remind (\S \ref{sec:Thegeometryofgradedmixedobjects}) that $\cH_0$ is the group stack $\Gm \ltimes B\Ga$, semi-direct product 
of $B\Ga$ with $\Gm$ (for the natural induced weight $1$ action of $\Gm$ on $\Ga$).

\begin{prop}\label{pII-5}
Let $X$ be a derived stack and $\F,\F' \in \Fol(X/k)$ two derived foliations on $X$. Then, there
exists a stack $\pi : \cM(\F,\F') \to B\cH_0$ with the following properties.

\begin{enumerate}
    \item The base change $\cM(\F,\F') \times_{B\cH_0}B\Gm$,  is equivalent to 
    $$\uMap^{gr}_{\LL_{X/k}/}(\LL_{\F'/k},\oplus_{i\geq 1}(\wedge^i_{\OO_X}\LL_{\F}^{big}[i-1])^{qcoh})$$
    the derived stack of graded morphisms from $\LL_{\F'/k}$ to $\oplus_{i\geq 0}\wedge^i_{\OO_X}\LL_{\F}^{big}[i-1]$
    under $\LL_{X/k}$, whith the natural gradings where $\LL_X$ and $\LL_{\F'}$ have weight one, 
    and $\wedge^i_{\OO_X}\LL_{\F}^{big}$ has weight $i$.
    \item There is a natural equivalence of spaces
    $$\HH(B\cH_0,\cM(\F,\F')):=p_*(\cM(\F,\F'))(k) \simeq \Map_{\Fol(X)}(\F,\F')$$
    where $p : B\cH_0 \to *$ is the structural morphism, and $p_* : \dSt_k/B\cH_0 \to \dSt_k$ is the
    direct image $\s$-functor.
\end{enumerate}
\end{prop}

\textit{Proof.} We use the equivalence 
$$\mathbf{CAlg}(\QCoh(B\cH_0)) \simeq \egrcdga_k$$
between commutative algebras in the symmetric monoidal $\s$-category $\QCoh(B\cH_0)$
and graded mixed cdga's (see \S \ref{sec:Thegeometryofgradedmixedobjects}). Using this equivalence, for two graded mixed cdga $D$ and $D'$ 
the mapping space $Map_{\egrcdga}(D,D')$ can be written as in the Proposition
$$Map_{\egrcdga_k}(D,D') \simeq p_*(\cM(D,D'))(k).$$
Here the derived stack $\cM(D,D')$ lives over $B\cH_0$ and can be defined functorially as
follows. For $u : \Spec\, A \to B\cH_0$, the space 
$\cM(D,D')(A)$ is given by $\Map_{\cdga_A}(u^*(D),u^*(D'))$, where $D$ and $D'$
are here considered as quasi-coherent cdga's over $B\cH_0$ via the equivalence mentioned above. 
Note that a similar formula remains valid when $D$ and $D'$ are $D_0$-graded mixed cdga, for 
a given graded mixed cdga $D_0$, with $\cM(D,D')$ being replace by 
$\cM_{D_0}(D,D')$ the derived stack of morphisms of $D_0$-cdga's.
Assuming further that $D$ (resp. $D'$ and $D_0$) are, as graded cdga's, of the form $Sym_A(E)$ 
(resp. $Sym_A(E')$ and $Sym_A(E_0)$) for some
$A$-dg-module $E$, we can instead of $\cM_{D_0}(D,D')$ consider its reduced version $$\cM^*_{D_0}(D,D'),$$
that is  the fiber of $\cM_{D_0}(D,D') \longrightarrow \cM_{D_0}(D,A)$, taken at the natural
augmentation $D \to A$.

Now, the mapping space $\Map_{\Fol(X/k)}(\F,\F')$ is given by the limit
$$\Map_{\Fol(X/k)}(\F,\F') \simeq \lim_{\Spec\, A \to X}\Map_{\DR(A/k)-\egrcdga}(\DR_{\F'}(A),\DR_{\F}(A)).$$
Using the previous identification for mapping spaces of graded mixed cdga's, we see that the
Proposition $(2)$ is true for a certain derived stack $\cM(\F,\F') \to B\cH_0$, namely
$$\cM(\F,\F') \simeq \underset{\Spec\, A \to X}{\mathrm{lim}}\cM^*_{\DR(A/k)-\egrcdga}(\DR_{\F'}(A),\DR_{\F}(A)).$$
The base change along $B\Gm \to B\cH_0$ is then given by the derived stack 
$$\lim_{\Spec\, A}\uMap^*_{\DR(A/k)-\grcdga}(\DR_{\F'}(A),\DR_{\F}(A)),$$
limits of the derived stacks $\uMap^*_{\DR(A/k)-\grcdga}(\DR_{\F'}(A),\DR_{\F}(A))$
of morphisms of $\DR(A/k)$-linear morphisms of graded cdga's from $\DR_{\F'}(A)$ to 
$\DR_{\F}(A)$ (obtained by forgetting the mixed structures) compatible with the augmentations to $A$.
Because these cdga's are graded free, we see that 
this base change can also be written in terms of graded morphisms of complexes
$$\lim_{\Spec\, A}\uMap^{gr}_{\LL_{A/k}/\grcdga}(\LL_{\F'}^{big}(A),\oplus_{i\geq 
1}\wedge^i_{A}\LL_{\F}^{big}(A)[i-1])) \simeq 
\uMap^{gr}_{\LL_X^{big}}(\LL_{\F'}^{big},\oplus_{i\geq 1}\wedge^i_{\OO_X}\LL_{\F}^{big}[i-1]).$$
(note that the sum does not involve the term $i=0$ because we use the reduced version
of the derived stacks of morphisms of graded cdga's).
To finish the proof of the Proposition it thus simply remains to prove that the quasi-coherator
construction induces for each $\OO_X$-dg-module $E$ with a morphism $\LL_{X}^{big} \to E$, an equivalence
$$\uMap_{\LL_X^{big}}(\LL_{\F'}^{big},E) \longrightarrow \uMap_{\LL_X}(\LL_{\F'},E^{qcoh}),$$
as the result will be obtained by the special cases where $E=\wedge^i\LL_\F^{big}$. 

To see this last step, by replacing $X$ with $X \times \Spec\, A$ we see that
it is enough to check this simply on mapping spaces rather than mapping derived stacks. To check this
last step, we consider the commutative diagram
$$\xymatrix{
\Map_{\LL_X^{big}}(\LL_{\F'}^{big},E) \ar[r] \ar[d] &  \Map(\LL_{\F'}^{big},E) \ar[r] \ar[d] &  \Map(\LL_{X}^{big},E) \ar[d] \\
\Map_{\LL_X}(\LL_{\F'},E^{qcoh}) \ar[r] & \Map(\LL_{\F'},E^{qcoh}) \ar[r] & \Map(\LL_X,E^{qcoh}).
}$$
Each row of this diagram is a fibration sequence, for the base point in $\Map(\LL_X^{big},E)$
corresponding to the fixed morphism $\LL_{X}^{big} \to E$. However, we know that the square
$$\xymatrix{
\LL_X^{big} \ar[r]  & \LL_{\F'}^{big}  \\
\LL_{X} \ar[r] \ar[u] & \LL_{\F'} \ar[u]
}$$
is cocartesian, and thus the square 
$$\xymatrix{
 \Map(\LL_{\F'}^{big},E) \ar[r] \ar[d] &  \Map(\LL_{X}^{big},E) \ar[d] \\
 \Map(\LL_{\F'},E^{qcoh}) \ar[r] & \Map(\LL_X,E^{qcoh}).
}$$
is cartesian. This implies that the natural morphism 
$\Map_{\LL_X^{big}}(\LL_{\F'}^{big},E) \to 
\Map_{\LL_X}(\LL_{\F'},E^{qcoh})$ is indeed an equivalence as wanted.
\hfill $\Box$ \\

An important corollary of the previous Proposition is obtained when $X$ is moreover a derived \emph{Artin} stack, as 
we know from \ref{tA-5-1} that $\wedge^i_{\OO_X}\LL_{\F}^{big}[i-1])^{qcoh} \simeq \wedge^i_{\OO_X}\LL_{\F}[i-1]$. 
In particular we get the following important consequence.

\begin{cor}\label{cpII-5}
Let $X$ be a derived Artin stack and let $\F$ be a derived foliation such that 
$\mathbb{L}_\F\simeq 0$. Then $\F$ is an initial object in $\Fol(X)$.
\end{cor}

\textit{Proof of Corollary.} In this case, for any $\F' \in \Fol(X)$, 
the natural projection $\cM(\F,\F') \to \times_{B\cH_0}B\Gm$
must be an equivalence and thus $p_*(\cM(\F,\F'))=*$ is the punctual derived stack. It follows
from the proposition \ref{pII-5} that $\Map_{\Fol(X)}(\F,\F') \simeq *$.
\hfill $\Box$ \\

\subsection{More examples}\label{subsec:Moreexamples}

We have given several examples of derived foliations on derived affine schemes in 
Section \ref{subsec:basicexamplesofderivedfoliations}. We now provide several global examples on 
derived schemes and stacks. \\

\begin{ex}\label{ex:tautologicalfoliationsglobal}
\emph{\textbf{Tautological derived foliations.} For any derived stack $X$, the
$\s$-category possesses a final and and sometimes an initial object. The final object
is denoted by $*_X$, whose cotangent complex is $\LL_{X/k}$. As an $\s$-functor
$(\dAff_k/X)^{op} \to \egrcdga$ it is given by sending $\Spec\, A \to X$
to $\DR(A/k)$. }

\emph{The initial object of $\Fol(X/k)$, when it exists, is denoted by $0_X$. 
When $X$ is representable by a derived Artin stack, such an initial object
always exists and moreover its cotangent complex is $0$.
As an $\s$-functor
$(\dAff_k/X)^{op} \to \egrcdga$, $0_X$ is given by sending $\Spec\, A \to X$
to $\DR(A/X)$, the relative de Rham algebra of $\Spec\, A$ over $X$. We need here
to provide more definitional details, as $X$ is not assumed to be affine
and thus $\DR(A/X)$ has not been defined yet. We define this notion 
by using descent over $X$, by the standard formula
$$\DR(A/X) = \llim_{\Spec\, B \to X}\DR(\Spec\, A \times_{\Spec, B}X/B).$$
In more general terms, for a morphism of derived stacks $f : Y \to Z$, we set 
$$\DR(Y/Z):=\llim \, \DR(A/B) \in \egrcdga_k,$$
where the limit is taken over the $\s$-category of commutative diagrams of the form
$$\xymatrix{
\Spec\, A \ar[r] \ar[d] & Y \ar[d] \\
\Spec\, B \ar[r] & Z.}$$
As the diagonal of $X$ is representable by a derived Artin stack, 
all the morphisms $\Spec\, A \to X$ are also representable, 
and thus the relative cotangent complex $\LL_{A/X}$ always exits as an $A$-dg-module.
Moreover, by the descent statement theorem 
\ref{tII-1} we have an equivalence of graded cdga's
$$\DR(A/X) \simeq Sym_A(\LL_{A/X}[1]).$$
This shows that $0_X$ is indeed a stack of graded mixed cdga's on $\dAff/X$ which 
corresponds to an object in $\Fol(X/k)$. Moreover, the big cotangent complex
of $0_X$ is the $\OO_X$-dg-module sending $\Spec\, A \to X$ 
to $\LL_{A/X}$, and we thus have an exact triangle of $\OO_X$-dg-modules
$$
\LL_{X/k} \to \LL_{X/k}^{big} \to \LL_{0_X}^{big}.$$
By our theorem \ref{tA-5-1} we have $\LL_{0_X}=(\LL_{0_X}^{big})^{qcoh}\simeq 0$ and thus
by Corollary \ref{cpII-5} $0_X$ is indeed an initial object.} \end{ex}

\begin{ex}\label{ex:integrablefoliations}
\emph{\textbf{Globally integrable derived foliations.}  
As in the affine case (see \ref{ex:derhamfoliations}, we will say that a derived foliation 
$\F \in \Fol(X/k)$ is} globally integrable (or simply integrable) \emph{if
there is a morphism $f : X \to Y$ with $Y$ a derived Artin stack such that $\F \simeq 
f^*(0_Y)$. We note that $Y$ is not uniquely defined, as any factorization
$f: \xymatrix{X \ar[r]^-{g} & Y' \ar[r]^-p & Y}$ with $p$ being étale has the property
that $g^*(0_Y')\simeq f^*(0_Y)$ (simply because $p^*(0_Y)=0_Y'$ by étaleness and the
corollary \ref{cpII-5}).} \end{ex}

\begin{ex}\label{ex:relativefoliations}
\emph{\textbf{Relative global derived foliations.}  Let $f : X \to Y$ be a morphism in $
\dSt_k$ with $Y$ being a derived Artin stack. 
 We define the $\s$-category of \emph{derived foliations on $X$ relative to $Y$}
to be the comma $\s$-category 
$$\Fol(X/Y):=\Fol(X/k)/f^*(0_Y),$$
of derived foliations over $f^*(0_Y)$ (we know $0_Y$ exists by our previous example).
Informally, a derived foliations on $X$ relative to $Y$ is a family of derived foliations in
the fibers of $f$. When $X=\Spec\, A \to Y=\Spec\, B$ is 
a a morphism of derived affine schemes, $\Fol(X/Y)$ is canonically 
equivalent to $\Fol(X/B)$. This follows simply from the observation that 
the $\s$-category of graded mixed $\DR(A/B)$-cdga's can be also expressed as the
comma $\s$-category of graded mixed cdga's under $\DR(A/B)$. By gluing we also deduce that 
$\Fol(X/Y)$ is equivalent to $\Fol(X/B)$ in the sense of Definition \ref{dII-4} for
any derived stack $X$ with a map to $Y=\Spec\, B$.} \end{ex}

\begin{ex}\label{ex:closed1forms}
\emph{\textbf{Closed 1-forms.} We recall from \cite{ptvv} the notion of closed $p$-forms on 
a derived stack $X$ (relative to $k$). This notion can be expressed easily in terms
of graded mixed complexes. For this, we consider 
the graded mixed complex $k(p)[p-n]$, which is $k$ sitting in pure weight $p$ and
cohomological deree $n-p$ (and both the cohomological differential as well as
the mixed structure being zero). The space of closed $p$-forms
on $X$ of degree $n$ is defined to be the mapping space of graded mixed complexes
$$\mathcal{A}^{p,cl}(X/k):=\Map_{\egrdg_k}(k(p)[p-n],\DR(X/k)).$$
We can use here the equivalence between graded mixed complexes and complete
filtered complexes in order to express the space of closed $p$-forms 
in a more traditional manner in terms of truncated de Rham complex 
$F^{-p}\CDR(X/k):=F^{-p}|\DR(X/k)|^t$ (see Definition \ref{dI-4}).}

\emph{Let now $\omega : k(1)[1-n] \to \DR(X/k)$ be a closed $1$-form on $X$ 
of degree $n$. As $\DR(X/k)$ is the graded mixed cdga of global sections
of $\DR_{X/k}$ (see theorem \ref{tII-1}), the form $\omega$ can also be considered
as a morphism of stacks of graded mixed cdgas on $\dAff/X$
$$Sym_{k}(k(1)[1-n]) \longrightarrow \DR_{X/k}.$$ 
We set 
$$\DR_\omega:=\DR_{X/k}\otimes_{Sym_{k}(k(1)[1-n])}k,$$
where we use the canonical augmentation morphism $Sym_{k}(k(1)[1-n])k \to k$.
Note that, as a stack of graded cdga's, $\DR_\omega$ is of the form
$Sym_{\OO_X}(\LL_\omega[1]^{big})$ where $\LL_\omega[1]^{big}$ is the
$\OO_X$-dg-module whose values on $\Spec\, A \to X$ is 
the cone of the morphism $A[-n] \to \LL_{A/k}$ given by the restriction of $\omega$ 
as a $1$-form of cohomological degree $n$ on $\Spec\, A$. It is clear from this description
that $\DR_{\omega}$ defines an object in $\Fol(X/k)$, denoted by $\F_\omega$, 
and called the \emph{derived foliation associated to $\omega$}. Note that the
cotangent complex of $\F_\omega$ is by construction $(\LL_\omega[1]^{big})^{qcoh}$
and thus sits in an exact triangle of quasi-coherent complexes
$$\OO_X[-n] \to \LL_{X/k} \to \LL_{\F_\omega/k}$$
where the first morphism is defined by $\omega$.
In particular, when $X$ is a smooth scheme over $k$ and $n=0$, a closed $1$-form (of degree
$0$) on $X$ simply is a section $\omega$ of $\Omega^1_{X/k}$ such that $dR(\omega)=0$. 
The derived foliation associated to $\omega$ has a cotangent complex given 
explicitly by a lengh two complex $\OO_X \to \Omega^1_{X/k}$, where 
$\Omega^1_{X/k}$ sits in degree $0$. We see in particular
that the derived foliation $\F_\omega$ is smooth in the sense of Definition
\ref{dII-8} if and only if $\omega$ is nowhere
vanishing, or in other words if $\omega$ identifies $\OO_X$ with a subbundle 
of $\Omega^1_{X/k}$. When $X$ is moreover proper over $k$, Hodge theory implies that 
any global section $\omega \in \Gamma(X,\Omega_{X/k}^1)$ is closed for the de Rham
differential, and thus defines a derived foliation $\F_\omega$. This provides many examples
of derived foliations on smooth proper schemes over $k$. } \end{ex}

\begin{ex}\label{ex:vectorfields}
\emph{\textbf{Vector fields.} Let $X=\Spec\, A$ be a derived affine scheme
over $k$ and $\nu : \LL_{X/k} \to A$ be a morphism of $A$-dg-module. Such a morphism
will be called a vector field on $X$ (relative to $k$). Note that this is the 
same thing as a $k$-linear derivation $A \to A$. We define a graded mixed
cdga $\DR_\nu$ by considering $Sym_A(A[1]) = A \oplus A[1]$ as a graded cdga, and 
define the mixed structure by the derivation $\nu : A \to A$. The weight zero part
is equal here to $A$, and thus we have a natural induced morphism of graded mixed
cdga's $\DR(A/k) \to \DR_\nu$. By construction $\DR_\nu$ defines
a derived foliation on $X$ relative to $k$ whose cotangent complex is $\OO_X$.
By \'etale gluing this contruction can be extended to any 
derived Deligne-Mumford stack $X$ and $\nu \in H^0(X,\TT_{X/k})$. Extension
to derived Artin stacks is trickier, can be done using formal localization techniques
from \cite{cptvv}, but will not be discussed in this book. We can moreover
extend this example by starting with a derivation $\nu : \LL_{X/k} \to \cL$, with $\cL$
an arbitrary line bundle. The previous construction provides 
a derived foliation $F_\nu$ whose cotangent complex is $\cL$ and 
whose sheaf of graded mixed cdga is $Sym_{\OO_X}(\cL[1]) \simeq \OO_X \oplus \cL[1]$ endowed
with the mixed structure defined by $\nu$.} \end{ex}

\begin{ex}\label{ex:groupactions}
\emph{\textbf{Group actions.}
Let $X$ be a derived Artin stack and $G$ a smooth group scheme acting
on $X$. We can form the quotient map $p : X \to Y=[X/G]$, which is a morphism
of derived Artin stacks. The foliation $p^*(0_Y) \in \Fol(X/Y)$ or
rather its image in $\Fol(X/k)$, is called the
\emph{foliation induced from the $G$-action}. Its big cotangent complex
has values on $u : \Spec\, A \to X$ the $A$-dg-module $\LL_{A/Y}$, and thus sits
in an exact triangle
$$\LL_{X/Y} \to \LL_{p^*(0_Y)}^{big} \to \LL_{0_X}^{big}.$$
Applying the quasi-coherator we find that the natural morphism
$\LL_{X/Y} \to \LL_{p^*(0_Y)}$ is an equivalence. Moreover, the cartesian square of
derived stacks
$$\xymatrix{
X \ar[r] \ar[d] & Y \ar[d] \\
\bullet \ar[r] & BG,}$$
implies that $\LL_{X/Y}$ is of the form $g^*\otimes_k \OO_X$, where $g^*=\LL_{*/BG}$ is the dual
Lie algebra of $G$. } \end{ex}

\section{Comparison with derived Lie algebroids}\label{sec:ComparisonwithderivedLiealgebroids}

For a commutative $k$-dg algebra $A$, $\mathsf{T}_A:= \mathsf{Der}_k(A, A)$ is the tangent $A$-dg module consisting of $k$-derivations $\partial: A \to A$. Note that the graded commutator $[\partial, \partial']:= \partial \circ \partial' + (-1)^{|\partial||\partial'|+1} \partial' \circ \partial$ endows $\mathsf{T}_A$ with the structure of a dg-Lie algebra over $k$.

\begin{df} Let $A$ be a commutative $k$-dg algebra. We denote by $\mathsf{dgLieAlgds}_A$ the category of \emph{dg-Lie algebroids over $A$},  whose objects are pairs $(\mathfrak{g}, \rho)$ where $\mathfrak{g}$ is an (a priori, unbounded) $A$-dg module, endowed with a further structure $[\,,\,]_{\mathfrak{g}}: \mathfrak{g}\otimes_k \mathfrak{g} \to \mathfrak{g}$ of $k$-dg Lie algebras, and $\rho:\mathfrak{g} \to \mathsf{T}_A$ is both a morphism of $A$-dg modules and of $k$-dg-Lie algebras, called the \emph{anchor map}, such that
\begin{equation}\label{oideq} 
\left[X, aY\right]_{\mathfrak{g}} = (-1)^{|X||a|}a[X,Y]_{\mathfrak{g}}+ \rho(X)(a)Y.
\end{equation}
Morphisms $(\mathfrak{g}, \rho) \to (\mathfrak{g}', \rho')$ in $\mathsf{dgLieAlgds}_A$ are morphisms $\mathfrak{g} \to \mathfrak{g}'$ of $A$-dg modules and $k$-dg-Lie algebras commuting with the anchor maps.
\end{df}

When no confusion is possible, we will simply write $[\,, \,]$ instead of $[\,, \,]_{\mathfrak{g}}$ for the Lie bracket of a dg-Lie algebroid.

\begin{rmk} \emph{Condition (\ref{oideq}) can be viewed (and remembered) as a sort of a Leibniz rule, relating the anchor map to the failure of the bracket $[\,, \,]_{\mathfrak{g}}$ to be an $A$-linear dg-Lie bracket on $\mathfrak{g}$: the anchor map gives the correction to $A$-linearity of $[\,, \,]_{\mathfrak{g}}$ in the Leibniz rule (\ref{oideq}). }
\end{rmk}

\begin{thm}\label{dglieNuiten}{\cite[Theorem 3.1]{nuit}} Let $A$ be a commutative $k$-dg algebra. The category $\mathsf{dgLieAlgds}_A$ of dg-Lie algebroids over $A$ admits a right proper, tractable, semi-model structure whose fibrations are degreewise surjections, and whose weak equivalences are quasi-isomorphisms.
\end{thm}

For the necessary background on semi-model structures, the reader is referred 
to \cite{nuit}.

\begin{df} Let $A$ be a commutative $k$-dg algebra, and $QA$ its cofibrant replacement. We denote by
\begin{itemize}
\item $\mathbf{dgLieAlgds}_A$ the $\s$-category obtained by Dwyer-Kan localization of $\mathsf{dgLieAlgds}_{QA}$ along its weak-equivalences of Theorem \ref{dglieNuiten}.
\item $\mathbf{dgLieAlgds}^{\emph{perf}, \geq 0}_A$ the full $\s$-subcategory of $\mathbf{dgLieAlgds}_A$ consisting of dg-Lie algebroids $(\mathfrak{g}, \rho)$ where $\mathfrak{g}$ is a perfect $A$-dg-module cohomologically concentrated in degrees $[0, +\s)$.
\end{itemize}
\end{df}

For a dg-Lie algebroid $(\mathfrak{g}, \rho)$ over $A$, we will consider its Chevalley-Eilenberg complex $\mathsf{CE}(\mathfrak{g}, \rho)$, which we now define. Let $\mathrm{Sym}_{A}(\mathfrak{g}[-1])$ be the symmetric algebra over $A$ of $\mathfrak{g}[-1]$ viewed as an $A$-dg-module. 

\begin{df} For a dg-Lie algebroid $(\mathfrak{g}, \rho)$ over $A$, we will denote by $\mathsf{CE}^*(\mathfrak{g}, \rho)$ the $A$-linear dual of the $A$-dg-module $\mathrm{Sym}_{A}(\mathfrak{g}[-1])$, endowed with the following structures:
\begin{itemize}
\item the grading $\mathsf{CE}^*(\mathfrak{g}, \rho)^{(n)}:= \underline{\mathrm{Hom}}_A(\mathrm{Sym}^n_{A}(\mathfrak{g}[-1]), A)\,, \, n\geq 0$.
\item the graded commutative $A$-dg-algebra structure induced by being the $A$-dual of the graded cocommutative coalgebra $\mathrm{Sym}_{A}(\mathfrak{g}[-1])$.
\item the internal differential $$d: (\mathsf{CE}^*(\mathfrak{g}, \rho)^{(n)})^m \to (\mathsf{CE}^*(\mathfrak{g}, \rho)^{(n)})^{m+1}, \, \omega \longmapsto d_A \circ \omega - (-1)^m \omega \circ d_{\mathrm{Sym}^n_{A}(\mathfrak{g}[-1])} $$ i.e. the differential of the $A$-dg module $\underline{\mathrm{Hom}}_A(\mathrm{Sym}^n_{A}(\mathfrak{g}[-1]), A)$. 
\item the mixed structure $$\epsilon^{(n)}_{\mathsf{CE}}: \mathsf{CE}^*(\mathfrak{g}, \rho)^{(n)} \to (\mathsf{CE}^*(\mathfrak{g}, \rho)^{(n+1)})[-1]   $$ is given by (omitting Koszul signs)
$$\epsilon^{(n)}_{\mathsf{CE}}(\omega)(\xi_1, \ldots, \xi_{n+1})= \sum_i  \rho(\xi_i)(\omega(\xi_1, \ldots, \cancel{\xi_i}, \ldots, \xi_{n+1})) - \sum_{i<j} \omega([\xi_i, \xi_j], \xi_1, \ldots, \cancel{\xi_i}, \ldots, \cancel{\xi_j}, \ldots, \xi_{n+1}).$$
\end{itemize}
Endowed with the above structures, $\mathsf{CE}^*(\mathfrak{g}, \rho)$ becomes a graded mixed cdga over $k$.
\end{df}

\begin{rmk}\emph{Usually, $\mathsf{CE}^*(\mathfrak{g}, \rho)$ is regarded as a cdga by taking the total 
differential, given by the sum of the internal differential and the mixed one. However, as we stressed 
out in several other points, it is quite important for our purposes to keep the two differentials 
separated and view it as a graded mixed cdga.}
\end{rmk}

Since $\mathbf{DR}(-/k)$ is left adjoint to the weight $0$ part functor (see Section 
\ref{sec:GradedmixedcomplexesandderiveddeRhamtheory}), and $$\mathsf{CE}^*(\mathfrak{g}, \rho)^{(0)}= 
\underline{\mathrm{Hom}}_A(A, A),$$  the identity map $A \to A$ induces a morphism of graded mixed 
cdga's  $\mathbf{DR}(A/k) \to \mathsf{CE}^*(\mathfrak{g}, \rho)$. Therefore, if $\mathfrak{g}$ is a 
perfect $A$-dg-module, this yields the structure of a derived foliation on $\mathsf{CE}^*(\mathfrak{g}, 
\rho)$ (Definition \ref{dII-1}). \\

The following theorem expresses the precise comparison between perfect foliations and perfect dg-Lie algebroids.

\begin{thm}\label{tII-4}
Let $A$ be a connective $k$-linear cdga of finite presentation over $k$.
The $\s$-functor sending a dg-Lie algebroid $(\mathfrak{g}, \rho)$ over $A$ to its
Chevalley complex $\mathsf{CE}^*(\mathfrak{g}, \rho)$ induces a equivalence
$$\mathsf{CE}^* : (\mathbf{dgLieAlgds}^{\mathrm{perf}}_A)^{op} \simeq \Fol^{\mathrm{perf}}(A/k),$$
where $\mathbf{dgLieAlgds}^{\mathrm{perf}}_A$ is the full sub-$\s$-category of dg-Lie algebroids
which are perfect as an $A$-dg-module, and $\Fol^{\mathrm{perf}}(A/k)$ is the $\s$-category of perfect 
derived foliations over $A$.
\end{thm}

\begin{proof}
This is a consequence of \cite[Theorem 5.15]{jfuarxiv}
\end{proof}

As an immediate consequence, we get the following corollary which relates perfect derived foliations
over the point $\Spec\, k$ and dg-Lie algebras over $k$ which are perfect as $k$-dg-modules.
We note also that the above theorem extends to almost perfect derived foliations under the
condition that $A$ is furthermore cohomogically bounded. We can thus extract the most 
important case of application for us, which is the following corollary.

\begin{cor}\label{ctII-4}
The Chevalley complex construction induces an equivalence of $\s$-categories
$$\dglie_k^{dap} \simeq \Fol^{ap}(k/k)$$
between dual almost perfect dg-Lie algebras over $k$ and almost perfect derived foliations over $k$.
\end{cor}

In the previous corollary \emph{dual almost perfect} means that the $k$-linear dual is
almost perfect. As $k$ is a field this condition is equivalent to the condition 
of being cohomologically bounded on the left and with finite dimensional individual cohomology 
spaces.

\section{Formal structures}\label{sec:Formalstructures}

All along this section, we will assume that our base ring $k$ is a field (of characteric zero).
We will be interested in derived foliations restricted to formal 
derived stacks, and their interpretations in terms of dg-Lie algebras over $k$. These results are
not optimal and replacing $k$ be a more general base is also possible. However, the
general notion of formal derived stacks over general bases is out of the scope of this
book, and we refer the interested reader to 
\cite{brantner2025formalintegrationderivedfoliations} for recent and pretty optimal results in that
direction.

\subsection{Derived stacks and formal moduli problems}\label{subsec:Derivedstacksandformualmoduliproblems}

In this section we recall the fundamental results of \cite{lupmf,lupmf2} (see also
\cite[\S IV]{SAG} and \cite{tobour}), relating pointed derived formal stacks
and dg-lie algebras. For this, we first introduce
$\dgart_k$, the $\s$-category of commutative, connective, Artinian 
local dg-algebras over $k$. By definition, this is the full sub-$\s$-category
of $\cdga_k$ consisting of objects $A$ satisfying the following conditions:
\begin{itemize}
    \item $H^i(A)=0$ for all $i>0$
    \item $H^0(A)$ is a local Artinian $k$-algebra
    \item $H^*(A)$ is finite dimensional over $k$.
\end{itemize}

Note that the third condition implies in particular that
$H^i(A)$ is finite dimensional for all $i$, but also that 
$H^i(A)=0$ when $i<<0$. 
In a similar manner, we denote by $\dgart_k^*$ the $\s$-catgeory 
of pointed objects in $\dgart_k$, or equivalently the $\s$-category 
of augmented connective Artinian dg-algebras. Note that in particular, 
if $A \in \dgart_k^*$ the residue field of $H^0(A)$ is then automatically $k$ itself. 

\begin{df}\label{dII-7}
A \emph{formal moduli problem} is an $\s$-functor
$$F : \dgart_k^* \to \Top$$
satisfying the following conditions.
\begin{enumerate}
    \item The space $F(k)$ is contractible.
    
    \item For any cartesian square in $\dgart_k^*$
    $$\xymatrix{ B' \ar[r] \ar[d] & A' \ar[d] \\
    B \ar[r] & A}$$
    such that $H^0(A') \to H^0(A)$ is surjective, the 
    square 
    $$\xymatrix{ F(B') \ar[r] \ar[d] & F(A') \ar[d] \\
    F(B) \ar[r] & F(A)}$$
    is cartesian in $\Top$.

\end{enumerate}
Formal moduli problems over the field $k$
form an $\s$-category $\FMP_k$, which by definition is a full 
sub-$\s$-category of $Fun^\s(\dgart_k^*,\Top)$.
\end{df}

Formal moduli problems are related to pointed derived stacks by means 
of an adjunction, whose right adjoint is the formal completion functor. 
Let $\dSt_k^*$ be the $\s$-category of pointed derived stacks. An object 
in $\dSt_k^*$ will be denoted as $(X,x)$, where $X \in \dSt_k$ and 
$x : \Spec\, k \to X$ is the chosen base point. To any such object
$(X,x)$ we associated an $\s$-functor 
$$(X,x)^\sim : \dgart_k^* \to \Top,$$
by sending $A \in \dgart_k^*$ to $Map_{\dSt_k^*}(\Spec\, A,(X,x))$. 
Here, $\Spec\, A$ is considered as an object in $\dSt_k^*$ by its
augmentation $A \to k$. In general, $(X,x)^\sim$ is not 
a formal moduli problem (in general, it does not satisfy condition (2)). However, the natural inclusion 
$$\FMP_k \hookrightarrow Fun^\s(\dgart_k^*,\Top)$$
is the right adjoint (see \cite{lupmf,SAG,tobour}) of an adjunction, 
and the image of $(X,x)^\sim$ by the left adjoint will be denoted by
$(X,x)^\wedge \in \FMP_k$.

This construction defines an $\s$-functor $(-)^\wedge : \dSt_k^* \to \FMP_k$
called the \emph{formal completion}. It is easy to 
construction a left adjoint as follows. We start by the full embedding
$$\Spec\, : (\dgart_k^*)^{op} \hookrightarrow \dSt_k^*$$
and we left Kan extend it to get an $\s$-functor
$$j : \FMP_k \hookrightarrow Fun^\s(\dgart_k^*,\Top) \longrightarrow \dSt_k^*.$$
We thus have defined an adjunction of $\s$-categories
$$j : \FMP_k \rightleftarrows \dSt_k^* : (-)^\wedge.$$
As a side comment, the left adjoint of the inclusion 
$\FMP_k \hookrightarrow Fun^\s(\dgart_k^*,\Top)$ is very non-explicit, 
and thus so is the $\s$-functor $(-)^\wedge$ in general. However, 
this $\s$-functor becomes very simple when applied to 
derived Artin stacks. Indeed, with the notations above, 
when $(X,x)$ is a pointed derived Artin stack, the 
$\s$-functor $(X,x)^\sim : \dgart_k^* \to \Top$ is in fact already
a formal moduli problem, as this follows easily from the
fact that derived Artin stack possesses an obstruction theory 
(see for instance \cite[\S 1.4]{hagII}). Therefore,
$(X,x)^\sim \simeq (X,x)^\wedge$, and we thus have 
for any $A \in \dgart_k^*$, a cartesian square in $\Top$
$$\xymatrix{
(X,x)^\wedge(A) \ar[r] \ar[d] & X(A) \ar[d] \\
\bullet \ar[r]_-{x} & X(k).}$$

We can now combine the previous adjunction with the
main equivalence $\FMP_k \simeq \dglie_k $ of \cite{lupmf} in order to obtain an adjunction
between $\dSt_k^*$ and the $\s$-category of dg-lie algebras. We thus
have constructed an adjunction\footnote{The functor $\MC$ is so named after L. Maurer and E. Cartan.}
$$\MC : \dglie_k \leftrightarrows \dSt_k^* : \ell_-.$$
The $\s$-functor $\ell_-$ sends a pointed derived stack
$(X,x)$ to the dg-lie algebra $\ell_{(X,x)}$ controlling the formal
moduli problem $(X,x)^\wedge$, and $\MC$ is its left adjoint. 
We have then the following important observation.

\begin{prop}\label{pII-6}
The $\s$-functor
$$\MC : \dglie_k \longrightarrow \dSt_k^*$$
is fully faithful.
\end{prop}

\textit{Proof.} As a first reduction, we can 
start by noticing that $k \mapsto \dglie_k$ and
$k \mapsto \dSt_k^*$ are both stacks for the fpqc topology on $k$.
Therefore, the $\s$-functor $\MC$ being compatible with the change of base field $k$, 
we can assume that $k$ is algebraically closed. 

By construction, the $\s$-functor $\MC$ is equivalent to the composition
$$j : \xymatrix{ \FMP_k \ar[r] & Fun^\s(\dgart_k^*,\Top) \ar^-{\phi}[r] & \dSt_k^*.}$$
Here, the first $\s$-functor is the canonical embedding and the second
$\s$-functor $$\phi : Fun^\s(\dgart_k^*,\Top) \to \dSt_k^*$$ is the composition
of the left Kan extension
$$i_! : Fun^\s(\dgart_k^*,\Top) \hookrightarrow */Fun^\s(\cdga_k^{\leq 0},\Top)$$
along  the inclusion $i : \dgart_k \hookrightarrow \cdga_k^{\leq 0}$,
with the associated stack $\s$-functor $$*/Fun^\s(\cdga^{\leq 0}_k,\Top) \to \dSt_k^* .$$
It is therefore enough to show that $\phi$ is fully faithful. The right adjoint  
to $\phi$ is the restriction $\psi$ sending a pointed stack $F : \cdga_k^{\leq 0} \to \Top$ 
to its restriction to Artinian dg-algebras. However, 
as $k$ is algebraically closed, any Artinian local dg-algebra is 
trivial for the \'etale topology, and thus for any prestack $F$, with
associated stack $F'$, the natural map $F \to F'$ induces equivalences
$F(A) \simeq F'(A)$ for any Artinian local $k$-algebra $A$. As a consequence, 
for any $F \in Fun^\s(\dgart_k^*,\Top)$, the adjunction morphism
$F \to \psi \phi(F)$, is simply equivalent to the adjunction morphims
$F \to i^*i_!(F)$, which is an equivalence as the left Kan extension
$$i_! :  Fun^\s(\dgart_k^*,\Top) \hookrightarrow */Fun^\s(\cdga_k^{\leq 0},\Top)$$
is fully faithful. \hfill $\Box$ \\

We will mainly use the previous adjunction when restricted to derived Artin 
stacks. We note that for such a pointed derived Artin stack $(X,x)$, the
underlying complex of the dg-lie algebra $\ell_{(X,x)}$ is 
functorialy quasi-isomorphic to $\TT_{X,x}[-1]$, the shifted tangent complex of $X$ at $x$.
This follows directly from a simple unfolding of the various definitions.

Finally, because of Proposition \ref{pII-6} we will often consider formal
moduli problems as fully embedded in  pointed derived stacks. In particular, 
we will often use the notation $(X,x)^\wedge$ to denote also the
pointed derived stack $\MC((X,x)^\wedge)$ and call it the \emph{formal 
completion of $X$ at $x$}. This formal completion $(X,x)^\wedge$ can be
in fact described in the usual manner on smooth atlases of $X$. 
Assuming first that $X=\Spec\, A$ is a pointed affine derived scheme, the point 
provides an augmentation $e : A \to k$, and we can define a pro-object in 
$\dgart_k^*$, called the formal completion of $A$ along $e$, by
$$\widehat{A}:="\lim_{A \to A'}" A',$$
where $A'$ runs in the $\s$-category of morphism of augmented
cdga morphism $A \to A'$, where $A'$ is a local Artinian cdga. 
When $A$ is an underived $k$-algebra of finite type, then 
it can be shown that this pro-object is equivalent to the usual 
formal completion $"\lim_n" A/m^n$ where $m \subset A$ is the maximal ideal
corresponding to the augmentation $e$ (see Proposition \ref{pA-3-1}). The pro-object 
$\widehat{A}$ pro-represents $(X,x)^\wedge$ is the sense that one has
$$(X,x)^\wedge \simeq \Spf(\underline{\widehat{A}}) = \mathrm{colim}\, \Spec\, A'.$$
In general, up to base changing to a finite extension $k'/k$, 
we can find a smooth morphism $p : U \to X$ of pointed derived stacks for which
$U$ is affine and such that $x$ lies in the image of $p$. We can form
the nerve $U_* \to X$ of $p$. We then have 
$$(X,x)^\wedge \simeq \mathrm{colim}_n (U_n,x)^\wedge.$$
We can this way compute the formal completion of $(X,x)^\wedge$ inductively on 
the geometricity of $X$. We could also replace the nerve of $p$ here by 
any pointed smooth hypercovering (locally around $x$) $U_* \to X$ where all 
$U_i$ are affine, and get this way a formula for $(X,x)^\wedge$ as the colimit
of formal completions $(U_n,x)^\wedge$ which can individually be described as the formal
spectrum of a pro-Artinian cdga. \\

\subsection{Derived foliations on formal moduli problems}\label{subsec:Derivedfoliationsonformalmoduliproblems}

In this part we show how derived foliations can be defined over
formal moduli problems. For this we simply use the fully faithful $\s$-functor defined
in the previous section
$j : \FMP_k \longrightarrow \dSt_k^*$. 

\begin{df}\label{dII-9}
For a formal moduli problem $X \in \FMP_k$ we set 
$$\Fol(X/k):=\Fol(j(X)/k).$$
When the base field is clear we will simply use the notation $\Fol(X)$.
\end{df}

Note that $X \mapsto \Fol(X)$ naturally defines an $\s$-functor $\Fol^{op} : \FMP_k \to \scat$.
Unfolding definitions and constructions we can give the following explicit formula
$$\Fol(X)\simeq \lim_{\Spec\, A \to X}\Fol(A),$$
where $\Spec\, A$ runs over the $\s$-category of $A \in \dgart_k^*$ endowed with 
a morphism $\Spec\, A \to X$ (concretely the comma $\s$-category $\Spec\, /X$
where $\Spec : \dgart_k^{op} \to \FMP_k$ is the Yoneda embedding). \\

Note that, for any pointed
derived stack $(X,x) \in \dSt_k^*$, there exists a canonical $\s$-functor
$$\Fol(X) \longrightarrow \Fol((X,x)^\wedge).$$
Indeed, by construction of the adjunction between $j$ and $(-)^\wedge$, we
have a natural adjunction morphism of derived stack $i_x : (X,x)^\wedge \to X$, 
and thus an induced pull-back $\s$-functor
$i_x^* : \Fol(X) \to \Fol((X,x)^\wedge)$. 

\begin{df}\label{dII-10}
For a pointed derived stack $(X,x)$
the $\s$-functor $i^*_x : \Fol(X) \to \Fol((X,x)^\wedge)$ above is
the \emph{formalization at $x$}, or \emph{formal completion at $x$}. It is denoted by 
$\F \mapsto \hat{\F}_x$.
\end{df}

Proposition \ref{pII-9} below shows that 
derived foliations which are \emph{almost perfect}, can be described explicitly 
on the formal completion of an affine derived scheme. Before explaining this, let 
$A$ be a connective cdga with an augmentation $A \to k$, and assume that 
$A$ is almost of finite presentation (see appendix \S \ref{secapp:Derivedschemesandstacks}).
Wet set $X=\Spec\, A$ endowed with its natural point $x \in X(k)$ associated to the augmentation.
We denote by
$\underline{\hat{A}}_x$ (or simply by $\underline{\hat{A}}$)  
the formal completion of $A$ at the augmentation $x : A \to k$. Recall (see appendix \S \ref{secapp:Formalcompletion})
that this is the pro-object in $\dgart_k^*$ given by 
$"\lim_{A \to A'}A'$, where $A'$ runs in the $\s$-category of all 
local Artinian cdga with an augmentation preserving morphism $A \to A'$. 
We let $X=\Spec\, A$ and $(X,x)^\wedge:=\Spf(\underline{\hat{A}})$ be the corresponding 
derived affine schemes and formal completion. 

We define the graded mixed cdga $\hat{\DR}(\underline{\hat{A}}/k)$, called the \emph{formal
de Rham algebra of $\underline{\hat{A}}$}, by the following limit
$$\hat{\DR}(\underline{\hat{A}}/k) := \lim_{A \to A'}\DR(A'/k).$$
By Corollary \ref{cpA-3-1-2} we know that, as graded cdga's we have
$$\hat{\DR}(\underline{\hat{A}}/k) \simeq 
Sym_{\hat{A}}(\hat{\LL}_{A/k}[1]),$$
where $\hat{A}=\lim_{A \to A'}A'$ is the realization of the pro-object
$\underline{\hat{A}}$, and $\hat{\LL}_{A/k}=\lim_{A \to A'}\LL_{A'/k}$ is the
formal cotangent complex of $\underline{\hat{A}}$.
Now, the limit $\s$-functor provides a natural $\s$-functor
$$\psi : \Fol((X,x)^\wedge)=\lim_{A \to A'}\Fol(Spec\, A') \longrightarrow 
(\hat{\DR}(\hat{A}/k)-\egrcdga_k)^{op},$$
from derived foliations on the formal completion of $X$ at $x$ to graded
mixed commutative $\hat{\DR}(\hat{A}/k)$-algebras. Explicitly, the $\s$-functor $\psi$
sends an object $\{\F_A'\} \in \lim_{A \to A'}\Fol(Spec\, A')$ to 
$\lim_{A \to A'}\DR(\F_A')$, as an algebra over $\lim_{A \to A'}\DR(A'/k)$.

\begin{prop}\label{pII-9}
Let $(X,x) =\Spec\, A$ be a pointed derived affine scheme almost of finite presentation over $k$.
The $\s$-functor
$$\psi : \Fol((X,x)^\wedge)=\lim_{A \to A'}\Fol(Spec\, A') \longrightarrow 
(\hat{\DR}(\hat{A}/k)-\egrcdga_k)^{op}$$
constructed above, is fully faithful when restricted to 
almost perfect derived foliations. Its essential image consists of 
all graded mixed commutative $\hat{\DR}(\hat{A}/k)$-algebras whose underlying graded cdga'a are
of the form $Sym_{\hat{A}}(E)$, with $E$ an almost perfect
$\hat{A}$-dg-module.
\end{prop}

\textit{Proof.} The limit $\s$-functor 
$$\lim_{A \to A'}(\DR(A'/k)-\egrcdga_k) \longrightarrow 
\hat{\DR}(\hat{A}/k)-\egrcdga_k$$
admits a left adjoint $\phi$, which sends a graded mixed commutative 
$\hat{\DR}(\hat{A}/k)$-algbera $B$ to the family of objects 
$\hat{\DR}(A'/k)\otimes_{\hat{\DR}(\hat{A}/k)}B$, together with their natural transition morphisms. 

We denote by $\C \subset \hat{\DR}(\hat{A}/k)-\egrcdga$ be the full sub-$\s$-category 
consisting of $B$ which, as graded $\hat{\DR}(\hat{A}/k)$-cdga's, are of the form
$Sym_{\hat{A}}(E)$ for some $E\in \APerf(\hat{A})$. We first observe that the essential image 
of $\phi$ is contained in $\C$.  
Let $\F \in\lim_{A \to A'}\Fol^{ap}Spec\, A') $ be a derived foliation 
on $(X,x)^\wedge$ which is assumed to be almost perfect It consists of a family of $\DR(\F_{A'})$ of graded mixed
$\DR(A'/k)$-cdga, which, as graded $\DR(A'/k)$-algebras, are of the form 
$Sym_{A'}(\LL_{\F_{A'}}[1])$ for an almost perfect $A'$-dg-module
$\LL_{\F_{A'}} \in \APerf(A')$. We set $\LL_\F:=\lim_{A \to A'} \LL_{\F_A'}$, which 
is an $\hat{A}$-module, and let $\N^*_{\F}$ be the fiber of the natural morphism
$\hat{\LL}_{A}\otimes_{\hat{A}} A' \to \LL_\F$.

We consider the natural commutative diagram of exact triangles
of pro-$\underline{\hat{A}}$-dg-modules
$$\xymatrix{``\lim_{A\to A'}"(\hat{\LL}_{A}\otimes_{\hat{A}} A') \ar[r] \ar[d] & ``\lim_{A \to A'}"(\LL_{\F}\otimes_A A')\ar[d] \ar[r] & ``\lim_{A \to A'}" (\N^*_{\F}\otimes_A A')[1] \ar[d] \\
``\lim_{A\to A'}"\LL_{A'} \ar[r] & ``\lim_{A\to A'}"\LL_{\F_{A'}} \ar[r] & 
``\lim_{A\to A'}"\N^*_{\F_{A'}}[1].}$$
The vertical left arrow is an equivalence of pro-objects thanks to theorem \ref{cpA-3-1-2}. 
The right vertical arrow
is also an equivalence (and even a levelwise equivalence) thanks to the fact that the natural $\s$-functor
$$\lim_{A \to A'} \APerf(A') \longrightarrow \APerf(\hat{A})$$
is an equivalence of $\s$-categories (see appendix \ref{pA-4-1})
and its inverse is given by 
$E \mapsto \lim_{A \to A'}E\otimes_{\hat{A}} A'$. Therefore, the natural morphism
$``\lim_{A \to A'}"(\LL_{\F}\otimes_{\hat{A}} A') \to ``\lim_{A\to A'}"\LL_{\F_{A'}}$ is again an 
equivalence of pro-$\underline{\hat{A}}$-dg-modules. As a consequence, we see that for all $p$, 
$$``\lim_{A \to A'}"(\wedge_{\hat{A}}^p(\LL_{\F})\otimes_{\hat{A}} A') \to ``\lim_{A\to 
A'}"\wedge^p_{A'}\LL_{\F_{A'}}$$
is an equivalence of pro-objects. This implies in particular that $\psi(\F)$ lies indeed in the
sub-$\s$-category $\C$. 

We thus have seen that the adjunction
$$\phi : \xymatrix{\hat{\DR}(\hat{A}/k)-\egrcdga_k  \ar@<2pt>[r]_{}  & 
\lim_{A \to A'}\DR(A'/k)-\egrcdga_k  \ar@<2pt>[l]_{}} : \psi$$
restricts to an adjunction
$$\phi : \xymatrix{\C \ar@<2pt>[r]_{}  & 
\lim_{A \to A'}(\Fol^{ap}(\Spec\, A'/k)) \ar@<2pt>[l]_{}} : \psi.$$
To finish the proof of the proposition we need to check that the 
unit and counit of this adjunction are equivalence. Let $B \in \C$, 
and $\LL:=B^{(1)}[1]$ be its weight one part. It is an object in 
$\APerf(\hat{A})$ by assumption, and thus if we denote
by $\N_\F^*$ the fiber of $\widehat{\LL_A} \to \LL$, we have that 
$\N_\F^*$ is an almost perfect $\hat{A}$-dg-module. The adjunction
morphism $B \psi(\phi(B))$, on the weight one part, is equivalent to the canonical
morphism
$\LL \to \lim_{A\to A'}\LL \otimes_{\hat{A}}A'$. This is an equivalence thanks to 
proposition \ref{pA-4-1}. This shows
that $\phi$ is fully faithful. But $\psi$ is also conservative. Indeed, if 
a morphism $f : \F \to \F'$ in $\lim_{A \to A'}\Fol^{ap}(Spec\, A')$ induces an equivalence
after applying $\phi$, then in particular the induced morphism on conormal complexes
$\lim_{A\to A'}\N^*_{\F'_{A'}} \to \lim_{A\to A'}\N^*_{\F_{A'}}$ is an equivalence.
The limit $\s$-functor $\lim_{A\to A'}\APerf(A') \longrightarrow \APerf(\hat{A})$
is an equivalence, and is thus conservative, we have therefore that for each $A'$ the induced morphism
$\N^*_{\F'_{A'}} \to \N^*_{\F_{A'}}$ is an equivalence, and thus that $\F_{A'} \to \F'_{A'}$
is an equivalence of derived foliations on $\Spec\, A'$. This shows
that $\psi$ is conservative and thus the proof of the Proposition.
\hfill $\Box$ \\

Finally, we have the formal version of Theorem \ref{tII-1}. For this, recall that $X^{art}$ is the
$\s$-site of pointed morphisms $\Spec\, A \to X$ with $A$ a local Artinian connective cdga
(see \S \ref{secapp:Cotangentcomplexesofformalcompletion}). It comes equipped with the canonical stack of cdga's $\OO_{X^{art}}$, as well 
as its $\s$-category of $\OO_{X^{art}}$-dg-modules. We have the full sub-$\s$-category
$\QCoh(X^{art}) \subset \OO_{X^{art}}-\dg$ of quasi-coherent complexes, which can canonically be
identified with $\QCoh((X,x)^{\wedge})$. Any derived foliation $\F \in \Fol((X,x)^\wedge)$
possesses a big cotangent complex $\LL_{\F/k}^{big} \in \OO_{X^{art}}-\dg$, and 
a cotangent complex $\LL_{F/k} \in \QCoh(X^{art})$. By definition, they are related by the formula
$\LL_{\F/k} = (\LL_{\F/k}^{big})^{qcoh},$
where $(-)^{qcoh} : \OO_{X^{art}}-\dg \longrightarrow \QCoh(X^{art})$ is the right adjoint
to the natural inclusion $\s$-functor.

\begin{prop}\label{pII-7}
Let $(X,x)$ be a pointed derived stack locally almost of finite presentation and $(X,x)^{\wedge}$ 
its formal completion. Then, for any almost perfect derived foliation 
$\F \in \Fol^{ap}((X,x)^{\wedge})$ such that $\N_\F^*$ is perfect, 
and all $i\geq 0$ the canonical morphism
$$\wedge^i_{\OO_{X^{art}}}\LL_{\F} \longrightarrow 
(\wedge^i_{\OO_{X^{art}}}\LL_{\F}^{big})^{qcoh}$$
is a quasi-isomorphism.
\end{prop}

\textit{Proof.} This is exactly the same proof as that of Theorem \ref{tII-1}, using Theorem
\ref{tA-6-1} instead of Theorem \ref{tA-5-1}. We leave the details to the reader. \hfill $\Box$ \\

\subsection{Formal integrability}\label{subsec:Formalintegrability}

We will now 
show that almost perfect derived foliations
on the formal completion of a derived Artin stack locally almost of finite presentation can be
described purely in algebraic terms, using dg-Lie algebras. This classification 
result also implies that these derived foliations are always
formally integrable at each point.

Recall that the conormal complex of a derived foliation $\F$ is defined by the exact triangle
$$\xymatrix{\mathcal{N}^*_\F \ar[r] & \LL_X \ar[r] & \LL_\F. }$$
An important remark is that 
$\mathcal{N}^*_\F$ can be defined using either 
the cotangent quasi-coherent complexes or the big cotangent 
complexes. Indeed, the commutative square
$$\xymatrix{
\LL_X \ar[r] \ar[d] & \LL_\F \ar[d] \\
\LL_X^{big} \ar[r] & \LL_\F^{big}}$$
is always cartesian (see Lemma \ref{ltII-1}) so that the fibers of the horizontal arrows agree. As a result, $\mathcal{N}_\F^* $ always makes sense as an object in $\mathbf{QCoh}(X)$, 
even if $X$ is merely a derived stack (i.e. not necessarily Artin).
Another important fact is that conormal complexes of derived foliations
are stable by pull-backs (as opposed to cotangent complexes themselves), and therefore
almost perfect derived foliations form a full substack $\Fol^{ap} \subset 
\Fol$ of the stack of derived foliations. In particular, for any formal
moduli problem $X \in \FMP$, we get a full sub-$\s$-category
$\Fol^{ap}(X) \subset \Fol(X)$, simply defined by $\Fol^{ap}(j(X))$ (Definition \ref{dII-9}).  \\

Let now $(X,x)$ be a pointed derived Artin stack locally almost of finite presentation 
and $\F \in \Fol((X,x)^\wedge)$
a derived foliation on the formal completion of $X$ at $x$. The natural base point 
$x : \Spec\, k \to (X,x)^\wedge$ provides $x^*(\F) \in \Fol(*)$. Moreover, 
the tautological derived foliation $0_X$ can be formally completed to
$\hat{0}_{X_x} \in \Fol((X,x)^\wedge$. 

\begin{lem}\label{ltII-2}
The object $\hat{0}_{X_x}$ is the initial object in $\Fol((X,x)^\wedge$.
\end{lem}

\textit{Proof.} The formal completion construction $\F \mapsto \hat{\F}_x$ (Definition \ref{dII-10}) being obtained by restriction, 
the conormal complex $\cN_{\hat{F}_x}$ is simply the restriction of $\cN_\F$. Applied to 
the case of $\F=0_{X}$, we see that the conormal complex of $\hat{0}_{X_x}$ is identified
with $\LL_{(X,x)^{\wedge}}$, via the canonical morphism. Equivalently, 
$\LL_{\hat{0}_{X,x}} \simeq 0$. We then finish the proof using the exact same argument 
as in the non-formal case (see example \ref{ex:tautologicalfoliationsglobal}, by replacing the 
big site $\dAff/X$ with 
$\dgart^*_k/X$, and using Proposition \ref{pII-5} 
(and Proposition \ref{pII-7} instead of Theorem \ref{tII-1}).
\hfill $\Box$ \\

Lemma \ref{ltII-2} provides 
a canonical morphism $x^*(0_X) \to x^*(\F)$, of derived foliations over
$\Spec\, k$. This defines an $\s$-functor
$\psi : \Fol((X,x)^{\wedge}) \to x^*(0_X)/\Fol(*)$, which can be 
restricted to perfect derived foliations 
$$\psi : \Fol^{ap}((X,x)^{\wedge}) \to x^*(0_X)/\Fol^{ap}(*).$$

\begin{thm}\label{tII-2}
With the above notations and assumptions,
the $\s$-functor
$$\psi : \Fol^{ap}((X,x)^{\wedge}) \to x^*(0_X)/\Fol^{ap}(*)$$
is an equivalence of $\s$-categories.
\end{thm}

Before proving this theorem, we state its most important corollary. 
The target $\s$-category of the equivalence of Theorem \ref{tII-2} can thus
be identified, according to corollary \ref{ctII-4}, 
with the comma $\s$-category of dual almost perfect dg-Lie algebras
$t_X / \dglie^{dap}_k$. Here $t_X$ is the dg-Lie algebra corresponding to
$x^*(0_X)$ according to the equivalence of \ref{ctII-4}. As this equivalence is functorial, 
it is easy to check that $x^*(0_X)$ is the dg-Lie algebra
obtained by pulling back the dg-Lie algebroid $\TT_X$ to the point $x$, 
which is nothing else that $\ell_{(X,x)}$, corresponding to the
formal moduli problem $(X,x)^\wedge$. The main corollary of Theorem \ref{tII-2} thus reads
as follows.

\begin{cor}\label{ctII-2}
Let $(X,x)$ be a pointed derived Artin stack locally almost of finite presentation. Then, 
the $\s$-functor $\psi$ combined with the equivalence of Corollary \ref{ctII-4} induces 
an equivalence of $\s$-categories
$$\Fol^{ap}((X,x)^\wedge) \simeq \ell_{(X,x)}/\dglie_k^{dap},$$
between almost perfect derived foliations on the formal completion 
$(X,x)^\wedge$ and the comma $\s$-category of dual almost perfect dg-Lie algebras under
$\ell_{(X,x)}$. This equivalence will be denoted by 
$$\F \mapsto \ell_{\F}.$$
\end{cor}

The above corollary states formal integrability 
of derived foliations. Indeed, for any $\F \in \Fol^{ap}(X)$, its
formal completion at a $k$-point $x \in X(k)$ provides by the equivalence
of the corollary, a morphism of dg-Lie algebras $\ell_{(X,x)} \to \ell_{\hat{\F}_x}$.
By applying the $\s$-functor $\MC : \dglie_k \to \dSt_k^*$ we get a morphism
$$u : \MC(\ell_{(X,x)}) \simeq (X,x)^\wedge \to \MC(\ell_{\hat{\F}_x}).$$
The pointed derived stack $\MC(\ell_{\hat{\F}_x})$ is, by definition, the
\emph{formal leaf space} of $\F$ at $x$. It is possible to show that 
the morphism $u$ integrates $\hat{\F}_x$ in the sense that 
$\hat{F}_x \simeq u^*(0)$, but we will not do this here. We will however 
prove this in the next section,  
in the specific case where $X$ is a smooth variety and 
$\F$ is rigid and quasi-smooth (see \S \ref{subsec:Transverselysmoothandrigidderivedfoliations}). 
For now, we gather some terminology 
related to Corollary \ref{ctII-2}. \\

\begin{df}\label{dII-11}
Let $(X,x)$ be a pointed derived Artin stack locally almost of finite presentation
and $\F \in \Fol^{ap}(X)$ be an almost perfect derived foliation. 

\begin{enumerate}
    \item The \emph{formal leaf space of $\F$ at $x$} is the formal
    moduli problem corresponding to the dg-Lie algebra $\ell_{\hat{\F}_x}$.
    It is denote by $\hat{\cM}_x(\F)$.
    
    \item The \emph{formal leaf of $\F$ passing through $x$} is the formal
    moduli problem corresponding to the dg-Lie algebra which is the
    fiber of the canonical morphism $\ell_{(X,x)} \to \ell_{\hat{\F}_x}$. It is denoted by 
    $\hat{\cL}_x(\F)$.
\end{enumerate}
\end{df}

We finish this section by proving Theorem \ref{tII-2}. \\

\textit{Proof of Theorem \ref{tII-2}.} As a first reduction, by a standard Galois descent argument, 
we may assume that $k$ is algebraically closed. We then proceed by induction on the geometricity 
of $X$. To start with, assume that $X$ is affine, and thus of the form $X=\Spec\, A$
for $A$ a connective cdga almost of finite presentation over $k$. We let $\underline{\hat{A}}$
be the formal completion of $A$ at the point $x : \Spec\, A \to k$, and
consider $\hat{\DR}(\hat{A}/k)$ the associated formal de Rham algebra. By Proposition \ref{pII-9}
$\Fol^{ap}((X,x)^{\wedge})$ is naturally equivalent to the (opposite) $\s$-category
of graded mixed $\hat{\DR}(\hat{A}/k)$-cdga's which are equivalent, as
graded $\hat{\DR}(\hat{A}/k)$-cdga, to $Sym_{\hat{A}}(E)$ for $E$ an almost perfect
$\hat{A}$-dg-module. 
Tracing back the various equivalences of $\s$-categories, we obtain a canonical 
equivalence of graded mixed cdga's
$$\DR(x^*(0_X)) \simeq \hat{A} \otimes_{\hat{\DR}(\hat{A}/k)}k,$$
and the $\s$-functor of the theorem $\psi$ is then induced by
$$\psi : \hat{\DR}(\hat{A}/k)-\egrcdga \longrightarrow \egrcdga/(A \otimes_{\hat{\DR}(\hat{A}/k)}k)$$
which sends $\hat{\DR}(\hat{A}/k) \to B$ to $B\otimes_{\hat{\DR}(\hat{A}/k)}k$
together with its natural augmentation to $\hat{A} \otimes_{\hat{\DR}(\hat{A}/k)}k$ induced
by the natural morphism $B \to \hat{A}$. The $\s$-functor $\psi$
is thus subsumed by the following commutative diagrams of graded mixed cdga's with cocartesian squares
$$\xymatrix{
\hat{\DR}(\hat{A}/k) \ar[r] \ar[d] & B \ar[d] \ar[r] & \hat{A} \ar[d] \\
k \ar[r] & \psi(B) \ar[r] & \hat{A}\otimes_{\hat{\DR}(\hat{A}/k)}k.
}$$

We construct an $\s$-functor $\phi$ in the other direction by using de Rham
cohomology of the first kind explained in Section \S \ref{secapp:InternalderiveddeRhamtheory}. 
For this, we first prove the following lemma.

\begin{lem}\label{ltII-2-2}
The canonical morphism 
$k \to \hat{\DR}(\hat{A}/k)$ induces an equivalence after taking realizations
$$k \simeq |\hat{\DR}(\hat{A}/k)|.$$
\end{lem}

\textit{Proof of the lemma.} By definition 
$|\hat{\DR}(\hat{A}/k)|$ is the derived de Rham cohomology of the formal
completion of $\Spec\, A$ at $x$, and the lemma simply states that derived de Rham
cohomology is trivial on formal completion. 
By definition $\hat{\DR}(\hat{A}/k) \simeq  \lim_{A \to A'}\DR(A'/k)$. As $|-|$ is a right 
adjoint it commutes with limits, and thus we are reduced to show that, for any 
local augmented dg-Artinian algebra $A$, the natural morphism $k \to \CDR(A/k)$ is
a quasi-isomorphism. 
The lemma then follows from proposition \ref{pA-8-1}.
\hfill $\Box$

\begin{lem}\label{ltII-2-3}
There exists a canonical equivalence of graded mixed cdga's
$$\hat{\DR}(\hat{A}/k) \simeq \CDR^{I}(\hat{A}\otimes_{\hat{\DR}(\hat{A}/k)}k/k).$$
\end{lem}

\textit{Proof of the lemma.} Let $D:=\hat{\DR}(\hat{A}/k)$. 
By base change, we have a canonical equivalence of double
graded mixed cdga's
$$\DR^{int}(\hat{A}/D)\otimes_D k \simeq \DR^{int}(\hat{A}\otimes_D k/k).$$
We therefore have a natural morphism of double graded mixed cdga's
$$\DR^{int}(\hat{A}/D) \to \DR^{int}(\hat{A}/D)\otimes_D k \simeq \DR^{int}(\hat{A}\otimes_D k/k),$$
where $D$ is endowed with the trivial mixed structure as its second graded mixed structure. By applying
the first realization functor we get a morphism of graded mixed cdga's 
$$\CDR^I(\hat{A}/D) \to \CDR^{I}(\hat{A}/D)\widehat{\otimes}_{|D|_1} k \simeq \CDR^I(\hat{A}\otimes_D k/k),$$
where $\widehat{\otimes}$ is the completed tensor product of complete filtered complexes (see \S \ref{subsec:Thecompletedmonoidalstructure}). 
Now, Lemma \ref{ltII-2-2} precisely states that $|D|_1 \simeq k$, and thus the above morphism
turns out to be an equivalence of graded mixed cdga's
$$\CDR^I(\hat{A}/D) \simeq \CDR^I(\hat{A}\otimes_D k/k).$$
To finish the proof of the lemma we apply proposition \ref{pA-8-2} of Appendix \ref{secapp:InternalderiveddeRhamtheory} showing the existence
of a canonical equivalence of graded mixed cdga's
$D \simeq \CDR^I(\hat{A}/D)$.  \hfill $\Box$ \\

Using Lemma \ref{ltII-2-3}, we can define an $\s$-functor 
$$\phi : \egrcdga/(\hat{A} \otimes_{\hat{\DR}(\hat{A}/k)}k) \longrightarrow \hat{\DR}(\hat{A}/k)-\egrcdga$$
by sending an object $B' \to (\hat{A} \otimes_{\hat{\DR}(\hat{A}/k)}k)$ to the induced morphism in 
de Rham cohomology of the first kind
$$\CDR^{I}(\hat{A} \otimes_{\hat{\DR}(\hat{A}/k)}k/k) \longrightarrow 
\CDR^{I}(\hat{A}\otimes_{\hat{\DR}(\hat{A}/k)}k/B').$$
We claim that $\psi$ and $\phi$ are inverse to each others when restricted to the appropriate
full sub-$\s$-categories. 

We start by showing that $\phi\circ\psi$ is equivalent to the identity. Let $\hat{\DR}(\hat{A}/k) \to B$
be an object in $\hat{\DR}(\hat{A}/k)-\egrcdga$ corresponding to an almost perfect derived foliation on $(X,x)^\wedge$.
By definition, the underlying graded cdga of $B$ is of the form $Sym_{\hat{A}}(\LL[1])$, for 
an almost perfect $\hat{A}$-dg-module $\LL$. The object $\phi\circ\psi(B)$ is 
$\CDR^{I}(\hat{A}\otimes_D k / B\otimes_{D}k)$, 
where we keep using the notation $D:=\hat{\DR}(\hat{A}/k)$. We have to produce an 
equivalence of graded mixed cdga's $u : B \simeq \CDR^{I}(\hat{A}\otimes_D k/B\otimes_D k)$, 
functorial in $B$. For this,
we consider the commutative diagram with push-out squares in graded mixed cdga's
$$\xymatrix{
\hat{\DR}(\hat{A}/k) \ar[r] & \DR^{int}(\hat{A}/B) \ar[r] & \DR^{int}(\hat{A}\otimes_D k/B\otimes_D k) \\
\hat{\DR}^{int}(B/k) \ar[r] \ar[u] & B \ar[u] \ar[r] & k \ar[u]
}$$
where $\hat{\DR}^{int}(B/k)$ is by definition $\lim_{A\to A'}\DR^{int}(B\otimes_D \DR(A'/k)/k)$.
By passing to the realizations along the first mixed structure, and applying 
proposition \ref{pA-8-2}, we obtain a natural morphism
\begin{equation}\label{questaeq}B\simeq \CDR^I(\hat{A}/B) \to \CDR^I(\hat{A}\otimes_D k / B\otimes_D k).\end{equation}
It is then easy to check that this morphism is an equivalence, simply by considering 
the induced morphism on each graded pieces for the natural (complete) filtrations. 
Assume that $B$ is of the form $Sym_{\hat{A}}(E)$ 
as a graded cdga, then in weight $1$, the above morphism is the natural projection
$$E \longrightarrow E \otimes_{\hat{A}}|\hat{A}\otimes_D k|.$$
But, by considering the lemma \ref{ltII-2-3} for the weight $0$ pieces,  
we have $|\hat{A}\otimes_D k|\simeq \hat{A}$, and thus this morphism is
an equivalence of complexes. The same argument shows that (\ref{questaeq}) is  
an equivalence also in arbitrary weights $i>0$, being the canonical morphism
$$Sym^i_{\hat{A}}(E) \longrightarrow Sym^i_{\hat{A}}(E)  \otimes_{\hat{A}}|\hat{A}\otimes_D k|.$$

In the other direction, to see that $\psi\circ\phi \simeq id$, we simply observe that, 
for any diagram of graded mixed cdga's $k \to B' \to \hat{A}\otimes_D k$, with $B'$ of the form 
$Sym_k(L[1])$ as a graded cdga, the following square 
$$\xymatrix{D \ar[r] \ar[d] & \CDR^I(\hat{A}\otimes_D k/B') \ar[d] \\
k \ar[r] & B'}$$
is a push-out diagram of graded cdga's: this can be easily checked by considering the various
cotangent complexes.  \\

We have therefore finished the proof of the theorem for $X$ \emph{affine}, and it remains to extend
the statement to the  case of an arbitrary derived Artin stack $X$. For this we proceed on induction
on the geometricity of $X$, as we have already seen that 
the theorem is true when $X$ is affine, or more generally a disjoint 
union of affines. We fix a smooth atlas $U \to X$, where $U$ is a disjoint union of
of affine derived schemes locally of finite presentation. As $k$ is algebraically closed
we can lift the point $x : \Spec\, k \to X$ to a point $x_0 : \Spec\, k \to U$. We consider
the nerve $U_* \to X$ of $U\to X$, which comes equipped with a $k$-point $x_* : \Spec\, k \to U_*$
compatible with the point $x$. The theorem is true for all $U_i=U\times_X \times \dots \times_X U$
by our induction hypothesis.
In order to deduce the result for $X$, we consider the following commutative square of $\s$-categories
$$\xymatrix{\Fol^{ap}((X,x)^\wedge) \ar[r] \ar[d] & \ell_{X,x}/\dglie^{dap}_k \ar[d] \\
\lim_{i} \Fol^{ap}((U_i,x_i)^\wedge) \ar[r] & \lim_i \ell_{U_i,x_i}/\dglie^{dad}_k.}$$
The bottom horizontal $\s$-functor is an equivalence by our induction assumption. The vertical $\s$-functor
on the left is an equivalence, because $(U_*,u_*)^\wedge \to (X,x)^\wedge$ is a hypercovering
of derived stacks. Finally, the induced augmented simplicial object
$\ell_{U_*,u_*} \to \ell_{X,x}$ is a smooth hypercovering of dg-Lie algebras in the
sense of \cite[Prop. 1.5.8]{lupmf2}, and thus induces an equivalence of dg-Lie algebras
$$\mathrm{colim}_i \ell_{U_*,u_*} \simeq \ell_{X,x}.$$
This implies that also the vertical $\s$-functor on the right is an equivalence of $\s$-categories. 
\hfill $\Box$ \\

\subsection{Transversely smooth and rigid derived foliations}\label{subsec:Transverselysmoothandrigidderivedfoliations}

In this section we prove a stronger formal integrability statement
for quasi-smooth and rigid derived foliations, which is not only valid
formally at each point but globally on $X$. This result provides 
an important connection between derived foliations and classical notions
such as Gobdillon-Vey sequences that will be explored later on. \\

For this section we fix $X$ a smooth DM-stack over $k$. The results of these sections
can certainly be generalized to the more general setting, but we restrict to this case for
simplicity.

\begin{df}\label{dII-12}
A derived foliation $\F \in \Fol(X)$ is \emph{transversely smooth and rigid} if
its normal complex $\N_\F$ is a vector bundle.
\end{df}

Note that because $X$ is itself a smooth DM-stack
then $\LL_\F$ is quasi-isomorphic to the 2-term complex of vector bundles $\N_\F^* \to \Omega_X^1$, 
and that $\F$ is automatically quasi-smooth in the sense of the general Definition \ref{dII-8}. 
The term "rigid" in Definition
\ref{dII-12} refers here to the fact that $\Omega_X^1 \to H^0(\LL_\F)$ must be an epimorphism of sheaves, 
so that $\F$ does not have any "infinitesimal automorphisms" or no "fixed points". As an example, 
the foliation induced by a group scheme $G$ acting on $X$ (see Example \ref{ex:groupactions}) is transversely smooth
and rigid if and only if the stabilizers of the action are all discrete.

We start by constructing a canonical formal thickening $X \hookrightarrow \cJ_\F$
associated to a transversely smooth and rigid derived foliation $\F$. 
For this, we consider the push-out of graded mixed cdga's over $X$
$$\xymatrix{\DR(X/k) \ar[r] \ar[d] & \DR(\F/k) \ar[d] \\
\OO_X \ar[r] & B.}$$
The graded mixed cdga $B$ is $\OO_X$-linear via its canonical morphism 
$\OO_X \to B$, and as a graded $\OO_X$-linear cdga it is equivalent to 
$Sym_{\OO_X}(\N_\F^*[2])$. Using the equivalence between graded mixed complexes
and complete filtered complexes of Corollary \ref{clI-4}, we associate to $B$ a complete
filtered $\OO_X$-linear cdga $\B_\F := |B|^t$. The associated graded
of $B$ is canonical equivalent to $Sym_{\OO_X}(\N_\F^*)$. Moreover, as 
$B$ comes equipped with an augmentation $B \to \OO_X$, the filtered cdga $\B_\F$
is itself canonically augmented $\B_\F \to \OO_X$.

Let us denote the filtration on $\B_\F$ by
$$\xymatrix{
\dots  \ar[r] & F^i\B_\F \ar[r] & F^{i-1}\B_\F \ar[r] & \dots & F^0\B_\F=\B_\F,}$$
and note that $F^i\B_\F=F^{i+1}\B_\F=\B_\F$ for all $i\geq 0$.
We thus have
a pro-objects in the $\s$-category of quasi-coherent $\OO_X$-linear cdga's
$$\hat{\B}_\F = ``\lim_i"\,\B_\F/(F^i\B_\F).$$
This pro-object possesses a formal spectrum
$$\cJ_\F := \Spf(\hat{\B}_\F)=\mathrm{colim}_i \, \Spec\, (\B_\F/(F^i\B_\F)) \to X,$$
which is an affine formal scheme over $X$. The formal scheme $\cJ_\F$ comes equipped with a canonical 
closed embedding $X \to \cJ_\F$ induced from the augmentation $\B_\F \to \OO_X$, in such a way that 
$X$ is a retract of $\cJ_\F$ because of the $\OO_X$-linear structure
on $\hat{\B}_\F$. Note also that $\cJ_\F \to X$ is locally on $X$ of the form
$X \times \hat{\mathbb{A}}^d$, where $d$ is the codimension of $\F$ (i.e. the rank of
the normal bundle $\N_\F$). Indeed, as $\cJ_\F$ is a relative formal spectrum, its
construction is stable by base change, so we may assume that $X$ is affine of the form 
$\Spec\, A$. In this case, the filtered $A$-cdga $\B_\F$ is such that $Gr^*(\B_\F)\simeq
Sym_A(\N_\F^*)$, where $\N_\F*$ is of weight $1$. Since 
$\N_\F^*$ is a vector bundle and $X$ is affine, there are no obstructions to lift 
the canonical morphism of $A$-dg-modules
$$\N_\F^* \to \B_\F/F^2 \simeq A \oplus \N_\F^*$$
to a morphism of $A$-dg-modules $\N_\F^* \to \B_\F$. This extends by multiplicativity to
a morphism of completed $A$-algebras
$\hat{Sym}_A(\N_\F^*) \to \B_\F$ which is clearly an equivalence. This shows that $\cJ_\F$ is thus
equivalent, as a formal affine scheme over $X$, to $\hat{\VV(\N_\F)} \to X$, the
formal completion at the zero section of the total space of the vector bundle $\N_\F$ over $X$.
Note that however this is not a functorial description, and therefore
the two formal schemes $\hat{\VV(\N_\F)}$ and $\cJ_\F$ are only locally (not globally) equivalent, for 
the étale topology on $X$, and the difference between these two objects is understood globally 
via the filtration on $\hat{\B}_\F$ interpolating between them.

So we have constructed, for any transversely smooth and rigid derived foliation $\F$ on $X$, 
a formal scheme $\cJ_\F \to X$ together with a section $j : X \to \cJ_\F$, which is
moreover obviously functorial in the pair $(X, \F)$. 

\begin{df}\label{dII-13}
For a transversely smooth and rigid $\F \in \Fol(X)$, the relative formal scheme $\cJ_\F \to X$
constructed above is called the \emph{(transversal) jet space of $\F$}.
\end{df}

The importance of the transversal jet space lies in the following integrability result.

\begin{prop}\label{pII-8}
Let $X$ be a smooth DM-stack, and $\F \in \Fol(X)$ transversely smooth and rigid with transverse 
jet space $\cJ_\F \to X$. There exists a canonical smooth, 
transversely smooth and rigid derived foliation
$\F' \in \Fol(\cJ_\F)$ such that $j^*(\F')$ is canonically equivalent to $\F$ in $\Fol(X)$.
\end{prop}

\textit{Proof.} It is enough to prove the Proposition on $X=\Spec\, A$ in such a way that 
the construction of $\F'$, and the equivalence $j^*(\F')\simeq \F$, are both 
stable by base change along any \'etale morphism $Y \to X$ between smooth affine schemes. We are thus
essentially reduced to the case where $X$ is smooth and affine. 

Let $\DR(\F/k)$ be the graded mixed $\DR(A/k)$-cdga corresponding to 
$\F$. We have an exact triangle of complexes $\N_\F^* \to \Omega_{A/k}^1 \to \LL_\F$. Equivalently,
$\LL_\F$ identifies canonically with the two terms complex $\N_\F^* \to \Omega_{A/k}^1$
where $\Omega_{A/k}^1$ sits in degree $0$. The formal scheme $\cJ_\F$ is here the formal
spectrum of the filtered $A$-algebra $\B_\F=|B|^t$, where $B$ is the $A$-linear graded mixed
cdga defined by $B:=\DR(\F/k)\otimes_{\DR(A/k)}A$. Note that $\B_\F$ is non-canonically
isomorphic to $\hat{Sym}_{A}(\N_\F)$.

We let $\hat{\DR}(\cJ_\F/k)$ be defined by 
$$\hat{\DR}(\cJ_\F/k):=\lim_{n} \DR((\B_\F/F^n)/k)$$
be the completed graded mixed cdga associated to the formal scheme $\cJ_\F$. 
We have the
following statement concerning perfect derived foliations on $\cJ_\F$, whose proof is 
very much similar (actually simpler) than proposition \ref{pII-9}, and is left to the reader
as an exercise.

\begin{lem}\label{lpII-8}
There is a natural equivalence between the following two $\s$-categories:
\begin{itemize}
    \item $\Fol^{ap}(\cJ_\F)$, the $\s$-category of perfect derived foliations on $\cJ_\F$
    \item the $\s$-category of graded mixed $\hat{\DR}(\cJ_\F/k)$-cdga $D$, of the
    form $Sym_{\B_\F}(V[1])$ as graded cdga where $V$ is aan almost perfect $\B_\F$-dg-module.
\end{itemize}
\end{lem}

To achieve the proof proposition \ref{pII-8}, it is enough to 
construct a smooth derived foliation $\F' \in \Fol(\cJ_\F)$ together with an equivalence
$j^*(\F') \simeq \F$, which is functorial under \'etale base change on $X$. Indeed,  
as $\F$ is rigid, such an $\F'$ must automatically be rigid as well: the conormal
complex $\N^*_{\F'}$ being perfect on $\cJ_\F$ is such that 
$j^*(\N^*_{\F'}) \simeq \N_\F^*$ is vector bundle, and thus $\N_{\F'}^*$
must already be a vector bundle on the formal thickening $\cJ_\F$. As a result, 
such a smooth $\F'$ with $j^*(\F') \simeq \F$ must be transversally smooth, and thus
is automatically rigid because the short exact sequence of bundles 
$0 \to \N^*_{\F'} \to \Omega_{\cJ_\F}^1 \to \LL_{\F'} \to 0$
locally splits.

In order to construct $\F'$, we use our notion of de Rham cohomology of the first kind 
explained in \S \ref{secapp:InternalderiveddeRhamtheory}. We recall the push-out square of graded mixed cdga's
$$\xymatrix{
\DR(X/k) \ar[r] \ar[d] & \DR(\F/k) \ar[d] \\
\OO_X \ar[r] & B,}$$
and we consider $\CDR^I(B/\DR(\F/k))$, the de Rham cohomology of the first kind.
A direct computation of
cotangent complexes shows that, as a graded algebra, $\CDR^I(B/\DR(\F/k))$ is of the form 
$Sym_{\B_{\F}}(\Omega_X^1 \otimes_A \B_\F[1])$. Therefore, by the equivalence of $\s$-categories
describing $\Fol(\cJ_\F)$, $\CDR^I(B/\DR(\F/k))$ defines a perfect derived foliation
$\F' \in \Fol(\cJ_\F/k)$. Moreover, the cotangent complex of $\F'$ is $p^*(\Omega_{X/k}^1)$, 
where $p : \cJ_\F \to X$ is the natural projection, and thus $\F'$ is smooth.

We consider now the canonical morphism $\CDR^I(B/\DR(\F/k)) \to \CDR^I(A/\DR(\F/k))$, via
the augmentation $B \to A$. By 
Proposition \ref{pA-8-2}, we know that the right hand side is naturally equivalent
to $\DR(\F/k)$, and we thus have a morphism of graded mixed cdga's
$\CDR^I(B/\DR(\F/k)) \to \DR(\F/k)$, making the following diagram commutative
$$\xymatrix{
\CDR^I(B/\DR(\F/k)) \ar[r] & \DR(\F/k) \\
\hat{\DR}(\cJ_\F/k) \ar[r] \ar[u] & \DR(A/k). \ar[u]
}$$
This produces a morphism
$$\CDR^I(B/\DR(\F/k)) \otimes_{\hat{\DR}(\cJ_\F/k)} \DR(A/k) \longrightarrow \DR(\F/k)$$
which is easily seen to be an equivalence by computing the cotangent complexes. By construction, this
is the required equivalence $j^*(\F') \simeq \F$. \hfill $\Box$

\begin{rmk}\label{rII-1}
\begin{enumerate}
    \item \emph{An important point of proposition \ref{pII-8} is that the
association $\F \mapsto (\cJ_\F,\F')$ is canonical. The foliation $\F'$ will be called
the \emph{canonical spread} of the derived foliation $\F$. We will see below that the 
pair $(\cJ_\F,\F')$ is closely related to the classical notion of Godbillon-Vey sequences (see 
for instance \cite{MR1955577} and \cite[$(3.4)$]{mal2} for classical references).}
\item \emph{In the special case where, in the statement of the proposition \ref{pII-8} $\F$ is \emph{integrable}, 
the objects $\cJ_\F$ and $\F'$ can be described explicitly as follows. Let $f : X \to Y$
be a moprhism of smooth separated schemes such that $f^*(0_Y)\simeq \F$. Then 
$\cJ_\F$ is identified with the formal completion of $X \times_k Y$ along the graph of $f$, 
and $\F'$ is $q^*(0_Y)$ where $q : X \times_k Y \to Y$ is the second projection. 
}
\item \emph{Being at the same time smooth and transverally smooth, the derived foliation
$\F'$ is a genuine foliation on $\cJ_\F$, given by an integrable sub-bundle 
$\TT_{\F'} \subset \T_{\cJ_\F}$ of the tangent sheaf.}
\end{enumerate}
\end{rmk}

One important corollary of proposition \ref{pII-8} is the following formal integrability result, 
that could also be obtained from the more general theorems \ref{tII-2} and \ref{ctII-2}. 

\begin{cor}\label{cpII-8}
Let $X$ be a smooth DM-stack and $x : \Spec\, k \to X$ a global point. Let $\F \in \Fol(X/k)$ be a transversely
smooth and rigid derived foliation on $X$. There exists a morphism of formal schemes
$f : (X,x)^\wedge \to \hat{\mathbb{A}}^d$, such that $\hat{F}_x \simeq f^*(0)$ in $\Fol((X,x)^\wedge)$. 
In other words, $\F$ is formally integrable at each point on $X$.
\end{cor}

\textit{Proof of the Corollary.} The statement is local at $x$, so we can assume that $X$ is affine, 
and moreover that the normal bundle $\N_\F\simeq \OO_X^d$ is trivial, as well 
as $\Omega_{X/k}^1$.
We use the canonical spread $j : X \hookrightarrow X\times \hat{\mathbb{A}}^d$,
with $\F' \in \Fol(X\times \hat{\mathbb{A}}^d)$, $j^*(\F')\simeq \F$. We formally complete at $x$ to 
get 
$$\hat{j} : (X,x)^{\wedge} \longrightarrow (X\times \mathbb{A}^d,(x,0))^{\wedge}.$$
The formal completion at $(x,0)$ defines a foliation 
$\hat{\F'} \in \Fol((X\times \mathbb{A}^d,(x,0))^{\wedge})$, which is here
given by a family of formal vector fields $\nu_1,\dots,\nu_n$ on the formal scheme 
$(X\times \mathbb{A}^d,(x,0))^{\wedge}\simeq \hat{\mathbb{A}}^{n+d}$, where $n=Dim_k X$,
satisfying the classical 
Frobenius condition: each $[\nu_i,\nu_j]$ can be written as a linear combination
of $\nu_k$'s. The formal Frobenius lemma (see e.g. \cite{MR2588634}) then tells us 
that there is a morphism of formal schemes
$p : (X\times \mathbb{A}^d,(x,0))^{\wedge} \to \hat{\mathbb{A}}^d$, such that 
$\hat{\F'}\simeq p^*(0)$. We can thus take $f=p\circ \hat{j}$.
\hfill $\Box$ \\

\section{Existence of derived enhancements}\label{sec:Existenceofderivedenhancements}

In this section we study the relations between the notion of derived foliation 
on a derived scheme $X$, and the classical notion of differential ideal on the truncation
$\tau_0(X)$. We show that any derived foliation underlies a classical differential ideal. However, 
not every differential ideal arises this way, and in the case where $X$ is a smooth variety, we provide
a criterion for the existence of a derived foliation enhancing a given classical 
differential ideal, based on the classical notion of Gobdillon-Vey sequences (see 
for instance \cite{MR1955577,mal2}). 

\subsection{Derived foliations and differential ideals}

Let $X$ be a derived scheme locally almost 
of finite presentation over $k$, and $\F \in \Fol(X/k)$ be a derived
foliation on $X$. We assume that $\F$ is connective in the sense of 
Definition \ref{dII-1}, i.e. $H^i(\LL_{\F/k})\simeq 0$ for $i>0$.
We consider the corresponding sheaf of graded mixed $\DR_{X/k}$-cdga $\DR_{\F}$.
We have the canonical morphism of quasi-coherent complexes
$$a : \LL_{X/k} \to \LL_{\F/k}$$
from which we form the kernel of the induced morphism on $H^0$
$$K_{\F}:=Ker(H^0(\LL_{X/k}) \to H^0(\LL_{\F/k})).$$
The sheaf $K_\F$ is a coherent subsheaf of $H^0(\LL_{X/k})\simeq \Omega_{X_0/k}^1$, where
$X_0=\tau_0(X)$ is the underived scheme obtained from $X$ by truncation (see Appendix \ref{secapp:Derivedschemesandstacks}). 

\begin{lem}\label{lII-2}
The subsheaf $K_\F \subset \Omega_{X_0/k}^1$ is a coherent differential ideal: 
$$dR(D_\F) \subset Im(D_\F\otimes \Omega_{X_0/k}^1 \to \Omega^2_{X_0/k}).$$
\end{lem}

\textit{Proof.} Since the graded mixed
structure on $\DR_\F\simeq Sym_{\OO_X}(\LL_{\F/k}[1])$ is compatible with the 
de Rham differential on $Sym_{\OO_X}(\LL_{X/k})$, we have a commutative diagram of sheaves
of $k$-modules on $X$
$$\xymatrix{
H^0(\LL_{X_0}/k) \ar[d]_-{a} \ar[r]^-{dR} & H^0(\wedge_{\OO_X}^2\LL_{X_0}/k) \ar[d]^-{\wedge^2 a} \\
H^0(\LL_{\F/k}) \ar[r]_-{\epsilon} & H^0(\wedge^2_{\OO_X}\LL_{\F/k}).
}$$
Because $\LL_{X_0/k}$ and $\LL_{\F/k}$ are both connective, the above diagram is isomorphic to 
$$\xymatrix{
\Omega_{X_0/k}^1 \ar[d]_-{a} \ar[r]^-{dR} & \Omega_{X_0/k}^2 \ar[d]^-{\wedge^2 a} \\
H^0(\LL_{\F/k}) \ar[r]_-{\epsilon} & \wedge^2_{\OO_X}H^0(\LL_{\F/k}).
}$$
The lemma thus follows by considering the induced morphism on the kernels of the vertical
morphisms. \hfill $\Box$ \\

\begin{df}\label{dII-14}
\begin{enumerate}
    \item With the notations and conditions above, the differential ideal $K_\F \subset \Omega_{X_0/k}^1$
is called the \emph{classical truncation} of the derived foliation $\F \in \Fol(X/k)$.

    \item For a coherent differential ideal $K\subset \Omega_{X_0/k}^1$, a \emph{derived enhancement
of $D$} is a connective derived foliations $\F \in \Fol(X/k)$ such that $K=K_\F$.
\end{enumerate}

\end{df}

Note that, if $\F \to \F'$ is a morphism of connective derived foliations, then we have a natural
inclusion of differential ideals $K_\F' \subset K_\F \subset \Omega^1_{X_0/k}$. Therefore, 
if we denote by $IDiff(X_0)$ the category of coherent differential ideals on $X_0$ (with morphisms
given by inclusions), the classical truncation is an $\s$-functor
$$K_{-} : \Fol^c(X/k)^{op} \longrightarrow IDiff(X_0),$$
from connective derived foliations on $X$ to differential ideals on $X_0$.
As we will now explain, there are no reasons to expect that the above $\s$-functor is essentially surjective, even in the case
where $X$ is a nice smooth scheme over $k$. Indeed, if $K\subset \Omega_{X/k}^1$ is 
a coherent differential ideal, we can certainly form a graded mixed cdga by considering
$Sym^{st}_{\OO_X}(E[1])$, where $E$ is the coherent sheaf $\Omega_{X/k}^1/K$, and where
the mixed structure is defined by the unique possible extension of the de Rham differential
along the quotient map $\DR(X/k) \to Sym^{st}_{\OO_X}(E[1])$. 

It is important to note that here
$Sym^{st}$, as opposed to notation $Sym$, denotes the strict or naive symmetric algebra, 
defined in the tensor abelian category of coherent sheaves on $X$. This makes a difference as 
$E$ is in general not a flat $\OO_X$-module, and thus the canonical morphism 
$Sym_{\OO_X}(E[1]) \to Sym^{st}_{\OO_X}(E[1])$, from the derived $Sym$-algeba to the naive $Sym$-algebra is not 
a quasi-isomorphism. In particular, the graded mixed cdga $Sym^{st}_{\OO_X}(E[1])$ does not 
satisfy the condition appearing in Definition \ref{dII-1}, and thus does not define a 
derived foliation on $X$. In order to promote $Sym^{st}_{\OO_X}(E[1])$ into a derived foliation,
we would have to lift the de Rham differential along the natural morphism $Sym_{\OO_X}(E[1]) \to 
Sym^{st}_{\OO_X}(E[1])$. We will see in the next section that there are obstructions to the
existence of such a lift in general, and that these obstructions can be expressed explicitly in simple terms. 

The main reason for expecting that the classical truncation $K_{-}$ is not 
essentially surjective comes from integrability results for derived foliations, such 
as Theorem \ref{ctII-2} or proposition \ref{pII-8}, which indicate that derived foliations
are always formally integrable at each point. There are however example of differential
ideals, or even analytic germs of differential ideals, which are not 
formally integrable (see remark \ref{rdII-14} $(1)$ below). Of course,
the notion of derived enhancement of Definition
\ref{dII-14} is extremely general as we do not restrict the type of derived foliation $\F$, 
so  theorem \ref{ctII-2} or proposition \ref{pII-8} are not strictly speaking counterexamples
to the existence of derived enhancements, and possibly the questions of their existence
is complicated and not well posed. However, when restricting the type of derived
foliations allowed, the question of existence of derived enhancements can be 
treated in great detail. This is what we will do in the next section for \emph{transversely smooth
and rigid} derived foliations. \\

\begin{rmk}\label{rdII-14}
\begin{itemize}
    \item \emph{Examples of not formally integrable differential ideals already exists
    in dimension $2$. Indeed, for a germ $\omega$ of holomorphic 1-differential with an isolated zero 
    at $0 \in \CC^2$ and satisfying $d\omega \wedge \omega = 0$, formal integrability and holomorphic 
    integrability are equivalent conditions by the main result of \cite{mal1}. Moreover, it is well 
    known that
    holomorphic integrability is restricted by topological obstructions, see for instance 
    \cite[Theorem B]{zbMATH03717678}. This provides examples of (analytic germs of) 
    differential ideal which do not admit derived enhancements. }
    \item \emph{It is worth mentioning that the truncation functor $K_{-}$ of Definition \ref{dII-14}
can also be understood in terms of a natural $t$-structure on graded mixed complexes. 
Indeed, we can introduce a $t$-structure on graded mixed complexes by 
declaring that $E \in \egrdg_k$ is \emph{connective} if for each $n\in \ZZ$, the part $E^{(n)}$
of weight $n$ is $n$-connective: $H^i(E^{(n)})=0$ for $i<-n$. This $t$-structure is compatible
with the symmetric monoidal structure, and thus t-truncations functors extends to 
commutative algebra objects. A connective derived foliation $\F \in \Fol^c(X/k)$ is such that 
$\DR(\F)$ is a connective graded mixed cdga, and its $0$-truncation for the above mentioned
$t$-structure is precisely $Sym^{st}_{\OO_X}(E[1])$ where $E=\Omega_{X_0/k}^1/K_\F$.}
\end{itemize}

\end{rmk}

\subsection{Derived enhancements and Godbillon-Vey sequences}

Let $X=\Spec\, A$ be a smooth affine scheme over $k$, and $K\subset \Omega_{X/k}^1$ be a
differential ideal and we assume that $K$ is a free an $A$-module (but we do not assume that 
it is a sub-bunld eof $\Omega_{X/k}^1$). 
We chose a set of generators $(\omega_1,\dots,\omega_d)$ of $K$ as an $A$-module, 
which is a basis. 
We recall from \cite[(3.4)]{mal2} that  a \emph{Godbillon-Vey (\emph{or simply} GV) sequence of forms
adapted to $(\omega_1,\dots,\omega_d)$} consists of a sequence of elements 
$\omega_{i,\alpha} \in \Omega_{X/k}^1,$
where $\alpha \in \NN^p$, satisfying the equations

\begin{equation}\label{eq:GV}
dR(\omega_{i,\alpha}) =  \sum_{j,\beta}\binom{\alpha}{\beta} \omega_{j,\beta} \wedge
\omega_{i,\alpha-\beta+e_j} \qquad \omega_{i,0}=\omega_i
\end{equation}
where $e_j=(0,\dots,0,1,0,\dots,0) \in \NN^p$ is the $j$-th element of the standard basis.
The GV sequences can also be interpreted as follows.

\begin{lem}\label{lII-3}
With the notations above, there is a bijection between the following sets.

\begin{itemize}
    \item The set of GV sequences 
adapted to $(\omega_1,\dots,\omega_d)$.

\item Families of differential forms $(\widetilde{\omega_1},\dots,\widetilde{\omega_d})$ on the
formal scheme $X \times \widehat{\mathbb{A}}^{d}$ with the following conditions
satisfied.
\begin{enumerate}
    \item $j^*(\widetilde{\omega_i})=\omega_i$, where $j : X \hookrightarrow X \times \widehat{\mathbb{A}}^{d}$
    is the embedding at $0$.
    \item For each $i$, the form $\widetilde{\omega_i}$ can be written as 
    $$\widetilde{\omega_i}=dR(t_i)+\sum_k a_k.p^*(\lambda_k)$$
    with $\lambda_k \in \Omega_{X/k}^1$, $p : X \times \widehat{\mathbb{A}}^{d} \to X$ being the projection 
    onto the first factor, $(t_1,\dots, t_p)$ the
    standard coordinates on $\widehat{\mathbb{A}}^{d}$, and $a_k \in A[[t_1,\dots,t_p]]$ formal functions on  $X \times \widehat{\mathbb{A}}^{d}$.
    \item The forms $(\widetilde{\omega_1},\dots,\widetilde{\omega_d})$ define a differential 
    ideal in $\widehat{\Omega}_{X \times \widehat{\mathbb{A}}^{d}}^1$.
\end{enumerate}
\end{itemize}
\end{lem}

\textit{Proof.} This is almost straightforward: to a GV sequence $\omega_{i,\alpha}$ we
associate the forms
$$\widetilde{\omega_i}:=dR(t_i)+\sum_{i,\alpha} p^*(\omega_{i,\alpha}).\frac{t^{\alpha}}{\alpha !}$$
We refer to \cite[Prop. 3.2]{mal2} for more details. \hfill $\Box$ \\

It is important to note that because of condition $(2)$, the 
forms $(\widetilde{\omega_1},\dots,\widetilde{\omega_d})$ are automatically linearly independent
everywhere, and thus the differential ideal they generate is a free subbundle of rank $p$ in 
$\widehat{\Omega}_{X \times \widehat{\mathbb{A}}^{d}}^1$. Therefore, they define a smooth, 
transversely smooth and rigid foliation (i.e. 
a classical foliation in the usual sense) on the formal scheme $X \times \widehat{\mathbb{A}}^{d}$
exactly as in proposition \ref{pII-8}. In fact, proposition \ref{pII-8}
can be used in order to prove the following result, which expresses that 
GV sequences and transversely smooth and rigid derived derived enhancements are essentially the same
thing.

\begin{prop}\label{pII-10}
Let $X$ be a smooth scheme over $k$ and $K \subset \Omega_{X/k}^1$ be a coherent differential ideal
such that $K$ is a vector bundle.
The following two conditions are equivalent:
\begin{enumerate}
    \item Locally for the Zariski topology on $X$, there exists a Godbillon-Vey sequence (adapated to 
    some choice of generators $(\omega_1,\dots,\omega_p)$ of $K$).
    
    \item Locally for the Zariski topology on $X$ there exists a transversely smooth and
    rigid derived foliation $\F \in \Fol(X/k)$ such that $K=K_\F$.
    
\end{enumerate}
\end{prop}

\textit{Proof.} We have already seen that $(1) \Rightarrow (2)$. Indeed, 
a GV sequence for $K \in IDiff(X)$ defines a smooth, transversely smooth and rigid derived foliation 
$\F'$ on $X \times \widehat{\mathbb{A}}^{d}$. We let $\F:=j^*(\F') \in \Fol(X/k)$. 
As pull-backs of derived foliations preserve conormal complexes, we see that 
$\F$ is transversely smooth and rigid. Moreover, by construction, the differential 
ideal $K_\F$ is generated by $(\omega_1,\dots,\omega_d)$ and thus coincides with $K$. 

Conversely, suppose that $K$ is of the form $K_\F$. We apply proposition \ref{pII-9}
to $\F$. As the statement is local, we can even assume that the normal bundle 
of $\F$ is trivialized, and the proposition gives us $\F' \in \Fol(X\times \widehat{\mathbb{A}}^d/k)$
such that $j^*(\F')=\F$. Moreover, by observing the proof of proposition \ref{pII-9}, 
it is easy to see that $\F'$ is given by a differential ideal 
$(\widetilde{\omega_1},\dots,\widetilde{\omega_d}) \subset \widehat{\Omega}_{X \times 
\widehat{\mathbb{A}}^{d}}^1$, where each $\widetilde{\omega_i}$ is of the form
$$\widetilde{\omega_i}=\sum_j f_{i,j}.dR(t_j)+\sum_k a_k.p^*(\lambda_k)$$
with $f_{i,j} \in A[[t]]=A[[t_1,\dots,t_p]]$ a formal function such that 
the the $p \times p$ matrix $M:=(f_{i,j})$ is invertible. Indeed, 
it is easy to see by construction that the cotangent complex
of $\F'$ is a vector bundle of rank $dim 
X$, and that the anchor map
$$\widehat{\Omega}_{X \times \widehat{\mathbb{A}}^{d}}^1 \longrightarrow \LL_{\F'},$$
precomposed with $p^*(\Omega_{X/k}^1) \to \widehat{\Omega}_{X \times \widehat{\mathbb{A}}^{d}}^1$
induces an isomorphism $p^*(\Omega_{X/k}^1) \simeq \LL_{\F'}$. As 
$$\widehat{\Omega}_{X \times \widehat{\mathbb{A}}^{d}}^1 \simeq p^*(\Omega_{X/k}^1) \bigoplus \oplus_{i}
A[[t]].dR(t_i),$$
we see that the kernel of the anchor map, which is the subsheaf generated by 
the $\widetilde{\omega_i}$'s, becomes isomorphic to $\oplus_i A[[t]].dR(t_i)$
via the projection $\widehat{\Omega}_{X \times \widehat{\mathbb{A}}^{d}}^1 \to \oplus_i A[[t]].dR(t_i)$.
We can thus 
change the basis $(\widetilde{\omega_1},\dots,\widetilde{\omega_d})$
by $M^{-1}(\widetilde{\omega_1},\dots,\widetilde{\omega_d})$, which defines 
the same foliation $\F'$, and thus simply assume that $f_{i,j}=\delta_{i,j}$.
By Lemma \ref{lII-3}, the family $(\widetilde{\omega_1},\dots,\widetilde{\omega_d})$
defines a GV sequence for the differential ideal $K$. \hfill $\Box$ \\

The proposition \ref{pII-10} has an important consequence concerning existence of derived
enhancement under appropriate codimension conditions. The question of uniqueness of a derived 
enhancement is more
complicated and will be studied later using analytic methods (see \S \ref{section_anintderenh}).

\begin{cor}\label{cpII-10}
Let $X=\Spec\, A$ be a smooth affine scheme over $k$ and $K \subset \Omega_{X/k}^1$ be a 
differential ideal. We assume that the coherent sheaf $\Omega_{X/k}^1/K$ is a vector bundle
on a dense open $U \subset X$, whose closed complement $Z$ is of codimension at least $3$
in $X$. 
Then, locally for the Zariski topology on $X$, $K$ is of the form $K=K_\F$ for 
a transversely smooth and rigid derived foliation $\F \in \Fol(X/k)$.
\end{cor}

\textit{Proof of the Corollary.} We first notice that $K$ must be a vector bundle on $U$, and thus
on $X$ by Hartogs.
Since the statement is local for the Zariski topology on $X$, we can assume the existence of 
an isomorphism of coherent sheaves $\OO_X^d \simeq K$. We consider
the induced injective morphism $\OO_X^d \hookrightarrow \Omega_{X/k}^1$. This defines
a complex of vector bundles concentrated in cohomological degress $[-1,0]$ that we
denote by $\LL$. We then use the following well-known vanishing lemma.

\begin{lem}\label{lcpII-10}
Under the conditions above, for all $n\geq 0$ we have
$H^i(\wedge^{n+2}_{\OO_X}(\LL))\simeq 0$
for all $i \leq -n$.
\end{lem}

\textit{Proof of the lemma.} The complex $\wedge^{n+2}_{\OO_X}(\LL)$ is a perfect complex which 
is a direct factor of $\LL^{\otimes n+2}$, which is itself perfect and represented by the
genuine complex of vector bundles $(\OO_X^d \to \Omega_{X/k}^1)^{\otimes n+2}$. This is
a complex of vector bundles sitting in cohomological degrees $[-n-2,0]$. We form
the dual perfect complex $\LL^\vee$ which is a bounded complex of vector bundles
concentrated in degrees $[0,n+2]$. Moreover, on the open $U$ the complex $\LL$ is equivalent to 
a vector bundle, and this all the cohomology sheaves $H_i := \underline{H}^i(\LL^\vee)$ vanishes
on $U$ as soon as $i>0$. In particular, by Grothendieck duality 
the ext-sheaves $\underline{Ext}^i_{\OO_X}(H_i,\OO_X)$
are zero for $i=0$, $i=1$ and $i=2$. This clearly implies, by postnikov induction on $\LL^\vee$,
that $\LL \simeq \mathbb{R}\underline{Hom}_{\OO_X}(\LL^\vee,\OO_X)$ has cohomology sheaves
concentrated in cohomological degrees $[-n,0]$ as wanted.
\hfill $\Box$ \\

Coming back to the proof of the Corollary, by proposition \ref{pII-10} it is
enough to show the existence of a GV sequence $\omega_{i,\alpha}$ adapted to the basis
$(\omega_1,\dots,\omega_d)$ of $K$ coming from our choice of isomorphism $\OO_X^d \simeq K$. 
We construct the $\omega_{i,\alpha}$ by induction on the multi-index $\alpha$. Suppose
that we have already constructed $\omega_{i,\alpha}$ for all $|\alpha|\leq  n$, such that 
equations (\ref{eq:GV}) are satisfied for all $|\alpha| < n$. For all $i$, and a multi-index
$\alpha$ with $|\alpha|=n$, consider the equation (\ref{eq:GV}) in the following form
\begin{equation}\label{eq_GVused}\sum_{j}\omega_{i,0}\wedge \omega_{i,\alpha+e_j} =  dR(\omega_{i,\alpha}) - \sum_{j,\beta\geq 1}\binom{\alpha}{\beta} \omega_{j,\beta} \wedge
\omega_{i,\alpha-\beta+e_j},\end{equation}
where $\omega_{i,\alpha+e_j}$ are to be determined. The lemma \ref{lcpII-10} tells us that the complex $\wedge^{n+2}_{\OO_X}(\LL)$ is exact in degree $[n-2,n]$, 
and thus, in particular, that the sequence
$$Sym^{n+1}(\OO_X^d) \otimes \Omega_{X/k}^1 \to Sym^{n}(\OO_X^d) \otimes \Omega_{X/k}^2
\to Sym^{n-1}(\OO_X^d) \otimes \Omega_{X/k}^3$$
is exact. The right hand side of (\ref{eq_GVused}) can be considered as cocycle in $Sym^{n}(\OO_X^d) \otimes \Omega_{X/k}^2$, 
because equations (\ref{eq:GV}) are satisfied for $|\alpha|<n$. The condition that this
cocycle is in fact a coboundary precisely says that there exists forms $\omega_{i,\alpha+e_j}$
such that the above equation (\ref{eq_GVused}) is satisfied. This finishes the proof of the Corollary. 
\hfill $\Box$ \\

We finish by emphasising a important consequence of lemma \ref{lcpII-10} that will be used
later on. Let $\F$ be a derived foliation which is rigid and quasi-smooth, and 
$K_\F$ the corresponding ideal. We can compare the derived de Rham 
cohomology of $\F$ with the \emph{naive de Rham cohomology} of $K_\F$ by means of the natural 
projection
$$p : \DR_\F \to \DR^{st}_\F=Sym_{\OO_X}^{st}(\Omega_{\F}^1[1]).$$
Here, $\Omega_\F^1 = H^0(\LL_\F)$ is the coherent sheaf of $1$-forms along $\F$, which 
is quasi-isomorphic to $\LL_\F$ via the natural truncation map $\LL_\F \to H^0(\LL_\F)\simeq 
\Omega_X^1/\cN_\F^*$. 
The graded mixed cdga $Sym_{\OO_X}^{st}(\Omega_{\F}^1[1])$ is the naive or \emph{strict} 
symmetric product algebra over the coherent sheaf $\Omega_\F[1]$, where the mixed
structure is induced by the de Rham differential $\OO_X \to \Omega_\F^1$. The morphism
$p$ is also the natural projection of the connective object $\DR_\F$ to its
$0$-th truncation for the Beilinson's t-structure on graded mixed objects (see \ref{dBeilinsontstructure}).
By passing to the non-Tate realization we have an induced morphism of the corresponding sheaves of 
de Rham cohomologies
$$|\DR_\F| \simeq \CDR^*(\F)^u \to |\DR^{st}_\F|=:C_{DR,naive}^*(\F).$$
Passing on hyper cohomology groups over $X$ we get an induced morphism
$$\widehat{H}^i_{DR}(\F) \to H^i_{DR,naive}(\F)$$
relating the Hodge completed derived de Rham cohomology with naive de Rham cohomology of $\F$.

\begin{cor}\label{clcpII-10}
With the notations above, assume that $\LL_\F$ is a vector bundle outside of a closed
subset $Z \in X$ of codimension $d$. Then, the induced morphism 
$$\widehat{H}^i_{DR}(\F) \to H^i_{DR,naive}(\F)$$
is an isomorphism for all $i<d-1$.
\end{cor}

\textit{Proof of the corollary.} This is a direct application of the lemma
\ref{lcpII-10} and the compatibility of the morphism $|\DR_\F| \simeq \CDR^*(\F)^u \to |\DR^{st}_\F|=:C_{DR,naive}^*(\F)$ with the Hodge filtrations on both sides.
\hfill $\Box$ \\

\section{Direct images}\label{sec:Directimages}

Let $f : X \to Y$ be a morphism of derived stacks. The pull-back $\s$-functor
$f^* : \dSt/Y \to \dSt/X$, sending $(Z \to Y)$ to $(X\times_Y Z \to X)$, commutes with 
arbitrary colimits. In particular, it admits a right adjoint
$$f_* : \dSt/X \to \dSt/Y,$$
called the \emph{direct image}, or \emph{Weil restriction} along $f$. 

For any $Z \to X$, the adjunction morphism is a commutative diagram of derived stacks
$$\xymatrix{
X\times_Y f_*(Z) \ar[r] \ar[dr] & Z \ar[d] \\
 & X.}$$
The adjunction counit $X\times_Y f_*(Z) \to Z$ will be denoted by $ev$ and will be called the
\emph{evaluation morphism}. The following theorem shows that, under certain properness conditions, 
derived foliations can be push-forward. This provides a very powerful tool to construct new derived
foliations from existing ones.

\begin{thm}\label{tII-3}
Assume that $f : X \to Y$ is a morphism of derived Artin stacks (locally of finite presentation) 
which is quasi-smooth and representable by a proper derived algebraic space. 
Let $Z \to X$ be a derived Artin 
stack locally of finite presentation over $X$.
Then, 
there exists an $\s$-functor
$$f_* : \Fol^{p}(Z/X) \to \Fol^{p}(f_*(Z)/Y)$$
from perfect derived foliations on $Z$ relative to $X$ to 
perfect derived foliations on $f_*(Z)$ relative to $Y$ (see Example \ref{ex:relativefoliations}). Moreover, the cotangent 
complex of $f_*(\F)$ is given by 
$$\LL_{f_*(\F)}\simeq q_+(ev^*(\LL_\F))$$
where $q : X\times_Y f_*(Z) \to f_*(Z)$ is the second projection and 
$q_+ : \Parf(X\times_Y f_*(Z)) \to \Parf(f_*(Z))$ is defined on perfect complexes by 
$q_+(E):=(q_*(E^\vee))^\vee$.
\end{thm}

The above theorem is proven in \cite{alf} and its proof will not be
covered in this book. In the next two subsections we will simply explain how the $\s$-functor $f_*$ is constructed, the reader will find details in \cite{alf}.

\subsection{Derived linear stacks}\label{subsec:Derivedlinearstacks}

For an affine derived scheme $X=\Spec\, A$, and an $A$-dg-module $E$, 
we can form the linear derived stack associated to $E$, denoted by $\VV(E)$, as follows. 
We consider the $\s$-functor
$$\VV(E) : \ncdga_A \to \Top$$
sending a commutative connective $A$-linear cdga $B$ to 
$\Map_{\dg_A}(E,B)$, the mapping space of $A$-dg-modules from $E$ to $B$. Note that, by adjunction, $\Map_{\dg_A}(E,B) \simeq \Map_{\dg_B}(E\otimes_A^{L}B,B)$. The group
scheme $\Gm$ acts on $\VV(E)$ by acting on $B$ by its weight $1$ action. More concretely, 
$\Gm(B)$ can be identified with the simplicial set of self autoequivalences of $B$
as a $B$-dg-modules. This is, by definition, the simplicial sub-monoid of 
$\Map_{\dg_B}(B,B)$ of endomorphisms of $B$ as a $B$-dg-modules, and thus it acts
naturally on $B$ as an object in the $\s$-category $\dg_B$.

The linear derived stack associated to $E$, $\VV(E)$, will be considered as 
an object in $\dSt_X^{\Gm}$, the $\s$-category of derived stacks over $X$ endowed with a 
$\Gm$-action. These local constructions 
are obviously compatible with base change on $X$, so they can be glued to obtain, for any derived stack $X$, an $\s$-functor
$$\VV : \QCoh(X)^{op} \longrightarrow \dSt_X^{\Gm}.$$

It is shown in \cite{mon} that the $\s$-functor $\VV$ has  
a left adjoint sending an $\Gm$-equivariant derived stack $Y \to X$ to 
its sheaf of functions of weight $1$. Moreover, when restricted to bounded above
quasi-coherent complexes, the counit of this adjunction is an equivalence, so that 
the induced $\s$-functor
$$\VV : \QCoh^{<0}(X)^{op} \longrightarrow \dSt_X^{\Gm}$$
is a full embedding. An object in the image of
$\VV$ will be called a \emph{derived linear stack over $X$}. 
We will mainly use this notion for the more specific case of \emph{perfect} 
complexes on $X$.

On a final note, when $E \in \Perf(X)$ is a perfect complex on $X$, then the
linear derived stack $\VV(E) \to X$ is a relative derived Artin stack which is
strongly finitely presented (see \cite[Prop. 2.23]{tova}). Such $\Gm$-equivariant
derived stacks will be called \emph{perfect derived linear stacks over $X$}. They form
an $\s$-category naturally equivalent to $\Parf(X)^{op}$ by means of the construction 
$E \mapsto \VV(E)$.
Note that
$\VV(E)$ is a derived Artin stack locally of finite presentation over $k$ if $X$
is so. Moreover, the
projection $\VV(E) \to X$ is smooth if and only if $E$ is of positive Tor amplitude, 
and is affine if and only if $E$ is connective.

\subsection{Derived foliations as equivariant derived linear stacks}
\label{subsec:Derivedfoliationsasequivariantderivedlinearstacks}

The importance of derived linear stacks for us is the fact that 
they can be used in order to provide a geometric interpretation 
for derived foliations. This geometric interpretation has its own interests, and
we will use it for instance in order to construct the Hodge filtration 
for general derived foliations. It will also be useful for the notion of derived foliations
in non-zero characteristic discussed in Chapter \ref{chapter:nonzerocharacteristics}.

We recall the group stack $\cH_0$, defined as the semi-direct product of $B\Ga$ and $\Gm$ for the
$\Gm$-action of weight $-1$ on $\Ga$ (see \S \ref{sec:Thegeometryofgradedmixedobjects}). For a derived affine scheme $X=\Spec\, A$, 
we have its \emph{graded loop space} $\cL^{gr}(X/k)$, defined as the derived mapping stack
$$\cL^{gr}(X/k):=\uMap(B\Ga,X) \in \dSt_k.$$
The group $\cH_0$ acts on $B\Ga$, where the $\Gm$ component acts by weight $-1$ and 
$B\Ga$ acts by translation on itself. As a result, $\cH_0$ acts on $\cL^{gr}(X/k)$, and
we will thus consider the graded loop stack as a $\cH_0$-equivariant derived stack. 
Note also that as $X$ is affine, so is $\cL^{gr}(X/k)$, and it is explicitly given by 
$$\cL^{gr}(X/k) \simeq \Spec\, (Sym_{A}(\LL_{A/k}[1])) \simeq \VV(\LL_{A/k}[1]),$$
and as a result we have a natural equivalence of graded cdga's
$$\OO(\cL^{gr}(X/k)) \simeq \DR(X/k).$$
Moreover, as shown in \cite{mrt}, the $\cH_0$-action on $Sym_{A}(\LL_{A/k}[1])=\DR(X/k)$
corresponds, by means of the above identification and our corollary \ref{cpI-4}, 
to its usual graded mixed structure given by the de Rham differential.

When $X$ is a general derived stack, we define the completed graded loop stack 
by imposing descent
$$\ccL^{gr}(X/k):=\mathrm{colim}_{\Spec\, A \to X}\cL^{gr}(\Spec\, A/k) \in \dSt_k^{\cH_0}.$$
Note that $\ccL^{gr}(X/k)$ can not always be identified with the
derived mapping stack $\uMap(B\Ga,X)$. There exists a natural descent morphism
$\uMap(B\Ga,X) \to \ccL^{gr}(X/k)$, but it is not an equivalence
in general (it is so when $X$ is a derived scheme, or more generally
a derived Deligne-Mumford stack). Finally, for an arbitrary morphism
of derived stacks $X \to Y$, we can define the relative graded loop stack
$$\ccL^{gr}(X/Y):=\ccL^{gr}(X/k) \times_{\ccL^{gr}(Y/k)}Y$$
which is now an $\cH_0$-equivariant derived stack \emph{over} $Y$.

We are now ready to define a notion of \emph{geometric derived foliations}, 
based on linear stacks and graded loop stacks. This definition is shown in \cite{alf} to be 
equivalent to Definition \ref{dII-4}. For any derived
stack $X$, we consider its derived graded loop stack $\ccL^{gr}(X/k)$ endowed with 
its natural $\cH_0$-action. It comes equipped with a projection $\pi : \ccL^{gr}(X/k) \to X$
by evaluating at the natural base point of $B\Ga$. We note here that 
even though $\pi$ is not $\cH_0$-equivariant (for the trivial action of $\cH_0$ on $X$), 
it is still endowed with a canonical $\Gm$-equivariant structure, given by considering the induced
action via the canonical section $\Gm \to \cH_0$. 

\begin{df}\label{defgeoderfol}
With the notations above, a \emph{geometric derived foliation $\F$ on $X$ (relative to $k$)} consists of
the datum of a $\cH_0$-equivariant derived stack $\VV(\F)$ over $\ccL^{gr}(X/k)$, 
such the following condition is satisfied:
for any derived affine scheme $S$ and morphism $S \to X$, the $\Gm$-equivariant morphism
$$\VV(\F)\times_{\ccL^{gr}(X/k)}\ccL^{gr}(S/k) \to S$$
makes $\VV(\F)$ into a perfect derived linear stack over $S$.

By definition, 
geometric derived foliations $\F$ on $X$ form a full sub-$\s$-category 
of $\dSt^{\cH_0}_{/\ccL^{gr}(X/k)}$, the $\s$-category of $\cH_0$-equivariant derived 
stacks over $\ccL^{gr}(X/k)$.
\end{df}

As a first comment concerning the above definition, it is important to note that 
$\VV(\F)$ is not a derived linear stack globally over $X$, in the exact same manner as the
big cotangent complex of a derived foliation is not a quasi-coherent complex (see Definition \ref{dII-5}). 
This is expressed
by the strange looking condition in Definition \ref{defgeoderfol}. When $X$ is an affine derived scheme, taking
$S=X$, this condition simply states that $\VV(\F)$ is a derived linear stack over $X$. The same holds, more generally, whenever $X$ is a derived DM stack, since $\cL^{gr}(-/k)$ has \'etale descent.
It is shown in \cite{alf} that the above notion of geometric derived foliations 
is equivalent to our original notion of derived foliations. The equivalence between the two notions
passes through the interpretation of graded mixed dg-modules as quasi-coherent complexes
with $\cH_0$-action (see \S \ref{sec:Thegeometryofgradedmixedobjects}). We will not go into too many details here, and we will limit ourselves to mention 
some basic facts. As a first comment, we can always restrict to the case where $X$ is an affine
derived scheme, and then proceed by gluing for general derived Artin stacks.
For a derived foliation $\F$, on an affine derived scheme $X$, 
the corresponding geometric derived foliation will be given by $\VV(\LL_{\F/k}[1])$, 
the derived linear stack associated to the shifted tangent $\LL_{\F/k}[1]$. One of the result
of \cite{alf} states that the graded mixed structure on $\DR(\F/k)$ induces
a $\cH_0$-action on $\VV(\LL_{\F/k}[1])$ which makes it into a geometric derived foliation in the 
sense of 
Definition \ref{defgeoderfol}. Conversely, if $\F$ is a geometric derived foliation on $X$, 
the corresponding derived foliation can be obtained by considering the push-forward
$$\pi_*(\OO) \in \QCoh(X \times B\Gm),$$
of the structure sheaf along the projection $\VV(\F)\to X$. This quasi-coherent sheaf is endowed 
with a commutative multiplication, as well as a graded mixed structure coming from the action of $\cH_0$ 
(see Corollary \ref{cpI-4}), and is therefore a graded mixed cdga.
The projection $\VV(\F) \to \ccL^{gr}(X/k)$ induces a structure morphism
$\DR_{X/k} \to \pi_*(\OO)$, which makes $\pi_*(\OO)$ into a derived foliation over $X$ (relative to 
$k$).\\

Definition \ref{defgeoderfol} will be reconsidered (with an important variation) in Chapter \ref{chapter:nonzerocharacteristics} in order to define derived foliations in non-zero characteristic.

\chapter{De Rham cohomology of derived foliations}

In this third chapter we study the foliated de Rham cohomology of derived foliations.
For this we introduce
the notion of crystals along derived foliations. Crystals must be thought as
quasi-coherent complexes endowed with partial flat connections along the foliation, 
and are
natural coefficients for foliated de Rham cohomology. In the third section of the chapter we
construct the Hodge filtration associated to a derived foliation $\F \in \Fol(X/k)$, which 
is a filtration on the de Rham cohomology of $X$. Its associated graded are identified with 
the foliated de Rham cohomology with coefficients in exterior powers of the conormal
complex. Finally, the conormal complex of a derived foliation possesses very rich structures, 
and defines a Lie algebra objects inside the $\s$-category of crystals. This can be used to 
study the transversal geometry of a derived foliations, such as transversal jets of functions, 
but also algebraic invariants related to the classical notion of holonomy. \\

All along this section $k$ is a fixed $\QQ$-algebra.

\section{Crystals along derived foliations}\label{sec:crystalsdefinition}

Let $\F \in \Fol(X/k)$ be a derived foliation on a 
derived stack $X$ over $k$,  
$\DR_\F$ the corresponding sheaf of graded mixed cdga's on $X$, 
and 
$$\DR_\F-\egrdg_X:= \lim_{\Spec\, A \to X}\DR_\F(A)-\egrdg_k$$
the $\s$-category of sheaves of graded mixed $\DR_\F$-dg-modules on $X$ (see Definition \ref{dII-6} and 
corollary \ref{cdII-6}).
By definition of derived foliations, 
the weight $0$ part of $\DR_\F^{(0)}$ comes equipped with a canonical identification
to $\OO_X$. Therefore, for any sheaf $E$ of graded mixed $\DR_\F$-dg-modules, 
the weight $n$ part $E^{(n)}$ is endowed with a natural $\OO_X$-dg-module structure.

\begin{df}\label{defqcohcrys}
With the notations above, a \emph{quasi-coherent crystal over $\F$}
is a graded mixed $\DR_\F$-dg-module $E$ satisfying the following two
conditions.

\begin{itemize}
\item The weight $0$ dg-module $E^{(0)}$, considered as a $\DR_\F^{(0)}=\OO_X$-dg-module, 
is quasi-coherent over $X$.

\item For any $n \in \ZZ$, the natural morphism
$$E(0)\otimes_{\OO_X}\DR_\F^{(n)} \longrightarrow E^{(n)}$$
is a quasi-isomorphisms of $\OO_X$-dg-modules.

\end{itemize}
The $\s$-category of quasi-coherent crystals over $\F$ is 
the full sub-$\s$-category $\QCoh(\F)$ of $\DR_\F-\egrdg_X$
consisting of quasi-coherent crystals.
\end{df}

\begin{rmk}\label{rdefqcohcrys}
\emph{It should be noted that quasi-coherent crystals on $X$
are not quasi-coherent sheaves on $X$ in the sense that, except for $n=0$, 
the weight $n$ pieces $E^{(n)}$ are $\OO_X$-dg-modules which are not quasi-coherent over $\OO_X$. 
However, a quasi-coherent crystal is always quasi-coherent over the sheaf of graded mixed
cdga's $\DR_{F/k}$.}
\end{rmk}

Various operations on crystals will be studied in Chapter \ref{ch:operationsoncrystals}, but we will mention here
 the existence of pull-backs.
The $\s$-category $\QCoh(\F)$ is contravariantly functorial in the pair $(X,\F)$ in 
the following sense. Suppose that we have two pairs $(X,\F)$ and $(Y,\G)$ consisting
of derived schemes endowed with derived foliations. A morphism
$f : (X,\F) \longrightarrow (Y,\G)$ will be a pair
$(g,u)$ consisting of 
\begin{itemize}
\item a morphism $g : X \to Y$ of derived schemes,
\item a morphism $u : \F \to g^*(\G)$
of derived foliations over $X$ 
(i.e. a morphism $\DR_{g^*\G}:=  \to \DR_\F$ of graded mixed cdga's over $X$).
\end{itemize}

Associated to such a morphism $f=(g,u)$ there is a pull-back $\s$-functor
$$f^! : \QCoh(\G) \longrightarrow \QCoh(\F)$$
constructed as follows. By definition of pull-backs of derived foliations, we have
an equivalence of graded mixed cdga's on $X$ 
$$\DR_{g^*(\G)}\simeq \DR_X \otimes_{g^{-1}(\DR_Y)}g^{-1}(\DR_\G).$$
The morphism $u$ thus corresponds to a morphism of sheaves of graded mixed cdgas over $X$
under $\DR_X$
$$u : \DR_X \otimes_{g^{-1}(\DR_Y)}g^{-1}(\DR_\G) \longrightarrow \DR_\F$$
or equivalently, to a morphism of sheaves of graded mixed cdga's under $g^{-1}(\DR_Y)$
$$u : g^{-1}(\DR_\G) \longrightarrow \DR_\F.$$
The pull-back $\s$-functor $f^!$ on quasi-coherent crystals is thus simply defined by the following formula
(for $E$ a graded mixed $\DR_Y$-module)
$$f^!(E):=g^{-1}(E) \otimes_{g^{-1}(\DR_\G)}\DR_\F.$$

Clearly, the rule $((X,\F) \mapsto \QCoh(\F), f \mapsto f^!)$ can be 
promoted to an $\s$-functor $\QCoh^{!}$ from the $\s$-category of pairs $(X,\F)$ to the
$\s$-category of $\s$-categories. Note also that 
$\QCoh(\F)$ comes equipped with a natural symmetric monoidal structure (induced by 
the tensor product of graded mixed $\DR_\F$-modules), and that the pull-back $f^!$ has a natural symmetric
monoidal structure as well. 

Moreover, $f^!$ is compatible with the pull-back of quasi-coherent sheaves on derived stacks 
in the following sense.
By definition of $\QCoh(\F)$, we have a forgetful $\s$-functor
$$\QCoh(\F) \longrightarrow \QCoh(X)$$
which sends a graded mixed $\DR_\F$-module $E$ to its weight zero part $E^{(0)} \in \QCoh(X)$. For
a morphism $f=(g,u) : (X,\F) \longrightarrow (Y,\G)$ as above,the compatibility of $f^!$ with the pull-back of quasi-coherent sheaves is expressed by the commutativity of the following square
$$\xymatrix{
\QCoh(\G) \ar[r]^-{f^!} \ar[d] & \QCoh(\F) \ar[d] \\
\QCoh(X) \ar[r]_-{g^*} & \QCoh(Y),
}$$
as this can be easily seen using the explicit formula
$f^!(E):=g^{-1}(E) \otimes_{g^{-1}(\DR_\G)}\DR_\F$ and the condition $(2)$ of Definition \ref{defqcohcrys}. 

The following definition fixes some terminology.

\begin{df}\label{dIII-2}
Let $X\in \dSt_k$ and $\F \in \Fol(X/k)$.
\begin{enumerate}
    \item 
For a quasi-coherent crystal
$E \in \QCoh(\F)$, the \emph{underlying quasi-coherent complex of $E$} is $E^{(0)}$. 

\item For $E' \in \QCoh(X)$, a \emph{crystal structure on $E'$ along $\F$} consists
of $E \in \QCoh(X)$ together with an equivalence $E^{(0)} \simeq E'$. 
     
\end{enumerate}
\end{df}

We warn the reader that we will often refer to objects in $\QCoh(\F)$ by their
underlying objects in $\QCoh(X)$. The typical example is the unit crystal $\OO_X$,
that is the structure sheaf endowed with its canonical crystal structure. 
Strictly speaking this object should be written $\DR_\F$, as it corresponds to 
$\DR_\F$ as a module over itself. However, we will keep referring to this object
with the symbol $\OO_X$, to keep in mind that crystals are somehow $\D$-module
structures (see Chapter \ref{ch:operationsoncrystals} for a precise statement). \\

\subsection{Examples}\label{subsec:crytalsexamples} We conclude this section by listing some basic examples 
of the notion quasi-coherent crystals. \\

\noindent \textbf{Crystals over the trivial foliation.} Assume that $X$ is
a derived Artin stack, so the initial derived foliation exists (see Example 
\ref{ex:tautologicalfoliationsglobal}). Now, for $\F=0_X$ the initial foliation,
the $\s$-category $\QCoh(\F)$ clearly is equivalent to 
$\QCoh(X)$, the $\s$-category of quasi-coherent complexes over $X$. This
equivalence is realized by sending $E\in \QCoh(\F)$ to its weight $0$ part $E(0) \in \QCoh(X)$. It
can be promoted to an equivalence of symmetric monoidal $\s$-categories. \\

\noindent \textbf{Crystals over the tautological foliation.} Assume that $\F=*_X$ is the final 
foliation. If $X$ is a smooth variety, then there is a canonical equivalence of 
$\s$-categories
$$\QCoh(*_X)\simeq \D_X-\dg,$$
between quasi-coherent crystals along $*_X$ and 
quasi-coherent complexes of left $\D_X$-modules. This equivalence will be reviewed and generalized in the next Section \ref{section_foliatedDRcoh}. This equivalence
is again compatible with the natural symmetric monoidal structures involved. When $X$ is not smooth
anymore, and for instance is a derived scheme, we believe that our notion of quasi-coherent crystals
on $*_X$ coincides with the notion of derived $\D_X$-modules of \cite{beraldo_derD}. \\

\noindent \textbf{Crystals over the Dolbeault foliation.} Let $\F=*_{Dol}$ be the Dolbeault foliation 
on
a smooth variety $X$, defined by 
$\DR(\F):=Sym_{\OO_X}(\Omega_X^1[1]))$ where the graded cdga $Sym_{\OO_X}(\Omega_X^1[1]))$ is endowed with the 
\emph{zero} mixed structure. Then $\QCoh(*_{Dol})$ is naturally equivalent to 
the derived $\s$-category of complexes of quasi-coherent Higgs sheaves on $X$
in the sense of \cite{sim}. Dolbeault foliations and the tautological foliations can be combined
together into a \emph{Hodge foliation} denoted by $*_{Hod}$. It is a derived foliation
on $X \times \mathbb{A}^1$ interpolating between $*_{DR}$ (fiber at $1$) and $*_{Dol}$ (fiber at $0$). 
This can be constructed using the deformation to the normal bundle of our appendix \ref{secapp:Deformationtothenormalcone} and will be reviewed in our next section
concerning the Hodge filtration (see \S \ref{sec:Hodgefiltration}).
 \\

\noindent \textbf{Crystals over integrable foliations.} Let $f : X \longrightarrow Y$ be a morphism
of derived schemes and $\F:=f^*(0_Y)$ be the corresponding derived foliation 
(see example \ref{ex:relativefoliations}).
Then $\QCoh(\F)$ is, by definition, the $\s$-category of \emph{relative
$\D$-modules on $X/Y$}. When $f$ is a smooth morphism between smooth varieties 
these relative $\D$-modules can be written as complexes of modules over $\D_{X/Y}$, the sheaf (of algebras) of
relative differential operators along the fibers of $f$. When $f$ is no more supposed to be smooth, 
we will see that $\D_{X/Y}$ only exists as a sheaf of \emph{dg-algebras} on $X$ (see 
Definition \ref{dIV-1}). \\

\noindent \textbf{Representations of Lie algebras as crystals.} Let $X=\Spec\, k$ and 
$\mathfrak{g}$ be a dual almost perfect dg-Lie algebra over a field $k$. By theorem \ref{tII-4}, $L$ 
determines a derived foliation 
$\F_\mathfrak{g} \in \Fol(*)$, for which the graded mixed cdga is the Chevalley-Eilenberg complex
$\DR(\F):=\mathsf{CE}^*(\mathfrak{g})$. There exists a natural $\s$-functor
$$\phi : \mathfrak{g}-\dg_k \longrightarrow \DR(\F)-\egrdg,$$
from $\mathfrak{g}$-dg-modules to the $\s$-category of graded mixed $\DR(\F)$-dg-modules, 
given by the Chevalley-Eilenberg complex for dg-modules. As for theorem \ref{tII-4}, 
it can be shown that the $\s$-functor $\phi$ is an equivalence of $\s$-categories. 
This example of course can be generalized to the case of dg-Lie algebroids. 

\section{Foliated de Rham cohomology}\label{section_foliatedDRcoh}

The notion of quasi-coherent crystals along derived foliations can be used in order to 
define a general notion of \emph{foliated de Rham cohomology with coefficients}. 
Let $X\in \dSt_k$ be any derived stack and $\F \in \Fol(X/k)$ be a derived foliation over $X$.
The symmetric monoidal $\s$-category $\QCoh(\F)$ sits inside $\DR_\F-\egrdg$. We have
a global section $\s$-functor
$$\Gamma : \DR_\F-\egrdg \longrightarrow \egrdg_k,$$
which is a lax monoidal $\s$-functor to graded mixed complexes. It can be composed with 
our Tate realization functor of \S \ref{subsec:Gradedmixedcomplexesandfilteredcomplexes},
in order to get a lax monoidal global section 
$\s$-functor
$$\Gamma : \DR_\F-\egrdg \longrightarrow \cfdg_k$$
to complete filtered complexes. 

\begin{df}\label{def-foldRwithCryscoeff}
The above lax monoidal $\s$-functor $\Gamma$, when restricted to 
$\QCoh(\F)$ is denoted by 
$$\CDR(\F,-) : \QCoh(\F) \longrightarrow \cfdg_k.$$
For any object $E \in \QCoh(\F)$, 
the complex $\CDR(\F,E)$ is called the \emph{foliated de Rham cohomology along $\F$
with coefficients in $E$}.
\end{df}

By construction, $\CDR(\F,E)$ is a complete filtered complex, which is 
functorial in $E$, but also in $\F$ and $X$. The lax monoidal
structure provides natural morphisms of complete filtered complexes
$$\CDR(\F,E) \widehat{\otimes} \CDR(\F,E') \to \CDR(\F,E\otimes E').$$
In particular, if $E$ is any (commutative) algebra object (e.g the unit crystal $\OO_X$),
we get a complete filtered (commutative) dg-algebra of cohomology
$\CDR(\F,E)$. Typically, the de Rham cohomology of the unit crystal along $\F$
provides a complete filtered cdga $\CDR(\F,\OO_X)$, where $\OO_X$ is the structure sheaf
endowed with its canonical crystal structure along $\F$. \\

Coming back to the basic examples of crystals given in \S \ref{subsec:crytalsexamples}, 
we can for each of them relate the foliated de Rham cohomology with more
classical notions. \\

\textbf{Initial and final foliations.} We have $\QCoh(0_X)\simeq \QCoh(X)$ and in this case
foliated de Rham cohomology simply coincides with quasi-coherent cohomologies of $X$ with
coefficients in quasi-coherent complexes. For the final foliation $*_X$ on a smooth
variety, $\QCoh(*_X)$ identifies with the $\s$-category of quasi-coherent $\D$-modules, 
and foliated de Rham cohomology simply is the usual de Rham cohomology of $X$ relative
to $\Spec\, k$. \\

\textbf{Dolbeault foliation.} In this case $\QCoh((*_{Dol})$
is identified with quasi-coherent Higgs complexes over $X$, and foliated de Rham cohomology coincide
here with the \emph{Hodge completed derived Higgs cohomology}. 
This is defined the exact same manner as the 
Hodge complete derived de Rham cohomology. \\

\textbf{Lie algebras representations.} For a derived foliation $\F$ on $\Spec\, k$ determines by 
a dual almost perfect dg-Lie algebra $\mathfrak{g}$, $\QCoh(\F)$ identifies with the
$\s$-category of $\mathfrak{g}$-dg-modules and foliated de Rham cohomology 
is equivalent to the usual Lie algebra cohomology defined via the Chevalley complex.

\section{The Hodge filtration}\label{sec:Hodgefiltration}

In this section we show that a derived foliation $\F$ on a derived stack $X$ induces 
a canonical filtration on the completed de Rham complex $\CDR(X/k)$, which is moreover multiplicative 
and functorial in $X$. We call it the \emph{Hodge filtration associated to $\F$}, as the usual Hodge filtration
on $\CDR(X/k)$ will be recovered for the special case where $\F=0_X$ is the initial derived foliation. When 
the derived foliation $\F$ is integrable by a morphism $f : X \to Y$, the associated Hodge filtration is
intimately related to the Gauss-Manin connection with respect to $f$.

As usual, we start by the case of an affine derived scheme $X=\Spec\, A$, and obtain the globalization by usual gluing. 
Let $\F \in \Fol(X/k)$ be a derived foliation on $X$, which we assume to be almost 
perfect. We consider the
corresponding $\cH_0$-equivariant derived linear stack $p : \VV(\F) \to \cL^{gr}(X/k)$ (see \S \ref{subsec:Derivedfoliationsasequivariantderivedlinearstacks}), an we 
perform the deformation to the normal bundle for the morphism $p$, as recalled in our appendix \ref{secapp:Deformationtothenormalcone}. 
This provides a morphism of $\cH_0$-equivariant derived stacks over $\A=[\mathbb{A}^1/\Gm]$
$$j : \VV(\F) \times \A \to Def(p).$$
Over the open point $[\Gm/\Gm]$, the morphism $j$ recovers the original morphism $p$, whereas 
over the closed point we get the zero section of the shifted relative tangent stack $T(p)[1]$, 
or equivalently $\VV(\LL_p[-1])$, the derived linear stack associated to the shifted cotangent 
complex of the morphism $p$. More generally, $Def(p) \to \A$ is itself a relative derived linear
stack, associated to $\LL_{X/k}[1] \to \LL_{\F}[1]$, considered as a one step filtered 
quasi-coherent complexes on $X$, and thus as a quasi-coherent complexes on $\A$ (see the end of
\S \ref{secapp:Deformationtothenormalcone}).
The cotangent complex $\LL_p$ can furthermore be described 
thanks to the following lemma.

\begin{lem}\label{lpIII-1}
Let $X=\Spec\, A$ be a derived affine $k$-scheme, $E$ be a  bounded above
quasi-coherent complex on $X$ and $\pi : \VV(E) \to X$ the induced
derived linear stack. Then, we have a functorial (in $E$ and $X$) equivalence in $\QCoh(\VV(E))$
$$\LL_{\pi} \simeq \pi^*(E).$$
\end{lem}

\textit{Proof of the lemma.} For any derived linear stack $\pi : \VV(E) \to X$
the relative cotangent complex $\LL_\pi$, by definition, satisfies the following property: for
any $A \to B$ of connective cdga's, andy 
$E \to B$ morphism of $A$-dg-modules corresponding to $u : \Spec\, B \to \VV(E)$, 
and any connective $M \in \dg_B$, we have functorial 
(in $B$ and $M$) equivalences
$$\Map_{\dg_B}(u^*(\LL_\pi),M)\simeq \Map_{\dg_A / B}(E,B\oplus M) \simeq \Map_{\dg_A}(E,M).$$
This shows that $u^*(\LL_\pi) \simeq B\otimes_A E$ functorially in $B$, and thus
that $\LL_\pi \simeq \pi^*(E)$. This identification is furthermore functorial
in $E$ in an obvious manner.
\hfill $\Box$ \\

By passing to $\cH_0$-equivariant global functions on $Def(p)$, we get this way 
a cdga endowed with a complete filtration.
The underlying object of the filtration is given by $\cH_0$-equivariant functions on $\cL^{gr}(X/k)$, 
i.e. 
by the completed de Rham complex $\CDR(X/k)$. The associated graded for this filtration is given by 
$\cH_0$-equivariant functions on $\VV(\LL_p[-1])$, and can be written as 
foliated de Rham cohomology with coefficient in some crystals as follows.

The graded cdga of graded functions on 
$\VV(\LL_p[-1])$, thanks to Lemma \ref{lpIII-1} is given, as a graded 
$\DR(\F)$-cdga, by
$$Sym_{\DR(\F)}(\LL_p[-1])\simeq Sym_{A}(\N_{\F}^*[1])\otimes_{A}\DR(\F),$$
where $\N_\F^*$ is the conormal complex of $\F$ defined as usual as being the fiber of $\LL_{X/k} \to \LL_{\F/k}$. We need to be slightly more precise here. The cdga 
$Sym_{\DR(\F)}(\LL_p[-1]) \simeq Sym_{A}(\N_{\F}^*[1])\otimes_{A}\DR(\F)$
has an action of $\cH_0$, its natural graded mixed structure, and a second action of $\Gm$, as $\VV(\LL_p[-1])$ lives over $B\Gm$. These two actions combined in an action of $\Gm\times \cH_0$, and
defines a graded objects inside the $\s$-category of graded mixed $\DR(\F)$-dg-modules. 
As such we have a decomposition of graded mixed $\DR(\F)$-dg-modules
$$Sym_{\DR(\F)}(\LL_p[-1]) \simeq \bigoplus_{i \geq 0}((\wedge^i \N_\F^*) \otimes_A \DR(\F))[i],$$
but for which $(\wedge^i \N_\F^*)$ is of weight $i$. In other words, the graded mixed structure on 
each piece $(\wedge^i \N_\F^*) \otimes_A \DR(\F)$ is not quite a crystal in the sense of 
Definition \ref{defqcohcrys} because it is not of weight $0$. Rather, 
the graded mixed $\DR(\F)$-dg-modules $(\wedge^i \N_\F^*) \otimes_A \DR(\F)$
are such that their weight shifts ("Tate twists")
$$((\wedge^i \N_\F^*) \otimes_A \DR(\F))(-i)$$
belongs to $\QCoh(\F)$.

Therefore, the $\cH_0$-action on each piece $((\wedge^i \N_\F^*) \otimes_A \DR(\F))(-i)$ 
induces a canonical crystal structure on the quasi-coherent complexes
$\wedge^i\N_{\F}^*$.
The associated graded of the filtration on $\CDR(X/k)$ can thus be written as 
$$|Sym_{\DR(\F)}(\LL_p[-1])|^t \simeq \bigoplus_{i\geq 0} |(\wedge^i \N_\F^*) \otimes_A \DR(\F)|^t
\simeq \oplus_{\geq 0}\CDR(\F,\wedge^i\N_\F^*)[-i],$$
where the extra shifts $[-i]$ comes from the fact that $(\wedge^i \N_\F^*) \otimes_A \DR(\F)$
is of weight $i$ rather than weight $0$ and thus its Tate realization is shifted by $-2i$
according its definition (see \S \ref{subsec:Gradedmixedcomplexesandfilteredcomplexes}).
We gather the above construction into the following proposition.

\begin{prop}\label{pIII-1}
For a perfect derived foliation $\F$ on a derived Artin stack $X$. 
\begin{enumerate}
    \item The 
wedge powers of the cotangent complex $\wedge^i_{\OO_X}\N_\F^*$ admit
canonical crystal structure along $\F$.
\item 
There exists a canonical (descending)
filtration $F^*_\F$ on the (complete filtered) cdga $\CDR(X/k)$ of completed de Rham cohomology on $X$. It possesses moreover the
following properties.
\begin{enumerate}
    \item The filtration $F^*_\F$ is functorial in $X$ and $\F$.
    \item The associated graded, as a graded cdga, is given by for $i\leq 0$
    $$Gr^{i}_{F^*_\F}(\CDR(X/k)) \simeq \CDR(\F,\wedge^{-i}\N_\F^*)[i].$$
\end{enumerate}
\end{enumerate}
\end{prop}

In the above proposition, we have two extreme cases: when $\F=0_X$ is the initial 
derived foliation the filtration $F_\F^*$ is the usual Hodge filtration, and
when $\F=*_X$ is the final foliation then the filtration $F_\F^*$ is the trivial filtration.\\

Interesting examples of the Hodge filtration already appear for derived foliations over the point $*=\Spec\, k$ for a field $k$. 
By corollary \ref{ctII-4}, we know that an almost perfect derived foliation $\F$ on $\Spec\, k$ (relative to $k$) is given by 
a dg-Lie algebra $\mathfrak{g}^\vee$ which is dual almost 
perfect over $k$. The $\s$-category of quasi-coherent 
crystals on $\F$ is equivalent to the $\s$-category of $\mathfrak{g}$ dg-modules, and
de Rham cohomology is given by the Lie algebra cohomology.
The completed de Rham complex of $X$ reduces here to $k$.
The conormal complex $\N_\F^*$ naturally identifies with $\mathfrak{g}[-1]$, and its crystal structure
is given by the coadjoint action. The Hodge filtration 
is then a filtration on $k$ whose associated graded are given by $C^*(\mathfrak{g},Sym^i\mathfrak{g}^\vee)[-2i]$, 
the Lie algebra cohomology with coefficients in the wedge powers of the coadjoint representation.
One way to see this filtration in more geometric terms, is to consider the formal derived stack
$\MC(\mathfrak{g}) \in \dSt_k^*$ associated to $\mathfrak{g}$, as explained in \S \ref{subsec:Derivedstacksandformualmoduliproblems}. 
By corollary \ref{pA-8-1}, we have 
$\CDR(\MC(\mathfrak{g}))\simeq k$, and the usual Hodge filtration on the completed de Rham cohomology provides
the Hodge filtration for the derived foliation $\F$. \\

We finish this part by mentioning another construction of a Hodge filtration, which 
we will call the \emph{longitudinal Hodge filtration} in order to distinguish it from
what we have seen in proposition \ref{pIII-1}.  For any graded mixed cdga $D$, 
we define a new graded mixed cdga $D_{Dol}$ as follows. 
The underlying graded cdga of $D_{Dol}$
will be $D \otimes_k k[T]$ where the new variable $T$ has weight $0$, and the
new mixed structure will be $k[T]$-linear and given by 
$T.\epsilon$, where $\epsilon$ is the mixed structure on 
$D$. This new graded mixed structure lives in $\mecdga_{k[T]}$ and is equipped with
a canonical compatible $\Gm$-action for which $T^n.D \subset D \otimes_k k[T]$ is of weight $n$.
Finally, the weight $0$-part of $D_{Dol}$ is $D^{(0)}[T]$. 

We now apply this to a derived foliation $\F$
on a derived affine $k$-scheme $X=\Spec\, A$ and its corresponding 
graded mixed cdga $\DR(\F/k)$. The graded mixed cdga $\DR(\F/k)_{Dol}$ is, as a graded cdga,
equivalent to $Sym_{A[T]}((\LL_{\F}\otimes_A A[T])[1])$, and thus defines 
a new derived foliation on $X\times \mathbb{A}^1$. Moreover, the $\Gm$-action 
defines a descent data as a derived foliation on $X\times \A$. This new derived foliation will be 
denoted by $\F_{Hod}$. 

\begin{df}\label{dIII-1}
The \emph{Dolbeault foliation associated to $\F$} is $\F_{Hod} \in \Fol(X\times \A)$
defined above.
\end{df}

Note that the pull-back of $\F_{Hod}$ at $X = [X\times \Gm/\Gm] \hookrightarrow X\times \A$
recovers the original $\F$. On ther other hand, the pull-back at the origin $X \to X\times \A$
is given by the graded mixed cdga $Sym_A(\LL_{\F}[1])$ with the zero mixed structure, 
with the induced morphism $\DR(X/k) \to Sym_A(\LL_{\F}[1])$ being the factorization
$\DR(X/k) \to \OO_X \to Sym_A(\LL_{\F}[1])$. This derived foliation will be denoted
by $\F_{Dol}$ and is called the Dolbeaul foliation associated to $\F$. When 
$\F=*_{DR}$, the Dolbeault foliation $\F_{Dol}$ is naturally equivalent to $*_{Dol}$
described in our previous section. Note that $(*_{DR})_{Hod}$, which we will denote by 
$*_{Hod}$, provides an interpolation between $*_{DR}$ and $*_{Dol}$.

Finally, on the level of foliated cohomology, $\CDR(\F_{Hod},\OO)$ now has
an induced filtration, called the \emph{longitudinal Hodge filtration}, whose associated
graded in weight $i$ is $C^*(X,\wedge^{i}\LL_{\F})[-i]$.  We gather the previous
constructions in the following proposition.

\begin{prop}\label{pIII-2}
For any derived foliation $\F$ on a derived stack $X$, we can associated functorially 
$\F_{Hod} \in \Fol(X\times \A)$. The pull-back of $\F_{Hod}$ on 
$X = [X\times \Gm/\Gm] \hookrightarrow X\times \A$ canonically identifies with $\F$, whereas
the pull-back of $\F_{Hod}$ at the origin $X \to X\times \A$ is $\F_{Dol}$. 

The foliated de Rham cohomology of $F_{Hod}$ relative to $\A$ defines 
a (descending) filtration on $\CDR(\F,\OO)$ whose associated graded in weight $i\leq 0$ is
$C^*(X,\wedge^{-i}\mathbb{L}_{\F})[i]$.
\end{prop}

As a final comment, the two Hodge filtrations given by propositions \ref{pIII-1} and \ref{pIII-2}
are different by are produced using similar constructions based on the exact triangle
$$\N_\F^* \to \LL_X \to \LL_\F$$
relating the cotangent complexes and the conormal complex. The filtration 
of \ref{pIII-1} is closely related by the natural $\mathbb{A}^1$-homotopy $\lambda.b$, where
$b : \LL_\F \to \N_\F^*[1]$ is the boundary map. Therefore, at the linear level, this filtration 
is simply induced by the one step filtration $\N_\F^* \to \LL_X$ on $\LL_X$. On the other side
the filtration \ref{pIII-2} is related to the $\mathbb{A}^1$-homotopy $\lambda.a$
where now $a : \LL_X \to \LL_\F$ is the anchor map. Therefore at the linear level
it is given by the one step filtration $\LL_\F[-1] \to \N_\F^*$ on the conormal
complex. These filtrations are both on some independent interest, depending if
one wants to study foliated de Rham cohomology, or rather absolute de Rham cohomology. 
The Hodge filtration will be used later in the study of characteristic classes
associated to foliations (see \S \ref{sec:characteristicclasses}).

\section{Transversal geometry}\label{sec:transversal}

In the previous section we have seen that conormal complex $\N_\F^*$ of a 
perfect derived foliation $\F$ on a derived Artin stack $X$ carries a canonical 
crystal structure along $\F$ and thus lifts to an object in $\QCoh(\F)$. 
This object reflects the tranversal geometry of $\F$, as when $\F$ is integrable
by a morphism $f : X \to Y$ the conormal complex $\N_\F^*$
is given by $f^*(\LL_{Y})$, the pull-back of the cotangent complex of $Y$, equipped with 
its natural crystal structure relative to $f$ (see \ref{sec:crystalsdefinition}). 
In this section we will show moreover
that the object $\N_\F[-1]$ comes equipped with a natural Lie bracket, encoded in 
a graded mixed structure on $Sym_{\OO_X}(\N^*_\F[2])$. More importantly this
Lie bracket is compatible with the crystal structure and makes $\N_\F[-1]$ into 
a Lie algebra object in $\QCoh(\F)$. The Koszul dual to this Lie object 
will be by definition the sheaf of transversal jets, which should be understood as the sheaf 
of formal functions in a tranversal direction along $\F$. 

To start with let us assume as usual that $X=\Spec\, A$ is a derived affine scheme
of finite presentation, and that $\F$ is a perfect derived foliation on $X$ (relative to $k$).
We consider the morphism of graded mixed cdga's $\DR(X) \to \DR(\F)$
and its internal de Rham algebra $\DR^{gr}(\DR(\F)/\DR(X))$ as introduced in 
our Appendix \ref{secapp:InternalderiveddeRhamtheory}. This is a double graded mixed cdga whose underlying 
bi-graded cdga is given by 
$$\DR^{gr}(\DR(\F)/\DR(X)) \simeq Sym_{\DR(\F)}(\DR(\F)\otimes_{A}\N^*_\F[2]) \simeq 
\DR(\F) \otimes_A Sym_{A}(\N^*_\F[2]).$$
This bi-graded cdga is endowed with two different graded mixed structure. This first one
is induced by the mixed structure on $\DR(X) \to \DR(F)$, whereas the second one
is the de Rham differential in $\DR^{gr}(\DR(\F)/\DR(X))$. When considered 
as a graded mixed structure for the former, $\DR(\F) \otimes_A Sym_{A}(\N^*_\F[2])$
defines a commutative algebra object inside $\QCoh(\F)$, the symmetric monoidal
$\s$-category of quasi-coherent crystals along $\F$. The underlying 
quasi-coherent complex is here $Sym_{A}(\N^*_\F[2])$. The second mixed structure tells
us that moreover this commutative algebra object in $\QCoh(\F)$ is endowed with 
an internal graded mixed structure. We can therefore consider its de Rham 
cohomology of the second kind as in Definition \ref{dA-8-1}, and define 
$$\widehat{\OO}_{\N_\F}:=\CDR^{II}(\DR(\F)/\DR(X))=|\DR^{gr}(\DR(\F)/\DR(X))|_2$$
which is a complete filtered commutative algebra object in $\QCoh(\F)$. 

When $X$ is a more general derived Artin stack locally of finite presentation, and
$\F$ is a perfect derived foliation on $X$, we obtain (by gluing the above 
constructions over all affines over $X$) a globally defined object 
$\widehat{\OO}_{\N_\F}$ 
which is a complete filtered commutative algebra object in $\QCoh(\F)$.

\begin{df}
The sheaf of \emph{transversal formal functions} along $\F$ is 
$\widehat{\OO}_{\N_\F}$ defined above. It is a complete filtered commutative algebra object 
in $\QCoh(\F)$.
\end{df}

By construction, the underlying object of $\widehat{\OO}_{\N_\F}$ is a sheaf 
of complete filtered $\OO_X$-linear cdga's on $X$ whose associated graded is
equivalent to $Sym_{\OO_X}(\N_\F^*)$ where $\N_\F^*$ is of weight $-1$ (remind from 
(\ref{eq:associatedgraded}) before proposition \ref{pI-2} 
that the associated graded of the realization $|E|$ of
a graded mixed complexes is of the form $\oplus_p E^{-p}[-2p]$).

To conclude this part, we mention here that 
the object $\widehat{\OO}_{\N_\F}$ can be used in order to define a very general
notion of algebraic holonomy of $\F$, by using the Tannakian formalism. This algebraic
holonomy will be compared to the more classical notion of holonomy in the analytic setting
later in \S \ref{sec:RRHTannakian}. 
Note however that the definition below makes sense over any general bases 
(including non-zero characteristics bases, see \S \ref{chapter:nonzerocharacteristics}). 
It provides 
a notion of 
algebraic holonomy defined over smaller fields, for instance number fields, if $\F$
is itself defined over smaller fields.

Let $\F$ be a perfect derived foliation over a derived Artin stack $X$ locally of finite
presentation over $k$. Let $K$ be a $k$-field and $x \in X(K)$ be a global 
$K$-point. We can consider the fiber functor
$$\omega_x : \Perf(\F) \longrightarrow \Perf(K),$$
sending a perfect crystal along $\F$ to its fiber at $x$. We let 
$<\N_\F^*>^{\otimes} \subset \Perf(\F)$ be the smallest thick rigid tensor triangulated $\s$-category
containing the object $\N_\F^*$, or in other words the thick triangulated sub-$\s$-category
generated by all tensor powers of $\N_\F^*$ and its duals. We restrict the 
fiber functor to a symmetric monoidal $\s$-functor
$$\omega_x : <\N_\F^*>^{\otimes} \to \Perf(K).$$

As explained in \cite{tan}, we can form the stack of $k$-linear tensor automorphisms
$aut^{\otimes}(\omega_x)$ of $\omega_x$. This is the (underived) group stack
$$aut^\otimes(\omega_x) : \Aff_k^{op} \longrightarrow Gp(\T),$$
sending a (genuine) commutative $k$-algebra $A$ to the group-like simplicial monoid
of self-equivalences of $A\otimes_k\omega_x(-)$ considered as an object in the $\s$-category
$$Fun_k^{\s,\otimes}(<\N_\F^*>^{\otimes},\Perf(A)),$$
of $k$-linear symmetric monoidal $\s$-functors from $<\N_\F^*>^{\otimes}$ to $\Perf(A)$.
Its classifying stack $Baut^{\otimes}(\omega_x)$ can be shown to be
a schematic homotopy type over $K$ in the sense of \cite{chaff}. For the sake 
of completeness we include a proof.

\begin{lem}\label{ldef_algholtype}
The (underived) stack $Baut^{\otimes}(\omega_x)$ above is a schematic homotopy type over the field 
$K$.
\end{lem}

\textit{Proof of the lemma.} By Definition \cite[3.1.2]{chaff} we have to prove the following 
two assertions.
\begin{enumerate}
    \item The stack $auu^{\otimes}(\omega_x)$ is affine.
    \item The stack $Baut^\otimes(\omega_x)$ is $P$-local. 
\end{enumerate}

As $aut^\otimes(\omega_x)$ coincide with 
$end^{\otimes}(\omega_x)$ (by the natural inclusion), because $\omega_x$
is a symmetric monoidal $\s$-functor between rigid symmetric monoidal $\s$-categories.
But $end^{\otimes}(\omega_x)$ is always equivalent to a certain limit of stacks of the form
$\underline{End}_K(E)$ for $E$ a perfect complex over $K$, where
$\underline{End}_K(E)$ is the stack of endomorphisms of $E$. Clearly 
this last stack is linear and thus affine, and explicitly given by the spectrum 
of $Sym_K(E \otimes_K E^\vee)$. As affine stacks are stable by limits the first assertion follows. 
Moreover, each stack $\underline{End}_K(E)$ is $P$-local in the sense of \cite[Def. 3.1.1]{chaff}, 
as it classifies endomorphisms of $E$. In other words, for any other $K$-stack $F$, we 
have 
$$\Map(F,\underline{End}_K(E)) \simeq \Map_{\Parf(F)}(p^*(E),p^*(E))$$
for $p : F \to \Spec\, K$ the canonical projection. As $P$-local stacks are also stable
by limits this finishes the proof of the lemma.
\hfill  $\Box$ \\

By the lemma \ref{ldef_algholtype}, we can apply \cite[Thm. 3.2.9]{chaff}, and show that 
the sheaves
$\pi_i(Baut^{\otimes}(\omega_x))\simeq \pi_{i-1}(aut^\otimes(\omega_x))$ 
are all representable by affine group schemes
over $K$, which are furthermore products of additive groups when $i>1$.

\begin{df}\label{def_algholtype}
The \emph{algebraic holonomy type of $\F$ at $x$} is the schematic homotopy type 
$Baut^{\otimes}(\omega_x)$ defined above, and is denoted by $\widehat{Hol}^{alg}_\F(x)$. 
Its homotopy groups are called the 
\emph{abstract algebraic holonomy groups of $\F$ at $x$}.
\end{df}

\begin{rmk}
\begin{enumerate}
    \item \emph{As a first comment, we have used the expression} abstract algebraic holonomy 
    groups \emph{in order to emphasis the existence of a more classical/topological version
    of algebraic holonomy in the analytic setting, which will be discussed in \S \ref{sec:leafspaces}.
    For $\pi_1$ a more concrete definition of algebraic holonomy will be also discussed below
    in Definition \ref{def_alghol}.} 

    \item \emph{The above definition depends on the choice of a global 
    base point $x$ in $X(K)$. There exists a base-point free version by considering 
    instead the whole stack of $k$-linear symmetric monoidal 
    $\s$-functors $<\N_\F^*>  \longrightarrow \Perf$, to the stack of
    perfect complexes. This leads to a stack that deserves the name of}
    foliated algebraic holonomy type of $\F$. \emph{The sheaf of connected components
    of this algebraic homotopy type is an algebraic version of an hypothetical \emph{space of leaves} 
    and can be very degenerate in general. A better behaved leaf space will be introduced
    in the analytic setting in \S \ref{sec:leafspaces}}
    
    \item \emph{Definition \ref{def_algholtype} can also be extended by considering 
    instead the tensor automorphisms of the fiber functor defined now on the whole
    $\s$-category $\Parf(\F)$ rather than just on $<\N_\F^*>$. This leads to a similar schematic
    homotopy type which is called the} algebraic homotopy (or monodromy) type of $\F$ at $x$, 
    denoted by $\widehat{Mon}^{alg}_\F(x)$.
    \emph{By restriction, there is a canonical projection morphism from the algebraic
    monodromy to the algebraic holonomy, which corresponds to some form of a maximal 
    effective quotient, in the same spirit as the maximal effective quotient of 
    a DM stack discussed in \ref{sec:leafspaces}. When $X$ is a smooth variety over
    a field $k$, and $\F=*_{X}$ is the totaulogical foliation, this 
    algebraic holonomy is the differential-Galois homotopy types, and its fundamental
    group is the usual (affine part) of the differential Galois group of $X$. In this case
    the algebraic holonomy is easily seen to be trivial as $\N_\F^*= 0$. Finally, in the
    other extreme case where $\F=0_X$, the algebraic holonomy is trivial
    because of the Tannakian reconstruction theorem (see for instance 
    \cite{zbMATH06755315,zbMATH06902492})}.
\end{enumerate}
\end{rmk}

We want to discuss a more concrete realization of the abstract holonomy groups defined above, in 
particular for $\pi_1(\widehat{Hol}^{alg}_\F(x))$, in order to relate it to 
the more usual notion of holonomy as the action of fundamental groups of leaves
on transversals. For this, we first notice that $\widehat{\OO}_{\N_\F}$
can be considered as a complete filtered algebra objects in the $\s$-category
$Pro(<\N^*_\F>)$ of pro-objects. Indeed, as the filtration on $\widehat{\OO}_{\N_\F}$
is concentrated in non-negative degrees, we can extra from it a pro-object
$$"\lim_{i \leq 0}" \widehat{\OO}_{\N_\F}/F^i.$$
Each individual piece $\widehat{\OO}_{\N_\F}/F^i$ is here the realization of 
the truncated internal de Rham algebra $\DR^{gr,\leq i}(\DR(\F)/\DR(X))$, 
which is itself considered as a graded mixed commutative algebra object in $\QCoh(\F)$
as explained earlier in this section. Each $\widehat{\OO}_{\N_\F}/F^i$ 
is endowed with its natural filtration (as being a realization of a graded mixed
object) whose associated graded is $Sym^{\leq i}_{\OO_X}(\N_\F^*)$. This
implies that indeed each piece $\widehat{\OO}_{\N_\F}/F^i$ is a commutative
algebra object in $<\N^*_\F>$, and thus defines $\widehat{\OO}_{\N_\F} = "\lim_{i \leq 0}" 
\widehat{\OO}_{\N_\F}/F^i$
as a pro-object in commutative algebras in $<\N^*_\F>$.

The symmetric monoidal $\s$-functor $\omega_x$ induces an $\s$-functor
on the corresponding $\s$-categories of commutative algebra objects, and we
can thus consider the action of the $aut^\otimes(\omega_x)$ on 
$\omega_x(\widehat{\OO}_{\N_\F})$ by automorphisms of complete cdga's. We get this way 
a natural morphism of group stacks
$$aut^\otimes(\omega_x) \longrightarrow
aut_{\cdga^{fil}}(\widehat{\OO}_{\N_\F,x})$$
where $\widehat{\OO}_{\N_\F,x}$ is the fiber of $\widehat{\OO}_{\N_\F}$
at $x$, and $aut_{\cdga^{fil}}$ denotes the automorphism group stack 
of self-equivalences of filtered cdga's. Note that $\widehat{\OO}_{\N_\F,x}$
can be morally identified here with the cdga of formal functions at $x$ along
in the transversal direction.
This defines an action of $aut^{\otimes}(\omega_x)$
on the formal transversal at $x$, which is closer to the usual 
notion of a holonomy action. 

To go a tiny bit further, we assume now $\F$ is transversally smooth,
or in other words that $\N_F^*$ is a vector bundle on $X$.
the filtered cdga $\widehat{\OO}_{\N_\F,x}$ is thus a discrete
filtered cdga,  as its associated graded of 
is $Sym_K(\N_{\F,x}^*)$, the symmetric algebra over the fiber $\N_{\F,x}^*$ of $\N_F^*$ at $x$.
In this case, the group stack $aut_{\cdga^{fil}}(\widehat{\OO}_{\N_\F,x})$
is a genuine sheaf of groups, which is moreover representable by an affine group scheme
over $K$ which is abstractly isomorphic to the group scheme $\widehat{G}^d$ of (origin preserving) 
automorphisms
of a formal affine space $\widehat{\mathbb{A}}^d$
where $d$ is the rank of $\N_\F^*$. The morphism
$aut^\otimes(\omega_x) \longrightarrow
aut_{\cdga^{fil}}(\widehat{\OO}_{\N_\F,x})$ thus factors as 
$$\pi_1(\widehat{Hol}^{alg}_\F(x)))=\pi_0(aut^\otimes(\omega_x)) \longrightarrow 
aut_{\cdga^{fil}}(\widehat{\OO}_{\N_\F,x})\simeq \widehat{G}^d.$$

\begin{df}\label{def_alghol}
With the notations and convention above, if $\F$ is a transversally smooth
derived foliation on $X$, and $x \in X(K)$, the \emph{algebraic holonomy
group at $x$} is the image of the morphism of affine group schemes over $K$
$$\pi_1(\widehat{Hol}^{alg}_\F(x))) \longrightarrow 
aut_{\cdga^{fil}}(\widehat{\OO}_{\N_\F,x}).$$
It is denoted by $Hol^{alg}_\F(x)$.
\end{df}

Later in this book, the topological version of holonomy will be introduced
in the complex analytic setting. We will show that $Hol^{alg}_\F(x)$ always
contains the topological holonomy as a sub-group. When $X$ and $\F$ are defined
over smaller fields, for instance number fields, this can provide interesting
informations about the topological holonomy.

\chapter{Operations on crystals}\label{ch:operationsoncrystals}

In this chapter we study furthermore crystals and their functorial properties.
We present in particular a new description of crystals, 
as dg-modules over a certain sheaf of differential operators which is closer to the usual
notion of $\D$-modules. We study the standard functorialities of crystals, pull-backs and push-forwards.
We use the $\D$-module approach in order to introduce the notion of good filtration on crystals and use
it to define characteristic cycles. The Grothendieck-Riemann-Roch theorem then states that 
characteristic cycles are compatible with push-forward along proper and quasi-smooth morphisms. \\

All along this section $k$ is a fixed $\QQ$-algebra.

\section{Sheaves of differential operators}

In this Section we associate to a derived foliation $\F$ on a derived DM-stack $X$, 
 a sheaf of filtered dg-algebras $\D_\F^{\mathrm{fil}}$ on $X$, called the sheaf of \emph{differential operators along} 
 $\F$, and prove that
 dg-modules over the underlying dg-algebra 
 $\D_\F:=(\D_\F^{\mathrm{fil}})^u$ (i.e. $\D_\F^{\mathrm{fil}}$ with filtration forgotten) corresponds to quasi-coherent crystals along $\F$. \\

\subsection{Sheaf of differential operators} 

Let $X=\Spec\, A$ be a derived affine $k$-scheme,  
and $\F$ be a derived foliation on $X$. We have a corresponding 
graded mixed cdga $\DR(\F)$, and a  canonical augmentation of graded mixed cdga's 
$\DR(\F) \to A$. This makes $A$ into a graded mixed $\DR(\F)$-dg-module. 

\begin{rmk}\label{rIV-1}
\emph{The graded mixed $\DR(\F)$-dg-module $A$ is not a crystal along $\F$, as it 
is obviously not free as a graded $\DR(\F)$-module (unless $X$ is étale over $\Spec\, k$). It should
not be confused with the canonical crystal structure on the structure sheaf $\OO_X$ in
the sense of  Definition \ref{dIII-2}, which 
corresponds to $\DR(\F)$ as a graded mixed dg-module over itself.
In this section we will use at several places non-free graded mixed modules, which thus do not
correspond to crystals along $\F$.
}
\end{rmk}

We consider $\rch_{\DR(\F)}(A,A)$, 
the internal Hom object 
of endomorphisms of $A \in \DR(\F)-\medg_X$. This is a graded mixed associative algebra
in the $\s$-category $\DR(\F)-\egrdg_k$ of graded mixed $\DR(\F)$-dg-modules
\begin{equation}\label{Diff}
\rch_{\DR(\F)}(A,A) \in \mathbf{Alg}(\medg_k).
\end{equation} 
Its underlying graded dg-algebra, obtained by forgetting the mixed structure, is explicitly given
by
$$\rch_{\DR(\F)}(A,A)^{\cancel{\epsilon}} \simeq 
\rch_{Sym_A(\LL_\F[1])}(A,A) \simeq \bigoplus_{i\geq 0}(Sym^i_A(\LL_\F))^\vee[-2i],$$
where $(-)^\vee$ denotes the $\OO_X$-linear dual. When $\F$ is a perfect derived foliation, 
$\LL_\F$ is perfect and thus duality commutes with symmetric powers, and in this case we can 
write furthermore
$$\rch_{\DR(\F)}(A,A)^{\cancel{\epsilon}} \simeq Sym_A(\T_\F[-2])$$
where the tangent complex $\mathbb{T}_\F=\LL_\F^{\vee}$ is here of weight $-1$. Note that, in 
general, the mixed structure induced on the right hand 
side is non-trivial. For instance, it induces a morphism in $\dg_k$
from the piece of weight $-2$ to the (shift by $-1$ of the) piece of weight $-1$
$Sym^2_A(\mathbb{T}_\F)[-4] \to \mathbb{T}_\F[-3]$, which defines a extension class
$$Sym^2_A(\mathbb{T}_\F) \to \mathbb{T}_\F[1].$$
The corresponding triangle 
$$\xymatrix{\mathbb{T}_\F \ar[r] & \D_\F^{[1,2]} \ar[r] & Sym^2_A(\mathbb{T}_\F)}$$
which defines the complex of $\D_\F^{[1,2]}$ differential operators of degree bigger than 
$1$ less than $2$ along $\F$.

\begin{df}\label{dIV-1}
The \emph{filtered ring of differential operators of} $\F$ is 
defined to be the realization
$$\D_\F^{\mathrm{fil}}:=|\rch_{\DR(\F)}(A,A)|^t \in \mathbf{Alg}(\fdg_k).$$
The \emph{ring of differential operators along $\F$} is 
the underlying $E_1$-algebra over $X$ obtained by forgetting the filtration
and is denoted by
$$\D_\F:=(\D_\F^{\mathrm{fil}})^u = colim_{i\geq 0}F^i\D_\F^{\mathrm{fil}} \in \mathbf{Alg}(\dg_k).$$
\end{df}

\begin{rmk}\label{rdIV-1}
\emph{The filtered dg-algebra $\D_\F^{\mathrm{fil}}$ depends on the base ring $k$, and should rather
be denoted by $\D_{\F/k}^{\mathrm{fil}}$. We will keep the base ring implicit and continue
to use the notation $\D_\F^{\mathrm{fil}}$ unless necessary.
}
\end{rmk}

By construction, $\D_\F^{\mathrm{fil}}$ is $k$-linear associative dg-algebras. 
As usual, we set 
$$\D_\F^{\leq i}:=F^i\D_\F,$$
and call $\D_F^{\leq i}$ the \emph{ring of differential operators along $\F$ of order $\leq i$}.
Using the explicit description of the $\s$-functor $|-|^t$ (see \S \ref{subsec:Gradedmixedcomplexesandfilteredcomplexes} (\ref{eq:associatedgraded})) 
it is straightforward 
to verify that, when $\F$ is perfect
the associated graded $\mathrm{Gr}(\D_\F^{\mathrm{fil}})$ is naturally equivalent to 
$Sym_{\OO_X}(\mathbb{T}_\F)$, where $\mathbb{T}_\F$ has weight $1$. When $\F$ is not perfect
a more complicated formula holds
$$\mathrm{Gr}(\D_\F^{\mathrm{fil}}) \simeq \bigoplus_{i\geq 0}(Sym^i_A(\LL_\F))^\vee[-2i].$$

Suppose that $\F$ is a smooth derived foliation, and $X$ is a non-derived scheme (e.g. a smooth
variety), 
$\mathbb{T}_\F$ is a vector bundle on $X$, and thus
$Gr^*(\D_\F^{\mathrm{fil}})$ is automatically cohomologically concentrated degree $0$.
As it is a complete filterted object, this implies that  
$\D_\F^{\leq i}$ is also in degree $0$ for all $i$. Therefore, in this case, $\D_\F^{\mathrm{fil}}$
is a genuine filtered $k$-algebra. For more general $X$, $\D_\F$ is a dg-algebra
that might have non-trivial cohomologies in an infinite number of degrees. 
As $\F$ is assumed to be perfect, 
 $\mathbb{T}_\F$ is a perfect $A$-module, as well as 
 its symmetric powers $Sym^i_A(\mathbb{T}_\F)$. 
As a result each $\D_\F^{\leq i}$ is a perfect $A$-module.
It is
moreover concentrated in non-negative degrees when $X$ is a non-derived scheme. 
When $\F$ and $X$ are both quasi-smooth (i.e. $\LL_\F$ and $\LL_X$ are perfect of tor-amplitude $[-1,0]$),
the dg-algebra $\D_F$ is moreover cohomologically bounded. This can be checked, for instance, using the exact triangles 
$$\xymatrix{\D_\F^{\leq i} \ar[r] & \D_\F^{\leq i+1} \ar[r] & Sym^{i+1}_{\OO_X}(\T_\F)
}$$
and induction on $i$. More is true: if $X$ is a non-derived scheme and 
$\T_\F$ can be represented by a two term complex of vector bundles $V \to W$, then $\D_\F$ is
cohomologically concentrated in degrees $[0,rk(W)]$. \\

Definition \ref{dIV-1} can be generalized to non-affine objects by sheafification. However, 
we have to restrict to Deligne-Mumford stacks, as the construction of $\D_\F$ is not 
stable by pull-backs unless we restrict to étale morphisms. It is possible to define
a non-quasi-coherent sheaf of dg-algebras over general derived Artin stacks, but we do not find it
particularly helpful so will not discuss this here.

Let $f : X=\Spec\, A \to Y=\Spec\, B$ be a formally étale morphism between derived 
affine schemes over $k$, and $\F \in \Fol^p(Y/k)$ be a perfect 
derived foliation on $Y$. We contemplate the
commutative diagram of graded mixed cdga's
$$\xymatrix{
\DR(B) \ar[r] \ar[d] & \DR(A) \ar[d] \\
\DR(\F) \ar[r] \ar[d] & \DR(f^*(\F)) \ar[d] \\
B \ar[r] & A
}$$
and realize that the top square is cartesian by Definition of $f^*(\F)$, and the bottom
square is cartesian to because $f$ is formally étale. As a result, the exterior square is
cartesian. Therefore, the base change $\s$-functor on graded mixed dg-modules
$$\DR(f^*(\F))\otimes_{\DR(\F)} - :  \DR(\F)-\egrdg_k \to  \DR(f^*(\F))-\egrdg_k$$
sends $B$ to $A$ via the canonical morphism $\DR(f^*(\F))\otimes_{\DR(\F)} B \to A$.
This implies that the formation of $\D_\F$ is functorial in $f$
and in particular it can be sheafified over the small étale 
site of any derived Deligne-Mumford stack. 

For any derived Deligne-Mumford stack $X$ over $k$ and $\F \in \Fol(X)$, we 
first have a sheaf of graded mixed associative algebras
$\rch_{\DR_{\F,et}}(\OO_X,\OO_X)$, and its realization
$$\D_\F^{\mathrm{fil}} := |\rch_{\DR_{\F,et}}(\OO_X,\OO_X)| \in \Alg(\edg_k(X)),$$
where $\edg_k(X)$ is the $\s$-category of sheaves of graded mixed $k$-modules
over $X_{\text{ét}}$ the small étale site of $X$.

\begin{df}\label{dIV-2}
For a derived DM-stack $X$ over $k$, the \emph{filtered sheaf of differential operators of} a 
derived foliation $\F \in \Fol(X/k)$ is 
$$\D_\F^{\mathrm{fil}}:=|\rch_{\DR_{\F,et}}(\OO_X,\OO_X)|^t \in \mathbf{Alg}(\fdg_k(X))$$
defined above.
The \emph{sheaf of differential operators along $\F$} is 
the underlying $E_1$-algebra over $X$ obtained by forgetting the filtration
and is denoted by
$$\D_\F:=(\D_\F^{\mathrm{fil}})^u = colim_{i\geq 0}F^i\D_\F^{\mathrm{fil}} \in \mathbf{Alg}(\dg_k(X)).$$
\end{df}

When $X$ is a smooth DM-stack over $k$ and $\F=*_{X}$ is the tautological foliation, 
$\D_\F^{\mathrm{fil}}$ is a sheaf of filtered associative 
$k$-algebras  whose associated graded
in $Sym_{\OO_X}(\mathbb{T}_{X/k})$. The following result identifies it 
with the usual ring of differential operators on $X$ relative to $k$.

\begin{prop}\label{pIV-1}
Let $X$ be a smooth Deligne-Mumford stack over $\Spec\, k$. There is a natural isomorphism
of sheaf of filtered $k$-algebras on $X_{\text{ét}}$
$$Diff_{X} \simeq \D_{*_X}^{\mathrm{fil}}$$
where $Diff_X$ is the usual sheaf of rings of differential operators (relative to $k$) 
filtered by order of differential operators.
\end{prop}

\textit{Proof.} We use corollary \ref{clI-4} in order to 
compute $|\rch_{\DR_{\F,et}}(\OO_X,\OO_X)|^t$. It can be identified 
as the internal derived endomorphism objects inside $\CDR(-)-\cfdg_k(X)$, the $\s$-category 
of sheaves of complete filtered complexes of $\CDR(-)$-modules on $X_{\text{ét}}$
$$|\rch_{\DR_{\F,et}}(\OO_X,\OO_X)|^t \simeq \rch_{\CDR(-)-\cfdg_k(X)}(\OO_X,\OO_X).$$
Here $\CDR(-)$ is the sheaf of filtered dg-algebras $(U \to X) \mapsto \CDR(U,\OO_U)$, 
defined on the small \'etale site of $X$ and given by the Hodge-completed derived de Rham 
cohomology of Definition \ref{dI-4}.
As $X$ is smooth this is 
nothing else than the genuine de Rham complex $(\Omega_{X/k}^*,dR)$, with $\Omega_{X/k}^i$
in cohomological degrees $i$, and endowed with its 
Hodge filtration $\Omega^{\geq p}_{X/k} \subset \Omega_{X/k}^*$. The statement of the proposition 
is then that the usual ring of differential operators $Diff_X$ appears as
$$Diff_X \simeq \rch_{\CDR(-)-\cfdg_k(X)}(\OO_X,\OO_X),$$
which is a consequence of the classical Koszul duality between the ring
of differential operators and the de Rham complex. One manner to see this is to use the standard
resolution of $\OO_X$ by locally free module over $Diff_X$ given by the Spencer complex
$Diff_X \otimes_{\OO_X} \wedge^*\T_{X/k}$
$$\xymatrix{
0 \ar[r] & Diff_X \otimes_{\OO_X}\omega_{X/k} \ar[r] &  Diff_X \otimes_{\OO_X}\wedge^{d-1}\T_{X/k}
 \dots \ar[r] & Diff_X \otimes_{\OO_X}\T_{X/k} \ar[r] & Diff_X \ar[r] & \OO_X
}$$
where $d$ is the rank of $\T_{X/k}$ (we can assume $X$ connected). This
resolution implies the existence of a natural quasi-isomorphism of sheaves of filtered dg-algebras 
$$\rch_{Diff_X-\cfdg_k(X)}(\OO_X,\OO_X) \simeq \CDR(-).$$
By adjunction, this produces a canonical morphism of sheaves of complete filtered dg-algebras on $X$
$$\phi : Diff_X \longrightarrow \rch_{\CDR(-)-\cfdg_k(X)}(\OO_X,\OO_X).$$
This last morphism can then be checked to be an equivalence by looking at the associated graded, 
which reduces to the standard equivalence of graded rings
$$Gr^*(\phi) : Sym_{\OO_X}(\T_{X/k}) \simeq \rch_{\Omega^*_{X/k}-\grdg_k(X)}(\OO_X,\OO_X).$$
\hfill $\Box$ \\

Here are some basic examples of Rings of differential operators.

\begin{ex}\label{Dex}
\begin{enumerate}
\emph{\item When $\F$ is the \emph{zero foliation} (i.e. $\mathbb{L}_\F=0$)
then $\D_\F=\OO_X$ endowed with the trivial filtration.
\item When $\F$ is the \emph{tautological foliation} on a smooth variety $X$
(i.e. $\DR(\F)=\DR(X)$), then $\D_\F=\D_X$ is the usual ring
of differential operators endowed with its usual filtration by the order of
operators. This is the proposition \ref{pVI-1} above.
\item When $\F$ is the \emph{Dolbeault foliation} on a smooth variety, that 
is $\DR(\F)=Sym_{\OO_X}(\Omega^1_X[1])$ with trivial mixed structure, 
then $\D_\F=Sym_{\OO_X}(\mathbb{T}_X)$ with the split filtration. More generally,
when $\DR(\F)=Sym_{\OO_X}(\LL_\F[1])$ with trivial mixed structure (abelian 
derived foliation), then $\D_\F\simeq Sym_{\OO_X}(\T_\F)$. 
\item If the foliation $\F$ is \emph{smooth}, i.e it is a Lie algebroid (see Theorem \ref{tII-4} or \cite[Section 1.3.1.]{ToVe-RH}), induced
by a smooth groupoid $G$ acting on $X$, then 
$\D_\F$ is the ring of distributions on $G$, i.e. the $\OO_X$-linear 
dual of formal functions of $G$, endowed with the convolution product.
This coincides with the universal enveloping algebra of the Lie algebroid.
\item When $\F$ is \emph{globally integrable} by a flat 
and generically smooth
morphism of smooth varieties $f : X \longrightarrow Y$, 
$\D_\F$ is called the dg-algebra of \emph{relative differential operators}. 
The reason for this name comes from the fact that $H^0(\D_\F)$ is indeed a subring of 
$\D_X$ consisting of differential operators stabilizing the
fibers of $f$. }
\end{enumerate}
\end{ex}

\begin{rmk}\label{weyl2} 
\begin{enumerate}
    \item \emph{It is interesting to note the following basic example. Let $X=\mathrm{Spec}( k[u])$ with $u$ in degree $-2$, and $\F=*_X$ be the tautological (i.e. final) derived foliation on $X$. Then
$\D_\F$ is here the Weyl dg-algebra over one generator $u$ in degree $-2$. In other
words, it is the dg-algebra freely generated by two cocycle $u$ and $\frac{\partial}{\partial u}$, 
respectively in degrees $-2$ and $2$, with the usual commutation relation
$[u,\frac{\partial}{\partial u}]=-1.$ Note that $\D_\F$ is \emph{not} Morita equivalent to $k$, as opposed to the case
when $\deg u$ is odd (see \cite[Proof of Cor. 4.3.13]{beraldo_derD}). The reader will find more details about $\D_{*_X}$ in \cite{beraldo_derD}, even for more general $X$ than the ones considered here.}
\item \emph{The previous example shows that our notion of differential operators, and thus
of crystals (see below) differs from the notion of $\D$-modules use in \cite{garo}. 
For more precise comparison we refer to \cite{carlob}.}
\end{enumerate}
\end{rmk}

For $X$ a fixed derived DM-stack, 
the construction $\F \mapsto \D_{\F}$ defines an $\s$-functor $\Fol(X) \to \mathbf{Alg}(\fdg(X))$: 
if $\F \to \F'$ 
is a morphism in $\Fol(X)$, given by a morphism $\varphi: \DR_\F' \to \DR_\F$ 
of sheaves of graded mixed 
cdga's over $X$, $\varphi$ induces a morphism on internal endomorphism objects 
$\rch_{\DR_\F}(\OO_X,\OO_X) \to \rch_{\DR_\F'}(\OO_X,\OO_X)$,
and thus a morphism $\D_\F^{\mathrm{fil}} \to \D_{\F'}^{\mathrm{fil}}$ in $\mathbf{Alg}(\fdg(X))$. In particular,
the maps $0_X \to \F \to *_X$, from the initial and to the final foliations, provide
maps of filtered dg-algebras over $X$
$$\xymatrix{
\OO_X \ar[r] & \D_\F^{\mathrm{fil}} \ar[r] & \D_X^{\mathrm{fil}},}$$
where $\OO_X$ has the trivial filtration, and $\D_X^{\mathrm{fil}}$ is by definition
the ring of differential operators on $X$ with its natural filtration by order of differential 
operators. Moreover, the first morphism exhibits $\OO_X$ as $F^0\D_\F^{\mathrm{fil}}$.
By induction on the filtration it is easy to see, when $\F$ is assumed to be perfect, that
$\D_\F$ is a quasi-coherent complex on $X$, for the left $\OO_X$-module structure given by the morphism
above $\OO_X \to \D_\F$

\subsection{Quasi coherent $\F$-crystals and $\D_\F$-modules}

As recalled in \S \ref{clI-4}, the realization functor on 
$X$ provides an equivalence of  
symmetric monoidal $\s$-categories
$$|-|^t : \medg_k \simeq \cfdg_k,$$
where the right hand side is endowed with the complete tensor product of
complete filtered complexes. By passing to sheaves on the small étale sites, 
we have a symmetric monoidal
equivalence
$$|-|^t : \medg_k(X) \simeq \cfdg_k(X),$$
for any derived DM-stack $X$ over $k$.

Let us fix a derived DM-stack $X$ and $\F\in \Fol(X)$ a derived foliation on $X$.
We will construct an $\s$-functor
$$\theta_X : \QCoh(\F) \longrightarrow \D_\F-\dg,$$
from quasi-coherent crystals along $\F$ to sheaves of (left) $\D_\F$-dg-modules on $X_{et}$. We start, as in formula (\ref{Diff}) above, by letting 
$$S_\F:=\rch_{\DR_{\F,et}}(\OO_X,\OO_X) \in \Alg(\medg(X)),$$
the internal Hom object 
of endomorphisms of $\OO_X$ as a sheaf of graded mixed $\DR_{\F,et}$-modules over $X_{et}$.

The object $\OO_X$ can thus be considered as a graded mixed bi-module with 
its right action by $\DR_{\F,et}$ and left action by $S_\F$. Therefore it can be used in order
to produce the following $\s$-functor (between categories of left modules)
$$\OO_X \otimes_{\DR_{\F,et}} - : \DR_{\F,et}-\medg \longrightarrow S_\F-\medg.$$
Note that this $\s$-functor obviously sends  
$\DR_{\F,et}$-dg-modules on $X$ which are quasi-coherent on $X$ and graded free on weight zero, i.e. quasi-coherent crystals along $\F$, to 
graded mixed $S_\F$-dg-modules whose underlying graded modules are pure of weight $0$ and $\mathcal{O}_X$-quasi-coherent. We thus get 
an induced $\s$-functor
$$\theta'_X: \QCoh(\F) \longrightarrow S_\F-\medg_{qcoh,w=0},$$
where $S_\F-\medg_{qcoh,w=0} \subset S_\F-\medg$ is the full sub-$\s$-category of objects
whose underlying graded $S_\F$-module are pure of weight $0$ and quasi-coherent over 
$\OO_X$.

We now compose the previous $\s$-functor with the Tate realization $|-|^t$ in order to get 
an $\s$-functor to dg-modules over $\D_{\F}$
$$\theta_X : \xymatrix{\QCoh(\F) \ar[r]^-{\theta'_X} & S_\F-\medg_{qcoh,w=0} \ar[r]^-{|-|^t} & \D_{\F}-\dg.}$$ 

\begin{thm}\label{pdmod}
Let $X$ be a derived DM-stack over $k$ and $\F\in \Fol^p(X)$ a perfect derived foliation.
The $\s$-functor $\theta_X: \QCoh(\F) \to \D_\F-\dg$ defined above is fully faithful. Its essential image, $\D_\F-\dg^{qcoh}_{X}$, consists of 
all $\D_\F$-modules over $X$ which are quasi-coherent as $\OO_X$-modules
$$\theta_X : \QCoh(\F) \simeq \D_\F-\dg^{qcoh}_{X}.$$
\end{thm}

\noindent \textit{Proof.} The key lemma to prove the theorem is the following, which 
is a well known version of Koszul duality in the graded context.

\begin{lem}\label{lpdmod}
Let $A$ be a cdga and $M$ a perfect $A$-dg-module. We let $C:=Sym_A(E)$. Then, the following holds.
\begin{enumerate}
    \item There exists a natural identification of graded dg-algebras
    $$\rch_{C}(A,A)\simeq Sym_A(E^\vee[-1]) =: D.$$
    
    \item Using the identification above makes $A$ into a bi-dg-module over $(C,D)$, and the
    induced $\s$-functor on $\s$-categories of graded modules
    $$A\otimes_C - : C-\grdg \longrightarrow D-\grdg$$
    induces an equivalence from the full sub-$\s$-category of graded $C$-modules which are graded
    free over their part of weight $0$, and the graded $D$-modules which are pure of weight $0$
    as graded modules over $k$ (i.e. comes from graded $A$-dg-modules via the augmentation $D \to A$).
\end{enumerate}
\end{lem}

\textit{Proof of the lemma.} This is pretty standard. The $\s$-functor $A\otimes_C -$
graded modules of the form $C\otimes_A M(0)$ to $M(0)$ with the trivial $D$-module
structure induced by the augmentation $D \to A$. Conversely, the right adjoint 
to $A\otimes_C -$ is $\rch_{D}(A,-)$. We first prove that for any $A$-dg-module $M(0)$, considered
as graded pure of weight $0$, the 
adjunction morphism
$$C\otimes_A M(0) \to \rch_{D}(A,M(0))$$
is an equivalence. For this, we use the explicit cofibrant replacement $QA \twoheadrightarrow A$
in $D-\grdg$ given explicitly by 
$QA := Sym_A(E^\vee\oplus E^\vee[-1])$ where the differential on 
$E^\vee \oplus E^\vee[-1]$ is the usual differential
making it contractible $\begin{pmatrix} d & 0 \\ id & -d \end{pmatrix}$.
This implies that $\rch_D(A,M(0))$ is equivalent to 
$\rch_A(Sym_A(E^\vee),M(0))$, the graded hom complexes over $A$ of morphisms from $Sym_A(E^\vee)$ 
to $M(0)$.
But, as $E$ is perfect, and as these are computed in the graded modules $\s$-categories, we have
$$\rch_A(Sym_A(E^\vee),M(0)) \simeq Sym(E)\otimes_A M(0)$$
as required. This shows that the adjunction $(A\otimes -, \rch_D(A,-))$ restricts to an 
equivalence on the $\s$-categories stated in the lemma.
\hfill $\Box$ \\

We now prove the theorem, and 
we prove first that the $\s$-functor
$$\theta'_X: \QCoh(\F) \longrightarrow S(\F)-\medg_{qcoh,w=0}$$
is an equivalence of $\s$-categories, and then identify the essential image of $\theta_X$. 

For this, we consider the $\s$-functor that forgets the mixed structures and only retain
the graded structures. We have obtain this way a commutative diagram of $\s$-categories
$$\xymatrix{
\DR_{\F,et}-\egrdg \ar[r]^-{\theta'_X} \ar[d] & S_\F-\egrdg \ar[d] \\
\DR_{\F,et}-\grdg \ar[r] & S_\F-\grdg
}$$
where the vertical $\s$-functors consists of forgetting the mixed structures. The horizontal
bottom $\s$-funcor is the one of the lemma $\OO_X \otimes_{\DR_{\F,et}}-$. Moreover, 
the above square is right adjointable in the sense of see \cite[Def. 4.7.4.13]{HA}, and the
right adjoints are given by $\rch_{S_\F}(\OO_X,-)$. The lemma \ref{lpdmod} implies that 
$\theta'_X$ is an equivalence of $\s$-categories, as the unit of conunit of the involved adjunction
can be checked on the level of the underlying graded objects without mixed structures. \\

To finish the proof of the theorem, we now consider the realization
$$|-|^t : S_\F-\medg_{X} \longrightarrow \D_\F^{\mathrm{fil}}-\fdg,$$
and observe that, by definition, $\D_\F^{\mathrm{fil}}:=|S_\F|^t$. 
We know that this $\s$-functor is fully faithful, and its
image consists of complete filtered $\D_\F^{\mathrm{fil}}$-modules over $X$. We restrict this
to $S_\F-\medg_{qcoh,w=0}$, the full sub-$\s$-category of graded mixed
module which are quasi-coherent and of weight $0$. Its image by $|-|^t$ is then easily checked to consist of 
all filtered $\D_\F^{\mathrm{fil}}$-modules $E$ satisfying the following two conditions

\begin{enumerate}
\item The filtration on $E$ is tautological: $F^i(E)=E$ if $i\geq 0$ and $0$ if $i<0$.
\item The underlying $\OO_X$-module of $E$ is quasi-coherent.
\end{enumerate}

Now, these two conditions define a full sub-$\s$-category of $\D_\F^{\mathrm{fil}}-\fdg$
which is equivalent, via the underlying object 
$\s$-functor $(-)^u : \fdg_k(X) \to \dg_k(X)$, to the $\s$-category
$\D_\F-\dg_{qcoh}$, of \emph{unfiltered} $\D_\F$-modules that are 
$\mathcal{O}_X$-quasi-coherent. \hfill $\Box$ \\

The following definition/notation will be used throughout the rest of the Chapter.

\begin{df}\label{dqcoh}
\emph{Let $\F \in \Fol(X)$ be a perfect derived foliation over a derived scheme $X$. We denote by 
$\QCoh^{\mathrm{fil}}(\F)$ the full sub-$\s$-category of $\D_\F^{\mathrm{fil}}-\fdg$ consisting of
all filtered modules which are quasi-coherent as filtered $\OO_X$-modules via restriction of scalars along 
the natural
morphism of filtered dg-algebras $\OO_X \longrightarrow \D_\F^{\mathrm{fil}}$. We will call 
$\QCoh^{\mathrm{fil}}(\F)$  
the \emph{$\s$-category of filtered quasi-coherent crystals along $\F$}.}
\end{df}

\begin{rmk}
\emph{Without the perfect condition the theorem \ref{pdmod} does not hold. In fact, 
in general the sheaf $\D_\F$ has no reasons of being quasi-coherent over $\OO_X$, even 
if $\F$ is almost perfect.}
\end{rmk}

\subsection{Inverse image and induction $\s$-functors}\label{s-ind}

Let $u : \F \to \G$ be a morphism of derived foliations on a derived DM-stack $X$, which 
we assume to be perfect in order to be able to use our theorem \ref{pdmod}. We have seen that it 
induces a morphism of sheaves of filtered dg-algebras over $X$
$$\D_\F^{\mathrm{fil}} \longrightarrow \D_\G^{\mathrm{fil}}.$$
Associated to this, we have  the usual forgetful and base-change adjunction on the corresponding
$\s$-categories of sheaves of filtered modules
\begin{equation}\label{for-filt}
u_!:=\D_\G^{\mathrm{fil}}\otimes_{\D_\F^{\mathrm{fil}}}- : \D_\F^{\mathrm{fil}}-\fdg \rightleftarrows \D_\G^{\mathrm{fil}}-\fdg : u^!.\end{equation}
The right adjoint $u^!$ is called the \emph{inverse image $\s$-functor}. 
The left adjoint $u_!$ is called the \emph{induction along $u$} or \emph{direct image $\s$-functor}. 
By forgetting the filtrations, we have a corresponding non-filtered adjunction
\begin{equation}\label{for-nofilt} u_!:=\D_\G\otimes_{\D_\F}- : \D_\F-\dg \rightleftarrows \D_\G-\dg : u^!.
\end{equation}
Both $\s$-functors $u_!$ and $u^!$ preserve quasi-coherence, and thus induce
an adjunction on quasi-coherent modules. By Theorem \ref{pdmod}, this can also be interpreted
as an adjunction on the $\s$-category of quasi-coherent crystals
$$u_! : \QCoh(\F) \rightleftarrows \QCoh(\G) : u^!$$
where again, $u^!$ is called the inverse image $\s$-functor, and $u_!$ the induction or direct image
$\s$-functor. Through Definition \ref{dqcoh}, the filtered adjunction (\ref{for-filt}) can be regarded as an adjunction on filtered crystals, as well
$$u_! : \QCoh^{\mathrm{fil}}(\F) \rightleftarrows \QCoh^{\mathrm{fil}}(\G) : u^!.$$
The filtered and unfiltered versions of $u_{!}$ and $u^{!}$ are of course compatible with the underlying object $\s$-functor, i.e. the
following squares naturally commutes
\begin{equation}\label{dunno1}\xymatrix{
\QCoh^{\mathrm{fil}}(\F) \ar[r]^{u_!} \ar[d]_-{(-)^u} & \QCoh^{\mathrm{fil}}(\G) \ar[d]^-{(-)^u} &  & \QCoh^{\mathrm{fil}}(\G) 
\ar[r]^{u^!} \ar[d]_-{(-)^u} & \QCoh^{\mathrm{fil}}(\F) \ar[d]^-{(-)^u} \\
\QCoh(\F) \ar[r]_-{u_!} & \QCoh(\G) &  & \QCoh(\G) \ar[r]_-{u^!} & \QCoh(\F)
}
\end{equation}
The same is true when the underlying object $\s$-functor $(-)^u$ is replaced with the associated graded $\s$-functor $\mathrm{Gr}$. 
\begin{equation}\label{dunno}\xymatrix{
\QCoh^{\mathrm{fil}}(\F) \ar[r]^{u_!} \ar[d]_-{(-)^u} & \QCoh^{\mathrm{fil}}(\G) \ar[d]^-{Gr} &  & \QCoh^{\mathrm{fil}}(\G) 
\ar[r]^{u^!} \ar[d]_-{(-)^u} & \QCoh^{\mathrm{fil}}(\F) \ar[d]^-{Gr} \\
Gr(\D_\F^{\mathrm{fil}})-\grdg \ar[r]_-{u_!} & Gr(\D_\G^{\mathrm{fil}})-\grdg &  & Gr(\D_\G^{\mathrm{fil}})-\grdg \ar[r]_-{u^!} &  Gr(\D_\F^{\mathrm{fil}})-\grdg
}
\end{equation}

\begin{rmk}
\emph{Tracking back the equivalence of Theorem \ref{pdmod} it is easy to see that the inverse image functor $u^! : \QCoh(\F) 
\longrightarrow \QCoh(\G)$
may also be identified  with the base change at the level of graded mixed dg-modules
$$\DR(\F) \otimes_{\DR(\G)} - : \DR(\G)-\medg_X \longrightarrow \DR(\F)-\medg_X$$
for the morphism $\DR(\G) \longrightarrow \DR(\F)$ corresponding to $u : \F \to \G$ in $\Fol(X)$.
In other words, it does coincide with the $!$-pull-back of quasi-coherent crystals previously defined 
in Section \ref{sec:crystalsdefinition} along the morphism of pairs $(id_X,u) : (X,\F) \to (X,\G)$.}
\end{rmk}

We will now examine two specific important cases of adjunctions (\ref{for-filt}) and (\ref{for-nofilt}):
when either $\F$ is the initial foliation, or $\G$ is the final foliation (so that the morphism $u$ is 
uniquely defined in either cases). 

Let us first consider the morphism $u: 0_X \to \F$. We know that 
$\QCoh(0_X)$ is naturally equivalent to $\QCoh(X)$. The corresponding induction $\s$-functor
$$u_! : \QCoh(X) \longrightarrow \QCoh(\F)$$
is then simply induced by $\D_\F \otimes_{\OO_X}-$. We warn the reader that 
this $\s$-functor does not have an easy description on the level of graded mixed modules, and this shows a particular instance of the usefulness of Theorem \ref{pdmod}. For example, $u_!$ 
sends $\OO_X$ to $\D_\F$ which is a rather big and complicated object inside $\QCoh(\F)$, not 
concentrated in degree $0$ (except if $\F$ and $X$ are both assumed to be smooth). The $\s$-functor
$u_!$ will play an important role for us later and will be referred to as \emph{the induction $\s$-functor
for $\F$}. There is also a corresponding filtered version 
$$u_! : \QCoh^{\mathrm{fil}}(X) \longrightarrow \QCoh^{\mathrm{fil}}(\F),$$
which takes a filtered object $E$ in $\QCoh(X)$ and sends it to 
$\D_\F^{\mathrm{fil}}\otimes_{\OO_X}E$.

\begin{df}\label{dinduc}
Let $\F \in \Fol(X)$ be a perfect derived foliation and $u : 0_X \to \F$ the canonical morphism. 
The \emph{induction $\s$-functor for $\F$} is the $\s$-functor
$$\mathsf{Ind}_\F:=u_! : \QCoh(X) \longrightarrow \QCoh(\F).$$
The \emph{filtered induction $\s$-functor for $\F$} is the $\s$-functor
$$\mathsf{Ind}^{\mathrm{fil}}_\F:=u_! : \QCoh^{\mathrm{fil}}(X) \longrightarrow \QCoh^{\mathrm{fil}}(\F).$$
\end{df}

A direct consequence of the existence of the induction $\s$-functor is the following important
fact.

\begin{cor}\label{ccompact}
Let $X$ be a derived DM-stack over $lk$. If the $\s$-category $\QCoh(X)$ admits a compact generator,
then so are the $\s$-categories $\QCoh^{\mathrm{fil}}(\F)$ and $\QCoh(\F)$. 
\end{cor}

\noindent \textit{Proof.} Pick a family of compact generators
$\{E_i\}_{i\in I}$ in $\QCoh(X)$. It is formal to check that the family of induced objects
$\mathsf{Ind}_\F(E_i) \in \QCoh(\F)$ is a family of compact generators.
In the filtered case, we first consider the family of objects 
$(i_n)_!(E_i) \in \QCoh^{\mathrm{fil}}(X)$, for $\in \ZZ$ and $i\in I$. Here 
$(i_n)_!(E_i)$ is the filtered object in $\QCoh(X)$ whose filtration layers are all zero 
below $n$ and are equal to $E_i$ in layers above $n$. In other words, we have
$$\Map_{\QCoh^{\mathrm{fil}}(X)}((i_n)_!(E_i),E') \simeq \Map_{\QCoh(X)}(E_i,F^nE').$$
The family of objects $\{(i_n)_!E_i\}_{n,i}$ is a family of compact generators
of $\QCoh^{\mathrm{fil}}(X)$, and it is formal to check that the induced family
$\mathsf{Ind}^{\mathrm{fil}}_\F(\{(i_n)_!E_i)\}_{n,i}$ is a again a family of compact
generators of $\QCoh^{\mathrm{fil}}(\F)$.
\hfill $\Box$ \\

\begin{cor}\label{ccompact2}
Let $X$ be either a quasi-compact and quasi-separated derived $k$-scheme, or 
a separated and quasi-compact derived Deligne-Mumford $k$-stack, and let $\F \in \Fol^p(X/k)$
a perfect derived foliation on $X$. Then 
the $\s$-categories $\QCoh^{\mathrm{fil}}(\F)$ and $\QCoh(\F)$ are both compactly generated.
\end{cor}

\textit{Proof.} If $X$ is a quasi-compact and quasi-separated derived scheme
then $\QCoh(X)$ is compactly generated by \cite[Thm. 3.7]{deraz}. When $X$ is a separated
and quasi-compact derived DM-stack it possesses a coarse moduli space $\pi : X \to M$
where $M$ is a separated and quasi-compact derived algebraic space (see for instance
\cite[Thm. 3.12]{ahlqvist2023goodmodulispacesderived}). We can write 
$\QCoh(X) \simeq \Gamma(M_{et},\pi_*(\QCoh))$, the global sections of the 
quasi-coherent stack $\pi_*(\QCoh)$ on $M_{et}$, the small étale site of $M$. Locally on 
$M_{et}$ the derived stack $X$ is equivalent to a quotient stack $[V/G]$ with $V$ affine and
$G$ a finite group, and thus $\pi_*(\QCoh)$ is locally compacty generated as a quasi-coherent
stack on $M_{et}$.
We can thus apply \cite[Thm. 10.3.2.1]{SAG}, which implies that 
$\Gamma(M_{et},\pi_*(\QCoh))\simeq \QCoh(X)$ is compactly generated.

We thus have seen that in both cases we can apply the corollary \ref{ccompact} and therefore have
proven the result.
\hfill $\Box$ \\

\begin{rmk}\label{rccompact2}
\emph{Slightly more is true for corollary \ref{ccompact2}, as we also know that the compact 
objects in $\QCoh(X)$ are the perfect complexes, are both results \cite[Thm. 3.7]{deraz} and
\cite[Thm. 10.3.2.1]{SAG} identify the compact objects in the global section $\s$-categories as
the objects being compact locally.}
\end{rmk}

The other interesting special case of (\ref{for-nofilt}) is that of the unique morphism 
$u : \F \to *_X$ to the final foliation $*_X$. By definition $\D_{*_X} = \D_X$ 
is the sheaf of differential operators on $X$, which 
coincides with the usual sheaf of differential operators when $X$ is smooth 
(see proposition \ref{pIV-1}). The induction $\s$-functor, in this
situation, produces a $\s$-functor
$$u_ ! : \QCoh(\F) \longrightarrow \D_X-\dg^{qcoh},$$
from quasi-coherent crystals along $\F$ to quasi-coherent $\D_X$-modules on $X$. The $\D_X$-modules 
of the form $u_!(E)$ will be called \emph{induced from the foliation $\F$}. One important example  
is the induced $\D_X$-module $u_!(\OO_X)$. This is a canonical $\D_X$-module on $X$ 
associated to the derived foliation $\F$ which contains interesting informations about $\F$. For instance, 
when $\F$ is smooth, then $u_!(\OO_X)$ is a coherent $\D_X$-module. However, this is not true anymore
for non-smooth derived foliations $\F$. Measuring the defect of coherence of $u_!(\D_X)$ is
a very interesting question related to invariants of singularities of derived foliations, generalizing 
classical invariants such as Milnor numbers. We will touch the general study of singularities of
derived foliations in this book and details will appear in a future work.

\section{Direct and inverse images of filtered crystals}
In this Section we define (filtered and unfiltered) direct image functors between quasi-coherent crystals, along proper and quasi-smooth maps. We also prove the important result that filtered direct and filtered inverse images commutes both with the underlying and the associated graded objects functors.

\subsection{Direct images}\label{s-dirim} 
We have seen in \S \ref{sec:crystalsdefinition} that $(X,\F) \mapsto \QCoh(\F)$ is a contraviant
$\s$-functor by using $!$-pull-backs. This functoriality can be extended to the
filtered case as follows. Let $f=(g,u) : (X,\F) \longrightarrow (Y,\G)$
be a morphism of pairs consisting of derived DM-stacks endowed with perfect derived foliations. 
Thus $f$ is given by a morphism
$g : X \rightarrow Y$ of derived stacks, 
and a morphism $u : \F \rightarrow g^*(\G)$ of derived foliations on $X$
 (i.e. a morphism $\DR_{g^*(\G)} \to \DR_\F$
of graded mixed cdga's over $X$). 
This pair 
induces a morphism on the corresponding 
sheaves of differential operators
$$\D_\F^{\mathrm{fil}} \to \D^{\mathrm{fil}}_{g^*(\G)}.$$
We consider the sheaf of filtered $k$-modules on $X$ defined by 
$$\D_{f}^{\mathrm{fil}} := |\rch_{g^{-1}(\DR_\G)-\egrdg}(g^{-1}(\OO_Y),\OO_X)|^t.$$
This obviously carries a right filtered $g^{-1}(\D_\G^{\mathrm{fil}})$-module structure by 
 via the canonical morphism of sheaves filtered dg-algebras
$$g^{-1}(\D_\G^{\mathrm{fil}})=g^{-1}|(\rch_{\DR_\G-\egrdg}(\OO_Y,\OO_Y)|^t \to 
|\rch_{g^{-1}(\DR_\G)-\egrdg}(g^{-1}(\OO_Y),g^{-1}(\OO_Y))|^t.$$
The filtered complex $\D_f^{\mathrm{fil}}$ also carries a left filtered 
by $\D_{g^*(\G)}^{\mathrm{fil}}$ by means of the natural morphism
of sheaves of filtered dg-algebras
$$\D^{\mathrm{fil}}_{g^*(\G)} = |\rch_{(\DR_X \otimes_{g^{-1}(\DR_Y)}g^{-1}(\DR_\G))-\egrdg}(\OO_X,\OO_X)|^t
\longrightarrow |\rch_{(g^{-1}(\DR_\G))-\egrdg}(\OO_X,\OO_X)|^t.$$
The filtered bi $(\D_{g^*(\G)}^{\mathrm{fil}},g^{-1}(\D_\G^\mathrm{fil})$-module
$\D_f^{\mathrm{fil}}$ defines a pull-back
$\s$-functor on filtered modules 
$$\phi : \D_{\G}^{\mathrm{fil}}-\fdg \longrightarrow \D_{g^*(\G)}^{\mathrm{fil}}-\fdg$$
by  the standard formula
$$\phi(E):=\D_{f}^{\mathrm{fil}}\otimes_{g^{-1}(\D_\G^{\mathrm{fil}})}g^{-1}(E).$$
Finally, we can compose this with the forgetful $\s$-functor for the morphism 
of filtered dg-algebras $u : \D_{g^{*}(\G)}^{\mathrm{fil}} \to \D_\F^{\mathrm{fil}}$
in order to get the requires pull-back $\s$-functor
$$f^! : \D_{\G}^{\mathrm{fil}}-\fdg \longrightarrow \D_{\F}^{\mathrm{fil}}-\fdg.$$

On the underlying $\OO$-modules, the $\s$-functor $f^!$ acts as the usual pull-back, 
as shown by the following lemma.

\begin{lem}\label{ldIV-1}
With the notations and assumptions above, we have a natural equivalence of filtered $\OO_X$-modules
$$\D_f^{\mathrm{fil}} \simeq \OO_X\otimes_{g^{-1}(\OO_Y)}g^{-1}(\D_{\G}^{\mathrm{fil}}).$$
In particular, the underlying $\OO_X$-module of $f^!(E)$ is naturally equivalent to 
$f^*(E)$ and is therefore quasi-coherent. As a consequence we have a commutative square
of $\s$-functors 
$$\xymatrix{
\QCoh^{\mathrm{fil}}(\G) \ar[r]^-{f^!} \ar[d] & \ar[d]  \QCoh^{\mathrm{fil}}(\F) \\
\QCoh^{\mathrm{fil}}(Y) \ar[r]_-{g^* } & \QCoh^{\mathrm{fil}}(X),
}$$
where the vertical $\s$-functors are the forgetful $\s$-functors to filtered quasi-coherent
complexes.
\end{lem}

\textit{Proof of the lemma.}  This computations have already been seen in the lemma \ref{lpdmod}, 
because $\G$ is perfect we have an equivalence of sheaves of graded mixed complexes
$$\rch_{g^{-1}(\DR_\G)-\egrdg}(g^{-1}(\OO_Y),\OO_X) \simeq 
\rch_{g^{-1}(\DR_\G)-\egrdg}(g^{-1}(\OO_Y),g^{-1}(\OO_Y)) \otimes_{g^{-1}(\OO_Y)} \OO_X.$$
Passing to the realization provides the wanted result
$$\D_f^{\mathrm{fil}} \simeq \OO_X\otimes_{g^{-1}(\OO_Y)}g^{-1}(\D_{\G}^{\mathrm{fil}}).$$
The rest of the lemma follows directly from the first assumption as the quasi-coherent
pull-back is given by the standard formula 
$$g^*(E) \simeq g^{-1}(E)\otimes_{g^{-1}(\OO_Y)}\OO_X.$$
\hfill $\Box$ \\

Acoording to the previous lemma we do have constructed the pull-back of filtered crystals
$$f^! : \QCoh^{\mathrm{fil}}(\G) \longrightarrow \QCoh^{\mathrm{fil}}(\F).$$

\begin{df}\label{dIV-3}
The $\s$-functor constructed above is called the \emph{filtered pull-back of filtered crystals}.
\end{df}

A first important property is the existence of a left adjoint to the $!$-pull back, given by 
the $!$-pushfoward. It will exists under sufficient properness and quasi-smoothness 
conditions to insure
that the standard direct image of quasi-coherent sheaves preserves perfect complexes.

\begin{prop}\label{ldirect}
Let $f = (g,u) : (X,\F) \to (Y,\G)$ be a morphism between derived DM-stacks endowed 
with perfect derived foliations. We assume that $g : X \to Y$ is proper and quasi-smooth, 
and that
$Y$ is either a quasi-compact and quasi-separated derived scheme, or a separated
and quasi-compact derived DM-stack. 
Then, the $\s$-functor
$$f^! : \QCoh^{\mathrm{fil}}(\G) \longrightarrow \QCoh^{\mathrm{fil}}(\F)$$
admits a left adjoint 
$$f_! : \QCoh^{\mathrm{fil}}(\F) \longrightarrow \QCoh^{\mathrm{fil}}(\G).$$
\end{prop}

\textit{Proof.} We already know from Corollary \ref{ccompact} that 
both $\s$-categories are compactly generated. For the existence of 
$f_!$ it is thus enough to check that 
$f^!$ commutes with limits. For this, we use that $f^!$ is compatible
with the usual pull-backs of filtered $\OO$-modules along $g$ (see Lemma \ref{ldIV-1})
$$\xymatrix{
\QCoh^{\mathrm{fil}}(\G) \ar[r] \ar[d]_-{f^!} & \QCoh^{\mathrm{fil}}(Y) \ar[d]^-{(g^*)^{\mathrm{fil}}} \\
\QCoh^{\mathrm{fil}}(\F) \ar[r] & \QCoh^{\mathrm{fil}}(X).}$$
The horizontal $\s$-functors are clearly conservative and commute with 
limits and colimits. Therefore, to check that $f^!$ 
preserves limits, it is enough to show that $(g^*)^{\mathrm{fil}}$ does so. Again, 
as $(g^*)^{\mathrm{fil}}$ is obtained by extending the 
usual pull-back $g^* : \QCoh(Y) \longrightarrow \QCoh(X)$ to filtered objects, we are 
reduced to show that $g^*$ commutes with limits.
 
This last step follows easily from the assumption that $g$ is proper and quasi-smooth. Indeed, 
let $E \in \QCoh(X)$ be a compact generator (thus a perfect 
complex on $X$), and $\{F_i\}_{i\in I}$ a diagram in $\QCoh(Y)$. 
Using that $E$ is dualizable and the projection formula, we have
$$Map(E,g^*(\lim_i F_i)) \simeq Map(\OO_X,E^\vee \otimes_{\OO_X} g^*(\lim_i F_i)) \simeq
Map(\OO_Y,g_*(E^\vee)\otimes_{\OO_Y} (\lim_i F_i)).$$
Now we use that $g$ is proper and quasi-smooth, so that $g_*(E^\vee)$ is again perfect and thus dualizable, 
and therefore the functor $g_*(E^\vee)\otimes_{\OO_Y} -$ commutes with limits. So we have
$$Map(\OO_Y,g_*(E^\vee)\otimes_{\OO_Y} (\lim_i F_i)) \simeq \lim_i Map(\OO_Y,g_*(E^\vee)\otimes_{\OO_Y} F_i) \simeq
\lim_i Map(E,g^*(F_i)).$$
This shows that the canonical
map
$$Map(E,g^*(\lim_i F_i)) \longrightarrow \lim_i Map(E,g^*(F_i))\simeq Map(E,\lim_i \, g^*(F_i))$$
is indeed in equivalence, so that $g^*$ preserves limits (since $E$ is a generator of $\QCoh(Y)$). \hfill $\Box$ \\

The proposition \ref{ldirect} allows us to give the following definition.

\begin{df}\label{ddirect}
Let $f=(g,u) : (X,\F) \longrightarrow (Y,\G)$ be a morphism with 
$g : X \to Y$ proper and quasi-smooth and satisfying the condition of
the proposition \ref{ldirect}. The \emph{filtered direct image} is the $\s$-functor 
$$f_! : \QCoh^{\mathrm{fil}}(\F) \longrightarrow \QCoh^{\mathrm{fil}}(\G),$$
left adjoint to the pull-back $\s$-functor $f^!$.
\end{df}

Note that for $f=(\mathrm{id}_X, u : 0_X \to \F)$, we clear have $f_!= \mathsf{Ind}^{\mathrm{fil}}_\F$ 
(see Definition \ref{dinduc}).

\begin{rmk}
\emph{Note that, although it is possible, we do not try to define direct images 
for non-proper or non-quasi-smooth morphisms (this would require working with Ind-coherent crystals rather than quasi-coherent ones). In the rest of the book, we will only need direct images for proper and quasi-smooth morphisms.}
\end{rmk}

We also define the unfiltered direct image
$$f_! : \QCoh(\F) \longrightarrow \QCoh(\G)$$
as being the left adjoint to the unfiltered version of $f^!$ whose existence is proved similarly. \\

The direct image functors satisfy the usual pseudo-functoriality properties, $(ff')_!\simeq f_!f_!'$ (for composeable, proper and quasi-smooth $f$ and $f'$). An important 
consequence of this property is the following result.

\begin{prop}\label{cdirectinduct}
Let $f=(g,u) : (X,\F) \longrightarrow (Y,\G)$ be a morphism with $g$ proper and quasi-smooth, 
$\F$ and $\G$ perfect, and $Y$ either a qcqs derived scheme or a separated and quasi-compact 
derived DM-stack. 
Then, 
the following diagram naturally commutes
$$\xymatrix{
\QCoh^{\mathrm{fil}}(X) \ar[r]^-{\mathsf{Ind}^{\mathrm{fil}}_\F} \ar[d]_-{g_*(\omega_{X/Y}\otimes -)[d]} & 
\QCoh^{\mathrm{fil}}(\F) \ar[d]^-{f_!} \\
\QCoh^{\mathrm{fil}}(Y) \ar[r]_-{\mathsf{Ind}^{\mathrm{fil}}_\G} & \QCoh^{\mathrm{fil}}(\G),
}$$
where $\omega_{X/Y}$ is the relative canonical line bundle of $X$ over $Y$, and 
$d$ the relative dimension of $X$ over $Y$.
\end{prop}

\textit{Proof.} We have a commutative diagram of pairs
$$\xymatrix{
(X,0_X) \ar[r]^-{u} \ar[d]_-{g} & (X,\F) \ar[d]^-{f} \\
(Y,0_Y) \ar[r]_-{v} & (Y,\G)
}$$
where $u$ and $v$ are the unique morphisms from the initial foliation. We get from this a natural isomorphism
of $\s$-functors 
$$f_!u_!\simeq v_!g_!.$$
The proposition then follows from the fact that, by definition $u_!$ and $v_!$ are the induction $\s$-functors, 
and from the explicit formula for $g_! : \QCoh(X) \longrightarrow \QCoh(Y)$
in terms of relative Serre duality. \hfill $\Box$ \\

An interesting application of proposition \ref{cdirectinduct} is to the direct image of $\D_\F^{\mathrm{fil}}$ itself. Keeping the same notations we obtain 
$$f_!(\D_\F^{\mathrm{fil}}) \simeq \D_\G^{\mathrm{fil}}\otimes_{\OO_Y}g_*(\omega_{X/Y})[d].$$

When $X$ is proper and quasi-smooth over $k$,
and $f : (X,\F) \longrightarrow (\mathrm{Spec}\, k,0)$ is the projection to the point (endowed with its trivial
foliation), 
we see in particular that $f_!(\D_\F)$ computes $H^{*+d}(X,\omega_X)$, that is \emph{coherent homology 
of $X$ with coefficients in $\OO_X$.} \\

The following result is a direct consequence of the fact that 
$f^!$ commutes with colimits (and can also be deduced  
by  using proposition \ref{cdirectinduct}
and the proof of corollary \ref{ccompact}). 

\begin{cor}\label{cproper}
Let $f=(g,u) : (X,\F) \to (Y,\G)$ be a morphism of derived schemes endowed with perfect derived
foliations, with $g$ proper and quasi-smooth and $Y$ being either
a qcqs derived scheme or a separated quasi-compact derived DM-stack. 
Then, the direct image $\s$-functors
$$f_! : \QCoh^{\mathrm{fil}}(\F) \longrightarrow \QCoh^{\mathrm{fil}}(\G) \qquad 
f_! : \QCoh(\F) \longrightarrow \QCoh(\G)$$
preserve compact objects.
\end{cor}

\subsection{Compatibility with underlying and associated graded objects} We conclude this section with the important result that 
filtered direct images commute with both taking underlying and associated graded objects.
For this, let us consider a morphism $f=(g,u) : (X,\F) \longrightarrow (Y,\G)$ 
with $g$ proper and quasi-smooth, and assume the usual assumptions of proposition \ref{ldirect}. We have the corresponding  adjunction on filtered crystals
$$f^{\mathrm{fil}}_! : \QCoh^{\mathrm{fil}}(\F) \rightleftarrows \QCoh^{\mathrm{fil}}(\G) : f_{\mathrm{fil}}^! \, ,$$
and its unfiltered version
$$f_! : \QCoh(\F) \rightleftarrows \QCoh(\G) : f^!.$$
We may also consider $\F^{\epsilon=0}$ which is the derived foliation 
whose underlying graded mixed cdga is $\DR(\F)^{\epsilon=0}$ (i.e. with trivial mixed structure $\epsilon=0$), 
and the same for $\G^{\epsilon=0}$. For such derived foliations, we have now a graded push-forward
$$f^{\mathrm{gr}}_! : \QCoh^{\mathrm{gr}}(\F^{\epsilon=0}) \longrightarrow \QCoh^{\mathrm{gr}}(\G^{\epsilon=0})$$
defined as the left adjoint to the graded pull-back $f^!_{\mathrm{gr}} : \QCoh^{\mathrm{gr}}(\G^{\epsilon=0}) \rightarrow 
\QCoh^{\mathrm{gr}}(\F^{\epsilon=0})$.
Note that the $\s$-categories $\QCoh^{\mathrm{gr}}(\F^{\epsilon=0})$ and $\QCoh^{\mathrm{gr}}(\G^{\epsilon=0})$
can also be identified with the categories of graded $Sym_{\OO_X}(\T_\F)$-modules over $X$ and of 
graded $Sym_{\OO_Y}(\T_\G)$-modules over $Y$, respectively (where $\T_\F$ and $\T_\G$ both sit in weight $1$). We leave to the reader the precise construction of these graded functorialities, 
which are obtained using similarly as the filtered case.

By putting all these functors together, we may write the following diagram of vertical adjunctions
\begin{equation}\label{compGr}
\xymatrix{
\QCoh^{\mathrm{gr}}(\F^{\epsilon=0}) \ar@<-1.0ex>[d]_-{f^{\mathrm{gr}}_!} & \QCoh^{\mathrm{fil}}(\F) \ar[l]_-{\mathrm{Gr}} \ar[r]^-{(-)^u}
\ar@<-1.0ex>[d]_-{f^{\mathrm{fil}}_!} & \QCoh(\F) \ar@<-1.0ex>[d]_-{f_!} \\
\QCoh^{\mathrm{gr}}(\G^{\epsilon=0}) \ar@<-1.0ex>[u]_-{f_{\mathrm{gr}}^!} & 
\QCoh^{\mathrm{fil}}(\G) \ar[l]^-{\mathrm{Gr}} \ar[r]_-{(-)^u}
\ar@<-1.0ex>[u]_-{f_{\mathrm{fil}}^!} & \QCoh(\G) \ar@<-1.0ex>[u]_-{f^!}
}
\end{equation}

We already know that this diagram, when restricted to the inverse images only, naturally commutes. This
implies that the diagram restricted to direct images is naturally lax commutative. In fact, the induced natural
transformations
$$f_!^{\mathrm{gr}}\circ \mathrm{Gr} \Rightarrow \mathrm{Gr}\circ f^{\mathrm{fil}}_! 
\qquad
f_!\circ (-)^u \Rightarrow (-)^u\circ f_!^{\mathrm{fil}}$$
turn out to be equivalences. Indeed, as all $\s$-functors involved commute with
colimits, it is enough to check this property on compact generators, and  we can thus use
proposition \ref{cdirectinduct}
and the proof of corollary \ref{ccompact} to conclude. As a consequence of the commutativity of (\ref{compGr}), we are allowed to (and will from now on) simply
write $f_!$ and $f^!$ \emph{without} any decorations $(-)^{\mathrm{fil}}$ or $(-)^{\mathrm{gr}}$. 
Because of its importance, and for later reference, we state this result in the following corollary.

\begin{cor}\label{ccommute}
With the notations and assumptions above, the  
filtered direct and filtered inverse images of quasi-coherent crystals commute with 
taking the underlying or the associated graded objects, i.e. either considering its $f^!$ part or its $f_!$ part, the following diagram
naturally commutes
\begin{equation}
\xymatrix{
\QCoh^{\mathrm{gr}}(\F^{\epsilon=0}) \ar@<-1.0ex>[d]_-{f_!} & \QCoh^{\mathrm{fil}}(\F) \ar[l]_-{\mathrm{Gr}} 
\ar[r]^-{(-)^u}
\ar@<-1.0ex>[d]_-{f_!} & \QCoh(\F) \ar@<-1.0ex>[d]_-{f_!} \\
\QCoh^{\mathrm{gr}}(\G^{\epsilon=0}) \ar@<-1.0ex>[u]_-{f^!} & 
\QCoh^{\mathrm{fil}}(\G) \ar[l]^-{\mathrm{Gr}} \ar[r]_-{(-)^u}
\ar@<-1.0ex>[u]_-{f^!} & \QCoh(\G) \ar@<-1.0ex>[u]_-{f^!}
}
\end{equation}
\end{cor}

Here is an equivalent way of looking at the leftmost adjunction $f^{\mathrm{gr}}_! 
:\QCoh^{\mathrm{gr}}(\F^{\epsilon=0}) \leftrightarrows \QCoh^{\mathrm{gr}}(\G^{\epsilon=0}) : 
f^{\mathrm{gr}\, !}$ of (\ref{compGr}) that avoids introducing the auxiliary foliations 
$\F^{\epsilon=0}$ and $\G^{\epsilon=0}$. First of all, we have the associated graded object $\s$-functor
$$\mathrm{Gr}: \QCoh^{\mathrm{fil}}(\F)= \D_\F^{\mathrm{fil}}-\fdg_{qcoh} \longrightarrow 
Gr(\D_\F^{\mathrm{fil}})-\grdg_{qcoh},$$ 
from $\OO_X$-quasi-coherent sheaves of filtered $\D_\F^{\mathrm{fil}}$-dg-modules
on $X$ to $\OO_X$-quasi-coherent sheaves of graded $Gr(\D_\F^{\mathrm{fil}})$-dg-modules
on $X$ (and similarly for $\G$ and $Y$). Note that  $\mathrm{Gr}(\D_\F^{\mathrm{fil}}) \simeq
Sym_{\OO_X}(\mathbb{T}_\F)$ and $\mathrm{Gr}(\D_\G^{\mathrm{fil}}) \simeq 
Sym_{\OO_Y}(\mathbb{T}_\G)$ in the $\s$-category of sheaves of graded associative dg-algebras $\mathrm{Alg}(\grdg(X))$ and  $\mathrm{Alg}(\grdg(Y))$
(with $\mathbb{T}_\F=\LL_\F^\vee$ and $\mathbb{T}_\G=\LL_\G^\vee$ of pure weight $1$). 
Now, we proceed as in \S \ref{s-dirim} and consider
$$\mathrm{Gr}(\D_f^{\mathrm{fil}})\simeq g^*(\mathrm{Gr}(\D_\G^{\mathrm{fil}}))=\OO_X 
\otimes_{g^{-1}(\OO_Y)}g^{-1}(\mathrm{Gr}(\D_\G^{\mathrm{fil}})).$$ 
which is a graded $(\mathrm{Gr}(\D_\F^{\mathrm{fil}}),g^{-1}(\D_\G^{\mathrm{fil}}))$-bi-module over 
$X$, and as such, it defines a pull-back
$\s$-functor on graded modules
$$f_{\mathrm{gr}}^! : Gr(\D_{\G}^{\mathrm{fil}})-\grdg 
\longrightarrow Gr(\D_\F^{\mathrm{fil}})-\grdg$$
by  
$$f_{\mathrm{gr}}^!(E):=\mathrm{Gr}(\D_f^{\mathrm{fil}})
\otimes_{g^{-1}(\mathrm{Gr}(\D_\G^{\mathrm{fil}}))}g^{-1}(E).$$
This $\s$-functor respects the property of being quasi-coherent over $X$ and $Y$, so it induces a 
$\s$-functor
$$f_{\mathrm{gr}}^! : Gr(\D_\G^{\mathrm{fil}})-\grdg_{qcoh} \longrightarrow 
Gr(\D_\F^{\mathrm{fil}})-\grdg_{qcoh}$$ 
which has a left adjoint 
$$f^{\mathrm{gr}}_!: Gr(\D_\F^{\mathrm{fil}})-\grdg_{qcoh} \leftrightarrows  Gr(\D_\G^{\mathrm{fil}})-\grdg_{qcoh} : f_{\mathrm{gr}}^!,$$ 
that coincides with the leftmost adjunction  of (\ref{compGr}), via the natural equivalenes 
$\QCoh^{\mathrm{gr}}(\F^{\epsilon=0}) \simeq 
Gr(\D_\F^{\mathrm{fil}})-\grdg_{qcoh}$, 
and $\QCoh^{\mathrm{gr}}(\G^{\epsilon=0}) \simeq 
Gr(\D_\G^{\mathrm{fil}})-\grdg_{qcoh}$.

\section{Characteristic cycles}

In this section we introduce the notion of \emph{characteristic cycle} of
a quasi-coherent crystal (or of a $\D_\F$-module thanks to the equivalence
of Theorem \ref{pdmod}) along a derived foliation $\F$. For this we first introduce the global cotangent stack $T^*\F$ of a derived
foliation $\F$, which is a derived Artin $n$-stack, where $n$ is the tor-amplitude of the perfect complex
$\LL_\F$. We will then discuss
the notion of \emph{bounded
coherent crystals} and of \emph{good filtrations} on them. By definition, the associated graded 
to a good filtration will be a $\Gm$-equivariant perfect complex on $T^*\F$, that will be used to define 
 \emph{characteristic cycles}. We investigate the 
existence of good filtrations and prove some independence (of good filtrations) results for 
characteristic cycles.

\subsection{Cotangent stacks of derived foliations}\label{subsec-ctgtofderivedfol}

Let $X$ be a derived scheme and $E$ be a perfect complex on $X$ of tor-amplitude 
contained in $[0,n]$ for some non-negative integer $n$ (see
\cite{tova} for the notion of tor-amplitude of perfect complexes). 
Recall from \S \ref{subsec:Derivedlinearstacks} that the linear stack $\VV(E)$
associated to $E$ comes equipped with its natural 
$\Gm$-action, covering the projection $\pi : \VV(E) \to X$. We pass to
the corresponding morphism on quotient stacks
$$\pi^{\mathrm{gr}} : [\VV(E)/\Gm] \longrightarrow X \times B\Gm.$$
The direct image along this morphism is a symmetric lax monoidal $\s$-functor
\begin{equation}\label{for-this}\pi_*^{\mathrm{gr}} : \QCoh^{\mathrm{gr}}(\VV(E)):=\QCoh([\VV(E)/\Gm]) \longrightarrow \QCoh(X \times B\Gm)=:\QCoh^{\mathrm{gr}}(X).\end{equation}
Since the structure sheaf $\OO:= \OO_{[\VV(E)/\Gm]}$ is the monoidal unit in $\QCoh^{\mathrm{gr}}(\VV(E))$, the lax-monoidal the $\s$-functor $\pi_*^{\mathrm{gr}}$ in (\ref{for-this}) factors via a $\s$-functor (denoted by the same symbol)
$$\pi_*^{\mathrm{gr}} : \QCoh^{\mathrm{gr}}(\VV(E)) \longrightarrow \pi^{\mathrm{gr}}_*(\OO)-\grdg_{qcoh},$$
from graded quasi-coherent complexes on $\VV(E)$ to 
graded $\pi^{\mathrm{gr}}_*(\OO)$-modules over $X$ which are quasi-coherent as $\OO_X$-modules.

\begin{prop}\label{plinear}
There exists a fully faithful $\s$-functor
$$Sym_{\OO_X}(E)-\grdg_{\mathrm{perf}} \longrightarrow \QCoh^{\mathrm{gr}}(\VV(E))$$
from 
perfect graded $Sym_{\OO_X}(E)$-modules over $X$, sending
the i-th weight twist $Sym_{\OO_X}(E)(i)$ of the tautological graded module $Sym_{\OO_X}(E)$, to the i-th weight twist  $\OO_{\VV(E)}(i)$ of the structure sheaf of $\VV(E)$.
\end{prop}

\textit{Proof.} We start by assuming that $X=Spec\, A$ is affine, with $A$ a
connective cdga, and by \cite{mon}
we known that there is a natural equivalence graded cdga's over $X$
$$\pi^{\mathrm{gr}}_*(\OO)\simeq Sym_{\OO_X}(E).$$
The direct image $\s$-functor
$$\pi_*^{\mathrm{gr}} : \QCoh^{\mathrm{gr}}(\VV(E)) \longrightarrow \pi^{\mathrm{gr}}_*(\OO)-\grdg_{qcoh,X},$$
from graded quasi-coherent complexes on $\VV(E)$ to 
graded $\pi_*(\OO)$-modules over $X$ which are quasi-coherent as $\OO_X$-modules possesses a left
adjoint $\pi^*_{gr}$
$$\pi^*_{gr} : \pi^{\mathrm{gr}}_*(\OO)-\grdg_{qcoh} \simeq Sym_{\OO_X}(E)-\grdg_{perf} \longrightarrow 
\QCoh^{\mathrm{gr}}(\VV(E)).$$
As the unit of this adjunction
$$Sym_{\OO_{X}}(E) \to \pi_*^{gr}(\pi_{gr}^*(Sym_{\OO_X}(E)) = \pi_*^{gr}(\OO)$$
is an equivalence, the $\s$-functor $\pi^*_{gr}$, when 
restricted to the thick triangulated sub-$\s$-category generated by the objects
$\OO(i)$, induces an equivalences with perfect graded $Sym_{\OO_X}(E)$-modules over $X$. This
proves the proposition when $X$ is affine.

The case where $X$ is a general derived Artin stack easily follows by descent from the affine case
thanks to the commutative diagram 
$$\xymatrix{Sym_{\OO_X}(E)-\grdg_{perf} \ar[r] \ar[d]_-{\simeq} & \QCoh^{\mathrm{gr}}(\VV(E)) 
\ar[d]^-{\simeq} \\
\lim_{u : Spec\, A \to X}\left(Sym_{A}(u^*(E))-\grdg_{perf}\right) \ar[r] & \lim_{u : Spec\, A \to X}\QCoh^{\mathrm{rm}}(\VV(u^*(E))).
}$$
\hfill $\Box$ \\

The proposition \ref{plinear} will be applied later to the particular
case where $E=\T_\F$, the tangent complex of a derived foliation $\F$ on $X$. 
The derived stack $V(\T_\F)$ will be denoted by $T^*\F$, and will be called
the \emph{global derived cotangent stack} of $\F$.

\begin{df}\label{dcotangent}
For a derived scheme $X$ and a perfect derived foliation $\F \in \Fol^p(X)$, 
the \emph{(derived) cotangent stack 
of $\F$} is defined by 
$$T^*\F:=V(\T_\F).$$
equipped with is natural $\Gm$-action.
\end{df}

In the next section, we will use the functor of the proposition \ref{plinear} in order to define \emph{characteristic cycles} of crystals on $\F$ admitting a good filtration.

\subsection{Coherent crystals and good filtrations}

Let $X$ be a separated and quasi-compact derived Deligne-Mumford stack 
and $\F \in \Fol^p(X)$ be a perfect derived foliation on $X$. We let 
$E \in \QCoh(\F)$ be a quasi-coherent crystal. Via the equivalence of our theorem 
\ref{pdmod}, we will freely identify 
$E$ with a $\D_\F$-module quasi-coherent over $X$.

\begin{df}\label{dgood}
Let $X$, $\F$ and $E$ as above. 
\begin{enumerate}
\item The crystal $E$ is called \emph{coherent} if it is a compact object
in $\QCoh(\F)$. The full sub-$\s$-category of coherent crystals is denoted by 
$\Coh(\F) \subset \QCoh(\F)$.

\item For $E \in \Coh(\F)$, a \emph{good filtration on $E$} is the data of
a compact object $E^{\mathrm{fil}} \in \QCoh^{\mathrm{fil}}(\F)$ together with an equivalence in $\QCoh(\F)$
$$(E^{\mathrm{fil}})^u \simeq E.$$
\end{enumerate}
\end{df}

Compact objects in $\QCoh(\F)$ and in $\QCoh^{\mathrm{fil}}(\F)$ can be easily characterized, either
through the induction $\s$-functors or as sheaves of $\D_\F$-modules.

\begin{prop}\label{pcompact}
Let $X$ and $\F$ be as above. An object $E \in \QCoh(\F)$ (resp.  
$E^{\mathrm{fil}} \in \QCoh^{\mathrm{fil}}(\F)$) is compact if and only if it
satisfies one of the following two equivalent conditions.

\begin{enumerate}
\item The object $E$ belongs to the thick triangulated sub-category generated
by objects of the form $\mathsf{Ind}_{\F}(E_0)$ (resp. $\mathsf{Ind}^{\mathrm{fil}}_{\F}(E^{\mathrm{fil}}_0)$)
for a perfect complex of $\OO_X$-modules $E_0$ (resp. a filtered perfect
complex of $\OO_X$-modules $E_0^{\mathrm{fil}}$).

\item $E$ (resp. $E^{\mathrm{fil}}$) is a perfect $\D_\F$-module (resp. a perfect $\D_\F^{\mathrm{fil}}$-module) i.e. locally on $X_{\text{ét}}$, $E$ (resp. $E^{\mathrm{fil}}$) is
a retract of a finite cell $\D_\F$-module (resp. of  a finite cell filtered $\D_\F^{\mathrm{fil}}$-module).

\end{enumerate}
\end{prop}

\textit{Proof.} Condition $(1)$ have been already considered in the proof of Corollary \ref{ccompact}. 
Clearly, condition $(1)$ implies condition $(2)$. Finally condition $(2)$ clearly implies 
compactness when $X$ is affine. The general case follows from
the fact that $X$ is assumed to be separated and quasi-compact and therefore 
\ref{ccompact2} holds and from the fact that any object satisfying condition (2)
is a compact object (see our remark \ref{rccompact2}). \hfill $\Box$ \\

Suppose that $E$ is a coherent crystal along $\F$, equipped with a good filtration 
$E^{\mathrm{fil}}$ in the sense of definition \ref{dgood} above. By proposition \ref{pcompact}, 
$\mathrm{Gr}(E^{\mathrm{fil}})$ is a perfect graded $\mathrm{Gr}(\D_\F^{\mathrm{fil}})=Sym_{\OO_X}(\T_\F)$-module over $X$. 
Now, proposition \ref{plinear} implies that $\mathrm{Gr}(E^{\mathrm{fil}})$ defines a graded perfect 
complex on the stack $T^*\F$, and by forgetting the $\Gm$-action 
we get a perfect complex on $T^*\F$. We consider the $K$-theory spectrum
$K(T^*\F)$, defined as the $K$-theory of the $\s$-category of perfect
(not graded) modules over $T^*\F$. The perfect complex $\mathrm{Gr}(E^{\mathrm{fil}})$
on $T^*\F$ thus defines a class
$$Ch(E^{\mathrm{fil}}) \in K_0(T^*\F).$$

\begin{df}\label{dsing1}
The \emph{characteristic cycle} of $E^{\mathrm{fil}}$ is the element
$$Ch(E^{\mathrm{fil}}) \in K_0(T^*\F)$$
defined above.
\end{df}

Note that the above definition depends a priori on the choice of $E^{\mathrm{fil}}$. We will see later
that, in fact, modulo a class of \emph{phantoms objects}, it does not (see \S \ref{sect-inde}, and 
proposition \ref{pinv}). When $X$ and the derived foliation $\F$ are both smooth these phantoms
are all trivial, but in general we have not been able to show that $Ch(E^{fil})$ only depends on $E$ 
and not of the choice of the good filtration.\\

\subsubsection{Existence of good filtrations.} The existence of good filtrations in general
seems to be a complicated question, and the authors do not know if good filtrations always
exist for arbitrary coherent crystals, as it is the case for usual $\D$-modules on smooth varieties. 
It can be shown that they do exist for \emph{smooth} foliations on \emph{smooth} varieties, but the 
fact that, for  a general foliation $\F$, its ring of differential operators $\D_\F$ is not concentrated
in degree $0$, creates complications in constructing good filtrations. The following result is 
therefore very 
useful in practice.

\begin{prop}\label{pgooddirect}
Let $f=(g,u) : (X,\F) \longrightarrow (Y,\G)$ be a morphism of separated and quasi-compact derivbed
DM-stacks 
endowed with perfect derived foliations and assume that $g$ is proper and quasi-smooth. If 
$E^{\mathrm{fil}}$ is a good filtration
on a coherent crystal $E \in \QCoh(\F)$, then $f_!(E^{\mathrm{fil}})$ is a good
filtration on $f_!(E)$.
\end{prop}

\textit{Proof.} This follows from Corollary \ref{cproper} together with the fact that 
direct images commutes with taking the underlying object (see diagram (\ref{dunno})). \hfill $\Box$ \\

We can also single out a large class of coherent crystals for which existence
of good filtrations is guaranteed: the \emph{finite cell} crystals. For
a derived scheme $X$ and $\F \in \Fol(X)$, we define a non-thick 
triangulated sub-$\s$-category $\Coh^{cell}(\F) \subset \Coh(\F)$, as being generated
(by finite limits and shifts)
by the objects of the form $Ind_\F(E)$ for $E$ a perfect complex on $X$. An object 
of $\Coh^{cell}(\F)$ will be called  a \emph{finite cell crystal}. There is an obvious filtered
version, too.
Finite cell crystals are
obviously coherent, but the converse has no reason to be true in general. Finite cell crystals are
however useful because of the following result.

\begin{prop}\label{pcell}
With the notation above, any object $E \in \Coh^{cell}(\F)$ admits 
a good filtration. Moreover, a good filtration $E^{\mathrm{fil}}$ can be chosen to
be a finite cell object in $\QCoh^{\mathrm{fil}}(\F)$.
\end{prop}

\textit{Proof.} Let $E \in \Coh^{cell}(\F)$. By definition of a finite cell object, there is a finite sequence of morphisms
in $\Coh(\F)$
$$\xymatrix{
E_{-1}=0 \ar[r] & E_0 \ar[r] & \dots \ar[r] & E_i \ar[r] & E_{i+1} \ar[r] & \dots \ar[r] & 
E_n=E,}$$
with the following property: for all $i$ there is a perfect complex $K_i$ on 
$X$ and a cartesian square in $\Coh(\F)$
$$\xymatrix{
E_i \ar[r] & E_{i+1} \\
Ind_\F(K_i) \ar[u]^-{u_i} \ar[r] & 0.  \ar[u]}$$
We will show, by induction, that each $E_i$ has a good filtration. For this, assume
that $E_i$ has a good filtration $E_i^{\mathrm{fil}}$. The morphism $u_i$ is given, by adjunction, by a morphism 
$v_i : K_i \longrightarrow E_i=(E_i^{\mathrm{fil}})^u$ in $\QCoh(\F)$.
The quasi-coherent sheaf $(E_i^{\mathrm{fil}})^u$ is the filtered colimit $\mathrm{colim}_k F^k(E_i^{\mathrm{fil}})$, and
as $K_i$ is a compact object in $\QCoh(X)$, $v_i$ can be factored as
$\xymatrix{
K_i \ar[r]^-{w_i} & F^k(E_i^{\mathrm{fil}}) \ar[r] & E_i
},$
for some index $k$. By using the left adjoint
$Ind_\F^{\mathrm{fil}}$, the morphism $w_i$ corresponds to a morphism of filtered $\D_\F^{\mathrm{fil}}$-modules
$$\alpha_i :  Ind_{\F}^{\mathrm{fil}}(K_i)(-k) \longrightarrow E_i,$$
where $(-k)$ denotes the endofunctor of $\QCoh^{\mathrm{fil}}(\F)$ that shifts by $-k$ the filtration. The cone 
of $\alpha_i$ clearly defines a good filtration on $E_{i+1}$.
\hfill $\Box$ \\

\subsubsection{Independence of the good filtration.}\label{sect-inde} We already noticed that the characteristic 
cycle of definition 
\ref{dsing1} \emph{depends}, a priori, on the good filtration $E^{\mathrm{fil}}$. In order to solve this
problem, we introduce a reduced $K$-group, and prove that the image of $Ch(E^{\mathrm{fil}})$ in this reduced $K$-group, only
depends on the object $E$. \\

Let $\F$ be a derived foliation on a derived scheme $X$. 
We consider $\Coh^{\mathrm{fil}}(\F)$, the $\s$-category of compact objects $\QCoh^{\mathrm{fil}}(\F)$, and
the underlying object $\s$-functor $(-)^u : \Coh^{\mathrm{fil}}(\F) \longrightarrow \Coh(\F)$.
An object $E \in \Coh^{\mathrm{fil}}(\F)$ will be called a \emph{phantom} if 
$E^u\simeq 0$. We then give the following definition.

\begin{df}\label{reduced}
With the notations above, the reduced $K$-group $K^{red}_0(T^*\F)$ is the 
quotient of $K_0(T^*\F)$ by the subgroup generated by the classes of $Gr(E)$ for 
$E \in \Coh^{\mathrm{fil}}(\F)$ a phantom.
\end{df}

We have the following independence result.

\begin{prop}\label{pinv}
Let $\F$ be a derived foliation on a derived scheme $X$, and
$E \in \QCoh(\F)$ be a coherent crystal.
Let $E_1^{\mathrm{fil}}$ and $E_2^{\mathrm{fil}}$
be two good filtrations on $E$ in the sense of Definition \ref{dgood}.
Then we have $Ch(E_1^{\mathrm{fil}})=Ch(E_2^{\mathrm{fil}})$ in $K^{red}_0(T^*\F)$.
\end{prop}

\textit{Proof.} We start by the following lifting lemma.

\begin{lem}\label{llift}
Let $M,N \in \QCoh^{\mathrm{fil}}(\F)$ with $M$ compact. Then 
any morphism $u : M^u \longrightarrow 
N^u$ in $\QCoh(\F)$ can be lifted, via the $\s$-functor $(-)^u$, to a morphism
$v : M \longrightarrow N(k)$ in $\QCoh^{\mathrm{fil}}(\F)$, where $N(k)$ the $k$-weight shift of $N$, for some integer $k$.
\end{lem}

\textit{Proof of the lemma.} Recall that $N(k)$ denotes the filtered crystal 
with the filtration shifted by $k$, i.e.
$$F^i(N(k)):=F^{i+k}N.$$
As $M$ is compact, we know by the proof of Corollary \ref{ccompact}, that 
$M$ is a retract of a finite cell object in $\QCoh^{\mathrm{fil}}(\F)$ in the sense of \ref{pcell}. 
Clearly, if the lemma is true for $M$ (and any $N$) it is also true for any of its retracts. We 
may therefore assume that $M$ is a finite cell object. By induction on the number of cells, we 
reduce ourselves to the following statement. Assume that the lemma is true for $M$ (and any $N$), and let 
us consider a push-out
$$\xymatrix{
M \ar[r] & M' \\
Ind_{\F}(K) \ar[u]-^{\alpha} \ar[r] & 0 \ar[u]
}$$
with $K$ compact in $\QCoh(X)$.
We must prove that the lemma remains true for $M'$ (and any $N$). Let $v : (M')^u \longrightarrow
N^u$ be a morphism. It consists of the data of a morphism $u : M^u \longrightarrow N^u$
and a homotopy to zero of $(u\circ \alpha ) : Ind_{\F}(K) \longrightarrow N^u$, i.e., by adjunction, 
a homotopy to zero $h$ of the induced morphism $\beta : K \longrightarrow N^u$. 
By assumption on $u$ we can lift $u$ to $w : M \longrightarrow N(k)$ for some $k$. Moreover, 
 $h$ defines a homotopy to zero of $v=(w)^u : M^u \longrightarrow N(k)^u \simeq N^u$. 
As $K$ is compact and $N^u=\mathrm{colim}_k F^iN$, the homotopy $h$ factors as a homotopy to zero
$h'$ of $K \longrightarrow F^{k'}N \longrightarrow N^u$ for some $k'\geq k$. This pair $(w,h')$ then 
defines the required lift $M' \longrightarrow N(k')$. 
\hfill $\Box$ \\

Let us go back to the proof of Proposition \ref{pinv}. As $E_i^{\mathrm{fil}}$, $i=1, 2$ are a good filtrations on the same crystal $E$, we have a natural 
equivalence in $\QCoh(F)$
$$u : (E_1^{\mathrm{fil}})^u \simeq (E_2^{\mathrm{fil}})^u.$$
By Lemma \ref{llift}, $u$ can be lifted to a morphism of filtered crystals
$v : E_1^{\mathrm{fil}} \longrightarrow E_2^{\mathrm{fil}}$. We set 
$M$ to be the cone of $v$ in $\QCoh^{\mathrm{fil}}(\F)$. This is a compact object which is obviously a
phantom, i.e. $M^u\simeq 0$. We thus have an $Ch(E_1^{\mathrm{fil}})=Ch(E_2^{\mathrm{fil}})$ in $K^{red}_0(T^*\F)$.
\hfill $\Box$ \\

Proposition \ref{pinv} shows that the following notion is well defined.

\begin{df}
Let $X$ be a quasi-compact and separated derived DM-stack and $\F \in \Fol^p(X)$
a perfect derived foliation. If $E\in \Coh(\F)$ admits a good filtration $E^{\mathrm{fil}}$, 
then its characteristic cycle is defined as
$$Ch(E)=Ch(E^{\mathrm{fil}}) \in K^{red}_0(T^*\F).$$
\end{df}

It is possible to show that when $X$ and $\F$ are both smooth, then 
the natural projection $K_0(T^*\F) \longrightarrow K_0^{red}(T^*\F)$
is bijective, or, in other words, that the class of $Gr(E)$ is trivial
in $K_0(T^*\F)$ for any phantom $E \in \Coh^{\mathrm{fil}}(\F)$. This relies on 
using regularity and Quillen's devissage techniques, that do not work in 
our general setting. We will not need this isomorphism $K_0(T^*\F) \simeq K_0^{red}(T^*\F)$ in the rest of the book, so we omit its proof. 
However, the following very simple particular case will be useful later,
in order to get numerical formulas out of our general Grothendieck-Riemann-Roch statement (Section \ref{sec:GRR}).

\begin{prop}\label{ppoint}
Let $X$ be a quasi-compact and separated derived DM-stacks and $\F=0_X$ be the initial foliation,
so that $T^*\F\simeq X$. Then, the natural projection
$$K_0(X) \longrightarrow K_0^{red}(X)$$
is bijective.
\end{prop}

\textit{Proof.} In this case $\QCoh^{\mathrm{fil}}(\F)\simeq \QCoh^{\mathrm{fil}}(X)$
is the $\s$-category of filtered quasi-coherent complexes. The compact objects
in $\QCoh^{\mathrm{fil}}(X)$ are therefore finite filtered perfect complexes, and thus
$Gr(E)$ has only a finite number of non-trivial summand for any compact object
$E \in \QCoh^{\mathrm{fil}}(X)$. As a result, if $E$ is a phantom, we have
$0=[Gr(E)]$ in $K_0(X)$. \hfill $\Box$ \\

\begin{rmk}\label{rem:K0redofSpeck}
\emph{A particularly important case of proposition \ref{ppoint} is when $X=\mathrm{Spec}\, k$, for which
we find $K_0^{red}(\mathrm{Spec}\, k)\simeq \ZZ$.}
\end{rmk}

\section{Grothendieck-Riemann-Roch for derived foliations}\label{sec:GRR}

In this Section we will state and prove the \emph{Grothendieck-Riemann-Roch} (GRR) \emph{formula} for proper maps between
derived schemes endowed with derived foliations. 

\subsection{The GRR formula}

Let 
$f=(g,u) : (X,\F) \longrightarrow (Y,\G)$ be a morphism of derived stack endowed with  
derived foliaions, with $g$ proper and quasi-smooth. We will make the standard
assumption that $Y$ is either a qcqs derived scheme, or a separated and quasi-compact derived
$DM$-stack.
Associated to $f$ is the 
so-called ``Japanese correspondence''
$$\xymatrix{
T^*\F & T^*\G \times_Y X \ar[r]^-{p} \ar[l]_-{q} & T^*\G.}$$
Here $q$ is induced by the morphism of $g^*(\LL_\G) \to \LL_\F$ perfect complexes on $X$ 
induced by $u$ while the morphism $p$ is the first projection.

Define the push-forward on K-groups
$$f_!:=p_!(q^*(-)) : 
K_0(T^*\F) \longrightarrow K_0(T^*\G).$$
where $p_!$ is the left adjoint of the pull-back of quasi-coherent sheaves
$p^*  : \QCoh(T^*\G\times_Y X) \longrightarrow \QCoh(T^*\G)$. This left adjoint exists 
because $p$ is representable, proper and quasi-smooth, and it is given explicitly by 
$$p_!(E)=p_*(E\otimes \omega_p[d])$$ for $E\in \Parf(T^*\G\times_Y X)$, where the integer
$d$ is the relative dimension of $X$ over $Y$, and $\omega_p$ is the relative canonical sheaf of the 
morphism $p$. Note that $\omega_p$ coincides with the
pull-back of $\omega_{X/Y}$ along the projection $T^*\G\times_Y X \to X$. 

We start by noticing that $f_!$ is compatible with the quotient defining 
reduced $K$-groups. Indeed, if $E \in \Coh^{\mathrm{fil}}(\F)$ is a phantom, then so is $f_!(E)$, and
the image by $f_!$ of $Gr(E)$ is $Gr(f_!E)$
(see Corollary \ref{ccommute}). Therefore, $f_{!}$ induces a well defined map
$$f_! : K_0^{red}(T^*\F) \longrightarrow K_0^{red}(T^*\G).$$

\begin{thm}\label{tgrr}
Let $f=(g,u): (X,\F) \longrightarrow (Y,\G)$ be a morphism of 
derived DM-stack schemes endowed with perfect derived foliations with $g$ proper and quasi-smooth.
We assume that either $Y$ is separated and quasi-compact, or a qcqs derived scheme.
For any coherent crystal $E \in \Coh(\F)$ along $\F$ admitting a good filtration (e.g. 
a finite cell object by proposition \ref{pcell}), we have
$$Ch(f_!(E))=f_!(Ch(E))$$
in the group $K_0^{red}(T^*\G)$.
\end{thm}

\textit{Proof.} This is a direct consequence of the fact that direct images
commutes with taking associated graded for good filtrations (corollary \ref{ccommute}). \hfill $\Box$ \\

The following consequence of Theorem \ref{tgrr} is obtained when $Y=\mathrm{Spec}\, k$ and $\G=0_Y$. 
Let $s : X \to T^*\F$ denotes the zero section of the
cotangent stack of $\F$. When $X$ is a proper and quasi-smooth derived DM-stack, we denote by 
$p : X \to Spec\, k$ the canonical map, and by
$$p_*=:\int_X : K_0(X) \longrightarrow K_0(k)\simeq \ZZ$$
the induced push-forward on perfect complexes.

\begin{cor}\label{chrr}
Let $(X,\F)$ be a quasi-smooth and proper derived DM-stack endowed with a perfect 
derived foliation $\F\in \Fol^p(X)$. 
Let $f : (X,\F) \to (Spec\, k,0_{\mathrm{Spec}\, k})$ be the canonical morphism.
For any coherent crystal $E$ that admits a good filtration we have
$$\chi(f_!(E))=\int_{X}s^*(Ch(E))\otimes \omega_X.$$
\end{cor}

This corollary is a \emph{Hirzeburch-Riemann-Roch} (HRR) \emph{formula} for crystals along the foliation $\F$. The complex
$f_!(E)$ is what should be called the \emph{foliated cohomology of $(X,\F)$ with coefficients in $E$}.
If we denote this cohomology by $H^*_{\F}(X,E)$, the HRR formula reads
$$\chi(H^*_{\F}(X,E))=\int_{X}s^*(Ch(E))\otimes \omega_X.$$

\subsection{The non-proper case: Fredholm crystals}\label{s-ell}

We explain here briefly how to extend theorem \ref{tgrr} to \emph{non-proper} maps, 
by introducing the notion of \emph{Fredholm crystal}. \\

We start by extending direct images to the case of \emph{compactifiable morphisms}. 
Assume that $f=(g,u) : (X,\F) \to (Y,\G)$ is a morphism of derived DM-stack endowed with perfect 
derived foliations, such that 
$g: X \to Y$ is quasi-smooth. A \emph{quasi-smooth compactification} of $f$ is the datum of a factorization
\begin{equation}\label{fact}
f : \xymatrix{
(X,\F) \ar[r]^-{(j,v)} & (\bar{X},\bar{\F}) \ar[r]^-{\bar{f}} & (Y,\G)}
\end{equation}
where $j$ is an open embedding, $v$ is an equivalence $\F \simeq j^*(\bar{\F})$, 
and $\bar{f}$ is a proper and quasi-smooth morphism. For such a compactification, 
we set 
$$j_* : \QCoh(\F) \to \QCoh(\bar{\F})$$
to be the right adjoint to the $\s$-functor $j^!$. We note here that $j_*$ exists because
$j^!$ commutes with colimits. Moreover, $j_*$ is compatible with the usual push-forward of quasi-coherent
sheaves through the forgeful $\s$-functors, i.e. the following diagram naturally commutes
$$\xymatrix{
\QCoh(\F) \ar[r]^-{j_*} \ar[d] & \QCoh(\bar{\F}) \ar[d] \\
\QCoh(X) \ar[r]_-{j_*} & \QCoh(\bar{X}),
}$$
where the vertical $\s$-functors are the natural forgetful functors. Since 
$j$ is an open immersion and $v$ is an equivalence, we have a natural 
equivalence $j^{-1}\D_{\bar{\F}} \simeq \D_\F$ of sheaves of dg-algebras on $X$,
and a canonical adjunction 
morphism on $\bar{X}$ 
$$a : \D_{\bar{\F}} \longrightarrow j_*(\D_\F).$$
The $\s$-functor $j_* : \QCoh(\F) \to \QCoh(\bar{\F})$ is then simply induced by the usual push-forward of 
sheaves along $j : X \hookrightarrow \bar{X}$: it sends a
$\D_\F$-module $E$ to $j_*(E)$, viewed as a $\D_{\bar{\F}}$-module via the map
$a$ above. This description shows that we also have a commutative square 
involving the induction $\s$-functors
$$\xymatrix{
\QCoh(\F) \ar[r]^-{j_*}  & \QCoh(\bar{\F}) \ \\
\QCoh(X) \ar[r]_-{j_*} \ar[u]^-{Ind_{\F}} & \QCoh(\bar{X}). \ar[u]_-{Ind_{\bar{\F}}}
}$$

We define the \emph{$*$-direct image along $f$} as
$$f_*=\bar{f}_!\circ j_* : \QCoh(\F) \longrightarrow \QCoh(\G).$$
A standard argument, by using the product embedding, proves that $f_*$, defined as above, does not depend on the
choice of the factorization (\ref{fact}) (see e.g. \cite[Exp XVII]{SGA4.3} or \cite[\S 8]{FrKi}). 
Indeed, $f_*$ clearly commutes with colimits, so 
the independence of the factorization in its definition can be checked on compact $\D_\F$-modules
of the form $Ind_\F(E)$, where it is easily deduced from the compatibility of push-forwards
for quasi-coherent sheaves, and from the explicit formula of proposition \ref{cdirectinduct}. 

\begin{df}\label{dell}
Let $f=(g,u) : (X,\F) \to (Y,\G)$ be a morphism of derived DM-stack 
endowed with perfect derived foliations. Suppose that $f$ admits a quasi-smooth compactification. An object $E \in \Coh(\F)$ is called f-\emph{Fredholm}
if it admits a good filtration $E^{\mathrm{fil}} \in \Coh^{\mathrm{fil}}(\F)$ such that $j_*(E^{\mathrm{fil}})$ is a compact object in $\QCoh^{\mathrm{fil}}(\F)$ for some quasi-smooth compactification $j : (X,\F) \hookrightarrow
(\bar{X},\bar{\F})$ of the morphism $f$.
\end{df}

Suppose $E \in \Coh(\F)$ is $f$-Fredholm; pick a quasi-smooth compactification $(\bar{X},\bar{\F})$ 
and a filtration $E^{\mathrm{fil}}$ as in the definition above. 
We have the associated graded $Gr(E^{\mathrm{fil}}) \in \Parf(T^*\F)$, and its direct image by $j$ is by 
again a perfect complex on $T^*\bar{\F}$, so that, in particular, the support 
of $Gr(E^{\mathrm{fil}})$, $SS(E) \subset T^*\F$ is closed in $T^*\bar{\F}$. If we pull-back 
this support by the canonical map $T^*\G \times_Y X \to \T^*\F$,  
we get a closed subset in $T^*\G \times_Y X$ which is proper over $T^*\G$. 
In other words, the morphism $f$ is proper when restricted to the support of $Gr(E^{\mathrm{fil}})$. 
In our situation, the notion of Fredholm is a priori stronger, it is unclear that 
properness of $f$ on the support of $Gr(E^{\mathrm{fil}})$ is enough to recover the fact that 
$j_*(E^{\mathrm{fil}})$ remains a compact object. Therefore, our definition of being $f$-Fredholm might
be hard to check in practice, already because it involves the existence of 
a compactification of $f$ which seems a delicate question even when $X$ and $Y$ are both 
smooth quasi-projective varieties. \\

We can now state the Grothendieck-Riemann-Roch formula for possibly \emph{non-proper} maps, 
and \emph{Fredholm} coefficients. 

\begin{thm}\label{tgrrell}
Let $f=(g,u): (X,\F) \longrightarrow (Y,\G)$ be a morphism of 
derived DM-stacks endowed with perfect derived foliations. We assume as usual
that either $Y$ is a qcqs derived schemes, or a separated and quasi-compact DM-stack.
If $f$ admits a quasi-smooth
compactification, and $E \in \Coh(\F)$ is an $f$-Fredholm crystal,
then we have an equality 
$$Ch(f_*(E))=f_*(Ch(E))$$
of elements in $K_0^{red}(T^*\G)$.
\end{thm}

\textit{Proof.} Simply apply Theorem \ref{tgrr} to the morphism $\bar{f}$ and to the 
object $j_*(E)$ over $(\bar{X},\bar{\F})$,
which has a good filtration given by $j_*(E^{\mathrm{fil}})$ for $E^{\mathrm{fil}}$ given by 
Definition \ref{dell}. \hfill $\Box$ \\

\section{Examples and applications}\label{GRRapplications}

We finish this Chapter by giving a sample of examples and applications of the general Grothendieck-Riemann-Roch formula of Theorem \ref{tgrr}. \\

\subsection{Grothendieck-Riemann-Roch for derived $\D$-modules.}\label{s-grrdmod}

The very first special case of Theorem \ref{tgrr} is when $\F=*_X$ and $\G=*_Y$ are the
final derived foliations on $X$ and $Y$ respectively. The $\s$-categories
$\QCoh(*_X)$ and $\QCoh(*_Y)$ are then called the \emph{$\s$-categories of derived 
$\D$-modules on $X$ and $Y$}, and denoted by $\D_X-\dg_{qcoh}$ and $\D_Y-\dg_{qcoh}$, respectively. 
We think they coincide, for all $X$, with the derived $\D$-modules introduced and studied in 
\cite{beraldo_derD}, where they are denoted by $\D^{der}(X)$ and $\D^{der}(Y)$. We
refer the reader to \cite{carlob} for precise comparison statements. However, we remark that, contrary 
to the notion of crystals introduced in \cite{garo}, these $\s$-categories of derived $\D$-modules  
are sensitive to derived structures (see remark \ref{weyl2}).

In this situation, $T^*\F=T^*X$ and $T^*\G=T^*Y$ are the global derived cotangent stacks of $X$ and $Y$, and the GRR formula is the equality  
\begin{equation}\label{ddmodform}Ch(f_!(E))=f_!Ch(E) \,\,\, \textrm{in} \,\, K_0^{red}(T^*Y) \end{equation}
for a compact object $E \in \D_X-\dg_{qcoh}$ admitting a good filtration, and $f : X \to Y$ 
a proper and quasi-smooth morphism. When $X$ and $Y$ are smooth varieties, this recovers the well known 
GRR formula for $\D$-modules of \cite{la} (see also \cite{sab}). However, already when $X$ and $Y$ are
just underived $k$-schemes, formula (\ref{ddmodform}) is new and provides a Grothendieck-Riemann-Roch formula for $\D$-modules on \emph{possibly singular} schemes. 

An interesting feature of this situation is that any morphism $f=(g,u) : (X,\F) \to (Y,\G)$ fits in a commutative
diagram
$$\xymatrix{
(X,\F) \ar[r]^-{f} \ar[d]_-{u} & (Y,\G) \ar[d]^-{v} \\
(X,*_X) \ar[r]_-g & (Y,*_Y).
}$$
As explained at the end of \S \ref{s-ind}, the push-fowards $u_!$ and $v_!$ define
induction $\s$-functors  
$$u_! : \QCoh(\F) \to \D_X-\dg_{qcoh} \qquad
v_! : \QCoh(\G) \to \D_Y-\dg_{qcoh}.$$
The GRR formula for $f$, which is an equality in $K^{red}(T^*\G)$, can be
pushed forward along $v$ to a formula in $K_0^{red}(T^*Y)$:
$$v_!Ch(f_!(E))=g_!Ch(u_!E))$$
where $E$ is a coherent crystal along $\F$ admitting a good filtration, 
$u_!(E)$ the induced coherent $\D_X$-module, and $v_! : K^{red}_0(T^*\G) \to K^{red}_0(T^*Y)$
the direct image along the canonical morphism of foliations $\G \to *_Y$. This gives a formula for the push-forward of
the characteristic cycle of a $\D_X$-module \emph{induced from} the derived foliation $\F$. When $f=id$ the formula reads
$$v_!Ch(E)=Ch(u_!E)$$
and should be understood as a formula for the characteristic cycle on a $\D_X$-module induced
from the derived foliation $\F$. \\

\subsection{Smooth Lie algebroids} Let $X$ be a smooth variety and $\F \in \Fol(X)$ a smooth foliation on $X$ (i.e. $\LL_{\F}$ is a vector bundle on $X$). As explained in theorem
\ref{tII-4}, $\F$ is determined by a Lie algebroid $\T_\F \to \T_X$. The sheaf $\D^{\mathrm{fil}}_\F$ is then isomorphic to the universal
enveloping algebra $U(\T_\F)$ of $\T_\F$, equipped with its natural PBW filtration. 
As already noticed $K_0^{red}(T^*\F)\simeq K_0(T^*\F)$, in this case (of a smooth foliation on a smooth variety). It can also be shown that any coherent crystal along $\F$
admits a good filtration.
Moreover, $\D_\F$ satisfies all the 
conditions of Quillen theorem (see \cite{quillen}) and we thus have natural isomorphisms of K-groups
$$\tau_X : K_0(\Coh(\F)) \simeq K_0(T^*\F) \simeq K_0(X).$$
The first of this isomorphisms is precisely given by $E \mapsto Gr(E^{\mathrm{fil}})$ for $E^{\mathrm{fil}}$ a good filtration on $E$, and
the second isomorphism is the pull-back along the zero section $s : X \to T^*\F$.

The GRR formula \ref{tgrr} tells us that the isomorphism $\tau_X : K_0(\Coh(F)) \simeq K_0(X)$ is 
covariantly functorial in $(X,\F)$
in the following sense: given a morphism $f=(g,u) : (X,\F) \to (Y,\G)$ of smooth varieties endowed
with smooth foliations, with $g$ proper, the following diagram
$$\xymatrix{
K_0(\Coh(\F)) \ar[r]^-{\tau_X} \ar[d]_-{f_!} & K_0(X) \ar[d]^-{g_*(-\otimes \omega_{X/Y})[d]} \\
K_0(\Coh(\G)) \ar[r]_-{\tau_Y} & K_0(Y)
}$$
commutes. This recovers the well known GRR formula for $\D$-modules on smooth varieties 
(see \cite{la}), and its natural extension
to Lie algebroids. This extension to Lie algebroids is probably a folklore result (as its proof is word by word the same as for $\D$-modules), but we have not been able to find an appropriate reference
in the literature. \\

\subsection{Shifted Poisson structures} An important class of derived foliations are provided by 
\emph{shifted Poisson structures}
in the sense of \cite{cptvv}. Indeed, let $X$ be a derived scheme endowed with a 
shifted Poisson structure of degree $n$. The Poisson bracket defines 
a morphism of perfect complexes on $X$
$$a : \LL_X[-n] \to \T_X$$
making $\LL_X[-n]$ into a dg-Lie algebroid over $X$. This defines ($\S$ 
\ref{sec:ComparisonwithderivedLiealgebroids}) a derived foliation $\F$ 
whose underlying graded mixed cdga is $Sym_{\OO_X}(\T_X[n])$, where the mixed structure is
induced the bracket $[-,p]$, $p$ being the bivector defining the Poisson structure. 
We refer to 
\cite{tomić2025shiftedlagrangianthickeningsshifted} where this construction is explained in details, 
as several technical complications related to homotopy coherences have to be overcome.

This derived foliation associated to a shifted Poisson structure is the derived analogue of the 
foliation by symplectic leaves of a classical
Poisson strcuture on a smooth variety. Its \emph{leaves}, in the sense of Chapter \ref{ch:leafspaces}, 
are by definition
the symplectic leaves of the shifted Poisson structure. When the Poisson structure is non-degenerate, 
then 
the foliation $\F$ is the final foliation. In general $\F$ is a very interesting derived foliation 
containing 
information about the Poisson structure. For instance, our notion $\QCoh(\F)$ of quasi-coherent 
crystals provides a useful setting to study various versions of Poisson cohomology. 
As an example, we may define \emph{derived Poisson cohomology} to be 
$$H_{der}^{\mathrm{Poiss}}(X):=\rh_{\QCoh(\F)}(\OO_X,\OO_X),$$ 
which coincides with derived de Rham cohomology when the Poisson structure is a symplectic structure. 
It is also
possible to consider the induced $\D_X$-module $u_!(\OO_X) \in \D_X-\dg_{X,qcoh}$, where
$u : \F \to *_X$ is the unique morphism to the final foliation. Then, we may use the $\D_X$-module
$u_!(\OO_X)$  in order to define a version of Poisson \emph{homology} by 
$$H_{\mathrm{Poiss}}(X):=p_*(u_!(\OO_X)) \in \dg$$
where $p : X \to Spec\, k$ is the structure map (assuming here that $X$ is quasi-smooth and  $p$ admits 
a quasi-smooth compactification as defined in \S \ref{s-ell}). This offers a criterion for 
finiteness of Poisson homology
by requiring $u_!(\OO_X)$ to be Fredholm as a $\D_X$-module on $X$. When $X$ is smooth and 
the Poisson structure is a classical Poisson structure of degree $0$, the $\D_X$-module
$u_!(\OO_X)$ has been considered in \cite{brtr}, where it was used to define  the notion 
of \emph{holonomic Poisson varieties}
and to get finiteness results for Poisson cohomology. The GRR formula \ref{tgrrell}, and more generally 
the general formalism of crystals along derived foliations, provides a way to extend these notions 
and results  
to shifted Poisson structures. \\

\subsection{A foliated index formula} As a last application of Theorem \ref{tgrrell} we propose an \emph{index formula} for
a weakly Fredholm differential operator along the leaves of a derived foliation on a quasi-smooth derived
scheme. We like to think of it as an algebraic version of the longitudinal index theorem of \cite{cosk}, 
possibly valid outside the smooth setting, i.e. for non-smooth derived schemes and non-smooth
derived foliation over them. 

Let $(X,\F)$, where $X$ is a derived DM-stack and $\F\in \Fol^p(X)$ a perfect derived foliation. 
We assume that 
$X$ has a quasi-smooth compactification $j : X \hookrightarrow \bar{X}$ (i.e.  the projection $X\to 
Spec\, k$ admits a quasi-smooth compactification as in \S \ref{s-ell}), which 
implies in particular that $X$ is separated and quasi-compact. Note that we only
assume the existence of $\bar{X}$, and we do not assume that $\F$ can be extended to $\bar{X}$.
Let $E$ and $E'$ be two perfect complexes on $X$. A \emph{differential operator $P$ along $\F$
from $E$ to $E'$} is by definition a morphism of quasi-coherent sheaves on $X$
$$P : E \longrightarrow Ind_\F(E')=\D_\F \otimes_{\OO_X}E',$$
or, equivalently, a morphism of coherent crystals along $\F$
$$a_P : Ind_\F(E) \longrightarrow Ind_\F(E').$$
As $E$ is compact in $\QCoh(X)$, we note that $P$ must factor through a morphism
$E \to \D_\F^{\leq i}\otimes_{\OO_X}E'$ for some integer $i$. In this case, we say that 
$P$ is of \emph{order} $\leq i$.

We let $\D(P)$ be the cone of the morphism $a_P$ inside $\Coh(\F)$, and we want to apply  Theorem \ref{tgrrell} to the coherent crystal $\D_X \otimes_{\D_\F} \D(P)$. Note that 
this is also the cone of the morphism $\D_X \otimes_{\D_\F} E \to \D_X \otimes_{\D_\F} E'$, induced by the composition
$$P : \xymatrix{
E \ar[r] & \D_\F \otimes_{\OO_X}E' \ar[r] & \D_X \otimes_{\OO_X}E'. }$$
In order to apply Theorem \ref{tgrrell},  we need to impose a condition on $P$ ensuring that $\D_X \otimes_{\D_\F} \D(P)$ is a Fredholm object (Definition \ref{dell}).
We consider the least integer $i$ such that $P$ factors as
$$E \longrightarrow \D^{\leq i}_\F \otimes_{\OO_X}E'.$$
Associated to this, we have a morphism of filtered $\D_\F^{\mathrm{fil}}$-modules
$$Ind_\F^{\mathrm{fil}}(E)(-i) \to Ind_\F^{\mathrm{fil}}(E),$$
whose cone defines a good filtration $\D^{\mathrm{fil}}(P)$ on $\D(P)$, and thus, by 
base change, a good filtration $\D^{\mathrm{fil}}_X \otimes_{\D^{\mathrm{fil}}_\F} \D^{\mathrm{fil}}(P)$ on 
$\D_X \otimes_{\D_\F} D(P)$. 

\begin{df}\label{dellop}
The operator $P$ along $\F$ is \emph{weakly Fredholm} if 
$j_*(\D^{\mathrm{fil}}_X \otimes_{\D^{\mathrm{fil}}_\F} \D^{\mathrm{fil}}(P))$ is a compact object 
in $\QCoh^{\mathrm{fil}}(\bar{X})$, for a quasi-smooth compactification $j : X \hookrightarrow \bar{X}$.
\end{df}

Assuming that $P$ is weakly Fredholm, we are thus able to apply Theorem \ref{tgrrell}
to get our index formula, that has values in $K_0(k)\simeq \ZZ$, and thus is
an equality of two numbers we are now going to describe.\\
Let $f : X \to Spec\, k$ the projection, that we factor
as $\xymatrix{X \ar[r]^-j & \bar{X} \ar[r]^-{\bar{f}} & Spec\, k.}$
We first describe $f_*(\D(P))\simeq \bar{f_!}(j_*(\D(P)))\simeq \bar{f}_!(\D_X \otimes_{\D_\F} D(P))$.
This is a perfect complex of $k$-modules, which by definition, is the cone of the morphism induced by $P$
$$\Gamma(P) : \Gamma(X,E\otimes \omega_X)[d] \to \Gamma(X,E'\otimes \omega_X)[d].$$
Note that none of the two complexes $\Gamma(X,E\otimes \omega_X)$ or $\Gamma(X,E'\otimes \omega_X)$
is perfect, and only the cone of $\Gamma(P)$ is so. This is the effect of the weakly Fredholm property, implying that 
$\Gamma(P)$ is indeed a Fredlhom operator. Therefore, we have a first well defined number, called
the \emph{algebraic index of $P$} 
$Ind(P)$, which is the Euler characteristic of the cone of $\Gamma(P)$
$$Ind(P):=(-1)^d.\chi(cone \left(\Gamma(P) : \Gamma(X,E\otimes \omega_X) \to \Gamma(X,E'\otimes \omega_X)\right)).$$

On the other hand, the object $\D_X \otimes_{\D_\F} \D(P)$ is endowed with the good filtration
$\D^{\mathrm{fil}}_X \otimes_{\D^{\mathrm{fil}}_\F} \D^{\mathrm{fil}}(P)$, and by assumption its associated graded
defines a perfect complex of $T^*X$ which remains perfect on $T^*\bar{X}$. 
This associated graded can be described as follows. Recall that 
$i$ is the least integer such that $P$ defines a morphism in $\QCoh(X)$
$P : E \to \D_X^{\leq i} \otimes_{\OO_X}E'.$
By projection $\D_X^{\leq i} \to Sym_X^i(\T_X)$, we find 
$$\sigma(P)  : E \to Sym^i_X(\T_X) \otimes_{\OO_X}E'$$
which is called the \emph{symbol of $P$}. This extends to $$Sym_X(\T_X) \otimes_{\OO_X}E \to Sym_X(\T_X) \otimes_{\OO_X}E$$
as a morphism of graded $Sym_{\OO_X}(\T_X)$-modules 

The cone of this morphism defines a perfect complex, still denoted by $\sigma(P)$, on 
$T^*X$, which remains perfect on $T^*\bar{X}$ by the weakly Fredholm assumption.
If we denote by $s : X \to T^*X$ the zero section map, we thus get a well defined perfect complex of
$k$-modules
$\Gamma(X,s^*(\sigma(P))\otimes \omega_X)[d]$. The Euler characteristic of this complex is called
the \emph{K-theoretical index of the operator $P$}, and denoted by 
$$Ind_\sigma(P):=(-1)^d\chi(\Gamma(X,s^*(\sigma(P))\otimes \omega_X)).$$

Theorem \ref{tgrrell} then implies the following

\begin{cor}\label{cindex}
With the notations above, we have
$$Ind(P) = Ind_\sigma(P).$$
\end{cor}

In plain words, Corollary \ref{cindex}, says that the \emph{algebraic index of $P$}, computing the Fredholm index of $P$ acting on global sections of $E$ and $E'$, equals 
the \emph{$K$-theoretic index of $P$}, which should be understood as an intersection number between $X$ and the cycle defined by the symbol of $P$ inside
the total cotangent stack $T^*X$.

\chapter{Analytic aspects}\label{chapter-analyticaspects}

This chapter is devoted to introduce the notion of derived foliation in the holomorphic context. 
We remind some basic definitions on holomorphic rings and holomorphic 
cdga, and the corresponding notion of derived analytic stacks. We then 
introduce sheaves of holomorphic graded mixed cdga and use them to define 
holomorphic derived foliations on general derived analytic stacks in a very similar
fashion than the in the algebraic setting. We study analytic integrability of quasi-smooth and rigid
derived foliations, locally in the analytic topology. Finally we construct an analytification
functor and prove a verion of GAGA relating derived algebraic foliations and derived analytic
foliations on a proper derived Deligne-Mumford algebraic stack. \\

Throughout this chapter the base ring $k$ is the field of complex numbers $\CC$. 

\section{Quick reminders on derived analytic geometry}\label{sec_quickremonanalytic}

We recall 
the definition of an 
algebraic theory $hol$, of holomorphic rings as follows (see \cite{MR3337959,MR4036665,porta}
and \cite{MR3121621,nuiten} for the $\Cs$-version). 
The algebraic theory $hol$ is by definition the full sub-category of
$\CC Man$, of complex manifolds and holomorphic maps, 
consisting of all the objects $\CC^n$ for $n\geq 0$. 
We will denote by $hol(n)$ the set of holomorphic maps $\CC^n \to \CC$, so that 
in the category $hol$ the maps from $\CC^n$ to $\CC^m$ are in natural bijection with 
$hol(n)^m$.

By definition, a \emph{holomorphic ring} $A$
is an algebra over the algebraic theory $hol$, that is a product preserving
functor $A : hol \to Sets$. In more concrete terms, such a holomorphic ring 
consists
of a commutative $\CC$-algebra $A$, 
endowed with extra operations: for any 
holomorphic function $f : \CC^n \to \CC$
we have an application
$$\gamma_f : A^{n} \to A,$$
satisfying some natural properties with respect to compositions of holomorphic maps. 
Morally, $\gamma_f(a_1,\dots,a_n)$ stands for
"$f(a_1,\dots,a_n)$", the "evaluation of $f$ at the elements $a_i$". 
When $f$ is restricted to polynomials maps, the operations $\gamma_f$ determines
a commutative $\CC$-algebra structure on $A$ (the sum and multiplication 
being induced by the two holomorphic maps $\mathbb{C}^2 \to \mathbb{C}$, $(x,y) \mapsto x+y$ and $(x,y) \mapsto xy$).
The typical example of a holomorphic ring is the set $\OO_X(X)$ of (globally defined) holomorphic
functions on a given complex analytic space $X$, where the operations $\gamma_f$
are actually defined by $\gamma_f(a_1,\dots,a_n) = f(a_1,\dots,a_n)$, for $a_i \in \OO_X(X)$.

Morphisms between holomorphic rings are maps commuting with all the operations
$\gamma_f$ (for all $f$ as above). Holomorphic rings thus form a category denoted by $\CR^h$, which comes equipped with 
with a forgetful functor $\CR^h \to \CC-\CR$  to 
commutative $\CC$-algebras. We will often use the same notations for 
a holomorphic ring $A$ and its underlying $\CC$-algebra, but in situations where
the difference is important the latter will be denoted by $N(A)$. The forgetful 
functor $\CR^h \to \CC-\CR$ commutes with limits and filtered colimits, and thus 
possesses a left adjoint, denoted by $A \mapsto A^h$, sending
a $\CC$-algebra to its "holomorphication". This left adjoint
$A \mapsto A^h$ can be understood as a left Kan extension
as follows. We start by the functor sending the polynomial algebra $\CC[x_1,\dots,x_n]$ to the
ring $hol(n)$ of holomorphic functions on $\CC^n$ (with its natural
holomorphic structure), and any $\CC$-algebra morphism $\CC[X_1,\dots,X_n] \to \CC[Y_1,\dots,Y_m]$, 
defined by polynomials $P_1(Y), \dots, P_n(Y)$, to the map $hol(n) \to hol(m)$
obtained by composition with $(P_1,\dots, P_n) : \CC^m \to \CC^n$. We then 
obtain $(-)^h$ as the left Kan extension of this functor along the
canonical inclusion of the category of polynomial algebras into the whole
category of commutative $\CC$-algebras. 
For any $\CC$-algebra $A$
the holomorphic ring $A^h$ can thus be described by choosing a presentation 
of $A$ by polynomial algebras, whose image by 
$(-)^h$ provides a presentation by holomorphic rings $hol(n)$. More explicitely, 
if $A$ is finite generated, we can write $A$ as a co-equalizer 
$$\left(\CC[X_1,\dots,X_n] \rightrightarrows \CC[Y_1,\dots,Y_m]\right) \longrightarrow A$$
and obtain $A^h$ as the corresponding coequalizer in $\CR^h$
$$\left(hol(n) \rightrightarrows hol(m)\right) \longrightarrow A^h.$$
When $A$ is not finitely presented we write $A=colim_i A_i$ as a filtered colimit of finite presented
$\CC$-algebras and we have $A^h = colim_i (A_i)^h$. 

The forgetful functor $\CR^h \to \CR$ commutes with filtered colimits (and even 
sifted colimits), but does \emph{not} commute with push-outs in general. In order to avoid
confusions we will then use $\otimes^h$ to denote the categorical sum in $\CR^h$. For instance
the push-out $hol(n)\otimes^h hol(m)$ is naturally isomorphic to $hol(n+m)$, simply because
$hol(n)$ turns out to be the free holomorphic ring over $n$ generators
$$Hom_{\CR^h}(hol(n),A) \simeq N(A)^n.$$
However, the $\CC$-alegbras $hol(n+m)$ is certainly not isomorphic, 
by the canonical map of $\CC$-algebras, to the algebraic tensor product $hol(n) \otimes hol(m)$.
However, as $(-)^h$ is left adjoint, we always have a canonical isomorphism of holomorphic rings
$$A^h \otimes^h B^h \simeq (A\otimes_\CC B)^h$$
for any two commutative $\CC$-algebras $A$ and $B$.

Any holomorphic ring $A$ possesses a module of (holomorphic) differential forms $\Omega_A^1$, which is 
an $A$-module with the following universal property:  
$A$-modules maps $\Omega_A^1 \to M$ are in one-to-one
correspondence with \emph{holomorphic} derivations $\delta : A \to M$. Recall
here that a derivation $\delta : A \to M$ is holomorphic 
if for all $f \in hol(n)$ we have the chain rule
$$\delta(\gamma_f(a_1,\dots,a_n)) = \sum_i \frac{\partial f}{\partial z_i}(a_1,\dots,a_n).\delta(a_i).$$
It is important here to make a difference between 
a holomorphic ring $A$ and its underlying $\CC$-algebra $N(A)$. Indeed, 
$\Omega_A^1$ is different from the usual model of Kahler differentials $\Omega_{N(A)/\CC}^1$  of $N(A)$.
There is an obvious morphism $\Omega_{N(A)/\CC}^1 \to \Omega_A^1$ coming from the fact that the
universal holomorphic derivation $A \to \Omega_A^1$ is a derivation in the usual algebraic
sense, but this is not an isomorphism in general. 
On the other hand, if we start with a commutative $\CC$-algebra $R$, 
then there exists a natural isomorphism of $R^h$-modules (which follows from the fact 
that the forgetful functor $\CR^h \to \CC-\CR$ preserves trivial square zero extensions)
$$R^h\otimes_R \Omega_{R/\CC}^1 \simeq \Omega_{R^h}^1.$$
When $A$ is a holomorphic ring $\Omega_A^1$ will always mean 
the holomorphic differentials, while we will use the notation $\Omega_{N(A)}^1$ for the 
(non-holomorphic) Kahler differentials. \\

\subsection{The $\s$-category of holomorphic cdga's}

The next definition is essentially similar 
to the notions used in \cite{MR3337959,MR4036665}. Note however that in the \emph{non-connective} case, our
notion of holomorphic cdga differs slightly, as the holomorphic structure will exist
on the cohomological degree $0$ subalgebra (and not merely on the algebra of $0$-cocyles). Note
also the extra condition on the differential
which does not appear in the above mentioned references. However, for connective cdga's it coincides
with the notions of \cite{MR3337959,MR4036665} (note that for connective cdga's, a holomorphic structure, as in Definition \ref{def_holcdga} below, is just a holomorphic ring structure on its degree $0$ part). We will non-connective holomorphic
cdga's in an essential manner, for instance in order to consider 
the holomorphic \emph{de Rham complex} as a holomorphic cdga.

\begin{df}\label{def_holcdga}
\begin{enumerate}
    \item 
Let $A \in \cdga_\CC$ be a possibly non-connective 
$\CC$-linear cdga. A \emph{holomorphic structure on $A$} 
consists of a holomorphic ring structure on the $\CC$-algebra $A^0$ of elements of cohomological
degree $0$, such that the cohomological differential $d : A^0 \to A^{1}$ is
a holomorphic derivation. A \emph{holomorphic cdga} is a $\CC$-linear cdga endowed
with a holomorphic structure.
\item 
A \emph{morphism $A\to B$, between two holomorphic cdga's},
is a morphism of (underlying) cdga's for which the induced morphism $A^0 \to B^0$ 
is a morphism of holomorphic rings.
\end{enumerate}
\end{df}

Holomorphic cdga's form a category $cdga^h$, endowed with a forgetful functor
$cdga^h \to cdga_\CC$, to the category of $\CC$-linear cdga's. This functor is the right adjoint
of an adjunction, whose left adjoint is denoted by $A \to A^h$. Explicitly, 
for $A \in cdga_\CC$, $A^h$ is defined as follows. We start by the graded
$\CC$-algebra $(A^0)^h\otimes_{A^0}A$, by extending $A^0$ to its
holomorphication. The cohomological differential $d : A^0 \to A^1$ extends uniquely
into a holomorphic derivation $d^h : (A^0)^h \to (A^0)^h\otimes_{A^0}A^1$. We then 
define the cohomological differential $d^h$ on the whole graded
algebra $(A^0)^h\otimes_{A^0}A$ by the formula
$$d^h(a\otimes x) = d^h(a).x  + a.d(x),$$
for any homogenuous $x \in A$ and $a \in (A^0)^h$. This defines a structure of cdga
on $A^h$. Moreover, $(A^h)^0=(A^0)^h$ comes equipped with its natural holomorphic ring structure
making $A^h$ into a holomorphic cdga in the sense of the above definition.

More generally, it is possible to describe the coproducts in the category $cdga^h$ 
as follows. For this we start by understanding the 
case of $\ZZ$-graded holomorphic rings $(\CR^h)^\ZZ$. These are, by definition,  
graded commutative $\C$-algebras $A$, with the usual sign rules $xy = (-1)^|x||y| = yx$, 
together with a holomorphic structure on $A^0$. We have a corresponding
notion of graded holomorphic derivations and a graded module of holomorphic
differential forms $\Omega_{A}^1$ for $A\in (\CR^h)^\ZZ$. In particular, 
any holomorphic cdga $A$ can be understood via its underlying 
graded holomorphic algebra $A^{gr} \in (\CR^h)^\ZZ$, together with 
its differential given by a morphism of graded $A^{gr}$-modules
$$d : \Omega_{A^{gr}}^{1} \longrightarrow A^{gr}(-1),$$
where $A^{gr}(-1)$ is the weight twist by $1$ of the unit graded module $A^{gr}$. With these notations
coproducts in $(\CR^{h})^h$ are easily identified by the formula
$$(A \otimes^h B)^n = \bigoplus_{p+q = n}(A^p\otimes B^q)\otimes_{A^0\otimes B^0}(A^0 \otimes^h B^0),$$
with the usual commutative multiplicative structure. Note that by construction 
$(A \otimes^h B)^0 = A^0 \otimes^h B^0$ is the coproduct of holomorphic rings in weight $0$.

Suppose now that $A$ and $B$ are two objects in $cdga^h$, which we represent as 
graded holomorphic rings $A^{gr}$ and $B^{gr}$ together with graded holomorphic derivations
$$d_A : \Omega_{A^{gr}}^1 \to A^{gr}(-1) \qquad B_A : \Omega_{B^{gr}}^1 \to B^{gr}(-1).$$
We can form the coproduct in graded rings $A^{gr}\otimes^h B^{gr}$, and the two 
differentials $d_A$ and $d_B$ provide for us 
a canonical morphism of graded $A^{gr}\otimes^h B^{gr}$-modules
$$d_A \otimes id + id \otimes d_B : 
\Omega_{A^{gr} \otimes^h B^{gr}}^1 \simeq (\Omega_{A^{gr}}^1 \otimes_{A^{gr}} 
A^{gr}\otimes^h B^{gr}) \oplus 
(A^{gr} \otimes^h B^{gr} \otimes_{B^{gr}} \Omega_{B^{gr}}^1) \longrightarrow (A^{gr}\otimes^h 
B^{gr})(-1).$$
This defines a differential on the graded holomorphic ring 
$A^{gr}\otimes^h B^{gr}$ which is holomorphic in degree $0$, and thus
defines a holomorphic cdga. This holomorphic cdga is the coproduct of $A$ and $B$ computed in 
the category $cdga^h$. The important point here is that 
the forgetful functor $cdga^h \to (\CR^h)^\ZZ$, from holomorphic cdga's to graded holomorphic rings, 
preserves coproducts. \\

The adjunction $cdga^h \leftrightarrows cdga_\CC$ can be used in order to lift the model
structure on $cdga_\CC$ to a model structure on $cdga^h$, for which fibrations and
equivalences are defined via the forgetful functor, as shown in the proposition below.

\begin{prop}\label{modstructureonholcdgas}
The above notions of equivalences and fibrations make $cdga^h$ into a model category for which the
forgetful functor $cdga^h \to cdga$ is a right Quillen functor.
\end{prop}

\textit{Proof.} We apply the criterion of transferring model structures along a right adjoint 
using the path object argument (see \cite[2.6]{MR2016697}). The only non-trivial statement is 
to prove that any object $A \in \cdga^h$ admits a path object. We construct
this path object explicitly has follows. 

We introduce the object $\CC[\Delta^1]^h \in \cdga^h$ which is the image
by $(-)^h$ of the $\cdga$ freely generated by the 
acyclic complex $\left(\xymatrix{\CC \ar[r]^-{id} & \CC}\right)$
concentrated in cohomological degrees $[0,1]$. It can be explicitly described as
the holomorphic cdga with only two non-trivial degrees 
$$(\CC[\Delta^1]^h)^0 = (\CC[\Delta^1]^h)^1 = hol(1) \qquad (\CC[\Delta^1]^h)^i = 0 \, for \, i\neq 0,1$$
and whose differential in degree $0$ is given by 
$$hol(1) \to hol(1)$$
sending a holomorphic function $f$ to its derivative $f'=\frac{\partial f}{\partial z}$.
The holomorphic cdga $\CC[\Delta^1]^h$ admits two augmentation to $\CC$, by 
simply evaluating a function $f$ at $0$ and at $1$
$$ev_0,ev_1 : \CC[\Delta^1]^h \rightrightarrows \CC.$$
We claim that for any $A\in \cdga^h$ the induced diagram 
$$\xymatrix{
A \ar[r]^-j & A \otimes^h \CC[\Delta^1]^h \ar@<-.5ex>[r]_-{ev_1} \ar@<.5ex>[r]^-{ev_0} &  A
}$$
is a path object for $A$: $j$ is a weak equivalence, 
$ev_0 \times ev_1 : A \otimes^h \CC[\Delta^1]^h \to A \times A$ is a fibration 
and $ev_i \circ j = id$. The last equality $ev_i\circ j =id$ is obvious and we focus
on the two other properties.

As a start, we want to describe $A \otimes^h \CC[\Delta^1]^h$ explicitly. 
We know that coproducts in $cdga^h$ can be described on the level of the underlying
graded holomorphic rings by forgetting the differential. Therefore, in degree $i$ we have 
$$(A \otimes^h \CC[\Delta^1]^h)^i \simeq (A \otimes^h hol(1))^i \oplus (A\otimes^h hol(1))^{i-1},$$
and the differential is given by 
$$\begin{pmatrix} d  & \partial \\ 0 & d \end{pmatrix} : (A \otimes^h hol(1))^{i} \oplus (A \otimes^h 
hol(1))^{i-1} \longrightarrow 
(A\otimes^h hol(1))^{i+1} \oplus (A \otimes^h hol(1))^i,$$
where $d$ the differential of $A\otimes^h hol(1)$, and $\partial : (A \otimes^h hol(1))^{i-1} \to 
(A \otimes^h hol(1))^{i}$ is the natural morphism
induced by the holomorphic derivation $\partial/\partial z : hol(1) \to hol(1)$
by base change to $A \otimes^h hol(1)$. In other words, 
the underlying complex of $\CC$-modules of $A \otimes^h \CC[\Delta^1]^h$ is the cocone of the morphism
$$\partial : A \otimes^h hol(1) \to A \otimes^h hol(1).$$
In order to prove that $j$ is a weak equivalence, it is therefore enough to show that 
the natural inclusion $A \to A \otimes^h hol(1)$ induces a short exact sequence of complexes
\begin{equation}\label{eq:shortexact}
    \xymatrix{
0 \ar[r] & A \ar[r] & A \otimes^h hol(1) \ar[r]^-{\partial} & A \otimes^h hol(1) \ar[r]& 0.}
\end{equation}
The statement we want to prove here is idenpendant of the differential of $A$ and we can thus
work entirely with the underlying graded holomorphic rings. We start by the degree $0$ part. 
The holomorphic ring $A^0$ can be written as a filtered colimit in $\CR^h$ of finitely presented
holomorphic rings $A^0 = colim_i B_i$. Each $B_i$  is thus the holomorphic ring of 
closed analytic space $X_i \subset \CC^n$ (for some $n$ depending on $i$), given by 
a finite number of global holomorphic equations. For each $i$, we can 
identify $B_i\otimes^h hol(1) \simeq \OO(X_i \times C)$. The diagram above is then isomorphic to the
following
$$\xymatrix{
0 \ar[r] & colim_i \OO(X_i) \ar[r] & colim_i \OO(X_i \times \CC) \ar[r]^-{\partial} & 
colim_i \OO(X_i\times \CC) \ar[r]& 0.}$$
Moreover, $\partial : \OO(X_i \times \CC) \to \OO(X_i \times \CC)$ turns out to be the
partial derivative along the canonical coordinate on $\CC$. It is then clear that the diagram
above is a filtered colimit of short exact sequences and thus is an exact sequence. This shows that
\ref{eq:shortexact} is indeed exact in degree $0$. Exactness in non-zero degrees is proven similarly.
For $i$ fixed, let $M=A^i$ considered as an $A^0$-module, and we can replace $A$ by the graded 
holomorphic ring $A^0 \oplus M[-i]$ in order to prove that \ref{eq:shortexact} is exact in degree $i$.
We first reduce to the case where
$A^0$ is of the form $\OO(X)$ for some $X \subset \CC^n$ as above. We then write $M$ as a filtered
colimit of finite presented $\OO(X)$-modules, and thus reduce to the case where 
$M$ is finite presented. 
In this case it corresponds to a well defined coherent $\cM$ sheaf
on the Stein space $X$ with a global presentation $\OO_X^q \to \OO_X^p \to \cM \to 0$. 
The diagram \ref{eq:shortexact} in degree $i$ then can be identified with
$$\xymatrix{
0 \ar[r] & \Gamma(X,\cM) \ar[r] & \Gamma(X \times \CC,\pi^*(\cM)) \ar[r]^-{\partial} & 
\Gamma(X\times \CC,\pi^*(\cM)) \ar[r]& 0,}$$
where $\pi : X \times \CC \to X$ is the first projection and $\partial$ is again 
the partial derivative along the coordinate on $\CC$, which exists canonically on $\pi^*(\cM)$. 
The above diagram is again a short exact sequence. 

We have thus shown that the sequence \ref{eq:shortexact} is exact, and thus that 
the morphism $j : A \to A \otimes^h \CC[\Delta^1]^h$ is indeed a weak equivalence of holomorphic 
cdga's. It remains to show that 
$$ev_0 \times ev_1 : A \otimes^h \CC[\Delta^1]^h \to A\times A$$
is a fibration, or equivalently that for any $i$ the induced map
$(A \otimes^h \CC[\Delta^1]^h)^i \to A^i\times A^i$ is surjective. This last map is the composition
$$(A \otimes^h hol(1))^i \oplus (A\otimes^h hol(1))^{i-1} \to 
(A \otimes^h hol(1))^i \to A^i \times A^i,$$
where the first morphism is the canonical projection on the first factor, and
the second map is induced by $ev_{0,1} : hol(1) \to \CC \times \CC$
given by evaluating a function $f$ at the subset $\{0,1\} \subset \CC$. As the projection is obviously 
surjective it remains to show that 
$$(A \otimes^h hol(1))^i \to A^i \times A^i,$$
is always a surjective morphism for any $A \in cdga^h$. Again, this is a statement idenpendant
of the differential of $A$ and thus can be proven more generally for any 
graded holomorphic ring $A$. As above, when $i=0$ we reduce the question to the case
where $A^0 = \OO(X)$ for $X \subset \CC^n$ globally defined by a finite number of holomorphic equations.
Then, the map under consideration is the evaluation at $X\times \{0,1\}$
$$\OO(X\times \CC) \to \OO(X)\times \OO(X).$$
This is surjective as a couple of functions $(f,g)$ will be the image
of $z.f + (1-z).g$ where $z$ is the global coordinate on $\CC$. The case where $i\neq 0$
is treated similarly as before, by reducing to the case where $A^0=\OO(X)$ and furthermore
$A^i$ is a finitely presented $A^0$-module corresponding to a coherent sheaf $\cM$ on 
$X$ which admits a global free resolution. The map under consideration in degree $i$ is then
$$\Gamma(X\times \CC,\cM) \to \Gamma(X,\cM) \times \Gamma(X,\cM).$$
To prove the surjectivity of this last map we can write $\cM$ as a quotient $\OO_X^p \to \cM$.
As $X$ is Stein, the induced morphism $\Gamma(X,\OO_X^p) \to \Gamma(X,\cM)$ is surjective, and thus
we are reduced to the case where $\cM$ is free which follows from the case already treated $\cM =\OO_X$.

This finishes the proof that $A\otimes^h \CC[\Delta^1]^h$ is a path object in $cdga^h$, and thus
by \cite[2.6]{MR2016697} of the proposition. \hfill $\Box$ \\

The $\s$-category
obtained from $cdga^h$ by inverting the weak equivalences (the morphisms inducing quasi-isomorphisms
on the underlying cdga's) will be denoted by $\cdga^h$. The Quillen adjunction of the above proposition 
produces
an adjunction of $\s$-categories
$$\cdga_\CC \leftrightarrows \cdga^h,$$
whose left adjoint will be denoted by $A \mapsto A^h$. The forgetful right adjoint, if necessary, 
will be denoted by $A \mapsto N(A)$, but most of the time it will not be written explicitly (so that 
we'll simply write again $A$ for the underlying cdga). \\

\subsection{Cotangent complexes}

As a first comment, for a holomorphic cdga $A \in \cdga^h$ an $A$-dg-module will be, by definition, 
a dg-module over $N(A)$ the underlying cdga of $A$. Theire $\s$-category is denoted as usual
$$\dg_A := \dg_{N(A)}.$$
Any holomorphic cdga $A$ possesses a module of \emph{holomorphic differential forms} $\Omega_A^1$. This
is an $A$-dg-module such that morphisms of dg-modules $\Omega_A^1 \to M$
are in one-to-one correspondence with holomorphic derivations $A \to M$. In order
to be more precise, we introduce the holomorphic trivial square zero extension $A \oplus M$. This
is a holomorphic cdga whose underlying cdga is the trivial square zero extension
of $N(A)$ by $M$. The holomorphic structure on $A^0 \oplus M^0$ is defined by 
formula
$$\gamma_f(a_1+m_1,\dots,a_n+m_n) = \gamma_f(a_1,\dots,a_n) + 
\sum_{i}\frac{\partial f}{\partial z_i}(a_1,\dots,a_n).m_i.$$
The holomorphic derivations from $A$ to $M$ are then defined as sections
of the natural projection $A \oplus M \to A$. They can be naturally identified with
derivations $d : N(A) \to M$ on the underlying cdga which satisfies furthermore the
condition that $d^0 : A^0 \to M^0$ is a holomorphic derivation from the holomophic ring $A^0$
to the $A^0$-module $M^0$.

As in the algebraic case, the construction $A \mapsto \Omega_A^1$ can be
left derived in order to produce a cotangent complex $\LL_A$ for any 
$A \in \cdga^h$. This cotangent complex satisfies the usual properties
of stability by base change and functoriality inside $\cdga^h$. For
$A \in \cdga_\CC$ we have a natural equivalence
$$A^h \otimes_A \LL_{A} \simeq \LL_{A^h},$$
as this follows by adjunction from the straight forward observation that $N(A\oplus M)\simeq 
N(A) \oplus M$ (i.e. trivial square zero extensions are preserved by the forgteful $\s$-functor
$\cdga^h \to \cdga_\CC$).

Any object $A \in \cdga^h$ possesses a connective cover $\tau_{\leq 0}A$. This is the right adjoint
of the inclusion $\s$-functor $\cdga^h_{c} \hookrightarrow \cdga^h$, where 
$\cdga^h_{c}$ stands for the full sub-$\s$-category spanned by holomorphic cdga's 
whose underlying cdga are connective. The connective cover construction 
$A \mapsto \tau_{\leq 0}A$ is moreover compatible with the forgetful $\s$-functor to $\cdga$. 
Indeed, we simply have to notice that if $A$ is a holomorphic cdga, the kernel
$Ker(d : A^0 \to A^{-1}) \subset A^0$ is a sub-holomorphic ring of $A^0$ (because 
$d$ is a holomorphic derivation). 

A connective holomorphic cdga $A \in  \cdga^h_c$ is called a 
\emph{finite cell object} if there is a finite sequence of morphisms
in $\cdga_c^h$
$$\xymatrix{
0=A_0 \ar[r] & A_1 \ar[r] & \dots \ar[r] & A_{n-1} \ar[r] & A_n=A,}$$
where each $A_{i+1}$ is obtained from $A_i$ by choosing some elements $a_i \in A_i^{d_i}$
of degree $n_i$, and freely adding a finite number of
variables $y_i$ of degrees $n_{i-1}$ with $d(y_i)=a_i$.  
In a slightly weaker manner, 
$A \in \cdga^h_c$ is an \emph{almost finite cell object} if there is a countable sequence of morphisms
$$\xymatrix{
0=A_0 \ar[r] & A_1 \ar[r] & \dots \ar[r] & A_{n-1} \ar[r] & A_n \ar[r] & \dots, }$$
with $A=\mathrm{colim}_i A_i$ and 
where each $A_i \to A_{i+1}$ is as above but with the condition that the sequence
of integers $n_i$ tends to $\s$. In other words, an almost finite cell
object is a cell object with a finite number of cell in each dimension. 

\begin{df}
\begin{enumerate}
    \item A connective holomorphic cdga $A$ is \emph{of finite presentation} if 
it is a retract, in the $\s$-category $\cdga_c^{h}$ of a finite cell object.
    \item A connective holomorphic cdga $A$ is \emph{almost of finite presentation} if it is a retract of
an almost finite cell object.
\end{enumerate}
\end{df}

Using \cite[Prop. 2.2]{tova}, it can be shown that finitely presented holomorphic cdga's 
are the compact objects in $\cdga_c^h$. Similarly, almost finitely presented
holomorphic cdga's are the object $A$ 
for which $Map(A,-)$ commutes with filtered colimits of uniformly truncated objects: for
any filtered system $\{A_i\}$, with $H^j(A_i)=0$ for all $j\leq n$ (for some fixed
integer $n$), we have 
$colim_i Map(A,A_i) \simeq Map(A,colim_i A_i)$. 

\subsection{Derived complex analytic spaces and stacks}

\begin{df}
The \emph{$\s$-category of affine derived analytic spaces} is the full 
sub-$\s$-category of the opposite $\s$-category of $\cdga_c^{h}$ formed by 
holomorphic cdga's almost of finite presentation. 
It is denoted by $\dAff^{h}$. The object in $\dAff^{h}$ corresponding
to an object $A \in \cdga_c^{h}$ will be denoted symbolically by $\Spec^{h}\, A$.
\end{df}

We note that for an affine derived analytic space 
$X=\Spec^h\, A$, the holomorphic ring $H^0(A)$ can be canonically identified with  
the ring of holomorphic functions on 
a closed analytic subspace $X$ in $\CC^n$. More precisely, 
we can write $A$ as a connective cell object with finitely many cells in each dimension. 
For such a cell object, $A^0$ is a free holomorphic ring on $n$ generators (for some
$n$), and thus $H^0(A)$ becomes a quotient of the holomorphic ring $hol(n)=\OO(\CC^n)$ 
by a finitely number of relations. 
These relations are given by elements $f_i \in hol(n)$ and thus they define
a complex analytic subspace $\tau_0(X)$ inside $\CC^n$. The holomorphic ring $H^0(A)$ is then 
naturally isomorphic to the holomorphic ring $\OO_X(X)$ of holomorphic functions on $\tau_0(X)$.
The complex space $\tau_0(X)$ is called the \emph{truncation} of $X$.
In the same manner, 
the cohomology groups $H^i(A)$ are finitely generated $H^0(A)$-modules, and therefore they
correspond to globally generated coherent sheaves on $\tau_0(X)$ 
(see \cite[Prop. 1.24, Lem. 1.25]{MR4036665}). \\

The natural transcendent topology on $\CC$ defines a Grothendieck topology
on the $\s$-categories $\dAff^{h}$. If $i : \Spec^h\, A\to \Spec^h\, B$ 
is a morphism of affine derived analytic spaces, 
we say that $i$ is an open immersion if the induced morphism $i_0 : \Spec^h\, H^0(A) \to 
\Spec^h\, H^0(B)$ is an open immersion of analytic spaces, and if furthermore the
natural morphisms $H^i(A)\otimes_{H^0(A)}H^0(B) \to H^i(B)$ are isomorphisms
(i.e. the coherent sheaves $H^i(A)$ on $\Spec^h\, H^0(A)$ restrict 
to $H^i(B)$ along the open immersion $i_0$). An open covering of affine derived analytic spaces 
is then a family of open immersion $\{U_i \to X\}_i$ such that $\coprod_i \tau_0(U_i) \to \tau_0(X)$
is a surjective map of complex spaces.

We can also define smooth and \'etale morphisms, simply by using the holomorphic cotangent complexes. 
We say that $\Spec^h\, A \to \Spec^h\, B$ is smooth (resp. étale) if the relative cotangent complex 
$\LL_{B/A}$ is a projective $B$-module of finite rank (resp. $\LL_{B/A}\simeq 0$). We refer
to \cite{MR3337959,MR4036665,porta} for more details on these notions. \\

We are now in a situation in which we have a Grothendieck $\s$-site $\dAff^{h}$
together with a notion of smooth morphisms. We can thus
define Artin stacks by using the formal geometricity procedure (see \cite[Ch. 1.3]{hagII}, \cite[3.3]{toenems} and \cite[\S 4]{pv}). In particular, 
we have the $\s$-topos of stacks on this $\s$-site, which will be denoted by 
$\dSt^h$ and whose objects are called \emph{derived analytic stacks}. We also have the full sub-$\s$-category of Artin (aka geometric) stacks inside $\dSt^h$.

\begin{rmk} \emph{Note that, by definition, the site 
$\dAff^h$ consists only of almost finitely presented objects, and thus
Artin stacks in these settings will be automatically locally almost 
finitely presented. In particular, for any Artin derived analytic stack $X$, 
its truncation $\tau_0 X$ is an (underived) Artin analytic stack 
locally of finite presentation, and the homotopy groups sheaves $H^i(\OO_X)$ 
always define coherent sheaves on $\tau_0(X)$.} 
\end{rmk}

\section{De Rham algebras in the  holomorphic context}

We introduce here a category of \emph{holomorphic} $\CC$-linear \emph{graded mixed} cdga's. 
Its objects are graded mixed cdga's $B$ (see \S  \ref{sec:GradedmixedcomplexesandderiveddeRhamtheory})  
together with a holomorphic structure on the weight $0$ cdga $B^{(0)}$, and such that 
the mixed structure $\epsilon : B^{(0)} \to B^{(1)}[-1]$ is a holomorphic derivation
on the holomorphic cdga $B^{(0)}$.
We can endow this category with a model category structure simply by defining equivalences
and fibrations on the underlying complexes $B^{(i)}$ (by forgetting the
holomorphic structures). 

\begin{prop}
The category $(\epsilon-CAlg(k)^{gr})^{h}$ of holomorphic graded mixed cdgas, with the above notions of fibrations 
and equivalences, defines a model category structure such that the forgetful functor
to graded mixed cdga's is right Quillen.
\end{prop}

\textit{Proof.} This is similar to the proof of Proposition \ref{modstructureonholcdgas} (with model structure on $\epsilon-CAlg(k)^{gr}$ as in Definition \ref{dI-3}), and in fact easier, since the forgetful functor
$$(\epsilon-CAlg(k)^{gr})^{h} \longrightarrow \epsilon-CAlg(k)^{gr}$$
commutes with push-outs. We leave the details to the reader. \hfill $\Box$ \\

The $\s$-category obtained from the model category of holomorphic graded mixed cdga's  will be denoted
by $(\egrcdga)^{h}$, and will be called the $\s$-category of \emph{holomorphic
graded mixed cdga's}. For any 
holomorphic cdga $A$ its holomorphic de Rham algebra is defined using the
holomorphic cotangent complex
$$\DR(A)=Sym_A(\LL_A[1]).$$ 
It is first a graded mixed cdga where the mixed structure is the de Rham differential. 
But, the de Rham differential being obviously holomorphic, $\DR(A)$ is an object
in $(\egrcdga)^{h}$. We can state things differently, by considering the $\s$-functor
$(\egrcdga)^h \to \cdga^h$, sending $B$ to $B^{(0)}$, and observe that it has a left adjoint
given by $B \to \DR(B)$. This adjunction of $\s$-categories is moreover induced
by a Quillen adjunction on the level of model categories. 
Of this adjunction we retain the following property

\begin{prop}\label{pA-10-1}
Let $A$ be a holomorphic cdga, and $\DR(A)$ its holomorphic de Rham complex
as an object in $(\egrcdga)^{h}$. For any $B \in (\egrcdga)^{h}$, the morphism induced on 
mapping spaces
$$\Map_{(\egrcdga)^{h}}(\DR(A),B) \longrightarrow \Map_{\cdga^h}(A,B^{(0)})$$
is an equivalence.
\end{prop}

\begin{rmk}
\emph{There is of course a big difference between $\DR(A)$ and $\DR(N(A))$, as already 
the holomorphic cotangent complex $\LL_A$ differs from the algebraic cotangent 
complex $\LL_{N(A)}$.}
\end{rmk}

We will also make use of a second adjunction between holomorphic graded mixed cdga's and holomorphic cdga's. 
Namely, we consider the $\s$-functor $A \mapsto A(0)$, sending a holomorphic cdga $A$ to
the holomorphic graded mixed cdga $A(0)$ consisting of $A$ sitting purely in weight $0$ and endowed with the trivial
mixed structure. This defines an $\s$-functor $\cdga^h \to (\egrcdga)^h$ which commutes with 
colimits, and so it admits a right adjoint. This right adjoint is called the \emph{holomorphic
realization} and is denoted by 
$$|-| : (\egrcdga)^h \to \cdga^h.$$

\begin{prop}
For any $B \in (\egrcdga)^h$, the underlying cdga $N(|B|)$ associated to $|B|$
is canonically equivalent to $|N(B)|$, the realization (as in Remark  \ref{nonTaterealizationforgradedmixedalgebras}) of the underlying graded mixed
cdga $N(B)$ associated to $B$.
\end{prop}

\textit{Proof.} This is clear by adjunction and from the fact that $A(0)^h$ obviously identifies with 
$A^h(0)$. For $A \in \cdga_\CC$, we have
a chain of natural equivalences
$$Map_{\cdga_\CC}(A,N(|B|)) \simeq Map_{\cdga^h}(A^h,|B|) \simeq
Map_{(\egrcdga)^h}(A(0)^h,B)$$
$$\simeq Map_{\egrcdga}(A(0),N(B)) \simeq Map_{\cdga_\CC}(A,|N(B)|).$$
\hfill $\Box$ \\

\begin{df}\label{dholdRcohomology}
Let $A$ be a (connective) holomorphic cdga and $\DR(A)$ its de Rham holomorphic graded mixed cdga.
The \emph{holomorphic derived de Rham complex of $A$} is defined by 
$\CDR^*(A) := |\DR(A)| \in \cdga^h$.
\end{df}

\begin{rmk}
\begin{enumerate}
\item \emph{There is discrepancy between the notation of Definition \ref{dholdRcohomology}
and Definition \ref{dI-4}, as the former involved also a filtration. To be perfectly
coherent we should have used instead the notation $\CDR^*(A)^u$. However,
we will not use the Hodge
filtration in the holomorphic context, so will keep the notation $\CDR^*(A)$ as above.  
It is clear that Definition \ref{dholdRcohomology} also has a refined filtered version, 
using the natural notion of a} filtered holomorphic cdga. \emph{We will not pursue is that direction
and leave these details to the interested readers.}
    \item \emph{In most cases we will only consider the underlying cdga of $\CDR^*(A)$ and will not use
often its natural holomorphic structure. However, this structure will be useful for some specific arguments (see 
for instance the proof of Theorem \ref{tV-1}). Because of this, 
we will still denote by $\CDR^*(A)$ the underlying cdga and simply 
call it the derived de Rham complex of $A$.}
\end{enumerate}

\end{rmk}

\section{Analytic derived foliations}

For any affine derived analytic space $\Spec^h\, A$, 
we remind that we have a holomorphic graded mixed cdga $\DR(A)$ whose underlying
graded cdga is $Sym_A(\LL_A[1])$, and for which the mixed structure
in given by the de Rham differential.
The construction $A \mapsto \DR(A)$ defines an $\s$-functor
$$\DR(-) : (\dAff^h)^{op} \longrightarrow (\egrcdga)^h.$$
We then consider the $\s$-functor
$$\Fol^{pr} : (\dAff^h)^{op} \longrightarrow \scat,$$
sending $X=\Spec^h \, A \in \dAff^h$ to the $\s$-category of holomorphic graded mixed $\DR(A)$-cdga's 
$B$ such that the following two conditions are satisfied
\begin{enumerate}
    \item the natural morphism $A \to B^{(0)}$ is a quasi-isomorphism
    
    \item $B^{(1)}$ is an almost perfect complex of $A$-dg-modules
    
    \item the natural morphism of graded cdga's
    $$Sym_A(B^{(1)}) \to B$$
    is a graded quasi-isomorphism.
\end{enumerate}

As opposed to the algebraic situation, the $\s$-functor $\Fol^{pr} : (\dAff^h)^{op} \longrightarrow \scat$
is only a prestack i.e. it does not satisfy descent. 
Indeed, if $X\subset \CC^n$ is a closed analytic subspace, 
globally defined by a finite number of equations, and with
holomorphic ring of functions $A$, then almost perfect $A$-modules correspond to 
almost perfect complexes on $X$ than can be written as bounded complex of globally generated
vector bundles on $X$. Of course, these are not \emph{all} the almost perfect complexes on $X$. 
The associated stack $\Fol$ to $\Fol^{pr}$ can however be computed easily. For this,
we sheafify everything on $X$. We first consider the small
site of affine opens $U \subset X$. The holomorphic ring $A$ 
sheafifies to a sheaf of holomorphic cdga's $\OO_X$.
We next consider $\DR_X$ the sheaf $(U=\Spec^h\, B \subset X) \mapsto \DR(B)$ of 
holomorphic graded mixed cdga's on the underlying topological space of $X$. 
The objects in $\Fol(X)$ can then be described as
sheaves of holomorphic graded mixed cdga's $B$, together with a morphism 
$\DR_X \to B$ in $(\egrcdga)^h(X)$ such that 

\begin{enumerate}[label=({\alph*})]
    \item the natural morphism $\OO_X \to B^{(0)}$ is a quasi-isomorphism
    of sheaves of cdga's on $X$
    
    \item $B^{(1)}$ is a almost perfect complex of $\OO_X$-dg-modules
    
    \item the natural morphism of sheaves of graded cdga's on $X$
    $$Sym_{\OO_X}(B^{(1)}) \to B$$
    is a quasi-isomorphisms of sheaves of graded complexes.
\end{enumerate}

By definition $\Fol(X)$ is a full sub-$\s$-category of the opposite of
$\DR_X/(\egrcdga)^h(X)$ of sheaves of holomorphic graded mixed cdga's under $\DR_X$. 
We thus get an $\s$-functor $\Fol : (\dAff^h)^{op} \to \scat$ satisfying the descent 
property can be left Kan extended to a colimit preserving $\s$-functor on the $\s$-category of all
derived analytic stacks
$$\Fol : (\dSt^h)^{op} \longrightarrow \scat.$$

\begin{df}\label{dV-1}
For a derived analytic stack $X\in \dSt^h$, the \emph{$\s$-category of derived 
foliations on $X$} is $\Fol(X)$. For a morphism of derived analytic stacks $f : X \to Y$
the induced $\s$-functor $\Fol(Y) \to \Fol(X)$ is called the \emph{pull-back of derived foliations}
and is denoted by $f^*$.
\end{df}

\begin{rmk}
\emph{By definition derived foliations in the analytic settings are always almost perfect by 
assumptions. This is only because we want to avoid technical complications related to the notion
of quasi-coherent sheaves in the analytic setting, and moreover all the examples treated in this
book will be almost perfect.}
\end{rmk}

We finish this section  by mentioning that most of the general results and notions we have seen 
for derived foliations on derived stacks remains valid in the analytic setting, possibly with mild modifications.
As a start, if $\F \in \Fol(X)$ is a derived foliation on a
derived analytic stack $X \in \dSt^h$, we can define its fake (or big) cotangent 
complex $\LL_\F$, which is an $\OO_X$-module, simply by sending $u : S:=\Spec^h\,A \to X$
to the $\OO_S$-module $\LL_{u^*(\F)}$. Here, $u^*(\F) \in \Fol(\Spec^h\,A)$ and can thus
be represented by a sheaf of graded mixed $\DR_S$-algebras $\DR_{u^*(\F)}$
satisfying conditions (a), (b), (c) above. The $\OO_S$-module $\LL_{u^*(\F)}$ is then by 
definition $\DR_{u^*(\F)}^{(1)}[-1]$, the (shifted) weight one piece of 
$\DR_{u^*(\F)}$, which by assumption is an almost perfect $\OO_S$-module. We get this
way a family of $\OO_S$-modules, functorial in $S \to X$, or in other
words an $\OO_X$-module denoted by $\LL_\F^{big}$. It comes equipped with a canonical
morphism $\LL_X^{big} \to \LL_\F^{big}$. 

A first observation is that the cone of this morphism is an almost perfect 
$\OO_X$-modules. Indeed, this is due to the fact that conormal complexes
of derived foliations are stable by pull-backs. We denote this cone by $\N_\F^*$.
Suppose now that $X$ is a derived analytic Artin stack, and let us consider
its cotangent complex $\LL_X$, which comes equipped with its canonical 
morphism $\LL_X^{big} \to \LL_X$. 

\begin{df}
Let $X$ be an analytic derived Artin stack and $\F \in \Fol(X)$ be a derived foliation on $X$.
The \emph{cotangent complex of $\F$}, $\LL_\F$, is defined by the cartesian square of $\OO_X$-modules
$$\xymatrix{
\LL_X^{big} \ar[d] \ar[r] & \LL_\F^{big} \ar[d] \\
\LL_X \ar[r] & \LL_\F.
}$$
\end{df}

We note that the cone of $\LL_X \to \LL_\F$ is again $\N_\F^*$, that we have seen to be
almost perfect. Therefore, $\LL_\F$ is itself an almost perfect $\OO_X$-module. Also, from the definition
it is easy to deduce the functorial nature of cotangent complexes. For a morphism
$f : X \to Y$ of analytic derived Artin stacks, and $\F \in \Fol(Y)$, we have a Cartesian 
square of almost perfect $\OO_Y$-modules
$$\xymatrix{
f^*(\LL_Y) \ar[r]\ar[d] & f^*(\LL_\F) \ar[d] \\
\LL_X \ar[r] & \LL_{f^*(\F)}.
}$$

Finally, we write $\DR_\F$ for the sheaf of holomorphic graded mixed cdga's associated to $\F \in \Fol(X)$, which is defined on the big site $\dAff^h/X$ of affine derived analytic spaces
over $X$. We denote by $\CDR^*(\F)=|\DR_\F| $ its realization
in the sense of Definition \ref{dholdRcohomology}, which is a sheaf of holomorphic cdga's on $X$.  

\section{Analytic integrability and existence of derived enhancements}\label{section_anintderenh}

The following local integrability result is specific to the analytic setting, and is 
a version of the classical Frobenius theorem for transversely smooth and rigid derived foliations. 
It is a direct consequence of a result of Malgrange combined with the analytic version of 
Corollary \ref{cpII-8}. However, we still display the result as a theorem to insist on its importance.

As in the algebraic setting, we will say that a derived foliation $\F \in \Fol(X)$ on 
a derived analytic space $X$ is \emph{analytically integrable}, if there
exists a morphism $f : X \to Y$ to another derived analytic space $Y$ such that 
$\F \simeq f^*(0_Y)$, where $0_Y \in \Fol(Y)$ is the initial foliation on $Y$. When such an
$f$ only exists locally on $X$ we say that $\F$ \emph{locally analytically integrable}.

\begin{thm}\label{tV-1}
Let $X$ be smooth Deligne-Mumford analytic stack, and $\F \in \Fol(X)$
be derived foliation. We suppose that the following two conditions are satisfied.
\begin{enumerate}
    \item The derived foliation $\F$ is transversely smooth and rigid.
    \item The derived foliation is smooth in codimension $2$.
\end{enumerate}
The, locally on $X$, $\F$ is integrable.
\end{thm}

\textbf{Proof.} The result being local on $X$ we can assume that $X$ is affine and we will
allow to shrink it furthermore when necessary.
For any point $x\in X$, we consider the formal completion $\hat{X}_x$ of $X$ at $x$, and
the restriction $\hat{\F}_x$ of $\F$ on $\hat{X}_x$. We know by our results on 
formal integrability that $\hat{\F}_x$ is formally integrable (see Corollary \ref{cpII-8}). 

We now consider the truncation $K$ of $\F$ in the sense of Definition \ref{dII-14}. This 
truncation is a sheaf of holomorphic graded mixed cdga of the form 
$Sym^{naive}_{\OO_X}((\Omega_{X}^1/K)[1])$, 
where $K \subset \Omega_{X}^1$ is the differential ideal defined as the kernel of the
morphism $\Omega_X^1 \to H^0(\LL_\F)$. The sub-sheaf $K$ is also the image of 
$\N_\F^* \to \Omega_X^1$. By generic smoothness the moprhism between 
vector bundles $\N_\F^* \to \Omega_X^1$ must be a monomorphism 
as it is a monomorphism on the dense open $U$. Therefore $K$ is, as an abstract $\OO_X$-module 
isomorphic to $\N_\F^*$. Shrinking $X$ is necessary
we can assume that $\N_\F^*$ is a free vector bundle, and thus that $K \simeq \OO_X^d$. 
The image of the basis elements of $\OO_X^d$ in $\Omega_X^1$ defines 
differential forms $(\omega_1,\dots,\omega_d)$ on $X$ which are linearly independent
over $\OO_X$. We can thus apply \cite[Thm. 3.1, B]{mal2} to the forms
$(\omega_1,\dots,\omega_d)$ and get that (possibly after localizing on $X$), 
there exists a holomorphic map $f : X \to S$, where $S$ is another complex manifold (that can be
chosen to be an open in $\CC^d$), 
such that $f$ integrates $K$: $f^*(\Omega_S^1)=K \subset \Omega_X^1$. We now claim that 
$f$ also integrates $\F$ in the sense that $\F$ is isomorphic to $f^*(0_S)$ as derived
foliations on $X$. 

For this, we consider $\OO_\F \subset \OO_X$ the subring of
functions which are killed by the relative de Rham differential $d : \OO_X \to \Omega_X^1/K$.
This is also the $0$-th naive de Rham cohomology group of the differential ideal 
$K$, denoted by $H^0_{dR,naive}(K)$.
By our assumptions on the codimension of the singularities of $\F$,  
we know that the natural morphism from derived de Rham cohomology along $\F$
to naive de Rham cohomology of $K$ induces an isomorphism of sheaves of holomorphic rings
(see \ref{clcpII-10})
$$H^0_{dR}(\F) \longrightarrow H^0_{dR,naive}(K)=\OO_\F.$$
In particular, as $\CDR^*(\F)=|\DR_\F|$ is the realization of 
$\DR_\F$, and because quasi-smoothness implies that $H^i(|\DR_\F|)\simeq 0$ for all $i<0$, 
we have a canonical morphism of sheaves of
holomorphic cdga's
$$\OO_\F\simeq H^0_{dR}(\F) \longrightarrow |\DR_\F|,$$
which we can compose with the adjunction morphism
$|\DR_\F| \to \DR_\F$ where the sheaf of holomorphic cdga 
$|\DR_\F|$ is considered with its trivial graded mixed structure (remind that $|-|$ is right adjoint 
to the $\s$-functor remembering only the weight $0$ piece, see after Proposition \ref{pA-10-1}).

\begin{lem}\label{ltV-1}
There exists a commutative square of sheaves of holomorphic graded mixed cdga's on $X$
$$\xymatrix{\DR_X \ar[r] & \DR_\F  \\
\DR_{\OO_\F} \ar[u] \ar[r] & \OO_\F, \ar[u]}$$
where $\OO_\F \to \DR_\F$ is the natural morphism described above, and 
$\DR_{\OO_\F} \to \DR_{X}$ the morphism induced from the inclusion 
of sheaves of holomorphic rings $\OO_\F \to \OO_X$.
\end{lem}

\textit{Proof of the lemma.} The square of the lemma provides two possible composed morphisms
$\DR_{\OO_\F} \to \DR_\F$. These are morphisms of holomorphic graded mixed cdga's. 
By the universal 
property of $\DR_{\OO_\F}$ (see Proposition \ref{pA-10-1}), in order construct a homotopy between 
these two morphisms it is enough to see that they induce the homotopic morphisms of holomorphic
rings in weight $0$. But, by construction, these two morphisms in weight $0$ both coincide with the 
canonical inclusion $\OO_\F \to \OO_X$ and thus are homotopic in a uniquely possible manner
(the space of self homotopies of this inclusion is contractible). \hfill $\Box$ \\

By the previous lemma, we have constructed a natural morphism of sheaves of holomoprhic graded mixed cdga's on $X$
$$\DR_{\OO_X/\OO_\F}=\DR_X\otimes_{\DR(\OO_\F)}{\OO_\F} \longrightarrow \DR_\F.$$
We claim that this morphism is an equivalence. Indeed, it is enough to check that 
it induces an equivalence on the associated graded, or equivalently (by multiplicativity) on weight $1$. Now, in weight $1$ the above morphism is the canonical morphism
$\LL_{\OO_X/\OO_\F} \to \LL_\F,$
where $\LL_{\OO_X/\OO_\F}$ is the holomorphic relative cotangent complex of $\OO_\F \subset \OO_X$, 
and the morphism above is induced by the holomorphic derivation
$\OO_X \to \DR_\F^{(1)}$ which vanishes on $\OO_\F$. To see that this morphism
$\LL_{\OO_X/\OO_\F} \to \LL_\F$ is an equivalence we can work locally on $X$, and thus
assume that the differential ideal $K$ is integrable by a morphism $f : X \to S$
with $S$ a smooth manifold. In this case, \cite[Thm. 2.1.1]{mal2} shows that 
$\OO_\F\simeq f^{-1}(\OO_S)$, and thus $\LL_{\OO_X/\OO_\F}$ is identified with 
the coherent sheaf $\Omega_f^1$ of relative $1$-forms on $X$ over $S$. The fact that 
$\Omega_f^1 \simeq \LL_\F$ then follows from the very definition of the fact that 
$f$ integrates the differential ideal $K$. 

To summarize, we have seen that under the hypothesis of the theorem $\DR_\F$ is
of the form $\DR_{\OO_X/\OO_\F}$. Moreover, locally on $X$, $\DR_{\OO_X/\OO_\F}$
is the relative de Rham algebra $\DR_{X/S}$ for a morphism $f : X \to S$ integrating 
the differential ideal $K$. This shows, in particular, that $f$ also integrates $\F$
as claimed. \hfill $\Box$ \\

Also the following two existence and uniqueness results are specific to the analytic setting
and does not have an algebraic counterpart (see however Corollary \ref{cGAGAperfectfol} in 
the proper case).
They are very close to Theorem \ref{tV-1}, and they are again consequences
of the results of \cite{mal2}.

\begin{prop}\label{pV-1}
Let $X$ be a smooth complex manifold and $K \subset \Omega_X^1$ be a differential ideal 
such that $K$ is a vector bundle.
\begin{enumerate}
    \item If the coherent sheaf $\Omega_X^1/K$ is a vector bundle in codimension $2$, 
    then there exists at most one, up to equivalence, 
    transversely smooth and rigid derived enhancement for $K$.
    \item If the coherent sheaf $\Omega_X^1/K$ is a vector bundle in codimension $3$, 
    then there exists a unique, up to equivalence, transversely smooth and rigid derived enhancement for $K$.
\end{enumerate}
\end{prop}

\textbf{Proof.} $(1)$ This is a consequence of the proof of Theorem \ref{tV-1}. Indeed,
in that proof we have shown that if $\F \in \Fol(X)$ is a transversely smooth and 
rigid derived enhancement of $K$, then we must have $\DR_\F \simeq \DR_{\OO_X/\OO_\F}$, 
for $\OO_\F \subset \OO_X$ the holomorphic subring of functions killed by
the de Rham differential $d : \OO_X \to \Omega^1_X/K$. This shows indeed that 
$K$ determines $\F$ and thus implies the uniqueness of $\F$. 

$(2)$ The point $(1)$ already insures unicity, so it remains to prove existence. By unicity
it is enven enough to prove existence locally on $X$.
The argument is then the 
same as in the proof of Theorem \ref{tV-1}. We let $\OO_\F \subset \OO_X$ be the
holomorphic subring of functions killed by the de Rham differential $d : \OO_X \to \Omega_X^1/K$. 
We set $\DR_\F:=\DR_{\OO_X/\OO_\F}$ be the relative holomorphic de Rham algebra of the
inclusion $\OO_\F \subset \OO_X$. It is endowed with the canonical
morphism $\DR_X \to \DR_\F$, and we claim that as such it is a transversely smooth
and rigid derived enhancement of $K$. For this we use case (A) of \cite[Thm. 3.1]{mal2}
which insures that $K$ is locally integrable. During the proof of Theorem \ref{tV-1} we have
seen that $\DR_\F$ is locally of the form $\DR_{X/S}$ for some morphism $f : X \to S$ with $S$ smooth. 
This clearly implies that the sheaf of graded mixed $\DR_X$-cdga $\DR_\F$ satisfies
the conditions of being a transversely smooth and rigid derived foliation and 
that the corresponding differential ideal is equal to $K$. \hfill $\Box$ \\

\section{Analytification of derived foliations: the GAGA theorem}\label{sec_GAGAderfol}

We consider the full sub-$\s$-category $\dSt_\CC^{aft} \subset \dSt_\CC$, consisting
of derived stacks which are locally almost of finite presentation over $\CC$. These are the 
derived stacks that can be obtained as colimits of objects of the form $\Spec\, A$ with
$A$ a connective cdga which is almost of finite presentation over $\CC$. The canonical 
embedding $\dSt_\CC^{aft} \subset \dSt_\CC$ can also be identified with the left Kan 
extension of the Yoneda embedding restricted to derived affine schemes locally almost of finite
presentation $\dAff_\CC^{aft} \subset \dAff_\CC \subset \dSt_\CC$. In particular, 
an alternative description of $\dSt_\CC^{aft}$ is as the $\s$-category of (hypercomplete)
stacks on the \'etale $\s$-site $\dAff_\CC^{aft}$. 

The analytification $\s$-functor $\cdga^{aft}_{c, \mathbb{C}} \to \cdga_c^{h}$  of \S \ref{sec_quickremonanalytic} induces an analytification $\s$-functor
$$(-)^h : \dSt^{aft}_\CC \longrightarrow \dSt^h$$
from the $\s$-category of derived stacks locally almost of finite presentation to the $\s$-category
of derived analytic stacks. It is obtained as the left Kan extension of the analytification $\s$-functor
followed by the Yoneda embedding
$$\xymatrix{\dAff_\CC^{aft} \ar[r] & \dAff^h \ar[r] & \dSt^h}$$
where the leftmost $\s$-functor is the (opposite of the) analytification functor of \S \ref{sec_quickremonanalytic}, and it sends $\Spec\, A$
to $\Spec^h \, A^h$, where $A^h$ is the holomorphic cdga generated by $A$. The analytification $\s$-functor $(-)^h : \dSt^{aft}_\CC \to \dSt^h$
is the left adjoint of the restriction $\s$-functor, sending $X \in \dSt^h$
to the derived stack defined by $A \mapsto X(A^h)$. This restriction will be denoted by $r$ so that we have 
an adjunction
$$(-)^h : \dSt^{aft}_\CC \leftrightarrows \dSt^h : r.$$
Moreover, the right adjoint $r$ is a geometric morphism of $\s$-topoi, or equivalently 
$(-)^h$ commutes with finite limits. To avoid confusions with the various $(-)^h$ notations involved
(for instance for graded mixed cdga's), 
we will denote this adjunction by 
$$u^{-1} : \dSt^{aft}_\CC \leftrightarrows \dSt^h : u_*.$$

Remind that to both $\s$-topoi $\dSt^{aft}_\CC$ and $\dSt^h$  are associated canonical stacks
(in $\s$-categories) of \emph{almost perfect} derived foliations, 
denoted by $\Fol^{ap}$ in the algebraic setting (see Definition \ref{dII-8}), and simply by $\Fol$ in the holomorphic setting. 
These can be compared by means of
the analytification $\s$-functor as shown in the following proposition.

\begin{prop}
There exists a canonical morphism of derived analytic stacks
$$u^{-1}(\Fol^{ap}) \longrightarrow \Fol.$$
\end{prop}

\textit{Proof.} By adjunction it is enough to construct a morphism $\Fol^{ap} \to u_*(\Fol)$. 
For an object $\Spec\, A \in \dAff_\CC^{aft}$, we consider the analytification construction
$(-)^h : \Fol^{ap}(X) \to \Fol(X^h)$, sending a graded mixed cdga $B$ with $B^{(0)} \simeq A$
to $B^h$, which is a holomorphic graded mixed cdga with $(B^h)^{(0)}\simeq A^h$
satisfying the conditions to define a holomorphic derived foliations on 
$X^h$.
This is clearly
functorial in $A$, and thus defines a morphism $\Fol^{ap} \to u_*(\Fol)$ as required. \hfill $\Box$ \\

The above proposition is needed in order to construct the analytification $\s$-functor
for derived foliations. Indeed, for $X \in \dSt_\CC^{aft}$, we define
an $\s$-functor
$$(-)^h : \Fol^{ap}(X) \to \Fol(X^h)$$
by means of the analytification $\s$-functor followed by the morphism of proposition
$$\Fol^{ap}(X)\simeq Map(X,\Fol^{ap}) \to Map(X^h,u^{-1}(\Fol)) \to Map(X^h,\Fol) \simeq \Fol(X^h).$$

\begin{df}\label{def_analytificationforfoliations}
For $X \in \dSt_\CC^{aft}$, the $\s$-functor defined above
$$(-)^h : \Fol^{ap}(X) \to \Fol(X^h)$$
is called the \emph{analytification $\s$-functor for derived foliations}.
\end{df}

It is easy to check that $(-)^h$ extends to a morphism of stacks $ \Fol^{ap}(-) \to \Fol(-) \circ (-)^h$. i.e. it is functorial for pullbacks of derived foliations: for any $f: X \to Y$ in $\dSt_\CC^{aft}$, there is an induced $\s$-functor $f^* : \Fol(Y^h) \to \Fol (X^h)$, compatible with composition.\\

For the next proposition, remind from \cite{porta}
the existence of an analytification $\s$-functor for 
almost perfect complexes
$$(-)^h : \APerf(X) \to \APerf(X^h),$$
functorially in $X \in \dSt_k^{aft}$.

\begin{prop}\label{prop_anpreservesctgtcplxforfoliation} For $X \in \dSt_\CC^{aft}$
a derived Artin stack locally of finite presentation, 
the analytification functor $(-)^h : \APerf(X) \to \APerf(X^h)$ preserves cotangent complexes of derived
foliations, i.e. for any $\F \in \Fol^{ap}(X)$, we have a canonical isomorphism 
$$(\LL_\F)^{h} \simeq \LL_{\F^{h}}$$ 
in $\APerf(X^h)$.
\end{prop}
\textit{Proof.} This is an easy consequence of \cite[Thm. 5.21]{PoYu_repr}. \hfill $\Box$ \\

The following result is our version of GAGA theorem for (almost perfect) derived foliations.

\begin{thm}\label{GAGAperfectfol}
Let $X$ be a proper algebraic Deligne-Mumford stack. Then the $\s$-functor
$$(-)^h : \Fol^{ap}(X) \to \Fol(X^h)$$
from almost perfect algebraic derived foliations on $X$ to almost perfect holomorphic derived foliations
on $X^h$, is an equivalence of $\s$-categories.
\end{thm}

\textit{Proof.} We start by enlarging quite a bit the $\s$-categories $\Fol^{ap}(X)$
and $\Fol(X^h)$. For this, let $C_X$ be the $\s$-category of sheaves of graded mixed
cdga's $B$ on $X_{\acute{e}t}$ together with an augmentation $B \to \OO_X$ (of graded mixed cdga's for the trivial
graded mixed structure on $\OO_X$) and satisfying the following conditions.

\begin{enumerate}
    \item The augmentation $B \to \OO_X$ induces an equivalence on weight 
    zero parts $B^{(0)} \simeq \OO_X$.
    
    \item The negative weight pieces of $B$ vanish: $B^{(i)}\simeq 0$ for all $i<0$.
    
    \item For all $i$, the $\OO_X$-module $B^{(i)}$ is almost perfect.
    
\end{enumerate}

Similarly, we define $C_{X^h}$ the $\s$-category of sheaves of holomorphic graded mixed
cdga's $B$ on $X^h_{et}$ together with an augmentation $B \to \OO_{X^h}$ (of 
holomorphic graded mixed cdga's for the trivial
graded mixed structure on $\OO_{X^h}$) and satisfying the following conditions.

\begin{enumerate}
    \item The augmentation $B \to \OO_{X^h}$ induces an equivalence on weight 
    zero parts $B^{(0)} \simeq \OO_{X^h}$.
    
    \item The negative weight pieces of $B$ vanish: $B^{(i)}\simeq 0$ for all $i<0$.
    
    \item For all $i$, the $\OO_{X^h}$-module $B^{(i)}$ is almost perfect.
    
\end{enumerate}

As a first step, we observe the following finiteness result for cotangent complexes.

\begin{lem}\label{auxlemmayes}
Let $B \in C_X$ (resp. $B \in C_{X^h}$), then the graded cotangent complex 
$\LL_B$ is such that each individual weight piece $\LL_B^{(i)}$ is 
almost perfect over $\OO_X$ and is zero for $i<0$.
\end{lem}

\textit{Proof of Lemma \ref{auxlemmayes}.} We give the proof for objects in $C_X$, 
the holomorphic case can be treated similarly.

First of all, this is  a statement about the underlying 
graded cdga so we can forget the mixed structure. 
We construct a "graded cell decomposition" of $B$
in the usual way, for which cells in a given dimension will be parametrized
by an almost perfect $\OO_X$-modules.
We construct, by induction, a sequence of sheaves of graded $\OO_X$-cdga's
$$
    \xymatrix{
C(0)=\OO_X \ar[r] & C(1) \ar[r] & \dots \ar[r] & C(n) \ar[r] & C(n+1) \ar[r] & \dots \ar[r] & B}
$$
such that $C(n) \to B$ induces an equivalence in weights less or equal to $n$, and for each $n$
there exists a push-out square of sheaves of graded cdga's
\begin{equation}\label{eq:pushoutC(n)}\xymatrix{
C(n) \ar[r] & C(n+1) \\
Sym_{\OO_X}(E(n)) \ar[u] \ar[r] & \OO_X \ar[u]
}\end{equation}
with $E(n)$ an almost perfect $\OO_X$-module pure of weight $n+1$. 
Assuming that $C(n) \to B$ has been constructed, we consider the induced morphism 
$C(n)^{(n+1)} \to B^{(n+1)}$ on weight $n+1$ pieces
and we denote by $E(n)$ is homotopy fiber, which 
is a graded $\OO_X$-module pure of weight $(n+1)$. The graded $\OO_X$-cdga $C(n+1)$ is then
defined by the push-out above. Because the morphism $K(n) \to B$ is
canonically homotopic to zero, the morphism $C(n) \to B$ factors canonically 
through $C(n+1) \to B$. Finally, $C(n)$ and $C(n+1)$ do coincide in weight less than $n$, so
$C(n+1) \to B$ induces an equivalence in weights less than $n$. Moreover, the weight $n$ part 
of $C(n+1)$ is by definition equivalent to the cone of $E(n) \to C(n)^{(n+1)}$, and thus
is sent by an equivalence to $B^{(n+1)}$ as required. 

The existence of the sequence $C$ does imply the lemma, as $\LL_B$ will coincide with 
$\LL_{C(n)}$ in weight less than $n$. Moreover, by induction 
on $n$ and from the push-out square (\ref{eq:pushoutC(n)}), $\LL_{C(n)}$ is seen to be 
a graded perfect $C(n)$-module for all $n$, and therefore 
is graded almost perfect over $\OO_X$. \hfill $\Box$ \\

A second observation is that $C_X$ (resp. $C_{X^h}$) contains $\Fol^{ap}(X)$ (resp.
$\Fol(X^h)$) as a full sub-$\s$-category.

\begin{lem}\label{auxlemmayesyes}
The natural forgetful $\s$-functors 
$$\Fol^{ap}(X) \longrightarrow C_X \qquad \Fol(X^h) \longrightarrow C_{X^h}$$
are fully faithful. Their essential images consists of $B \in C_X$ (resp. 
$B \in C_{X^h}$) such that $Sym_{\OO_X}(B^{(1)}) \simeq B$ (resp.
$Sym_{\OO_{X^h}}(B^{(1)}) \simeq B$) as graded cdga's.
\end{lem}

\textit{Proof of Lemma \ref{auxlemmayesyes}.} This follows by simply unraveling the various definitions,
and the universal property of the de Rham algebras. We give the argument for $C_X$, the case
of $C_{X^h}$ being completely analogous.

For two derived foliations $\F$ and $\F'$ 
in $\Fol^{ap}(X)$, corresponding to sheaves of graded mixed $\DR_X$-cdga's $\DR_\F$ and $\DR_{\F'}$, 
the mapping space $Map(\F,\F')$ is by definition the homotopy fiber, taken at the identity 
of $\OO_X$, of the natural morphism
$$Map_{\egrcdga(X)}(\DR_{\F'},\DR_\F) \longrightarrow Map_{\egrcdga(X)}(\DR_X,\DR_\F) \simeq $$ $$ \simeq
\Map_{\cdga(X)}(\OO_X,\DR_\F^{(0)}) \simeq
Map_{\cdga(X)}(\OO_X,\OO_X)$$
(where we have denoted $\egrcdga_X$ and $\cdga_X$ the $\s$-categories of sheaves of graded mixed
cdga's and of cdga's on $X_{\acute{e}t}$).
Analogously, when $\DR_\F$ are considered as augmented towards $\OO_X$ via their
natural augmentation to $\DR_\F^{(0)}\simeq \OO_X \simeq \DR_{\F'}^{(0)}$, their mapping space in 
$C_X$ is the homotopy fiber, taken at the identity, of the natural morphism
$$Map_{\egrcdga(X)}(\DR_{\F'},\DR_\F) \longrightarrow Map_{\egrcdga(X)}(\DR_{\F'},\DR_\F^{(0)}) \simeq$$
$$ \simeq Map_{\cdga(X)}(\DR_{\F'}^{(0)},\DR_\F^{(0)}) \simeq Map_{\cdga(X)}(\OO_X,\OO_X).$$
These two morphisms are easily seen to be equivalent, and thus so are their homotopy fibers. 
This shows that the $\s$-functors of the lemma are fully faithful. Similarly, 
if $B \to \OO_X$ is an object in $C_X$, choosing an inverse of the equivalence $B^{(0)} \simeq \OO_X$
provides a morphism $\OO_X \to B^{(0)}$, and thus a morphism $\DR_X \to B$ by the universal
property of de Rham alegbras. The graded mixed cdga $B$, together with this morphism from $\DR_X$ clearly
defines an object in $\Fol^{ap}(X)$, whose image by the forgetful $\s$-functor is equivalent to $B$, showing
the essential surjectivity. \hfill $\Box$ \\

The analytification
$\s$-functor $(-)^h : \Fol^{ap}(X) \to \Fol(X^h)$ extends to 
$(-)^h : C_X \to C_{X^h}$, simply by sending a graded mixed cdga to the corresponding 
holomorphic graded mixed cdga by the analytification construction. Concretely, 
this sends a sheaf of graded mixed cdga $B$ on $X_{\acute{e}t}$ 
to the sheaf (on $X^{h}_{\acute{e}t}$) of graded holomorphic 
cdga $\OO_{X^h}\otimes_{u^{-1}(\OO_X)}u^{-1}(B)$ endowed
with its canonical mixed structure together with its natural augmentation to $\OO_{X^h}$. To prove the theorem.  
It is therefore enough to prove that this extended $\s$-functor $(-)^h : C_X \to C_{X^h}$ is an 
equivalence of $\s$-categories. 

In order to show this, we consider the forgetful $\s$-functor
to graded cdga's by forgetting the mixed structures. We let $C^o_X$ be the $\s$-category of sheaves 
of graded
cdga's $B$ on $X_{\acute{e}t}$ together with an augmentation $B \to \OO_X$ and satisfying the
graded analogues of the previous three conditions: $B^{(0)} \simeq \OO_X$, $B^{(i)}$ are
zero for $i<0$ and are almost perfect for $i>0$.  Similarly we have 
$C^o_{X^h}$ consisting of sheaves 
of graded
cdga's $B$ on $X^h_{\acute{e}t}$ together with an augmentation $B \to \OO_{X^h}$ and satisfying the
graded analogues of the previous three conditions: $B^{(0)} \simeq \OO_{X^h}$, $B^{(i)}$ are
zero for $i<0$ and are almost perfect for $i>0$. We have forgetful $\s$-functors 
$C_X \to C^o_X$ and $C_{X^h} \to C^o_{X^h}$ which commute with analytification functors
\begin{equation}\label{eq:forgetanalyt}
    \xymatrix{
C_X \ar[r]^-{(-)^h} \ar[d] & C_{X^h} \ar[d] \\
C^o_X \ar[r]_-{(-)^h} & C^o_{X^h}.
}
\end{equation}

\begin{lem}\label{lemauxYES}
The analytification $\s$-functor
$$(-)^h : C^o_X \to C^o_{X^h}$$
is an equivalence of $\s$-categories.
\end{lem}

\textit{Proof of Lemma \ref{lemauxYES}.} This is a simple application of the GAGA theorem for almost 
perfect complexes of proper derived Artin stack of \cite{porta}. Indeed, GAGA implies that 
the analytification $(-)^h : \APerf(X) \to \APerf(X^h)$ is an equivalence (note that 
$\APerf$ is denoted by $Coh^-$ in \cite{porta}). We deduce that it also induces
an equivalence on the $\s$-categories of $\NN$-graded objects $\APerf^\NN(X) \to \APerf^\NN(X^h)$.
As this equivalence is moreover an equivalence of symmetric monoidal $\s$-categories, we get 
an induced equivalence on commutative algebra objects. Finally, considering commutative algebra
augmented towards $\OO_X$ and $\OO_{X^h}$ we get the statement of the lemma. \hfill $\Box$ \\

In order to deduce the theorem from the lemma we use the Barr-Beck theorem \cite[Theorem 4.7.3.5]{HA}. 
Indeed, the forgetful $\s$-functors $C_X \to C_X^o$ and $C_X \to C_{X^h}$ are both
monadics. The left adjoint to $C_X \to C_X^o$ sends an augmented sheaf of graded algebras $B \to \OO_X$
to $Sym_B(\LL_B[1])$, the graded de Rham algebra of $B$ endowed with its total graduation and 
its natural mixed structure induced by the de Rham differential. The weight $0$ part
of $Sym_B(\LL_B[1])$ is clearly $\OO_X$, and thus we do get an object in $C_X$. Similarly, 
the forgetful $\s$-functor $C_{X^h} \to C^o_{X^h}$ has a left adjoint given by the 
graded holomorphic de Rham algebra construction. Because cotangent complexes commutes with 
analytification (see Proposition \ref{prop_anpreservesctgtcplxforfoliation}), we have that the commutative square \ref{eq:forgetanalyt} is adjointable and
the adjoint square 
$$ \xymatrix{
C_X \ar[r]^-{(-)^h} & C_{X^h}  \\
C^o_X \ar[u] \ar[r]_-{(-)^h} & C^o_{X^h} \ar[u]
}$$
naturally commutes (the canonical natural transformation between the two composition
is an equivalence). Therefore, the equivalence $(-)^h : C_X^o \simeq C_{X^h}^o$ preserves
the monads defining $C_X$ and $C_{X^h}$, and we thus have that the induced
$\s$-functor $(-)^h : C_X \to C_{X^h}$ is an equivalence as required. \hfill $\Box$ \\

Combining Proposition \ref{pV-1} and Theorem \ref{GAGAperfectfol} we get the following important 
corollary.

\begin{cor}\label{cGAGAperfectfol}
Let $X$ be a smooth and proper algebraic complex variety 
and $K \subset \Omega_X^1$ be an algebraic differential ideal
such that $K$ is a vector bundle.
If the coherent sheaf $\Omega_X^1/K$ is a vector bundle in codimension $3$, 
then there exists a unique, up to equivalence, transversely smooth and rigid derived enhancement for 
$K$ in the sense of Definition \ref{dII-14}.
\end{cor}

The theorem \ref{GAGAperfectfol} also has a linear version, for crystals which are 
almost perfect as $\OO_X$-modules. The proof of the theorem below can 
be formally deduced from Theorem \ref{GAGAperfectfol} by interpreting graded mixed 
$\DR_\F$-modules (and in particular thus crystals) as trivial square zero 
extensions of $\DR_\F$. Its proof will thus be only sketched.

\begin{thm}\label{GAGAperfectcrys}
Let $X$ be a proper algebraic Deligne-Mumford stack, $\F \in \Fol^{ap}(X)$
and almost perfect derived foliation on $X$ and $\F^h \in \Fol(X^h)$ its analytification. 
The analytification $\s$-functor
$$(-)^h : \DR_\F-\egrdg(X) \to \DR_{\F^h}-\egrdg(X^h)$$
induces an equivalence on the full sub-$\s$-category of crystals which are almost perfect
as $\OO_X$-modules and $\OO_{X^h}$-modules respectively.
\end{thm}

\textit{Proof.} We leave to the reader the precise construction of the analytification $\s$-functor
$$(-)^h : \QCoh(\F) \to \QCoh(\F^h)$$
along the same lines as the construction $\F \mapsto \F^h$ explained at the beginning of this section.
It is such that the following diagram of $\s$-categories commutes
$$\xymatrix{
\DR_\F-\egrdg(X) \ar[r] \ar[d] & \DR_{\F^h}-\egrdg(X^) \ar[d] \\
\OO_X-\grdg(X) \ar[r] & \OO_{X^h}-\grdg(X^h),
}$$
where the vertical $\s$-functor forget the module structures and only remember the graded $\OO$-module
structures.
The commutative square above is adjointable, the monads being here 
given by $\DR_\F\otimes_{\OO_X}-$ and $\DR_\F\otimes_{\OO_{X^h}}-$ respectively. Using 
the GAGA theorem for almost perfect $\OO$-modules (see \cite{porta}) we easily deduce 
the theorem. \hfill $\Box$ \\

\chapter{Leaf spaces and holonomy}\label{ch:leafspaces}

In this chapter we pursue our study of holomorphic derived foliations focusing now of the
notion of analytic leaf and analytic leaf spaces. We start by 
discussing the notion of leaves in a very general setting. In the specific case of 
quasi-smooth and rigid derived foliation we state natural conditions, principally 
related to the singularities of the derived foliation, insuring that a leaf space
exists in the analytic category. The existence of this leaf space is used to define the holonomy, 
as well as the relative homotopy type describing the change of topology of the leaves.
We provide two main applications of the existence of the leaf space, a Riemann-Hilbert
correspondence relating algebraic crystals and relative local systems, and 
an algebraic integrability result based on notion of Reed stability.  \\

We set $k=\CC$ throughout the whole chapter.

\section{Formal, algebraic and analytic leaves}

In Definition \ref{dII-11} we have seen the notion of formal leaves at a given global point. We remind
that, for any derived Artin stack $X$ locally of finite presentation, any
perfect derived foliation $\F \in \Fol^p(X/\CC)$, and any $x\in X(\CC)$, the
formal completion $\hat{\F}_x$ of $\F$ at $x$ is determined by a morphism of 
perfect dg-Lie algebras
$p : \ell_{X,x} \to \ell_{\F_x}$. The \emph{formal leaf} of $\F$ at $x$ is defined to be the
formal moduli problem determined by the $\CC$-dg-Lie algebra which is the fiber 
of $p$. We denote by $\hat{\cL}_x$ this formal moduli problem (or the corresponding
derived stack by the full embedding of Proposition \ref{pII-6}). In this section 
we will extend this formal notion to more global algebraic and analytic versions.  \\

As a start, we introduce the $\s$-category $\dSt^{fol}_\CC$, of \emph{foliated derived stacks}. 
Its objects are pairs $(X,\F)$, where $X$ is a derived stack and $\F$ is a derived foliation
on $X$. A morphism $(X,\F_X) \to (Y,\F_Y)$ is defined to be a pair 
$(f,u)$, where $f : X \to Y$ is a morphism of derived stacks, and $u : \F_X \to f^*(\F_Y)$ is
a morphism in $\Fol(X/\CC)$. Formally, $\dSt^{fol}_\CC$ is obtained as the Grothendieck construction
for the $\s$-functor of (\ref{functorFolonderivedstacks})
$$\Fol(-/\CC) : \dSt_\CC^{op} \longrightarrow \scat.$$

We let $X$ be a derived Artin stack locally
of finite presentation over $\CC$, and $\F \in \Fol^p(X/\CC)$
a perfect derived foliation over $X$, and consider $(X,\F)$ as an object in $\dSt^{fol}_\CC$. 
For any derived Artin stack $L$, we consider a morphism of foliated derived stacks
$$(f,u) : (L,*_L) \longrightarrow (X,\F).$$
Such a morphism can be referred to as a \emph{preleaf of $\F$}. Considering the induced morphisms
on tangent complexes, we find the following cartesian square of perfect complexes on $L$ (see Proposition \ref{cltII-1})
$$\xymatrix{
f^*(\mathbb{T}_\F) \ar[r] & f^*(\mathbb{T}_X) \\
\mathbb{T}_{f^*(\F)} \ar[r] \ar[u] & \mathbb{T}_L. \ar[u]
}$$

The morphism $u$ determines a section of $\mathbb{T}_{f^*(\F)}  \to \mathbb{T}_L$, and thus
a morphism of perfect complexes $\theta_u : \mathbb{T}_L \to f^*(\mathbb{T}_\F)$.

\begin{df}\label{dVI-1}
Let $X$ be a derived Artin stack locally 
of finite presentation over $\CC$ and $\F \in \Fol^{p}(X/\CC)$
a perfect derived foliation over $X$.
With the notations and assumptions as above, a preleaf $(f,u) : (L,*_L) \longrightarrow (X,\F)$
is called an \emph{(algebraic) leaf of $\F$}, if $L$ is non-empty and the 
induced morphism
$\theta_u : \mathbb{T}_L \to f^*(\mathbb{T}_\F)$ is an equivalence.
\end{df}

We can organise the leaves of $\F$ into an $\s$-category, by considering the 
full sub-$\s$-category in $\dSt_\CC^{fol}/(X,\F)$ of objects of the forms $(L,*_L)$,
with $L$ a derived Artin stack, 
and such that the structure morphism $(L,*_L) \to (X,\F)$ satisfies the condition of the previous
definition. A first important observation is that a morphism $(L,*_L) \to (L',*_{L'})$ between 
two leaves of $\F$ induces automatically an \'etale morphism $L \to L'$. As opposed to the 
formal notion of leaves of Definition \ref{dII-11}, \emph{the} leaf passing through  
a given point $x \in X(\CC)$ is not a well defined notion, as we could always modify
$L$ by an \'etale cover to obtain another, non-equivalent, leaf. Under some extra conditions
(e.g. for transversely smooth and rigid derived foliations),
there is a canonical choice, namely the ``maximal leaf" passing through a point, which will be
studied later with the introduction of the \emph{leaf space}. 
We have however the following comparison between the formal and global versions, 
showing that formally at each point $L$ is equivalent to the formal leaf defined
previously.

\begin{prop}\label{pVI-1}
We let $X$ be a derived Artin stack locally of finite presentation over $\CC$ and $\F \in \Fol(X/\CC)$
a perfect derived foliation over $X$. Let $(L,*_L) \to (X,\F)$ be a leaf of $\F$, and 
$x \in L(\CC)$ be a global point. Then, the formal completion $(L,x)^\wedge$ of $L$ at $x$
is, as a formal moduli problem, naturally equivalent to the formal leaf $\hat{\cL}_x$ as 
in Definition \ref{dII-11}.
\end{prop}

\textit{Proof.} First of all, any morphism $f : L \to X$, induces 
a morphism of derived foliations over $\Spec\, k$ 
(or of dg-Lie algebras according to Corollary \ref{ctII-2})
$x^*(0_L) \to x^*(0_X)$. The derived foliation $\F$ itself determines a second morphism
$x^*(0_X) \to x^*(\F)$. Now, a lift of $f$ to a morphism of foliated derived stacks
$(L,*_L) \to (X,\F)$ defines a commutative diagram in $\Fol(*)$
$$\xymatrix{
x^*(0_L) \ar[r] \ar[d] & x^*(0_X) \ar[d] \\
0 \ar[r] & x^*(\F).
}$$
Finally, if we require the lift to be a leaf in the sense of Definition \ref{dVII-1}, then 
the square above is moreover cartesian. This implies that $x^*(0_L)$ is equivalent to the
fiber of $x^*(0_X) \to x^*(\F)$. This proves the Proposition, as via the equivalence of Corollary \ref{ctII-2}, 
$x^*(0_L)$ corresponds to the formal completion of $L$ at $x$, and 
the fiber of $x^*(0_X) \to x^*(\F)$ is by definition the formal leaf at $x$.
\hfill $\Box$ \\

Analytic leaves are defined in exactly the same way. 

\begin{df}\label{dVI-2}
We let $X$ be a derived analytic Artin stack locally of finite presentation
over $\CC$ and $\F \in \Fol(X/\CC)$
a perfect derived foliation over $X$.
A morphism $(f,u) : (L,*_L) \longrightarrow (X,\F)$
of foliated derived analytic Artin stacks is called an 
\emph{analytic leaf of $\F$}, if $L$ is non-empty and the 
induced morphism
$\theta_u : \mathbb{T}_L \to f^*(\mathbb{T}_\F)$ is an equivalence of perfect complexes on $L$
\end{df}

Of course, by analytification, an algebraic leaf $(L,*_L) \to (X,\F)$ produces
an analytic leaf $(L^{h},*_{L^{h}}) \to (X^{h},\F^{h})$.  
Algebraic leaves exist very rarely, as opposed to analytic leaves. \\

\begin{ex}
\emph{Let $f: X \to Y$ be a morphism of derived Artin stacks locally of finite presentation. 
Then, for any global point $y \in Y(\mathbb{C})$, the canonical morphism 
$(f^{-1}(y), *_{f^{-1}(y)}) \to (X, f^*(0_Y)\simeq *_{X/Y})$ identifies the fibers of $f$ with 
leaves of the foliation $\F=f^*(0_Y) \in \Fol(X)$. The same argument applies in the analytic setting. 
In other words, when the derived 
foliation is integrable by the morphisms $f$, then the fibers of $f$ are naturally leaves of $\F$, as
expected.
} 
\end{ex}

\section{Existence of leaf space and holonomy}\label{sec:leafspaces}

In this section we will restrict ourselves to the following specific setting. We let $X$ be
a separated and connected (smooth) 
complex manifold. Let $\F \in \Fol(X/\CC)$ be a derived foliation on $X$
which is transversely smooth and rigid (i.e. satisfies Definition \ref{dII-12} in the analytic setting). We will study the existence of an analytic leaf space
for $\F$, namely we would like to contract all the leaves to points and produce the quotient space. 
This quotient does not exist in general, and we start by studying extra conditions on $\F$ that will ensure
its existence.

\subsection{Some conditions}

We start by discussing some conditions on $\F$ that will later on ensure the existence of a leaf space. \\

\textbf{Codimension $2$.} Because $\F$ is assumed transversely smooth and rigid, 
the cotangent complex $\LL_\F$ is a length two complex of vector bundles on $X$
of the form $\cN^{*}_\F \to \Omega_{X/\CC}^1$, where $\cN^*_\F$ is the conormal bundle of $\F$. 
The \emph{codimension $2$ condition} for $\F$ states that there exists a closed analytic 
subset $Z \subset X$, of codimension $2$ or higher, such that 
$\LL_\F$ restricts on $U=X-Z$ to a vector bundle. Equivalently, 
$\cN^*_\F \to \Omega_{X/\CC}^1$ is injective, and is a sub-bundle when restricted to $U$.
In particular, the perfect complex $\LL_\F$ is equivalent to a single coherent sheaf in degree 
zero, namely $\Omega_{X/\CC}^1/\cN^*_\F$, which is a vector bundle outside of $Z$.
The closed subset $Z$ is called the \emph{singular set of $\F$}. 

When the codimension $2$ condition is satisfied for $\F$, we will also say that 
\emph{$\F$ is smooth in codimension $2$}. \\

\textbf{Flatness.} For any point $x \in X$, we have seen that there exists a formal leaf 
$\hat{\cL}_x(\F)$ passing through $x$ (Definition \ref{dII-11}). This formal leaf is a
formal moduli problem given by a dg-Lie algebra whose underlying complex is $\TT_{\F,x}[-1]$, the 
fiber of the tangent complex of $\F$ at $x$, shifted by $-1$. As $\F$ is transversely smooth 
and rigid, $\TT_{\F,x}[-1]$ is cohomologicaly concentrated in degrees $[1,2]$, and therefore
corresponds to a pro-representable formal moduli problem (see \cite[Prop. 13.3.3.1]{SAG}): 
there exists
a pro-Artinian connective local cdga $\underline{A}="\lim" A_i$ such that 
$\hat{\cL}_x(\F)\simeq \Spf \underline{A}$, where $\Spf \underline{A}$ is the formal spectrum of $A$. 

The \emph{flatness condition} states that, if $A=\lim_i A_i$ is the realization of $\underline{A}$, then 
$H^i(A)=0$ for all $i<0$. Equivalently, this means that $\hat{\cL}_x(\F)$ is flat over $\CC$. 
As $A$ is naturally quasi-isomorphic to the Chevalley complex of the dg-Lie algebra $\TT_{\F,x}[-1]$, 
the flatness condition can also be stated as the following  vanishing for dg-Lie algebra cohomology
$$H^i(\TT_{\F,x}[-1],\CC)\simeq 0 \qquad \forall \, i<0.$$

\textbf{Local connectedness.} The local connectedness condition only makes sense under the codimension $2$
condition. Indeed, we have already seen that when $\F$ is smooth in codimension $2$, 
it is locally (analytically) integrable on $X$ (Theorem \ref{tV-1}). Therefore, there exists a basis $\B$ for the topology of $X$, 
and for all open $U\in \B$, a holomorphic map $f : U \to \CC^d$, such that $\F_{|U} \simeq f^*(0)$. 
The \emph{local connectedness coinditon} states that $B$ and the holomorphic maps $f$ above can 
be moreover chosen so that the non-empty fibers of $f$ are all connected. This condition is of
topological nature and is not always satisfied. When $d=1$, it is always satisfied as $f$ is smooth
in codimension $2$ and thus it is well known that the Milnor fibers of $f$ are always 
connected (by \cite{zbMATH03484390} the Milnor fiber is $(n-2-k)$ connected, where $n$ is the dimension
of $X$ and $k$ the dimension of the singular locus of $\F$).
When $d>1$, there are examples of holomorphic maps $f$ such that the fibers close to a singular fibers
are connected or not depending on the direction we approach the singular values. A simple example
is given in \cite{zbMATH03861545}, explictly given by  
$f : \CC^3 \to \CC^2$ given by $f(x,y,z)=(x^2-y^2z,y)$. \\

Our existence theorem below requires all three conditions, codimension $2$, flatness and local connectedness. Before going further we would like to make some comments on these conditions.

\begin{rmk}\label{rVI-1}
\begin{enumerate}
    \item \emph{As already reminded before the codimension $2$ condition implies that $\F$ is locally 
    integrable (Theorem \ref{tV-1}). So, for any given point $x\in X$, there is an open $x\in U$ and a holomorphic
    map $f : U \to \CC^d$ such that $f^*(0)\simeq \F_{|U}$. This implies that the formal 
leaf $\hat{\cL}_x(\F)$ is the formal completion of $f^{-1}(f(x))$ at the point $x$, where 
$f^{-1}(f(x))$ denotes here the derived fiber of $f$ passing through $x$. As formal completion are faithfully flat (see 
Corollary \ref{cpA-3-1-3}), the flatness of the formal moduli space $\hat{\cL}_x(\F)$ implies that $f^-{1}(f(x))$ is flat over $\Spec\, \CC$
locally at $x$. Therefore, flatness implies that the derived fibers
are all flat over $\CC$, and thus that $f$ must be a flat morphism. The converse is obviously true. 
Therefore, we see that under the codimension $2$ condition, flatness is equivalent to the fact that the
all the local holomorphic maps integrating $\F$ are actually flat holomorphic maps.}
    \item \emph{If $\F$ is moreover smooth, and thus is given by an integrable subbundle of $\TT_X$, then 
all three conditions above are all automatic. Indeed, the local holomorphic maps $f$ must be smooth
and thus their fibers will be themselves smooth. In particular, the local connectedness is true because the maps 
$f$ are submersions. Finally, the formal leaves are the formal completions of smooth sub-varieties
(the fibers of $f$), and are thus all equivalent to formal affine spaces $\widehat{\mathbb{A}}^{n-d}$
(where $n=dim X$), and so they are flat.} 
\end{enumerate}
\end{rmk}

\subsection{Existence theorem}\label{subsec_existencethm}
The following result establishes, under appropriate hypotheses, the existence of a \emph{leaf space}. 

\begin{thm}\label{tVI-1}
Let $X$ be a smooth, connected and separated 
complex manifold and $\F \in \Fol(X/\CC)$ be a derived foliation
which is flat, smooth in codimension $2$ and locally connected. 
There exists a smooth and $1$-truncated Deligne-Mumford analytic stack $X\sslash \F$, together with 
a holomorphic morphism $\pi : X \to X\sslash\F$ with the following properties.

\begin{enumerate}
    \item The morphism $\pi$ is flat, surjective, and we have $\pi^*(0)\simeq \F$.
    \item The Deligne-Mumford stack $X\sslash \F$ is quasi-separated, effective 
    and is of dimension $d$ the codimension of $\F$.
    \item For any global point $x : * \to X\sslash \F$ the fiber of $\pi$ at $x$
    $$\pi^{-1}(x):=* \times_{X\sslash \F}X$$
    is a connected and separated analytic space.
    \item For any global point $x : * \to X\sslash \F$, the isotropy group $H_x$
    of $X\sslash\F$ at $x$ acts freely and properly on $\pi^{-1}(x)$.
    \item The morphism $\pi$ is relatively connected in the sense of topos theory: the pull-back
    $\s$-functor $\pi^{-1} : St(X\sslash \F) \to St(X)$ is fully faithful when restricted to 
    $0$-truncated objects (i.e. sheaves of sets).
\end{enumerate}
\end{thm}

Before giving a proof of Theorem \ref{tVI-1} some of the terms in $(2)$ require some  explanations. 
First of all, \emph{quasi-separated} here means that the diagonal of $X\sslash \F$ is a separated 
morphism. In more concrete terms, if the stack $X\sslash \F$ is presented as an \'etale groupoid
object $G \rightrightarrows U$, with $U$ a disjoint union of Stein manifolds, then 
the smooth analytic space $G$ must be separated. 
\emph{Effectiveness} of $X\sslash \F$ means that the action of $G$ on $U$, by germs of biholomorphic
maps, is faithful 
(see \cite[Def. 2.4]{zbMATH07082639}). Equivalently, the groupoid $G$ is a sub-groupoid of the
Haefliger’s groupoid of germs of biholomorphic isomorphisms.  As shown in 
\cite[Def. 2.9]{zbMATH07082639}, any 
Deligne-Mumford stack $Y$
possesses a universal effective quotient $Y^{eff}$. We will see later than 
$X\sslash \F$ is itself the effective quotient of a bigger stack which is pro-\'etale over it, given
as the relative pro-homotopy type of $X$ over $X\sslash\F$ (see 
\S \ref{sec:relhomotopy}). \\

Finally, Theorem \ref{tVI-1} allows to set the following useful terminology.

\begin{df}\label{dVI-3}
Let $X$ be a smooth complex manifold and $\F$ be a derived foliation on $X$, and 
assume the conditions of Theorem \ref{tVI-1} hold.

\begin{enumerate}
    \item The stack 
$X\sslash \F$ is called the \emph{leaf space of the derived foliation $\F$}.
    \item For any $x\in X$, the quotient analytic space
    $\pi^{-1}(\pi(x))/H_{\pi(x)}$ is called the \emph{maximal leaf passing through $x$}.
    \item The $H_{\pi(x)}$-covering $\pi^{-1}(\pi(x)) \to \pi^{-1}(\pi(x))/H_{\pi(x)}$ is called
    the \emph{holonomy covering}, and the group $H_{\pi(x)}$ the \emph{holonomy group} of $\F$ at $x$. 
\end{enumerate}
The maximal leaf will be denoted by $\cL_{x}^{max}(\F)$, and the holonomy group $Hol_\F(x)$.
The holonomy covering of $\cL_{x}^{max}(\F)$ will be denoted by $\widetilde{\cL}_{x}^{max}(\F)$.
\end{df}

We now prove Theorem \ref{tVI-1}, except for $(3)$: the connectedness of $\pi^{-1}(x)$ will be 
proved at the very end of the section, as a consequence of an explicit description of the underlying 
groupoid of $X\sslash \F$. \\

\textbf{Construction of $X\sslash \F$.} We let $\B$ be the set of all open subsets $U \subset X$
such that there exists a flat holomorphic function $f_U : U \to \CC^d$ with connected fibers, 
and with $f_U^*(0)\simeq \F_{|U}$. As $f_U$ is flat, it has an open image $V \subset \CC^d$, 
and thus factors $U \to V \subset \CC^d$.
For each $U \in \B$ we fix once for all $f_U : U \to V \subset \CC^d$, 
where $V$ is open, $f_U$ is flat, surjective with connected fibers, and
such that $f_U^*(0)\simeq \F_{|U}$. Our conditions on $\F$ ensure that the set $\B$
forms a basis for the topology of $X$. 

Before going further, we note that $\B$ consists of all opens $U \subset X$ with 
the required properties, so the only choices that have been made here are the choices of the functions $f_U$. 
In order to understand the impact of these choices we will need the following lemma. By the codimension $2$ condition
we know that the perfect complex $\LL_\F$ is a coherent sheaf, denoted by $\Omega^1_\F \simeq H^0(\LL_\F)$.
The de Rham differential for $\F$ therefore induces a morphism 
$dR_\F : \OO_U \to \Omega_\F^1$ of sheaves on $U$.
Because $f_U^*(0)\simeq \F$, this de Rham differential can be identified with 
the relative de Rham differential $dR_{U/V} : \OO_U \to \Omega_{U/V}^1$.
In particular, the composite morphism $f^{-1}(\OO_V) \to \OO_U \to \Omega_\F^1$
is zero.

\begin{lem}\label{ltVI-1}
The natural morphism of sheaves of rings on $U$
$$f^{-1}(\OO_V) \to Ker(dR_\F)$$
is an isomorphism.
\end{lem}

\textit{Proof of Lemma \ref{ltVI-1}.} Because $f$ is flat and surjective the induced morphism 
$f^{-1}(\OO_V) \to \OO_U$ is a monomorphism of sheaves. Therefore, $f^{-1}(\OO_V) \to Ker(dR_\F)$
is also mono. It remains to show that a local section of $Ker(dR_\F)$, locally descend to a local 
holomorphic function on $V$. 

\begin{slem}\label{sltVI-1}
Let $f : X \to Y$ be a flat surjective morphism of complex manifolds which is smooth
in codimension $2$. If an application $u : Y \to \CC$
is such that $f\circ u : X \to \CC$ is holomorphic then $u$ is holomorphic.
\end{slem}

\textit{Proof of sublemma \ref{sltVI-1}.} Let $U_0 \subset X$ be the open on which 
$f$ is smooth, so that $Z=X-U_0$ is a closed analytic subset of codimension 2 or higher.
The image of $U_0$ by $f$ is an open $V_0 \subset Y$, because $f$ is flat, which is dense
by Sard's theorem. As $f\circ u : U_0 \to V_0$ is a holomorphic submersion, it has local sections, and
thus we deduce that $u$ is holomorphic when restricted to $V_0$. Moreover, $f$ being flat and surjective
implies that the topology of $Y$ is the quotient of the topology on $X$ induced by the map $f$. This shows that 
$u$ is also continuous on $Y$. The application $u$ is continuous on $Y$ and 
holomorphic on a dense open, and so is holomorphic on the whole $Y$. \hfill $\Box$ \\

Now, let $s$ a be local section of $Ker(dR_\F)$ defined on an open $U' \subset U$. 
By possibly shrinking $U'$ and using the local connectedness condition, we can assume that the restriction of $f_U$
on $U'$ has connected fibers. The codimension $2$ condition implies that the
canonical map from derived de Rham cohomology of $U$ over $V$ coincides with the naive
de Rham cohomology in degree $0$ (see Corollary \ref{clcpII-10}). 
In particular, $dR_\F(s)=0$ implies that $s$ determines
a natural element in $H^0_{dR}(U'/V)$, the derived de Rham cohomology of $U'$ relative to $V$.
We can write this as $s\in |\DR(U/V)|(U')$, where $\DR(U/V)$ is the sheaf on $U$ of graded mixed
cdga's of relative de Rham theory of $U$ over $V$. The sheaf $|\DR(U/V)|$ is naturally a 
module over $f^{-1}(\OO_V)$. Moreover, if $y \in V$ is a point, corresponding to a morphism
of sheaves of rings $f^{-1}(\OO_V) \to \CC$,
then we have an equivalence of sheaves of cdga's
$$|\DR(U/V)|\otimes_{f^{-1}(\OO_V)}\CC \simeq |\DR(U/V)\otimes_{f^{-1}(\OO_V)}\CC|\simeq
j_*(|\DR(U_y/\CC)|),$$
where $j: U_y \hookrightarrow U$ is the inclusion of the fiber $U_y=f^{-1}(y)$, and where
the first equivalence follows from the fact that $\CC$ is a perfect $f^{-1}(\OO_V)$-module (so that 
the tensor operation $\otimes_{f^{-1}(\OO_V)}\CC$ commutes with limits and with the realization).
To summarize, the restriction of $s$ on $U'_y$, lies in derived de Rham cohomology $H^0_{dR}(U'_y/\CC)$. 
Now $U'_y$ is connected, so that by the comparison between derived de Rham cohomology and
Betti cohomology (see Proposition \ref{pI-6} which also holds in the
complex analytic setting), we have that $H^0_{dR}(U'_y/\CC) \simeq \CC$, and thus 
the restriction of $s$ on $U'_y$ is a constant function. As $f : U' \to V$ is flat with connected fibers, 
the function $s$ descends to an application defined on the open $V'=f(U')\subset V$, which is moreover
holomorphic by 
the sublemma \ref{sltVI-1}. This shows that $s$ comes from a local section of $f^{-1}(\OO_V)$ 
as required. \hfill $\Box$ \\

Lemma \ref{ltVI-1} has the following important corollary, showing that the choices of the
local morphisms $f : U \to V$ are essentially unique.

\begin{cor}\label{cltVI-1}
Let $f : U \to V$ and $g : U \to V'$ 
be two flat and surjective holomorphic morphisms of smooth manifolds
with $f^*(0) \simeq g^*(0)$. Assume that $V'$ is isomorphic to an open in $\CC^d$.
If the fibers of $f$ are connected Then, there exists a 
unique \'etale holomorphic morphism $\alpha : V \to V'$
such that $\alpha f = g$.
\end{cor}

\textit{Proof of Corollary \ref{cltVI-1}.} 
We apply Lemma \ref{ltVI-1} to local coordinate functions of $\CC^d$ restricted on $V'$,
and get a unique factorization $\alpha : V \to V'$ such that $\alpha f = g$. Because $f$ is
surjective, $\alpha$ must be unique. Finally, as $f^*(0)\simeq g^*(0)$, comparing cotangent 
complexes show that $\alpha$ must be \'etale. \hfill $\Box$ \\

We now come back to the base $\B$ of the topology. The set $\B$ is ordered by the inclusions and 
as such will be considered as a category. We define a functor
$$\B \to \CC Man^d_{et},$$
from $\B$ to the category of complex manifolds of dimension $d$ and \'etale morphisms between them. 
On objects, the 
functor $\phi$ sends $U$ to the space $V$, the base of the morphism $f_U : U \to V$. If $U' \subset U$
is an inclusion of opens in $\B$, with morphisms $f_U : U \to V$ and $f_{U'} \to V'$.
Corollary \ref{cltVI-1} implies that there exists a unique \'etale factorization
$\alpha : V' \to V$ such that $f_{U'}\alpha =(f_{U})_{|U'}$. This provides for each open inclusion
$U' \subset U$ in $\B$, a natural morphism $\phi(U')=V' \to \phi(U)=V$ in $\CC Man^d_{et}$. This defines
the functor $\phi$.

Smooth Deligne-Mumford stacks are stable arbitrary colimits of \'etale morphisms. 
We therefore consider the colimit of the functor $\phi$, and take its universal effective 
quotient defined in \cite[Def. 2.9]{zbMATH07082639}
$$X\sslash \F := (\mathrm{colim}_{U \in \B}V)^{eff}.$$
By construction $X\sslash \F$ is an effective smooth Deligne-Mumford stack of dimension $d$. 
Moreover, as $\B$ is a basis
for the topology of $X$, we have that $X \simeq \mathrm{colim}_{U \in B}U$. The local morphisms $f_U : U \to V$
therefore define the projection
$$\pi : X \simeq \mathrm{colim}_{U \in B}U \to \mathrm{colim}_{U \in \B}V \to (\mathrm{colim}_{U \in \B}V)^{eff}= X\sslash \F.$$

\bigskip

\textbf{The morphism $\pi$ is flat surjective and $\pi^*(0_{X\sslash \F})\simeq \F$.}  By construction, the morphism 
$\pi : X \to X\sslash \F$, restricted to an open $U \in \B$ factorizes as
$\pi_{|U} : \xymatrix{U \ar[r]^-{f_U} & V \ar[r]^-{p_V} & X\sslash \F}$, where 
$p_V$ is the composition of the natural morphism $V \to \mathrm{colim}_{U \in B}V$, 
with the canonical projection to the effective quotient $\mathrm{colim}_{U \in B}V\to X\sslash \F$.
The morphism $p_V$ is etale, and $f_U$ is flat, so we see that $\pi_{|U}$ is a flat morphism
for all $U \in \B$. As $\B$ is a basis for the topology of $X$, this
shows that $\pi$ is flat. 

To see that $\pi$ is surjective, it is enough to notice that the universal effective
quotient map is always surjective (because it does not the change the space of objects), and
moreover that $\coprod_{U \in \B} V \to \mathrm{colim}_{U \in \B}V$ is surjective. As a result, we see that 
the composition
$\coprod_{U \in \B}U \to X \to X\sslash/\F$ is a surjective morphism, 
and thus $\pi$ is surjective.

Finally, let us consider the pullback foliation $\pi^*(0_{X\sslash \F})$. We thus have two 
derived foliations on $X$, $\pi^*(0_{X\sslash \F})$ and $\F$. For all $U\in \B$, we know that $f_U^*(0)\simeq \F_{|U}$, 
but as $p_V : V \to X\sslash \F$ is \'etale, we have that $f_U^*(0)\simeq \pi^*(0)_{|U}$. Therefore,
the two derived foliations $\pi^*(0)$ and $\F$ coincide on each open $U$. This implies in particular
that the corresponding differential ideals $K_\F$ and $K_{\pi^*(0)}$, obtained as their truncations
(Definition \ref{dII-14}), agree on each $U$, and thus globally agree on $X$. This means that $\F$ and $\pi^*(0)$ are two
transversely smooth and rigid derived enhancements of the same differential ideal $K_\F$. It follows
from the uniqueness statement of derived enhancements in codimension $2$ (Proposition \ref{pV-1}) that 
$\F \simeq \pi^*(0)$ as required. \\

\textbf{Quasi-separatedness.} The quasi-separatedness property for $X\sslash \F$ follows from the general lemma 
below.

\begin{lem}\label{lVI-2}
Any smooth effective Deligne-Mumford stack is quasi-separated.
\end{lem}

\textit{Proof.} Let $Y$ be an effective Deligne-Mumford stack and $U \to Y$ an \'etale atlas
with $U$ a disjoint union of Stein manifolds (hence separated). The nerve of $U\to Y$ defines
an \'etale groupoid $G=U\times_Y U$ acting on $U$. The quasi-separatedness of $Y$ is equivalent to the
fact that the analytic space $G$ is separated. 

Let $Haf(U)$ be the \emph{Haefliger groupoid} on $U$. As a set, $Haf(U)$ consists of triplets 
$(x,y,u)$, with $x$ and $y$ points in $U$ and $u$ is a germ of biholomorphic isomorphism 
from $x$ to $y$. A basis for the topology on $Haf(U)$ is defined as follows. 
Fix $V \subset U$ and $W \subset U$ two opens, and $u : V \simeq W$ and isomorphism. 
Associated to $V,W$ and $u$ we have a subset $Hal(U)_{V,W,u} \subset Haf(U)$
consisting of all points $(x,u(x),u_x)$, where $x \in V$, $u(x)\in W$ its image and
$u_x$ the germ of isomorphism at $x$ from $x$ to $u(x)$. By definition, a basis for the topology
consists of all $Hal(U)_{V,W,u}$ where $(V,W,u)$ varies as above. 
The groupoid structure on $Haf(U)$ is the obvious one: the source (resp. target) 
map sends $(x,y,u)$ to $x$ (resp. $y$), and the composition is given by composing germs of 
isomorphisms. Since we are working in the holomorphic context, the space $Haf(U)$ is automatically 
separated, 
as two holomorphic maps agreeing on an open agree globally (on the connected components
meeting this open). 

Finally, there is a canonical morphism of groupoids over $U$, $\phi : G \to Haf(U)$. 
To a point $u \in 
G$, 
with source $x$ and target $y$, we associate $\phi(u)$ the germs of isomorphisms
defined by the diagram 
$$\xymatrix{U & \ar[r] \ar[l] G & U},$$
which, $G$ being an \'etale groupoid, induces an isomorphism on germs
$$\xymatrix{U_x & \ar[r]^-{\simeq} \ar[l]_-{\simeq} G_u & U_y}.$$
The isomorphism $U_x \simeq U_y$ obtained this way is $\phi(u)$. Now, $Y$ being effective
is equivalent to the fact that the morphism $\phi$ is a monomorphism 
(see \cite[Def. 2.4]{zbMATH07082639}). So 
$G \hookrightarrow Haf(U)$ is here an injective holomorphic map. As $Haf(U)$ is separated, this
implies that $G$ is also separated, as claimed. \hfill $\Box$ \\

\textbf{Actions of isotropy groups.} Let $x : * \to X\sslash \F$ be a point in 
$X\sslash \F$, and $H_x=aut(x)$ be the isotropy group of the DM-stack $X\sslash\F$ at $x$. We have a commutative diagram with cartesian squares
$$\xymatrix{ \pi^{-1}(x) \ar[d] \ar[r] & \pi^{-1}(BH_x)\simeq [\pi^{-1}(x)/H_x] \ar[r] \ar[d] & X \ar[d] \\
\bullet \ar[r] & BH_x \ar[r] & X\sslash \F.}$$
The canonical morphism $BH_x \to X\sslash \F$ is a monomorphism, and thus $[\pi^{-1}(x)/H_x] \to X$
is also a monomorphism. This implies that the stack $[\pi^{-1}(x)/H_x]$ is homotopically 
$0$-truncated and thus
is an analytic space. Moreover, as $X$ is separated, so is $[\pi^{-1}(x)/H_x]$. Therefore, 
$H_x$ acts freely and properly on $\pi^{-1}(x)$. \\

\textbf{Relative connectedness.} Let us denote by $Sh(Z)\subset St(Z):=\lim_{S \to X}St(S)$ the full sub-category of $0$-truncated
objects (for some stack $Z$), i.e. the category of sheaves of sets. By construction, the
functor $\pi^{-1} : Sh(X\sslash \F) \to Sh(X)$ is obtained as a functor on limit categories
$$\lim_{U\in \B} f_U^{-1} : Sh(X\sslash \F)\simeq \lim_{U\in B}Sh(V) \to \lim_{U\in B}Sh(U)\simeq Sh(X).$$
Therefore, in order to prove that $\pi^{-1}$ is fully faithful, it is enough to prove that 
$f_U^{-1} : Sh(V) \to Sh(U)$ is fully faithful. But this follows easily from the fact that 
$f_U$ is a flat surjective morphism with connected fibers, as shown by the next lemma.

\begin{lem}\label{lVI-3}
Let $f : U \to V$ be an open and surjective morphism such that the topology of $V$ is the
quotient of the topology of $U$ (i.e. a subset $T \subset V$ is open if and only if $f^{-1}(T) \subset U$
is open). Then, if $f$ has connected fibers, for any sheaf of sets $E$ on $V$ the canonical morphism
$$E \to f_*f^{-1}(E)$$
is an isomorphism.
\end{lem}

\textbf{Proof.} Let us represent $E$ by its espace étalé $p : E \to V$. As the statement is local on $V$, 
it is enough to show that the natural map $E(V) \to f^{-1}(E)(U)$ is bijective. Now, $E(V)$ 
is the set of continuous sections of $E \to V$, and $ f^{-1}(E)(U)$ is the set 
of continuous maps $s : U \to E$ such that $ps=f$. The map $E(V) \to f^{-1}(E)(U)$
simply sends a section $s : V \to E$ to $sf : U \to E$. As $f$ is surjective this map is clearly
injective. Moreover, if $s : U \to E$ is an element $f^{-1}(E)(U)$, 
we can restrict $s$ to the fiber of $f$ at a point $v\in V$. This provides a continuous map
$f^{-1}(v) \to E_v$, where $E_v$ is the fiber of $E$ at $v$ endowed with the discrete topology. 
As $f^{-1}(v)$ is connected we see that the restriction of $s$ along all the fibers 
of $f$ is a constant morphism. Because $f$ is a topological quotient, this implies that 
$s : U \to E$ descend to a continuous map $V \to E$, and thus $s$ is the image of an element in 
$E(V)$. \hfill $\Box$ \\

\textbf{An \'etale groupoid description.} It is possible to describe the stack $X\sslash \F$
by a (more or less explicit) \'etale groupoid acting on $W:=\coprod_{U\in \B}V$. Recall that we have chosen,
for any $U \in \B$, a smooth surjective holomoprhic map with connected fibers $f_U : U \to V$
which integrates $\F_{|U}$. We consider all pairs of embeddings
$$U_1 \supset U_0 \subset U_2$$
of opens $U_i \in \B$. Because of Corollary \ref{cltVI-1}, there is a corresponding pair of 
\'etale morphisms
$$\xymatrix{V_1 & V_0 \ar[r]^-s  \ar[l]_-b & V_2.}$$
Therefore, for any $v\in V_0$, we deduce a germ of isomorphism $\alpha(v)$
between $(V_1,s(v))$ and $(V_2,b(v))$. When $v$ runs in $V_0$, we 
get this way a family of elements
$\{(s(v),b(v),\alpha(v))\}_{v\in V_0}$ in $Haf(W)$, the Haefliger groupoid on $W=\coprod_{U\in \B}V$ (see the proof of \ref{lVI-2} for the definition of $Haf$).

\begin{lem}\label{lVI-4}
Let $G \rightarrow W$ be the \'etale groupoid obtained as the nerve of the canonical
morphism $W \to X\sslash \F$. Then, $G$ is isomorphic, as a groupoid over $W$, 
to the subgroupoid of $Haf(W)$ generated by the elements $(s(v),b(v),\alpha(v))$
for all possible choices of $U_1 \supset U_0 \subset U_2$ in $\B$ and of $v\in V_0$.
\end{lem}

\textit{Proof.} This follows from a more general formula for the effective part of 
a colimit of \'etale morphisms. Let $V_* : I \to \CC Man^d_{et}$ be a diagram of smooth 
complex manifolds of dimension $d$ with \'etale transition map. Let $X=(\mathrm{colim} U_i) ^{eff}$ and
$W=\coprod V_i \to X$ be the canonical map. Then, the nerve $W \to X$ is isomorphic, 
as a groupoid over $W$, to the subgroupoid of $Haf(W)$ generated by all the germs of isomorphisms generated
by diagrams of the form $\xymatrix{V_j & V_i \ar[r] \ar[l] & V_k}$, image of 
$\xymatrix{j & i \ar[r] \ar[l] & k}$ in $I$ by $V_*$. This last general fact is obvious
according to the universal properties of colimits and of the effective part. \hfill $^\Box$ \\

\textbf{Connectedness of the leaves and their holonomy coverings.} We prove here Theorem \ref{tVI-1} (3). We first prove that the leaves are connected. For this, 
let $x$ and $y$ be two points in $X$ such that $\pi(x)$ and $\pi(y)$ are isomorphic in $X\sslash \F$.
By Lemma \ref{lVI-4} we know that there exists a finite sequence of open inclusions in $\B$
$$U_1 \supset U_{1,2} \subset U_2 \supset U_{2,3} \subset U_3 \dots U_{n-1} \supset U_{n-1,n} \subset U_n,$$
with points $x_{i,j} \in U_{i,j}$, and such that 
\begin{enumerate}
    \item $x\in U_1$ and $y\in U_n$
    \item $f_1(x_{1,2}) = f_1(x)$ and $f_n(x_{n-1,n})=f_n(y)$
    \item $f_i(x_{i-1,i})=f_i(x_{i,i+1})$
\end{enumerate}
where $f_i : U_i \to V_i$ is our chosen map integrating $\F_{|U_i}$. As the fibers of $f_i$ are all connected,
it is possible to find a path $\gamma_{i}$ joining $x_{i-1,i}$ with $x_{i,i+1}$ and such that 
$f_i(\gamma_{i})$ is constant (for all $i\neq 1,n$). 
In the same manner, we can choose paths $\alpha$ from $x$ to $x_{1,2}$ with $f_1(\alpha)$ constant, 
and $\beta$ from $x_{n-1,n}$ to $y$ with $f_n(\beta)$ constant. The concatenation 
of the $\gamma_i$ with $\alpha$ and $\beta$ provides a path $\gamma : [0,1] \to X$, from
$x$ to $y$ and such that for all $t \in [0,1]$ the object $\pi(\gamma(t))$ are all isomorphic. 
By definition, this means that $\gamma$ factors through the image of the monomorphism
$\cL^{max}_x(\F) \hookrightarrow X$, and thus it is a continuous path in 
$\cL^{max}_x(\F)$ joining $x$ and $y$. This implies that $\cL^{max}_x(\F)$ is path connected. \\

It remains to show that the $H_x$-covering $\widetilde{\cL}^{max}_x(\F) \to \cL_x^{max}(\F)$
is connected or, equivalently, that the corresponding classifying morphism 
$\pi_1(\cL_x^{max}(\F),x) \to H_x$ is a surjective morphism of groups. This is proven using the
same argument as before with $x=y$. Indeed, an element $h \in H_x$ is determined by a 
finite sequence of open inclusions in $\B$
$$U_1 \supset U_{2,3} \subset U_2 \supset U_{2,3} \subset U_3 \dots U_{n-1} \supset U_{n-1,n} \subset U_n,$$
with points $x_{i,j} \in U_{i,j}$, and such that 
\begin{enumerate}
    \item $x \in U_1\cap U_n$
    \item $f_1(x_{1,2}) = f_1(x)$ and $f_n(x_{n-1,n})=f_n(x)$
    \item $f_i(x_{i-1,i})=f_i(x_{i,i+1})$.
\end{enumerate}

We have seen that there exists a loop $\gamma$ in $X$, pointed at $x \in X$, and
such that $\pi(\gamma(t))$ is independent of $t \in[0,1]$, up to isomorphism, in $X\sslash \F$. 
The image of $\gamma$ by the morphism $\pi_1(\cL_x^{max}(\F),x) \to H_x$
is $h$ by construction. \\

We conclude this \S \, with the following statement, concerning the functoriality of the 
leaf space construction, and thus of its characterization by a universal property.

\begin{prop}\label{pVI-2}
Let $f : X \to Y$ be a morphism between smooth, connected and separated complex manifolds
and $\F_Y \in \Fol(Y)$ and $\F_X:=f^*(\F_Y)$. We assume that both derived foliations
$\F_Y$ and $\F_X$ are transversely smooth, flat, smooth in codimension $2$, and locally 
connected. Then, there exists a morphism, unique up to a unique isomorphism,
$u : X\sslash \F_X \to Y\sslash \F_Y$ such that the following diagram 
commutes up to isomorphism
$$\xymatrix{
X \ar[r]^f \ar[d]_-{\pi_X} & Y \ar[d]^-{\pi_Y} \\
X\sslash \F_X \ar[r]_-{u} & Y \sslash \F_Y.
}$$
\end{prop}

\textit{Proof.} Left to the reader. \hfill $\Box$ \\

\section{Relative homotopy theory}\label{sec:relhomotopy}

Let $f : \cT \to \cT'$ be a geometric morphism of $\s$-topoi. The inverse image $f^{-1} : \cT' \to \cT$
is exact by definition, and thus admits a pro-left adjoint\footnote{We recall that this means  that $Pro(f^{-1}):Pro(\cT') \to Pro(\cT)$ has a genuine left adjoint, and $f_\sharp$ is then the composition of $\cT \to Pro(\cT)$ with such a genuine left adjoint. } 
$$f_\sharp : \cT \to Pro(\cT').$$

\begin{df}\label{dVI-4}
The pro-object $f_\sharp(*) \in Pro(\cT')$ is called the \emph{relative homotopy type of $\cT$ over $\cT'$}.
It is denoted either by $\cH_{f}$ or $\cH_{\cT/\cT'}$.
\end{df}

By definition, the relative homotopy type satisfies the following universal property: for any 
object $X \in \cT'$, we have 
$Map_{Pro(\cT')}(\cH_f,X) \simeq Map_{\cT}(*,f^{-1}(X))$. This must be understood as the fact that 
$\cH_f$ pro-represents non-abelian cohomology of $\cT$ with coefficients coming from $\cT'$. 
When $\cT'=\widehat{*}=\Top$ is the punctual $\s$-topos, the pro-object $\cH_f$ has been introduced in \cite{TOVESegal}
as the \emph{shape of $\cT$} (see also \cite[\S 7.1.6]{htt}, \cite[\S 2]{HoyoisHGT} and \cite[\S 3]{mvolpe}). \\

We will restrict ourselves to the case where $\cT$ and $\cT'$ are $\s$-topoi of stacks on
nice topological spaces. In particular, all these spaces will be locally 
CW complexes of finite dimension, and thus locally contractible and locally of finite homotopical 
dimension. Therefore, 
all our $\s$-topoi appearing below will be hypercomplete, making it easier to check
that morphisms are equivalences by just looking at fibers at each points.\\

Let $f : X \to Y$ be a 
holomorphic morphism between analytic Stein varieties. In this case, $X$ and $Y$ are
CW complexes of finite dimension by Hironaka (see \cite{HironakaTriang}), and thus are locally contractible, and
locally of finite homotopical dimension.
We still denote by $f : \St(X) \to \St(Y)$ the induced geometric morphism between the corresponding $\infty$-topoi of stacks on $X$ and $Y$.
The pro-object $\cH_f \in Pro(\St(Y))$ is in general non-constant, but has the following 
remarkable properties.

\begin{prop}\label{pVI-3}
With the notation above, the following assertions hold.
\begin{enumerate}
 \item Let $y\in Y$ be a point, then the fiber of $\cH_f$ at $y$ is pro-constant
    as a pro-object in $\Top$, equivalent to $Sing(f^{-1}(y))$, the homotopy type 
    of the fiber of $f$ at $y$.
    \item If $f$ is a smooth morphism, then the pro-object 
    $\cH_f$ is pro-constant.
    \item Suppose that $f$ is topologically a locally trivial fibration 
    with fiber $F_0$, then $\cH_f$ is pro-constant, and a locally constant stack
    stack with fiber the homotopy type $Sing(F_0)$.
\end{enumerate}
\end{prop}

\textit{Proof.} $(1)$ Let $i : \{y\} \hookrightarrow Y$ the inclusion map. Let $X_y:=f^{-1}(y)$ be the
fiber of $f$ at $y$, $j : X_y \hookrightarrow X$ the inclusion and $f_y : X_y \to *$ the structure map. 
For any space $F \in \St(*)=\Top$, we have
$$Map_{Pro(\Top)}(i^*(\cH_f),F) \simeq Map_{Pro(\St(Y))}(\cH_f,i_*(F)) \simeq Map_{\St(X)}(*,f^{-1}i_*(F)).$$
But, by proper base change, we have a canonical equivalence $f^{-1}i_*(F) \simeq (f_y)_*j^{-1}(F)$. 
Therefore, we have
$$Map_{\St(X)}(*,f^{-1}i_*(F)) \simeq Map_{\St(X_y)}(*,j^{-1}(F)) \simeq Map_{Pro(\Top)}(\cH_{f_y},F).$$
However, $X_y$ is Stein analytic space, and thus is a CW complex of finite dimension by \cite{HironakaTriang}. 
Using \cite[Cor. 3.2]{mvolpe} combined with \cite[Cor. 7.2.3.7]{htt}, we have that 
$\cH_{f_y}$ is pro-constant and equivalent to the usual homotopy type $Sing(X_y)$. As a final result, we have
functorial equivalences
$$Map_{Pro(\Top)}(i^*(\cH_f),F) \simeq Map_{\Top}(Sing(X_y),F),$$
which shows $(1)$ by Yoneda. \\

$(2)$ This is proven in \cite[Lem. 3.2]{mvolpe}, keeping in mind, as in the previous point,
that Stein analytic spaces are essential in the sense of  \cite{mvolpe} (as they 
are of finite covering dimension and locally contractible). \\

$(3)$ The statement being local on $Y$, we may assume that $X=Y\times F_0$ with $F_0$ another Stein 
space, and $f : Y \times F_0 \to Y$ the natural projection. Let $K=Sing(F_0)$ and
$\underline{K}_Y \in \St(Y)$ be the constant stack with fiber $K$ on $Y$. There exists a canonical
morphism $\cH_f \to \underline{K}_Y$ defined as follows. Such a morphism, by previous $(2)$,  is in fact 
a morphism of stacks on $Y$ and not of pro-stacks. By definition of $\cH_f$, it is thus  determined
by a global section $* \to \underline{K}_X$ of the constant stack $\underline{K}_X$ with fiber K on $X$.
As $X$ is a CW complex, we know that the stack $\underline{K}_X$ simply sends
an open $U \subset X$ to $Map_{\Top}(Sing(U),K)$. Therefore, a global section 
of $\underline{K}_X$ is determined by a morphism of simplicial sets $Sing(X) \to Sing(F_0)$, 
which we take as being the natural projection. This defines our morphism $\cH_f \to \underline{K}_Y$. 
To 
check this is an equivalence we use that $\St(Y)$ is hypercomplete, and thus it is enough
to check it is a fiberwise equivalence. The point $(1)$ therefore implies the result. \hfill $\Box$ \\

The relative homotopy theory mentionned above as several interesting applications, and one of them
is the definition of a \emph{monodromy groupoid}, bigger than the leave space $X\sslash\F$ 
of Definition \ref{dVI-3}. It is gather in the following definition, and will be 
an essential object for our Riemann-Hilbert correspondence studied in the next section.

\begin{df}\label{dVI-6}
Let $X$ be a smooth complex manifold and $\F \in \Fol(X)$ a derived foliation
satisfying the condition of Theorem \ref{tVI-1}. The \emph{monodromic leaf space}
of $X$ with respect to $\F$ is the relative homotopy type
$$\cH_\pi \in Pro(X\sslash\F),$$
where $\pi : X \to X\sslash\F$ is the natural projection. It is denoted by 
$\widetilde{X\sslash\F}$.
\end{df}

The monodromic leaf space $\widetilde{X\sslash\F}$ can be understood as a projective system of étale
morphism of analytic stacks over $X\sslash\F$, and thus as a smooth pro-Deligne-Mumford stack.
The major difference between $\widetilde{X\sslash\F}$ and $X\sslash\F$ is that the former retains the variation
of the homotopy types of the leaves of $\F$, whereas in $X\sslash\F$ these leaves are contracted to
stack points $BH_x$. By the base change formula Proposition \ref{pVI-3} $(1)$
the fiber of $\widetilde{X\sslash\F} \to X\sslash\F$ at a point $x$ is pro-constant and equivalent to
the usual homotopy type of the complex space $\widetilde{\cL^{max}_x(\F)}$. \\

\begin{rmk}
\emph{The expression} monodromic leaf space \emph{is not standard but refers to the classical notion
of the} monodromy groupoid  \emph{of a genuine foliation, which is an extension of the
holonomy groupoids. If we were to recover these groupoid notions we would simply have to consider 
the nerve groupoids of the projections $\widetilde{\pi} : X \to \widetilde{X\sslash\F}$ and 
$\pi : X \to X\sslash\F$, which should be called
the} monodromy groupoid and the holonomy groupoid of $\F$. \emph{The new feature here is that the 
presence of possible singularities in $\F$ first implies that these groupoids are flat groupoids but are not smooth in general (except when $\F$ is a genuine smooth foliation). Moreover, the monodromy
groupoid does not even exist as a genuine groupoid but only as a pro-object in groupoids.}
\end{rmk}

\section{The Riemann-Hilbert correspondence}\label{sec:RRHTannakian}

In this section we propose a first application of the existence of the leaf spaces and monodromy
leaf spaces discussed previously, by proving a Riemann-Hilbert correspondence relating on the
one side analytic crystals along a derived foliation $\F$, and a certain $\s$-category of
relative local systems along $\F$ that we will define. For integrable smooth foliation
this recovers Deligne's relative Riemann-Hilbert correspondence of \cite{del}, but our result
also includes a relative correspondence for flat, possibly non-smooth, morphisms
between complex manifolds, as well as its generalization to the case of derived foliations. \\

We fix $X$ a separated smooth complex manifold and $\F$ be rigid and quasi-smooth derived foliation
on $X$. We assume that the assumptions of Theorem \ref{tVI-1} are satisfied so 
we have a leaf space $X\sslash\F$ and a monodromic leaf space $\widetilde{X\sslash\F}$, related by
natural projections
$$\xymatrix{
X \ar[rd]_-{\pi} \ar[r]^-{\widetilde{\pi}} & \widetilde{X\sslash\F} \ar[d]^-p \\
 & X\sslash\F.
}$$
Recall that $\pi$ is flat surjective and integrates $\F$, $\pi^*(0) \simeq \F$, 
and similarly for $\widetilde{\pi}$. The holonomy leaf space  $\widetilde{X\sslash\F}$ is only 
a pro-Deligne Mumford stack étale over $X\sslash\F$ and which is relatively 
connected. \\

\begin{df}\label{dlocalalongF}
With the notations and assumptions above,
the \emph{$\s$-category of sheaves locally constant along $\F$} (also 
called \emph{relative local systems along $\F$} is defined by 
$$\D_{lc}(\F) := \Gamma(X,\pi^{-1}(\Parf))$$
where $\Parf$ is the stack of perfect complexes of $\OO$-modules on $X\sslash\F$.
\end{df}

Suppose that $\F=*_X$ is the tautological foliation on $X$, then 
$X\sslash\F=*$ and therefore $\pi^{-1}(\Parf) = \pi^{-1}(\Parf(k))$ is the stack of locally constant 
complexes of sheaves of $k$-modules on $X$. Therefore, $\D_{lc}(\F)$ can be identified in this case
with the $\s$-category of complexes of sheaves on $X$ with locally constant 
cohomology sheaves and perfect fibers. \\

Definition \ref{dlocalalongF} possesses an interesting interpretation in terms 
of the monodromic leaf space. Indeed, by definition we have and equivalence of $\s$-categories
$$\D_{ct}(\F) \simeq \Map_{Pro(St(X\sslash\F))}(\widetilde{X\sslash\F},\Parf).$$
Therefore, objects in $\D_{ct}(\F)$ can be, and probably should be, understood 
as representations of the connected pro-stack $\widetilde{X\sslash\F}$ in the stack 
of perfect complexes $\Parf$. Morally speaking these are 
analytic families, parameterised by $X\sslash\F$, of locally constant perfect complexes on 
the maximal leaves $\cL^{max}_x(\F)$, where $x$ varies in $X\sslash\F$. The possible pathologies
related to singularities of $\F$ are taken into account in the pro-structure of the
object $\widetilde{X\sslash\F}$. \\

We let $X$ and $\F \in \Fol(X)$ as in Theorem \ref{tVI-1}. We 
define a sheaf of cdga's on $X$ as
$$\cA_\F:=|\DR_\F| \in \cdga(X),$$
where $|-|$ is the non-Tate realization functor of \ref{nonTaterealizationforgradedmixedalgebras}.
In other words $\cA_\F$ is the sheaf sending $U \subset X$ to the Hodge completed
derived de Rham cohomology along the leaves $\CDR^*(\F_U)$. This is a non-negatively graded
sheaf of cdga's (because $\F$ is quasi-smooth) and moreover we know that our conditions
imply that $H^0(\cA_\F)$ is naturally isomorphic to $\OO_\F = \pi^{-1}(\OO_{X\sslash\F})$ (because of 
Corollary  \ref{clcpII-10} and Lemma \ref{ltVI-1}), the sheaf of locally constant holomorphic functions
along $\F$.

\begin{df}\label{dalmostlocalalongF}
With the notations and assumptions above,
the \emph{$\s$-category of sheaves almost locally constant along $\F$} is defined by 
$$\D_{alc}(\F) := \Parf(\cA_\F)$$
the $\s$-category of sheaves of perfect $\cA_\F$-modules on $X$.
\end{df}

Finally, we also have the perfect crystals along $\F$ which, by definition, are
sheaves $E$ of graded mixed $\DR_\F$-modules on $X$ satisfying the following two conditions.
\begin{enumerate}
    \item As a sheaf of graded $\DR_\F$-module $E$ is free of weight $0$
$$\DR_\F \otimes_{\OO_X}E^{(0)} \simeq E.$$
\item The $\OO_X$-module $E^{(0)}$ is perfect.
\end{enumerate}
This $\s$-category will be denoted by $\Parf(\F)$ as in the algebraic situation. 

We can now consider two natural $\s$-functors
$$\xymatrix{
\Parf(\F) \ar[rr]^-{Sol} & & \cA_\F-\dg(X) & &  \D_{lc}(\F). \ar[ll]_-{\cA_\F \otimes_{\OO_\F}}
}$$
The $\s$-functor $Sol$ (for \emph{solutions}) sends a crystal $E$ to 
the sheaf of morphisms $\rch_{\Parf(\F)}(\OO_X,E)$, where $\OO_X$ is the usual unit 
crystal along $\F$. It is obviouly equipped with a left $\cA_\F$-module structure,
as $\cA_\F \simeq \rch_{\Parf(\F)}(\OO_X,\OO_X)$. The $\s$-functor on the right hand side
simply is the base change from $\OO_\F$-module to $\cA_\F$ via the canonical morphism
$$\OO_\F = H^0(\cA_\F) \to \cA_\F$$
whose existence is insured by the fact that $\cA_\F$ is cohomologically concentrated in 
non-negative degrees. By construction, this second $\s$-functor factors throught the full 
sub-$\s$-category $\Parf(\cA_\F) \subset\cA_\F-\dg(X)$. However, we de not 
if $Sol$ has its image contained in perfect $\cA_\F$-modules as well, and we are thus forced
to introduce the notion of unipotent objects as follows.

\begin{df}
An object $E \in \Parf(\F)$ is \emph{unipotent}, if locally on $X$, 
$E$ belongs to the thick triangulated sub-$\s$-category of $\Parf(\F)$
generated by the unit crystal $\OO_X$. The full sub-$\s$-category of unipotent objects will be 
denoted by 
$$\Parf^{uni}(\F) \subset \Parf(\F).$$
\end{df}

By the definition above, it is now clear that if $E$ is unipotent, then $Sol(E)$ is itself locally
in the thick triangulated sub-$\s$-category of $\cA_\F-\dg(X)$ generated by $\cA_\F$. In other
words, $Sol$ sends unipotent crystals to perfect $\cA_\F$-modules. 

\begin{thm}\label{tVI-2}
The $\s$-functors constructed above
$$\xymatrix{
\Parf^{uni}(\F) \ar[rr]^-{Sol} & & \Parf(\cA_\F) & &  \D_{lc}(\F) \ar[ll]_-{\cA_\F \otimes_{\OO_\F}}
}$$
satisfy the following assertions.
\begin{enumerate}
    \item The $\s$-functor $Sol$ induces an equivalence
    $\Parf^{uni}(\F) \simeq \Parf(\cA_\F)$.
    \item The $\s$-functor $\cA_\F \otimes_{\OO_\F}$ induces an equivalence 
    when restricted to the genuine sub-category $\Vect(\OO_\F) \subset \D_{lc}(\F)$, consisting of 
    locally free $\OO_\F$-modules of finite rank. 
    \item Any object $E \in \Parf^{uni}(\F)$ which is a vector bundle as an $\OO_X$-module is
    unipotent.
\end{enumerate}
\end{thm}

\textit{Proof.} $(1)$ By the definition of unipotent objects $Sol$ is clearly fully faithful.
It is also locally essentially surjective, as perfect $\cA_\F$-modules are locally 
build out from $\cA_\F$. As both $\s$-categories are global sections of stacks over $X$, 
this shows that in fact $Sol$ is an equivalence of the $\s$-categories as wanted. 

$(2)$ This is formally implied by the fact that $\cA_\F$ is cohomologically non-negatively 
graded, and the fact that the morphism $\OO_\F \to \cA_\F$ identifies
$\OO_F$ with $H^0(\cA_\F)$. The image of $\Vect(\OO_F)$ by $\cA_\F \otimes_{\OO_\F}$ clearly 
consists of all locally free $\cA_\F$-modules, which we can denote by $\Vect(\cA_\F)$.

$(3)$ This is a consequence of \cite{mal2} and is proved in 
\cite[Thm. 3.2.3]{zbMATH07612808}. We do not reproduce the poof here.
\hfill $\Box$ \\

\begin{rmk}
\emph{Part $(1)$ of the theorem above is true under much weaker assumptions on $\F$. In particular
it also applies to some situations where $\F$ is not rigid. We refer to \cite{zbMATH07612808} for more
on this.}
\end{rmk}

An important consequence of Theorem \ref{tVI-2} is the following corollary, which is
the most interesting version of our Riemann-Hilbert correspondence.

\begin{cor}\label{ctVI-2}
Under the same assumptions as in Theorem \ref{ctVI-2}, the functor $Sol$ induces an equivalence
of categories
$$\Vect(\F) \simeq \Vect(\OO_\F) \simeq \Gamma(X,\pi^{-1}(\Vect)).$$
\end{cor}

Combine with GAGA, this gives the following corollary.

\begin{cor}\label{ctVI-3}
Let $X$ be a smooth and proper algebraic complex variety and $\F \in \Fol(X)$. We assume that 
$\F^h \in \Fol(X^h)$ satisfies the conditions of Theorem \ref{tVI-1}. Then, $Sol$
induces an equivalences of categories
$$\Vect(\F) \simeq \Vect(\OO_{\F^h}) \simeq \Gamma(X^h,\pi^{-1}(\Vect)).$$
\end{cor}

Tracking back the various functors, it is easy to see that the equivalences of the above corollary
are equivalence of tensor categories, and are further functorial in the pair $(X,\F)$. In 
particular, fixing a point $x\in X(\CC)$ provides a fiber functor
$$\omega_x : \Vect(\F) \to \Vect(\CC)$$
as explained in \ref{sec:transversal}. The abstract algebraic holonomy group 
$\widehat{Hol}^{alg}_\F(x)$ is defined in Definition \ref{ldef_algholtype} as the 
group of tensor automorphisms of $\omega_x$ restricted to the tensor
sub-category generated by the conormal bundle $\N_\F^*$. By 
corollary \ref{ctVI-3}, when $X$ is smooth proper and $\F$ satisfies the assumptions
of \ref{tVI-1}, this can also be interpreted as the group of tensor automorphisms of the fiber
functor
$$\omega_x : \Vect(\OO_{\F^h}) \to \Vect(\CC)$$
now restricted to the object $\pi^{-1}(\Omega^1_{X\sslash\F}) \in  \Vect(\OO_{\F^h})$. Now, the holonomy 
group $H_x$, in the sense of Definition \ref{dVI-3}, obviously acts on the
fiber of $\Omega_{X\sslash\F}^1$ at $x$, and thus there is a natural morphism
$$H_x \to \widehat{Hol}^{alg}_\F(x)$$
from the holonomy in the sense of \ref{dVI-3} to the abstract algebraic holonomy in the sense
of \ref{ldef_algholtype}. Both of these groups acts faithfully on the 
formal completion of $X\sslash\F$ at $\pi(x)$ and thus the morphism $H_x \to \widehat{Hol}^{alg}_\F(x)$
is injective. We believe it can also be Zariski dense in many situations but we do not adress
this question here. Finally, we can now compare the two notions of holonomy we have defined
in this book.

\begin{cor}\label{ctVI-4}
With the notations and assumptions of Corollary \ref{ctVI-3}, 
the holonomy group $H_x$ of \ref{dVI-3} identifies 
with a sub-group of the algebraic holonomy $Hol^{alg}_{\F}(x)$ of \ref{sec:transversal}.
\end{cor}

An important aspect of the corollary above is that 
the group $Hol^{alg}_{\F}(x)$ admits a purely algebraic description whereas the group $H_x$ is of 
transcendental nature.

\section{Reeb stability and algebraic integrability}

In this section we give a second application of the existence of a leaf space, in order to 
prove a fully algebraic integrability result. This theorem is on the one side an algebraic 
version of \cite{MR1913291} provde in the Khaler situation and for smooth foliations, but also
a generalization to the quasi-smooth and thus singular setting. Our proof is mainly inspired from 
\cite{MR1913291} with the use of Gromov compacity being replaced by the properness of Hilbert schemes.

\begin{thm}\label{tVI-3}
Let $X$ be a smooth, proper and connected
complex algebraic variety and $\F \in \Fol(X/\CC)$ be a transversely smooth
derived foliation on $X$. We suppose that $\F$ is flat, smooth in codimension $2$, and locally 
connected. We denote by $X^h$ and $\F^h$ the analytifications of $X$ and $\F$.
Then, the following three conditions are equivalent. 
\begin{enumerate}
    \item The analytic Deligne-Mumford stack $X^h\sslash \F^h$ is algebraizable and proper. 
    \item There exists a smooth and proper effective Deligne-Mumford stack $M$, and 
    a flat surjective morphism $p : X \to M$ with connected fibers such that 
    $p^{*}(0)\simeq \F$.
    \item There exists one point $x \in X^h$, such that $\cL_x^{max}(\F^h)$,
    the maximal leaf passing through $x$, is compact, and 
    its holonomy group $Hol_{\F^h}(x)$ is finite (or equivalently
    the holonomy covering $\widetilde{\cL_x^{max}(\F^h)}$ is compact).
\end{enumerate}
\end{thm}

\textbf{Proof.} $(1) \Rightarrow (2)$ Let $M$ be a smooth and proper Deligne-Mumford
stack with $M^h \simeq X^h\sslash \F^h$. By GAGA, the morphism $\pi : X^h \to M^h$
algebraizes to a morphism $p : X \to M$. Moreover, the GAGA theorem for perfect derived foliations
(Theorem \ref{GAGAperfectfol}), we know that the analytification $\s$-functor $\Fol(X) \to \Fol(X^h)$
is an equivalence of categories. From $\pi^*(0)\simeq \F^h$, and from the compatibility of
the analytification with pull-backs (see \S \ref{sec_GAGAderfol}), we deduce that $p^*(0)\simeq \F$. 
Moreover, as $\pi=p^h$ has connected fibers, so does $p$. \\

$(2) \Rightarrow (3)$ By the universal property of the leaf spaces (Proposition \ref{pVI-2}) we have
a canonical identification $M^h\simeq X^h\sslash \F^h$. In particular, 
the maximal leaves of $\F^h$ are automatically the analytification of the fibers
of $p$, and thus, by hypothesis, they are compact. In the same way, the holonomy groups
of $\F^h$ are the isotropy groups of $M$, and so they are finite because $M$ is proper. \\

$(3) \Rightarrow (1)$ This implication is the main content of the theorem, and its proof
will be somehow long. Let $x \in X$ be a closed point such that $\cL^{max}_{x}(\F^h)$
is compact. The monomorphism $\cL^{max}_{x}(\F^h) \hookrightarrow X^h$ is thus
a closed immersion, and by GAGA $\cL^{max}_{x}(\F^h)$ arises as the analytification
of a closed subscheme $X_x \subset X$. We denote by 
$\widehat{X_x}$ the formal completion of $X$ along $X_x$. 

In the same manner, the monomorphism $BH_x \hookrightarrow X^h\sslash /\F^h$ is 
a closed immersion. Indeed, $X^h \longrightarrow X^h\sslash \F^h$ being flat and surjective, 
it is enough (by flat descent) to show that the pull-back 
$BH_x \times_{X^h\sslash \F^h}X^h \hookrightarrow X^h$ is a closed immersion. But 
$BH_x \times_{X^h\sslash \F^h}X^h \simeq X_x^h\simeq \cL^{max}_x(\F^h)$, and we have seen 
that it is closed in $X^h$, by assumption. We can therefore consider $\widehat{BH_x}$, the formal
completion of the stack $X^h\sslash \F^h$ along the closed substack $BH_x$. We thus get a diagram 
of formal stacks

\begin{equation}\label{loccc}
 \xymatrix{
\widehat{X_x} \ar[r] \ar[d] & X \\
\widehat{BH_x}.
}
\end{equation}

The formal stack $\widehat{BH_x}$ is moreover the form $[\widehat{A}^d/H_x]$, and
as $H_x$ is finite 
we can even assume that the action of $H_x$ is induced by a linear action on $\mathbb{A}^d$. The above diagram can be 
consider as a formal $\widehat{BH_x}$-point of the stack $Fin(X)$, of finite schemes over $X$, 
whose reduced $BH_x$-point is given by $X_x \to X \times BH_x$. Since the stack $Fin(X)$ is 
an algebraic stack locally of finite presentation, this formal point can be algebraized by Artin
algebraization theorem (\cite{ArtinAlgI}). Therefore we can find 
a pointed smooth variety $(S,s)$, which is an \'etale neighborhood of $0$ in $\mathbb{A}^d$, 
with a compatible $H_x$-action (fixing $s$), 
a scheme $Z$ with a flat and proper map $Z \to [S/H_x]$, and 
a diagram

\begin{equation}\label{globbb}
  \xymatrix{
Z\ar[r]^f \ar[d]^-q & X \\
[S/H_x].
}  
\end{equation}

The diagram (\ref{globbb}) is such that the fiber $Z$ at $BH_x=[s/H_x] \hookrightarrow [S/H_x]$ is isomorphic
to $X_x$, as a subscheme in $X$. Moreover, the formal completion of (\ref{globbb})
at $X_x \subset Z$ recovers the diagram (\ref{loccc}) above.

\begin{lem}
By possibly shrinking $(S,s)$ to a smaller \'etale neighborhood of $0$ in $\mathbb{A}^d$, 
we have $f^*(\F^h)\simeq q^*(0)$ as derived analytic foliations on $Z$ (i.e. in $\Fol(Z^h)$). 
\end{lem}

\textit{Proof of the lemma.} To prove this lemma we use the uniqueness of derived enhancement 
of Proposition \ref{pV-1}. For this, we need to prove that $f^*(\F^h)$ and $q^*(0)$ are both smooth in codimension $2$, and
moreover that the corresponding differential ideals in $\Omega_{Z^h}^1$ coincide. 

We first notice that, by construction, the morpshim $f : Z \to X$
is formally \'etale around $X_x \hookrightarrow Z$, and thus, by shrinking $S$ if necessary,
we may assume that $f$ is an \'etale morphism. Therefore, it is clear that $f^*(\F^h)$ is
smooth in codimension $2$, as it is locally isomorphic to $\F$. Concerning $q^*(0)$, in order  to prove that it is smooth in codimension $2$ it is enough to show that 
$\LL_{Z/[S/H_x]}$ is a vector bundle outside of a codimension $2$ closed subset. However, 
the perfect complexes $\LL_{Z/[S/H_x]}$ and $f^*(\LL_\F)$ are quasi-isomorphic on the formal 
completion $\widehat{X_x}$. Because the stack of quasi-isomorphisms between 
two perfect complex is an Artin stack of finite presentation (see \cite{tova}), we deduce that these two perfect complexes must be quasi-isomorphic locally around $X_x \subset Z$. 
Therefore, by shrinking $S$ if necessary, we can assume that they are quasi-isomorphic on 
$Z$. This implies that also $q^*(0)$ is smooth in codimension $2$. Finally, 
we can use the same argument but now for the stack of quasi-isomorphisms between 
$\LL_{Z/[S/H_x]}$ and $f^*(\LL_\F)$ compatible with the canonical morphism
$$\xymatrix{\LL_{Z/[S/H_x]} & \Omega_Z^1 \ar[r] \ar[l] & f^*(\LL_\F).}$$
This shows that the two quotients of coherent sheaves
$$\Omega_Z^1 \to H^0(\LL_{Z/[S/H_x]}) \qquad \Omega_Z^1 \to H^0(f^*(\LL_\F))$$
are isomorphic, and thus their kernel must be equal. These kernels being precisely the
differential ideals $K_{q^*(0)}$ and $K_{f^*(\F)}$, we see that these must coincide. 
Therefore, all the conditions are met to apply Proposition \ref{pV-1} and we deduce the lemma.
\hfill $\Box$ \\

Since the morphism $f : Z \to X$ is \'etale, it is easy to see that $f^*(\F)^h$ 
is not only smooth in codimension $2$, but it is also flat and locally connected. We can thus
apply functoriality of leaf spaces (Proposition \ref{pVI-2}) in order to get a commutative diagram 
in the analytic category
$$\xymatrix{
Z^h \ar[r]^-f \ar[d]_-p & X^h \ar[d]^-{\pi} \\
Z^h\sslash f^*(\F)^h \ar[d]_-r \ar[r]_-{g} & X^g \sslash \F^h \\
[S^h/H_x]. & }$$
The morphism $r$ is here a finite \'etale cover, which can be identified with the
relative connected components of $Z \to [S/H_x]$. Since $X_x$ is connected, 
we see that this finite \'etale cover must be trivial, and thus $r$ is an equivalence. 
We therefore conclude the existence of commutative diagram
$$\xymatrix{
Z^h \ar[r]^-f \ar[d]_-p & X^h \ar[d]^-{\pi} \\
[S^h/H_x]\ar[r]_-{g} & X^g \sslash \F^h. }$$
Because $f$ is \'etale, and $q^*(0)\simeq \pi^*(0)$,  a comparison of the cotangent complexes shows that 
$g$ must be \'etale, too.

As a result, for any $y \in S^h$, we have an \'etale morphism of connected analytic spaces
$$\{y\}\times_{[S^h/H_x]}Z^h \longrightarrow \widetilde{\cL}^{max}_{f(y)}(\F).$$
However, because $p$ is proper we deduce that the left hand side is a proper analytic space, 
and that the above morphism is a finite etale cover, and thus a proper surjective morphism. 
In particular, $\widetilde{\cL}^{max}_{f(y)}(\F)$ must be compact. We conclude that 
for all $y \in X$ in the image of $f$, the maximal leaf $\cL^{max}_y(\F) \subset X^h$ is compact
and its holonomy $H_y$ is finite. As $f$ is an \'etale morphism, its image is a dense Zariski open 
subset in $X$, and we denote it by $W \subset X$. Let $y\in X-W$ be a closed point, and let 
$y_i$ be a sequence of closed points in $X$ converging to $y$ in the analytic topology. We can 
moreover assume that the holonomy groups of the $y_i$ are all trivial, as these form
a Zariski dense open substack in $[S/H_x]$.
We then consider the
corresponding sequence of closed subschemes $\cL^{max}_{y_i}(\F^h) \subset X$. By compactness 
of the Hilbert scheme of $X$, this sequence of subschemes has  
a limit $L \subset X$. Clearly $\pi$ sends $L$ to $BH_y \subset X^h\sslash \F^h$, and thus
$L$ is contained in the maximal leaf $\cL^{max}_{y}(\F^h)$. Moreover, $y$ being the limit of $y_i$, we have
$y \in L$. 
More generally, any point $y' \in \cL^{max}_{y}(\F)$ is a limit of points $y'_i \in \cL^{max}_{y_i}(\F^h)$, 
so that $y' \in L$. Therefore $L$ and $\cL^{max}_{y}(\F^h)$ coincide set theoretically. This
implies in particular that $\cL^{max}_{y}(\F^h)$ is compact, that the monomorphism
$\cL^{max}_{y}(\F^h) \hookrightarrow X^h$ is a closed immersion, and by GAGA that 
$\cL^{max}_{y}(\F^h)$ is algebraizable. We thus conclude that the maximal leaves
$\cL^{max}_{y}(\F^h)$ are compact and algebraizable for all $y \in X^h$.

It remains to show that the holonomy group $H_y=Hol_y(\F^h)$ is finite for all $y$. For this, 
we use the same techniques. We write the maximal leaf $\cL_y^{max}(\F^h)$ as the limit, in the
Hilbert scheme of $X$, of a sequence $\cL_{y_i}^{max}(\F^h)$ of compact leaves with trivial holonomy. We can fit this limit 
in a diagram of analytic spaces
$$\xymatrix{Z \ar[d] \ar[r] & X \\
S, & }$$
where $S$ is a small holomorphic disk, $Z$ is flat and proper with connected fibers over $S$, $Z \to X \times S$ is a closed immersion, and $Z_0=\cL_y^{max}(\F^h) \subset X$ is the
central fiber at the point $0 \in S$. We let $Z_t \subset X$ denotes the "generic" or Milnor fiber.
The nearby cycles in dimension $0$ form a constructible sheaf on the central fiber $X_0$, and therefore
any open cover $U_i$ of $X_0$ in $X$ can be refined to a cover such that each $U_i\cap X_t$
has a number of connected components which is uniformly bounded. By covering $X_0$ with $U_i$ belonging
to our basis $\B$ for the topology (see the construction of $X\sslash \F$ in \S \ref{subsec_existencethm}), we deduce the existence of an 
\'etale morphism $V \to X^h\sslash \F^h$ covering the point $y$, and such that $\pi^{-1}(V)\cap 
\cL_{y_i}^{max}(\F^h)$ possesses at most $m$ connected components for some fixed integer $m$. 
Equivalently, let $G \rightrightarrows V$ be the \'etale groupoid induced on $V$, 
so that $[V/G]$ is an open substack in $X^h\sslash \F^h$. Then, all the orbits of the points
$y_i$ are finite of order at most $m$ in $V$. As a result, if $h \in G_y$ is a point 
in the stabilizer of $y \in V$, since all the $y_i$'s have trivial holonomies, 
we must have that $h^m$ is the identity near $y_i$, for all $i$ (for all $i$ big enough 
such that $h^m$ is defined at $y_i$). By analytic continuation, we have that 
$h^m=id$. Now, let $H_y$ be the holonomy group at a point $y$. 
We have seen that $H_y$ is of finite exponent $m$. Moreover, as it is a subgroup of the
group of germs of holonomorphic isomorphisms of $\CC^d$ at $0$, it must be finite 
(see \cite[Lem. 2]{MR1913291}). 

To summarize, we have seen that all the maximal leaves are compact and all the holonomy groups
are finite. We still need to show that this implies that $X^h\sslash \F^h$ is \emph{proper} and \emph{algebraizable}. 
We start by proving that $X^h \sslash \F^h$ is a separated Deligne-Mumford stack, 
i.e. that its diagonal is a finite morphism. For this, we come back to the very beginning of the proof
of $(3) \Rightarrow (1)$. Now that we know thar all the maximal leaves are compact and all
holonomy groups are finite, we can run exactly the same argument starting with any point $x \in X$.
We have seen that there exists a commutative diagram 
$$\xymatrix{
Z^h \ar[r]^-f \ar[d]_-{p} & X^h \ar[d]^-{\pi} \\
[S^h/H_x] \ar[r]_-g & X^h\sslash \F^h.
}$$
In this diagram $p$ is a proper morphism, and moreover the natural 
morphism $Z^h \to [S^h/H_x]\times_{X^h\sslash \F^h} X^h$ is \'etale. 
Properness of $p$ implies that this \'etale morphism is moreover 
a finite covering. As a consequence $[S^h/H_x]\times_{X^h\sslash \F^h} X^h \to [S^h/H_x]$
is proper. As the \'etale morphisms of the form $[S^h/H_x] \to X^h\sslash \F^h$ cover
the whole stack $X^h\sslash \F^h$, we deduce that the morphism $\pi : X^h \to X^h\sslash \F^h$
is a proper morphism.

We now consider the nerve of $\pi$, which is a proper and flat groupoid $G \rightrightarrows X^h$. 
The diagonal morphism $G \to X^h\times X^h$ is thus a proper and unramified morphism, hence 
a finite morphism. As $X^h \to X^h\sslash \F^h$ is a flat covering, this implies that the
diagonal of the stack $X^h\sslash \F^h$ is a finite morphism, and thus that 
$X^h\sslash \F^h$ is a separated Deligne-Mumford stack. Finally, as $X^h \to X^h\sslash \F^h$
is proper and surjective we deduce that $X^h\sslash \F^h$ is also proper. 

To finish the proof, we need to show that $X^h\sslash \F^h$ is algebraizable. But this follows
by considering the groupoid $G \rightrightarrows X^h$ above (i.e. the nerve of the morphism $\pi$). 
Since $G \to X^h\times X^h$ is finite, and $X^h$ and $G$ are both proper algebraic spaces, 
we know by GAGA that the groupoid $G \rightrightarrows X^h$ is the analytification 
of a proper flat algebraic groupoid $G' \rightrightarrows X$. We can then consider the 
quotient stack $[X/G]$ of this flat groupoid: it is an (algebraic) Deligne-Mumford stack
whose analytification is $X^h\sslash \F^h$.
\hfill $\Box$ \\

\begin{rmk}
\emph{We can combine the above theorem with Corollary \ref{ctVI-4} which shows that 
finiteness of $H_x$ is implied by the finiteness of the algebraic holonomy group
$Hol_\F^{alg}(x)$, which can be described purely in algebraic terms with using the transcendent
topology.}
\end{rmk}

\chapter{Derived foliations on schemes and stacks in non-zero 
characteristic}\label{chapter:nonzerocharacteristics}

In this last chapter we introduce derived foliation over arbitrary ground ring. We present
two different non-equivalent notion; derived foliations and infinitesimal derived foliations.
We show in particular how derived foliations can be useful to study characteristic
classes, and prove in that direction vanishing theorem in crystalline cohomology. 
Finally we show that infinitesimal derived foliations have remarkable formal
integrability properties, and explain how they share relations with infinitesimal
cohomology as opposed to de Rham cohomology. \\

All along this chapter we denote by $k$ an \emph{arbitrary} commutative ring.

\section{The graded circle and graded derived loop spaces}\label{sec-gradedcircle}

We start by defining a group stack called the \emph{graded circle} that we will use in order to define derived foliations in non-zero characteristic. Let $k\!<\!x\!>$ the free commutative $k$-algebra with divided powers over one generator. It can be presented
as the quotient of the free commutative $k$-algebra over generators $x^{(n)}$, for $n\geq 0$, with relations
$$x^{(n)}.x^{(m)}= \binom{n+m}{n}x^{(n+m)}$$
for all $n \geq 0$ (in particular $x^{(0)}=1$). The commutative ring $k\!<\!x\!>$ is endowed
with a natural comultiplication
$$\Delta : k\!<\!x\!> \longrightarrow k\!<\!x\!> \otimes_k k\!<\!x\!>$$
which is the morphism of $k$-algebras defined on generators by
$$\Delta(x^{(n)}) = \sum_{i+j=n}x^{(i)} \otimes x^{(j)}.$$ 
This makes $k\!<\!x\!>$ into a commutative and cocommutative
Hopf $k$-algebra, whose corresponding commutative group scheme $\Spec\, k\!<\!x\!>$ 
can be naturally identified with the PD-hull at $0$ of the additive group $\Ga$. Using the notations
from \cite{drinfeld2024prismatization}, this group scheme will be denoted by $\Ga^{\sharp}$
(it was denoted by $Ker$ in \cite{mrt}). 
As the group $\Ga^{\sharp}$
is commutative, its classifying stack $B\Ga^\sharp$ inherits a canonical group structure and will be
considered as a group stack.

\begin{df}\label{dVII-1}
The \emph{graded circle} (over $k$) is the group stack 
$$S^{1}_{gr} := B\Ga^{\sharp}$$
\end{df}

\begin{rmk}\label{rVII-1}
\emph{The terminology} graded circle \emph{suggests that $S^{1}_{gr}$ is the associated graded of 
a filtered object. This} filtered circle \emph{is introduced in \cite{mrt} where it is
identified with the universal Hochschild-Kostant-Rosenberg filtration. }
\end{rmk}

We can give several interpretations of the group $\Ga^\sharp$ and of the group stack $S^{1}_{gr}$. We will mainly use
the following cohomological interpretation, showing that $S^1_{gr}$ behaves like a \emph{formal circle}.

The natural morphism $\Ga^\sharp \to \Ga$ induces a canonical morphism on the corresponding
classifying stacks $S^1_{gr} \to B\Ga$. This defines a canonical degree $1$ cohomology class
$$\epsilon^{(1)} \in H^1(S^1_{gr},\OO).$$
This class can be shown to be a generator for the cohomology ring 
and thus provides a natural graded ring isomorphism 
$$H^*(S^{1}_{gr},\OO) \simeq k\oplus k.\epsilon^{(1)}.$$
We can say more by considering the following commutative diagram of stacks
$$\xymatrix{
\Spec\, k[\epsilon] \ar[r] \ar[d] & Spec\, k \ar[d] \\
\Spec\, k \ar[r] & S^1_{gr},
}$$
for which the natural isomorphism making this diagram commutative is given 
by the element $\epsilon$ in $\Ga^{\sharp}(k[\epsilon])$: this is the element 
given by the morphism $k\!<\!x\!> \to k[\epsilon]$ sending $x^{(1)}$ to $\epsilon$ and
all the $x^{(n)}$ for $n>1$ to zero. The commutative diagram above then provides a canonical 
morphism on $\mathcal{O}$-cohomology complexes 
$$C^*(S^1_{gr},\OO) \to k \times_{k[\epsilon]}k$$
where the fiber product on the right hand side is taken in the $\s$-category $\dg_k$. 

\begin{prop}\label{pVII-1}
The morphism above induces an equivalence in $\dg_k$
$$C^*(S^1_{gr},\OO) \to k \times_{k[\epsilon]}k \simeq k\oplus k[-1].$$
\end{prop}

We will not provide the proof of the Proposition above and will refer to 
(see \cite[Prop. 3.4.3]{mrt} for details) for the details. 
This equivalence lives in $\dg_k$ but, by construction, can be lifted to include algebra 
structures. For us it is important to realize it as an equivalence of commutative cosimplicial 
algebras (see \cite[Thm. 3.4.17]{mrt}).

\begin{rmk}
\emph{Proposition \ref{pVII-1} is one of the main reason to introduce the group stack $S^1_{gr}$. 
We warn the reader that it badly fails if one tries to use $B\Ga$ as in the characteristic zero case, 
as already $H^2(B\Ga,\Ga)\neq 0$ (as length $2$ Witt vectors provides a non-trivial extension
of $\Ga$ by itself)}. 
\end{rmk}

The multiplicative group $\Gm$ acts on the group scheme $\Ga^{\sharp}$ by dilations. Functorially, this
action is described as follows. For a commutative $k$-algebra $A$, $\Gm(A)$ is the group
of invertible elements $\lambda$ in $A$ and $\Ga^{\sharp}(A)$ is the additive group of sequences
$\{a^{(n)}\}$ of elements in $A$ satisfying $a^{(n)}.a^{(m)}= \binom{n+m}{n}a^{(n+m)}$. The action 
of $\lambda \in A^*$ is then given by
$$\lambda\{a^{(n)}\} = \{\lambda^{n}a^{(n)}\}$$
This provides an action of $\Gm$ on the group stack $S^1_{gr}$. The multiplicative group
$\Gm$ also acts naturally on $\Spec\, k[\epsilon]$ by dilations, and the equivalence of  Proposition \ref{pVII-1}
is then naturally compatible with the $\Gm$-actions on both sides. \\

The importance of Proposition \ref{pVII-1} (when taking into account the algebras' structures) is that the graded circle acts as the spectrum of a
square zero trivial extension by $k[-1]$. An important consequence of this fact is the following
corollary.

\begin{cor}\label{cVII-1}
For any derived affine $k$-scheme $\Spec\, A$, we have a canonical $\Gm$-equivariant identification
of derived stacks
$$\uMap(S^1_{gr},X) \simeq \VV(\LL_X[1]) = \Spec\, (\mathbf{Sym}_A(\LL_{A/k}[1])).$$
\end{cor}

\textit{Sketch of proof.} The commutative diagram 
$$\xymatrix{
\Spec\, k[\epsilon] \ar[r] \ar[d] & Spec\, k \ar[d] \\
\Spec\, k \ar[r] & S^1_{gr},
}$$
defines a commutative diagram of derived stacks
$$\uMap(S^1_{gr},X) \longrightarrow X \times_{\VV(\LL_X)} X \simeq \VV(\LL_X[1]).$$
The fact that this morphism is an equivalence can be easily proven by the Whitehead theorem: 
a morphism of derived Artin stacks is an equivalence if and only if it is formally \'etale and
an equivalence on the classical truncations (see \cite[2.2.2]{hagII}). The formally \'etale assumption 
can be checked using Proposition \ref{pVII-1} as the tangent complex of the right hand side
can be understood as a cohomology complex of $S^1_{gr}$. We refer to \cite[Thm. 5.1.3]{mrt} 
for more details.
\hfill $\Box$ \\

As a consequence of Corollary \ref{cVII-1}, we see that, for any derived affine $k$-scheme $X$, $\VV(\LL_X[1])$
carries a canonical action of the semi-direct product group stack $$\cH :=\Gm \ltimes S^1_{gr}.$$ This semi-direct 
group stack $\cH$ is the pertinent generalization of our group $\cH_0 = \Gm \ltimes B\Ga$
(Definition \ref{dI-5}) outside of characteristic zero situations. The representations of $\cH$
are indeed graded mixed complexes, as we will see in the next sections. For the moment we use
$\cH$ in order to introduce the notion of graded loop space.

\begin{df}\label{dVII-2}
For a derived affine $k$-scheme $X$, the \emph{graded loop space} is the derived
affine scheme 
$$\cL^{gr}(X/k):=\uMap(S^1_{gr},X).$$
It comes naturally equipped with a canonical
action of the group $\cH$.
\end{df}

For non-affine derived stacks, we use the left Kan extension of the functor $X \mapsto \cL^{gr}(X/k)$. To 
avoid any confusion we will denote by $\ccL^{gr}(-/k)$ this extension. 

\begin{df}\label{dVII-3}
The derived graded loop space of a derived stack $X$ is $\ccL^{gr}(X/k)$
defined above. It comes with a canonical
action of the group $\cH$.
\end{df}

Note that for general derived Artin stack $X$, Corollary \ref{cVII-1} fails. There is always a 
natural descent morphism $\ccL^{gr}(X) \to \VV(\LL_X[1])$ but this is not 
an equivalence in general. 
For a general derived Artin stack things get even more complicated and $\ccL^{gr}(X)$ involves
some form of formal completion of $\VV(\LL_X[1])$ along the zero section. 

The induced morphism
$\ccL^{gr}(X/k) \to \VV(\LL_X[1])$ is however an equivalence when 
$X$ is more specifically a derived Deligne-Mumford stack.

\begin{cor}\label{cVII-2}
For a derived Deligne-Mumford stack $X$ the natural morphism
$$\ccL^{gr}(X/k) \longrightarrow \VV(\LL_X[1])$$
is an equivalence.
\end{cor}

\textit{Proof.} Indeed, both functors $\ccL^{gr}(-/k)$ and $\VV(\LL_{-}[1])$
are stacks on the small \'etale site of $X$, and Corollary \ref{cVII-1} implies that 
the morphism in consideration defines a local equivalence between these stacks.
\hfill $\Box$ \\

\begin{rmk}\label{rVII-2}
\emph{For a general Artin stack $X$ there is a close connection between $\ccL^{gr}(X/k)$ and
$\VV(\LL_X[1])$. It is for instance possible to prove that these two derived stacks
have equivalent algebras of $\Gm$-equivariant functions. This can be proved 
using some standard descent property as well as extensions of 
the results of \S\ref{secapp:CotangentcomplexesofderivedArtinstacks} to the case of arbitrary characteristics (replacing
$\Lambda$-powers with their left derived versions). We will not pursue this direction.}
\end{rmk}

The geometric interpretation of graded mixed complexes in characteristic zero mentioned in \S \ref{sec:Thegeometryofgradedmixedobjects}
remains valid in any characteristic using the group stack $\cH$ defined above when working over
general base rings. The construction of \S \ref{sec:Thegeometryofgradedmixedobjects} remains similar, and produces a canonical
$\s$-functor
$$\phi: \epsilon-\dg_k \longrightarrow \QCoh(B\cH).$$

\begin{prop}\label{pVII-2}
The $\s$-functor $\phi$ induces an equivalence of 
symmetric monoidal $\s$-categories
$$\phi : \edg_k \simeq \QCoh(B\cH).$$
\end{prop}

\textit{Proof.} This is proven in \cite[Prop. 4.2.3 (ii)]{mrt}.
\hfill $\Box$ \\

\section{Derived foliations via graded derived loop spaces}\label{derfolgradedloop}

The geometric interpretation of derived foliations given in Section \ref{subsec:Derivedfoliationsasequivariantderivedlinearstacks} (see 
Definition \ref{defgeoderfol}) can be generalized
in arbitrary characteristic by replacing $\cH_0$ with the group stack $\cH = \Gm \ltimes S^1_{gr}$.

\begin{df}\label{dVII-4}
Let $X$ be a derived stack over $k$. A 
\emph{derived foliation $\F$ on $X$ (relative to $k$)} consists of
the data of a $\cH$-equivariant derived stack $\VV(\F)$ over $\ccL^{gr}(X/k)$, 
such the following condition is satisfied:
for any derived affine scheme $S$ and morphism $S \to X$, the $\Gm$-equivariant morphism
$$\VV(\F)\times_{\ccL^{gr}(X/k)}\ccL^{gr}(S/k) \to S$$
makes $\VV(\F)$ into a perfect derived linear stack over $S$.

By definition, derived foliations $\F$ on $X$ form a full sub-$\s$-category 
of $\dSt^{\cH}_{/\ccL^{gr}(X/k)}$, the $\s$-category of $\cH$-equivariant derived 
stacks over $\ccL^{gr}(X/k)$.
\end{df}

\begin{rmk}\emph{
Outside of the case where $k$ is a $\QQ$-algebra, the above definition does not have an easy
purely algebraic description. Indeed graded mixed cdga's must be replaced by the more
complex notion of \emph{graded mixed LSym-algebra}, where $LSym$ is the $\s$-monad given by 
the left derived symmetric powers. We refer to \cite{raksit2020hochschildhomologyderivedrham} 
for such a purely algebraic approach, and warn the
reader that the comparison with our definition above is not straightforward. In what 
follows we will thus avoid any tentative to provide an algebraic version of derived foliations
over general rings and will stick to the geometric definition above. }
\end{rmk}

The derived foliations defined above share the same formal properties that we have already seen
in characteristic zero, such as pull-backs and existence of finite limits. For a morphism
of derived stacks $X \to Y$ we have an $\s$-category $\Fol(X/Y)$ of relative
derived foliations over $X$ relative to $Y$. This $\s$-category is functorial by 
pull-back for commutative diagrams of derived stacks
$$\xymatrix{
X' \ar[r]^-{f} \ar[d] & X \ar[d] \\
Y' \ar[r]_-{s} & Y.}
$$

Explicitly, an object $\F \in\Fol(X/Y)$ is represented by an $\cH$-equivariant derived
stack $\VV(\F) \to \ccL^{gr}(X/Y)$. Its pull-back along the commutative diagram above is
given by 
$$(f,s)^{*}(\VV(\F)) = \VV(\F) \times_{\ccL^{gr}(X/Y)}\ccL^{gr}(X'/Y').$$

Cotangent complexes of derived foliations are defined as in \S \ref{subsec-ctgtofderivedfol}. Namely, for $\VV(\F) \to \ccL^{gr}(X)$ a derived foliation over $X$,  
and a derived affine $U=\Spec\, A$ with a morphism $U \to X$, the $\Gm$-equivariant derived affine stack
object $\VV(\F) \times_{\ccL^{gr}(X)} U \to U$ corresponds, by \cite{mon}, to a 
well defined $A$-dg-module $\LL^{big}_\F(U)$ such that $\VV(\F) \times_{\ccL^{gr}(X)} U$ and $\VV(\LL^{big}_\F(U))$ are equivalent in $\dSt_k/U$.  When $U$ varies in $\dAff_k/X$, these glue to 
a (non-quasi-coherent) $\OO_X$-dg-module $\LL_\F^{big}$ called the \emph{big
cotangent complex of $\F$}. Its associated quasi-coherent complex (see \S \ref{secapp:OXmodulesquasi-coherentandperfectmodules}) is called \emph{the
cotangent complex of $\F$}, and is denoted by 
$$\LL_\F := (\LL_{\F}^{big})^{qcoh} \in \QCoh(X).$$
A derived foliation $\F$ on $X$ is called \emph{perfect} if its cotangent complex $\LL_\F$ is perfect.\\

For any derived stack $X$ over $k$, we have a final object $*_X \in \Fol(X/k)$ which is given by
$\VV(\F) = \ccL^{gr}(X/k)$. Similarly, when $X$ is a derived Artin stack locally of finite presentation 
over $k$,
the $\s$-category $\Fol(X/k)$ has an initial object $0_X \in \Fol(X/k)$. It can be described
as an explicit derived stack $\VV(\F) \to \ccL^{gr}(X)$ derived stack as follows. We write
$X$ as a colimit of derived affine opschemes $colim_i\, U_i$, so that $\ccL^{gr}(X) \simeq colim_i 
\,\cL^{gr}(U_i)$ ($\ccL^{gr}$ commutes with colimits has it is defined by a left Kan extension). 
For each index $i$, we form the relative graded loop stack
$$\ccL^{gr}(U_i/X) := \cL^{gr}(U_i) \times_{\ccL^{gr}(X)} X \to \cL^{gr}(U_i),$$
which comes equipped with its natural $\cH$-action. The $\cH$-equivariant derived stacks
$\ccL^{gr}(U_i/X)$ certainly glue to a global object over $\ccL^{gr}(X)$, however
in general $\ccL^{gr}(U_i/X)$ is not a derived affine stack over $U_i$. We thus replace
$\ccL^{gr}(U_i/F)$ by its graded affinization (see \S \ref{ss1.2})
$$\ccL^{gr}(U_i/F)^{aff} := \RR\Spec^{\Delta,gr}\, (\OO(\ccL^{gr}(U_i/F)^{aff})) \to U_i.$$
The group stack $\cH$ being affine over $B\Gm$, it continues to act on 
$\ccL^{gr}(U_i/F)^{aff}$. Finally, by our descent result for cotangent complexes (see \S 
\ref{secapp:CotangentcomplexesofderivedArtinstacks})
we have that $\ccL^{gr}(U_i/F)^{aff} \simeq \VV(\LL_{U_i/F}[1])$ as graded derived affine stacks.
The family of objects $\ccL^{gr}(U_i/F)^{aff}$ then glue as a global $\cH$-equivariant derived
stack $\VV(\F) \to \ccL^{gr}(F)$, which thus defines an object in $\Fol(F)$. It can be checked, as 
in \S \ref{ex:tautologicalfoliationsglobal}, that this defines an initial object in $\Fol(F)$. \\

We finish this general part by stating that integrable sub-bundles of the tangent
bundle of a smooth scheme canonically define derived foliations in the sense of Definition \ref{dVII-4}.
This result is taken from \cite{mon2}, whose proof is reproduced below.

\begin{prop}\label{pVII-3}
Let $f : X \to S$ be a smooth morphism of derived Deligne-Mumford stacks. 
Let $E \subset \TT_{X/S}$ be a sub-bundle of the relative tangent bundle of $X$ over $S$
such that $E$ is stable by the Lie bracket of relative vector fields. Then, 
there exists a unique action of $\cH$ on $\VV(E^{\vee}[1])$, covering the
$\cH$-action on $\VV(\LL_{X/S}[1]) \simeq \cL^{gr}(X/S)$.
\end{prop}

\textit{Proof.} The purpose of this proof is to show that the space of
$\cH$-action on $\VV(E^\vee[-1])$ compatible, via the map
$\VV(\LL_{X/S}[1]) \simeq \cL^{gr}(X/S)=\VV(\LL_X[-1])$, with the $\cH$-action on 
$\VV(\LL_X[1])$, is contractible. This is a local statement on the small étale site
of $X$ and thus we can assume that $X$ is an affine smooth scheme over $k$.

We first consider $\VV(E^\vee[1])$ as a $\Gm$-equivariant derived stack, and
we study the $\Gm$-equivariant group stack $G$ of self equivalence of $\VV(E^\vee[1])$.
More precisely, $G$ is the group stack over $B\Gm$ of self-equivalences of
$\VV(E^\vee[1]) \to B\Gm$. The group stack $G$ a priory comes equipped with a natural
derived structure, but we will only be interested in the underlying underived stack 
$G \to B\Gm$. Also, we will not need the full stack $G$ but only $G_0 \subset G$
the full sub-stack of the connected component of the identity, defined as the pull-back
of stacks relative to $B\Gm$
$$\xymatrix{
G_0 \ar[r] \ar[d] & G \ar[d] \\
\star \ar[r]_-{id} & \pi_0(G).
}$$
The group stack $G_0$ is connected, $\pi_0(G_0)\simeq *$ as a sheaf over $B\Gm$.

Because $\VV(E^\vee[1])=\Spec\, Sym_A(E^\vee[1])$, it is easy to compute the homotopy sheaves
$\pi_i(G_0)$, as sheaves over $B\Gm$, or equivalently as $\Gm$-equivariant sheaves of groups
on $\Spec\, k$.
We have short exact sequence of sheaves on the big fpqc site of affine $k$-schemes
$$\xymatrix{
0 \ar[r] & Hom_A(E^\vee,\wedge^i(E^\vee))^{\sim} \ar[r] & \pi_{i}(G_0) \ar[r] & 
Der_k(A,\wedge^{i-1}E^\vee)^{\sim} \ar[r] & 0.
}$$
The left hand side is the quasi-coherent sheaf on $\Spec\, k$ 
associated with the $k$-module of $A$-linear maps $E^\vee \to \wedge^i(E^\vee)$, 
whereas the right hand side is the quasi-coherent sheaf on $\Spec\,k $ associated 
with the $k$-module of $k$-linear derivations from $A$ to 
the $A$-module $\wedge^{i-1}E^\vee$. The $\Gm$-action on both of these quasi-coherent sheaves
is induced by  the canonical homothetic action on $E^\vee$ and is thus 
pure of weight $i-1$ on both side.

Now, an $\cH$-action on $\VV(E^\vee[-1])$, compatible with the already existing $\Gm$-action, 
is given by a $\Gm$-equivariant morphism 
$$\rho : B\cH \to BG_0.$$
The stack $B\cH$ is an affine stack (see Appendix \ref{app:derivedaffine}) 
of the form $\Spec^{\Delta}\, \ZZ[u]$, where
$u$ is a free variable in degree $2$, and $\ZZ[u]$ is the free commutative cosimplicial-simplicial 
ring freely generated by $u$, as this follows
easily from the fact that $\cH = \Spec^{\Delta}\, H^*(S^1)$ (see \ref{pVII-1} 
and \cite[Thm. 3.4.17]{mrt}). 
Moreover, the $\Gm$-action on $B\cH$ is such that $u$ is of weight $1$. We can thus
understand the space of $\Gm$-equivariant morphisms
$\rho : B\cH \to BG_0$ by Postnikov induction on $BG_0$, with in mind the above description of
the homotopy sheaves $\pi_i(BG_0) \simeq \pi_{i-1}(G_0)$. For this we use the standard
cellular decomposition, inside the $\s$-category of the $\Gm$-equivariant affine stacks,
$B\cH \simeq colim_k (B\cH)_{\leq 2k}$
$$\xymatrix{
B\cH_{\leq 2k} \ar[r] & B\cH_{\leq 2k+2} \\
S^{2k+1}_{gr}(k+1) \ar[r] \ar[u]^-{h_k} & \star \ar[u]
}$$
where $S^{p}_{gr}(q)$ denotes the \emph{affine graded sphere} $\Spec^{\Delta}\, \ZZ\oplus \ZZ[-p](q)$,
and $h_k$ is induced by the graded versions of the standard Hopf maps 
$S^{2k+1}_{gr}(k+1) \to B\cH_{\leq 2k}$.

Using this cellular decomposition, and by weight considerations, we see in 
particular that the projection on the first non-trivial Postnikov layer
$$BG_0 \to \tau_{\leq 2}(BG_0)=K(\pi_2,2)$$
induces an equivalences of mapping spaces between $\Gm$-equivariant stacks
$$\Map_{\Gm-\St}(B\cH,BG_0) \simeq \Map_{-\Gm-\St}(B\cH,K(\pi_2,2)) \simeq
Hom_{\Gm-Gp}(\Ga^{\sharp},\pi_2),$$
where $\pi_2 = \pi_1(G_0)$.
Using the fact that group morphisms $\Ga^{\sharp} \to \Ga$ are in one-to-one correspondence
with elements in $k$ (because $H^1(B\Ga^{\sharp},\OO)\simeq k$, see Proposition \ref{pVII-1}), 
we see that $Hom_{\Gm-Gp}(\Ga^{\sharp},\pi_2)\simeq \pi_2(\Spec\, k)$ can be identified
with the global sections of the sheaf $\pi_2$. 

We have thus seen that the space of $\Gm$-compatible actions on $\VV(E^\vee[1])$
is discrete and equivalent to the underlying set of $\pi_2(\Spec\, k)$, which sits
in a short exact sequence
$$\xymatrix{
0 \ar[r] & Hom_A(E^\vee,\wedge^2(E^\vee)) \ar[r] & \pi_{2}(\Spec\, k) \ar[r] & 
Der_k(A,E^\vee) \ar[r] & 0.
}$$
The elements in the set $\pi_{2}(\Spec\, k)$ can be described as pairs
$(d,u)$, where $d : A \to E^\vee$ is a $k$-linear derivation, and 
$u : E^\vee \to \wedge^2 (E^\vee)$ is a $d$-semi-linear morphism 
$$u(am) = d(a)\wedge m + a.u(m).$$
We can therefore conclude the proof of the proposition, 
as we can apply this result on both $\VV(\LL_X[-1])$ and $\VV(E^\vee[1])$ and see that the
only possible compatible $\cH$-action on $\VV(E^\vee[1]))$
must be induced by the de Rham differential on the quotient 
$$Sym_A(\Omega_A^1[1]) \twoheadrightarrow Sym_A(E^\vee)$$
and is therefore given by the pair $(d,u)$, where $d : A \to E^\vee$ is
the dual of the anchor map 
$a : E \to \T_X$, and $u$ is dual to the restriction of the Lie bracket on $E$
(i.e. the Chevalley Eilenberg differential)
$$u(\phi)(x\wedge y)=a(x)(\phi(y)) - a(y)(\phi(x)) - \phi([a(x),a(y)])$$
for $\phi : E \to A$ an $A$-linear form and $x,y \in E$.
\hfill $\Box$ \\

\section{Foliated de Rham theory}

Recall that any morphism of derived stacks $f : X \to Y$ induces and
adjunction given by push-forward and pull-back on quasi-coherent sheaves
$$f^* : \QCoh(X) \leftrightarrows \QCoh(Y) : f_*.$$
The left adjoint $p^*$ is a symmetric monoidal $\s$-functor, and therefore this 
adjunction induces an adjunction on the corresponding $\s$-categories of
commutative algebras in the sense of \cite[Def. 2.1.3.1]{HA}
$$p^* : CAlg(X) \leftrightarrows CAlg(Y) : p_*,$$
and this adjunction covers the adjunction on quasi-coherent sheaves
(i.e. $p_*$ and $p^*$ commute with the forgetful $\s$-functor 
forgetting the commutative algebra structure). To void confusions objects in 
$CAlg(\QCoh(X)))$ will be called \emph{quasi-coherent $E_\s$-algebras on $X$}, and 
we will use the general notation
$$E_\s-\Alg(\mathcal{C}):=CAlg(\mathcal{C}),$$
for any symmetric monoidal $\s$-category $\mathcal{C}$.
In characteristic zero these coincide with quasi-coherent sheaves of $\OO_X$-cdga's over $X$.

Let $X \in \dSt_k$ be a derived Artin stack over $k$, and 
$\F \in \Fol(X/k)$ be a derived foliation (relative to $k$). By definition, $\F$ is given 
by an $\cH$-equivariant derived stack $\VV(\F)$ together with an 
$\cH$-equivariant morphism
$$\VV(\F) \to \ccL^{gr}(X).$$
We can pass to the corresponding quotient stacks, and obtain this way 
a morphism of derived stacks over $B\cH$
$$\xymatrix{
[\VV(\F)/\cH] \ar[r] \ar[rd]_-{p} & [\ccL^{gr}(X)/\cH] \ar[d]^-{q} \\
 & B\cH.
}$$
By considering the adjunctions induced by this morphisms we get
a canonical morphism of quasi-coherent $E_\s$-algebras 
$q_*(\OO) \to p_*(\OO)$ on $B\cH$. We can compose with the symmetric monoidal 
equivalence given by the Tate realization (see Corollary \ref{clI-4})
$$|-|^t : \QCoh(B\cH) \simeq \cfdg_k,$$
and get this way a morphism of filtered and complete $E_\s$-algebras
$|q_*(\OO)|^t \to |p_*(\OO)|^t$.

\begin{df}\label{dVII-5}
Let Let $X \in \dSt_k$, and 
$\F \in \Fol(X/k)$. With the above notations, we define 
\begin{enumerate}
    \item the \emph{completed derived de Rham cohomology of $X$ (relative to $k$)}
    is the complete filtered $E_\s$-algebra
    $$\CDR(X,\OO) := |q_*(\OO)|^t \in E_\s-\Alg(\cfdg_k).$$
    \item The \emph{completed foliated derived de Rham cohomology of 
    the derived foliation $\F$ (relative to $k$)} is
    the complete filtered $E_\s$-algebra 
    $$\CDR(\F,\OO) := |p_*(\OO)|^t \in E_\s-\Alg(\cfdg_k).$$
\end{enumerate}
The corresponding cohomology groups will be denoted by
$$\widehat{H}^*_{DR}(X):=H^i(\CDR^*(X,\OO)) \qquad \widehat{H}^*_{DR}(\F):=H^i(\CDR^*(\F,\OO)).$$
\end{df}

We also have a definition \emph{with coefficients in crystals}, as 
in Definition \ref{def-foldRwithCryscoeff}. To get this, we introduce the $\s$-category of (quasi-coherent) crystals along 
$\F$ as a certain full sub-$\s$-category $\QCoh(\F)$ of $\QCoh([\VV(\F)/\cH])$. 
To explain the definition of $\QCoh(\F)$ ,
we start with a general derived stack $X$ with an action
of $\Gm$ together with a $\Gm$-equivariant morphism $X \to Y$, 
considered as a morphism $p : [X/\Gm] \to Y$.
A quasi-coherent complex $E \in \QCoh([X/\Gm])$
is said to be \emph{graded free relative to $Y$}, if the
canonical morphism
$$p^*(p_*(E)) \to E$$
is an equivalence in $\QCoh([X/\Gm])$. With this notion, 
by definition an object 
$E \in \QCoh([\VV(\F)/\cH])$ lies in $\QCoh(\F)$
if its pull-back to $\QCoh([\VV(\F)/\mathbb{G}_m])$ is graded free
with respect to the $\Gm$-equivariant projection $\VV(\F) \to X$. 

For any crystal $E \in \QCoh(\F)$ on $\F$, we thus have 
a \emph{complete foliated derived de Rham cohomology with coefficients in} $E$,
defined by
$$\CDR(\F,E):=|p_*(E)|^t.$$
The object $\CDR(\F,E)$ is canonically a complete filtered module over the complete filtered 
$E_\s$-algebra
$\CDR(\F,\OO)$. \\

One can prove, using the descent results of \S \ref{secapp:CotangentcomplexesofderivedArtinstacks}, 
that 
the underlying graded complex of $p_*(\OO)$ is equivalent
to 
$$p_*(\OO) \simeq \bigoplus_{p}\HH(F,\LL \Lambda^p\LL_\F)[p].$$

The Hodge filtration discussed in \S \ref{sec:Hodgefiltration} also makes sense in the present 
setting. The only difference 
is here that the \emph{left derived} functors of the exterior product must be used. We gather the 
results in the
following Proposition, and leave the details to the reader.

\begin{prop}\label{pVII-4}
For a perfect derived foliation $\F$ on a derived Artin stack $X$, there exists a canonical 
filtration $F^*_\F$ on the $E_\s$-algebra $\CDR(X/k)$ of completed de Rham cohomology on $X$, with the 
following properties.
\begin{enumerate}
    \item The filtration $F^*_\F$ is functorial in $X$ and $\F$.
    \item The associated graded, as a graded $E_\s$-algebra, is given by 
    $$Gr^{i}_{F^*_\F}(\CDR(X/k)) \simeq \CDR(\F,\LL\Lambda^{-i}\N_\F^*)[i].$$
\end{enumerate}
\end{prop}

In the above proposition, $\LL\Lambda^i\N^{*}_\F$ is the $i$-th derived exterior product of $\N^{*}_\F$ (defined as usual as being the fiber of $\LL_{X/k} \to \LL_{\F/k}$),
equipped with its natural structure of a perfect crystal along $\F$ 
(see Proposition \ref{pIII-1}).

\section{Characteristic classes}\label{sec:characteristicclasses}

In this section we construct \emph{characteristic classes} of derived foliations
as cohomology classes in the foliated de Rham coholomogy introduced in the previous 
section. When the base is a perfect field, we then use the comparison between 
derived de Rham cohomology and crystalline cohomology to prove the existence
of residues for foliations, in the sense of \cite{babo}, in the positive characteristic
setting. \\

We start by general results and constructions concerning the de Rham 
cohomology of the classifying stacks $B\Gl_n$. The completed derived de Rham
complex $\CDR^*(B\Gl_n,\OO)$ is a complete filtered $E_\s$-algebra whose
associated graded is equivalent to
$$\bigoplus_p H^*(B\Gl_n,S^p(\mathfrak{g}^*))[-2p],$$
where $\mathfrak{g}=Lie(\Gl_n)$ is the Lie algebra of the group scheme
$\Gl_n$ and $S^p(\mathfrak{g}^*)$ denotes the (underived) $p$-th symmetric power
of its linear dual (as a $k$-module). The group $\Gl_n$ acts by the adjoint action on $\mathfrak{g}$
and thus on $S^p(\mathfrak{g}^*)$, this action defines $S^p(\mathfrak{g}^*)$
as a quasi-coherent sheaf on the stack $B\Gl_n$. 

By construction the complete filtered $E_\s$-algebra $\CDR^*(B\Gl_n)$
can be written as $|q_*(\OO)|^t$ where $q_*(\OO)$ is a quasi-coherent
$E_\s$-algebra over $B\cH$. The symmetric monoidal $\s$-category
$\QCoh(B\cH)$ can be endowed with its canonical t-structure, for which the pull-back
along the canonical point $\Spec\, k \to \QCoh(\B\cH)$ is t-exact and conservative.
Therefore, $q_*(\OO)$ possesses a connective cover $B \to q_*(\OO)$ inside
$E_\s-\Alg(\QCoh(B\cH))$, which is simply given by the non-derived direct image
of $\OO$ by $q$. The underlying graded complex of this connective cover $B$
is concentrated in degree $0$ and is
isomorphic to 
$\oplus_p(S^p(\mathfrak{g}^*))^{\Gl_n}$. As $B$ is concentrated in cohomological
degree $0$ the $B\cH$-action is automatically trivial, and thus the induced morphism
on the Tate realizations provides a morphism of complete filtered $E_\s$-algebras
We thus obtain a natural morphism of filtered $E_\s$-algebras
$$|B|^t := \A_n \to \CDR^*(B\Gl_n,\OO)=|q_*(\OO)|^t.$$
By the definition of the Tate realization and because the mixed structure is trivial
on $B$, we have
$$\A_n \simeq \oplus_p S^p(\mathfrak{g}^*)^{\Gl_n}[-2p]$$ 
where the filtration is the split filtration associated with the
graduation for which $(\mathfrak{g}^*)^{\Gl_n}[-2]$ sits in weight $-1$.

The inclusion of the maximal torus in $\Gl_n$ induces, by restrictions,
a morphism of graded $E_\s$-algebras
$$\A_n \to k[x_1,\dots,x_n]^{\Sigma_n}$$
which is an isomorphism (and each $x_i$ sits in degree $2$ and weight $-1$).
We thus have an isomorphism of graded $E_\s$-algebras
$$k[c_1,\dots,c_n] \simeq \A_n$$
where $c_i$ corresponds to the $i$-th elementary symmetric polynomial
in the $x_1, \dots, x_n$, and is therefore of cohomological degree $2i$ and weight $-i$. 
Via the map $\A_n \to \CDR^*(B\Gl_n,\OO)$ these classes defines natural morphism of complexes
of $k$-modules
$$c_i : k \to F^{-i}\CDR^*(B\Gl_n)[2i].$$

\begin{rmk}
\emph{In \cite[Thm. 10.2]{zbMATH06896953} it is proven that the higher cohomology groups
of $\HH(B\Gl_n,S^p(\mathfrak{g}^*))$ all vanish, and therefore the natural morphism induces
an isomorphism}
$$\A_n \to \CDR^*(B\Gl_n)$$
\emph{However, the construction of the Chern classes $c_i$ in $\CDR^*(B\Gl_n)$ does not require
the results of \cite{zbMATH06896953}, as the mere existence of the morphism from $\A_n$
to $\CDR^*(B\Gl_n)$ is sufficient.}
\end{rmk}

We now extend the universal Chern class to a 
class in derived completed de Rham cohomology of the derived stack 
of perfect complexes, using the standard methods of \cite{zbMATH03751119}.

Each Chern class $c_i$ can be considered as a mophism of 
underived stacks on the big etale site of affine $k$-schemes
$$c_i : \ZZ \times B\Gl_\s \to F^{-i}\CDR^*(-,\OO)[2i],$$
where $\Gl_\s$ is defined as the colimit of the $\Gl_n$ for the natural
inclusion $\Gl_n \hookrightarrow \Gl_{n+1}$ given by the direct sum wuth a trivial line bundle.
By the universal property of Quillen "+ construction",
this factors as a morphism of stacks
$c_i : (\ZZ \times B\Gl_\s)^+ \to 
F^{-i}\CDR^*(-,\OO)[2i],$
or in other words as a morphism 
$$c_i : \underline{K}^{\mathrm{cl}} \to F^{-i}\CDR^*(-,\OO)[2i]$$
where $\underline{K}^{\mathrm{cl}}$ is the $K$-theory stack (i.e. the \'etale
stackification of the $K$-theory functor sending a commutative
$k$-algebra $A$ to its $K$-theory space $K(A)$). If we denote by 
$i^* : \dSt_k \to \St_k$ the restriction functor from the $\s$-category of derived stacks to the $\s$-category of stacks, 
the above morphism can be written as
$$c_i : i^*(\underline{K}) \to i^*(F^{-i}\CDR^*(-,\OO)[2i])$$
where now $\underline{K}$ is the derived stack sending 
a simplicial commutative $k$-algebra $A$ to the Waldhausen
$K$-theory space of the $\s$-category of perfect $A$-modules.
The functor $i^*$ has a left adjoint $i_! : \St_k \to \dSt_k$
obtained by the left Kan extension induced by
 the embedding of $k$-algebras into simplicial $k$-algebras.
By adjunction, the morphism $c_i$ induces a morphism of derived
stacks
$$c_i : i_!i^*(\underline{K}) \to F^{-i}\CDR^*(-,\OO)[2i],$$
However, it is well known that the adjunction morphism
$i_!i^*(\underline{K}) \to \underline{K}$ is an equivalence, and this
can be seen easily from the fact that the derived stack of vector bundles
$\Vect$ is such that the ajdunction morphism
$i_!^*(\Vect) \to \Vect$ is an equivalence (simply because $B\Gl_n$ is a smooth
underived stack). We therefore obtain a morphism of derived stacks
$$c_i : \underline{K} \to F^{-i}\CDR^*(-,\OO)[2i].$$
Composed with the canonical morphism $\Parf \to \underline{K}$
we obtain the desired extensions
$$c_i : \Parf \to F^{-i}\CDR^*(-,\OO)[2i],$$
and thus canonical classes
$$c_i \in H^{2i}(F^{-i}\CDR^*(\Parf,\OO))).$$

\begin{df}\label{dVI-5}
The \emph{universal $i$-th Chern class} is 
$$c_i \in H^{2i}(F^{-i}\CDR^*(\Parf,\OO))$$
constructed above.
\end{df}

We are now ready to define characteristic classes of
perfect crystals over derived foliations. Let $X$ be a derived Artin stack (locally of finite
presentation over $k$) and
$\F \in \Fol(X/k)$ a derived foliation on $X$ (relative to $k$).

Let $E \in \Parf(\F)$ be a perfect crystal along $\F$.
The object $E$ can be represented by 
a morphism between pairs, consisting 
of derived stacks with derived foliations
$$\phi_E : (X,\F) \to (\Parf,0_\Parf)$$
where $0_\Parf$ is the initial foliation on $\Parf$. 
Therefore, the morphism $\phi_E$ induces a morphism on completed derived de Rham cohomologies
$$\CDR(\Parf,\OO) \to \CDR(X,\OO)$$
which is a morphism of filtered $E_\s$-algebras for the Hodge filtrations
on both sides (Proposition \ref{pVII-4}). The Hodge filtration $F_{0_\Parf}^{*}$
on the left hand side is the usual Hodge filtration $F$ as the derived foliation on $\Parf$ is here
the initial foliation. We therefore obtain for any integer $i$ an induced morphism
$$F^{i}\CDR^*(\Parf,\OO) \longrightarrow F_\F^i\CDR^*(X,\OO).$$
As a consequence, we see that the image of the universal Chern classes $c_i$
in $\widehat{H}^{2i}_{DR}(X)$ has natural lifts through 
$\widehat{H}^{2i}(F_\F ^{-i}\CDR^*(X,\OO)) \to \widehat{H}^{2i}_{DR}(X)$.

A first direct consequence of the previous considerations is the following vanishing result
of Chern classes for perfect crystals along derived foliations.

\begin{cor}\label{cor-chernvanishing}
Let $X$ be a derived Artin stack locally of finite presentation over $k$, and
$\F \in \Fol(X/k)$ be a perfect derived foliation. Then, for any $E \in \Parf(\F)$, the chern 
classes $c_i(E) \in \widehat{H}^{2i}_{DR}(X)$ of the underlying perfect complex $E$ on $X$ possess
canonical lifts to $H^{2i}(F^{-i}_\F\CDR^*(X,\OO))$. It particular, the image of $c_i(E)$
by the natural moprhism
$$H^{2i}_{DR}(X) \to H^{2i}_{DR}(\F,\OO) \simeq H^{2i}(Gr_{F_\F}^0\CDR^*(X,\OO))$$
is zero as soon as $i>0$.
\end{cor}

A slightly more sophisticated consequence is the following result on characteristic classes
of the conormal complex of a derived foliation. Recall from Proposition \ref{pVII-4} that, for $X$ and $\F$ as above, 
the Hodge filtration $F^{*}_\F$ is a filtration on $\CDR^*(X,\OO)$ whose associated graded
is given, as a graded $E_\s$-algebra, by foliated de Rham cohomology with coefficients in the derived external products
of the conormal complex
$$Gr_{F^{*}_\F}^i\CDR^*(X,\OO) \simeq \CDR^*(\F,\LL\Lambda^{-i}\N^{*}_\F)[i]$$
where $\LL\Lambda^i\N^{*}_\F$ is the $i$-th derived exterior product of $\N^{*}_\F$, 
equipped with its natural structure of a perfect crystal along $\F$.
In particular, when $\N_\F$ is a vector bundle of rank $r$,  $\LL\Lambda^i\N^{*}_\F \simeq \Lambda^i\N^{*}_\F\simeq 0$ 
for all $i>r$. As a consequence, using the multiplicative nature of the Hodge filtration, we obtain
the following important corollary.

\begin{cor}\label{ccor-chernvanishing}
Let $X$ be a derived Artin stack locally of finite presentation over $k$, and $\F \in \Fol(X)$
be a perfect derived foliation whose conormal complex $\N^{*}_\F$ is a vector bundle of rank $r$. 
Then, for any homogenuous polynomial $P(x_1,\dots,x_r)$, with $deg(x_i)=-2i$, of total
degree $deg(P)=q < -r$, the element
$$P(c_1(\N^{*}_\F),\dots,c_r(\N^{*}_\F)) : k \to \CDR^*(X,\OO)$$
is canonically homotopic to zero.
\end{cor}

\textit{Proof.} The vector bundle $\N^{*}_\F$ has canonical structure of 
a crystal along $\F$, and thus, by Corollary \ref{cor-chernvanishing}, all the morphisms of dg-modules over $k$
$$c_i(\N^{*}_\F) : k \to \CDR^*(X,\OO)$$
possess canonical factorisations
$$c_i(\N^{*}_\F) \to F^{-i}_\F \CDR^*(X,\OO).$$
As the Hodge filtration is a filtration of $E_\s$-algebras, we deduce that 
the induced morphism
$P(c_1(\N^{*}_\F),\dots,c_r(\N^{*}_\F))$ 
has a canonical factorisation as
$$P(c_1(\N^{*}_\F),\dots,c_r(\N^{*}_\F)): k \to F^{-r-1}_\F\CDR^*(X,\OO) \to \CDR^*(X,\OO).$$
But, by assumption $\LL\Lambda^i\N^{*}_\F = 0$ for $i>r$, and thus $F^{-r-1}_\F\CDR^*(X,\OO) \simeq 0$.
\hfill $\Box$ \\

\section{Residues of singular foliations in positive characteristics}

In this section we exploit the results on characteristic classes defined and studied in the previous section, in the specific case
of a smooth variety over a perfect field of characteristic $p\neq 2$. Our objective is to prove
a crystalline version of the Baum-Bott residue theorem of \cite{baumbott}. \\

We let $k$ be a perfect field of characteristic $p>0$, and we assume $p\neq 2$.
We start by explaining the well-known relations between
the crystalline cohomology over $\WW(k)$ and derived de Rham cohomology.  
We let $S_n=\Spec\, \WW_n(k)$ be the spectrum of the ring of Witt vectors of length
$n$ over $k$, $\WW(k)=lim_n \WW_n(k)$ the corresponding filtered complete algebra, 
and $S:=\Spf\, \WW(k) = "colim"_n S_n$ its formal spectrum. In the sequel, in order to avoid confusions
with our notation $\DR(A/\WW_n(k))$, we will write
$dR(A/\WW_n(k))$ to denote the usual (underived) de Rham complex of a commutative
$\WW_n(k)$-algebra $A$
$$dR(A/\WW_n(k)) := \left ( A \to \Omega^1_{A/\WW_n(k)} \to \Omega^2_{A/\WW_n(k)} \to \dots \right )$$
where $A$ sits in degree $0$.

For a commutative simplicial $\WW_n(k)$-algebra $A$, which we assume to be 
a finite cell $\WW_n(k)$-algebra (so each $A_m$ is a polynomial algebra over $\WW_n(k)$), 
we consider the projection $A \to \pi_0(A)$ and the corresponding augmentation simplicial 
ideal $I \subset A$. We let $D(I)$ be the divided power envelope of $A$ along $I$ relative to $\WW_n(k)$
(endowed with its canonical divided power structure), so we have a commutative diagram 
of commutative simplicial $\WW_n(k)$-algebras
$$\xymatrix{
A \ar[r] \ar[rd] & D(I) \ar[d] \\
 & \WW_n(k).
}$$

The commutative $A$-algebra $D(I)$ carries a canonical flat connection,
and thus we can form its usual de Rham complex, which is the simplicial object in $\fdg_k$
$$[m] \mapsto D(I) \otimes_{A_m}dR(A_m/\WW_n(k)),$$
where the filtrations are here given by the $PD$-Hodge filtration of \cite[\S V-2.3.1]{MR384804} (convolution
of the standard Hodge filtration with the $PD$-adic filtration). We remind that the
$(-i)$-layer of this filtration is the subcomplex (for a fixed given $m$)
$$I_m^{[i]} \to I_m^{[i-1]}\otimes \Omega^1_{A_m/\WW_n(k)} \to I_m^{[i-2]}\otimes \Omega^2_{A_m/\WW_n(k)} \to \dots,$$
where $I_m^{[q]}$ denotes the $q$-th $PD$-powers of the ideal $I_m$ in $D(I_m)$. By taking the
geometric realization of this simplicial object (i.e. its colimit in the $\s$-category $\fdg_k$)
we obtain a well defined object $|dR(D(I)/\WW_n(k))| \in \fdg_{WW(k)}$ which is functorial in $A$. It defines
an $\s$-functor $A \mapsto |dR(D(I))|$, which by \cite[\S V-2.3.2]{MR384804}
is naturally equivalent to the functor sending $A$ to the complex of crystalline cohomology of 
the $\WW_n(k)$-algebra $\pi_0(A)$.

For any $A$ as above, the morphism $A \to D(I)$ induces a morphism on the de Rham complexes
$|dR(A)| \to |dR(D(I))|$
where the left hand side is the geometric realization of the simplicial object 
$[m] \mapsto dR(A_m/\WW_n(k))$. We can complete this morphism with respect to the filtrations on both
sides, and thus get an induced morphism
$$\widehat{|dR(A)|} \to \widehat{|dR(D(I))|},$$
which is functorial in $A$. The left hand side is naturally equivalent to the complex  $\CDR^*(X/S_n,\OO)$
computing derived de Rham cohomology of $X=\Spec\, A$ relative to $\WW_n(k)$, whereas the right hand
side is the $PD$-completed crystalline cohomology of $t_0(X)=\Spec\,\pi_0(A)$ relative to $\WW_n(k)$ denoted by
$$\widehat{\mathcal{C}}^*_{crys}(t_0(X)/\WW_n(k)) := \widehat{|dR(D(I))|}.$$
We thus have constructed a functorial diagram of objects in $\fdg_{\WW(k)}$
$$\xymatrix{
\CDR^*(X/S_n,\OO) \ar[r] &  \widehat{\mathcal{C}}^*_{crys}(t_0(X)/S_n) & \ar[l] 
\mathcal{C}^*_{crys}(t_0(X)/S_n) }$$
where the map on the right hand side is the completion morphism with respect to the $PD$-Hodge filtration.
This diagram of morphisms is not only functorial in $X$, but also functorial with respect to the
natural projections $\WW_{n+1}(k) \to \WW_n(k)$, so that we can pass to the limit over $n$
$$\xymatrix{
\CDR^*(X/S) \ar[r]^-{\phi} & \widehat{\mathcal{C}}^*_{crys}(t_0(X)/S)  & \ar[l] 
 \mathcal{C}^*_{crys}(t_0(X)/S). }$$
 
Finally, the functoriality in $X$ can be used to glue these constructions to get a global 
version valid for any derived stack $X$ which is locally of finite presentation over $S$ (i.e. 
a compatible system of derived Artin stacks $X_n$ of finite presentation over $S_n$). Moreover, 
we could also keep track of the multiplicative structures involved and contemplate the natural
lift of the previous diagram to the $\s$-category of filtered $E_\s$-algebras over $\WW(k)$.
Note also that, as $\widehat{\mathcal{C}}^*_{crys}(s/S)\simeq \WW(k)$, the
morphism $\phi$ enters in a commutative diagram of complete filtered $E_\s$-algebras
$$\xymatrix{
\CDR^*(X/S) \ar[r] & \widehat{\mathcal{C}}^*_{crys}(X/S) \\
\CDR^*(s/S) \ar[r] \ar[u] & \WW(k) \ar[u]
}$$
and thus produces a morphism
$$\CDR^*(X/S) \otimes_{\CDR^*(s/S)}\WW(k) \to \mathcal{C}^*_{crys}(X/S)$$
(here we have noted $s=\Spec\, k$).

\begin{prop}\label{compareDRCrys}
Let $X$ be a smooth variety over $k$, considered as a derived scheme of finite presentation 
over $S$ via the natural morphism $s = \Spec\, k \to S$. 
\begin{enumerate}
    \item The completion morphism with respect to the $PD$-Hodge filtration
$$\mathcal{C}_{crys}^*(X/S) \to \widehat{\mathcal{C}}^*_{crys}(X/S)$$    
is an equivalence. 

    \item The morphism $\phi$ constructed above 
    $$\CDR^*(X/S) \otimes_{\CDR^*(s/S)}\WW(k) \to  \widehat{\mathcal{C}}^*_{crys}(X/S)$$
is an equivalence.
\end{enumerate}
\end{prop}

\textit{Proof.} Assertion $(1)$ follows from the fact that the $PD$-Hodge filtration
is already complete (in fact it is finite) when $X$ is smooth over $k$. Indeed, by descent
we can assume that $X$ is affine to prove $(1)$. We can then compute the crystalline 
cohomology of $X$ over $S_n$ by choosing a smooth lift $X_n$ over $S_n$.
The scheme $X_n$ is the $PD$-hull of $X$ relative to $S_n$, and the $PD$-adic
filtration is then finite as this can be seen using 
\cite[\S I-3.2.4]{MR384804} (we use here that $p\neq 2$). 
We thus have that the $PD$-Hodge filtration is finite
and thus complete. 

It remains to prove $(2)$.
By construction 
$$\CDR^*(X/S) \otimes_{\CDR^*(s/S)}\WW(k) \to  \widehat{\mathcal{C}}^*_{crys}(X/S)$$
is a morphism of complete $\WW(k)$-dg-modules, and thus to check that it is
an equivalence it is enough to check that it is so after base change along the augmentation
morphism $\WW(k) \to k$. After such a base change, the morphism $\psi$ can be written as
$$\CDR^*(X \times s \times_S s/s) \otimes_{\CDR^*(s \times_S s/s)} k \longrightarrow \CDR^*(X/s).$$
By K\"unneth, this is also written as the canonical morphism
$$\CDR^*(X/s) \otimes_k \CDR^*(s \times_S s/s) 
\otimes_{\CDR^*(s \times_S s/s)} k \longrightarrow  \CDR^*(X/s)$$
which is obviously an equivalence.
\hfill $\Box$ \\

By combining Corollary \ref{ccor-chernvanishing} and Proposition 
\ref{compareDRCrys} we finally find the following corollary, which is 
our version of the Baum-Bott's residue theorem in characteristic $p>0$.
We let $X$ be a smooth variety over $k$. For a derived foliation $\F$ on $X$ relative to $s=\Spec\, k$,
we will call an \emph{$S$-structure on $\F$} the data of a derived foliation $\F' \in \Fol(X/S)$
lifting $\F$ via the pull-back $i : \Fol(X/S) \to \Fol(X/s)$ along $i : s \hookrightarrow S$.

\begin{cor}
Let $X$ be a smooth variety over $k$, and $D \subset \Omega^1_{X/k}$ be a coherent subsheaf 
which is a differential ideal: $dR : D \to \Omega^2_{X/k}$ factors through 
the image of $D \otimes \Omega^1_{X/k} \to \Omega^2_{X/k}$. Let $U \subset X$ be a dense open 
on which $D$ is a sub-bundle of rank $q$ of $\Omega_{X/k}^1$, and $Z$ its reduced closed complement. 
Let us fix an $S$-structure on the derived foliation defined by $D_{|U} \subset \Omega_{U/k}^1$
(by Proposition \ref{pVII-3}).

Then, for any homogeneous polynomial $P \in \ZZ[c_1,\dots,c_q]$ of degree $d<-q$ (with $deg(c_i)=-2i$), 
there exists a canonical element in local cohomology supported on $Z$
$$Res_{D,P} \in H^{-d}_{Z,crys}(X/S)$$
whose image in $H^{-d}_{crys}(X/S)$
coincide with $P(c_1(D),\dots,c_q(D))$.
\end{cor}

Note that the residue $Res_{D,P}$ does depend on the choice of the $S$-structure on 
the foliation defined by the differential ideal $D_{|U}$. Note also that when $U=X$, that is
when $D$ is a sub-bundle of rank $q$ in $\Omega_{X/k}^1$ then we find that $P(c_1(D),\dots,c_q(D))$
vanishes in cyrstalline cohomology $H^{d}_{crys}(X/S)$. Without the condition on the 
existence of an $S$-structure, we would only have that this class is divisible by $p$.

Finally, we note that if both $X$ and $D$ lifts to characteristic zero, then 
the existence of the $S$-structure is automatic (and the lift provides a canonical one). However, 
$S$-structures can exist without $X$ being liftable to $S$. For instance, the tautological
foliation $*_{X/s} \in \Fol(X/s)$ always possesses a canonical lift, namely $*_{X/S}$, even
if $X$ does not itself lift to $S$.
In a way, the existence of an $S$-structure should be thought as the existence of a lift 
of an hypothetical leaf space of the foliation from $s$ to $S$.

\section{Infinitesimal derived foliations}

In this section we present an alternative notion of derived foliations in arbitrary characteristics. 
In characteristic zero this notion is equivalent to our original notion of derived foliations, by 
means of the famous \emph{red shift self equivalence} of the $\s$-category of graded mixed
complexes. However, outside of characteristic zero this notion is different 
from the one presented previously in this chapter (Definition \ref{dVII-4}). We call this alternative notion \emph{infinitesimal
derived foliations}, as its relation with infinitesimal cohomology is similar to the
relation between derived foliation and de Rham cohomlogy. It is also very closely related to
the notion of \emph{partition Lie algebras and Lie algebras} and we refer to 
\cite{jfuarxiv} for a precise
comparison statement.

We will show that any infinitesimal derived foliation possesses an underlying derived foliation 
in the sense of \ref{dVII-4}. It is important to understand that infinitesimal derived foliations are
somehow more restricted and could be thought as derived foliations together with some
extra structures/properties. For instance, infinitesimal derived foliations should also be closely 
related to  p-restricted foliations of \cite{MR927978}. \\

We start by considering a derived group stack $\cH_\pi$, which will be playing a similar
role than  the group stack $\cH$ considered in Section \ref{sec-gradedcircle}. 
For this we consider $\Ga[-1]:=\Omega_0\Ga$, the
derived group scheme (over $k$) defined as the derived loop space of the additive group $\Ga$
taken at $0$. By definition, we have
$$\Ga[-1] \simeq \Spec\, k \times_{\Ga}\Spec\, k \simeq \Spec\, k\oplus k[1],$$
where $k\oplus k[1]$ is the simplicial commutative ring obtained as the 
trivial square zero extension of $k$ by $k[1]$. 
The natural $\Gm$-action on $\Ga$, which by 
our convention is of weight $-1$ (see right before Definition \ref{dI-5}), induces an action of $\Gm$ 
on the derived 
group scheme $\Ga[-1]$. Algebraically, this $\Gm$ action on $\OO(\Ga[-1]))=k\oplus k[1]$ is such that 
$k[1]$ is pure of weight $1$.

\begin{df}
The \emph{derived group scheme $\cH_\pi$} is defined to be the semi-direct product of $\Ga[-1]$ by $\Gm$
$$\cH_\pi:=\Gm \ltimes \Ga[-1].$$
\end{df}

According to the above definition, the underlying derived scheme  of $\cH_\pi$ is 
$\Spec\, (k\oplus k[1])[t,t^{-1}]$, i.e. 
the spectrum of the Laurent polynomials ring with coefficients in 
$k\oplus k[1]$. The group structure on $\cH_\pi$ provides $(k\oplus k[1])[t,t^{-1}]$
with a structure of a cogroup object inside the $\s$-category of simplicial commutative $k$-algebras, 
or in other words as a $k$-linear commutative Hopf simplicial algebra. If we write
$k\oplus k[1]$ as $k[\epsilon]$, with $\epsilon$ a free generator in homotopical degree $1$ 
(so cohomological degree $-1$), this makes $k[\epsilon][t,t^{-1}]$ into a 
a simplicial commutative Hopf algebra for the usual \emph{twisted} comultiplication
$$\Delta(\epsilon)=\epsilon\otimes 1 + t \otimes \epsilon \qquad
\Delta(t)=t\otimes t.$$

Similarly to Proposition \ref{pVII-2}, we have the following result, identifying 
graded mixed complexes with representations of the group $\cH_\pi$.

\begin{prop}\label{pVII-5}
There exists a natural equivalence of $k$-linear symmetric monoidal $\s$-categories
$$\egrdg_k \simeq \QCoh(B\cH_\pi).$$
\end{prop}

\textit{Proof.} The proof is similar to \cite[Prop. 4.2.3 (ii)]{mrt}, and in fact easier as the
group $\cH_\pi$ is an affine derived group scheme, which simplifies some of the arguments.
We start by writing $B\cH_\pi$ as the standard geometric realization of the
simplicial object $[n] \mapsto \cH_\pi^n$. By descent we have a canonical equivalence of 
symmetric monoidal $\s$-categories
$$\QCoh(B\cH_\pi) \simeq \lim_n \OO(\cH_\pi)^{\otimes n}-\dg_k.$$
We observe that the simplicial object $[n] \to \cH_\pi^n$ is the nerve of the
zero section $\Spec\, k \to [\mathbb{A}^1/\Gm]$, so we have, by pull-back, 
a symmetric monoidal $\s$-functor
$$\QCoh([\mathbb{A}^1/\Gm]) \longrightarrow \lim_n \OO(\cH_\pi)^{\otimes n}-\dg_k \simeq
\QCoh(B\cH_\pi).$$
By the comparison between graded mixed complexes and complete filtered 
complexes (see Corollary \ref{clI-4}) the Tate realization provides a 
lax symmetric monoidal full embedding
$$\egrdg_k \hookrightarrow \QCoh([\mathbb{A}^1/\Gm])$$
and thus by composition a lax symmetric monoidal $\s$-functor
$$\psi : \egrdg_k \longrightarrow \QCoh(B\cH_\pi).$$
We now prove that $\psi$ is a symmetric equivalence. For this, we start by noticing that 
the lax symmetric monoidal structure of $\psi$ is in fact strict. For this we consider the
natural projection $\cH_\pi \rightarrow \Gm$, and the induced morphism
$p : B\cH_\pi \to B\Gm$ on classifying derived stacks. The composition of $\psi$ by the
push-forward 
$p_* : \egrdg_k \to \QCoh(B\Gm)\simeq \grdg_k$, is easily seen to be equivalent to the
$\s$-functor sending a graded mixed objects $E$ to the graded
object $E^{(*)}[-2*]$. In particular, its lax symmetric monoidal
structure is strict. Moreover, as $p_* : \QCoh(B\cH_\pi) \to \QCoh(B\Gm)$ is clearly 
conservative (because $p$ is induced by an affine morphism), 
we deduce that the lax symmetric monoidal structure on $\psi$ is a symmetric monoidal structure. 

It remains to show that the $\s$-functor underlying $\psi$ is an equivalence of $\s$-categories.
For this we come back to the situation
$$\egrdg_k \hookrightarrow \QCoh([\mathbb{A}^1/\Gm]) \longrightarrow \QCoh(B\cH_\pi).$$
By Corollary \ref{clI-4}, the image of the first full embedding is the full sub-$\s$-category of objects
which are complete at $0$. We thus have to show that 
$\QCoh([\mathbb{A}^1/\Gm]) \longrightarrow \QCoh(B\cH_\pi)$ identifies the full sub-$\s$-category
of complete modules with $\QCoh(B\cH_\pi)$. This follows easily from the
fact that $B\cH_\pi$ is the realization of the nerve of the zero section 
$B\Gm \hookrightarrow [\mathbb{A}^1/\Gm]$ as shown by the following lemma.

\begin{lem}
Let $A$ be a $k$-algebra and $f\in A$ be a non-zero divisor. 
Let $B^*$ be the cosimplicial $E_\s$-algebra given by the derived conerve
of $A \to A/f$
$$B^n := A/f\otimes_A A/f \otimes_A \dots \otimes_A A/f.$$
Then the descent $\s$-functor
$$\dg_A \to \lim_n \dg_{B^n}$$
induces an equivalence 
$$\hat{\dg}_{A,f} \simeq \lim_n \dg_{B^n}$$
where $\hat{\dg}{A,f}$ is the $\s$-category of $f$-complete $A$-dg-modules ($E$ such that 
$E \simeq \lim_n (E/f^n)$).
\end{lem}

\textit{Proof of the lemma.} The descent $\s$-functor $\phi : \dg_A \to \lim_n \dg_{B^n}$
sends an $A$-dg-module to the pull-back $E\otimes_A B^*$. It has a right 
adjoint $\psi$ sending $E^* \in \lim_n \dg_{B^n}$ to $\lim_n E^n \in \dg_A$.

We can form the cocartesian square of commutative 
cosimplicial $E_\s$-algebras
$$\xymatrix{
A \ar[r] \ar[d] & A/f \ar[d] \\
B^* \ar[r] &C^*.}$$
Then, $A/f \to C^*$ is the derived conerve of 
the morphism $A/f \to A/f\otimes_A A/f$ which possesses a retraction given by the
multiplication in $A/f$. 
Therefore, 
the cosimplicial object $C^*$ is simplicially contractible. This implies that the descent 
$\s$-functor
$$\dg_{A/f} \to \lim_n(\dg_{C^n})$$
is in fact an equivalence of $\s$-categories. We consider now 
the commutative diagram of $\s$-categories obtained by base change from the above diagram
$$\xymatrix{
\dg_A \ar[d] \ar[r] &  \lim_n(\dg_{B^n})\ar[d] \\
 \dg_{A/f} \ar[r] & \lim_n(\dg_{C^n}).
}$$
In this diagram the vertical $\s$-functors are conservative and commutes with limits. Indeed, 
$A/f$ is a perfect $A$-dg-module and thus $(A/f)\otimes_A -$ commutes with limits. Similarly, 
each $C^n$ is perfect over $B^n$. As a result, the diagram is right adjointable in the
sense of \cite[Def. 4.7.4.13]{HA}, because of the explicit description of the right adjoints
of the horizontal $\s$-functors in terms of limits. As the bottom horizontal $\s$-functor is 
an equivalence, so is the top horizontal $\s$-functor, because the unit of counit of the top
adjunction can be checked to be equivalences after applying the vertical conservative
$\s$-functors (so that they become the unit and counit of the bottom adjunction).
\hfill $\Box$ \\

We can apply the previous lemma to the natural simplicial objects representing $B\Gm \hookrightarrow 
[\mathbb{A}^1/\Gm]$. This finishes the proof of the proposition.
\hfill $\Box$ \\

In parallel to the notion of graded loop space we introduce the \emph{infinitesimal} 
graded loop space as follows.

\begin{df}
Let $X$ be a derived stack. The \emph{infinitesimal graded loop space of $X$}
is defined by 
$$\cL^{gr}_{\pi}(X):=\uMap(\Ga[-1],X) \in \dSt.$$
If $X \to Y$ is a morphism of derived stacks, the \emph{relative infinitesimal graded loop space of $X$
over $Y$} is defined by
$$\cL^{gr}_{\pi}(X/Y):=\cL^{gr}_{\pi}(X)\times_{\cL^{gr}_{\pi}(Y)}Y,$$
for the natural morphism $Y \to \cL^{gr}_{\pi}(Y)$ corresponding to the projection $Y\times \Ga[-1] \to Y$.
\end{df}

The group $\Gm$ acts on $\Ga[-1]$, and  
the group $\Ga[-1]$ acts on itself by translation. Therefore we get a natural action of $\cH_\pi=\Gm \ltimes \Ga[-1]$ on the stack $\Ga[-1]$, and thus
on the infinitesimal graded loop space of any derived stack.  The object
$\cL^{gr}_{\pi}(X)$ will always be considered as an object in the $\s$-topos $\dSt^{\cH_\pi}$
of derived stacks endowed with a $\cH_\pi$-action (and $\cH_\pi$-equivariant maps). 

The following proposition identifies $\cL^{gr}_{\pi}(X)$ with its $\Gm$-action when $X$
is a derived Artin stack. For this we consider the canonical point $* \to \Ga[-1]$ given by the unit, 
and its induced projection $\cL^{gr}_\pi(X) \to X$. Note that this projection is $\Gm$-equivariant, but 
is not $\Ga[-1]$ equivariant as the $\Ga[-1]$-action by translation does not fix the base point.
The derived stack $\cL^{gr}_\pi(X) \to X$ thus lives in the $\s$-topos $\dSt_X^{\Gm}$ of $\Gm$-equivariant
objects over $X$. Recall from \cite{mon} that this category contains as a full sub-$\s$-category 
the $\s$-category of derived linear stacks over $X$, which is equivalent to $\QCoh^-(X)$ of 
eventually coconnective quasi-coherent complexes on $X$. The full-embedding
$\QCoh^-(X) \hookrightarrow \dSt_X^{\Gm}$ is given by the construction $\VV$, sending
$E \in \QCoh^-(X)$ to $\VV(E) \to X$ to the $\s$-functor sending an affine $u : \Spec\, A \to X$
to the space $Map_{\dg_A}(u^*(E),A)$ (see Definition \ref{dlinearstack}).

\begin{prop}
Let $X$ be a derived Artin stack. The $\Gm$-equivariant derived stack $\cL^{gr}_{\pi}(X) \to X$
is canonically equivalent to the linear stack $\VV(\LL_X[-1]) \to X$, associated to the
quasi-coherent complex $\LL_X[-1]$.
\end{prop}

\textit{Proof.} As $\Ga[-1]=\Spec\, k \oplus k[1]$, this is simply the definition 
of $\LL_X$ in terms of points with values in square zero extensions (see \cite[Def. 1.4.1.5]{hagII}). 
\hfill $\Box$ \\

The terminology \emph{infinitesimal graded loop space} has been chosen because of the theorem below, 
expressing that $\cH_\pi$-equivariant functions on $\cL^{gr}_\pi(X)$ computes infinitesimal 
cohomology of $X$, at least in the smooth case. Recall that if $X \to S=\Spec\, k$ is a smooth morphism 
of schemes, we can define $(X/S)_{inf}$ as the functor
sending an affine $\Spec\, A \to S$ to $X(A_{red})$, where $A_{red}$ is defined
to be the discrete ring $\pi_0(A)_{red}$. The functor $(X/S)_{inf}$ is a derived stack over $S$, and
by definition infinitesimal 
cohomology of $X$ relative to $S$ is defined by 
$$C^*_{inf}(X/S):=C^*((X/S)_{inf},\OO),$$
which will be considered as a $k$-linear $E_\s$-algebra.
Because $X$ is smooth over $S$, the canonical morphism $X \to (X/S)_{inf}$ is an epimorphism
of derived stacks (because of the infinitesimal lifting property). Therefore, $H^*_{inf}(X/S)$ can also be
computed as the limit
$$C^*_{inf}(X/S) \simeq \lim_{n}(C^*(X^{(n)},\OO)),$$
where $X^{(*)}$ is the the nerve of $X$ over $(X/S)_{inf}$, which can also be
described as the nerve of the formal groupoid $X^{(2)} \leftleftarrows X$ where
$X^{(2)}$ is the formal completion of $X \times_S X$ along the diagonal.

In order to state the theorem below, let us consider the quotient map
$$p : [\cL^{gr}_\pi(X/S)/\cH_\pi] \to S$$
The direct image of the structure sheaf $p_*(\OO)$ is by definition
the $\cH_\pi$-equivariant cohomology of $\cL^{gr}_\pi(X/S)$ with coefficients in the
structure sheaf $\OO$. It will be symbolicaly denoted by $C^*(\cL^{gr}_\pi(X/S))^{\cH_\pi}$.

We produce a natural morphism $C^*_{inf}(X/S) \to C^*(\cL^{gr}_\pi(X/S))^{\cH_\pi}$
as follows. The natural morphism $X \to (X/S)_{inf}$ induces a morphism on infinitesimal
graded loop spaces
$\cL^{gr}_\pi(X/S) \to \cL^{gr}_\pi((X/S)_{inf}/S)$. By the definition of $(X/S)_{inf}$, 
and from the fact that the reduced group $\cH_\pi$ is trivial, it is easy to see that 
the constant map morphism $(X/S)_{inf} \to \cL^{gr}_\pi((X/S)_{inf}/S)$ is an equivalence. 
In particular, we get this way a natural morphism of derived stacks with $\cH_\pi$-actions
$\cL^{gr}_\pi(X/S) \to (X/S)_{inf}$. By considering cohomology with coefficients in $\OO$ we 
finally get the required morphism
$$\phi_{X/S} : C^*_{inf}(X/S) \longrightarrow C^*(\cL^{gr}_\pi(X/S))^{\cH_\pi}.$$
By construction it is a morphism of $E_\s$-algebras functorial in $X/S$. Note that 
the above morphism always exists, for a general morphism of derived stacks $X \to S$ without 
assuming any extra conditions. 

\begin{thm}\label{tcompareinf}
Let $X \to S=\Spec\, k$ be a smooth scheme. Then, the
above morphism
$$\phi_{X/S} : C^*_{inf}(X/S) \simeq C^*(\cL^{gr}_\pi(X/S))^{\cH_\pi}$$
is an equivalence.
\end{thm}

\textit{Proof.} We consider the commutative diagram of derived stacks with $\cH_\pi$-actions
$$\xymatrix{
X \ar[r]^-{id} \ar[d] & X \ar[d] \\
\cL^{gr}_\pi(X/S) \ar[r]_-{\phi_{X/S}} & (X/S)_{inf}.
}$$
This diagram induces a morphism of simplicial objects on the corresponding nerves of the
vertical morphisms. The nerve of the morphism $X \to \cL^{gr}_\pi(X/S)$ can be described by the
simplicial object $[n] \mapsto \cL^{gr}_\pi(X/X^n)$, where $X$ sits inside $X^n$ by the diagonal embedding.
At the simplicial level $n$, the morpshim induced on the nerve is
thus an $\cH_\pi$-equivariant morphism of stacks over $S$
$$\psi_{X,n} : \cL^{gr}_\pi(X/X^n) \longrightarrow X^{(n)},$$
where $X^{(n)}$ is the formal completion of $X^n$ along its diagonal. 

\begin{lem}\label{lem-ausiliario}
The morphism $\psi_{X,n}$ induces by pull-backs an equivalence of $E_\s$-algebras
$$C^*(X^{(n)},\OO) \simeq C^*(\cL^{gr}_\pi(X/X^n))^{\cH_\pi}.$$
\end{lem}

\textit{Proof Lemma \ref{lem-ausiliario}.} By descent we can easily restrict to the case where
$X=\Spec\, A$ is a smooth affine scheme over $S=\Spec\, k$. In this case,
$C^*(X^{(n)},\OO)$ identifies canonically with $\hat{A}^{\otimes_k n}$, the formal
completion of $A^{\otimes_k n}$ along the canonical augmentation $A^{\otimes_k n} \to A$.
Similarly, $\cL^{gr}_\pi(X/X^n)$ is canonically identified with 
the smooth affine scheme $\Spec\, Sym_A((\Omega_{A/k}^1)^{n-1})$. Therefore, the morphism 
$\psi_{X,n}$ is a morphism of discrete commutative algebras
$$\psi_{A,n} : \hat{A}^{\otimes_k n} \longrightarrow |Sym_A((\Omega_{A/k}^1)^{n-1})|^t,$$
where the right hand side is the Tate realization of the graded mixed object $Sym_A((\Omega_{A/k}^1)^{n-1})$
induced by the $\cH_\pi$-action on $\cL^{gr}_\pi(X/X^n)$. The morphism $\psi_{A,n}$ is
compatible with the canonical augmentations to $A$, and thus is a morphism of complete
filtered commutative algebras. Therefore, in order to prove it is an isomorphism it is enough to show that it induces an isomorphism on the associated graded objects. This induced morphism is 
a morphism 
$$gr(\psi_{A,n}) : Gr^*(\hat{A}^{\otimes_k n})\longrightarrow Sym_A((\Omega_{A/k}^1)^{n-1}).$$
Moreover, because $A$ is a smooth $k$-algebra, the left hand side is canonically isomorphic
to $Sym_A(I_n/I_n^2)$, where $I_n$ is the augmentation ideal. Finally, 
we have the usual natural isomorphism of $A$-modules $I_n/I_n^2\simeq (\Omega_{A/k}^1)^{n-1}$, and thus
the morphism $gr\psi_{A,n}$ can be identified with a graded endomorphism
of the graded algebra $Sym_A((\Omega_{A/k}^1)^{n-1})$. As it is compatible with augmentations to $A$
it is thus determined by an $A$-linear endomorphism of $(\Omega_{A/k}^1)^{n-1}$. By functoriality in 
$A$ and in the simplicial direction $n$, we see moreover that it is determined by a functorial
endormorphism the $A$-module $\Omega_{A/k}^1$, and thus must be the multiplication by 
an element $\lambda \in k$. To conclude the lemma, we have to show that 
$\lambda=1$. 

By functoriality, it is easy to reduce to the case $A=k[T]$ is a polynomial ring, and by base change we
can even assume that $k=\ZZ$. In order to check that $\lambda=1$ we can base change from $\ZZ$ to $\QQ$, and
thus reduce to the case $A=\QQ[T]$ and $k=\QQ$. As we are now in characteristic zero, 
$C^*(\cL^{gr}_\pi(X/X^n))^{\cH_\pi}$ identifies canonically with the 
completed de Rham complex $\CDR^*(X/X^n)$. Indeed, the red shift equivalence sends 
the graded mixed cdga $\OO(\cL^{gr}_\pi(X/X^n))$ to $\OO(\cL^{gr}(X/X^n))$. Finally, 
the morphism $\psi_{X,n}$ is now the natural isomorphism relating 
derived de Rham cohomology of a closed embedding with functions on the formal 
completion (see for instance Proposition \ref{pI-7}). \hfill $\Box$ \\

The lemma implies the theorem by descent along the epimorphisms $X \to \cL^{gr}_\pi(X/S)$ and 
$X \to (X/S)_{inf}$. \hfill $\Box$ \\

\subsection{Infinitesimal derived foliations}

As before, we fix $k$ a base simplicial commutative ring. \\

\begin{df}\label{dinfol}
Let $X=\Spec\, A$ be a derived affine scheme over $k$. A (perfect) \emph{infinitesimal derived foliation
$\F$ on $X$ (relative to $k$)} consists of a $\cH_\pi$-equivariant derived stack $\VV(\F)$ 
together with a $\cH_\pi$-equivariant morphism $p : \VV(\F) \to \cL^{gr}_\pi(X/k)$, such that, 
as a $\Gm$-equivariant stack over $X$, $\VV(\F)$ is a perfect  linear derived affine stack (i.e. of the
form $\VV(E)$ for $E$ a perfect complex on $X$).
\end{df}

The infinitesimal derived foliations on $X$ (relative to $k$) form an 
$\s$-category $\Fol_\pi(X/k)$, which by definition is a full sub-$\s$-category 
of $(\dSt^{\cH_\pi}_k)/\cL^{gr}_\pi(X/k)$, the $\cH_\pi$-equivariant derived stacks over
$\cL^{gr}_\pi(X/k)$.  The construction $X \mapsto \Fol_\pi(X/k)$ can be promoted to 
an $\s$-functor $\Fol_\pi(-/k) : \dAff_k^{op} \to \scat$, by
defining pull-backs as follows: for a morphism $f : X \to Y$ of affine derived $k$-schemes, 
we define $f^*$ to be induced by the base change
$$(\dSt^{\cH_\pi}_k)/\cL^{gr}_\pi(Y/k) \longrightarrow (\dSt^{\cH_\pi}_k)/\cL^{gr}_\pi(X/k),$$
along the induced morphism $\cL^{gr}_\pi(f/k): \cL^{gr}_\pi(X/k) \to \cL^{gr}_\pi(Y/k)$. As the
$\s$-functor $\Fol_\pi(-/k)$ clearly is a hyperstack for the fpqc topology, it 
can be uniquely extended to a colimit preserving $\s$-functor
$$\Fol_\pi(-/k) : \dSt_k \to \scat^{op}.$$

Let $X$ be a derived affine $k$-scheme, and $\F \in \Fol_\pi(X/k)$ be an infinitesimal derived foliation
on $X$. By \cite{mon}, the linear stack $\VV(\F) \to X$ is of the form $\VV(\LL_{\F/k}[-1])$, for some uniquely 
defined perfect complex $\LL_{\F/k}$ on $X$. 

\begin{df}\label{dcot}
With the notations above, the complex $\LL_{\F/k}$ is called the \emph{cotangent complex
of $\F$ (relative to $k$)}. It is also simply denoted by $\LL_\F$ if $k$ is clear from the context.
\end{df}

When $X$ is a general derived stack and $\F \in \Fol_\pi(X/k)$, it is also 
possible to define the cotangent complex $\LL_{\F/k}$, as a quasi-coherent
complex over $X$ following the same lines as Definition \ref{dII-5}. We have first a big cotangent
complex $\LL_\F^{big}$, which is an $\OO_X$-dg-module over the big site of derived
affine schemes over $X$: it sends $u : \Spec\, A \to X$ to the $A$-dg-module
$\LL_{u^*(\F)/k}$. This $\OO_X$-dg-module is not quasi-coherent and we thus
apply the quasi-coherator construction of \S \ref{secapp:OXmodulesquasi-coherentandperfectmodules} to get the desired cotangent complex
$$\LL_{\F/k} := (\LL_{\F/k}^{big})^{qcoh}.$$

The special case where $X$ is a derived \emph{Deligne-Mumford} 
stack simplifies, as the cotangent complexes of infinitesimal foliations are
stable by \'etale pull-backs, and thus will glue globally on $X$. This applies in particular
to the case where $X$ is any derived scheme. In these cases, we have
the natural morphism $\LL_{\F/k}^{big} \to \LL_{\F/k}$
restricts to an equivalence on the small \'etale site of $X$. \\

Let $X$ be a derived affine scheme of finite presentation over $k$.
The $\s$-category $\Fol_\pi(X/k)$ admits a final
object, namely the identity morphism 
$\cL^{gr}_\pi(X/k) \to \cL^{gr}_\pi(X/k)$. It is denoted by $*_{X/k}$, and its cotangent 
complex is naturally equivalent to $\LL_{X/k}$, the cotangent complex of $X$. Indeed, 
it is easy to see, simply by contemplating the functors of points and using the
very definition of the cotangent complex, that $\uMap(G_0,X) \simeq \VV(\LL_{X/k}[-1])$.
Similarly, the initial object in $\Fol_\pi(X/k)$ exists, is denoted by $0_{X/k}$, and
its cotangent complex is $0$. It corresponds to the
constant-map morphism $\VV(0_{X/k}) = X \longrightarrow \cL^{gr}_\pi(X/k)$. 

\subsection{Infinitesimal cohomology of infinitesimal derived foliations}

Let $X$ be a derived Artin stack of finite presentation over $k$,  
and $\F \in \Fol_\pi(X/k)$ be an infinitesimal derived foliation on $X$. We define
the $\s$-category of quasi-coherent complexes along $\F$ as follows. 
We first consider $\QCoh([\VV(\F)/\cH_\pi])$, the $\s$-category of quasi-coherent complexes
on the quotient stack $[\VV(\F)/\cH_\pi]$ (which is also, by definition, the $\s$-category
of $\cH_\pi$-equivariant quasi-coherent complexes on $\VV(\F)$). The $\s$-category
$\QCoh(\F)$ is defined as the full sub-$\s$-category of $\QCoh([\VV(\F)/\cH_\pi])$
whose objects $E$ are graded free: the pull-back of $E$ on $[\VV(\F)/\Gm]$ is of the form
$p^*(E_0)$ where $p : [\VV(\F)/\Gm] \to X$ is the natural projection. 

In the definition above we use the identification $\QCoh(B\cH_\pi)$ with the
$\s$-category of graded mixed complexes (see Proposition \ref{pVII-5}).

\begin{df}\label{dinf}
With the notations above, 
the \emph{$\s$-category of quasi-coherent complexes along $\F$} is $\QCoh(\F)$. For $E \in 
\QCoh^\F(X)$, the \emph{(derived and Hodge completed) the infinitesimal cohomology of $X$ along $\F$ with coefficients in $E$} 
is defined to be
$$\Cinf^*(\F,E):=|q_*(E)|^t \in \cfdg_k$$
where $q : [\VV(\F)/\cH_\pi] \longrightarrow B\cH_\pi$ is the structural map.
\end{df}

The definition above is justified by our Theorem \ref{tcompareinf}, stating that 
when $X$ is smooth over $k$, and $\F=*_{X/k}$ and $E=\OO_X$, $\Cinf^*(\F,E)$ coincides with 
infinitesimal cohomology defined using the usual infinitesimal site. 
For a general
derived Artin stack $X$, $\Cinf(*_{X/k},\OO)$ computes the \emph{completed derived 
infinitesimal cohomology of $X$} (relative to $k$). Indeed, it is easy to see the $\s$-functor
$X \mapsto \OO(\cL^{gr}_\pi(X/k))$, from affine derived schemes to 
graded mixed complexes, is equivalent to its extension by sifted colimits from 
smooth $X$'s. For $A \in \sCR$, where each $A_n$ is a smooth algebra over $k$, we thus have
$$\Cinf(*_X,\OO_X) \simeq \mathrm{colim}_{[n]} C^*((Spec\, A_n)_{\mathrm{inf}},\OO),$$
where the colimit is computed in the $\s$-category of \emph{complete} filtered complexes. The right
hand side of this equivalence is, by definition, the completed derived infinitesimal 
cohomology of $X$, and is denoted by $\mathbb{R}\hat{C}^*(X_{\mathrm{inf}},\OO)$. This notion is
extended to any derived Artin stack locally of finite presentation over $k$ by the usual gluing formula 
$$\mathbb{R}\hat{C}^*(X_{\mathrm{inf}},\OO):=\lim_{U \to X} \mathbb{R}\hat{C}^*(U_{\mathrm{inf}},\OO)$$
where the limit runs over affine $U$ locally of finite presentation.

\begin{cor}\label{c1}
For a general derived Artin stack $X$ of finite presentation over $k$, 
$\Cinf(*_{X/k},\OO)$ recovers the Hodge completed derived infinitesimal 
cohomology of $X$ relative to $k$
$$\Cinf(*_{X/k},\OO) \simeq \mathbb{R}\hat{C}^*(X_{\mathrm{inf}},\OO).$$
\end{cor}

\subsection{Integrability by formal groupoids}

As opposed to derived foliations, infinitesimal derived foliations turn out to be 
\emph{formally integrable} in a very strong sense. The most general integrability statements can be
found in \cite{brantner2025formalintegrationderivedfoliations}. 
Below, we prove a special easy case by elementary means, 
showing that the $\s$-category of smooth infinitesimal derived foliations over a $k$-scheme $X$
is in fact \emph{equivalent} to the category of \emph{formally smooth formal groupoids} over $X$. \\

Let $X$ be a smooth $k$-scheme and $\F\in \Fol_\pi(X/k)$. We say that $\F$ is \emph{smooth} if 
$\LL_\F$ is a vector bundle over $X$. In this case, we consider $X \times_{\VV(\F)} X$, endowed
with its natural $\cH_\pi$-action and with its first projection to $X$. Note that, as a graded
derived stack, $X \times_{\VV(\F)} X$ is of the form $\VV(\LL_{F})$, and thus is a the total 
space of the tangent bundle $\T_\F$ of $\F$. Also, the natural projection
$X\times_{\VV(\F)} X \to X$ is $\cH_\pi$-equivariant. As a result, 
$X \times_{\VV(\F)} X$ is a new object in $\Fol_\pi(X)$, that will be simply denoted by $\Omega_X\F$, and
called the \emph{loop space} of $\F$. Being the fiber product of the canonical morphism
$0_X \to \F$, the loop space $\Omega_X\F$ comes equipped with a natural structure of a groupoid object in $\Fol_\pi(X)$, whose object of objects is $0_X \in \Fol_\pi(X/k)$. Note that as
$0_X$ is not the final object $\Omega_X\F$ is not a group object and merely a groupoid.

\begin{prop}\label{p2}
The full sub-$\s$-category of $\Fol_\pi(X/k)$ formed by objects $\F$ whose cotangent complex
is of the form $V[1]$ for $V$ a vector bundle on $X$, is naturally equivalent to the 
category of formal schemes $Z$ with an identification $Z_{red} \simeq X$, 
which are locally equivalent on $X$ to $X\times \widehat{\mathbb{A}}^d$ where $d=rank(V)$.
\end{prop}

\textit{Proof.} To each $\F$ as in the Proposition, we associate the sheaf 
of infinitesimal cohomology $\Cinf(\F)$ on $X$. By assumptions on the cotangent complex of $\F$, 
$\Cinf(\F)$ is a sheaf of discrete complete filtered commutative algebras, augmented 
to $\OO_X$ and 
its associated graded object is isomorphic to $\OO_X[[t_1,\dots,t_d]]$. Therefore, 
$Z_\F:=\Spf(\Cinf(\F))$ is a formal scheme with a canonical identification 
$Z_{red}\simeq X$ which is locally equivalent to $X\times \widehat{\mathbb{A}}^d$.

The fact that 
the construction $\F \mapsto Z_\F$ induces the equivalence of the Proposition simply follows from the
fact that the Tate realization induces an equivalence between graded mixed complexes and
complete filtered complexes (see \cite[Prop. 1.3.1]{derfol2}). \hfill $\Box$ \\

By Proposition \ref{p1}, for any smooth infinitesimal derived foliation $\F$, its loop space 
$\Omega_X \F$ can be realized as a groupoid object in formal schemes acting on $X$
$$Z_{\Omega_X\F} \longrightarrow X\times X.$$

The construction $\F \mapsto Z_{\Omega_X\F}$ provides an functor from the
category of smooth infinitesimal derived foliations on $X$, and the category of 
formally smooth formal groupoids over $X$. We can produce a functor in the other direction, as follows.
Starting from a formally smooth formal groupoid $G \rightarrow X\times X$, we can consider its nerve
$G_* \to X$, which is a simplicial object in formal schemes. Taking functions on $G_*$
provides a cosimplicial complete filtered commutative algebra $\OO(G_*)$ whose associated graded
is the cosimplicial algebra of functions on the simplicial scheme $BV$, where 
$V=e^*(\Omega_{G/X}^1)$ is the vector bundle of invariant relative forms on $G$.
Using the equivalence between graded mixed complexes and complete filtered complexes of 
\cite[Prop. 1.3.1]{derfol2}, 
we find that $\OO(G_*)$ can be realized as a cosimplicial object in the category of
graded commutative algebras endowed with a compatible action on $G_0$. Passing then to $Spec$, we get 
a simplicial diagram of affine schemes endowed with a $\cH_\pi$-action, whose underlying 
simplicial diagram (obtained by forgetting the action) is the usual simplicial scheme $BV$. 
Taking geometric realization we get a $\cH_\pi$-action on the derived stack $\VV(V[1])$, 
which defines a smooth infinitesimal derived foliations $\F$ on $X$. 

We can then subsume the previous discussion in the following corollary.

\begin{cor}\label{cint}
The category of smooth infinitesimal derived foliations on $X$ is equivalent to the
category of formally smooth formal groupoids on $X$. 
\end{cor}

We note that the above corollary marks a major difference between 
derived foliations of Section \S \ref{derfolgradedloop} and infinitesimal derived foliations. Indeed,
any finite dimensional Lie algebra $L$ over a field $k$, defines a 
smooth derived foliation over $k$, by considering $Sym(L^\vee[1])$ endowed with 
the mixed structure coming from the Chevalley-Eilenberg differential. However, 
it is not true that any Lie algebra can be integrated to a smooth formal
scheme in general. On the contrary, the above corollary, when applied to $X=\Spec\, k$ with $k$ a field, 
states that smooth infinitesimal derived foliations over $k$ form a category equivalent to 
smooth formal groups over $k$. Any such smooth infinitesimal derived foliation
is of the form $\VV(V[-1])$ for $V$ a finite dimensional $k$-vector space $V$. Therefore, the data
of an infinitesimal derived foliation structure on $\VV(V[-1])$ seems to be related to that of a \emph{partition
Lie algebra structure} in the sense of \cite{bracanu}. The precise comparison is out of the scope of this note, 
but has been carried out in \cite{jfuarxiv}.

\subsection{Comparison between infinitesimal derived foliations and derived foliations: the "de Rhamification" functor }

In this \S \, we explain how to construct a derived foliation from any infinitesimal
derived foliation. This construction will be called the \emph{de Rhamification} functor, 
as it is related to the natural "forgetful" morphism from (Hodge-completed and
derived) infinitesimal cohomology to de Rham cohomology. \\

We start by $X=\Spec\, A$ a derived affine scheme of finite presentation over $k$ and
$\F \in \Fol_\pi(X)$ an infinitesimal derived foliation. The natural 
projection $\VV(\F) \to X$ possesses a zero section $s : X \to \F$ which is moreover naturally
$\cH_\pi$-equivariant. We can then consider the relative derived loop space of $s$
$\ccL^{gr}(X/\VV(\F))$, which fits in a Cartesian square of $\cH_\pi$-equivariant derived stacks
$$\xymatrix{
\ccL^{gr}(X/\VV(\F)) \ar[r] \ar[d] & \cL^{gr}(X/k) \ar[d] \\
X \ar[r] & \ccL^{gr}(\VV(\F)/k).
}$$

The group stack $\cH$ acts by its standard action on $\ccL^{gr}(X/\VV(\F))$, and this
action is itself $\cH_\pi$-equivariant, or in other words is combined into a 
natural action of $\cH \times \cH_\pi$. The derived stack $\ccL^{gr}(X/\VV(\F))$
is a priori not affine, but we can consider its affinization (see \S \ref{app:derivedaffine})
$$\ccL^{gr}(X/\VV(\F))^{aff} := \RR\Spec^{\Delta}(\OO(\cL^{gr}(X/\VV(\F))).$$
Because the group stack $\cH \times \cH_\pi$ is affine and 
$\OO(\cH \times \cH_\pi)$ is perfect as a $k$-dg-module, it is easy to check that 
$\cH \times \cH_\pi$ continues to act on $\ccL^{gr}(X/\VV(\F))^{aff}$ in a natural manner (because
the affinization preserves taking products with $\cH \times \cH_\pi$). Moreover, 
using the descent theorem \ref{tII-1} we have a canonical identification
$$\ccL^{gr}(X/\VV(\F))^{aff} \simeq \RR\Spec^{\Delta}(Sym_A(\LL_{\F}[1])) = \VV(\LL_{\F}[1]).$$
We have therefore constructed a natural action of $\cH \times \cH_\pi$ on 
$\VV(\LL_\F[1])$. By construction, the induced action of $B\Gm \times B\Gm$ 
is the diagonal action, making $\LL_\F[1]$ pure of weight $(1,1)$.

\begin{lem}\label{actiontrivial}
With the above notations, the action of $\Ga[-1]$ on $\VV(\LL_\F[1])$ is canonically trivial.
\end{lem}

\textit{Proof of the lemma.} This is formally deduced from the fact that 
the action of $\Ga[-1]$, or more generally of $\cH_\pi$, respects the action
of $\cH$. In particular, $\Ga[-1]$ acts on $\VV(\LL_\F[1])$ by respecting the grading. 
Using \cite{mon} this action is then determined by a linear action of $\Ga[-1]$ on 
the perfect complex $\LL_{\F}$, compatible with the gradings, 
and thus by a graded mixed structure on $\LL_\F$, which is pure of weight $1$. 
Such a graded mixed $A$-dg-module is represented by an object $E \in \medg_k$ with 
a fixed quasi-isomorphism $E^{(1)} \to \LL_\F$ and $E^{(i)}$ being acyclic for all $i\neq 1$. 
Such an object is equivalent $E'$ with $(E')^{(1)} = \LL_\F$ and 
$(E')^{(i)}=0$ for $i\neq 1$ (and thus $\epsilon = 0$). \hfill $\Box$ \\

We now form the morphism $K := B\cH_\pi \times_{B\Gm} B\cH \to B\cH_\pi \times B\cH$, the corresponding
quotient derived stack
$[\VV(\LL_\F[1])/K] \to B\cH$, and its affinization
$$F:=[\VV(\LL_\F[1])/K]^{aff} \to B\cH.$$
Lemma \ref{actiontrivial} implies that $F\times_{B\cH}B\Gm$ is naturally equivalent to $\VV(\LL_\F[1])$, and the natural
projection $\VV(\LL_\F[1]) \to X$ makes $F$ into a derived foliation on $X$ whose cotangent complex
is identified with $\LL_\F$.

The above construction is clearly functorial in $\F$ and $X$, and thus defines for 
any derived Artin stack $X$ a well defined $\s$-functor
$$(-)^{DR} : \Fol_\pi(X/k) \to \Fol(X/k),$$
sending an infinitesimal derived foliation $\F$ to its de Rhamification $\F^{DR}$. 

\begin{df}
The above $\s$-functor $(-)^{DR}$ is called the \emph{de Rhamification} construction. For
an infinitesimal derived foliation $\F \in \Fol_\pi(X/k)$, the
derived foliation $\F^{DR}$ is called the \emph{underlying derived foliation} of $\F$.
\end{df}

We note in particular that the cotangent complexes are preserved by de Rhamification
$$\LL_{\F} \simeq \LL_{\F^{DR}}.$$

The behaviour of the de Rhamification construction is rather subtle. For instance, already the question of
which derived foliations arise as the de Rhamification of infinitesimal derived foliations is a non-trivial question.
Derived foliations are not formally integrable, so our formal integrability result Corollary \ref{cint}
provides a first obstruction for a derived foliation to come from an infinitesimal derived foliation. 
This can be seen already at the level of classical smooth foliations, as formally integrable
foliations on a smooth varieties in characteristic $p>0$ must always be stable by the $p$-th power 
operation on the tangent 
sheaf (see \cite{MR927978}). The precise relation between infinitesimal derived foliations and
an hypothetical notion $p$-restricted derived foliations is also not totally obvious. By the results
of \cite{jfuarxiv} the question translate, for derived foliations over a perfect field, 
to the comparison between $p$-restricted dg-Lie algebras and partition Lie algebras of \cite{bracanu}.
As far as we know this comparison has not been clarified yet.

Nevertheless the existence of $\F \mapsto \F^{DR}$ is already helpful in the following sense. Starting 
with an infinitesimal derived foliation $\F$ we can exploit its formal integrability, and use 
$\F^{DR}$ when we want to study its de Rham cohomology, which is much better behaved than 
infinitesimal cohomology. For instance, we believe that the $\s$-category of crystals on 
$\F^{DR}$ is more interesting for the study of cohomological properties than 
the $\s$-category of crystals on $\F$ (defined naturally using $\cH_\pi$-equivairant quasi-coherent complexes 
on $\VV(\F)$).

\begin{appendix}

\chapter{}

\section{Model categories and $\s$-categories}\label{secapp:Modelcategoriesandscategories}

\subsection{Model categories}\label{subsec_modelcats}
The notion of model category, due to D. Quillen \cite{quiHA}, is a very convenient axiomatization of a 
homotopy theoretical context. Modern standard references include  \cite{hov, hirschhorn}.

\begin{df}
Let $\mathcal{C}$ be a category. A \emph{model structure on} $\mathcal{C}$ consists of the choice of three classes 
$\mathsf{W}(\mathcal{C})$ (weak equivalences), $\mathsf{Fib}(\mathcal{C})$ (fibrations), and $\mathsf{Cofib}(\mathcal{C})$ 
(cofibrations)  of morphisms in $\mathcal{C}$, subject to the following properties.
\begin{enumerate}
    \item If $(f,g)$ is a pair of composeable morphisms in $\mathcal{C}$ (i.e. $f\circ g$ exists in $\mathcal{C}$), if two 
    among the morphisms $f, g, f\circ g$ is in $\mathsf{W}(\mathcal{C})$, then the same is true for the remaining one.
    \item Each class $\mathsf{W}(\mathcal{C})$, $\mathsf{Fib}(\mathcal{C})$, and $\mathsf{Cofib}(\mathcal{C})$ is stable 
    under composition and retracts\footnote{See \cite[Def. A.2.1.1 (3)]{htt} for a precise statement of stability under 
    retracts} in $\mathcal{C}$.
    \item Let $$\xymatrix{X \ar[r]^-{f} \ar[d]_-{i} & Y \ar[d]^-{p} \\
    Z \ar[r]_-{g} & T}$$ be a commutative diagram in $\mathcal{C}$, with $i \in \mathsf{Cofib}(\mathcal{C}) $ and $p \in 
    \mathsf{Fib}(\mathcal{C})$. If either $i$ or $p$ also belongs to $\mathsf{W}(\mathcal{C})$, then there exists a 
    morphism $h: Z \to Y$ such that $p \circ h = g$ and $h\circ i = f$.
    \item Any morphism $f: X \to Y$ in $\mathcal{C}$ can be factored as $f=p\circ i$ and as $f= q \circ j$ where $i \in 
    \mathsf{Cofib}(\mathcal{C})$, $p \in \mathsf{Fib}(\mathcal{C}) \cap \mathsf{W}(\mathcal{C})$, $j \in 
    \mathsf{Cofib}(\mathcal{C}) \cap \mathsf{W}(\mathcal{C})$, $q \in \mathsf{Fib}(\mathcal{C}) $. Moreover, these 
    factorizations are both required to be functorial in $f$.
\end{enumerate}
A \emph{model category} consists of a complete and cocomplete category $\mathcal{C}$ endowed with a model structure.
\end{df}

\noindent \textbf{Model categories' vocabulary.} Let $(\mathcal{C}, \mathsf{W}(\mathcal{C}), \mathsf{Fib}(\mathcal{C}), 
\mathsf{Cofib}(\mathcal{C}) )$ be a model category. Morphisms in $\mathsf{Fib}(\mathcal{C}) \cap \mathsf{W}(\mathcal{C})$ 
are called \emph{trivial fibrations}, while morphisms in $\mathsf{Cofib}(\mathcal{C}) \cap \mathsf{W}(\mathcal{C})$ are 
called \emph{trivial cofibrations}. An object $P$ (respectively, $Q$) in $\mathcal{C}$ is said to be \emph{fibrant} 
(respectively, \emph{cofibrant}) if the canonical map from $P$ to the final object $*_{\mathcal{C}}$ (resp., from the 
intial object $\emptyset_{\mathcal{C}}$ to $Q$) is a fibration (resp., a cofibration). By (4), any object $x\in 
\mathcal{C}$ admits a diagram, which is moreover \emph{functorial} in $x$, $$\xymatrix{ Qx \ar[r]^-{p} & x \ar[r]^-{i} & 
Rx}$$ where $Qx$ is cofbrant, $p$ is a trivial fibration (by applying (4) to $\emptyset_{\mathcal{C}} \to x$), $Rx$ is 
fibrant and $i$ is a trivial cofibration (by applying (4) to $x \to *_{\mathcal{C}}$). The \emph{homotopy category} 
$\mathsf{ho}(\mathcal{C})$ of the model category $\mathcal{C}$ is the Gabriel-Zisman localization 
$\mathcal{C}[\mathsf{W}(\mathcal{C})^{-1}]$ of $\mathcal{C}$ at weak equivalences. It admits also an equivalent 
description as the category whose objects are those objects in $\mathcal{C}$ that are both fibrant and cofibrant, with 
morphisms sets given by homotopy classes of morphisms between such objects (\cite[\S 1.2]{hov}).  

If $(\mathcal{C}, \mathsf{W}(\mathcal{C}), \mathsf{Fib}(\mathcal{C}), \mathsf{Cofib}(\mathcal{C}) )$ and $(\mathcal{D}, 
\mathsf{W}(\mathcal{D}), \mathsf{Fib}(\mathcal{D}), \mathsf{Cofib}(\mathcal{D}) )$ are model categories, an adjoint pair 
$(F,G)$, with left adjoint $F:\mathcal{C} \to \mathcal{D}$ is called a \emph{Quillen adjunction} or \emph{Quillen pair} 
(and $F$, respectively $G$, is called a \emph{left}, resp. \emph{right}, \emph{Quillen functor}) if $F$ preserves 
cofibrations and trivial cofibrations or if, equivalently, $G$ preserves fibrations and trivial fibrations. A Quillen 
adjunction induces an adjunction at the level of the corresponding homotopy categories (\cite[\S 1.3]{hov}), and it is 
said to be a Quillen equivalence if the induced adjoints are (mutually inverse) equivalences of categories.

\subsection{$\s$-categories} There are several, a priori different, definitions of $\s$-categories: $\s$-categories as 
simplicial (i.e. categories enriched over the category of simplicial sets) categories (\cite[Definition 1.1.4.1]{htt}), 
$\s$-categories as quasi-categories (\cite[Definition 1.1.2.4]{htt}), $\s$-categories as Segal categories 
(\cite{Simpson_2011}), $\s$-categories as complete Segal spaces (\cite{Rezk, JoyTie}), $\s$-categories as Segal spaces 
(\cite{MosNui}), $\s$-categories as topological (i.e. categories enriched in the category of compactly generated weakly 
Hausdorff topological spaces) (\cite[Definition 1.1.1.6]{htt}).
The common idea behind these definitions is that a $\s$-category should be a category (maybe weakly) enriched in 
topological spaces, or in any equivalent homotopical model for topological spaces (e.g. simplicial sets). For a precise  
comparison between these different models, the interested reader may consult \cite{htt}, \cite{Bergner2010}, 
and \cite{jberg}.\\

For example, any usual 1-category gives rise to an $\s$-category, to which it is usually identified. If for 
$\s$-categories we take the model of quasi-categories, then the $\s$-category associated to a $1$-category $\mathcal{C}$ 
is its \emph{nerve} $\mathsf{N}(\mathcal{C})$, i.e. the simplicial set such that $$\mathsf{N}(\mathcal{C})_n := 
\mathsf{Fun}([n], \mathcal{C}).$$ \\
All the standard concepts in usual (1-)category theory have appropriate extensions to $\s$-categories. For example, we 
have a notion of $\infty$\emph{-functor} $\mathbf{F}: \mathbf{C} \to \mathbf{D}$: in the model of $\infty$-categories 
given by quasi-categories (respectively, simplicial categories), an $\s$-functor is just a morphism between the 
corresponding simplicial sets (respectively, a simplicially enriched functor between the corresponding simplicially 
enriched categories). Given two $\s$-categories $\mathbf{C}$ and $\mathbf{D}$, there is, actually, an $\s$-category 
$\mathbf{Fun}_{\s}(\mathbf{C}, \mathbf{D})$ of $\infty$-functors from $\mathbf{C}$ to $\mathbf{D}$ (see \cite[1.2.7.2, 
1.2.7.3]{htt} for quasi-categories). There is a notion of \emph{adjoint $\s$-functors} between $\s$-categories (see 
\cite[Definition 5.2.2.1]{htt} for quasi-categories), notions of \emph{limit} and \emph{colimit} for diagrams in an 
$\s$-category (see \cite[\S 1.2.13]{htt} for these concepts inside a quasi-category), Kan extensions (see \cite[\S 
4.3]{htt} for quasi-categories), etc. For any pair $(x,y)$ of objects inside an $\s$-category $\mathbf{C}$, there is a 
\emph{mapping space} $Map_{\mathbf{C}}(x,y)$ between them, that can be identified with a simplicial set, and is defined up 
to homotopy (i.e. it is well defined in $ho(\mathsf{sSets})$). The \emph{homotopy category} of an $\s$-category 
$\mathbf{C}$, is the $1$-category $\mathsf{ho}(\mathbf{C})$ with the same objects as $\mathbf{C}$ and morphism from $x$ to 
$y$ given by $\pi_0(Map_{\mathbf{C}}(x,y))$. There is an obvious $\s$-functor $ \mathbf{C} \to \mathsf{ho}(\mathbf{C})$, 
and a morphism in $\mathbf{C}$ will be called an equivalence if it becomes an isomorphism in  $\mathsf{ho}(\mathbf{C})$.\\ 
Moreover, (small) $\s$-categories and $\s$-functors between them, organize themselves\footnote{With appropriate 
set-theoretic or universes's caveats.} again into a $\s$-category $\mathbf{Cat}_{\infty}$ (see \cite[Chapter 3]{htt} for 
quasi-categories). More is true: $\mathbf{Cat}_{\infty}$ is moreover a $(\infty, 2)$-category, i.e. an $\s$-categorical 
version of a bicategory, but we will mainly view it simply as an $\infty$-category. The nerve construction from 
1-categories to $\infty$-categories, recalled above, provides a fully faithful functor. \\

For an $\infty$-category $\mathbf{C}$ and a set $\mathbf{W}$ of morphisms, one can construct another $\s$-category denoted 
by either $L_{\mathbf{W}}(\mathbf{C})$ or $\mathbf{C}[\mathbf{W}^{-1}]$ where, informally speaking, the morphisms in 
$\mathbf{W}$ has been inverted. The construction also provides a canonical $\s$-functor $L: \mathbf{C} \to 
L_{\mathbf{W}}(\mathbf{C})$ such that, for any $\s$-category $\mathbf{D}$, the induced $\infty$-functor 
$$\mathbf{Fun}_{\infty}(L_{\mathbf{W}}(\mathbf{C}),\mathbf{D}) \longrightarrow 
\mathbf{Fun}_{\infty}(\mathbf{C},\mathbf{D})$$ is fully faithful with essential image consisting of those $\s$-functors 
$\mathbf{C}\to \mathbf{D}$ sending morphisms in $\mathbf{W}$ into equivalences in $\mathbf{D}$.  The $\s$-category 
$L_{\mathbf{W}}(\mathbf{C})$ is called the \emph{localization of $\mathbf{C}$ at (or along) $\mathbf{W}$}. 

\begin{df}\label{def_underlyinginfinitycattomodelcat}
Given a model category $(C, \mathsf{W}(\mathcal{C}), \mathsf{Fib}(\mathcal{C}), \mathsf{Cofib}(\mathcal{C}) )$, its 
\emph{underlying or associated $\s$-category} will be the localization 
$L_{\mathsf{W}(\mathcal{C})}(\mathsf{N}(\mathcal{C}^{\mathrm{cof}}))$, where $\mathcal{C}^{\mathrm{cof}}$ denotes the full 
subcategory of $\mathcal{C}$ consisting of cofibrant objects.
\end{df}

\begin{rmk} \emph{In the main text, if $\mathcal{C}$ is a model category, we have simply denoted by $\mathbf{C}$ its 
underlying $\infty$-category. }
\end{rmk}

A Quillen adjunction between model categories (\ref{subsec_modelcats}) induces an adjunction between $\s$-functors on the 
underlying $\s$-categories (\cite[\S 5.2.4]{htt}), which is an equivalence of $\s$-categories if the Quillen pair is a 
Quillen equivalence.

\section{Derived schemes and stacks: basic facts and notations}\label{secapp:Derivedschemesandstacks}

In this \S $\,$ we will mainly fix the notations used in the main text, we do \emph{not} mean to give a short introduction 
to derived algebraic geometry, for which we refer to \cite{toenseattle, toenems, pv}.\\

Let $k$ be a base discrete commutative ring (e.g. $k=\mathbb{Z}$), we denote by $\mathsf{sCAlg}_k$ the model 
category\footnote{This is, moreover, a cofibrantly generated simplicial model structure; see \cite{hirschhorn} for the 
meaning of these words.} of simplicial commutative $k$-algebras (where both weak equivalences and fibrations are detected 
by the forgetful functor $\mathsf{sCAlg}_k \to \mathsf{sSets}$), and by $\mathbf{sCAlg}_k$ the underlying $\s$-category 
(see \ref{subsec_modelcats}). Note that, for $A\in \mathsf{sCAlg}_k$,  $\pi_*(A)=\oplus_{i\geq 0}\pi_i (A)$ has a natural 
induced structure of a commutative graded $k$-algebra. In particular, $\pi_0(A)$ has an induced structure of usual 
commutative $k$-algebra, and each $\pi_i(A)$ is a $\pi_0(A)$-module. We say that $A$ is \emph{Noetherian} if $\pi_0 (A)$ 
is a Noetherian ring, and each $\pi_i (A)$ is a finitely generated (discrete) $\pi_0(A)$-module. We will say that $A$ is 
\emph{discrete} if $\pi_i (A)=0$ for any $i>0$ (i.e. $A$ is weakly equivalent to the constant simplicial $k$-algebra 
$\pi_0 (A)$). \\ 

For a simplicial commutative $k$-algebra $A$, the associated \emph{affine derived scheme} is the 
$\s$-functor
$$\mathbf{Spec}\, A := \mathsf{Map}_{\mathbf{sCAlg}_k} (A, -): \mathbf{sCAlg}_k \longrightarrow 
\mathbf{sSets}=\mathbf{Top}$$
We let $\mathbf{dAff}_k$, the $\s$-category of \emph{affine derived schemes} over $k$, be the full $\s$-subcategory of 
$\mathbf{Fun}_{\s}(\mathbf{sCAlg}_k, \mathbf{Top})$ essential image of the $\s$-functor $\mathbf{Spec}$. The $\s$-category 
$\mathbf{dAff}_k$ has arbitrary pullbacks: for $\xymatrix{B & A \ar[l] \ar[r] & C}$ a diagram in $\mathbf{sCAlg}_k$, the 
corresponding pullback in $\mathbf{dAff}_k$ is given by $\mathbf{Spec}(B\otimes_A^{\mathrm{L}}C)$, where 
$B\otimes_A^{\mathrm{L}}C$ is the derived tensor product of $B$ and $C$ in $A/\mathsf{sCAlg}_k$, obtained by replacing 
either $B$ or $C$ by a cofibrant replacement in the undercategory $A/\mathsf{sCAlg}_k$, before taking the usual (strict) 
tensor product.\\ 

The inclusion of usual discrete commutative $k$-algebras in $\mathsf{sCAlg}_k$, as constant simplicial 
object, induces a fully faithful functor $j:\mathbf{Aff}_k \to \mathbf{dAff}_k$ having as right adjoint $\mathbf{t}_0: 
\mathbf{Spec}\,A \mapsto \mathrm{Spec}(\pi_0 (A))$, called the \emph{truncation functor}. An important property of $j$ is 
that it does \emph{not} preserve general pullbacks: if $\xymatrix{S & R \ar[l] \ar[r] & T}$ is a diagram in 
$\mathsf{CAlg}_k$, then $S\otimes_R^{\mathrm{L}}T \neq S\otimes_R T$ unless $S$ or $T$ are flat over $R$.\\

If, as often in the main text, the fixed base ring $k$ is a commutative $\mathbb{Q}$-algebra, then we may replace 
$\mathsf{sCAlg}_k$ with the model category $\mathsf{cdga}_k$ of commutative differential non-positively graded (aka 
\emph{connective}) $k$-algebras, since we have an equivalence of the underlying $\s$-categories $\mathbf{sCAlg}_k \simeq 
\mathbf{cdga}_k$. If, under this equivalence, $A\in \mathbf{sCAlg}_k$ corresponds to $A'\in \mathbf{cdga}_k$, we have 
$\pi_{i}(A) \simeq H^{-i}(A')$ for any $i\geq 0$. For this reason, we will sometimes abuse notations and write $\pi_i(B)$ 
for $H^{-i}(B)$ even when $B \in \mathbf{cdga}_k$.\\

For $A\in \mathsf{sCAlg}_k$, its normalization $\mathsf{N}(A)$ is a commutative differential graded $k$-algebras, and we 
denote by $\mathsf{Mod}_A$ the model category of \emph{not necessarily bounded} $\mathsf{N}(A)$-dg modules 
(see \cite[Ch. 2.2]{hagII} for more details). As usual, $\mathbf{Mod}_A$ will be the underlying $\s$-category of 
$\mathsf{Mod}_A$. The same notation is used for $A \in \mathsf{cdga}_k$. For $M \in \mathsf{Mod}_A$ and 
$i \in \mathbb{Z}$, we will write $\pi_i(M)$ or $H^{-i}(M)$ interchangeably.\\

If $\mathbf{C}$ denotes either $\mathbf{sCAlg}_k$ or $\mathbf{cdga}_k$, a morphism $A \to B$ in $\mathbf{C}$ is called 
\emph{locally of finite presentation} if the functor $\mathbf{Map}_{A/\mathbf{C}}(B,-)$ commutes with filtered colimits 
(i.e. $B$ is a compact object inside $A/\mathbf{C}$). Since the forgetful $\s$-functor $A/\mathbf{C} \to \mathbf{Mod}_A$ 
preserves and reflects filtered colimits, this condition is equivalent to $\mathbf{Map}_{\mathbf{Mod}_A}(B,-): 
\mathbf{Mod}_A \to \mathbf{Top}$ preserving filtered colimits. A morphism $A \to B$ in $\mathbf{C}$ is called 
\emph{locally almost of finite presentation} (sometimes abbreviated to \emph{lafp}) $\mathbf{Map}_{\mathbf{Mod}_A}(B,-)$ 
preserves filtered colimits when restricted to each full $\s$-subcategory $\mathbf{Mod}^{\leq n}_A$ of \emph{truncated} 
$A$-modules (i.e. $A$-modules $M$ such that $\pi_i(M)=0$ for $i> n$. For example, if $A$ is Noetherian, then $A \to B$ is 
lafp iff $B$ is Noetherian, and $\pi_{0}(B)$ is a finitely presented $\pi_0 (A)$-algebra in the usual sense of commutative 
algebra.\\

One can consider different $\s$-topologies (\cite[6.2.2]{htt}, \cite[3.1]{hagI}) on the $\s$-category $\mathbf{dAff}_k$. 
For example, the \emph{\'etale} $\s$-topology is generated by the $\s$-pre-topology defined as follows: for 
$X=\mathbf{Spec}\, A \in \mathbf{dAff}_k$ where $A \in \mathsf{sCAlg}_k$ or $A \in \mathsf{cdga}_k$, \'etale covering 
families of $X$ are given by arbitrary families $\{A \to A_i\}_{i\in \mathrm{I}}$ such that the induced family ${\pi_0(A) 
\to \pi_0 (A_i)}_{i\in \mathrm{I}}$ is a usual \'etale covering family of $\mathbf{t}_0(\mathbf{Spec}\, A)= 
\mathrm{Spec}(\pi_0 (A))$, and moreover, the canonical maps $$\pi_i(A)\otimes_{\pi_0(A)}\pi_0 (A_i) \to \pi_i (A_i)$$ are 
isomorphisms\footnote{Recall our convention $\pi_i(A):= \mathrm{H}^{-i}(A)$ for $A \in \mathsf{cdga}_k$, and $i\geq 0$.} 
for all $i \geq 0$ (and we say, then, that each $A \to A_i$ is \emph{strongly \'etale}). Once an $\s$-topology $\tau$ on $\mathbf{dAff}_k$ is fixed, we define the $\s$-category of \emph{derived $\tau$-stacks} 
$$\mathbf{dSt}^{\tau}_k := \mathbf{Sh}((\mathbf{dAff}_k, \tau), \mathbf{Top})$$ the category of $\mathbf{Top}$-valued 
$\s$-sheaves on the $\s$-site $(\mathbf{dAff}_k, \tau)$ (\cite[Def. 6.2.2.6]{htt}, \cite[Def. 3.4.9]{hagI}). Therefore, a 
derived $\tau$-stack is an $\s$-functor $F: \mathbf{dAff}^{\mathrm{op}}_k \to \mathbf{Top}$ (i.e. a $\s$-preseheaf on 
$\mathbf{dAff}_k$) satisfying $\tau$-hyperdescent, i.e. the canonical morphism $F(X) \to F(U_{\bullet})$ is an equivalence 
in $\mathbf{Top}$, for any $\tau$-hypercover $U_{\bullet} \to X$.\\

In order to select geometrically meaningful and interesting derived stacks, one imposes recursively defined 
\emph{geometricity} conditions on atlases and representable morphisms (see \cite[Ch. 1.3]{hagII}, \cite[3.3]{toenems} and 
\cite[\S 4]{pv}). In particular, a derived \emph{Deligne-Mumford} stack is a derived \'etale stack admitting an \'etale 
atlas (by derived affine schemes), while a derived \emph{Artin} (or \emph{algebraic}) stack is a derived \'etale stack 
admitting a smooth atlas.\\

Any connective $k$-linear cdga $A$ possesses a cotangent complex $\LL_{A/k}$ relative to $k$. 
By definition $\LL_{A/k}$ is an $A$-dg-module such that derivations $A \to M$
are equivalent to $A$-dg-modules maps $\LL_{A/k} \to M$ (see \cite[\S 1.4.1]{hagII}). More precisely, 
for $M \in \dg_A$ an $A$-dg-module we can form the trivial square
zero extension $A \oplus M$, which is a new cdga augmented to $A$. We then have
functorial equivalences of mapping spaces
$$\Map_{\dg_A}(\LL_{A/k},M) \simeq \Map_{\cdga_k/A}(A,A\oplus M),$$
and the right hand side is, by definition, the space of derivations on $A$ 
with coefficients in $M$. In particular, the 
identity of $\LL_{A/k}$ defines a universal derivation $dR : A \to \LL_{A/k}$. 
The cotangent complex of non-affine derived schemes and more
generally of derived Artin stacks is discussed below in \S \ref{secapp:CotangentcomplexesofderivedArtinstacks}. \\

\section{Formal completion}\label{secapp:Formalcompletion}

Let $A$ be a connective cdga over a field $k$ of characteristic zero endowed with 
au augmentation $x : A \to k$. We can define a pro-Artinian cdga $\underline{\hat{A}}_x$
(or simply $\underline{\hat{A}}_x$), called the \emph{formal completion of $A$ at $x$}, as follows.
We consider $A/\dgart^*_k$ the full sub-$\s$-category of the $\s$-category of 
augmented connective cdga under $A$ consisting of all morphisms $A \to A'$ where
$A'$ is a local Artinian cdga. Recall that Artinian here means that 
$H^0(A')$ is a local Artinian $k$-algebra and moreover that $H^*(A')$ is finite dimensional
(so all $H^i(A')$ are finite dimensional and vanish for $i<<0$).
We note that $A/\dgart^*_k$ is a filtered $\s$-category as it admits finite limits. The following
definition therefore makes sense.

\begin{df}\label{dA-3-1}
With the above notations, \emph{the formal completion of $A$ at $x$} is the pro-object
in $\dgart_k^*$ defined by
$$\underline{\hat{A}}_x:="\lim"_{A' \in (A/\dgart^*_k)} A'.$$
\end{df}

The pro-object $\underline{\hat{A}}$ has a realization $\hat{A}=\lim_{A' \in (A/\dgart^*_k)} A'$
which is a cdga. We will see below that $\hat{A}$ is connective when 
$A$ is almost of finite presentation. This does not follow
formally from the definition (as connective cdga's are not stable by limits inside the
bigger $\s$-category of all cdga's). This will be a consequence of a more explicit 
description of the pro-object $\underline{\hat{A}}$. \\

Assume now that $A$ is a connective cdga wich is moreover almost of finite presentation.
We consider the ideal $I \subset H^0(A)$ of the induced augmentation $H^0(A) \to k$, and
pick $(f_1,\dots,f_p)$ a set of generators for $I$. For each integer $n>0$ we define
an object $A \to A_n$ of $A/\dgart_k^*$ as follows. Start with 
the Koszul cdga $K(A,f^n_*)$, which is the free $A$-linear cdga generated
by elements $h_1,\dots h_p$ in degree $-1$ with $dh_i=f^n_i$. We let $A_n:=\tau_{\leq n}
K(A,f^n_*)$ the n-truncation of $K(A,f^n_*)$. For $m>n$, we have a canonical 
morphism $p_{m,n} : A_m \to A_n$ which sends $h_i$ to $f_i^{m-n}.h_i$ making the $A_n$ into a projective
system in $A/\cdga$. Moreover, $H^0(A_n)\simeq H^0(A)/(f_1^n,\dots,f_p^n)$ is a local 
Artinian algebra, and each $H^i(A_n)$ is of finite type over $H^0(A_n)$. Indeed, 
we can see this by induction on the number of elements $f_i$, using that 
$K(A,f_*)\simeq K(K(A,f_{*< k}),f_k)$,
and the
long exact sequence 
$$\xymatrix{
H^{i-1}(A) \ar[r] & H^i(A) \ar[r]^-{f} & H^i(A) \ar[r] & H^i(K(A,f)) \ar[r] & 
H^{i+1}(A) \ar[r]^-{f} & \dots}$$
for the Koszul algebra $K(A,f)$ of a single element $f \in H^0(A).$ Finally, 
$A_n$ are truncated by construction.
The family of objects $A_n$ thus form a pro-object $"\lim_n" A_n$ in $A/\dgart_k^*$. There is an obvious
morphism of pro-object $\underline{\hat{A}} \to "\lim_{n}"A_n$, as the diagram given by the
canonical morphisms $A \to A_n$.

\begin{prop}\label{pA-3-1}
With the notations above, the natural morphism $\underline{\hat{A}} \to "\lim_{n}"A_n$ is an equivalence
of pro-objects in $A/\dgart_k^*$.
\end{prop}

\textit{Proof.} Let $\N$ be the category of natural numbers with its natural order. The family of
objects $A_n$ defines an $\s$-functor $\phi : \N^{op} \to A/\dgart_k^*$. The Proposition will follow
from the fact that $\phi$ is cofinal. 
As a first observation, for any object $A \to A'$ in $A/\dgart_k^*$, the image of 
all $f_i$ in $H^0(A')$ are nilpotent. Therefore, we can chose $n$ such that 
for each $i$, the image of $f_i^n$ is homotopic to zero in $A'$. By the universal property of
the $A$-cdga $A_n=K(A,f^n_*)$, there exists a factorization
$A \to A_n \to A''$.

Let $A' \in A/\dgart_k^*$, and let $\phi/A'$ be the slice $\s$-category of $\phi$. It is enough
to show that $\phi/A'$ is cofiltered. For this, let $B_* : J \to \phi/A'$ be an $\s$-functor
$j \mapsto B_j$
with $J$ a finite $\s$-category. We can for the limit along $J$ and get this way an dg-Artinian 
local algebra $A_J=\lim_{j\in J} B_j$ and a factorization $A \to A_J \to A'$. By our previous observation, 
there exists a factorization $A \to A_p \to A_J$. 

For any $j \in J$, $B_j$ is of the form $A_{n_j}$ for some integer $n_j$. The morphism
$A_p \to A_{n_j}$ is possibly not of the standard form, but it can be made so by increasing
$p$ thanks to the following lemma.

\begin{lem}\label{lpA-3-1}
For any integer $p\geq n$, and any morphism of $A$-cdga's
$u : A_p \to A_n$, there exists $q\geq p$ such that the composition 
$$\xymatrix{A_q \ar[r]^-{p_{q,p}} & A_p \ar[r]^-{u} & A_n}$$
is homotopic to the standard morphism $p_{q,n}$.
\end{lem}

\textit{Proof of the lemma.} This follows easily from the universal property of Koszul 
algebras. Indeed, for any $A$-cdga $B$, the space of $A$-cdga morphisms
$Map_{A-\cdga}(A_p,B)$ is $\Omega_{0,f^p}B$, the space of homotopies between $0$ and
the image of $f^p$ in $B$. It is therefore either empty or non-canonically 
equivalent to the loop space $\Omega_0B$, and a choice of a homotopy $0 \sim f^p$ induces
an equivalence $\Omega_{0,f^p}B \simeq \Omega_0B$. As a result,
for any $p\geq n$, the set of homotopy classes of $A$-cdga morphisms 
$\pi_0(Map_{A-\cdga}(A_p,A_n))$ is isomorphic to $H^{-1}(A_n)$, and moreover this
isomorphism is canonical as the universal homotopy $0 \sim f^n$ in $A_n$ induces
a homotopy $0 \sim f^p$ by multiplication by $f^{p-n}$. 

Under the bijection $\pi_0(Map_{A-\cdga}(A_p,A_n))\simeq H^{-1}(A_n)$, the
canonical morphism $p_{p,n}$ corresponds to $0$ in $H^{-1}(A_n)$. Therefore, the lemma simply
follows by observing that for $q>>p$ (and precisely $q\geq 2p$), 
the morphism $p_{q,p}$ induces the zero map $H^{-1}(A_q) \to H^{-1}(A_p)$. \hfill $\Box$ \\

Coming back to the proof of the Proposition. The lemma tells us that if $p$ is chosen 
big enough, then all the morphisms $A_p \to B_j$ are the standard morphisms. This means that 
the $A_p$ defines indeed an object in $\phi/A'$ that admits a morphism to 
the diagram $B_*$. \hfill $\Box$ \\

An important consequence of Proposition \ref{pA-3-1} is the description of 
the pro-$\underline{\hat{A}}$-modules $\wedge^p\LL_{\underline{\hat{A}}}:="\lim_{A\to A'}"\wedge^p \LL_{A'}$
of $p$-forms. 

\begin{cor}\label{cpA-3-1}
The natural morphism
$$"\lim_{A\to A'}" \LL_{A}\otimes_{A}A' \longrightarrow \LL_{\underline{\hat{A}}}$$
is an equivalence of pro-$\underline{\hat{A}}$-modules.
\end{cor}

\textit{Proof of Corollary \ref{cpA-3-1}.} For all $A \to A'$ we have an exact triangle of 
pro-$\underline{\hat{A}}$-modules
$$\LL_{A}\otimes_{A}A' \to \LL_{A'} \to \LL_{A'/A}.$$
Therefore, in order to prove the statement we need to prove that the pro-object
$"\lim_{A \to A'}" \LL_{A'/A}$ is zero. To see this, it is enough by Proposition 
\ref{pA-3-1} to see that $"\lim_{n}" \LL_{A_n/A}$ is the zero pro-$"\lim_nA_n"$-module. For each 
$n$, $\LL_{A_n/A}\simeq A_n[1]^p$, a basis being given by the de Rham differential $dR(h_i)$ (relative to $A$)
of $h_i$, of the canonical elements $h_i$ of degree $-1$ in $A_n$. Moreover, for
$p\geq 2n$, the canonical morphism $A_p \to A_n$ sends $dR(h_i)$ to $f^{p-n}.dR(h_i)$, and thus
is homotopic to the zero map because $f^{p-n}$ is homotopic to zero for $p\geq 2n$. 
This implies that the pro-$"\lim_nA_n"$-module $"\lim_{A \to A'}" \LL_{A'/A}$ is zero
as wanted. \hfill $\Box$ \\

By taking $p$-th wedge powers we find the following consequence.

\begin{cor}\label{cpA-3-1-2}
For any $p\geq 0$, the natural morphism
$$"\lim_{A\to A'}" (\wedge^p_A\LL_{A})\otimes_{A}A' \longrightarrow 
\wedge_{\underline{\hat{A}}}^p\LL_{\underline{\hat{A}}}\simeq "\lim_{A \to A'}" \wedge^p_{A'}\LL_{A'}$$
is an equivalence of pro-$\underline{\hat{A}}$-modules.
\end{cor}

Another important consequence is the following flatness result, stating the 
cohomology groups of the formal completion are the formal completion of the cohomology groups.

\begin{cor}\label{cpA-3-1-3}
Let $x : A\to k$ be an augmented connective cdga wich is almost of 
finite presentation and $\hat{A}=\lim_{A \to A'}A'$ be the realization of its formal completion at $x$.
Then, we have
$$H^0(\hat{A})\simeq \widehat{H^0(A)} \qquad H^i(\hat{A}) \simeq \widehat{H^i(A)} \simeq
H^i(A)\otimes_{H^0(A)}\widehat{H^0(A)}.$$
In particular, the morphism $A \to \hat{A}$ is a flat morphism of connected cdga's, which 
is moreover faithfully flat locally at $x$.
\end{cor}

\textit{Proof of Corollary \ref{cpA-3-1-3}.} As $A$ is almost of finite presentation over $k$, 
we can find a strict quasi-free model as a cdga of the form 
$$\xymatrix{ \dots \ar[r] & A^{-i} \ar[r] & A^{-i+1} \ar[r] & \dots \ar[r] & A^0}$$
where $A^0\simeq k[x_1,\dots,x_d]$ as a $k$-algebra, and each $A^i$ is a free 
$A^0$-module of finite rank. We can moreover assume that the augmentation $A \to k$ precomposed
with $A^0 \to A$ is given by the augmentation ideal $I=(x_1,\dots,x_q)$. Therefore, 
the $f_i$ are here represented by the $x_i$.

By definition of the Koszul algebra $A_n=K(A,f_*^{n})$ is thus given explicitly by 
$$K(A,f_*^{n})\simeq A\otimes_{A^{0}}K(A^0,x_*^{n}).$$
However, the $x_i^n$ form a regular sequence and thus $K(A^0,x_*^{n})$ is naturally 
quasi-isomorphic to the $k$-algebra $k[x_1,\dots,x_q]/(x_1^n,\dots,x_q^n)$. 
As $A$ is cofibrant as an $A^0$-module, we thus have 
$$K(A,f_*^{n}):=A\otimes_{A^{0}}K(A^0,x_*^{n}) \simeq A\otimes_{A^0}(A^0/(x_1^n,\dots,x_q^n)).$$
Therefore, by the Proposition \ref{pA-3-1}, the realization $\hat{A}$
is simply given by 
$$\hat{A}\simeq A\otimes_{A^0}\hat{A^0},$$
where $\hat{A^0}$ is the formal completion of $A^0$ at the ideal $(x_1,\dots,x_0)$. As formal
completions of noetherian rings are flat, we deduce the Corollary. \hfill $\Box$ \\

\section{$\OO_X$-modules, quasi-coherent and perfect modules}
\label{secapp:OXmodulesquasi-coherentandperfectmodules}

The $\s$-topos $\dSt_k$, of derived stacks over $k$, admits a tautological sheaf of cdga's, 
namely the left Kan extension of the tautological $\s$-functor sending $\Spec\, A$ to $A$.
It is denoted by $\OO_k$, and its restriction to $\dSt_k/X$ will be denoted by $\OO_X$. 
Associated to this sheaf of cdga's we have $\s$-categories of sheaves of modules, 
simply denoted by $\OO_k-\dg$ and $\OO_X-\dg$. In concrete terms, for $X \in\dSt_k$, and object $E$
in $\OO_X-\dg$ consists of the following data:

\begin{itemize}
    \item for all $u : \Spec\, A \to X$ an $A$-dg-module $E(u)$,
    \item for all $\xymatrix{\Spec\, A \ar[r]^-{f} & \Spec\, B \ar[r]^-{u} &  X}$
    a morphism of $B$-dg-modules $\phi_f : E(u\circ f) \to E(u)$,
\end{itemize}

together with the usual coherence conditions. In different terms, the $\s$-category 
$\OO_X-\dg$ can also be described as a lax limit
$$\OO_X-\dg \simeq lim_{\Spec\, A \to X}^{lax}A-\dg.$$

An object $E \in \OO_X-\dg$ will be called \emph{quasi-coherent} if all the transition 
morphisms $\phi_f : E(u\circ f) \to E(u)$ induces equivalences $A\otimes_B E(u\circ f) \simeq E(u)$.
The full sub-$\s$-category of quasi-coherent $\OO_X$-dg-modules (also called 
quasi-coherent complexes) will be denoted by 
$$j : \QCoh(X) \subset \OO_X-\dg.$$
This is a full embedding of presentable $\s$-catgeories, and obviously the inclusion $\s$-functor
commutes with arbitrary colimits. Therefore, the inclusion $j$ admits a right adjoint 
$$(-)^{qcoh} : \OO_X-\dg \longrightarrow \QCoh(X),$$
called the \emph{quasi-coherator} construction. 

We remind also two more full sub-$\s$-categories $\Perf(X) \subset \APerf(X) \subset \QCoh(X)$.
As a start let $A$ be a connective cdga over $k$. We recall that 
an $A$-dg-module $E$ is perfect if it is a compact object of $A-\dg$. It is
equivalent to as for $E$ being dualizable with respect to the natural
symmetric monoidal structure on $A-\dg$, and also equivalent to the fact that 
$E$ is a retract of a finite cell $A$-dg-module (see \cite[Prop. 2.2]{tova}). Similarly, 
an $A$-dg-module $E$ is \emph{almost perfect} if it can be approximated
by perfect objects: for all $n$, there exists a perfect $A$-dg-module $E_n$ and
a morphism $E_n \to E$ such that the induced morphism $H^i(E_n) \to H^i(E)$ is
an isomorphism for all $i>n$. As $A$ is connective, perfect $A$-dg-modules are bounded
above, and thus so are almost perfect $A$-dg-modules. Moreover, if $A$ is 
almost of finite presentation (i.e. $H^0(A)$ is of finite type over $k$ and
$H^i(A)$ is a finite type $H^0(A)$-module for all $i$) then 
$E \in A-\dg$ is almost perfect if and only if it satisfies the following two conditions:
\begin{itemize}
    \item $E$ is bounded above,
    \item $H^i(E)$ is a $H^0(A)$-module of finite type for all $i$.
\end{itemize}

To finish this section, remind that $\OO_X-\dg$ comes equipped with a canonical 
symmetric monoidal structure $\otimes_{\OO_X}$, determined for any $u : \Spec\, A \to X$ by 
the rule $(E\otimes_{\OO_X}F)(u):=E(u)\otimes_A F(u)$ plus sheafification. 
This symmetric monoidal structure preserves the full sub-$\s$-categories
$\Perf(X) \subset \APerf(X) \subset \QCoh(X) \subset \OO_X-\dg,$
and thus induces symmetric monoidal structures on all of them. Note also that $\Perf(X)$
consists of all dualizable objects in $\QCoh(X)$, and thus any perfect complex
is also dualizable as an $\OO_X$-dg-module. \\

We finish with an important statements concerning quasi-coherent sheaves on formal derived completions
of derived affine schemes. Let $X=\Spec\, A$ be a derived affine scheme almost of finite presentation
over a field $k$ (so in particular $H^0(A)$ is noetherain and each $Hi(A)$ is a $H^0(A)$-module
of finite type). Let $x \in X(k)$ be a global point and $\hat{\underline{A}}_x$ be the formal completion
of $A$ at the point $x$. It is a pro-object in $\cdga_k$ whose limit is
denoted by $\hat{A}_x$. 

\begin{prop}\label{pA-4-1}
The natural $\s$-functor
$$\phi : \APerf(\hat{A}_x) \to \Perf(\hat{\underline{A}_x}) = \lim_{A \to A'}\APerf(A'),$$
sending $E$ to the pro-module $"\lim_{A\to A'}" (A'\otimes_{\hat{A}_x}E)$, is an equivalence of
$\s$-categories.
\end{prop}

\textit{Proof.} From proposition \ref{pA-3-1} $\hat{\underline{A}}_x$ can be represented
as pro-object of the form 
$$\xymatrix{\hat{A}_x \ar[r] & 
\dots \ar[r] & A_n \ar[r] & A_{n-1} \ar[r] & \dots \ar[r] & A_0  \ar[r] &  k}$$
where each $A_i$ is a connective local artinian cdga over $k$. Moreover each 
morphism $H^0(A_n) \to H^0(A_{n-1})$ is a surjective map of artinian local $k$-algebras. The proposition
then becomes of statement about the $\s$-functor
$$\phi : \APerf(\hat{A}_x) \to \lim_{n}\APerf(A_n).$$
The $\s$-functor $\phi$ is the restriction of a left adjoint defined on the level of quasi-coherent
complexes
$$\phi : \QCoh(\hat{A}_x) \leftrightarrows \lim_n\QCoh(A_n) : \psi.$$
The left adjoint has already been described, it sends $E$ to $"\lim_{n}"A_n\otimes_{\hat{A}_x}E$. 
The right adjoint sends a pro-object $"\lim_{n}E_n"$ to its realization $\lim_n E_n \in \QCoh(\lim_n A_n)$.

A first observation is that the unit of the adjunction
$\hat{A}_x \to \psi\phi(\hat{A}_x)$ is obviously an equivalence, and thus this is also the case
for any $E \in \Perf(\hat{A}_x)$ (as $\hat{A}_x$ generates $\Perf(\hat{A}_x)$
by finite colimits, shifts and retracts). 

\begin{lem}\label{lpA-4-1}
Suppose that $E \in \QCoh(\hat{A}_x)$ is $k$-connective ($H^i(E) \simeq 0$ for $i>-k$) and
$H^k(E)$ is of finite type as a $H^0(\hat{A}_x)$-module. Then, $\psi\phi(E)$ remains
$k$-connective and the adjunction morphism
$E \to \psi\phi(E)$ induces an isomorphism
$$H^k(E) \simeq H^k(\psi\phi(E)).$$
\end{lem}

\textit{Proof of the lemma.} We use the Milnor short exact sequences 
$$\xymatrix{ 0 \ar[r] & 
\lim_{n}^1 H^{j-1}(E_n) \to \ar[r] & H^j(\lim_n E_n) \ar[r] & \lim^0_n H^j(E_n) \ar[r] & 0}$$
for any $j \in \ZZ$ and any projective system of complexes
$$\xymatrix{
\dots \ar[r]  & E_n \ar[r] & E_{n-1} \ar[r] & \dots \ar[r] & E_1 \ar[r] & E_0.
}$$ 
Here the
$\lim_n^i$ denote the $i$-th right derived functor associated to the $\lim_n$ functor, and 
$\lim_n E_n$ is understood in the $\s$-category $\dg_k$ of complexes of $k$-modules. To see that 
$\psi\phi(E)$ remains $k$-connective we thus have to show that $\lim_n^1H^{k}(E_n)=0$. 
But this follows from the fact that the system $n\mapsto H^{k}(E_n)$ is of the form 
$n\mapsto H^0(A_n)\otimes_{H^{0}(\hat{A}_x))}H^k(E)$ and thus all maps $H^k(E_n) \to H^k(E_{n-1})$
are surjective. Once $\psi\phi(E)$ is known to be $k$-connective, the unit of the adjunction, 
on the $k$-th cohomology group can be written as
$$H^k(E) \to \lim_n H^0(A_n)\otimes_{H^{0}(\hat{A}_x))}H^k(E).$$
This is an isomorphism because $H^{0}(\hat{A}_x)) \simeq \lim_n H^0(A_n)$ and
$H^k(E)$ is a module of finite type.
\hfill $\Box$ \\

Let $E \in \APerf(\hat{A}_x)$. By definition, for any $j$, we can find $E_0 \in \Perf(\hat{A}_x)$
and a morphism $E_0 \to E$ whose cone $E'$ is $j$-connective. We have a commutative diagram
with exact rows
$$\xymatrix{
E_{0} \ar[r] \ar[d] & E \ar[d]  \ar[r] & E' \ar[d] \\
\psi \phi (E_{0}) \ \ar[r] & \psi \phi (E) \ar[r] & \psi \phi (E').}$$
The left vertical morphism is an equivalence because $E_0$ is perfect. Moreover, 
by the lemma \ref{lpA-4-1} both $E'$ and $\phi(E')$ are $j$-connective. In particular, considering the
long exact sequences on cohomology groups we see that the vertical morphism of the middle must 
induce isomorphisms $H^k$ for all $k>-j$. As $j$ is arbitrary this shows the unit of the
adjunction $E \psi\phi(E)$ is an equivalence, and thus that the $\s$-functor $\phi$
of the proposition is fully faithful. 

To finish the proof we have to show that the right adjoint $\psi : \lim_n\QCoh(A_n) \leftrightarrows 
\QCoh(\hat{A}_x)$ sends objects in $\lim_n \APerf(A_n)$ to almost perfect complexes
and that the induced $\s$-functor
$$\psi : \lim_n \APerf(A_n) \to \APerf(\hat{A}_x)$$
is furthermore conservative. As we have seen during the proof of the lemma \ref{lpA-4-1}, 
the Milnor short exact sequences implies that $\psi$ preserves the $j$-connective objects.
Any object $E=\{E_n\} \in \lim_n \APerf(A_n)$ is $j$-connective for some $j$ (uniform in $n$). Indeed,
if $E_n \in \APerf(A_n)$, then it is $l$-connective for some integer $l$. Moreover, 
$A_n$ is an artinian local connective cdga, and thus the Yoneda lemma implies that 
if $E_n \otimes_{A_n}k$ is $j$-connective for some $j$ then $E_n$ was already $j$-connective as well.
In particular, if $E_0$ is $j$-connective then all $E_n$ are $j$-connective as well (for the same $j$), as
$$E_n \otimes_{A_n}k \simeq (E_n \otimes_{A_n}A_0)\otimes_{A_0}k \simeq E_0\otimes_{A_0}k.$$
Therefore, $\psi(E)$ is also $j$-connective. 

We now proceed as before by approximation by perfect complexes as before. 
For any $E=\{E_n\} \in \APerf(\hat{A}_x)$, which is $j$-connective for some $j$, 
we can find a free $A_0$-module $D_0=A_0^r$ and a morphism
$$u_0 : D_0[j] \to E_0$$
whose cones $E'_0$ is $(j+1)$-connective. As $D_0$ is free we can lift 
$u_0$ to a morphism $u : D[j] \to E$ whose cone $E'$ is $(j+1)$-connective, and $D$ is here
of the form $\phi(\hat{A}_x^r)$ (i.e. free of rank $r$). Continuing the resolution 
on $E'$, we construct inductively an object $F_r \in \Perf(\hat{A}_x)$
and a morphism $\phi(F_r) \to E$ whose cone $E_r'$ is now $(j+r)$-connective. For a fixed $r$,
we have a commutative diagram with exact rows
$$\xymatrix{
F_r \ar[r] & E   \ar[r] & E_r'  \\
\phi \psi (F_r) \ar[u] \ar[r] & \phi \psi (E)  \ar[u] \ar[r] & \ar[u] \phi \psi (E_r').}$$
As $E_r'$ and $\phi \psi (E_r')$ are both $(j+r)$-connective, we see that the vertical
morphism in the middle induces bijection on $H^l$ for all $l\geq -r-j+1$. As $r$ is arbitray, 
we deduce that $\phi\psi(E) \to E$ is an equivalence. 

Finally, the same exact sequence $F_r \to E \to E'_r$ show that the
cohomology groups of $\psi(E)$ are of finite type as $H^0(\hat{A}_x)$-modules, and thus that 
$\psi(E)$ is almost perfect as required.
\hfill $\Box$ \\

\section{Cotangent complexes of derived Artin stacks}
\label{secapp:CotangentcomplexesofderivedArtinstacks}

For a derived Artin stack $X$ over $k$, we have its big or fake cotangent complex
$\LL_{X/k}^{big}$. This is a non-necessarily quasi-coherent $\OO_X$-dg-module. 
As an $\s$-functor, it sends $u : \Spec\, A \to X$ to the $A$-dg-module 
$\LL_{A/k}$. There is a canonical pull-back morphism
$u^*\LL_{X/k} \to \LL_{A/k}$, which defines a morphism of $\OO_X$-dg-modules
$$\LL_{X/k} \to \LL_{X/k}^{big}$$
where $\LL_{X/k}$ is the global cotangent complex of $X$ relative to $k$ as
constructed in \cite[\S 1.4]{hagII}. As $\LL_{X/k}$ is itself quasi-coherent, we get a canonical
morphism of quasi-coherent complexes
$\LL_{X/k} \to (\LL_{X/k}^{big})^{qcoh}$. Considering exterior product, we also
get this way natural morphisms for all $i\geq 0$
$$\wedge^i_{\OO_X}\LL_{X/k} \to (\wedge^i_{\OO_X}\LL_{X/k}^{big})^{qcoh}.$$

\begin{thm}\label{tA-5-1}
Let $X$ be a derived Artin stack locally of finite presentation over $k$. 
For all $i\geq 0$ the natural morphism 
$$\wedge^i_{\OO_X}\LL_{X/k} \to (\wedge^i_{\OO_X}\LL_{X/k}^{big})^{qcoh}$$
is an equivalence of quasi-coherent complexes on $X$.
\end{thm}

\textit{Proof.} The proof is essentially the same as in \cite[Prop. 1.14]{ptvv}, together with the additional
result from \cite{mon}. Note that the statement of the theorem implies
\cite[Prop. 1.14]{ptvv} by taking global sections. \\

The proof goes by induction on the geometricity of $X$. When $X$ is affine, the statement follow
formally from the shape of the adjunction
$$i : \QCoh(X) \leftrightarrows \OO_X-\dg : (-)^{qcoh}.$$
Indeed, the $\s$-functor $E \mapsto E^{qcoh}$ is explicitely given by 
sending $E$ to the $A$-dg-module of global sections $\Gamma(X,E) \in A-\dg\simeq \QCoh(X)$.
From this, it is obvious that the canonical morphisms $\wedge^i_{\OO_X}\LL_{X/k} \to 
(\wedge^i_{\OO_X}\LL_{X/k}^{big})^{qcoh}$ are indeed equivalence. In the same manner,
the result remains true for disjoint union of derived affine schemes.

We assume now that the statement has been proved for $(n-1)$-geometric derived Artin stack. 
We let $f : U \to X$ be a smooth atlas with $U$ a disjoint union of affine derived schemes. 
Let $U_* \to X$ be the nerve of $f$, so each $U_i$ is a $(n-1)$-geometric derived Artin stack.
For a derived stack $Y$ we denote by $TY[1]$ its shifted total tangent stack. It is formally
defined by $\Map(\Spec\, k[\epsilon],Y)$, where $\epsilon$ sits in cohomological degree
$-1$. The derived stack $TY[1]$ sits over $Y$, and as such 
is of the form $\VV(\LL_{Y/k}[-1])$, the linear stack associated to the perfect
complex $\LL_{Y/k}[-1]$. This is a $\Gm$-equivariant derived stack over $Y$, 
and according to \cite{mon}, the projection $\pi : [TY[1]/\Gm] \longrightarrow Y\times B\Gm$
is such that the quasi-coherent push-forward $\pi_(\OO)$ is given by
$$\pi_*(\OO) \simeq \bigoplus_{p} Sym^p_{\OO_Y}(\LL_{Y/k}[-1]),$$
where the right hand side is a graded quasi-coherent complex on $Y$ 
defined by the natural grading of the $Sym$. In other words, the weight $i$ piece
of $\pi_*(\OO)$ is $\wedge^i\LL_{Y/k}[-i]$. 

It is very easy to check, for instance using the infinitesimal lifting properties, 
that $TU_0[1] \to TX[1]$ remains a smooth epimorphism. In particular, 
$TU_*[1]$, which is equivalent to the nerve of $TU_0[1] \to TX[1]$, induces
an equivalence of derived $\Gm$-equivariant stacks over $X$
$$TX[1] \simeq \colim_{p}TU_p[1].$$
If we denote by $p_i : TU_p[1]\to X$ the natural projection, and similarly $p : TX[1] \to X$, 
we thus have an equivalence of $\Gm$-equivariant quasi-coherent complexes
$p_*(\OO) \simeq \llim_{p}(p_p)_*(\OO).$
Applying the main result of \cite{mon} and considering the weight $p$ pieces, we get 
a canonical equivalence of quasi-coherent complexes
$$\wedge^i\LL_{X/k} \simeq \llim_{p} (f_p)_*(\wedge^i\LL_{U_p/k})$$
where $f_i : U_i \to X$ is the natural projection.

Because the natural morphism $\colim_p U_p \to X$ is an equivalence, we get, by 
descent, an equivalence of $\OO_X$-dg-modules
$\wedge^i\LL_{X/k}^{big} \simeq \llim_p (f_p)_*(f_p^*)(\wedge^i\LL_{X/k}^{big}).$
But we have a tautological identification $(f_p^*)(\wedge^i\LL_{X/k}^{big})\simeq 
\wedge^i\LL_{U_p/k}^{big}$, 
and thus we get
$\wedge^i\LL_{X/k}^{big} \simeq \llim_p (f_p)_*(\wedge^i\LL_{U_p/k}^{big}).$
As the $\s$-functor $(-)^{qcoh}$ is right adjoint it commutes with limits. Also, 
as pull-backs of $\OO$-dg-modules preserves quasi-coherent complexes, $(-)^{qcoh}$ 
also commutes with direct images. Finally, as each $U_p$ is $(n-1)$-geometric, we know
that the theorem holds for the $U_p$'s.
We thus get 
a natural equivalence
$$(\wedge^i\LL_{X/k}^{big})^{qcoh} \simeq \llim_p (f_p)_*(\wedge^i\LL_{U_p/k}).$$
But we have already seen that the right hand side is naturally equivalent to 
$\wedge^i\LL_{X/k}$, and thus we have 
$$(\wedge^i\LL_{X/k}^{big})^{qcoh} \simeq \wedge^i\LL_{X/k}.$$
Unfolding the various equivalences we see that this equivalence is indeed induced by the
canonical morphism $\wedge^i\LL_{X/k} \to \wedge^i\LL_{X/k}^{big}$, which proves the theorem.
\hfill $\Box$ \\

\section{Cotangent complexes of formal completion}\label{secapp:Cotangentcomplexesofformalcompletion}

As in the previous section let $X$ be a derived Artin stack of finite presentation over $k$.
We assume here that $k$ is a field and that we have given ourselves a global
point $x : \Spec\, k \to X$ (over more general rings $k$, we can always reduce to this special case by base 
change). 
As explained in \S \ref{subsec:Derivedfoliationsonformalmoduliproblems} we have the formal completion 
$(X,x)^\wedge$. This formal completion 
can be described as a ringed site as follows. We let $X^{art}$ be the $\s$-category 
$\dgart^*_k/X$, of derived local Artinian affine schemes over $X$ (preserving the base point).
The $\s$-category is a ringed site, for the induced etale topology and for the canonical
stack of simplicial rings sending $\Spec\, A \to X$ to $A$. We let 
$\OO_{X^{art}}-\dg$ be the $\s$-category of all $\OO_{X^{art}}$-dg-modules. We have a full 
sub-$\s$-category $\QCoh(X^{art}) \subset\OO_{X^{art}}-\dg$ consisting of
objects $E$ such that for all morphism in $X^{art}$
$$\xymatrix{
\Spec\, A \ar[rr] \ar[rd] & & \Spec\, B \ar[dl] \\
 & X & }$$
such that the induced morphism $E(B)\otimes_B A \longrightarrow E(A)$ is a quasi-isomorphism. 
Note that by definition $\QCoh(X^{art})$ identifies canonically with 
$\QCoh((X,x)^{\wedge})$, the $\s$-category of quasi-coherent complexes on the formal
completion $(X,x)^{\wedge}$ of $X$ at $x$.

By Freyd representability the natural inlcusion $j : \QCoh(X^{art}) \subset\OO_{X^{art}}-\dg$
possesses a right adjoint denoted $E \mapsto E^{qcoh}$. We also have a canonical restriction $\s$-functor
$$\OO_X-\dg \longrightarrow \OO_{X^{art}}-\dg,$$
and thus the cotangent complex of $X$ (relative to $k$) 
provides a canonical object $\LL_{X^{art}/k} \in \OO_{X^{art}}-\dg$
called the \emph{formal cotangent complex of $X$ at $x$}. It is obviously 
a quasi-coherent complex over $X^{art}$. In the same manner, we have 
the big or fake cotangent complex $\LL_{X^{art}/k}^{big} \in \OO_{X^{art}}-\dg$ given by
$\LL_{X^{art}/k}(\Spec\, A):=\LL_{A/k}$. \\

The following theorem is the formal version of our previous result theorem \ref{tA-5-1}.

\begin{thm}\label{tA-6-1}
With the notations and conditions as above, for all $i\geq 0$ the canonical morphism
$$\wedge^i_{\OO_{X^{art}}}\LL_{X^{art}/k} \longrightarrow
(\wedge^i_{\OO_{X^{art}}}\LL_{X^{art}/k}^{big})^{qcoh}.$$
\end{thm}
 
\textit{Proof.} We proceed as in theorem by induction on the geometricity of $X$. We let $f_* : U_* \to X$
be the nerve of a smooth atlas for $X$ with $U_0$ a disjoint union of derived affine schemes. Moreover,
by galois descent we can replace $k$ by a algebraic separable extension $k'/k$, and thus
assume that $x$ lifts to a $k$-point $u$ in $U_*$.

As we have
already seen during the proof of \ref{tA-5-1} we have equivalences in $\OO_X-\dg$ and in $\QCoh(X)$
$$\wedge^i\LL_{X/k}^{big} \simeq \llim_p (f_p)_*(\wedge^i\LL_{U_p/k}^{big})$$
$$\wedge^i\LL_{X/k} \simeq \llim_p (f_p)_*\wedge^i\LL_{U_p/k}.$$
As the restriction $\s$-functor $\OO_X-\dg \to \OO_{X^{art}}-\dg$ commutes with limits, we see that 
we have a commutative square in $\OO_{X^{art}}-\dg$
$$\xymatrix{
\wedge^i_{\OO_{X^{art}}}\LL_{X^{art}/k} \ar[d] \ar[r] & 
\llim_p (\hat{f}_p)_*\wedge^i\LL_{U_p^{art}/k} \ar[d] \\
(\wedge^i_{\OO_{X^{art}}}\LL_{X^{art}/k}^{big})^{qcoh} \ar[r] & \llim_p (\hat{f}_p)_*(\wedge^i\LL_{U_p^{art}/k}^{big})^{qcoh},
}$$
where $\hat{f}_* : U_*^{art} \to X^{art}$ is the induced atlas on formal completions. Using 
induction on the geometricity of $X$ we see that we can reduce the theorem to 
the case $X$ is disjoint union of affines. Restricting to the component containing $x$ we can 
further assume that $X$ is an affine derived scheme $X=\Spec\, A$, where $A$ is of finite presentation over
$k$.

We now fix $i\geq 0$, and consider the exact triangle of $\OO_{X^{art}}$-dg-modules
$\wedge^i\LL_{X^{art}} \to \wedge^i\LL_{X^{art}}^{big} \to E$. By construction, the $\OO_{X^{art}}$-dg-module
$E$ sends $\Spec\, A' \to X$ to the cone of $\wedge^i\LL_{A}\otimes_A A' \to \wedge^i\LL_{A'}$.
For any quasi-coherent complex $E_0 \in \QCoh(X^{art})$ we have a natural equivalence
$$\Map_{\OO_{X^{art}}-\dg}(E_0,E) \simeq \llim_{A \to A'} \Map_{A'-\dg}(E_0(A'),E(A')).$$
In other words, this mapping space if also the mapping space of pro-$\underline{\hat{A}}$-modules
between $"\llim_{A \to A'}"E_0(A')$ and $"\llim_{A \to A'}"E(A')$. By Corollary \ref{cpA-3-1-2} we know that 
$E$ is equivalent to zero as a pro-object, and thus that $\Map_{\OO_{X^{art}}-\dg}(E_0,E) \simeq 0$
for any quasi-coherent complex $E_0$, or equivalently that $E^{qcoh}\simeq 0$.
This shows as required that 
$\wedge^i\LL_{X^{art}} \simeq (\wedge^i\LL_{X^{art}}^{big})^{qcoh}$.
\hfill $\Box$ \\

\section{Deformation to the normal bundle}\label{secapp:Deformationtothenormalcone}

In this section we denote by $\A$ the quotient stack $[\mathbb{A}^1/\Gm]$. 
As explained in \cite{Moulinos}, derived stacks over $\A$ should be considered as
filtered derived stacks. We remind in particular 
(see \S \ref{subsec:Gradedmixedcomplexesandfilteredcomplexes}), that 
there exists a natural symmetric monoidal equivalence of $\s$-categories
$$\QCoh(\A) \simeq \fdg_k,$$
between quasi-coherent complexes over $\AA$ and filtered complexes of $k$-modules.
This equivalence is induced by a natural (decreasing) filtered object $\OO(-*)$
in $\QCoh(\A)$, namely the families of $\OO(-i)$ for $i \in \ZZ$. Here $\OO(-i)$ stands
for the line bundle on $\A$ which is trivial on $\mathbb{A}^1$ together with its
canonical $\Gm$-action of weight $-i$. The families of $\OO(-i)$ is
endowed with canonical
morphisms $\OO(-i) \to \OO(-i+1)$ induced by multiplication by $X \in \Gamma(\AA,\OO(1))$.
The $\s$-functor $\rh(\OO(*),-) : \QCoh(\AA) \longrightarrow \fdg_k$ is the equivalence
mentionned above.

We consider the inclusion of the origin $\Spec\, k \to \mathbb{A}^1$, and the induced morphism
on quotient stacks 
$$j : B\Gm \hookrightarrow [\mathbb{A}^1/\Gm]=\A.$$
As an object over $\A$, $B\Gm$ will be denoted by $\cQ \in \dSt_k/\A$. This is a 
relative derived affine scheme, whose fiber at any non-zero point is empty, and whose
fiber at $0$ is $\Spec\, k[\epsilon_{-1}]$, where $\epsilon_{-1}$ is of degree $-1$
and $\epsilon_{-1}^2=0$. The induced $\Gm$-action acts by weight $1$ on $\epsilon_{-1}$.

We now let $X \to Y$ be any morphism of derived stacks. we consider
$Y \times \cQ$ and $X \times \A$ as 
objects over $Y \times \A$, and we form the relative 
derived mapping stack
$$Def(X/Y):=\uMap_{/Y\times \A}(Y \times \cQ,X \times \A) \in \dSt_k/Y\times \A.$$

\begin{df}
The \emph{deformation to the normal bundle of a morphism $f : X \to Y$} is the derived stack
defined above
$$Def(X/Y)=Def(f):=\uMap_{/Y\times \A}(Y \times \cQ,X \times \A) \in \dSt_k/Y\times \A.$$
\end{df}

The deformation to the normal bundle $Def(X/Y)$ has the following two important properties.
Its fiber over the open $Y=[\Gm/\Gm]\times Y \hookrightarrow \A\times Y$
is $\Map_{/Y}(\emptyset,X)\simeq Y$, whereas its fiber at zero
$Y \times B\Gm \hookrightarrow Y \times \A$ is $[\Map_{/Y}(Y\times \Spec\, k[\epsilon_{-1}],X)/\Gm]$.
Note that $[\Map_{/Y}(Y\times \Spec\, k[\epsilon_{-1}],X)$ is identified with the
shifted relative tangent bundle $T(X/Y)[1]$, where $\Gm$ acts by weight $-1$ dilatation
along the fibers of $T(X/Y)[1] \to X$. When $X$ and $Y$ are derived Artin stacks 
locally of finite presentation, $T(X/Y)[1]$ is also the linear stack $\VV(\LL_{X/Y}[-1])$, 
associated to the perfect complex $\LL_{X/Y}[-1]$ on $X$, endowed with its natural
$\Gm$-action. Therefore, $Def(X/Y)$ provides indeed 
an interpolation, by means of filtration, between $X$ and $T(X/Y)[1]$, which 
is at the origin of the relative Hodge filtration.

Note that the projection $Y \times \cQ \to Y \times \A$ induces on derived mapping stacks
a morphism $j : X\times \A \to Def(f)$. Over the open point $[\Gm/\Gm]$ of $\A$ the morphism 
$j$ recovers  $f : X \to Y$, whereas over $B\Gm$ it is the zero section $X \to T(X/Y)[1]$.

We finish this section by the following comment. When $f : X=\VV(E') \to Y=\VV(E)$ is a morphism
of derived linear stacks over $\Spec\, A$, induced by a morphism of bounded
above $A$-dg-modules
$u : E \to E'$, then $Def(f)$ is itself globally, as a derived stack over $\Spec\, A \times \A$,
a derived linear stack. Indeed, this can be seen directly using that for any derived stack 
$Z \to \Spec\, A$ we have functorial equivalences
$$\Map_{\dSt_A}(Z,\VV(E))\simeq \Map_{\dg_A}(E,\OO_Z(Z)).$$
In this case, the derived linear stack $Def(f)$, relative to $\Spec\, A\times \A$, 
can be described explicitly as
$\VV(u)$, where $u : E \to E'$ is now considered as a one step filtered quasi-coherent complex 
on $\Spec\, A$, and thus as a quasi-coherent complex over $\Spec\, A \times \A$.

\section{Internal derived de Rham theory}\label{secapp:InternalderiveddeRhamtheory}

We remind from \cite{mrt} the group stack (over $k$) $\cH_0:=B\Ga\ltimes \Gm$, semi-direct product 
of the multiplicative group acting, by its standard weight 1 action, on $B\Ga$.
As we have seen, the $\s$-category $\QCoh(B\cH_0)$ is canonically equivalent to the
$\s$-catgeory of graded mixed complexes 
$$\QCoh(B\cH_0) \simeq \egrdg_k.$$
In the same manner, there is an equivalence of $\s$-categories of commutative algebras
$$\OO_{B\cH_0}-\cdga_{qcoh} \simeq \egrcdga_k,$$
where $\OO_{B\cH_0}-\cdga_{qcoh}$ is the $\s$-category of quasi-coherent cdga's over the stack 
$B\cH_0$, defined
by the usual left Kan extension formula
$$\OO_{B\cH_0}-\cdga_{qcoh}=\lim_{\Spec\, A \to B\cH_0}\cdga_A.$$
For any derived stack $X$, we define the $\s$-category of 
quasi-coherent graded mixed $\OO_X$-cdga's
by 
$$\OO_X-\egrcdga_{qcoh}:=\OO_{X\times B\cH_0}-\cdga_{qcoh}.$$
Note that, by definition, an object in $\OO_X-\egrcdga_{qcoh}$ 
is a sheaf of graded mixed cdga's over $X$
which is moreover linear over $\OO_X$ and quasi-coherent. That is, it is given by a family of
graded mixed $A$-linear cdga $D_A$, for each $\Spec\, A \to X$, together with 
compatible by base change equivalences
$$D_A \otimes_A B \simeq D_B$$
for each $\Spec\, B \to \Spec\, A \to X$. 
\begin{rmk}
\emph{One should not confuse 
$\OO_X-\egrcdga_{qcoh}$ with the notion of graded mixed cdga's over $X$ used in the definition
of derived foliations over the derived stack $X$: the latter objects are not 
quasi-coherent by any means. We hope that the wording 
"quasi-coherent graded mixed $\OO_X$-cdga" will help avoiding the confusion. }
\end{rmk}

Taking the weight $0$ part defines an $\s$-functor
$\OO_X-\egrcdga_{qcoh} \longrightarrow \OO_X-\cdga_{qcoh},$
where $\OO_X-\cdga_{qcoh}$ is the $\s$-category of sheaves of quasi-coherent
$\OO_X$-linear cdga's over $X$
(defined by $\OO_X-\cdga_{qcoh}=\lim_{\Spec\, A \to F}\cdga_A$). This forgetful $\s$-functor
possesses a left adjoint called the relative or internal de Rham algebra construction
$$\DR^{int}(-/\OO_X) : \OO_X-\cdga_{qcoh} \longrightarrow \OO_X-\egrcdga_{qcoh}.$$
Concretely, if $B \in \OO_X-\cdga_{qcoh}$, then for 
any $\Spec\, A \to X$ we consider $\DR(B_A/A) \in \egrcdga_A$. The base change
equivalences for $B$ induces base change equivalences for $\DR(B_A/A)$ making 
it into a quasi-coherent $\OO_X$-linear graded mixed cdga. \\

We will now be interested in the specific case where $X=B\cH_0$. We get this way an $\s$-functor
$$\OO_{B\cH_0}-\cdga_{qcoh} \longrightarrow \OO_{B\cH_0}-\egrcdga_{qcoh}\simeq 
\OO_{B\cH_0 \times B\cH_0}\cdga_{qcoh}.$$
Note that, using the same argument as in the proof of the equivalence
$\OO_{B\cH_0}-\cdga_{qcoh} \simeq \egrcdga$, the $\s$-category 
$\OO_{B\cH_0 \times B\cH_0}-\cdga_{qcoh}$ is naturally equivalent to the $\s$-category of \emph{double
graded mixed cdga}: that a cdga endowed with two multiplicative $\ZZ$-graduations
and two compatible mixed structures $\epsilon_1$ and $\epsilon_2$ such that 
$\epsilon_1\epsilon_2 + \epsilon_2\epsilon_1 = 0$ ($\epsilon_1$ is of weight bigrading $(1,0)$ and
$\epsilon_2$ of bigrading $(0,1)$, and both are of cohomological grading $-1$). 
The internal de Rham construction 
is then explicitly given by sending an graded mixed cdga $A$ to its
graded de Rham algebra $\DR^{gr}(A/k)$, which is endowed with a double graded mixed structure:
$\epsilon_1$ is induced by the original mixed structure on $A$, and $\epsilon_2$ is the 
de Rham differential of $\DR^{gr}(A/k)$.

Double graded mixed cdga's have two different realizations constructions, depending if
we realize the first of the second mixed structure. These corresponds to considering direct 
images along the two projections
$$B\cH_0 \leftarrow B\cH_0 \times B\cH_0 \rightarrow B\cH_0$$
which induced the two (non-Tate) realization $\s$-functors
$$|-|_1 \,,  |-|_2 : \OO_{B\cH_0 \times B\cH_0}-\cdga_{qcoh} \longrightarrow \OO_{B\cH_0}-\cdga_{qcoh}.$$
Applied to $\DR^{int}(A/k)$, for $A$ a graded mixed cdga, we get this way two possible 
de Rham complexes of $A$
$$\CDR^{I}(A/k):=|\DR(A/k)|_1 \qquad \CDR^{II}(A/k):=|\DR(A/k)|_2,$$
which both of them are graded mixed cdga's. 

\begin{df}\label{dA-8-1}
Let $A$ be a graded mixed cdga (over $k$). 
\begin{enumerate}
\item The \emph{(derived) de Rham complex of the first kind 
of $A$
relative to $k$} is the graded mixed cdga given by
$$\CDR^{I}(A/k):=|\DR^{int}(A/k)|_1.$$
\item The \emph{(derived) de Rham complex of the second kind 
of $A$
relative to $k$} is the graded mixed cdga given by
$$\CDR^{II}(A/k):=|\DR^{int}(A/k)|_2.$$
\item The \emph{absolute (derived) de Rham complex
of $A$
relative to $k$} is the graded mixed cdga given by
$$\CDR^{abs}(A/k):=||\DR^{int}(A/k)|_1|_2 \simeq ||\DR^{int}(A/k)|_2|_1.$$
\end{enumerate}
\end{df}

Visually, the two versions of the de Rham complexes of definition above can be 
summarized as follows.

$$\CDR^{I}(A/k) = |A| \to |\LL_{A/k}| \to |\wedge_A^2\LL_{A/k}|
\to \dots |\wedge_A^i\LL_{A/k}| \to |\wedge_A^{i+1}\LL_{A/k}| \to \dots$$
$$\CDR^{II}(A/k) = |A \to \LL_{A/k} \to \wedge_A^2\LL_{A/k}
\to \dots \wedge_A^i\LL_{A/k} \to \wedge_A^{i+1}\LL_{A/k} \to \dots|$$
Note also that the two graded mixed cdga's $\CDR^{I}(A/k)$
and $\CDR^{II}(A/k)$ are different, but their realizations are canonically equivalent
as cdga's
$$|\CDR^{I}(A/k)| \simeq |\CDR^{II}(A/k)|.$$

In a more explicit form, $\CDR^{I}(A/k)$ is a graded mixed cdgas whose
piece of weight $i$ is $|\wedge^i \LL_{A/k}|$, where 
$\LL_{A/k}$ is the graded cotangent complexe of $A$, endowed with the induced
mixed structure $\epsilon(dR(a))= - dR(\epsilon(a))$, where $\epsilon$ is the original mixed
structure on $A$. The mixed structure is then induced by the de Rham
differential of $A$.

On the other hand, the graded piece of weight $i$ of
$\CDR^{II}(A/k)$ is given by the complex
$$\prod_{j\geq 0}(\wedge^j\LL_{A/k})^{(i)}[-j],$$
endowed with the total differential $d+dR$. The mixed structure
relating these graded pieces is 
then induced by the mixed structure $\epsilon$ of $A$.
Note that for each $i$, we have
a decomposition
$$(\wedge^j\LL_{A/k})^{(i)} \simeq \bigoplus_{n_1k_1+\dots+n_jk_j=i} \bigotimes_{l} 
\wedge^{k_l}(\LL_{A/k}^{(n_l)})[-j],$$
so the weight $i$ piece of $\CDR^{*,II}(A/k)$ can also be written
$$\prod_{j\geq 0}(\bigoplus_{n_1k_1+\dots+n_jk_j=i} \bigotimes_{l} 
\wedge^{k_l}(\LL_{A/k}^{(n_l)})[-j]).$$

As an example of applications of the notions of internal de Rham complex we prove the following
well known lemma.

\begin{prop}\label{pA-8-1}
Let $A \to k$ be a local augmented $k$-linear Artinian dg-algebra, almost of
finite presentation over $k$. Then, the canonical
morphism induces a quasi-isomorphism of cdga's
$$k\simeq \CDR(A/k).$$
\end{prop}

\textit{Proof.} Let $D:=\DR(k/A) \in \egrcdga_k$, and we consider the canonical morphism
$A \to \DR(k/A)$, where $A$ is considered as a graded mixed cdga of weight $0$ with trivial
mixed structure. This morphism induces a morphism on internal de Rham algebras
$$\phi : \DR(A/k) \simeq \DR^{int}(A/k) \to \DR^{int}(D/k).$$
This is a morphism of double graded mixed cdga's, which can be realized 
in the absolute sense to get a morphism of cdga's
$$\CDR(A/k) \to \CDR^{abs}(D/k).$$
This morphism can be described as first realizing along the first graded mixed structure and then
realizing along the second one. The morphism $|\phi|_1$ is given by
$\DR(A/k) \to \CDR^{I}(D/k)$, as the first graded mixed structure on $\DR(A/k)$ is trivial of weight $0$.
On the peice of weight $i$ the morphism $|\phi|_1$ is the canonical morphism
$$\wedge^i \LL_{A/k}[i] \longrightarrow |\wedge^i \LL_{D/k}[i]|_1.$$
For $i=0$, this morphism is the canonical morphism of cdga's
$A \to \CDR(k/A)$, which is a quasi-isomorphism as shown in \cite[Lem. 2.2.4]{cptvv}. The morphism 
$|\phi|_1$ in weights $i>0$ is the canonical morphism
$$\wedge^i \LL_{A/k} \longrightarrow \CDR(k/A,\wedge^i \LL_{A/k}\otimes_A k),$$
where the right hand side is the derived de Rham cohomology of $k$ relative 
to $A$, with coefficients in $\LL_{A/k}\otimes_A k$. As $\mathbb{L}_{A/k}$ is almost perfect
$A$-dg-module, we also know that this morphism is a quasi-isomorphism (see \cite[Cor. 2.2.5]{cptvv}
for the case of perfect modules which can be extended to almost perfect modules in a straightforward 
manner by postnikov approximation).

This shows that $|\phi|_1 : \DR(A/k) \to \CDR^{I}(D/k)$ is a quasi-isomorphism
of graded mixed cdga's, and thus induces an equivalence after realizations
$\CDR(A/k) \simeq \CDR^{abs}(D/k)$. Now, $\CDR^{abs}(D/k)$ can also be 
computed as$||\DR^{int}(D/k)|_2|_1$. As a graded cdga, we have
$D\simeq Sym_k(E[2])$, where $E:=\LL_{A/k} \otimes_A k$. Therefore, its de Rham
cohomology is trivial, and we thus have $|\DR^{int}(D/k)|_2 \simeq k$.
This shows that 
$\CDR^{abs}(D/k) \simeq k$ and finishes the proof of the Proposition.
\hfill $\Box$

\begin{rmk}\label{rem:pA-8-1}
\emph{Proposition \ref{pA-8-1} is valid in a more general context, with a similar proof. For any 
derived affine $k$-scheme $X=\Spec\, A$ almost of finite presentation over $k$, the natural
morphism}
$$\CDR(A/k) \to \CDR(A_{red}/k)$$
\emph{is an equivalence.}
\end{rmk}

Another useful statement is the following reconstruction result.

\begin{prop}\label{pA-8-2}
Let $B$ be a graded mixed cdga with $B^{(0)}=A$ connective and 
$Sym_{A}(B^{(1)}) \simeq B$ as graded cdgas. Then, there is a canonical equivalence of graded mixed
cdga's
$$B \simeq \CDR^I(A/B).$$
\end{prop}

\textit{Proof.} We first notice that the canonical morphism 
$B \to \DR^{int}(A/B)$ induces an equivalence of graded mixed cdga's
$B \simeq \CDR^{II}(A/B)$. Indeed, this statement is independent of the mixed structure on $B$, 
as it is enough to show that $B \to \CDR^{II}(A/B)$ is an equivalence of graded cdga's. To see this
it is enough to write down $B=Sym_A(E)$ for a perfect $A$-dg-module $E$, and write
$A$ as $Sym_A(E')$ where $E'$ is the cone of the identity $E \to E$. as a graded cdga, the de Rham 
cohomology of the second kind $\CDR^{II}(A/B)$ can be modeled by the (graded) de Rham complex
of $Sym_A(E')$ over $B$, which is explicitly given by $Sym_A(K[-1])$ where $K$ is the cone
of the natural projection $E' \to E[1]$. The $A$-dg-module $K$ is canonically quasi-isomorphic
to $E[1]$ and the corresponding morphism $B \to Sym_A(K[-1])$ sends $E$ to $K[-1]$ by its canonical
identification and thus is an equivalence.

To finish the proof of the Proposition it is thus enough to identify $\CDR^I(A/B)$
canonically with $\CDR^{II}(A/B)$. For this we use the following statement about 
double graded mixed complexes. 

\begin{lem}\label{lpA-8-2}
There exists a fully faithful symmetric monoidal $\s$-functor
$$\QCoh(B\cH_0) \hookrightarrow \QCoh(B\cH_0 \times B\cH_0)$$
from graded mixed complexes to double graded mixed complexes, whose essential
image consists of all object $E$ which, as a bigraded complex, is concentrated on the diagonal:
$E^{(i,j)} \simeq 0$ for all $i\neq j$.
\end{lem}

\textit{Proof of the lemma.} The most easy proof consists of using the equivalence between 
graded mixed complexes and complete filtered complexes (see Proposition \ref{pI-3} and Corollary 
\ref{clI-4}). Via this translation the lemma
becomes obvious: the right adjoint $d_*$ consists of sending 
a filtered complex to the corresponding "diagonal" bilfiltered complex. More explicitely, 
if $\{F^j\}_{j\in \ZZ}$ is a filtered complex, the corresponding bifiltered complex $X$
is given by $X^{i,j}:=F^{min(i,j)}$ with the natural transition morphisms.
In terms of functors category, if we denote by $\ZZ$ the poset of integers, 
this full embedding is given by 
$$d_! : Fun(\ZZ,\dg_k) \longrightarrow  Fun(\ZZ\times \ZZ,\dg_k)$$
given by left Kan extension along the diagonal $d : \ZZ \hookrightarrow \ZZ \times \ZZ$.
This is clearly a fully faithful construction
and its essential image correspond precisely to complete bifiltered complexes $E$ whose
associated bigraded are concentrated on the diagonal.
\hfill $\Box$ \\
 
The full embedding of the previous lemma being symmetric monoidal $\s$-functor and thus
extends to a similar statement concerning commutative algebras object. This means that 
in particular, that any double graded mixed cdga $D$, which is concentrated on the diagonal
as a bi-graded complex, is naturally symmetric with respect to the action of $\ZZ/2$ which exchanges the
two gradings and the two mixed structures. In particular, for any such object $D$, 
we have a canonical equivalence $|D|_1 \simeq |D|_2$, and both are equivalent to 
$|E|$, where $E$ is the graded mixed cdga obtained from the lemma \ref{lpA-8-2}. 
This can therefore be applied to $D=\DR^{int}(A/B)$, which as a bi-graded cdga 
is equivalent to $Sym_A(E[2])$ and thus is indeed concentrated on the diagonal.
\hfill $\Box$ \\

\section{Derived affine stacks}\label{app:derivedaffine}

\subsection{Cosimplicial-simplicial algebras}

We let denote by $sCR$ be the category of simplicial commutative rings. It is endowed with
its usual \emph{model category} structure for which equivalences and fibrations are defined on the
underlying simplicial sets.

Similarly, we denote by $csCR$ the category of cosimplicial objects in $sCR$, or equivalently 
of cosimplicial-simplicial commutative rings. Its objects will be denoted by $A^*_*$, where the
superscript will refer to the cosimplicial direction and subscript to the simplicial direction.
The category $csCR$ can be endowed with a \emph{tensored} and \emph{cotensored} structure over the category $sSet$
of simplicial sets. We warn the reader here that there are several natural such structures and that 
the choice we make here is important for the sequel. 
For a simplicial set $K$ and an object $A$
of $csCR$, we define the cotensored structure by 
$A^K$ by the formula
$$(A^K)^q_p:=(A^q_p)^{K_q}$$
for $[p] \in \Delta^{op}$ and $[q] \in \Delta$. 
The tensored structure $K \otimes A$, for $A \in csCR$ and $K \in sSets$, is defined by the usual adjunction formula
$Hom(K\otimes A,B) \simeq Hom(A,B^K)$, for arbitrary $B \in csCR$. Similarly, we have a canonical simplicial enrichment
for which the simplicial Homs are denoted by $\uHom$ and are defined by 
$$\uHom(A,B)_p:=Hom_{csCR}(A,B^{\Delta^p}).$$

By considering the natural embedding $Mod_\ZZ \subset C(\ZZ)$, from $\ZZ$-modules to cochain complexes, 
any object $A \in csCR$ can be considered, by forgetting the multiplicative structure, 
as a cosimplicial-simplicial object in $C(\ZZ)$. We denote by $Tot^{\pi}(A)$
the complex equals to $\prod_{q-p=n}A^q_p$ in cohomological degree $n$
and whose differential is given by the alternated sums of the faces and cofaces of the
cosimplicial-simplicial diagram. Equivalently, $A$ can be first turned into a 
bi-complex, by individual totalisations along the cosimplicial and the simplicial directions, and
$Tot^\pi(A)$ can then be realized as the product-total complex of this bicomplex. 
Note also that $Tot^\pi(A)$ is a model for the homotopy limit-colimit of $A$
$$Tot^\pi(A) \simeq \mathrm{holim}_{[p] \in \Delta^{op}}\left(\mathrm{hocolim}_{[q] \in \Delta} A^q_p \right) \in C(\ZZ).$$

\begin{df}
A morphism $A \to B$ in $csCR$ is a \emph{completed quasi-isomorphism}
if the induced morphism $Tot^\pi(A) \to Tot^\pi(B)$ is a quasi-isomorphism of complexes.
\end{df}

We warn the reader that even though the categories $csCR$ and $scCR$, of cosimplicial-simplicial
commutative rings and of simplicial-cosimplicial rings, are equivalent, the notion 
of completed quasi-isomorphisms is not the same as the notion of weak equivalences used
for instance in \cite{bracanu}.

\begin{thm}\label{tsome}
The category $csCR$ is endowed with a model category structure for which the weak equivalences
are the completed quasi-isomorphism, and the fibrations are the epimorphisms. This model 
category structure is moreover cofibrantly generated and is a simplicial model category for
the above mentioned simplicial enrichment.
\end{thm}

\textit{Proof.} As fibrations\footnote{Note that $f:A\to B$ in $csCR$ is surjective iff $Tot^\pi(f)$ is.} and equivalences are defined via the functor
$Tot^\pi$, the proof of the theorem consists of transferring the projective model structure 
on $C(\ZZ)$ to $csCR$ via the functor $Tot^\pi$. For this we apply the path object
argument (see for instance \cite[\S 2.6]{MR2016697}).

The functor $Tot^\pi$ is a right adjoint whose left adjoint sends a complex $E \in C(\ZZ)$ to 
the free commutative ring over the cosimplicial-simplicial module $\phi(E)$ given by
$$(p,q) \mapsto (C^*(\Delta^p) \otimes C_*(\Delta^q) \otimes E)^0/Im\, d^{-1},$$
where $C^*(\Delta^p)$ and $C_*(\Delta^q)$ are the cohomology and homology complexes
of the standard simplex $\Delta^p$ and $\Delta^q$. This left adjoint 
will be denoted by $L\phi$, as being the composition of 
the free commutative ring functor $L : csMod \to csCR$, with the
functor $\phi : C(\ZZ) \to csMod$ defined by the above formula.
Clearly, the functor $Tot^\pi$ preserves small objects (but does not 
commute with filtered colimits in general). In order to apply \cite[\S 2.6]{MR2016697} we 
thus simply have to prove that $csCR$ possesses a fibrant replacement functor
and that all fibrant object possesses a path object. But all objects
are fibrant by definition, and we are thus reduced to show the existence
of path objects. For this we use the following lemma.

\begin{lem}\label{l1}
For any simplicial set $K$ and any $A \in csCR$, there exists a canonical 
quasi-isomorphism
$$Tot^\pi(A^K) \simeq Tot^\pi(A)^K,$$
where for $E \in C(\ZZ)$ and $K \in sSet$ we define
$$E^K :=\mathrm{holim}_{[p] \in \Delta} E^{K_p}.$$
\end{lem}

\textit{Proof of the lemma.} We remind that for any functor $F : \Delta \times \Delta \to \dg$
(or more generally any bi-cosimplicial object in a complete $\s$-category), 
the homotopy limit of $F$ can be computed, up to a quasi-isomorphism, 
either component-wise, or diagonally
$$\mathrm{holim}_{[p] \in \Delta} F(p,p) \simeq \mathrm{holim}_{[p]\in \Delta} \mathrm{holim}_{[q] \in \Delta} F(p,q).$$
We apply this to the functor sending $([p],[q])$ to the complex $N(A^q_*)^{K_p}$, where 
$N : sMod \to dg$ is the normalization construction. Computing the homotopy limit
diagonally yields $Tot^\pi(A^K)$, whereas the double homotopy limit construction yields $Tot^\pi(A)^K$.
\hfill $\Box$ \\

Lemma \ref{l1} now implies the existence of path objects. Indeed, for $A \in csCR$ an explicit
path object is given by $A \to A^{\Delta^1} \to A^{\partial \Delta^1}= A \times A$. The fact that
the constant map morphism $A \to A^{\Delta^1}$ is a completed quasi-isomorphism follows from
the lemma \ref{l1}. 
The fact that $A^{\Delta^1} \to A^{\partial \Delta^1}= A \times A$ is a fibration simply follows from the fact that
it is a levelwise epimorphism as this can be checked easily using the explicit formula for 
$A^K$. 

In order to finish the proof of the theorem, we are left to showing that the simplicial enrichement
is compatible with the model category structure. In other words, we must show that for
any cofibration of simplicial sets $i : K \hookrightarrow L$ and any fibration $f : A \to B$ in $csCR$, 
the induced morphism
$$f^i : A^L \longrightarrow B^L \times_{B^K}A^K$$
is a fibration, which is also an equivalence if $i$ is a trivial cofibration or $f$ is a trivial
fibration. The fact that the morphism $f^i$ is a fibration is clear because fibrations are epimorphisms
and because of the explicit formula for the exponentiation by a simplicial set in $csCR$. Finally, 
the fact that $f^i$ is also a weak equivalence when $i$ or $f$ is so simply follows from 
Lemma \ref{l1}.
\hfill $\Box$ \\

\begin{df}
The \emph{$\s$-category of cs-rings} is the $\s$-category 
obtained from $csCR$ by inverting the completed quasi-isomorphisms. It is denoted by
$\csCR$ (or by $\csCR_\mathbb{U}$ if one wants to specify cosimplicial-simplicial rings belonging
to a given universe $\mathbb{U}$).
\end{df}

\begin{rmk}
\emph{We note also that the model structure on $csCR$ can be equivalently be
obtained by first constructing a model structure on $csMod_\ZZ$, the category of cosimplicial-simplicial 
abelian groups, where again the equivalences and fibrations are defined via
the $Tot^\pi$ functor. The model structure on $csCR$ can then be obtained by
transferring along the forgetful functor $csCR \to csMod_\ZZ$, with left adjoint
given by the free commutative cosimplicial-simplicial ring.}
\end{rmk}

\begin{rmk}\emph{It is possible to a construct an $\s$-functor
$Tot^\pi : \csCR \longrightarrow LSym-Alg$, from our $\s$-category of cosimplicial-simplicial
commutative rings to the $\s$-category of \emph{derived rings}, namely modules over the
$LSym$-monad. Here we denote by $LSym$ the monad
associated to the derived commutative operad as in \cite[Ex. 3.71]{bracanu}. This $\s$-functor simply sends $A$ to $\lim_{q}A^q_*$, where
the limit is taken in $LSym-Alg$ and $A^q_*$ is considered as a connective $LSym$-algebra and thus
as an object in $LSym-Alg$. This $\s$-functor is the right adjoint of an adjunction, and
preserves free objects over perfect complexes. It is likely that $Tot^\pi$ is an equivalence of $\s$-categories
when restricted to nice enough objects.}
\end{rmk}

\subsection{Spectrum of cosimplicial-simplicial rings}\label{ss1.2}

For a fixed commutative simplicial ring $k \in sCR$, we can 
work relatively over $k$ and define the model category $k-csCR$ of cosimplicial-simplicial 
commutative $k$-algebras, and its associated $\s$-category $k-\csCR$. As usual, we have a 
natural equivalence of $\s$-categories
$$k-\csCR \simeq k/\csCR,$$
where $k$ is considered as an object in $\csCR$ which is constant is the cosimplicial direction.

We consider $dAff_k:=(k-sCR)^{op}$, 
the model category of derived affine $k$-schemes, which is defined to be the opposite category of
that of simplicial commutative $k$-algebras. Remind from \cite{hagII} that it can be endowed
with the fpqc model topology, and that we can consider the model category of (hyper-complete) stacks
$dAff_k^{\sim,fpqc}$. There are here some set-theoretical issues, that can be solved,
as usual, by fixing two universes $\mathbb{U} \in \mathbb{V}$. By definition $k-sCR$ refers here
to $\mathbb{U}$-small simplicial $k$-algebras, and $dAff_k^{\sim,fpqc}$ is then the category of
functors $k-sCR \to sSet_\mathbb{V}$ to $\mathbb{V}$-small simplicial sets.

We then consider a functor
$$\Spec^\Delta : k-csCR_\mathbb{V}^{op} \longrightarrow dAff_k^{\sim,fpqc},$$
from $\mathbb{V}$-small cosimplicial-simplicial commutative $k$-algebras to $dAff_k^{\sim,fpqc}$. 
It is defined by 
sending $A \in k-csCR_\mathbb{V}$ to the functor $\Spec^\Delta\, A : sCR \to sSet_\mathbb{V}$ defined by 
$$\Spec^\Delta\, A : B \mapsto \uHom(A,B)= Hom_{csCR}(A,B^{\Delta^\bullet}),$$
where $\uHom$ are the simplicial $Hom$'s of the simplicial enrichement described in the 
previous section, and $B$ is considered as an object in $k-csCR_\mathbb{V}$ by viewing it as constant
in the cosimplicial direction.

\begin{prop}\label{p1}
The functor $\Spec^\Delta : k-csCR_\mathbb{V}^{op} \longrightarrow dAff_k^{\sim,fpqc}$
is right Quillen, and the induced functor on $\s$-categories
$$\RR\Spec^\Delta : k-\csCR_\mathbb{U} \longrightarrow \dSt_k$$
is fully faithful. The essential image of $\Spec^\Delta$ is the smallest
full sub-$\s$-category of $\dSt_k$ containing the objects
$K(\Ga,n)$ for various $n$, and which is stable by
 $\mathbb{U}$-small limits.
\end{prop}

\textit{Proof.} This is proven in a very similar manner than \cite[Cor. 2.2.3]{chaff}.
The left adjoint to $\Spec^\Delta$, denoted by $\OO$, sends 
a representable presheaf $h^B=Hom(B,-) : k-sCR \to sSet$ to 
$B \in k-csCR$, which is considered as constant in the cosimplicial direction. It also sends
an object of the form $K \times h^B$, for a simplicial set $K \in sSet$, to 
$B^K \in k-csCR$. Finally, it is uniquely defined by these properties together with the requirement that 
it sends colimits in $dAff_k^{\sim,fpqc}$ to limits in $k-csCR$. 

To prove that $\Spec^\Delta$ is right Quillen we use that $dAff_k^{\sim,fpqc}$ is a left Bousfield localization
of the levelwise projective model structure on the category of simplicial presheaves 
$Fun(k-sCR_\mathbb{U},sSet)$ obtained by inverting fpqc-local equivalences and 
equivalences in $sCR$ (see \cite{hagII}). Therefore, by general facts about Bousfield localizations,
it is enough to show that $\Spec^\Delta$ is a right Quillen functor for the
levelwise projective model structure, and moreover that for any cofibrant $A \in k-csCR$, 
$\Spec^\Delta\, (A)$ is a fibrant object in $dAff_k^{\sim,fpqc}$. The first of these statements 
easily follows from the fact that $k-csCR$ is a simplicial model category in which 
every object is fibrant, and thus
that if $A \to A'$ is a (trivial) cofibration in $k-csCR$, for any $C \in k-sCR$ the morphism
$\uHom(A',C) \to \uHom(A,C)$ 
is a (trivial) fibration of simplicial sets. For the second of these statements, let $A \in k-csCR$ be 
a cofibrant object. Any equivalence $B \to B'$ in $k-sCR$ obviously induces an equivalence
in $k-csCR$ when considered both $B$ and $B'$ as constant in the cosimplicial direction, and thus
$\uHom(A,B) \to \uHom(A,B')$ is an equivalence of simplicial sets. This shows that 
$\Spec^\Delta(A)$ is fibrant when the topology is the trivial topology. Finally, 
if $B=\mathrm{holim}_i B_i$ is a homotopy limit in $k-sCR$, it is also a homotopy limit in 
$k-csCR$, and thus the natural morphism $\uHom(A,B) \to \mathrm{holim}_i \uHom(A,B_i)$
is an equivalence. As the fpqc model topology is subcanonical, this implies that 
$\Spec^\Delta(A)$ is indeed a stack for the fpqc topology, and thus fibrant 
as an object in $dAff_k^{\sim,fpqc}$ (see \cite{hagII}).

To finish the proof of the Proposition, let $A \in k-csCR_\mathbb{U}$ and 
$X=\Spec^\Delta(A)$. By definition of the simplicial structure on $k-csCR$, the functor
$X$ is the geometric realization of the simplicial object $q \mapsto h^{A_*^q}$, where
$h^{A_*^q} : k-sCR \to sSet$ is corepresented by $A_*^q \in k-sCR$. In other words, we have 
$\Spec^\Delta(A) \simeq \mathrm{hocolim}_{q} h^{A^q_*}$. Therefore, the adjunction morphism
$$A \longrightarrow \OO(\Spec^\Delta(A))$$
is the canonical morphism $A \to \mathrm{holim}_{q}A^q_*$ in $Ho(k-csCR_\mathbb{U})$. This canonical 
morphism is obviously an equivalence, showing the fully faithfulness property in the Proposition.
It is then formal that the essential image of $\RR\Spec^\Delta$ is stable by limits. 
It also contains the objects $K(\Ga,n)$, as these are the images of $L\phi(\ZZ[-n])$, where
$L\phi : C(\ZZ) \to k-csCR$ is the left adjoint to $Tot^\pi$. Finally, 
as $k-csCR$ is cofibrantly generated, any object is equivalent to a $I$-cell object, where
$I$ is the image by $L\phi$ of the generating cofibrations in $C(\ZZ)$. The images of 
$I$ by $\Spec^\Delta$ are equivalent to morphisms of the form $* \to K(\Ga,n)$, so that 
 any object is the essential image of $\RR\Spec^\Delta$ lies in the smallest
sub-$\s$-category containing the $K(\Ga,n)$ and stable by limits.
\hfill $\Box$ \\

\begin{df}
The $\s$-category of \emph{($k$-linear) derived affine stacks} is the essential image of 
$\RR\Spec^\Delta : k-\csCR_\mathbb{U}^{op} \hookrightarrow \dSt_k$. It is denoted by 
$\dChAff_k$.
\end{df}

The above notions and results also have \emph{graded versions}, for which the details are left to the reader.
We denote by $k-csCR^{gr}$ the category of graded cosimplicial-simplicial commutative $k$-algebras. We endow
this category with the model category structure for which equivalences and fibrations
are defined by the graded Tot functor to graded complexes $Tot^\pi : k-csCR^{gr} \to C(\ZZ)^{gr}$, 
defined by taking the $Tot^\pi$ of each graded component. The proof of the existence of this model 
category structure follows the same lines as in the non-graded case. The corresponding 
$\s$-category will be denoted by $k-\csCR^{gr}$.

We consider the strict multiplicative group object $\Gm^{st}$ in the model category $ dAff_k^{\sim,fpqc}$ of fpqc (hyper-complete) stacks (see Section \ref{ss1.2}). $\Gm^{st}$ is defined as the
functor sending $B \in k-sCR$ to the group $B_0^*$ of invertible elements in the ring $B_0$
of $0$-simplices in $B$. This is a group object in the model category $dAff_k^{\sim,fpqc}$, which 
is not a fibrant object (it does not preservers equivalences in $B$), but its image
in $\dSt_k$ is the usual multiplicative group scheme $\Gm$. Indeed, 
$\Gm^{st}=h^{k[t,t^{-1}]}$ is corepresented by $k[t,t^{-1}]$, and,
as shown in \cite{hagII}, a fibrant model for $h^A$ is the derived scheme $\Spec\, A : B \mapsto Map(A,B)$.
In particular, we have that the model category 
$\Gm^{st}-dAff_k^{\sim,fpqc}$, of $\Gm^{st}$-equivariant objects in $dAff_k^{\sim,fpqc}$, is a model
for the $\s$-category $\dSt^{\Gm}_k$, of $\Gm$-equivariant derived stacks, called
\emph{graded derived $k$-stacks}.

The graded version of the spec functor is the functor
$$\Spec^{\Delta,gr} : (k-csCR^{gr})^{op} \longrightarrow \Gm^{st}-dAff_k^{\sim,fpqc},$$
sending $A \in k-csCR^{gr}$ to $\uHom(A,-)$, endowed with the natural
$\Gm^{st}$-action coming from the grading on $A$. More explicitly, 
for $B \in k-sCR$, the set of q-simplices of $\Spec^{\Delta,gr}(A)(B)$ is the set
of morphisms $Hom(A^q,B)$, which is endowed with a $B_0^*$-action as follows. 
For $b \in B_0^*$ and $f : A^q \to B$, we define $bf : A^q \to B$ by the formula
$(bf)(x)=b^n.f(x)$ for a homogeneous element $x$ of degree $n$. As in the non-graded
case, $\Spec^{\Delta,gr}$ is a right Quillen functor, and the induced $\s$-functor 
$$\RR\Spec^{\Delta,gr} : (k-\csCR^{gr}_\mathbb{U})^{op} \longrightarrow \dSt^{\Gm}_k$$
is fully faithful. Its essential image is the smallest sub-$\s$-category 
containing all the objects $K(\Ga^{\chi},n)$, for $n\geq 0$ and $\chi \in \ZZ$,  
where $\Gm$ acts on $\Ga$ with weight $\chi \in \ZZ$, 
and which is stable by $\mathbb{U}$-small limits. This is the graded analog of Proposition \ref{p1}.

\begin{df}\label{daffchgr}
The $\s$-category of \emph{graded derived affine $k$-stacks} is the essential image of 
$\RR\Spec^{\Delta,gr} : (k-\csCR_\mathbb{U}^{gr})^{op} \hookrightarrow \dSt^{\Gm}_k$. It is denoted by 
$\dChAff_k^{gr}$.
\end{df}

We finish this section by recalling the notion of \emph{linear derived stack}, already introduced and studied in \cite{mon}. 
We first notice that the Tot functor
$Tot^\pi : csMod_k \longrightarrow C(\ZZ)$ admits a natural lax monoidal structure. Therefore, 
it gives rise to an $\s$-functor
$$Tot^\pi : \csMod_k \longrightarrow \dg_k,$$
where $\csMod_k$ is the $\s$-category of cosimplicial-simplicial $k$-modules.
This $\s$-functor turns out to be an equivalence of $\s$-categories, and we denote by
$\phi : \dg_k \to \csMod_k$ its inverse.

For any object $E \in \dg_k,$ we consider
$\phi(E) \in \csMod_k$, and define
$$Sym^{\Delta}(E):=L\phi(E) \in csCR^{gr}$$
the free cosimplicial-simplicial commutative ring generated by $\phi(E)$ (here, $L\phi$ is the left adjoint to the $Tot^\pi$ functor, same notation as in the proof of Theorem \ref{tsome}). Being a free
commutative ring, $Sym^{\Delta}(E)$ comes naturally equipped with a graduation, and thus
is considered as an object in the $\s$-category $\csCR^{gr}$. 

\begin{df}\label{dlinearstack}
The \emph{$k$-linear derived stack associated to a complex $E \in \dg_k$}, 
is the graded derived affine $k$-stack defined by
$$\VV(E):=\mathbb{R}\Spec^{\Delta,gr}\,(Sym^{\Delta}(E)) \in \dSt^{\Gm}_k.$$
\end{df}

Note that by construction the functor of points of $\VV(E)$ sends $B \in k-sCR$
to $Map_{\dg_k}(E,N(B))$, where $N(B)$ is the normalisation of $B$ as a simplicial
$k$-module. The $\Gm$-action is then the natural action of invertible elements 
of $B$ on $B$ itself. 

It is proven in \cite{mon} that the $\s$-functor sending $E$ to $\VV(E)$ is fully faithful when
restricted to bounded above $k$-dg-modules $E$, and in particular for perfect
$k$-dg-modules $E$.

\end{appendix}

\bibliographystyle{alpha}
\bibliography{Biblio.bib}


\end{document}